\documentclass[div=12, reqno, twoside,fontsize=10.5pt]{scrbook}
\usepackage[T1]							{fontenc}
\usepackage[utf8]						{inputenc}
\usepackage							{ellipsis,textcomp,lmodern}
\usepackage							{amstext,amssymb,amsmath,amsthm}
\usepackage[fixamsmath,disallowspaces]				{mathtools}
\usepackage[shortlabels]					{enumitem}
\setlist[enumerate]{leftmargin=8mm ,labelindent=0cm}
\usepackage							{esint}

\usepackage[pdfpagelabels]{hyperref}
\hypersetup{bookmarksopen, pdfstartview={FitH}, pdfborder={0 0 1}, colorlinks=false, citecolor=blue, 
	linktoc = page}
\usepackage{fancyhdr}

\fancypagestyle{my-style}{ %
	\fancyheadoffset{0ex}
	\fancyhf{} 
	\fancyhead[RO]{\nouppercase{\rightmark}}
	\fancyhead[LE]{\nouppercase{\leftmark}}
	\fancyfoot[C]{\thepage}
}
\usepackage[subnum]{cases}
\theoremstyle{plain}
\newtheorem{theorem}		{Theorem}[chapter]
\newtheorem{proposition}	[theorem]	{Proposition}
\newtheorem{corollary}	[theorem]	{Corollary}
\newtheorem{lemma}	[theorem]	{Lemma}

\theoremstyle{definition}
\newtheorem{definition}	[theorem]	{Definition}
\newtheorem{remark}	[theorem]	{Remark}
\newtheorem{example}	[theorem]	{Example}
\newtheorem*{remark*}		{Remark}
\newtheorem*{example*}		{Example}
\DeclareMathOperator	{\R}		{\mathbb{R}}

\DeclareMathOperator	{\supp}		{supp}
\DeclareMathOperator 	{\dist}		{dist}

\newcommand		{\pv}		{\operatorname{p.v.}}

\renewcommand 		{\epsilon}	{\varepsilon}
\usepackage{xpatch}
\xpatchcmd{\proof}{\itshape}{\prooffont}{}{}

\newcommand{\prooffont}{\bfseries}

\usepackage{mathrsfs,dsfont}
\usepackage{makeidx}
\usepackage{frcursive}
\usepackage{graphicx}
\usepackage[toc, page]{appendix}
\usepackage{amssymb,amsmath,esint}
\usepackage[shortlabels]{enumitem}
\usepackage{amsthm}
\usepackage{amscd,mdwlist}
\usepackage[
paper=a4paper,
left=23mm,
right=18mm,
top=23mm,
bottom=22mm,
includefoot,
foot=8mm,
bindingoffset=0mm]{geometry}
\usepackage{relsize}
\usepackage{color}
\usepackage{bbm, mathrsfs,pgf,tikz}
 \usepackage{tcolorbox}
 \usepackage{thmtools}
 \usepackage[backgroundcolor=white, bordercolor=blue,
 linecolor=blue]{todonotes}
 \usetikzlibrary{arrows}
\usepackage{cleveref}

\numberwithin{equation}{chapter}

\newtheorem{counterexample}[theorem]{Counterexample}

\newtheorem*{theorem*}{Theorem}
\newcommand{\il}{\int\limits}
\newcommand{\iil}{\iint\limits}
\renewcommand{\d }{\,\mathrm{d}}

\newcommand{\WnuOm}{W_{\nu}^{p}(\Omega)}
\newcommand{\WnuOmR}{W_{\nu}^{p}(\Omega|\mathbb{R}^d)}
\newcommand{\WnuOmO}{W_{\nu,\Omega}^{p}(\Omega|\mathbb{R}^d)}
\newcommand{\VnuOmO}{V_{\nu}^{\Omega}(\Omega|\mathbb{R}^d)}
\newcommand{\VnuOm}{V_{\nu}(\Omega|\mathbb{R}^d)}
\newcommand{\VnuOma}{V_{\nu_\alpha}(\Omega|\mathbb{R}^d)}

\newcommand{\HnuOm}{H_{\nu}(\Omega)}
\newcommand{\HnuOma}{H_{\nu_\alpha}(\Omega)}
\newcommand{\TnuOm}{T_{\nu} (\Omega^c)}

\newcommand{\eps}{\varepsilon}
\newcommand{\bet}[1]{\left| #1 \right|}

\newcommand{\vertiii}[1]{{\left\vert\kern-0.25ex\left\vert\kern-0.25ex\left\vert #1 \right\vert\kern-0.25ex\right\vert\kern-0.25ex\right\vert}}

\usepackage{color}
\definecolor{GF}{rgb}{0.1,0.5,0.1}
\definecolor{M}{rgb}{0.1,0.2,0.7}
\definecolor{H}{rgb}{0.7,0.1,0.2}
\definecolor{Y}{rgb}{0.8,0.35,0.1}
\definecolor{MY}{rgb}{0.5,0,0.45}

\begin{document}
	\title{
\begin{Huge} $L^2$-Theory For Nonlocal Operators On Domains \end{Huge}\\
\url{https://doi.org/10.4119/unibi/2946033}
\vspace*{3cm} \\
\begin{Large}	
	\vspace*{0.5cm}
	Universität Bielefeld \\
	\vspace*{0.3em}
	Fakultät für Mathematik\\
\end{Large}	{\ }\\
\vspace*{2cm}
\begin{huge}Dissertation	            \end{huge}\\
\vspace*{0.7cm}
\begin{Large}zur Erlangung des akademischen Grades \\ {\ }  \vspace*{-0.7em} \\
	\textsc{Doktor der Mathematik (Dr. math.)}\end{Large} \\
\vfill
\begin{Large}	eingereicht von \\
	\vspace*{0.3cm}
	M.Sc. Foghem Gounoue Guy Fabrice\\
	\vspace*{0.3cm}
	\date{\textbf{am 2. Juni 2020}}
\end{Large}	
}
\author{}


	\maketitle
	
	\thispagestyle{empty}
	\tableofcontents
	
	

\newpage

\thispagestyle{empty}
\vspace*{\fill}

	\begin{quote}
		
	\begin{center}\cursive{	\textbf{\textit{If I had only one hour to solve a problem, I would spend fifty-five minutes defining the problem, and only five minutes thinking about solutions.}}}
		
		\vspace{0.5cm}
		Albert Einstein.\footnote{This quote and several other versions of it are commonly attributed to Albert Einstein, but there seems to be no evidence witnessing his ownership. It is not incorporated in the comprehensive collection “The Ultimate Quotable Einstein” from Princeton University Press".}
	\end{center}
	\end{quote}
\vspace*{\fill}

	
\chapter{Introduction}	
\vspace{-3.5ex}

\subsection*{\textbf{Motivation and main goals}} 
\vspace{-1ex}

In recent years, the study of \textit{nonlocal} operators and related fields has been of significant interest to many researchers in the areas of both analysis and probability theory. This thesis is devoted to developing an accessible and moderate level of $L^2$-theory for symmetric \textit{nonlocal} operators on bounded domains. Our approach will closely follow the standard $L^2$-theory for \textit{local} operators, especially that of  elliptic partial differential operators of second order. Formally, an operator $A$ defined on a space of  functions is  defined as \textit{local} if it preserves 
the support, i.e. it obeys the rule, $\supp (Au) \subset \supp (u)$ for every $u$ in the domain of $A$. Otherwise, $A$ is classified  as a \textit{nonlocal} operator. Among local operators, we have the Laplace operator $A=\Delta $, the gradient operator $A = \nabla $ and the divergence operator $A=\operatorname{div};$ each is defined on the space of smooth functions. In terms of  nonlocal operators, in this thesis, we focus our attention on a subclass of integrodifferential operators, i.e. operators acting on a smooth function $u\in C_c^\infty(\R^d)$ as follows
\vspace{-1.3ex}
\begin{align*}
\mathscr{L} u(x)= \pv\int_{\R^d} (u(x)-u(y))\, \mu(x, \d y),\quad \quad (x\in \R^d).
\end{align*} 

\vspace{-1.2ex}

\noindent Here, $(\mu(x,\d y))_{x\in \R^d}$ is a family of Borel measures on $\R^d$ also called jump interaction measures such that $\mu(x, \{x\})=0$ for all $x\in \R^d$. It is noted that the pointwise evaluation of $\mathscr{L} u(x)$ may fail to hold even for a \textit{bona fide} function $u\in C_c^\infty (\R^d)$. Therefore, it is often convenient to evaluate $\mathscr{L} u$ in the generalized sense, i.e.  in the sense of distributions or via the associated energy form. A similar observation occurs in the local setting for a partial differential operator of the form $\mathscr{A}=-\operatorname{div}(A\cdot \nabla)$ where $A:\R^d\to \R^{d^2}$ with $A(x) =(a_{ij}(x))_{1\leq i,j\leq d}$ is a matrix-valued measurable function. $\mathscr{A}$ can be as good or as bad as $\mathscr{L}$.

\vspace{1mm}

\noindent The nonlocal feature of integrodifferential operators requires that elliptic conditions of the associated nonlocal problems on a domain must be assigned on the whole complement. The terminology \textit{nonlocal complement value problems} is thus appropriate. 
This contrasts with the local situation (elliptic boundary value problems) where the data are usually prescribed on the boundary of associated domains.

\vspace{0.5mm}

 \noindent The majority of our research concerns the study of IntegroDifferential Equations (IDEs) on open bounded domains subject to Neumann, Dirichlet, Robin and mixed \textit{complement} type conditions.
 We will not deal with the physical interpretations of these problems and only do mathematics
\emph{l’art pour l’art }\footnote{Meaning that mathematics is applied to  mathematics.}. For physical applications, many references are quoted in \cite{Va09}.   The solvability of integrodifferential equations on bounded domains raises a natural conceptual problem regarding the function spaces related to such operators. Thus, the overreaching goal of this work is to set up some milestones for a method of dealing with the aforementioned complement value problems in the framework of Hilbert spaces. 

\vspace{1mm}

\noindent The main  motivations of this work  are threefold: (i) to set-up convenient function spaces encoding certain nonlocal problems, (ii) to study the well-posedness of the aforementioned problems and their applications within the frameworks of Hilbert spaces and (iii) to  explore the asymptotic behavior of solutions associated with integrodifferential
operators with collapsing jump kernels. The latter objective aims to bridge a transition from nonlocal objects to the corresponding local objects. For instance, we show that limits of solutions to elliptic IDEs are solutions to elliptic partial differential equations (PDEs) of second order. 

\vspace{1mm} 

\noindent Nonlocal Dirichlet problems and related topics have received a great deal of  research attention over the last few years. See  for instance \cite{FKV15} where the well-posedness of the variational formulation is covered for a large class of symmetric and nonsymmetric jump interaction kernels. We shall mention some supplementary references later on. In the meantime, we note that nonlocal Robin and mixed complement value problems appear to be combinations of Dirichlet and Neumann complement value problems. Accordingly, our approach within this work is designed to give much more attention and priority to nonlocal Neumann problems.  For illustration purposes,  we shall illuminate our exposition by briefly reviewing the Neumann problem for the 
simple case of the Laplace operator. This summary is intended to help the reader to become 
more familiar with the corresponding nonlocal formulation.  


\subsection*{\textbf{Illustration in the local case}}Let $\Omega\subset \mathbb{R}^d $ be a bounded open subset whose boundary $\partial\Omega$ is sufficiently regular. Let $f:\Omega\to\mathbb{R}$ and $g: \partial\Omega \to \mathbb{R}$ be measurable. The classical inhomogeneous Neumann problem for the Laplace operator associated with data $f$ and $g$ consists of finding a function $u:\Omega\to \mathbb{R}$ satisfying the following:
\begin{align}\label{eq:local-Neumann}
-\Delta u = f \quad\text{in}~~~ \Omega \quad\text{ and } \quad \frac{\partial u}{\partial n}= g ~~~ \text{on}~~~ \partial \Omega
\end{align}
where  $\frac{\partial u}{\partial n} $ stands for the outward normal derivative of $u$ on $\partial\Omega$. 
To be more precise if we assume $u$ is the restriction on $\overline{\Omega}$ of a smooth function $\widetilde{u}:\mathbb{R}^d \to \mathbb{R}$ then $\frac{\partial u}{\partial n} (x)= \nabla \widetilde{u}(x)\cdot n(x)$ where $n(x) $ is the outer normal vector at $x\in \partial\Omega$. It is a common approach to consider $ f\in L^2(\Omega)$ and $g\in H^{1/2}(\partial\Omega)$, where $H^{1/2}(\partial\Omega)$ is the trace space of the classical Sobolev space $H^1(\Omega)$. However, for the variational formulation, the regularity requirement on $f$ and $g$ can be relaxed. For example, it is possible to consider $f\in (H^1(\Omega))'$ (dual of $H^1(\Omega)$) and $g\in H^{-1/2}(\partial \Omega)$ (dual of $H^{1/2}(\partial \Omega)$). For the sake of simplicity, assume that $f\in L^2(\Omega)$ and $g\in L^2(\partial \Omega)$. A function $u: \Omega\to \R$ is said to be a weak solution, or a variational solution, to \eqref{eq:local-Neumann} if $u\in H^1(\Omega)$,  such that:
\begin{align}\label{eq:local-Neumann-var} 
	\int_\Omega \nabla u(x)\cdot \nabla v(x) \d x =\int_\Omega f(x) v(x) \d x+ \int_{\partial \Omega} g(x)\gamma_0   v(x) \d \sigma(x)\quad\text{for all $v\in H^1(\Omega)$} 
\end{align}
\noindent where, $\gamma_0 v$  stands for the trace of $v$ on $H^1(\Omega)$ and $\d\sigma$ stands for the 
restriction of the Lebesgue measure on $\partial\Omega$. It is worth  emphasizing that the trace space $ H^{1/2}(\partial\Omega) $ can be seen as the space  that models the restriction\footnote{The restriction here is 
understood in the sense of traces. Indeed, it does not make sense to restrict a measurable function on $\partial\Omega$ since $\partial\Omega$ is the smooth manifold of a dimension lower that than $d$, and thus  has measure zero.} to the boundary $\partial\Omega$ of functions in 
$H^1(\Omega)$. Most importantly, the study of trace spaces is mostly motivated by the study of problems such as 
boundary value problems with Neumann condition. If $u$, $f$ and $g$ are sufficiently regular, 
say $u\in H^2(\Omega)$, $ f\in L^2(\Omega)$ and $g\in H^{1/2}(\partial\Omega)$ then  \eqref{eq:local-Neumann} and \eqref{eq:local-Neumann-var} are equivalent. This can be accomplished through 
applying the classical Green-Gauss formula (see \cite[Chapter 3]{Necas67}, \cite[Theorem 2.20]{FSU19} or \cite[Appendix A.3]{Trie92}), 
\begin{align}\label{eq:local-green-gauss} 
\int_\Omega \hspace{-0.2ex} [-\Delta u(x)]v(x) \d x= \int_\Omega \hspace{-1ex}\nabla u(x)\cdot \nabla v(x) \d x + \int_{\partial \Omega} \hspace{-2ex}\gamma_1(x)\gamma_0v(x) \d \sigma(x),\quad\text{for $u\in H^2(\Omega),\,\, v\in H^1(\Omega)$}. 
\end{align}
\noindent Note that $\gamma_1 u= \gamma_0 \frac{\partial u}{\partial n}$.  In fact, the Green-Gauss formula \eqref{eq:local-green-gauss} is the keystone to the weak formulation \eqref{eq:local-Neumann-var}. We will return later to some observations related to problems $\eqref{eq:local-Neumann}-\eqref{eq:local-Neumann-var}$. Meanwhile, it is apparent that, in many studies of PDEs, Neumann boundary problems receive less attention than Dirichlet boundary problems. Notwithstanding, some aspects of Neumann boundary problems involving elliptic 
differential operators of second order are studied in \cite{Dorothee-Triebel, Taylor, Jost-Jurgen, DjVj09, Mikhailov78}. 
A rigorous treatment on topics including regularity up to the boundary, Schauder estimates, $L^p$ estimates and variational formulation related to the Neumann boundary problem for the Laplace operator can be found in Giovanni Leoni’s detailed, thorough lecture notes \cite{Giovanni13}. See also \cite{Dan00} where several topics are extended to the local Robin boundary value problems. 

\subsection*{Illustration in the nonlocal case}

\noindent One of the goals of this work is to set up the inhomogeneous Neumann problems for a certain class of symmetric nonlocal operators. We focus our attention on symmetric nonlocal L\'{e}vy operators which are integrodifferential operators of the form
\begin{align*}
L u(x)&= \pv \int_{\mathbb{R}^d}(u(x)-u(y))\nu(x-y)\mathrm{d}y,\quad (x\in \mathbb{R}^d),
\end{align*}
defined for a given measurable function $u:\mathbb{R}^d \to \mathbb{R}$ whenever the right hand side exists and makes sense. Here and henceforward, the function $\nu:\mathbb{R}^d\setminus\{0\}\to[0,\infty]$ is the density of a symmetric L\'{e}vy measure. In other words, $\nu$ is positive and measurable such that
\begin{align*}
\nu(-h) = \nu(h) ~\text{ for all $h\in \mathbb{R}^d$ and} ~\int_{\mathbb{R}^d}(1\land |h|^2)\nu(h)\d h<\infty\,.
\end{align*}
Notationally, for $a,b\in \R$ we write $a\land b$ to denote $\min(a,b)$. In addition, we will also assume that $\nu$ does not vanish on sets of positive measure. In short we say that $\nu$ has full support. If $\nu$ is radial, we write $\nu$ for the profile, too, i.e. $\nu(r) =\nu(x)$ if $|x|=r$. For $Lu(x)$ to be defined, it is sufficient for $u$ to possess a $C^2$ regularity in the neighborhood of the point $x$ and some weighted integrability for $|x|\to \infty$.


\noindent A prototypical example of an operator  $L$ is obtained by taking, $\nu(h) =C_{d,\alpha} |h|^{-d-\alpha}$ for $h\neq 0$ where $ \alpha \in (0,2)$ is fixed. The resulting operator is the so called fractional Laplace operator $(-\Delta)^{\alpha/2}$. The constant $C_{d, \alpha}$ given by
\begin{align*}
C_{d,\alpha}:= \left(\int_{\mathbb{R}^d} \frac{1-\cos(h_1)}{|h|^{d+\alpha}} \d h\right)^{-1},
\end{align*}
is chosen so that the Fourier relation $\widehat{(-\Delta)^{\alpha/2} u}(\xi)= |\xi|^\alpha \widehat{u}(\xi)$ holds for all $u\in C^\infty_c(\mathbb{R}^d)$. The fractional Laplacian $(-\Delta)^{\alpha/2}$ is one of the most heavily studied  integrodifferential operators. Asymptotically, we have  $C_{d, \alpha}\asymp \alpha (2-\alpha)$. This will play an important role for our analysis. Further details about the fractional Laplacian $(-\Delta)^{\alpha/2}$ and the constant $C_{d,\alpha}$ are presented in Chapter \ref{chap:basics-nonlocal-operator}. Next, we introduce the Neumann problem associated with the operator $L$. 

\medskip

\noindent Assume $\Omega\subset \mathbb{R}^d$ is an open set. Let $f:\Omega\to \mathbb{R}$ and $g: \mathbb{R}^d\setminus \Omega\to \mathbb{R}$ be measurable. The Neumann problem for the operator $L$ associated with data $f$ and $g$ is to find a measurable function $u:\mathbb{R}^d\to \mathbb{R}$ such that 
\begin{align}\label{eq:nonlocal-Neumann-intro}\tag{$N$}
L u = f \quad\text{in}~~~ \Omega \qquad\text{ and } \qquad \mathcal{N} u= g ~~~ \text{on}~~~ \mathbb{R}^d\setminus\Omega.
\end{align} 
Here, $\mathcal{N}$ is another the integrodifferential operator defined by 
\begin{align}\label{eq:nonlocal-derivative-intro}
\mathcal{N}u(y) = \int_{\Omega}(u(x)-u(y))\nu(x-y)\,\mathrm{d}x,\qquad (y\in \mathbb{R}^d\setminus\Omega). 
\end{align}
 and is also called  for obvious reasons the \emph{nonlocal normal derivative operator} across the boundary of $\Omega$ with respect to $\nu$. We  justify this terminology  progressively through our exposition. The problem \eqref{eq:nonlocal-Neumann-intro} is said to be homogeneous if $g=0$ and inhomogeneous otherwise. Let us emphasize that, in contrast to the local situation 
where $\Delta u(x) $ is evaluated via an arbitrarily small neighborhood of $x\in \Omega$,
the eventual evaluation of $L u(x) $ requires us to know $u(y)$ for almost all $y\in \mathbb{R}^d$. It is therefore reasonable to prescribe \textit{the Neumann condition} on the complement of $\Omega$. The terminology "\textit{nonlocal Neumann complement value problem}" is thus appropriate for the system \eqref{eq:nonlocal-Neumann-intro}.  It is also important to observe  that unlike in the local setting, the definition of the nonlocal normal derivative $\mathcal{N} u$ requires neither the regularity of the function $u$ nor the regularity of the domain $\Omega.$ Another observation is that the definition of $\mathcal{N}$  implies that if $u$ solves \eqref{eq:nonlocal-Neumann-intro},  then for almost all $y\in \mathbb{R}^d\setminus \Omega$ we have 
\begin{align*}
u(y)= n^{-1}_\nu(y) \Big(-g(y)+ \int_{\Omega}u(x) \nu(x-y)\mathrm{d}x\Big)\,\quad\text{with }\quad n_\nu(y) = \int_{\Omega}\nu(x-y)\d x.
\end{align*}
In other words, if $u$ solves \eqref{eq:nonlocal-Neumann-intro} then the values of $u$ on $\mathbb{R}^d\setminus\Omega$ solely depend on its values inside $\Omega$ and the function $g$. Meanwhile,
it can be observed that the integrodifferential operators $L$ and $\mathcal{N}$ are closely related from their common integrands and differ solely with regard to the domain of integration. Furthermore, with the nonlocal 
analog of the normal derivative  at hand, it makes sense to derive a formula that resembles the classical Green-Gauss for formula \eqref{eq:local-green-gauss}.  Indeed, a routine check shows that the operators $L$ and $\mathcal{N}$ satisfy the following \textit{nonlocal Green-Gauss formula}
\begin{align}\label{eq:green-gauss-nonlocal-intro}
\int_{\Omega} [Lu(x) ]v(x)\d x= \mathcal{E}(u,v)- \int_{\Omega^c} \mathcal{N} u(y)v(y)\d y,\qquad\text{for all} \quad u,v \in C_c^\infty(\mathbb{R}^d)\,
\end{align}
where, here and henceforth, $\mathcal{E}(\cdot,\cdot)$ stands for the bilinear form, 
\begin{align*}
\mathcal{E}(u,v)&:= \frac{1}{2}\iint\limits_{(\Omega^c\times \Omega^c)^c} (u(x)-u(y))(v(x)-v(y))\nu(x-y)\d x\d y.
\end{align*}
The Green-Gauss formula \eqref{eq:green-gauss-nonlocal-intro} is the cornerstone for deriving the variational formulation of the Neumann problem \eqref{eq:nonlocal-Neumann-intro}. To do this, we need to introduce adequate function spaces. 

\vspace{1mm}

\noindent To the best of our knowledge,  the problem \eqref{eq:nonlocal-Neumann-intro} and the operator $\mathcal{N}$ were introduced for the first time in \cite{DROV17} for the fractional Laplace operator. They successfully analysed the well-posedness of the homogeneous problem by applying the Fredholm alternative. Therein, the inhomogeneous problem is analysed on a function space that depends on the Neumann complement data $g$, which is somewhat unaccustomed. 

\noindent We will propose a slightly different framework that does not  depend on the Neumann data. In fact, our function spaces are adapted for certain types of nonlocal problems similar to \eqref{eq:nonlocal-Neumann-intro} such as the Robin complement value problem (see \eqref{eq:nonlocal-Robin-intro}).  As mentioned above, an important role in our study is played by function spaces.

\subsection*{Energy spaces and nonlocal trace spaces}  Let us introduce several function spaces with respect to $\nu$ that are fundamental for our work. 

\medskip 

\noindent $\bullet$ We define the space $\HnuOm$ by, $\HnuOm= \big\{u \in L^2(\Omega): \, |u|^2_{\HnuOm}<\infty \big \}$ equipped with the norm $\|u\|^2_{ \HnuOm}= \|u\|^2_{L^2(\Omega)} + |u|^2_{ \HnuOm}$ where 
%
%
\begin{align*}
|u|^2_{ \HnuOm}&=\iil_{\Omega\Omega} \big(u(x)-u(y) \big)^2 \, \nu (x-y) \d y\,\d x\\
&= \iil_{\R^d\R^d} \!\!\big(u(x)-u(y) \big)^2 \, \min(\mathds{1}_\Omega (x), \mathds{1}_\Omega(y))\, \nu (x-y) \d y\,\d x.
\end{align*}

\noindent $\bullet$ Following \cite{FKV15,SV14}, we introduce the space  $\VnuOm$ by 
$$ \VnuOm = \big\{ u: \R^d \to \R \text{ meas.}:\, \, u|_{\Omega}\in L^2(\Omega)\,\,\text{and}\,\, \mathcal{E}(u,u)<\infty \big\}.$$
 Recall  that the form $\mathcal{E}(\cdot,\cdot) $ is simultaneously given by 
%
%

\begin{align*}
\mathcal{E}(u,u) &=\frac{1}{2} \hspace{-1ex}\!\!\iil_{(\Omega^c\times \Omega^c)^c} \hspace{-1ex}\!\!\big(u(x)-u(y) \big)^2 \, \nu (x-y) \d y\,\d x= \frac{1}{2} \!\!\iil_{\R^d\R^d} \!\!\big(u(x)-u(y) \big)^2 \,\max(\mathds{1}_\Omega (x), \mathds{1}_\Omega(y))\,\nu (x-y) \d y\,\d x.
\end{align*}

\noindent Keep in mind  that $(\Omega^c\times \Omega^c)^c=(\mathbb{R}^d\times \mathbb{R}^d)\setminus (\Omega^c\times \Omega^c)= (\Omega\times \Omega) \cup( \Omega^c\times \Omega)\cup( \Omega\times \Omega^c)$. It is worth emphasizing  that in the integrand of $\mathcal{E}(\cdot,\cdot)$, only the increments from $\Omega^c$ to $\Omega^c$ are not allowed. Hence a function $u\in \VnuOm$ has certain regularity inside $\Omega$ and across the boundary $\partial\Omega$. Furthermore, the space $\VnuOm$ can be redefined as 
$$\VnuOm = \big\{ u: \R^d \to \R \text{ meas.}: \, \, u|_{\Omega}\in L^2(\Omega)\,\,\text{and}\,\, |u|^2_{\VnuOm}<\infty \big\}$$ where we have, 
\begin{align*}
|u|^2_{\VnuOm} &:=\!\!\iil_{\Omega\R^d} \!\!\big(u(x)-u(y) \big)^2 \, \nu (x-y) \d y\,\d x= \frac{1}{2} \!\!\iil_{\R^d\R^d} \!\!\big(u(x)-u(y) \big)^2 \,[\mathds{1}_\Omega (x)+ \mathds{1}_\Omega(y)]\,\nu (x-y) \d y\,\d x.
\end{align*}
\noindent Indeed, we have $\mathcal{E}(\cdot,\cdot)\leq |\cdot|^2_{\VnuOm}\leq 2\mathcal{E}(\cdot,\cdot)$ as for all $x,y\in \R^d$,
$$\max(\mathds{1}_\Omega (x), \mathds{1}_\Omega(y))\leq \mathds{1}_\Omega (x)+ \mathds{1}_\Omega(y)\leq 2\max(\mathds{1}_\Omega (x), \mathds{1}_\Omega(y)).$$

\noindent Throughout, we equip the space $\VnuOm$ with the norm
\begin{align*}
\|u\|^2_{\VnuOm}=\|u\|^2_{L^2(\Omega)}+|u|^2_{\VnuOm}\asymp \|u\|^2_{L^2(\Omega)}+\mathcal{E}(u,u). 
\end{align*}

\vspace{1mm}

\noindent $\bullet$ Assume  $L$  is the L\'evy operator associated with the measure $\nu$. Define the new space  $V^1_\nu(\Omega|\R^d)$ by
\begin{align*}
V^1_\nu(\Omega|\R^d)= \Big\lbrace u\in \VnuOm:\, \text{ $Lu$ exists weakly and } Lu\in L^2(\Omega)\,\Big\rbrace \,.
\end{align*}
Here, the weak integrodifferentiability of $Lu$ is understood in the sense of Definition \ref{def:weak-integro}.

\vspace{2mm}

\noindent $\bullet$ Let us also introduce the  subspace of functions in $\VnuOm$ that vanish on the complement of $\Omega$, i.e.
\begin{align*}
\VnuOmO= \{ u\in \VnuOm~: ~u=0~~\text{a.e. on } \mathbb{R}^d\setminus \Omega\}\,.
\end{align*} 

\noindent The space $\VnuOmO$ is clearly a closed subspace of $\VnuOm$. The space $ \VnuOmO$ encodes the nonlocal Dirichlet problems related to the operator $L$. 

\medskip

\noindent Note that if $\nu(h) = |h|^{-d-\alpha}$ for $h \in \R^d, h \ne 0$ with $\alpha \in (0,2)$, the space $ \HnuOm$ equals the classical fractional Sobolev-Slobodeckij space $H^{\alpha/2}(\Omega)$. For the same choice of $\nu$, we define $V^{\alpha/2} (\Omega|\mathbb{R}^d)$ as the space $ \VnuOm$.  Next, we need to introduce the spaces of functions defined on the complement of $\Omega$ which aim to incorporate the complement data. 

\medskip

\noindent $\bullet$ We define trace space $\TnuOm$: the space of restrictions to $\mathbb{R}^d\setminus \Omega$ of functions of $\VnuOm$. That is, 
\begin{align*}
\TnuOm = \{v: \Omega^c\to \mathbb{R}~\text{meas.} ~~\hbox{ such that }~~ v = u|_{\Omega^c} ~~\hbox{with }~~ u \in \VnuOm\}.
\end{align*}
We equip $\TnuOm $ with its natural norm, 
\begin{align*}
\|v\|_{\TnuOm } = \inf\{ \|u\|_{\VnuOm }: ~ u \in \VnuOm ~ \hbox{ with }~ v = u|_{\Omega^c} \}. 
\end{align*}

\noindent $\bullet$ We consider the weighted $L^2$-spaces on $\Omega^c$ denoted by $L^2(\Omega^c, \nu_K)$ and $L^2(\Omega^c, \mathring{\nu}_K)$ where for a given measurable set  $K\subset \Omega$ with 
$0<|K|<\infty$, we define 
\begin{align}
\nu_K(x)&:=\underset{y\in K}{\operatorname{essinf}} ~ \nu(x-y)\quad\text{and}\quad
\mathring{\nu}_K(x):=\int_K 1\land \nu(x-y)\d y.
\end{align}

\medskip

\noindent The aforementioned spaces are Hilbert spaces and are  connected. Indeed, the following hold 
\begin{enumerate}[$(a)$]
	\item The embeddings $ \VnuOmO \hookrightarrow \VnuOm \hookrightarrow\TnuOm \hookrightarrow L^2(\Omega^c, \nu_K)$ are continuous.
	\item The trace operator $\operatorname{Tr}: \VnuOm \to L^2(\Omega^c, \nu_K)$ with $u \mapsto \operatorname{Tr}(u) = u\mid_{\Omega^c}$ is linear and continuous.
	\item We have $\operatorname{Tr}(\VnuOm)= \TnuOm$ and $ \ker(\operatorname{Tr} ) =  \VnuOmO$. 
\end{enumerate}

\noindent  The same holds true with $\nu_K$ replaced by $\mathring{\nu}_K$ . These interactions between the spaces $\VnuOm$, $\TnuOm$, $  \VnuOmO$ and $L^2(\Omega^c, \nu_K)$ respectively, are analogous to the ones between the classical Sobolev spaces $H^1(\Omega)$, $H^{1/2}(\partial\Omega)$, $H^1_0(\Omega)$ and $L^2(\partial\Omega)$. Recall  that $H^1_0(\Omega)$ is by definition the closure of $C_c^\infty(\Omega)$ in $H^1(\Omega)$.  In fact, it is well known that the following are true. 

\begin{enumerate}[(a')]
	\item The embeddings $ H_0^1(\Omega)\hookrightarrow H^1(\Omega)\hookrightarrow H^{1/2}(\partial\Omega) \hookrightarrow L^2(\partial\Omega)$ are continuous.
	\item The  classical trace operator $\gamma_0 :H^1(\Omega)\to L^2(\partial\Omega)$ whenever it exists, is linear and continuous. 
\item We have  $\gamma_0(H^1(\Omega))= H^{1/2}(\partial\Omega)$ and $\ker(\gamma_0)= H^1_0(\Omega)$.
\end{enumerate}

\noindent  In this regard, it is fair to view the spaces $\VnuOm$, $  \VnuOmO$, $\TnuOm$ and   
$L^2(\Omega^c, \nu_K)$ as the nonlocal replacements of $H^1(\Omega)$, $H^1_0(\Omega),$ 
$H^{1/2}(\partial\Omega)$ and $L^2(\partial\Omega)$, respectively. 
Therefore, it is natural to think of the space $\TnuOm$ as the space encoding Dirichlet complement data
whereas the space $L^2(\Omega^c, \nu_K)$ or $L^2(\Omega^c, \mathring{\nu}_K)$  encodes Neumann complement data. Let us mention  that  for sufficiently smooth $\Omega$,  spaces similar to $\TnuOm$ have been  recently studied in \cite{BGPR17, DyKa18} wherein, the main motivation is the study of the nonlocal extension problem for the 
space $\VnuOm$ which is an analog of the classical Sobolev extension problem from the space $H^{1/2}(\partial \Omega)$ to $H^1(\Omega)$. We discuss and compare the aforementioned spaces and their 
connections with the classical Sobolev spaces in  Section \ref{sec:function-spaces}.  Another approach to nonlocal trace spaces on constrained domains is considered in \cite{TD17}. 

\medskip

\noindent Let us say a few words about the weights $\nu_K$ and $\mathring{\nu}_K$.  To avoid $ L^2(\Omega^c, \nu_K)$ being  a trivial space, the following condition on $\nu_K$ is implicitly required
\begin{align}\label{eq:non-degeneracy-condition}
\nu_K(x):=\underset{y\in K}{\operatorname{essinf}} ~ \nu(x-y)>0,\quad\text{ for almost every $x\in \Omega^c$.}
\end{align}

\noindent The condition \eqref{eq:non-degeneracy-condition} also ensures that the function $\nu$ does not 
decay (degenerate) too much for large increments across the boundary of $\Omega$. A suitable way to picture 
this condition is to assume that the function $h\mapsto\nu(h)$ is continuous on $\mathbb{R}^d\setminus\{0\}$. Choose $K$ to be a nonempty compact subset of $\Omega$ and fix $x\in\Omega^c $ so 
that $\nu_K(x)=\nu(x-z)>0$ for a suitable $z\in K$. The same observation holds for $\mathring{\nu}_K$. Another crucial  observation is that the weights $\nu_K$ and $\mathring{\nu}_K$  annihilate the eventual singularity of 
$\nu$ at the origin. Let us illuminate our argument by looking at a particular case of $\nu$. Assume that $\nu$ is 
a unimodal, i.e. $\nu$ is radial and almost decreasing in the sense that there exists a constant $c>0$ such that 
$c\nu(y)\leq \nu(x)$ whenever $|y|\geq |x|$. In addition, assume $\nu$ satisfies the doubling growth condition i.e. 
$\exists\,\,\kappa >0 $ such that 
$$\nu(h)\leq \kappa \nu(2h) \quad \text{ for all $|h|\geq 1$}.$$ 
 Let $B_1$ be the unit ball of  $\R^d$ then  (see Theorem \ref{thm:Vnu-in-l2}) we have $\mathring{\nu}_{B_1}, \nu_{B_1}\in  L^1(\R^d)$, $\nu_{B_1}\asymp 1\land \nu$ and $\mathring{\nu}_{B_1}\asymp 1\land \nu.$ Clearly, $1\land \nu $ does not have a singularity at the origin. For a concrete example, one could consider the standard example $\nu(h) = |h|^{-d-\alpha}, ~\alpha\in (0,2)$. Then we have $\nu_{B_1}(x) \asymp (1+|x|)^{-d-\alpha} \asymp 1\land \nu(x)$ and $\mathring{\nu}_{B_1}(x) \asymp (1+|x|)^{-d-\alpha} \asymp 1\land \nu(x)$. 

\vspace{1mm}
\noindent Although the class of almost decreasing unimodal L\'evy kernel is fairly large, there exist some radial L\'evy kernels that are not almost decreasing . For example, for $\beta\in [-1, 2)$ define 
\begin{align*}
\nu_\beta(h) = |h|^{-d-\beta}\,\,\rho(h)\text{ with }\,\, \rho(h) :=\Big(\frac{2+\cos |h|}{3}\Big)^{|h|^4}.
\end{align*}
Note that $\nu_\beta$ is not almost decreasing since  $\rho(2\pi n) =1$ and $\rho(\pi (2n-1)) =3^{-\pi^4 (2n-1)^4}$ for all $n\in \mathbb{N}$. Now if $\beta\in (0,2)$ it is clear that $\nu_\beta$ is L\'evy integrable  since $\rho$ is bounded. If $-1\leq \beta \leq 0$ then $\nu_\beta$  is also L\'evy integrable since the map  $r\mapsto \rho(r)$ is in  $L^1(\R)$.  

\medskip
\subsection*{Variational formulation of the Neumann problem }
\noindent Let us now return to our main problem of interest. In light of the relation \eqref{eq:green-gauss-nonlocal-intro} it is reasonable to define weak solutions of the Neumann problem under consideration as follows. A measurable function $u:\mathbb{R}^d\to \mathbb{R}$ is a weak solution or a variational solution of the inhomogeneous Neumann problem \eqref{eq:nonlocal-Neumann-intro} if $u \in \VnuOm$ and satisfies the relation
	\begin{align}\label{eq:var-nonlocal-Neumann-intro}\tag{$V$}
	\mathcal{E}(u,v) = \int_{\Omega} f(x)v(x)\mathrm{d}x +\int_{\Omega^c} g(y)v(y)\mathrm{d}y,\quad \mbox{for all}~~v \in \VnuOm\,.
	\end{align}

\noindent In particular if $ \Omega$ is bounded, taking $v=1$, \eqref{eq:var-nonlocal-Neumann-intro} becomes the so called compatibility condition
	\begin{align}\label{eq:compatible-nonlocal-intro}\tag{$C$}
	\int_{\Omega} f(x)\mathrm{d}x +\int_{\Omega^c} g(y)\mathrm{d}y=0.
	\end{align}

\medskip

\noindent In fact, the compatibility condition \eqref{eq:compatible-nonlocal-intro} appears to be an implicit necessary requirement that the data $f$ and $g$ must fulfil beforehand for any attempt of solving the problems \eqref{eq:var-nonlocal-Neumann-intro} and \eqref{eq:nonlocal-Neumann-intro}. The local counterpart of this compatibility condition where $g$ is defined on $\partial\Omega$ (see \eqref{eq:local-Neumann-var}) is given by 
	\begin{align}
	\int_{\Omega} f(x)\mathrm{d}x +\int_{\partial\Omega} g(y)\mathrm{d}\sigma(y)=0.
	\end{align} 
\noindent It is worth highlighting  that as opposed to \cite{DROV17}, our functional test space $\VnuOm$ in the weak formulation \eqref{eq:var-nonlocal-Neumann-intro} does not depend on the Neumann data $g$. Moreover for the existence of weak solutions to \eqref{eq:var-nonlocal-Neumann-intro} we essentially choose $f\in L^2(\Omega)$ and  $g= g'\nu_K$ with  $g'\in L^2(\Omega^c, \nu_K)$. 

\medskip

\noindent Another important remark is  that under the compatibility condition 
\eqref{eq:compatible-nonlocal-intro}, the problems \eqref{eq:nonlocal-Neumann-intro} and 
\eqref{eq:var-nonlocal-Neumann-intro} might possess multiple solutions. Indeed, both integrodifferential operators $L$ and $\mathcal{N}$ annihilate additive constants. Therefore, as long as $u$ is a solution to the system \eqref{eq:nonlocal-Neumann-intro} or to the variational problem \eqref{eq:var-nonlocal-Neumann-intro} so does the function $\widetilde{u}= u+c$ for any $c\in \mathbb{R}$. Accordingly, both problems are ill-posed in the sense of Hadamard. 
The situation is similar in the local setting for problems \eqref{eq:local-Neumann} and 
\eqref{eq:local-Neumann-var} with the operators $L$ and $\mathcal{N}$ respectively replaced by the operators $-\Delta$ and $\frac{\partial}{\partial n}$. To overcome this anomaly we assume that $\Omega$ is bounded and we introduce the functional Hilbert subspace $ \VnuOm^{\perp}$ consisting of functions in $ \VnuOm$ with zero mean over $\Omega$. To be more precise, 
\begin{align*}
\VnuOm^{\perp}:= \Big\{ u\in \VnuOm: \int_{\Omega}u(x)\mathrm{d}x=0\Big\}.
\end{align*}

\noindent The space $\VnuOm^\perp $  enables us to reformulate the problem \eqref{eq:var-nonlocal-Neumann-intro} as follows: find $u\in \VnuOm^\perp $ such that
\begin{align}\label{eq:var-nonlocal-Neumann-bis-intro}\tag{$V'$}
\mathcal{E}(u,v) = \int_{\Omega} f(x)v(x)\mathrm{d}x +\int_{\Omega^c} g(y)v(y)\mathrm{d}y,\quad \mbox{for all}~~v \in \VnuOm^{\perp}\,.
\end{align}

\noindent Let us emphasize  that in contrast to \eqref{eq:var-nonlocal-Neumann-intro}, the variational problem \eqref{eq:var-nonlocal-Neumann-bis-intro} possesses at most one solution since $\mathcal{E}(\cdot, \cdot)$ defines a scalar product on $\VnuOm^\perp$. Furthermore, if $u'$ solves \eqref{eq:var-nonlocal-Neumann-bis-intro} then all solutions of the variational problem \eqref{eq:var-nonlocal-Neumann} are of the form $u'+c$ with $c\in \R$. An analogous observation holds in the local setting for the problem 
\eqref{eq:local-Neumann-var}. There, one would need to introduce the space $H^1(\Omega)^\perp$ of 
the functions in $H^1(\Omega)$ whose mean over $\Omega$ vanishes.  

\medskip 
\noindent We point out  that if $f\in L^2(\Omega)$ and $g=g'\nu_K$ with $g'\in L^2(\Omega^c, \nu_K)$ then under some additional conditions on $\nu$ and $\Omega$ (which we mention below) we are able to prove that the problem  \eqref{eq:var-nonlocal-Neumann-bis-intro} has a unique solution in $\VnuOm$. From this perspective,  it is legitimate to say that the function spaces $\VnuOm$ and $L^2(\Omega^c, \nu_K)$ introduced earlier are of great  importance for the study of Neumann complement value problems. 

\medskip

\noindent In the same spirit, we are able to formulate the variational formulation of others IDEs. For example, let us introduce another variant of the problem \eqref{eq:nonlocal-Neumann-intro}. Let $b:\Omega^c\to \mathbb{R}$ 
be a measurable function. The perturbation of the complement Neumann condition by $\mathcal{N}u + b u =g$ 
gives rise to a new type of nonlocal complement value problem and allows  us to remedy the restriction to the space $\VnuOm^\perp$ in the weak formulation. This will be called the Robin complement condition. In the local setting, 
the Robin boundary condition is also known as the Fourier boundary condition or the third boundary condition. 
Given $ f:\Omega \to \mathbb{R}$ and $b, g:\Omega^c\to \mathbb{R}$ as measurable functions, the nonlocal Robin boundary problem related to the operator $L$ is to find a measurable function $u:\mathbb{R}^d\to \mathbb{R}$ such that 
	\begin{align}\label{eq:nonlocal-Robin-intro}
	L u = f \quad\text{in}~~ \Omega \quad\text{ and } \quad \mathcal{N} u+ b u= g ~~ \text{on}~~ \mathbb{R}^d\setminus\Omega.
	\end{align} 
	
	\noindent We say that $u$ is a weak solution of the problem \eqref{eq:nonlocal-Robin-intro} if $u\in \VnuOm$
 and $u$ satisfies the relation
	\begin{align}\label{eq:var-nonlocal-Robin_0}
	\mathcal{E}(u,v)+ \hspace*{-0.5ex}\int\limits_{\Omega^c} \hspace*{-0.5ex} b(x)u(x)v(x)\d x = \int\limits_{\Omega} \hspace*{-0.5ex} f(x)v(x)\mathrm{d}x +\int\limits_{\Omega^c} \hspace*{-0.5ex} g(y)v(y)\mathrm{d}y,\quad \mbox{for all}~~v \in \VnuOm\, .
	\end{align}

	\medskip
	
	\subsection*{Further nonlocal problems}
\noindent The function spaces defined above  also apply to the study of the following 
IDEs. 

\begin{itemize}
	\item Nonlocal Dirichlet problem: for $f\in L^2(\Omega)$ and $g\in \TnuOm$,
	\begin{align*}
	L u = f \quad\text{in}~~ \Omega \quad\text{ and } \quad  u= g ~~ \text{on}~~ \mathbb{R}^d\setminus\Omega.
	\end{align*} 
		\item Nonlocal Neumann problem: for $f\in L^2(\Omega)$ and $g\in L^2(\Omega^c, \nu_K)$,
	\begin{align*}
\quad	L u = f \quad\text{in}~~ \Omega \quad\text{ and } \quad \mathcal{N} u= g\nu_K ~~ \text{on}~~ \mathbb{R}^d\setminus\Omega.
	\end{align*} 
		\item Nonlocal Robin problem:  for $f\in L^2(\Omega)$ and $g\in L^2(\Omega^c, \nu_K)$
	\begin{align*}
\qquad\quad	L u = f \quad\text{in}~~ \Omega \quad\text{ and } \quad \mathcal{N} u+ b u= g\nu_K ~~ \text{on}~~ \mathbb{R}^d\setminus\Omega.
	\end{align*} 
		\item Nonlocal mixed problem: Assume $\Omega^c=D\cup N$, where $D$ and $N$ are measurable such that $|D\cap N|=0$. For $f\in L^2(\Omega)$ and $g_D\in \TnuOm$ and $g_N\in L^2(\Omega^c, \nu_K)$
	\begin{align*}
	L u = f \quad\text{in}~~ \Omega \quad\text{ and } \quad  u=g_D~ \text{on}~ D, ~~ \mathcal{N} u= g_N\nu_K~~ \text{on}~~ N.
	\end{align*} 
	
		\item Nonlocal Helmholtz problem: for  $\lambda\in \R$
	\begin{align*}
	L u -\lambda u= f \quad\text{in}~~ \Omega \quad\text{ and } \quad \mathcal{N} u=0 \text{  (or $u=0$)} ~~ \text{on}~~ \mathbb{R}^d\setminus\Omega.
	\end{align*} 
	\item Nonlocal heat equation: Let $f\in L^2(\Omega)$ and $T>0$.
\begin{align*}
\partial_t u + Lu = f ~~\text{in $\Omega\times [0, T)$},
\quad \mathcal{N} u =0  \text{  (or $u=0$)}~~\text{on $\Omega^c\times [0, T)$},
\quad u=u_0~~ \text{on  $\Omega\times \{0\}$}.
\end{align*}

\item Nonlocal Schr\"odinger equation: for  $u_0, f\in L^2(\Omega)$ and $T>0$
\begin{align*}
i\partial_t u + Lu = f ~~\text{in $\Omega\times [0, T)$},
\quad \mathcal{N} u =0  \text{ (or $u=0$)} ~~\text{on $\Omega^c\times [0, T)$},
\quad u=u_0~~ \text{on  $\Omega\times \{0\}$}.
\end{align*}
\item Nonlocal wave equation: for $u_1,u_0, f\in L^2(\Omega)$ and $T>0$,
\begin{align*}
\partial_{tt}^2 u + Lu = f ~~\text{in $\Omega\times [0, T)$},
\quad \mathcal{N} u =0  \text{ (or $u=0$)} ~~\text{on $\Omega^c\times [0, T)$},
~ \partial_t u=u_1, ~u=u_0~~ \text{on  $\Omega\times \{0\}$}.
\end{align*}
\end{itemize}

\noindent Our treatments of the above IDEs are motivated by many existing analogous concepts from the theory of elliptic PDEs of second order. We exploit follow the contents  of \cite{Dan00,Dorothee-Triebel,Ev10,Hunter14,LeDret16,Necas67,Ro14}.

\vspace{2mm}

 \noindent \subsection*{Main results and literature review} We wish to formulate our main results and provide some further references. 
The solvability of some of the above mentioned IDEs easily follows from  the Lax-Milgram Lemma if Poincar\'e type  inequalities hold. Frequently, the  Poincar\'{e} type inequalities can be derived via compactness arguments. Our first main result concerns the compactness of the embedding of $\HnuOm$ and $\VnuOm$ into $L^2(\Omega)$. Let us start with some basic observations and formulate some sufficient assumptions on $\nu$ and $\Omega$. To this end, we need to reinforce our general assumption on the function $\nu: \mathbb{R}^d\setminus\{0\}\to [0, \infty]$. Let us recall that $\nu$ satisfies
\begin{align}\tag{$I_1$}\label{eq:integrability-condition+even-intro}
\nu(-h) = \nu(h) ~\text{ for all $h\in \mathbb{R}^d$ and} ~\int_{\mathbb{R}^d}(1\land |h|^2)\nu(h)\d h<\infty\,.
\end{align}

\noindent First and foremost, observe  that for $\nu\in L^1(\mathbb{R}^d)$ the  space $\HnuOm$ coincides with $L^2(\Omega)$ and therefore, cannot be compactly embedded into $L^2(\Omega)$. Likewise, if $\nu\in L^1(\mathbb{R}^d)$, then the spaces $\VnuOm \cap L^2(\mathbb{R}^d)$ and $H_{\nu}(\mathbb{R}^d)$ coincide with $L^2(\mathbb{R}^d)$, which is not even locally compactly embedded in $L^2(\Omega)$. In other words, the 
least necessary condition for compact embeddings to hold is that $\nu$ is not integrable. Therefore, it is 
necessary to consider the following non-integrability condition 
\begin{align}\tag{$I_2$}\label{eq:non-integrability-condition-intro}
\int_{\mathbb{R}^d}\nu(h)\d h=+ \infty.
\end{align}

\noindent By having the condition \eqref{eq:integrability-condition+even-intro} at hand, it is possible to strengthen the condition \eqref{eq:non-integrability-condition-intro} by assuming  that 
\begin{align}\tag{$I'_2$}\label{eq:limit-at-0-explode-intro}
\lim_{|h|\to 0}|h|^d \nu(h)= \infty. 
\end{align}
\noindent That is, the  condition \eqref{eq:limit-at-0-explode-intro} clearly implies 
\eqref{eq:non-integrability-condition-intro}. It is important to point out  that if $\Omega$ is bounded then 
conditions \eqref{eq:integrability-condition+even-intro} and \eqref{eq:non-integrability-condition-intro} are 
sufficient to obtain the local compactness of $\HnuOm$ and $\VnuOm$ into $L^2(\Omega)$ (see 
Theorem \ref{thm:local-compactness}). In fact, this is reminiscent of the main result of \cite{JW19}, which shows  
that the embedding  $\VnuOmO\hookrightarrow L^2(\Omega)$  is compact. We will revisit this result under slightly modified assumptions.

\medskip

\noindent  The global compactness requires some extra  compatibility assumptions 
between  $\Omega$ and $\nu$.  We establish the global compactness by exploiting the recent results from \cite{JW19} and \cite{DMT18}. We intend to provide an alternative approach to the compactness result in \cite[Theorem 2.2]{CP18}. The technique therein is adapted from \cite[Theorem 7.1]{Hitchhiker} for fractional Sobolev spaces and uses the Sobolev extension property of the corresponding domain. However, the proof provided in \cite{CP18} only  seems to be valid for domains that can be written as a finite union of cubes; unless 
the corresponding nonlocal function space possesses the extension property. Our approach is rather standard and follows the idea used to prove  the classical  Rellich–Kondrachov theorem, i.e. the compactness of the embedding $W^{1,p}(\Omega)\hookrightarrow L^p(\Omega)$ for  sufficiently smooth $\Omega$. It involves applying the local compactness and using an approximation argument near the boundary of $\Omega$. 

\medskip

\noindent Let us  introduce some regimes relating $\Omega$ and $\nu$ under which the global compactness holds true. We will enumerate these assumptions on  $(\nu, \Omega)$ into different classes. We say that  the couple $(\nu, \Omega)$ is in the class $\mathscr{A}_i,~~i=1,2,3$ if $\Omega\subset \mathbb{R}^d$ is an open bounded set and $\nu: \mathbb{R}^d\setminus\{0\}\to [0, \infty]$ satisfies the conditions \eqref{eq:integrability-condition+even-intro} and \eqref{eq:non-integrability-condition-intro} and additionally $\nu$ and $\Omega$ satisfy:  

\vspace{1mm}

\noindent $\bullet$ The class $\mathbf{\mathscr{A}_1}$: there exists a $\WnuOm$-extension operator $E: \WnuOm\to W^p_{\nu}(\mathbb{R}^d)$, i.e.  there is a constant $C: C(\nu, \Omega, d)>0$ such that for every $u\in \WnuOm$, $\|Eu\|_{W^p_\nu(\mathbb{R}^d)}\leq C\|u\|_{\WnuOm}$ and $Eu|_\Omega =u$.

\vspace{1mm}

\noindent $\bullet$ The class $\mathbf{\mathscr{A}_2}$: $\Omega$ has a Lipschitz boundary, $\nu$ is radial and 
\begin{align}\label{eq:class-lipschitz-intro}
q(\delta):= \frac{1}{\delta^p}\int_{B_\delta(0)} |h|^p\nu(h)\d h \xrightarrow[]{\delta \to 0}\infty. 
\end{align}

\vspace{1mm}

\noindent $\bullet$ The class $\mathbf{\mathscr{A}_3}$: setting $\Omega_\delta= \{x\in \Omega: \operatorname{dist}(x,\partial\Omega)>\delta\}$ for $\delta>0$,  the following condition holds true

\begin{align}\label{eq:class-sing-boundary-intro}
\widetilde{q}(\delta): = \inf_{a\in \partial\Omega}\int_{\Omega_\delta}\nu(h-a)\d h \xrightarrow[]{\delta \to 0}\infty.
\end{align}

\medskip

\noindent Let us introduce a fourth class $\mathbf{\mathscr{A}_4}$ of interest. 

\noindent $\bullet$ The class $\mathbf{\mathscr{A}_4}:$ we say that the couple $(\nu, \Omega)$ is in the class $\mathbf{\mathscr{A}_4}$ if  $\Omega$ is any open bounded subset of $\mathbb{R}^d$ and $\nu:\mathbb{R}^d\setminus\{0\}\to [0, \infty]$ is a unimodal L\'{e}vy measure that is, $\nu$ is radial, almost decreasing and $\nu \in L^1(\R^d, 1\land |h|^2\d h)$. 
\medskip

\noindent Let us mention that the class $\mathscr{A}_2$ is inspired from \cite{DMT18}. Moreover, in the class 
$\mathscr{A}_4,$ $\nu$ is not necessarily singular near $0$. Note that if $\Omega$ is bounded and Lipschitz, 
then the couple $(|\cdot|^{-d-\alpha}, \Omega)$ with $\alpha\in (0,2)$ belongs to  each $\mathscr{A}_i, ~i=1,2,3,4$.  We discuss the classes $\mathscr{A}_i, ~i=1,2,3,4$ in Chapter \ref{chap:nonlocal-sobolev}.  Here is our compactness result. 

\begin{theorem*}[Theorem \ref{thm:embd-compactness}]
	Let $\Omega\subset \mathbb{R}^d$ be open and bounded and let $\nu:\mathbb{R}^d\setminus\{0\} \to [0,\infty]$ be a measurable function. If the couple $(\nu, \Omega)$ belongs to one of the classes
$\mathscr{A}_i,~ i=1,2,3$ then the embedding $\HnuOm \hookrightarrow L^2(\Omega)$ is compact. In particular, the embedding $\VnuOm \hookrightarrow L^2(\Omega)$ is compact.
\end{theorem*}
\noindent A noteworthy consequence of our compact embeddings result is the well-known Rellich-Kondrachov compactness theorem which implies the compactness of the embedding
 $H^1_0(\Omega)\hookrightarrow L^2(\Omega)$ (resp. $H^1(\Omega)\hookrightarrow L^2(\Omega)$ when $\Omega$ has a Lipschitz boundary). Indeed, the embeddings $H^1_0(\Omega)\hookrightarrow \VnuOmO \hookrightarrow L^2(\Omega) $ (resp. $H^1_0(\Omega)\hookrightarrow \VnuOmO \hookrightarrow L^2(\Omega) $) are continuous. Another crucial application of the above compactness theorem is the Poincar\'{e} inequality (cf. Theorem \ref{thm:poincare-inequality}) which also holds for the class $\mathscr{A}_4$. To be more precise,  
 we are able to show that if the couple $(\nu, \Omega)$ belongs to one of the classes $\mathscr{A}_i,~ i=1,2,3,4$ then there exists $C=C(\nu,\Omega,d)>0$  such that for every $u\in L^2(\Omega)$ we have 
\begin{align*}
\big\|u- \mbox{$\fint_\Omega u$}\big\|^2_{L^2(\Omega)} \leq C\iil_{\Omega\Omega}(u(x)-u(y))^2\nu(x-y)\d x\d y.
\end{align*}

\noindent  Accordingly, it follows that  for every  $u\in L^2(\Omega)$ we have 
 \begin{align*}
 \big\|u- \mbox{$\fint_\Omega u$}\big\|^2_{L^2(\Omega)} \leq C\mathcal{E}(u,u). 
 \end{align*}
 
\noindent In the same spirit, by establishing a Poincaré-Friedrichs type inequality, we provide an alternative short proof of \cite[Lemma 2.7]{FKV15}. That is, for every $u\in \VnuOmO$ we also have
 \begin{align*}
\big\|u\big\|^2_{L^2(\Omega)} \leq C\mathcal{E}(u,u). 
\end{align*}

\noindent As complementary results on function spaces, we establish the density of smooth functions in the spaces $\HnuOm,\VnuOm$ and $\VnuOmO$. In particular, we improve the density result from \cite[Lemma 2.12]{Voi17} and \cite[Theorem A.4]{BGPR17} whose proofs are incomplete. The results can be summarized as follows. 

\begin{theorem*}[Theorem \ref{thm:density}, \ref{thm:density-bis} \& \ref{thm:density-zero-outside}]
	Assume that $\Omega\subset \mathbb{R}^d$ is open. 
	\begin{itemize}
	\item  $\HnuOm\cap C^\infty(\Omega)$ is dense in $\HnuOm$.
	\item  If $\partial \Omega$ is Lipschitz and  compact or $\Omega=\R^d$ then $C_c^\infty(\R^d)$ is dense in $\VnuOm$ (\cite{Voi17,FGKV19}).
	\item If $\partial \Omega$ is continuous and  compact  then $C_c^\infty(\Omega)$ is dense in $\VnuOmO$  (\cite{BGPR17,FKV15}).
	\end{itemize}
\end{theorem*}

\noindent We established the well-posedness of  elliptic IDEs. In particular, we have the following result. 

\begin{theorem*}
		Assume that the couple $(\nu, \Omega)$ belongs to one of the classes $\mathscr{A}_i,~i=1,2,3,4$. 
	Then given $ f \in L^2 (\Omega)$ and $g=g'\nu_K$ with $g'\in L^2(\Omega^c, \nu_K)$, there exists a unique solution $u\in\VnuOm^{\perp} $ to the variational problem \eqref{eq:var-nonlocal-Neumann-bis-intro}.  If the compatibility \eqref{eq:compatible-nonlocal-intro} holds,  then all solutions to \eqref{eq:var-nonlocal-Neumann-intro} are of the form $w=u+c$ with $c\in \mathbb{R}$. Moreover, there exists a constant $C: = C(d,K, \Omega, \nu)>0$ independent of $f$ and $g$ such that any solution $w$ of $\eqref{eq:var-nonlocal-Neumann-intro}$ satisfies the following weak regularity estimate
	\begin{align*}
	\|w-\hbox{$\fint_{\Omega}w$} \|_{\VnuOm}\leq C \Big(\|f\|_{L^2(\Omega)}+\|g\|_{L^{2}(\Omega^c, \nu_K^{-1})}\Big).
	\end{align*}
	\end{theorem*}

\vspace{1mm}

\noindent Note  that  the study of the nonlocal Neumann problem \eqref{eq:nonlocal-Neumann-intro} 
was introduced in \cite{DROV17} for the fractional Laplace operator. The significant difference to our 
approach to studying \eqref{eq:var-nonlocal-Neumann-intro}, however, is that the test space introduced therein depends 
on the Neumann data $g$. Following the approach of \cite{DROV17}, the recent articles \cite{BMPS18,ML19} also 
study inhomogeneous nonlocal problems with the Neumann complement condition. We point out that \cite{ML19}
is a remake of \cite{DROV17} is the $L^p$-setting while considering the Neumann problem for the  so called $p$-fractional Laplacian operator.  Most importantly, the nonlocal normal derivative $\mathcal{N}$ is currently appearing in more and more works. 
For example, see \cite{AL18,Che18,CC20,FJS19} for the study of nonlocal semilinear problems with homogeneous 
Neumann complement condition, \cite{LMPPS18} for the study of the principal eigenvalue of the fractional Laplacian with a mixed complement condition, \cite{Zoran19} for the study of a probability interpretation of 
nonlocal quadratic forms, \cite{Aba20}  for a comparative study on different types of nonlocal Neumann conditions and \cite{GSU16} for the study of the Calder\'on problem for the fractional Laplacian. It is important to highlight  that some authors prefer to formulate nonlocal Neumann problems via the regional fractional Laplacian \cite{Warma2018,Warma2016,Warma2015, Gru16,CS16}. The homogeneous Neumann problem for nonlocal regional type operators is also studied in the area of peridynamic models see for instance \cite{TTD17}. However, other authors work on the Neumann problem by  defining the fractional power of the  Neumann Laplacian \cite{MBG12,SB15,DSV15}. In the latter contexts, the Neumann conditions are rather prescribed on the boundary of the underlying domains, therefore, our set-up does not apply. 

\vspace{2mm}

\noindent There exists a substantial amount of literature on nonlocal Dirichlet
problems. For example, for their solvability see \cite{FKV15} where the topic is extended to nonlocal operators with nonsymmetric kernels. See \cite{Ru18}  for a study of nonlocal Dirichlet problems involving symmetric nonlocal operators whose  driven jump interaction measure need not be absolutely continuous.See \cite{PR18} for an examination of the Dirichlet problem for the fractional Laplacian in perforated domains and also \cite{Ros15,hoh96} for complementary approaches. There are also several works on related subjects, for example, see \cite{DyKa20,COS17} for interior regularity of solutions, \cite{FK13} for interior regularity of parabolic problems, \cite{Ros15,Ros14} for regularity up to the boundary, \cite{CS09,CS11,Li14} for regularity for viscosity solution and \cite{JW19max} for the maximum principle. There are also several works on related subjects for example, for interior regularity of solutions see \cite{DyKa20,COS17}, for interior regularity of parabolic problems see \cite{FK13}, for regularity up to the boundary see \cite{Ros15,Ros14}, for regularity for viscosity solutions see \cite{CS09,CS11,Li14} and for the maximum principle see \cite{JW19max}.
 
\vspace{2mm}

\noindent Having a suitable set-up for nonlocal Dirichlet and Neumann problems to hand makes looking at some aspects of the corresponding Dirichlet-to-Neumann map legitimate (at least the definition and the spectrum).The Dirichlet-to-Neumann operator exposed here is largely inspired by \cite{WoRa07,WoRa12,Bet15,AR19}, where an analogous subject is treated for the Laplace operator. The Dirichlet-to-Neumann map plays a crucial role in the study of the so-called Calder\`on problem, see for instance \cite{FSU19}. Our approach leads to a slightly different Dirichlet-to-Neumann map than the one derived in \cite{GRSU20,RS20} for the fractional Laplacian.

\vspace{2mm}

\noindent Meanwhile, the compactness result constitutes a powerful tool for further investigations. Indeed, with the compactness of the embedding $\VnuOm\hookrightarrow L^2(\Omega)$  in force, we are able to analyse the
 following:

\begin{itemize}
\item The spectral decomposition of the integrodifferential operator $L$ subject to the Dirichlet, Neumann or Robin complement condition using the Rayleigh quotient.  Indeed, the operator $L$ turns out to have a discrete spectrum and a compact resolvent like  the Laplace operator.  
\item the spectrum of the nonlocal Dirichlet-to-Neumann map via the spectrum of the operator $L$ subject to the Robin complement condition. To do this, we closely rely on the approach from  \cite{WoRa07, WoRa12} where an analogous  characterization  is  derived for the Laplace operator.   

\item The essentially self-adjointness for the operator $L$ with  the Dirichlet, Neumann or Robin complement condition. Indeed, the operator $L$ turns out to be a symmetric unbounded operator on $L^2(\Omega)$ the same as for the Laplacian. Our survey on the essentially self-adjointness for the operator $L$ closely follow the material in \cite{Kow09,Dav96} where a similar study is carried out for the Laplace operator. Let us mention that a different study of the essentially self-adjointness for the fractional Laplacian is studied in \cite{HKM17}  in the context where $\Omega\subset \R^d$ is the complement of a compact set.

\item The profile solutions to evolution equations involving the operator $L$ which are  Initial Complement Value Problems (ICVP). The parabolic equation with the Neumann complement condition for the fractional Laplacian is also discussed in \cite{DROV17}. We  exploit some ideas from \cite{Ro14} which treats the parabolic problem elliptic with the Robin boundary condition for the Laplace operator.
\end{itemize}

\noindent We have seen several familiarities between the characteristics of the (\textit{local}) Laplace operator $-\Delta$ and the characteristics of the (\textit{nonlocal})  integrodifferential operator $L$. 
Next, we want to  bridge a connection from the nonlocal world to the local one. Strictly speaking, we try to understand local  objects as limits of nonlocal objects. Accordingly, let us introduce $(\nu_\alpha)_{0<	\alpha<2}$, a family of L\' evy radial functions  approximating the Dirac measure at the origin, i.e. for every $\alpha, \delta > 0$
\begin{align*}
\begin{split}
\nu_\alpha\geq 0\,\,\text{ is radial}, \quad \int_{\mathbb{R}^d}	(1\land |h|^2)\nu_\alpha (h)\d h=1, \quad \lim_{\alpha\to 2}\int_{|h|>\delta}	\nu_\alpha(h)\d h=0\,.
\end{split}
\end{align*}                      

\noindent For a family  $(J^\alpha)_{0<\alpha<2}$ of  positive symmetric kernels $J^\alpha: \mathbb{R}^d\times \mathbb{R}^d \setminus \operatorname{diag} \to [0, \infty]$ we set-up the following conditions: 
 
\begin{itemize}
%
	\item[(G-E)] There exists a constant $\Lambda\geq 1$ such that for every $\alpha\in  (0,2)$ and all $x,y \in \mathbb{R}^d$, with $x\neq y$, 
	\begin{align}\label{eq:global-elliptic-condition-intro}\tag{$G$-$E$}
	\Lambda^{-1} \nu_\alpha (x-y) \leq J^\alpha(x,y) \leq  \Lambda \nu_\alpha (x-y).
	\end{align}
\end{itemize}

\noindent Given $x \in \R^d$ and $\delta > 0$, we define the symmetric matrix $A(x) = (a_{ij}(x))_{1\leq i,j\leq d}$ by
\begin{align*}
a_{ij}(x) = \lim_{\alpha\to 2^{-}} \int_{B_\delta}  h_ih_j J^\alpha(x,x+h)dh.
\end{align*}
It is noteworthy to mention that the matrix $A= (a_{ij})_{ij}$  does not depend on the choice of $\delta>0$ and satisfies the elliptic condition
\begin{align*}
\Lambda^{-1} \, d^{-1}|\xi|^2\leq \langle A(x)\xi,\xi\rangle\leq \Lambda \, d^{-1}|\xi|^2, \quad \text{for all $x,\xi\in \R^d$}.
\end{align*}
Furthermore, assume we have  $J^\alpha(x,y) = C_{d,\alpha}|x-y|^{-d-\alpha}$ for a suitable choice of $\nu_\alpha(h)= c_\alpha|h|^{-d-\alpha}$ then $a_{ij}=0$ if $i\neq j$ and $a_{ii}=1$. Moreover, motivated by \cite{BBM01},  we show that, for sufficiently smooth $\Omega$, the nonlocal spaces $\HnuOma$ and $\VnuOma$ both converge to  the Sobolev space $H^1(\Omega)$ as $\alpha\to 2$. To be more precise, we have the following simplified result (see Section \ref{sec:charact-W1p}).
\begin{theorem*} Let $\Omega$ be  an $H^1(\Omega)$-extension domain.  Then for all $u \in H^1(\Omega)$ we have 
	\begin{align*}
		&\lim_{\alpha\to 2^-} \iil_{\Omega\Omega}(u(x)-u(y))^2\nu_\alpha(x-y)\d x\d yx= \lim_{\alpha\to 2^-}  \hspace{-2ex}\iil_{(\Omega^c\times \Omega^c)^c}\hspace{-3ex} (\overline{u}(x)-\overline{u}(y))^2\nu_\alpha(x-y)\d x\d y= K_{d,2}\int_{\Omega}|\nabla u(x)|^2\d x. 
	\end{align*}
	where $\overline{u}\in H^1(\R^d)$ is any extension of $u$ to $\R^d$ and $K_{d,2}= \frac{1}{d}$. 
	Moreover, a function $u\in L^2(\Omega)$, belongs to $H^1(\Omega)$ if and only if 
	\begin{align*}
	\liminf_{\alpha\to 2^-} \iil_{\Omega\Omega}(u(x)-u(y))^2\nu_\alpha(x-y)\d x\d y<\infty.
	\end{align*}
\end{theorem*}

\noindent Note that, the above convergence result can be found in \cite{BBM01}  when $\Omega$ is a bounded Lipschitz domain and in \cite{Brezis-const-function} for $\Omega=\R^d$,  and however include the case where $\Omega$ is an unbounded domain with an $H^1(\Omega)$-extension property. Next, let us define the  operators  $\mathscr{L}_\alpha$ and $\mathscr{N}_\alpha$ by
\begin{align*}
\mathscr{L}_\alpha u(x):= \pv  \int_{\R^d} (u(x)-u(y)) J^\alpha (x,y)\d y \quad \text{and}\quad \mathscr{N}_\alpha u(y):= \int_{\Omega} (u(y)-u(x)) J^\alpha (x,y)\d x \,.
\end{align*}	
\noindent Let us also introduce the outwards normal derivative of a function $v$ on $\partial \Omega$ with respect to $A$ defined for $x\in \partial \Omega$ by $ \frac{\partial v}{\partial n_A}(x) = A(x) \nabla v(x) \cdot n(x).$
\noindent We now formulate some convergence results  in a simplified setting.  
\begin{theorem*}[Theorem \ref{thm:convergence-solutioin-I}] Let $\Omega\subset \mathbb{R}^d$ be open bounded with a Lipschitz boundary and connected. Let $(f_\alpha)_\alpha$ be a family converging weakly  to some $f$ in $L^2(\Omega)$ as $\alpha\to2$. Define $g_\alpha=  \mathscr{N}_\alpha\varphi$ and $g=\frac{\partial \varphi}{\partial n_A}$ for $\varphi \in C^2_b(\mathbb{R}^d)$.  Assume that the condition \eqref{eq:global-elliptic-condition} holds and suppose $u_\alpha \in  \VnuOma^\perp$ is a weak solution of the nonlocal Neumann problem
 \begin{align*}
	 \mathscr{L}_\alpha u= f_\alpha\text{ on~}\Omega\quad \text{and}\quad\mathscr{N}_\alpha u= g_\alpha\text{ on~}\Omega^c. 
\end{align*}	
	
	\noindent  Let $u\in H^1(\Omega)^\perp$ be the unique weak solution in $H^1(\Omega)^\perp$ to the  Neumann problem 
\begin{align*}
	-\operatorname{div}(A(\cdot) \nabla ) u=f\, \, \text{ on ~$\Omega$}\quad\text{ and}\quad  \frac{\partial u}{\partial n_A} =g \text{ on ~$\partial\Omega$}.
	\end{align*}	
Then $(u_\alpha)_\alpha$ converges to $u$ in $L^2(\Omega)$  as $\alpha\to 2$, i.e. $\|u_\alpha-u\|_{L^2(\Omega)}\xrightarrow[]{\alpha\to 2}0.$
\end{theorem*}
\vspace{2mm}

\noindent Likewise, if we assume that $g\in H^1(\R^d)$ then under the condition \eqref{eq:global-elliptic-condition-intro} (see Theorem \ref{thm:convergence-solutioin-III})  weak solutions to nonlocal Dirichlet problems $\mathscr{L}_\alpha u= 
f_\alpha$ on $\Omega$ and $ u= g$  on $\Omega^c$  with $\alpha\in (0,2)$) converge in $L^2(\Omega)$ to  the weak solution of the local Dirichlet problem $-\operatorname{div}(A(\cdot)\nabla)u= f$ on $\Omega$ and $ u= g$  on $\partial\Omega$.

\vspace{2mm}

\noindent Furthermore, we also show the convergence of normalized eigenpairs. A couple $(\lambda, \phi)$, with $\lambda\in \R$ and $u\in L^2(\Omega)$  is called a normalized eigenpair of $\mathscr{L}_\alpha$ subject to the Neumann (resp. Dirichlet) complement condition if $\|u\|_{L^2(\Omega)} =1$  and $u$ is a weak solution to the nonlocal Neumann (resp. Dirichlet) problem $\mathscr{L}_\alpha u= \lambda u $ in $\Omega$ and $\mathscr{N}_\alpha u=0$ on $\Omega^c$ (resp. $ u=0$ on $\Omega^c$). In the same manner, one can define an eigenpair of the  elliptic operator $-\operatorname{div}(A(\cdot)\nabla)$.  
 
\begin{theorem*}[Theorem \ref{thm:convergence-eigenpair-I}--\ref{thm:convergence-eigenpair-II}]
Let $(\lambda_\alpha,u_\alpha)$ be a  normalized eigenpair of the nonlocal operator $\mathscr{L}_\alpha$. 
Then, up to a subsequence,  $(\lambda_\alpha,u_\alpha)_\alpha$ converge in $\R\times L^2(\Omega)$  
to couple a $(\lambda,u)$  where  the latter is  a normalized eigenpair of the local operator $-\operatorname{div} (A(\cdot) \nabla)$. 
\end{theorem*}

\noindent It is worth noting that if $J^\alpha(x,y)= C_{d,\alpha}|x-y|^{-d-\alpha}$ then $\mathscr{L}_\alpha=(-\Delta)^{\alpha/2}$ and $-\operatorname{div}(A(\cdot)\nabla) =-\Delta$.  
Correspondingly, weak solutions (resp. normalized eigenpairs) of the fractional Laplacian converge to weak 
solutions (resp. normalized eigenpairs) of the Laplacian as $\alpha\to 2$.

\vspace{2mm}

\noindent Let us point out that the  crucial tools for accomplishing the aforementioned convergence results are  robust Poincar\'e type inequalities and  Mosco convergence of nonlocal quadratic forms. First, by exploiting the  technique of \cite{Ponce2004},  we show the following robust Poincar\'e inequalities: there exists $\alpha_0\in(0,2)$  and a constant $C= C(d, \Omega)$ such that
\begin{align}\label{eq:robust-poincare-intro}
\big\|u- \mbox{$\fint_\Omega u$}\big\|^2_{L^2(\Omega)} \leq C\iil_{\Omega\Omega}(u(x)-u(y))^2\nu_\alpha(x-y)\d x\d y,\quad\text{for all  $\alpha\in (\alpha_0,2)$ and $u\in L^2(\Omega)$},
\intertext{and}
\big\|u\big\|^2_{L^2(\Omega)} \leq C\iil_{\R^d\R^d}(u(x)-u(y))^2\nu_\alpha(x-y)\d x\d y,\quad\text{for all  $\alpha\in (\alpha_0,2)$ and $u\in C_c^\infty(\Omega)$}.\label{eq:robust-poincare-friedrcih-intro}
\end{align}
 The inequality \eqref{eq:robust-poincare-intro} (resp. \eqref{eq:robust-poincare-friedrcih-intro})  is robust in the sense that  the constant $C$  does not depend on $\alpha$ and by letting, $\alpha\to 2$ gives the classical Poincar\'e (resp. Poincar\'e-Friedrichs) inequality
 \begin{align*}
 \big\|u- \mbox{$\fint_\Omega u$}\big\|^2_{L^2(\Omega)} & \leq C\int_{\Omega}|\nabla u(x)|^2\d x,\quad\text{for all $u\in L^2(\Omega)$},\\
 (\text{resp.}\quad \big\|u\big\|^2_{L^2(\Omega)} &\leq C\int_{\Omega}|\nabla u(x)|^2\d x,\quad\text{for all  $u\in C_c^\infty(\Omega)$}).
 \end{align*}

\noindent On the other hand, we establish \cite{FGKV19} the Mosco convergence of the nonlocal to local quadratics forms. 

\begin{theorem*}[Theorem \ref{thm:Mosco-convergence}]
Assume $\Omega\subset \R^d$ is open and bounded with a Lipschitz boundary. Assume the condition \eqref{eq:global-elliptic-condition-intro} holds. Then, as $\alpha\to 2$,  the family of nonlocal quadratic forms $(\mathcal{E}^\alpha, \VnuOma)_\alpha$ and $(\mathcal{E}^\alpha_\Omega, \HnuOma)_\alpha$ converge to the local quadratic form $(\mathcal{E}^A, H^1(\Omega))$ whereas, $(\mathcal{E}^\alpha, V_{\nu_\alpha}^\Omega(\Omega|\R^d))_\alpha$ converges to $(\mathcal{E}^A , H^1_0(\Omega))$. Here, we define
\begin{align*}
	 	\mathcal{E}^\alpha(u,u) &:=\iil_{(\Omega^c\times\Omega^c)^c} (u(x)-u(y))^2J^\alpha(x,y)\d x\d y,\\
\mathcal{E}^\alpha_\Omega(u,u)&:=\iil_{\Omega\Omega} (u(x)-u(y))^2J^\alpha(x,y)\d x\d y,\\
\mathcal{E}^A(u,u)&:=\int_{\Omega} (A(x)\nabla u(x)\cdot \nabla u(x))\d x.
\end{align*}
\end{theorem*}
\noindent In \cite{Mos94} (see also \cite{KS03}) it is shown that Mosco convergence of a sequence of symmetric closed forms is equivalent to the 
convergence of the sequence of associated semigroups (or of the associated resolvents) and implies the weak convergence of  the finite-dimensional distributions of the corresponding processes if any.  Note that several authors have studied the weak convergence of Markov processes with the help of Dirichlet forms, e.g., in \cite{LyZh96, KuUe97,SU16, MRZ98, Sun98, Kol05, Kol06, BBCK09, CKK13,CS02}. Most of the related results are concerned with situations where the type of the process does not change, i.e., the diffusions converge to a diffusion or jump processes converge to a jump process. The present work considers the cases  where a sequence of jump processes in bounded domains converges to a diffusion process.

\vspace{2mm}
\noindent  \subsection*{Outline} The thesis is organized as follows.  Chapter \ref{chap:basics-nonlocal-operator} is devoted to the introduction of the  basics of integrodifferential operators. First, we look at different  characterizations of the operator $L$. We end up with fourteen  characterizations of the fractional Laplacian 
$(-\Delta)^{\alpha/2}$. We compute and study the asymptotic behavior of  the normalization constant $C_{d,\alpha}$. Afterwards, we define the notion of nonlocal elliptic operators in divergence and non-divergence form,  and we show some correlations with local elliptic operators of second order. Finally we define nonlocal mixed
(anisotropic) operators. 

\vspace{2mm}
\noindent In  Chapter \ref{chap:nonlocal-sobolev} we introduce  nonlocal Sobolev-like spaces that are in a certain sense generalized Sobolev-Slobodeckij spaces which we often encounter. Roughly speaking, these are just some refinement of classical Lebesgue $L^p$-spaces with $1\leq p< \infty$ (like the classical Sobolev spaces $W^{1,p}(\Omega)$ and $W_0^{1,p}(\Omega) $) whose additional structures are of importance. As such spaces are less common, we will examine some of their rudimentary properties, e.g. their Banach structure, their relation with the classical Sobolev spaces, their embeddings, their approximation by smooth functions, their  extension property, their compact embeddings into $L^p(\Omega)$ and Poincar\'{e} type inequalities. A complement to this chapter is recorded in the Appendix \ref{chap:complemnt-Lebesgue}. 

\vspace{2mm}

\noindent Chapter \ref{chap:IDEs} is dedicated to the solvability of nonlocal IntegroDifferential Equations (IDEs) and some related problems such as the spectral decomposition, the essentially self-adjointness of the integrodifferential operator $L$ and the Dirichlet-to-Neumann map for $L$. We also look at some nonlocal evolution problems that are Initial Complement Value Problems (ICVP). 

\vspace{2mm}

\noindent Finally, in Chapter \ref{chap:Nonlocal-To- Local}, we  deal with convergence transitions from nonlocal to local.  We start with the convergence of spaces and characterize classical Sobolev spaces with the help of  nonlocal spaces. Afterwards, we establish robust Poincar\'e type inequalities. We conclude by proving the convergence of solutions and eigenpairs of nonlocal problems, to the local ones.


\chapter*{Acknowledgement}
This work would not have been possible without the support and collaboration of the many contributors. The financial support from the DeutscheForschungsgemeinschaft (DFG, German Research Foundation), which I obtained during my final years as a PhD student within the International Research Training Group (IRTG 2235): "Searching for the regular in the irregular: Analysis of singular and random system" is gratefully acknowledged. 

\noindent I express my deep and immense gratitude  to my advisor Prof.Dr Moritz Kassmann for his patience, guidance, support  and valuable suggestions throughout the development of this doctoral thesis. I would like to express a special thank to Prof Ki-Ahm Lee for his collaboration during my exchange stay at the Seoul National University (SNU) as his guest student. This collaboration has been an impetus among many others. I would to address special thanks to Dr Angkana R\"uland, Dr Vanja Wagner and  Dr Bart{\l}omiej Dyda  with whom I had several uplifting and stimulating discussions during their stay at the Bielefeld University and during our international meetings abroad. 

\noindent This work has been fuelled by my participation in conferences, workshops, seminars, schools and excursions, which would not have been possible without the financial support of the IRTG 2235. Special thanks goes to the IRTG administrators Nadine Brehme, Claudia K\"ohler, Anke Bodzin and Rebecca Reischuk for organizing my travel to all these international meetings. 
 
\noindent An immense thank you to all my IRTG fellows for making our leisure time and extra-academic activities more familiar, social and relaxed and to our IRTG fellows from SNU for their hospitality and their guidance in South Korea during our six-month stay. I am also grateful to my family and friends for being a source of encouragement and motivation. A big thank you to my officemate Fillip Bosni\'{c} for our stimulating daily discussions of mathematical questions. Finally, I wish to thank my friends who helped me with my other life problems while I was focused on this work. I am particularly indebted to Melissa Meinert, Robert Schippa, Andre Schenke and Marvin Weidner who helped me to proofread this thesis

\noindent I am highly indebted to my friend Frances Pairaudeau from the UK for her valuable help and support after we met at AIMS-Cameroon in 2016. Inviting me to the  UK granted my stay with a very nice vacation in Europe.  
A huge thank you to my landlady in Germany Petra Udelhoven for her warm hospitality and for opening her doors to me. 

\newpage

\noindent \textbf{Abgrenzung des eigenen Beitrags gemäß §10(2) der Promotionsordnung}

\noindent The proofs of Theorem \ref{thm:Mosco-convergence}, Theorem \ref{thm:Mosco-convergence-bis}  and the density of $C_c^\infty(\R^d)$ in $\VnuOm$ of Theorem \ref{thm:density}  (for $p=2$) are published in \cite{FGKV19} and were established in collaboration between the author,  his supervisor and Dr. Paul Voigt.

	\newpage
	
\chapter*{Notation}
We collect various notations that will be frequently used in this thesis. We only summarize here the most common notations.
 \begin{itemize} 
 	\item Throughout, $d\geq 1$ is an integer and $\mathbb{R}^d$ 
 	represents the $d$-dimensional Euclidean space furnished with the usual Euclidean inner product defined for two elements $x=( x_1,x_2,\cdots+ x_d) $ and $ y=(y_1,y_2,\cdots,y_d)$ by $x\cdot y = x_1y_1+x_2y_2\cdots+ x_dy_d$ and the norm shall be denoted by $|x|= \sqrt{x\cdot x}$. 
 	\item We shall assume $\mathbb{R}^d$ is automatically equipped with the topology induced by this norm and we denote an open ball of radius $r>0$ centered at $x$ by $B_r(x)$ or merely $B_r$ if $x=0$.
 	Further the space $\mathbb{R}^d$ will be furnished with the Borel $\sigma$-algebra and the Lebesgue measure $dx$. We will simply write "measurable" instead of "Borel measurable" for sets and functions. 
 	\item Given $A$ a subset of $\mathbb{R}^d$ we shall synonymously write $A^c$ or $\mathbb{R}^d\setminus A$ to designate the complement of $A$ and $\mathds{1}_A $ denotes the characteristic function of $A$. 
 	
 	\item The notation $\partial A $ shall stand for the boundary of $A$ that is $ \partial A= \overline{A}\setminus \overset{\circ }{A}$ where $\overline{A}$ and $\overset{\circ }{A}$ are respectively the closure and the interior of $A$. 
 	\item Given two sets $A$ and $B$ define $\operatorname{dist}(A, B) = \inf\big\{|a-b|:~a\in A \,, b\in B\big\}$ and for $x\in \R^d$, $\operatorname{dist}(x, A) = \operatorname{dist}(\{x\}, A) $.
 	
 	\item Besides, if $A$ is measurable we shall write $|A|$ to denote the Lebesgue measure of $A$ and also $|\partial A|$ to denote the Hausdorff measure of $\partial A$( which in this context is considered as the restriction to $\partial A$ of Lebesgue measure). Especially $ |\mathbb{S}^{d-1}|$ denotes the surface of the $d-1$-dimensional sphere of $\mathbb{R}^d$.
 	 \item For a measurable function $u$ and a measurable set such that $0<|A|<\infty$, we denote the mean value of $u$ over $A$ by  $\displaystyle\fint_A u= \frac1{|A|}\int_A u(x)\d x$. 
 	 
 	 \vspace{-2ex}
 	 
 	\item For $1\leq p\leq \infty$ we shall define the number $p'$ by $\frac{1}{p}+\frac{1}{p'} =1$ with the understanding that $p' = \infty$ if $p = 1$, and $p' = 1$ if $p = \infty$. 
 \item For $h\in \mathbb{R}^d $, $\tau_h$ denotes the shift function defined by $\tau_h u(x) = u(x+h)$ when $u$ is well understood. 

 \item 	We denote the support of a continuous function $u : \Omega\to \mathbb{R}^d $ by $\operatorname{supp} u = \overline{ \{x \in \Omega : u(x) \neq 0\}}.$
 %
 If $u$ is only measurable, then $\operatorname{supp} u =\R^d\setminus\mathcal{O}$ where $\mathcal{O}$ is largest open set on which $u$ vanishes. Ex. $u(x) =\mathds{1}_{\mathbb{Q}}(x) =0$ a.e. is not continuous but $ \operatorname{supp} u= \R$.
 \item For a multiindex $\alpha= (\alpha_1,\cdots, \alpha_d)\in \mathbb{N}_0$, we write $|\alpha|=\alpha_1+\cdots,+\alpha_d$ and $\partial^\alpha=\displaystyle \frac{\partial^{\alpha_d}\cdots\partial^{\alpha_1} }{ {\partial x_d^{\alpha_d} \cdots \partial x_1^{\alpha_1}}}$.
 
 \item Let $m\in \mathbb{N}$, the space $C_b^{m}(\Omega)$ is the collection of differentiable functions whose classical derivatives up to the order $ m$ are bounded.
 
 \item Let $m\in \mathbb{N}$ and $0<\sigma<1$. We define the H\"older space $C_b^{m+\sigma}(\Omega)$ also denote by $C_b^{m,\sigma}(\Omega)$ to be collection of functions in $C_b^{m}(\Omega)$ whose classical derivatives of order $|\alpha|=m$ belongs to $ C_b^{\sigma}(\Omega)$.The usual norm is the H\"older norm defined by 
 
 \vspace{-1cm}
 \begin{align*}
 \|u\|_{C^{m, \sigma }_b(\Omega)}= \sum_{|\alpha|\leq m-1}\sup_{x\in \Omega}|\partial^\alpha u(x)|+ \sum_{|\alpha|= m}\sup_{\stackrel{0<|x-y|<1}{x,y\in \Omega}}\frac{\big|\partial^\alpha u(x)-\partial^\alpha u(y)\big|}{|x-y|^\sigma}.
 \end{align*}
 \item The space $C_c^\infty(\Omega)$ is the collection of functions compactly supported in $\Omega$ whose classical derivatives $ \partial^\alpha u$ of every order exist. $C^m(\overline{\Omega})$ denotes the space  of  restrictions to $\overline{\Omega}$ of function of  $C^m(\mathbb{R}^d)$. 
 
 \item Given two comparable quantities $a$ and $b$ we denote $a\land b= \min(a,b)$, $a\lor b= \max(a,b)$  and $a\asymp b$ means there is a constant $C>0$ such that $C^{-1}a\leq b\leq C a$.
 \end{itemize}

\chapter{Basics On Nonlocal Operators }\label{chap:basics-nonlocal-operator}

It is the aim of this chapter to introduce some rudimentary notions on nonlocal operators. Let us recall that 
an operator $A$ defined on a function space is called local if it preserves the support, i.e. it obeys the rule that $\supp (Au) \subset \supp (u)$ for every function in the domain of $A$. Otherwise, $A$ will be called a nonlocal operator. Examples of local operators include, the Laplace operator $A=\Delta $, the gradient operator $A = \nabla $ and the divergence operator $A=\operatorname{div}$; each defined on the space $C^2(\mathbb{R}^d).$ The author admits that the title is sort of claptrap as it is conceived to attract attention. Indeed, in this chapter, we restrict ourselves on purely on nonlocal integrodifferential operators of L\'{e}vy type, which can be seen as a generalization of nonlocal operators acting on a smooth function $u\in C_c^\infty(\R^d)$ as follows
\begin{align*}
 L u(x) = \pv \int_{\mathbb{R}^d}(u(x)-u(y))\nu(x-y)\,\d y, \qquad (x\in \mathbb{R}^d)
\end{align*}
\noindent where the function $\nu:\mathbb{R}^d \setminus\{0\}\to [0,\infty] $ is even and satisfies the L\'{e}vy integrability condition that is $ \nu(-h)=\nu(h),~h\neq 0$ and $\nu\in L^1(\R^d,1\land|h|^2\d h)$.
 We will define the above operator from many different perspectives and later on, we provide a generalization of such operators viewed as integrodifferential operators. Their connections with elliptic differential operators of second order will be provided afterwards.

\section[Characterization of a L\'{e}vy operator]{Characterization of 
	a purely nonlocal symmetric L\'{e}vy operator} \label{sec:charac-levy-operator}

 We intend to define a purely nonlocal symmetric L\'{e}vy operator from several  perspectives. Our exposition here is certainly not exclusively original and is largely influenced by \cite{Mateusz2017} mainly treating the particular case of the fractional Laplace operator. We believe the content here expands upon this, however. Throughout this section, we assume that the function $\nu:\mathbb{R}^d \setminus\{0\}\to [0,\infty] $ is the density of a symmetric L\'{e}vy measure, i.e. $ \nu(h) = \nu(-h)$ and satisfies the integrability condition
\begin{align*}
 \int_{\mathbb{R}^d}(1\land |h|^2)\nu(h)\,\d h<\infty. 
\end{align*}
\noindent This integrability condition suggests that on the one hand  $\nu$ has some decay at infinity but is also allowed to have some singularity at the origin. 
Further generalization of such nonlocal operators will be considered in another section.

\medskip

\noindent\textbf{D.1: Singular integral} Given $\nu:\mathbb{R}^d \setminus\{0\}\to [0,\infty] $ as  the density of a symmetric L\'{e}vy measure we define a pure nonlocal symmetric L\'{e}vy operator $L$ acting on a smooth function $u:\mathbb{R}^d\to \mathbb{R}$ by
\begin{align}\label{eq:principal-val}
 L u(x) = \pv \int_{\mathbb{R}^d}(u(x)-u(y))\nu(x-y)\,\d y = \lim_{\varepsilon\to 0^+}L_\varepsilon u(x), \qquad (x\in \mathbb{R}^d),
\end{align}
with 
\begin{align*}
 L_\varepsilon u(x)= \il_{\mathbb{R}^d\setminus B_\varepsilon(x)}(u(x)-u(y))\nu(x-y)\,\d y\,.
\end{align*}
The notation $\pv $ stands for the common abbreviation of the Cauchy principal value. This makes sense within our context as the function $\nu$ might eventually have a singularity at the origin. 

\medskip 

\noindent \textbf{D.2: First order difference} One easily gets the following representations
\begin{align*}
 L u(x) &= \pv \int_{\mathbb{R}^d}(u(x)-u(x+h))\nu(h)\,\d h 
 = \pv  \int_{\mathbb{R}^d}(u(x)-u(x-h))\nu(h)\,\d h.
\end{align*}
\textbf{D.3: Second order difference} If $Lu(x)$ exists then we have, 
\begin{align}\label{eq:second-order-difference}
 L u(x) = -\frac{1}{2}\int_{\mathbb{R}^d}(u(x+h)-u(x-h)-2u(x))\nu(h)\,\d h .
\end{align}
We point out that the expression \eqref{eq:second-order-difference} may exist while the one in \eqref{eq:principal-val} does not. However, if \eqref{eq:principal-val} exists, then the above 
representations coincide. Indeed, for fixed $\varepsilon>0$, the mere change of variables $ y=x\pm h $ gives 
\begin{alignat*}{2}
 L_\varepsilon u(x)&= \il_{\mathbb{R}^d\setminus B_\varepsilon(x)}(u(x)-u(y))\nu(x-y)\,\d y
 &&= \il_{\mathbb{R}^d\setminus B_\varepsilon(0)}(u(x)-u(x+h))\nu(h)\,\d h\\
 &=\il_{\mathbb{R}^d\setminus B_\varepsilon(0)}(u(x)-u(x-h))\nu(h)\,\d h
 &&= -\frac{1}{2}\il_{\mathbb{R}^d\setminus B_\varepsilon(0)}(u(x+h)+u(x-h)-2u(x))\nu(h)\,\d h
\end{alignat*}
where the last equality is obtained by adding the two preceding ones. The below proposition explains why the second order difference allows us to get rid of the principal value. We show that $Lu(x)$ in \eqref{eq:second-order-difference} is well defined for bounded functions $u: \mathbb{R}^d\to \mathbb{R}$ which are sufficiently regular in a neighborhood of $x\in \mathbb{R}^d$. 
\medskip

\begin{proposition}\label{prop:well-defined}
Assume the function $u: \mathbb{R}^d\to \mathbb{R}$ is bounded and $C^2$ in a neighborhood of $x\in \mathbb{R}^d$. Then $L u(x) $ exists and one has
\begin{align*}
 L u(x)= -\frac{1}{2}\il_{\mathbb{R}^d}(u(x+h)+u(x-h)-2u(x))\nu(h)\,\d h.
\end{align*}
\end{proposition}

\medskip

\begin{proof} Assume $u$ is $C^2$ on a ball $B_{4\delta}(x)$ for $0<\delta<1$ small enough. For $|h|\leq \delta$ the fundamental Theorem of calculus yields, 
 \begin{align*}
 (u(x+h)+u(x-h)-2u(x)) &= \int_{0}^{1} \big[\nabla u(x+th) - \nabla u(x-th)\big]\cdot h\, \,\d t\\
 &= \int_{0}^{1} \int_{0}^{1} 2t \big[D^2 u(x-th + 2sth) \cdot h\big]\cdot h\, \,\d s\,\d t\, .
 \end{align*} 
Since $u$ is bounded on $\mathbb{R}^d $ and its Hessian $D^2 u $ on $B_{4\delta}(x)$, we get the estimate
 \begin{align}\label{eq:bound-second-difference}
 \left|u(x+h)+u(x-h)-2u(x)\right|\leq 4(\|u\|_{ C_b(\mathbb{R}^d)}+\|D^2u\|_{ C(B_{4\delta}(x) })(1 \land |h|^2).
 \end{align}
 Thereupon, it follows that for all $\varepsilon>0$ we have 
 \begin{align*}
 \frac{1}{2}\il_{\mathbb{R}^d\setminus B_\varepsilon(0)} \left|u(x+h)+u(x-h)-2u(x)\right|\nu(h)\,\d h\leq C\int_{\mathbb{R}^d} (1 \land |h|^2)\nu(h)\,\d h<\infty. 
 \end{align*}
\noindent This proves on one hand by the dominated convergence theorem that $L_\varepsilon u(x)$ converges to $Lu(x)$ as $\varepsilon \to 0^+$ and on the other hand that $Lu(x)$ exists and is worth the desired expression.
\end{proof}

\noindent Later on, we will show that the domain of $L$ can be extended to a larger space. For the moment let us evaluate the convergence of the family $(L_\varepsilon u)_\varepsilon$ to the operator $L u$ for smooth functions. 
\medskip

\begin{proposition}\label{prop:uniform-cont}
	Assume $\Omega \subset \mathbb{R}^d $ is an open bounded set and let $ u\in C^2_b(\mathbb{R}^d)$. The following properties then are satisfied.
	\begin{enumerate}[$(i)$]
		\item The map $x\mapsto L u(x)$ is uniformly continuous and bounded.

		\item For each $\varepsilon>0$, the map $x\mapsto L_\varepsilon u(x)$ is uniformly continuous.
		\item The family $(L_\varepsilon u(x))_\varepsilon$ is uniformly bounded and uniformly converges to $Lu$ as $\varepsilon\to 0$, i.e.
		\begin{align*}
		\|L_\eps u-Lu\|_{L^\infty(\mathbb{R}^d)}\xrightarrow[]{\eps \to 0}0.
		\end{align*}
	\end{enumerate}
\end{proposition}

\medskip

\begin{proof} 
	Here we use the second order difference representation \eqref{eq:second-order-difference}. Let $u \in C^2_c(\mathbb{R}^d)$, by a simple change of variables one gets
	\begin{align*}
	L_\varepsilon u(x)= -\frac{1}{2} \int\limits_{\mathbb{R}^d \setminus B_\varepsilon(0)}(u(x+h)+u(x-h)-2u(x)) \nu(h)\,\mathrm{d}h.
	\end{align*}
	Since $u$ and its Hessian $D^2 u $ are bounded functions, from \eqref{eq:bound-second-difference} we get the estimate
	\begin{align}\label{eq:second-difference}
	\left|u(x+h)+u(x-h)-2u(x)\right|\leq 4\|u\|_{ C^2_b(\mathbb{R}^d)}(1 \land |h|^2),\quad \quad x,h \in \mathbb{R}^d.
	\end{align}
	The integrability of the function $h\mapsto (1 \land |h|^2)\nu(h)$ entails the boundedness of $x\mapsto L u(x)$
	and the uniform boundedness of $x\mapsto L_\varepsilon u(x)$. Thereupon, it allows us to get rid of the principal value. We also get the uniform convergence of $(L_\varepsilon u)_\varepsilon$ to $Lu $ as follows
	\begin{align*}
	\|L_\eps u-Lu\|_{L^\infty(\mathbb{R}^d)} \leq 4\|u\|_{ C^2_b(\mathbb{R}^d)} \int_{B_\varepsilon(0)}(1 \land |h|^2)\nu(h)\,\d h \xrightarrow[]{\varepsilon\to0}0. 
	\end{align*}
	To prove the uniform continuity, we fix $x,z\in \mathbb{R}^d $ close enough, say $|x-z|\leq \delta$ with $0<\delta<1$,  then for all $h\in \mathbb{R}^d$ we have 
	
	\begin{align*}
	2|u(x)-u(z)| + |u(x+h)-u(z+h)|+ |u(x-h)-u(z-h)| \leq 4\delta \|u\|_{ C^2_b(\mathbb{R}^d)}.
	\end{align*} 
	This combined with \eqref{eq:second-difference} yields the uniform continuity via the integrability of $h\mapsto (1 \land |h|^2)\nu(h)$ as follows, 
	\begin{align*}
	\|L u(x)-L u(z) \|_{L^\infty(\mathbb{R}^d)} \leq 4\|u\|_{ C^2_b(\mathbb{R}^d)} \int_{\mathbb{R}^d}(\delta \land |h|^2)\nu(h)\,\d h \xrightarrow[]{\delta \to0}0. 
	\end{align*}
	The uniform continuity of $x\mapsto L_\varepsilon u(x)$ follows analogously. 
\end{proof}
\bigskip

\noindent It is possible to define $Lu$ with $u$ belonging to a space bigger than $C^2(\Omega)\cap L^\infty(\R^d)$. 
\begin{definition} The Blumenthal-Getoor index of $\nu$ is the real number $\beta_\nu \in [0,2]$ given by 
	$$\beta_\nu=\inf\{\gamma\geq0 :~\int_{B_1 (0)} |h|^\gamma\nu(h)\d h<\infty\}= \inf\{\gamma\geq 0 :~\int_{\R^d} (1\land |h|^\gamma)\nu(h)\d h<\infty\}. $$
	For  example, by considering $\nu(h) = |h|^{-d-2s}$ with $0<s<1$ one readily finds that $\beta_\nu= 2s$. 
	
\end{definition}

\begin{proposition}
	Let $0\leq \beta_\nu\leq 2 $ be the Blumenthal-Getoor index of $\nu$. Let $\delta>0$ then for every 
	$u\in C_b^{\beta_\nu+\delta}(\mathbb{R}^d)$. The map $x \mapsto L u(x) $ is well defined  and one has
	\begin{align*}
	L u(x)= -\frac{1}{2}\il_{\mathbb{R}^d}(u(x+h)+u(x-h)-2u(x))\nu(h)\,\d h.
	\end{align*}
	Moreover, the conclusions of Proposition \ref{prop:uniform-cont} remain true with the space $C_b^{2}(\mathbb{R}^d)$ replaces by $C_b^{\beta_\nu+\delta}(\mathbb{R}^d)$.
\end{proposition}

\begin{proof}
	Arguing on the two cases $0\leq \beta_\nu\leq 1$ and $1\leq \beta_\nu\leq 2$ then as for   the proof of \eqref{eq:bound-second-difference} one gets 
	\begin{align*}
	\left|u(x+h)+u(x-h)-2u(x)\right|\leq 4\|u\|_{ C_b^{\beta_\nu+\delta}(\mathbb{R}^d)}(1 \land |h|^{\beta_\nu+\delta}).
	\end{align*}
	The result follows.
\end{proof}

\noindent We point out that $\beta_\nu<2$ one can choose $\delta>0$ for which $\beta_\nu+\delta<2$. In this case we find that $C_b^{2}(\mathbb{R}^d) \subset C_b^{\beta_\nu+\delta}(\mathbb{R}^d)$.

\vspace{2mm}

\noindent \textbf{D.4: Integral form with gradient} Assume $u\in C^2(\mathbb{R}^d)\cap C_b(\mathbb{R}^d)$ just as for the preceding case, one has 
\begin{align*}
 Lu(x) = \il_{\mathbb{R}^d} [u(x)-u(x+h) -\mathds{1}_{B_1}(h)\nabla u(x)\cdot h]\,\nu(h)\,\d h,\qquad (x\in \mathbb{R}^d).
\end{align*}

\noindent This representation of $L$ is dear to the community from the area of probability analysis area as it represents the generator of a pure jump L\'{e}vy processes related to $\nu$. The expression on the right-hand side is justified since $h \mapsto [\nabla u(x)\cdot h]\nu(h) $ is odd and therefore has vanishing integral over $B_1\setminus B_\varepsilon $ for all $0<\varepsilon<1$.
Whence for all $0<\varepsilon<1$,
 \begin{align*}
 L_\varepsilon u(x) = \il_{\mathbb{R}^d\setminus B_\varepsilon } [u(x)-u(x+h) -\mathds{1}_{B_1}(h)\nabla u(x)\cdot h]\,\nu(h)\,\d h.
\end{align*}
\noindent As with the second-order difference \eqref{eq:second-order-difference}, one can easily get rid of the principal value 
since we have 
\begin{align*}
u(x)-u(x+h) -\mathds{1}_{B_1}(h)\nabla u(x)\cdot h &= \int_0^1[\nabla u(x+th)\cdot h -\nabla u(x)\cdot h]\,\d t\\
&= \int_0^1t \int_0^1 [D^2u(x+sth)\cdot h]\cdot h]\,\d s\,\d t
\end{align*}
for $|h|<1$ small enough. 
\bigskip

\noindent \textbf{D.5: Pseudo-differential operator} We now show that the integrodifferential operator $L$ can be realized as a pseudo-differential operator.  Let $\mathcal{S}^d(\mathbb{R}^d)$ denote the space of Schwartz functions on $\mathbb{R}^d$. Recall that the Fourier transform of a function $u\in \mathcal{S}^d(\mathbb{R}^d)$ is given by
\begin{align*}
 \widehat{u}(\xi) = (2\pi)^{-d/2}\int_{\mathbb{R}^d}e^{-i\xi \cdot x} u(x)\,\d x,\qquad(\xi\in \mathbb{R}^d). 
\end{align*}
The operator $L$ can be defined using the Fourier symbol as follows
\begin{align*}
 \widehat{Lu} (\xi) = \psi(\xi)\widehat{u}(\xi) \qquad(\xi\in \mathbb{R}^d)
\end{align*}
where, with the helping hand notation $m_\xi (h)=\cos{(\xi\cdot h)}$, the Fourier symbol $\psi(\xi)$ is given by
\begin{align*}
 \psi(\xi) = Lm_\xi (0)= \int_{\mathbb{R}^d}(1- \cos{(\xi\cdot h)})\nu(h)\,\d h.
\end{align*}
This relation is well defined since one readily observes that, for each $\xi\in \R^d$, the function $m_\xi (h)=\cos{(\xi\cdot h)}$ is in $C_b^2(\mathbb{R}^d)$. Hence, by Proposition \ref{prop:well-defined}, $Lm_\xi (0)= \psi(\xi)$ exists. Alternatively, we have $|1- \cos{(\xi\cdot h)}|\leq 1$ and by using the elementary inequality $|\sin t|\leq |t|$ true for all $t\in \mathbb{R}$, we get 
\begin{align*}
 |1- \cos{(\xi\cdot h)}|= |2\sin^2{\frac{\xi\cdot h}{2}}|\leq \frac{|\xi|^2|h|^2 }{2}.
\end{align*}
Thus, $|1- \cos{(\xi\cdot h)}|\leq (1+\frac{|\xi|^2 }{2}) (1\land |h|^2)$ which assures the existence of $\psi(\xi)$ since $\nu$ is L\'{e}vy integrable. Let us formally prove the above relation involving the Fourier transform. 
\begin{proposition}\label{prop:fourier-symbol}
For $u \in \mathcal{S}^d(\mathbb{R}^d)$ and a fixed frequency variable $\xi$ the following relation holds
\begin{align*}
 \widehat{Lu} (\xi) = \psi(\xi)\widehat{u}(\xi).
\end{align*}
with the Fourier symbol $\psi(\xi)$ given by
\begin{align*}
 \psi(\xi) = Lm_\xi (0)= \int_{\mathbb{R}^d}(1- \cos{(\xi\cdot h)})\nu(h)\,\d h.
\end{align*}
It is common to also name the Fourier symbol $\psi$ as the \textit{the characteristic exponent}. 
\end{proposition}

\bigskip

\begin{proof}
The key to the proof is to apply Fubini's theorem. 
To this, write 
$$|u(x+h)+ u(x-h) -2u(x)|\nu(h)= \Lambda_1(x,h)+ \Lambda_2(x,h)$$ for every $(x,h)\in \mathbb{R}^d\times \mathbb{R}^d $ with 
\begin{align*}
 \Lambda_1(x,h)&= \mathds{1}_{B_1}(h)\nu(h)|u(x+h)+ u(x-h) -2u(x)|,\\
 \Lambda_2(x,h)&= \mathds{1}_{B^c_1}(h)\nu(h)|u(x+h)+ u(x-h) -2u(x)|.
\end{align*}

\noindent First, observe that for each $h\in B_1$, 
 \begin{align*}
 \int_{\mathbb{R}^d}|u(x+h)+u(x-h)-2u(x)| \,\d x &= \int_{0}^{1} \int_{0}^{1} 2t \big[D^2 u(x-th + 2sth) \cdot h\big]\cdot h\, \,\d s\,\d t\,\\
 &\leq |h|^2 \int_{0}^{1} \int_{0}^{1} 2t \int_{\mathbb{R}^d} \big|D^2 u(x-th + 2sth)\big| \,\d s\,\d t\,\d x\\
 &=|h|^2 \int_{\mathbb{R}^d} \big|D^2 u(x)\big|\,\d x. 
 \end{align*}
Note that the mapping $x\mapsto |D^2u(x)|$ representing the norm in $ \mathbb{R}^{d^2}$ of the Hessian matrix of the function $u\in \mathcal{S}(\mathbb{R}^d)$ is integrable. 
Therefore, $\Lambda_1\in L^1(\mathbb{R}^d\times \mathbb{R}^d)$ since
\begin{align*}
 \iil_{\mathbb{R}^d \mathbb{R}^d} \Lambda_1(x,h) \,\d x\,\d h
 &=\int_{B_1}\nu(h) \int_{\mathbb{R}^d}|u(x+h)+u(x-h)-2u(x)| \,\d x\, \,\d h\\
 &\leq \int_{B_1}|h|^2\nu(h) \,\d h\cdot \int_{\mathbb{R}^d}\big|D^2 u(x)\big| \,\d x <\infty.
\end{align*}
Besides this, we have 
\begin{align*}
 \iil_{\mathbb{R}^d\mathbb{R}^d} \Lambda_2(x,h) \,\d x\,\d h
 &=\int_{B_1^c}\nu(h) \int_{\mathbb{R}^d}|u(x+h)+u(x-h)-2u(x)| \,\d x\, \,\d h\\
 &\leq 4\int_{B_1^c}\nu(h) \,\d h\cdot \int_{\mathbb{R}^d}|u(x)| \,\d x<\infty.
\end{align*}
 As a result, the function $(x,h)\mapsto \nu(h)|u(x+h)+ u(x-h) -2u(x)|$ is integrable on $\mathbb{R}^d\times \mathbb{R}^d$. Therefore, by using the translation rule of the Fourier transform $\widehat{u(\cdot+h )}(\xi) = \widehat{u}(\xi) e^{i\xi \cdot h}$ for all $\xi,h\in \mathbb{R}^d$ along with Fubini's theorem we get the desired result as follows
 \begin{align*}
 \widehat{Lu}(\xi) &= -\frac{1}{2} \,\d x\int_{\mathbb{R}^d} e^{-i\xi \cdot x} \int_{\mathbb{R}^d}(u(x+h)+u(x-h)-2u(x)) \nu(h)\,\,\d h\\
 &= -\frac{1}{2} \int_{\mathbb{R}^d} \nu(h)\,\d h\int_{\mathbb{R}^d}e^{-i\xi \cdot x} (u(x+h)+u(x-h)-2u(x)) \,\,\d x\\
 &=-\frac{1}{2} \widehat{u}(\xi)\int_{\mathbb{R}^d} (e^{i\xi \cdot h}+ e^{-i\xi \cdot h}-2) \nu(h)\,\d h\\
 &= \widehat{u}(\xi)\int_{\mathbb{R}^d} (1-\cos{(\xi \cdot h)}) \nu(h)\,\d h.
 \end{align*}
\end{proof}

\medskip

\begin{remark}
Assume that $\nu$ is radial. In this case we abuse the notation by writing $\nu(h)= \nu(|h|) $ for all $h\in \mathbb{R}^d\setminus\{0\}$. Let $O\in \mathcal{O}(d)$ be a rotation on $\mathbb{R}^d$ such that $O(|\xi|e_1) = \xi$ and enforce the change of variables $h'=O^Th$. Then $dh= dh'$ by the rotation invariance of the Lebesgue measure. Note that the rotation preserves the inner product and $O^T= O^{-1}$ hence $\xi\cdot h = O(|\xi|e_1) \cdot OO^T h= |\xi|e_1 \cdot O^T h = |\xi|e_1 \cdot h' $. In addition since $\nu$ is radial, we have $\nu(|h|) = \nu(|h'|)$. As a result, the Fourier symbol $\psi(\xi)$ becomes,
\begin{align}
 \psi(\xi)& = \int_{\mathbb{R}^d}(1- \cos{(\xi\cdot h)})\nu(|h|)\,\d h\notag \\
 &= \int_{\mathbb{R}^d}(1- \cos{(|\xi|e_1 \cdot h' )})\nu(|h'|)\,\d h'\notag\\
 &= \int_{\mathbb{R}^d}(1- \cos{(h_1})\nu(h/|\xi|)\frac{\,\d h}{|\xi|^{-d}}\quad(\xi\neq 0 )\label{eq:fourier-symbol-radial}\\
 &= \psi(|\xi|e_1)\notag.
\end{align}
\end{remark}

\vspace{2mm}

 \noindent\textbf{D.6: Generator of a symmetric L\'{e}vy process and of a semi-group} Assume that $\nu$ is radial, almost decreasing and $\int_{\R^d}\nu(h)\d h=\infty$. According to \cite[Lemma 2.5]{KR16} for each $t\geq0$, there exists a continuous function $p_t\geq0$ in $\R^d\setminus \{0\}$ such that in the Fourier space we have 
 \begin{align*}
 \widehat{p_t}(\xi)= \int_{\R^d}e^{-i\xi\cdot x}p_t(x)\d x = e^{-t\psi(\xi)},\qquad\qquad(\xi \in \R^d). 
 \end{align*}
 \noindent Obviously, $p_{t+s} = p_t p_s=p_sp_t$ for all $t,s\geq0$. Therefore, the family $(P_t)_t$ is defined by
 \begin{align*}
P_tu(x)= u*p_t= \int_{\R^d}u(y)p_t(x-y)\d y ,\qquad\qquad (x \in \R^d). 
\end{align*}
$(P_t)_t$ is clearly a strong continuous semigroup on $L^2(\R^d)$, whose generator is $-L$. We solely show that the operator $-L$ is the generator of $(P_t)_t$. Let $u\in \mathcal{S}(\R^d)$ we know that $\widehat{Lu} (\xi) = \widehat{u}(\xi) \psi(\xi) $. The Plancherel theorem implies,
\begin{align*}
\Big\|\frac{P_t u-u}{t} -(-Lu) \Big\|_{L^2(\R^d)} = \Big\|\frac{\widehat{p_t} \widehat{u}-\widehat{u}}{t} +\widehat{u}\psi \Big\|_{L^2(\R^d)} = \Big\|\widehat{u} \psi\frac{e^{-t\psi} -1+t\psi}{t\psi } \Big\|_{L^2(\R^d)}\xrightarrow{t\to 0}0.
\end{align*}
The rightmost term goes to $0$ as $t\to 0$ since the function $\zeta: s\mapsto \frac{e^{-s}-1+s}{s}$ with $\zeta(0)=0$ is continuous and bounded on $[0,\infty)$. Thus, applying the dominated convergence theorem suffices. 

\medskip 

\noindent 
The Kolmogorov extension theorem (see \cite{Sat13}) actually infers the existence of a stochastic process $(X_t)_t$ with the transition density is $p_t(x,y) = p_t(x-y)$, namely $\mathbb{P}^x(X_t\in A)=\mathbb{E}^x[\mathds{1}_A(X_t)]$. More generally 
\begin{align*}
\mathbb{E}^x[u(X_t)]= \int_{\R^d} u(y)p_t(x,y)\d y. 
\end{align*}
Here $\mathbb{P}^x$ (resp. $\mathbb{E}^x$) is the probability (resp. the expectation) corresponding to a process $(X_t)_t$ starting from the position $x$, i.e. $\mathbb{P}^x(X_0=x) =1$. 
The generator of such a stochastic process turns out to be $-L$. Indeed for a smooth function $u$,
\begin{align*}
\lim_{t\to 0}\frac{\mathbb{E}^x[u(X_t)]-u(x)}{t} =\lim_{t\to 0}\frac{P_t u(x)-u(x)}{t} =- Lu(x). 
\end{align*}
In fact, $(X_t)_t $ is a pure-jump isotropic unimodal L\'{e}vy process in $\R^d$, i.e.,  a stochastic process with stationary and independent increments and c\`{a}dl\`{a}g paths whose transition function $p_t(x)$ is isotropic and unimodal. We refer to \cite{Sat13} for a more extensive study on L\'{e}vy processes. 


\vspace{2mm}

\noindent\textbf{D.7: Energy form} We now show that the integrodifferential operator $L$ is intimately related with a Hilbert space of greatest interest in its own right. Let $H_\nu(\mathbb{R}^d)$ be the space of functions $u\in L^2(\mathbb{R}^d)$ such that $\mathcal{E}_{\mathbb{R}^d}(u,u)<\infty$ where we consider the bilinear form,
\begin{align*}
\mathcal{E}_{\mathbb{R}^d}(u,v) = \frac{1}{2} \iil_{\mathbb{R}^d \mathbb{R}^d}(u(x)-u(y))(v(x)-v(y))\nu(x-y)\,\d y\,\d x\qquad (u,v\in H_\nu(\mathbb{R}^d) ). 
\end{align*}

\noindent The space $H_\nu(\mathbb{R}^d)$ becomes a Hilbert space when furnished with the inner product
\begin{align*}
(u,v)_{H_\nu(\mathbb{R}^d)}:= (u,v)_{L^2(\mathbb{R}^d)}+\mathcal{E}_{\mathbb{R}^d}(u,v),\qquad u,v\in H_\nu(\mathbb{R}^d) .
\end{align*}
 Regarding the increments between the variables $x$ and $y$ involved in the integrand of $\mathcal{E}_{\mathbb{R}^d}(u,v)$ one legitimately suspects some close relation with integrodifferential operator $L$. This intuition is in fact correct, and  we show below that $\mathcal{E}_{\mathbb{R}^d}(\cdot, \cdot )$ can be viewed as the energy form associated with $L$. To do this, let us first observe that $\mathcal{S}(\mathbb{R}^d)\subset H_\nu(\mathbb{R}^d)$. Indeed, for $u\in \mathcal{S}(\mathbb{R}^d)$ we clearly
 have that $u, |\nabla u|\in L^2(\mathbb{R}^d) $ and arguing analogously as the proof of Proposition \ref{prop:fourier-symbol} ones readily arrives at
 \begin{align*}
 \hspace{-2ex}\iil_{\mathbb{R}^d \mathbb{R}^d}(u(x)-u(y))^2&\nu(x-y)\,\d y\,\d x 
\leq \int_{B_1}\hspace{-2ex}|h|^2\nu(h) \,\d h\cdot \int_{\mathbb{R}^d}\hspace{-1ex}|\nabla u(x)|^2 \,\d x + 2\int_{B_1^c}\hspace{-2ex}\nu(h) \,\d h\cdot \int_{\mathbb{R}^d}\hspace{-1ex}|u(x)|^2 \,\d x\hspace{-0.5ex}<\infty. 
 \end{align*}
 Noting that $|1-e^{-it}|^2 = 4\big|\tfrac{|e^{it/2}-e^{-it/2}}{2i}\big|^2= 4\sin^2\frac{t}{2} =2(1-\cos{t})$ for every $t\in \mathbb{R}$, it follows that for $u\in \mathcal{S}(\mathbb{R}^d)$, the Plancherel theorem yields,

 \begin{alignat*}{2}
\mathcal{E}_{\mathbb{R}^d}(u,u) &= \frac{1}{2} \iil_{\mathbb{R}^d \mathbb{R}^d}(u(x)-u(y))^2\nu(x-y)\,\d y\,\d x
& &= \frac{1}{2} \int_{\mathbb{R}^d} \nu(h) \,\d h\int_{\mathbb{R}^d}(u(x)-u(x+h))^2\,\d x\\
&= \frac{1}{2}\int_{\mathbb{R}^d} \nu(h) \,\d h\int_{\mathbb{R}^d} |\widehat{u}(\xi)|^2|1-e^{-i\xi \cdot h}|^2\,\d \xi
&&= \int_{\mathbb{R}^d} |\widehat{u}(\xi)|^2\int_{\mathbb{R}^d} (1-\cos{(\xi \cdot h)})\nu(h) \,\d h \,\d \xi\\
&= \int_{\mathbb{R}^d} |\widehat{u}(\xi)|^2 \psi(\xi) \,\d \xi =\int_{\mathbb{R}^d} \widehat{u}(\xi)\psi(\xi) \overline{\widehat{u}(\xi)} \,\d \xi&&= \int_{\mathbb{R}^d} \widehat{Lu}(\xi) \overline{\widehat{u}(\xi)} \,\d \xi.
\end{alignat*}

\noindent Here, the notation $\overline{z} $ denotes the conjugate of a complex number $z\in \mathbb{C}$. Employing the Plancherel theorem  again to the last expression produces the relation 
\begin{align*}
\mathcal{E}_{\mathbb{R}^d}(u,u)= \int_{\mathbb{R}^d} u(x) Lu(x) \,\d x.
\end{align*}
\noindent Replacing $u$ by $u+v$ leads to the relation,
\begin{align*}
\mathcal{E}_{\mathbb{R}^d}(u,v)= \int_{\mathbb{R}^d} v (x) Lu(x) \,\d x= \int_{\mathbb{R}^d} u(x) L v(x) \,\d x\qquad\text{for all} ~~u ,v \in \mathcal{S}(\mathbb{R}^d).
\end{align*}

\noindent Therefore, due to the density of $\mathcal{S}(\mathbb{R}^d)$, $Lu$ can be seen as a continuous linear form on $H_\nu(\mathbb{R}^d)$. Moreover,  through the dual pairing
 \begin{align*}
 \langle Lu, v\rangle= \mathcal{E}_{\mathbb{R}^d}(u,v)\qquad\text{for all} ~~ v \in H_\nu(\mathbb{R}^d). 
 \end{align*}

 \noindent The integrodifferential operator $L$ can be extended to functions $u$ in $ H_\nu(\mathbb{R}^d)$. Thereupon, $L$  can legitimately be regarded as a linear bounded operator from $H_\nu(\mathbb{R}^d)$ into its dual, i.e. $L: H_\nu(\mathbb{R}^d)\to \big(H_\nu(\mathbb{R}^d)\big)'$. In this case, we observe that $H_\nu(\mathbb{R}^d)$ is a fairly large domain for $L$. It is worthy nothing that through the correspondence $L: H_\nu(\mathbb{R}^d)\to \big(H_\nu(\mathbb{R}^d)\big)'$, $L$ may not always be evaluated in the classical sense. 
 
 \vspace{1mm}

\noindent\textbf{D.8: Distributions} It is natural to want to know  the largest possible functional space for which $Lu$ is defined.
 In an attempt to answer this question, assume in addition that $\nu:\mathbb{R}^d \setminus\{0\}\to [0,\infty] $ is unimodal, i.e. $\nu$ is radial and almost decreasing  and  there is a constant $c$ such that $\nu(y)\leq c \nu(x)$ whenever $|y|\geq |x|$. Let us define the function 
 $$\widehat{\nu}(x)=\nu(\tfrac{1}{2}(1+|x|)).$$ 
 It is not difficult to show that $\widehat{\nu}\in L^1(\mathbb{R}^d)\cap L^\infty(\mathbb{R}^d)$(c.f. the proof of Lemma \ref{lem:natural-norm-on-V}). Note that the space $L^1(\mathbb{R}^d,\widehat{\nu})$ is fairly large and contains the spaces $C_b^2(\mathbb{R}^d)$, $L^\infty(\mathbb{R}^d)$ and $ C^2(\mathbb{R}^d)\cap L^\infty(\mathbb{R}^d)$. Most importantly,  we show in  Chapter \ref{chap:nonlocal-sobolev} that $H_\nu(\R^d)\subset L^1(\mathbb{R}^d,\widehat{\nu})$.
 Furthermore, for $u\in C^2(\mathbb{R}^d)\cap L^1(\mathbb{R}^d,\widehat{\nu})$, we have that $Lu(x)$ exists in the classical sense for all $x\in \mathbb{R}^d$. Let us prove this formally.
 \begin{proposition}\label{prop:L1-def-L}
 	Assume  $\nu:\mathbb{R}^d \setminus\{0\}\to [0,\infty] $ is unimodal. Let $u\in C^2(\mathbb{R}^d)\cap L^1(\mathbb{R}^d,\widehat{\nu})$,  then $Lu (x)$ is well defined for all $x\in  \R^d$. 
 \end{proposition}
 \vspace{1mm}

 \begin{proof} For $x\in \R^d$, set $R=2|x|+1$, then for $y\in B^c_R(0)$ we have $|x-y|\geq \frac{|y|}{2}+ \frac{R}{2}-|x|\geq\frac{1}{2}(1+|y|). $
 	As $\nu$ is almost deceasing,  if $y\in B^c_R(0)$,  we have $\nu(x-y)\leq c\widehat{\nu}(y)$. Consequently, for $u \in L^1(\R^d, \widehat{\nu})$ we have 
 	\begin{align*}
 \Big|\int_{B^c_R(0)} (u(x)-u(y))\nu(x-y)\d y\Big|&\leq c\int_{\R^d} |u(x)-u(y)|\widehat{\nu}(y)\d y \\
 &\leq c|u(x)|\|\widehat{\nu}\|_{L^1(\R^d)}+ c\int_{B^c_R(0)} |u(y)\widehat{\nu}(y)\d y<\infty. 
 	\end{align*}
By exploiting \eqref{eq:second-difference} we get
 \begin{align*}
 \Big|\int_{B_R(0)} (u(x)-u(y))\nu(x-y)\d y\Big| &= \Big|\frac12 \int_{B_R(0)} (u(x+h)+ u(x-h)-2u(x))\nu(h)\d h\Big|\\
 &\leq 4\|u\|_{C^2(B_R(0))}\int_{\R^d} (1\land |h|^2)\d h<\infty. 
 \end{align*}
In conclusion, $Lu(x)$ exists since 
\begin{align*}
Lu(x)= \pv\int_{B_R(0)} (u(x)-u(y))\nu(x-y)\d y+ \int_{B^c_R(0)} (u(x)-u(y))\nu(x-y)\d y. 
\end{align*}
 \end{proof} 

\noindent A refinement of Proposition \ref{prop:L1-def-L}  for a possibly larger space, with an analog proof is given as follows: 
\begin{proposition}
	Assume  $\nu:\mathbb{R}^d \setminus\{0\}\to [0,\infty] $ is radial and almost decreasing. Let $\eps>0$ be sufficiently small and $\beta_\nu$ be the Blumenthal-Getoor index of $\nu$. Then, for $u\in C^{\beta_\nu+ \varepsilon}(\mathbb{R}^d)\cap L^1(\mathbb{R}^d,\widehat{\nu})$,  $Lu (x)$ is well defined for all $x\in  \R^d$. 
\end{proposition} 
\vspace{1mm}

\begin{remark} For the simple instance $\nu(h) = |h|^{-d-\alpha}$($\alpha\in (0,2)$) we have $\widehat{\nu}\asymp 1\land \nu\asymp (1+|h|)^{-d-\alpha}$. 
\end{remark}

\noindent We now deal with the situation where $u\in L^1(\mathbb{R}^d,\widehat{\nu})$ and show that $Lu$ is a  distribution. We recall that $\mathcal{D}(\R^d)$ is the space $C_c^\infty(\R^d)$ endowed its natural topology and a sequence $(\varphi_n)_n\subset \mathcal{D}(\R^d) $, 
 $\varphi_n\xrightarrow{n\to \infty}0$ in $ \mathcal{D}(\R^d)$ if and only if there exists a Compact set $K\subset \R^d$ such that $\supp \varphi_n\subset K$ for all $n\in \mathbb{N}_0$ and $\|\partial^\alpha \varphi_n\|_{C(K)}\xrightarrow{n\to \infty}0$ for all $\alpha\in \mathbb{N}^d_0$. The space $ \mathcal{D}'(\R^d)$ 
 whose elements are usually called \textit{distributions}, is the collection of continuous linear forms on $\mathcal{D}(\R^d)$. Furthermore, a sequence of linear forms $(T_n)_n\subset \mathcal{D}'(\R^d) $ converges to another one $T$ if and only if $ \langle T_n-T, \varphi \rangle \xrightarrow{n\to \infty}0$ for every $\varphi\in \mathcal{D}(\R^d) $. Here $\langle \cdot, \cdot \rangle$ is the dual pairing between $\mathcal{D}'(\R^d) $ and $\mathcal{D}(\R^d) $.

 \begin{proposition}\label{prop:distribution-L}
 For $u \in L^1(\mathbb{R}^d, \widehat{\nu})$, $Lu$ defines a distribution
 via the mapping $$\varphi\mapsto \langle Lu, \varphi \rangle := (u, L\varphi)_{L^2(\mathbb{R}^d)},\qquad\varphi\in C_c^\infty(\R^d).$$ 
 Moreover, by applying this procedure, the linear map $L:\, L^1(\mathbb{R}^d, \widehat{\nu})\to \mathcal{D}'(\R^d)$ with $u\mapsto Lu$ is continuous.
 \end{proposition}
 
 \medskip
 
 \begin{proof}
 Let $\varphi\in C_c^\infty(\mathbb{R}^d)$ be supported in $B_R(0)$ where $R\geq 1$ is sufficiently large.  
 First, we claim that 
\begin{align}\label{eq:estimate-test}
|L\varphi(x)|\leq C\|\varphi\|_{ C^2_b(\mathbb{R}^d)} \widehat{\nu}(x)\qquad \text{for all}~x\in \mathbb{R}^d.
\end{align}
Here the constant $C=C(R,d,\nu)$ depends only on $R,d$ and $\nu$. To be sure, suppose $|x|\geq 4R $, so that $\varphi(x)=0$. Since $|x-y|\geq \frac{|x|}{2}+ \frac{|x|}{2}-|y|\geq \frac{|x|}{2}+R\geq \frac{1}{2}(1+|x|)$ for $y \in B_R(0)$, the monotonicity of $\nu$ implies $\nu(x-y)\leq c \widehat{\nu}(x)$. Accordingly, 
\begin{align*}
|L\varphi(x) |\leq \int_{B_R(0)} |\varphi(y)|\nu(x-y)\,\d y \leq c|B_R(0)|\|\varphi\|_{C^2(\mathbb{R}^d)} \widehat{\nu}(x).
\end{align*}
 Whereas, if $|x|\leq 4R$ the proof of \eqref{eq:estimate-test} is completed by applying  \eqref{eq:second-difference} as follows: we have $\widehat{\nu}(x)\geq c_1$ for an appropriate constant $c_1>0$, depending on $R$ and $\nu$ since $\frac{1}{2}(1+ |x|)\leq 4R$. 
 \begin{align*}
 |L\varphi(x)|\leq 4\Theta \|\varphi\|_{C^2_b(\mathbb{R}^d)}\leq c_1^{-1}4\Theta \|\varphi\|_{C^2_b(\mathbb{R}^d)} \widehat{\nu}(x)
\end{align*} 
with $\Theta= \int_{\mathbb{R}^d} (1\land |h|^2)\nu(h)\,\d h$. Finally, \eqref{eq:estimate-test} yields 
\begin{align*}
|(u, L\varphi)_{L^2(\mathbb{R}^d)}|\leq C \|\varphi\|_{ C^2_b(\mathbb{R}^d)}\int_{\mathbb{R}^d}|u(x)|\widehat{\nu}(x)\,\d x.
\end{align*}
This spontaneously shows that $Lu $ is a distribution when $u\in L^1(\mathbb{R}^d, \widehat{\nu})$ and that $L : L^1(\mathbb{R}^d, \widehat{\nu})\to \mathcal{D}'(\R^d)$ is continuous, i.e. $Lu_n\xrightarrow{n\to \infty} Lu$ in $\mathcal{D}'(\R^d)$ if $u_n\xrightarrow{n\to \infty} u$ in $L^1(\mathbb{R}^d, \widehat{\nu})$. Indeed the above estimate implies
\begin{align*}
|\langle Lu_n-Lu, \varphi\rangle|\leq C \|\varphi\|_{ C^2_b(\mathbb{R}^d)}\int_{\mathbb{R}^d}|u_n(x)-u(x)|\widehat{\nu}(x)\,\d x\xrightarrow{n\to \infty}0\quad\text{for all}\,\, \varphi\in \mathcal{D}(\R^d).
\end{align*}
\end{proof}
\medskip

\noindent In light of Proposition \ref{prop:distribution-L}, we are forced to formulate the following definition. 
\begin{definition}\label{def:weak-integro}
A function $u\in L^1_{\operatorname{loc}}(\mathbb{R}^d)$ will be called to be weakly integrodifferentiable (or integrodifferentiable in the sense of distributions) with respect to $\nu$ on an open set $\Omega\subset \R^d$ (eventually $\Omega= \R^d$) if there exists a function $g\in L^1_{\operatorname{loc}}(\Omega) $ such that for every compact set $K\subset \Omega$ and for every $\varphi\in C_c^\infty(\mathbb{R}^d)$ with $\supp\varphi \subset K$,
\begin{align*}
\int_K g(x)\varphi(x)\,\d x=\int_K u(x)L\varphi(x)\,\d x. 
\end{align*}

\vspace{2mm}
 \noindent If $\nu$ is well understood, we briefly say that $Lu$ is the weak integrodifferential of $u$ or that $Lu$ is the integrodifferential of $u$ in the sense of distributions. In fact, as it can be easily shown that $ g$ is unique up to a set of measure zero, we shall merely write $g =Lu $ a.e on $\Omega$. 
\end{definition}

\medskip

\begin{remark}
Of course, if for some locally integrable function $u: \mathbb{R}^d\to \mathbb{R} $, the expression 
$Lu $ exists almost everywhere and belongs to $L^1_{\operatorname{loc}(\mathbb{R}^d)}$,  then $u$ is weakly integrodifferentiable and its weak integrodifferential coincides with $Lu$.
\end{remark}

\medskip
\noindent 
Let us  collect some simple facts involving the nonlocal operator $L$ under consideration.
\begin{proposition} \label{prop: properties-Levy-op} 
	The following assertions are true: 
\begin{enumerate}[$(i)$]
\item The operator $L$ commutes with translations. More generally, $L$ commutes with rigid motions\footnote{ A rigid motion is any transformation that can be obtained as a finite composition of translations and rotations.} if  in addition $\nu$ is radial. To be more precise, if $Lu(x)$ for all $x\in \mathbb{R}^d$ and $\tau :\mathbb{R}^d\to \mathbb{R}^d$ is a rigid motion then $[Lu]\circ \tau(x) = L[u \circ \tau](x) $ for all $x\in \mathbb{R}^d$.

\item  For any multiindex $\alpha\in \mathbb{N}_0^d$ we have $\partial^\alpha (Lu) = L\partial^\alpha u$ provided that $u$ and $Lu$ are  sufficiently smooth. 

\item  The convolution rule $L(u*\varphi)= Lu*\varphi$ holds for all $\varphi \in C^2_b(\mathbb{R}^d)$ and $u\in L^1(\mathbb{R}^d)$.

\item Assume $Lu$ is well defined in a distributional sense, then for all $\varphi \in C^2_b(\mathbb{R}^d)$ the equality $L(u*\varphi)= u*L\varphi$ holds as well in the distributional sense.

\item Given two functions $u,v: \mathbb{R}^d\to \mathbb{R}$, the relation $L[uv] =u[Lv] +[Lu]v- \Gamma(u,v)$ holds provided that $L[uv],$ $Lu$ and $Lv$ exist. Here $\Gamma(u,v)$ is the so called \textit{carr\'{e} du champs} operator associated with $L$ and is defined by 
\begin{align*}
\Gamma(u,v)(x) =\frac12 \int_{\R^d}(u(x)-u(y))(v(x) -v(y))\d y.
\end{align*}

\item Let $(\varphi_j)_{j\in \mathbb{N}}$ be a sequence of functions in $C_c^\infty(\mathbb{R}^d)$, such that for each $j\geq1$, $ \varphi_j=1$ on $B_j(0)$. Assume $u\in H_\nu(\R^d)$, such that $Lu\in L^12(\Omega)$, then we have  $\|L[\varphi_j u]- L u\|_{L^2(\R^d)}\xrightarrow{j\to \infty}0$. 
\end{enumerate}
\end{proposition}
\begin{proof} $(i)$ is a routine verification.  $(ii)$ follows by iterating the procedure from Proposition \ref{prop:uniform-cont}. $(iv)$ follows from $(iii)$. 
	To prove $(iii)$, observe that from  the estimate \eqref{eq:second-difference}, for all $x,h,z\in \R^d$ we have 
	\begin{align*}
	\big|u(z)\big(\varphi(x-z+h) + \varphi(x-z-h) -2\varphi(x-z) \big)\big|\nu(h)\leq C|u(z)|(1\land |h|^2)\nu(h)\in L^1(\R^d\times \R^d).
	\end{align*}
Thus, by applying Fubini's theorem, one easily arrives at $L(u*\varphi)= Lu*\varphi$. 
$(vi)$ is a consequence of $(v)$ whereas $(v)$ follows by integrating the identity
 $$u(x)v(x)-u(y)v(y) = u(x)(v(x)-v(y))- (u(x)-u(y))(v(x)-v(y)) + v(x)(u(x)-u(y)).$$ 
 
\end{proof}

\section{Case of the fractional Laplacian}

In this section we focus on the most studied integrodifferential operator which is of course interesting in its own right.  Before we define this object formally, let us observe that for $s\in \R$ the function $h\mapsto |h|^{-d-2s}$ is L\'{e}vy integrable if and only if $0<s<1$. Concretely,
\begin{align*}
	\int_{\mathbb{R}^d}(1\land |h|^2)|h|^{-d-2s}\d h<\infty\quad\text{if and only if}\quad 0<s<1.
\end{align*}
A prototypical example of a symmetric L\'{e}vy operator is obtained by putting, $\nu(h) =C_{d,s} |h|^{-d-2s}$  with $ s \in (0,1)$.  The resulting  L\'{e}vy operator is the so called fractional Laplace operator and is denoted by  $(-\Delta)^{s}$ or  $(-\Delta)^{\alpha/2}$ with $\alpha=2s$. The "$s$" notation "with a fractional tendency" usually suites to  PDEs  community whereas the probability community uses the "$\alpha$" notation as its remind the $\alpha$-stable processes.
Here, and in what follows depending on the context, we may simultaneously use the notation $\alpha=2s$. The constant $C_{d, s}\equiv C_{d, \alpha}$ is chosen so that the Fourier relation $\widehat{(-\Delta)^{\alpha/2} u}(\xi)= |\xi|^\alpha \widehat{u}(\xi)$ holds for all $u$ in $\mathcal{S}(\mathbb{R}^d)$. It follows from  \eqref{eq:fourier-symbol-radial} that  the Fourier symbol associated with $\nu(h)= C_{d,s}|h|^{-d-2s}$ is given by
\begin{align*}
\psi(\xi) = C_{d,s} \int_{\mathbb{R}^d}(1- \cos{(x_1})|x/|\xi||^{-d-2s}\frac{\,\d x}{|\xi|^{-d}} =C_{d,s} |\xi|^{-2s} \int_{\mathbb{R}^d} \frac{1-\cos(x_1)}{|x|^{d+2s}} \,dx.
\end{align*}
By identification it follows that 
\begin{align*}
C_{d,\alpha} \equiv C_{d,s}:= \left(\int_{\mathbb{R}^d} \frac{1-\cos(x_1)}{|x|^{d+2s}} \,dx\right)^{-1}.
\end{align*}
Note that within the Fourier symbol $\psi(\xi)$ we have already shown (see characterization D.5)  that $C_{d,s}$ is well defined. Later, we use a simple approach to compute this constant and study its asymptotic(c.f Section \ref{sec:cnst-cds}).  For the moment, it is important to keep in mind that 
\begin{align}\label{eq:pre-value-cds}
C_{d, \alpha}= \frac{2^\alpha \Gamma\big(\frac{d+\alpha}{2}\big)}{\pi^{d/2}\big|\Gamma\big(-\frac{\alpha}{2}\big)\big|}= \frac{2^{2s} \Gamma\big(s+\frac{d}{2}\big)}{\pi^{d/2}\big|\Gamma\big(-s\big)\big|}.
\end{align}

\noindent Patently from this Fourier characterization of $(-\Delta)^{\alpha/2}$, we can already glimpse that for $u\in \mathcal{S}(\R^d)$, $(-\Delta)^{\alpha/2} u(x) \xrightarrow{\alpha\to  2}  -\Delta u(x)$ and $(-\Delta)^{\alpha/2} u(x)\xrightarrow{\alpha\to  0} u(x)$. Here, we are reminded that $\xi \mapsto |\xi|^2$ corresponds to the Fourier multiplier of $-\Delta$.
This type of convergence extends to a more general context for continuous bounded functions that are $C^2$ in the vicinity of $x$ (see Section \ref{sec:nonlocal-elliptic}). It is worth emphasizing that the conclusions from Section \ref{sec:charac-levy-operator} apply to the fractional Laplacian $(-\Delta)^{\alpha/2}$. In particular with $\alpha=2s\in (0,2)$, from D.3  we have 
\begin{align}\label{eq:def-fractional-laplace}
\begin{split}
(-\Delta)^{\alpha/2}u(x) &:= C_{d,\alpha} \, \pv \int_{\R^d}\frac{(u(x)-u(y))}{|x-y|^{d+\alpha}}\d y\\
&= 
-\frac{C_{d,\alpha} }{2}\int_{\R^d}((u(x+h)+u(x-h)-2u(x) )\frac{\d h}{|h|^{d+\alpha}}.
\end{split}
\end{align}

\noindent The fractional operator $(-\Delta)^{\alpha/2}$ naturally appears as the generator of the rotationally symmetric $\alpha$-stable L\'{e}vy processes see \cite{Sat13}. More extensive works on the fractional Laplace operator can be found in \cite{Aba15,bucur2016,sil05,Sti19}. Let us proceed with some further representations of the fractional Laplace operator. Most of these representations are inspired by \cite{Mateusz2017}. 

\medskip

 \noindent\textbf{D.9: Inverse of Riesz’s potential} The fractional Laplacian can be realized as the inverse of the Riesz potential. The Riesz potential  of order $\alpha\in (0, d)$ which we denote as $I_\alpha$ is a pseudo-differential operator whose Fourier multiplier is $\xi\mapsto |\xi|^{-\alpha}$, i.e. for all $u\in \mathcal{S}(\R^d)$ we have  $\widehat{I_\alpha u}(\xi)= |\xi|^{-\alpha}\widehat{u}(\xi)$. It is easy to show that (see Remark \ref{rmk:rieszpotential-cnst}) for all  $u\in \mathcal{S}(\R^d)$ we have 
 \begin{align*}
 	I_\alpha u(x)=\frac{1}{\gamma_{d,\alpha}} \int_{\R^d}\frac{u(y)}{|x-y|^{d-\alpha}}\d y, \qquad\text{($x\in\R^d$)}, 
 \end{align*}
 \noindent where the constant $\gamma_{d,\alpha}$ (see Remark \ref{rmk:rieszpotential-cnst} for the computation) is given by 
 	\begin{align*}
 \gamma_{d,a}= \int_{\mathbb{R}^d}\frac{e^{-i x_1 }}{|x|^{d-a} }\d x= \pi^{\frac{d}{2}-a} \frac{\Gamma\left(\frac{d}{2}\right)}{\Gamma\left(\frac{d-a}{2}\right)}. 
 \end{align*}
 Using the Fourier characterization of $(-\Delta)^{\alpha/2}$ and $I_\alpha$ for $\alpha
 \in(0,2)$ we get  $(-\Delta)^{\alpha/2}\circ I_\alpha u=u$ and $ I_\alpha  \circ (-\Delta)^{\alpha/2} u =u$. That is,  $(-\Delta)^{\alpha/2}=I_\alpha^{-1}$ (the inverse of the operator $I_\alpha$). 
 
 \medskip 
 
 \noindent\textbf{D.10: Dynkin’s definition} For all $u\in C_b^2(\R^d)$ and for  $x\in \R^d$, we have 

\begin{align*}
(-\Delta)^{\alpha/2}u(x) &= C_{d,\alpha} \,\lim_{\eps\to 0^+}
\il_{ B^c_\eps(0)}\frac{(u(x)-u(x+h))}{|h|^2(|h|^2-\eps^2)^{\alpha/2}}\d h.
\end{align*}
The Dynkin definition of the fractional Laplacian is useful for studying  (see \cite{Aba15}) $\alpha$-harmonic functions, i.e. functions solving the equation $(-\Delta)^{\alpha/2} u=0$ on $\R^d$.  Let now us establish it. Fix, $0<\eps<1$. From the relation \eqref{eq:second-difference}, the following estimate holds 

\begin{align*}
\Big|\frac{(u(x+h)+u(x-h)-2u(x))}{|h|^2(|h|^2-\eps^2)^{\alpha/2}}\mathds{1}_{B^c_\eps(0)}(h)\Big|\leq 4\|u\|_{C_b^2(\R^d)}\frac{(1\land |h|^2)}{|h|^2\big||h|^2-1\big|^{\alpha/2}}.
\end{align*}
It is possible to show that the dominating function on the right side is integrable on $\R^d$. Thus the claim follows from the dominated convergence theorem, as we find  from \eqref{eq:def-fractional-laplace} that 

\begin{align*}
C_{d,\alpha} \,\lim_{\eps\to 0^+}\il_{ B^c_\eps(0)}\frac{(u(x)-u(x+h))}{|h|^2(|h|^2-\eps^2)^{\alpha/2}}\d h
&=-\frac{C_{d,\alpha} }{2}\,\lim_{\eps\to 0^+}\il_{ B^c_\eps(0)}\frac{(u(x+h)+u(x-h)-2u(x))}{|h|^2(|h|^2-\eps^2)^{\alpha/2}}\d h\\
&=-\frac{C_{d,\alpha} }{2}\,\int_{\R^d }\frac{(u(x+h)+u(x-h)-2u(x))}{|h|^{d+\alpha}}\d h= (-\Delta)^{\alpha/2} u(x). 
\end{align*}

\medskip

 \noindent\textbf{D.11: Bochner’s formula} The claim here is that for $\alpha\in (0,2)$ we have the identity $$B_\alpha u= (-\Delta)^{\alpha/2} u\qquad \text{for all $u\in \mathcal{S}(\R^d).$}$$ Here $B_\alpha$  is accomplished through the  following Bochner's integral formula
 \begin{align*}
 B_\alpha u(x): = \frac{1}{\big|\Gamma\big(-\frac{\alpha}{2}\big)\big|}\int_0^\infty \frac{(e^{-t\Delta} u(x)-u(x))}{t^{1+\alpha/2}}\d t.
 \end{align*}
 \noindent  For a rigorous proof of this identity (see for example \cite{ST10}) one would need some advance knowledge on Bochner integrals and functional calculus for unbounded operators.   Here we provide some intuitive approaches for  guessing the Bochner identity for $(-\Delta)^{\alpha/2}$. For $\lambda>0$, let us apply the integration by part to following the integral:  
\begin{align*}
\int_0^\infty(e^{-t\lambda}-1)t^{-\alpha/2-1}\d t&=  \lim_{\eps\to0^+}\Big[-\frac{2}{\alpha}(e^{-t\lambda}-1)t^{-\alpha/2}\Big]_\eps^\infty-\frac{2\lambda}{\alpha} \int_0^\infty e^{-t\lambda}t^{-\alpha/2}\d t\\
&=-\lambda^{-\alpha/2}\frac{2}{\alpha}\int_0^\infty e^{-s}s^{-\alpha/2-1}\d s=-\lambda^{-\alpha/2} \frac{2}{\alpha}\Gamma\big(1-\tfrac{\alpha}{2}\big). 
\end{align*}

\noindent Since $\alpha\in (0,2)$, by using  the duplication formula for the Gamma function (see the formula \eqref{eq:xAna-extension} ), we can  write, 
$\Gamma\big(1-\tfrac{\alpha}{2}\big)= -\tfrac{\alpha}{2}\Gamma\big(1-\tfrac{\alpha}{2}\big)= \tfrac{\alpha}{2}\big| \Gamma\big(-\tfrac{\alpha}{2}\big)\big|.$ Finally we have
 
\begin{align}\label{eq:bochner-formula-|ambda}
\lambda^{\alpha/2} u(x) = \frac{1}{\big|\Gamma\big(-\frac{\alpha}{2}\big)\big|}\int_0^\infty \frac{(e^{-t\lambda} -1)}{t^{1+\alpha/2}}\d t.
\end{align}

\noindent By roughly substituting $\lambda$ with $\Delta$ and $1$ with the identity operator, this gives 
\begin{align}\label{eq:bochner-formula}
(-\Delta)^{\alpha/2} u= \frac{1}{\big|\Gamma\big(-\frac{\alpha}{2}\big)\big|}\int_0^\infty \frac{(e^{-t\Delta} u -u)}{t^{1+\alpha/2}}\d t.
\end{align}

\noindent A more serious approach would  be to consider the Fourier multiplier. By knowing that, $\xi \mapsto|\xi|^2$ is the Fourier multiplier of $-\Delta$ it is not difficult to show that $\xi\mapsto e^{-t|\xi|^2}$ is the Fourier multiplier of $e^{-t\Delta}$. Therefore, using the formula \eqref{eq:bochner-formula-|ambda} we obtain
\begin{align*}
\widehat{B_\alpha u}(\xi)= \frac{\widehat{u}(\xi)}{\big|\Gamma\big(-\frac{\alpha}{2}\big)\big|}\int_0^\infty \frac{(e^{-t|\xi|^2} -1)}{t^{1+\alpha/2}}\d t= |\xi|^\alpha \widehat{u}(\xi)= \widehat{(-\Delta)^{\alpha/2} u}(\xi), \quad \text{for $u\in \mathcal{S}(\R^d)$ and $  \xi \in \R^d$}.
\end{align*}
The uniqueness of the Fourier transform, implies $B_\alpha u= (-\Delta)^{\alpha/2} u$, i.e. the identity \eqref{eq:bochner-formula} holds true.

 \vspace{2mm}
 
 \noindent\textbf{D.12: Balakrishnan’s formula} The idea here is the same as previously: fix $\lambda>0$, then using the Schwartz reflexion formula \eqref{eq:schwarz-reflexion} we get
 \begin{align*}
 \int_0^\infty\frac{\lambda t^{\alpha/2-1}}{t+\lambda}\d t= \lambda^{\alpha/2}\int_0^1 (1-t)^{-\alpha/2}t^{\alpha/2-1}\d t= \lambda^{\alpha/2} B\big(\tfrac{\alpha}{2}, 1-\tfrac{\alpha}{2}\big)=\frac{\lambda^{\alpha/2}\pi}{\sin\big(\frac{\alpha\pi}{2}\big)}.
 \end{align*}
 
 \noindent Therefore,  we have 
 
 \begin{align*}
 \lambda^{\alpha/2}=  \frac{\sin\big(\frac{\alpha\pi}{2}\big)}{\pi}\int_0^\infty\frac{\lambda t^{\alpha/2-1}}{t+\lambda}\d t.
 \end{align*}
 
\noindent Substituting $\lambda$ by $-\Delta $  yields the following Balakrishnan formula with $u\in \mathcal{S}(\R^d)$
\begin{align}
(-\Delta)^{\alpha/2}u=  \frac{\sin\big(\frac{\alpha\pi}{2}\big)}{\pi}\int_0^\infty-\Delta (tI-\Delta)^{-1} u\,\, t^{\alpha/2-1} \d t
\end{align}
\noindent It is worth noting that one is able to establish this formula using the Fourier multiplier as in D.11

\vspace{2mm}

\noindent\textbf{D.13: From long jump random walk} The aim here (see \cite{Va09}) is to show that the fractional Laplacian can be approximated by the generator of random walks. Define $K:\R^d\to [0,\infty]$ with   $K(h) = \theta_{d,\alpha}|h|^{-d-\alpha}$ for $h\neq 0$ and $K(0) =0$, where the constant $\theta_{d,\alpha}$ is chosen so that $$\sum\limits_{k\in \mathbb{Z}^d} K(k)=1.$$  For $\delta>0$, assume there is a random particle on the lattice $\delta\mathbb{Z}^d$ with the following properties. 

\begin{itemize}
	\item At any unit of time $\tau$ the particle jumps from any point of $\delta\mathbb{Z}^d$ to any other point. 
	
	\item The probability that the particle jumps from $x\in \delta\mathbb{Z}^d$ to $x+\delta k\in \delta\mathbb{Z}^d$ is given by $\theta_{d,\alpha}|k|^{-d-\alpha}$. This includes the fact that the particle can also experience a long jump with a small probability.

\item Let $u(x,t)$ be the probability that the particle sites at the point $x\in \delta\mathbb{Z}^d$ at the time $t\in \tau \mathbb{N}$.
\end{itemize} 

\noindent Patently, $u(x,t+\tau)$ (the probability that the particle sites at point $x$, at time $t+\tau$) corresponds to the probability that the particle sites at any other point $y=x+\delta k\in \delta\mathbb{Z}^d$ at time $t$ weighted with the jump rate  $\theta_{d,\alpha}|k|^{-d-\alpha}$. Together with the fact that $\sum\limits_{k\in \mathbb{Z}^d} \theta_{d,\alpha}|h|^{-d-\alpha}=1$, this yields
\begin{align}\label{eq:random-walk}
\begin{split}
u(x,t+\tau)-u(x,t)&= \theta_{d,\alpha}\sum_{k\in \mathbb{Z}^d}|k|^{-d-\alpha}u(x+\delta k,t) -u(x,t)\\
&= \theta_{d,\alpha}\sum_{k\in \mathbb{Z}^d}|k|^{-d-\alpha}\big(u(x+\delta k,t)-u(x,t)\big).
\end{split}
\end{align} 

\noindent Assume the unit of time is $\tau= \delta^\alpha$ then $\frac{|k|^{-d-\alpha}}{\tau}= \delta^d|\delta k|^{-d-\alpha}$. Combining  this with\eqref{eq:random-walk} and the Riemann sum of $h\mapsto (u(x+h,t) -u(x,t)) |h|^{-d-\alpha}$ over $\R^d$ yields the following 
\begin{align}
\partial_t u(x,t)&= \lim_{\tau\to 0}\frac{u(x,t+\tau)-u(x,t)}{\tau}\notag\\ 
&= \lim_{\delta\to 0}\delta^d \theta_{d,\alpha}\sum_{k\in \mathbb{Z}^d} |\delta k|^{-d-\alpha}\big(u(x+\delta k,t)-u(x,t)\big)\label{eq:riemann-sum}\\
&=\theta_{d,\alpha}\int_{\R^d}\frac{(u(x+h,t)-u(x,t))}{|h|^{d+\alpha}}\d h\label{eq:riemann-to-fractional}\\
&=-\theta_{d,\alpha}C^{-1}_{d,\alpha}(-\Delta)^{\alpha/2}u(x,t). \notag
\end{align} 

\noindent One should observe that the expression in \eqref{eq:riemann-sum} is the Riemann sum of \eqref{eq:riemann-to-fractional}. We have shown that $\partial_t u(x,t) =-\theta_{d,\alpha}C^{-1}_{d,\alpha}(-\Delta)^{\alpha/2}u(x,t)$. In other words, the generators of our random particle is $-(-\Delta)^{\alpha/2}$ up to a positive constant. 

\vspace{2mm}

\noindent\textbf{D.14: Caffarelli-Silvestre extension}
The Caffarelli-Silvestre  extension is probably the most skillful way to derive the fractional Laplacian. Indeed, Caffarelli and Silvestre showed in \cite{CS07} that the fractional  Laplacian can be determined as an operator that maps a Dirichlet boundary condition to a Neumann-type condition (Dirichlet-to-Neumann correspondence) via a PDE-extension problem. Concretely, for $u\in C_b^2(\R^d)$ and $\alpha\in (0,2)$, assume that a function $U:\R^d\times [0, \infty)\to \R$ satisfies the extension problem 
\begin{align}\label{eq:caffarelli-extension-problem}
	\begin{split}
\operatorname{div}(t^{1-\alpha}\nabla U) =0\quad\text{ in $\R^d\times [0,\infty)$}\quad\text{and}\quad  U(\cdot,0)= u\quad\text{on $\,\R^d$}. 
	\end{split}
\end{align}

\noindent Then,  for a constant $c_\alpha$ to be specified later,  we have 
\begin{align}\label{eq:caffarelli-extension}
\begin{split}
c_\alpha(-\Delta)^{\alpha/2} u(x)= -\frac{1}{\alpha}\lim_{t\to 0} t^{1-\alpha}\partial_tU(x,t) = -\lim_{t\to 0} \frac{U(x,t)-U(x,0)}{t^{\alpha}}. 
\end{split}
\end{align}

\noindent First, a routine check shows that we have $\operatorname{div}(t^{1-\alpha}\nabla )= t^{1-\alpha}\big(\Delta_x +\frac{1-\alpha}{t}\partial_t + \partial^2_{tt}\big)$, where $\Delta_x$ is the $d$-dimensional Laplacian with respect to the $x$ variable. We shall now introduce the Green function and the 
Poisson kernel associated with the operator $\operatorname{div}(t^{1-\alpha}\nabla ) $ which will be indispensable for establishing  the relations in \eqref{eq:caffarelli-extension}. 
\begin{definition}[\cite{CS07}]
The fundamental solution associated with the operator $\operatorname{div}(t^{1-\alpha}\nabla ) $ is the radial function $\Phi_{d,\alpha}:(\R^d\times \R)\setminus\{0\}\to (0,\infty)$ defined by 
\begin{align}
\Phi_{d,\alpha}(x)&:=-\frac{1}{2\pi}\ln\big(|x|^2+t^2\big)&&\text{ if $d+2-\alpha=2$\,, i.e. $(d=1,\alpha=1)$},\\\notag
\Phi_{d,\alpha}(x)&:=-\omega^{-1}_{d,\alpha}\big(|x|^2+ t^2\big)^{-\frac{d-\alpha}{2}}\, \,&& \text{ if $d+2-\alpha
>2$}.
\end{align}
 
 \noindent Assume that $d+2+\alpha>2$. The Poisson kernel associated with the operator $\operatorname{div}(t^{1-\alpha}\nabla ) $ is the function $P_{d,\alpha}:\R^d\times (0, \infty)\setminus\to (0,\infty)$ defined by 
 \begin{align}\label{eq:poisson-kernel}
 P_{d,\alpha}(x,t):=-t^{\alpha-1}\partial_t \Phi_{d,2-\alpha}(x,t)= \frac{t^\alpha (d-2+\alpha)  }{\omega_{d,2-\alpha}\big(|x|^2+ t^2\big)^{\frac{d+\alpha}{2}}}.
 \end{align}
 Note that the constant $\omega_{d,\alpha}$ is determined so that for every $t>0$ the function $P_{d,\alpha}(\cdot, t): \R^d\to(0, \infty)$ with $x\mapsto P_{d,\alpha}(x,t)$ has a unit mass, i.e. it satisfies the following integral condition 
 \begin{align}\label{eq:poisson-unit-mass}
 \int_{\R^d}P_{d,\alpha}(x,t)\d x=1.
 \end{align}
 
\end{definition}

\medskip

\begin{proposition}
Assume that $d+2-\alpha>2$ then we have 
\begin{align*}
\Phi_{d,\alpha}(x,t)&= \frac{1}{\omega_{d,\alpha}}\frac{ 1}{\big(|x|^2+ t^2\big)^{\frac{d-\alpha}{2}}} & &\quad\text{with}\quad \omega_{d,\alpha}= \frac{1}{(d-\alpha)}\frac{\Gamma\big(\frac{d+2-\alpha}{2}\big)}{\pi^{d/2}\Gamma\big(\frac{2-\alpha}{2}\big)},\\
P_{d,\alpha}(x,t)&= \frac{V_{d,\alpha}t^\alpha }{\big(|x|^2+ t^2\big)^{\frac{d+\alpha}{2}}}&& \quad\text{with}\quad V_{d,\alpha}= \frac{\Gamma\big(\frac{d+\alpha}{2}\big)}{\pi^{d/2}\Gamma\big(\frac{\alpha}{2}\big)}. 
\end{align*}
\end{proposition}

\begin{proof}
It is sufficient to compute $\omega_{d, 2-\alpha}$. The relation \eqref{eq:poisson-unit-mass} implies that 
\begin{align*}
	1= \omega^{-1}_{d,2-\alpha}(d+2-\alpha)\int_{\R^d}\frac{t^\alpha\d x}{\big(|x|^2+ t^2\big)^{\frac{d+\alpha}{2}}}&\overset{x=tz}{=} \omega^{-1}_{d,2-\alpha}(d-2+\alpha)\int_{\R^d}\frac{\d z}{\big(1+|z|^2\big)^{\frac{d+\alpha}{2}}}\\&= \omega^{-1}_{d,2-\alpha}(d-2+\alpha)A(d+1,\frac{\alpha-1}{2}).
\end{align*}

\noindent Here, $A(d+1,\frac{\alpha-1}{2})$ is given by the expressions \eqref{eq:Ans-expression} and \eqref{eq:value-Ans} below from which  the claim follows since,  
\begin{align*}
	 \omega_{d,2-\alpha}= \frac{1}{(d-2+\alpha)A(d+1,\frac{\alpha-1}{2})}= \frac{1}{(d-2+\alpha)}
	 \frac{\Gamma\big(\frac{d+\alpha}{2}\big)}{\pi^{d/2}\Gamma\big(\frac{\alpha}{2}\big)}.
\end{align*}

\end{proof}

\begin{remark} For the particular case $\alpha=1$, 	$\Phi_{d,1}$ and $	P_{d,1}$ are the  Green kernel  and the Poisson kernel of the Laplacian in the $d+1$ dimension space, respectively. Moreover, if $\omega_d$ denotes the area of the $d+1$-dimensional unit sphere of $\R^{d+1}$ then we have 
		\begin{align*}
		\frac{1}{\omega_{d,1}}= \frac{1}{(d-1)\omega_d}\quad 
		\text{with}\quad \omega_d= \frac{2\pi^{\frac{d+1}{2}}}{\Gamma\big(\tfrac{d+1}{2}\big)}.
		\end{align*}
	
\end{remark}

\medskip 

\noindent We omit  the computational details of the next proposition. Recall, $\operatorname{div}(t^{1-\alpha}\nabla )= t^{1-\alpha}\big(\Delta_x +\frac{1-\alpha}{t}\partial_t + \partial^2_{tt}\big)$. 
\begin{proposition}\label{prop:poisson-green} Let $\alpha\in (0,2)$ such that $d+2-\alpha>2$. The following assertions are true. 
	\begin{enumerate}[$(i)$] 
		\item For each $t>0$, $P_{d,\alpha}(\cdot, t) $ is radial positive,  satisfies $P_{d,\alpha}(\cdot, t) = t^{-d}P_{d,\alpha}(\frac{x}{t}, 1)$ and $\il_{\R^d}\hspace{-1ex}P_{d,\alpha}(x, t) \d x=1$.
		\item  The family $(P_{d,\alpha}(\cdot, t) )_t$ weakly converges to the Dirac mass $\delta_0$ as $t\to0$. To be more precise, 
		\begin{align*}
		\varphi(0)= \lim_{t\to0}\int_{\R^d} P_{d,\alpha}(x, t) \varphi(x)\d x\quad\text{ for all $\varphi\in C_c^\infty(\R^d)$}. 
		\end{align*} 
		\item $P_{d, \alpha}$ satisfies the equation $\operatorname{div}(t^{1-\alpha}\nabla)P_{d,\alpha}=0$ in $ \R^d\times (0,\infty)$.
		\item  $\Phi_{d, \alpha}$ satisfies the equation $\operatorname{div}(t^{1-\alpha}\nabla)\Phi_{d,\alpha}=\delta_0$ in $\R^{d+1}$.
	\end{enumerate}
\end{proposition}

\medskip 

\noindent We can now derive the fractional Laplace operator from the Dirichlet-to-Neumann map associated with the extension problem \eqref{eq:caffarelli-extension-problem}.
\begin{proposition}
Let $u\in C^2_b(\R^d)$ and $\alpha\in (0,2)$. Define the convolution $U(x,t) = [P_{d,\alpha}(\cdot, t)*u](x)$ then the following assertions are true. 

\begin{enumerate}[$(i)$]
	\item $U$ solves the equation \eqref{eq:caffarelli-extension} in the sense that we have $\operatorname{div}(t^{1-\alpha}\nabla) U(x,t) =0$ for every $x\in \R^d$ and every $t>0$ and 
	for every $x\in \R^d$ we have 
	\begin{align*}
	U(x,0)= \lim_{t\to0}[P_{d,\alpha}(\cdot, t)*u](x)= u(x). 
	\end{align*}
	\item Consider the constant $c_{\alpha}= \frac{V_{d,\alpha}}{C_{d,\alpha} }= \frac{\big|\Gamma\big(-\frac{\alpha}{2}\big)\big|}{2^\alpha\Gamma\big(\frac{\alpha}{2}\big)}$ then for every $x\in \R^d$ we have 
	\begin{align*}
	 -\lim_{t\to 0} t^{1-\alpha}\partial_t U(x,t)= -\frac{1}{\alpha}\lim_{t\to0}	\frac{U(x,t)- U(x,0)}{t^\alpha}
	 =  c_{\alpha}(-\Delta)^{\alpha/2}u(x).
	\end{align*}
	
	\item We define  $U^*(x,t) =U(x, \alpha t^{1/\alpha})$, for $x\in \R^d$ and $t>0$. Let $\beta= -\frac{2(1-\alpha)}{\alpha}$. Then $U^*$ verifies 
	\begin{align}\label{eq:cafferalli-extension-problem-bis}
		\big(\Delta_x+ t^\beta \partial^2_{tt}\big) U^*= 0\quad\text{ in \, $\R^d\times [0,\infty)$}\quad \text{and} \quad U^*(\cdot,0)= u\quad\text{on $\,\R^d$}. 
	\end{align}

\noindent Moreover, for all $x\in \R^d$ we have 
\begin{align*}
-\partial_t U^*(x,0)= -\lim_{t\to0}\frac{U^*(x,t)-U^*(x,0)}{t}= \alpha^{\alpha} c_{\alpha}(-\Delta)^{\alpha/2} u(x). 
\end{align*}

%
\end{enumerate}
\end{proposition}

\medskip 

\begin{proof} The differentiation rule under the integral sign and Proposition \ref{prop:poisson-green} $(iv)$ imply that  

$$\operatorname{div}(t^{1-\alpha}\nabla) U(x,t) = [\operatorname{div}(t^{1-\alpha}\nabla) P_{d,\alpha}(\cdot,t) ]*u(x) =0\quad\text{for all $x\in \R^d$ and all $t>0$}.$$ 
 Next, for $x\in \R^d$ and $0<t<1$, using the condition \eqref{eq:poisson-unit-mass}, we have 
	\begin{align*}
U(x,t)-u(x)=&	[P_{d,\alpha}(\cdot, t)*u](x)- u(x)
	= \int_{\R^d}(v(x-h) -v(x)) P_{d,\alpha}(h, t)\d h\\
	&=\frac12\int_{\R^d}(u(x+h) +u(x-h)-2u(x)) P_{d,\alpha}(h, t)\d h\\
	&= \frac{V_{d,\alpha}}{2}\int_{\R^d}(u(x+h) +u(x-h)-2u(x))\frac{t^\alpha\d h}{\big(|h|^2+ t^2\big)^{\frac{d+\alpha}{2}}}.
	\end{align*}
In view of the estimate \eqref{eq:second-difference}, we find that for every $h\in\R^d$,
\begin{align*}
	\Big|(u(x+h) +v(x-h)-2u(x))\frac{t^\alpha}{\big(|h|^2+ t^2\big)^{\frac{d+\alpha}{2}}}\Big|\leq 4\|u\|_{C^2_b(\R^d)} (1\land|h|^2)|h|^{-d-\alpha}. 
	\end{align*}
The dominant $h\mapsto (1\land|h|^2)|h|^{-d-\alpha}$ is integrable over $\R^d$. For each $h\in \R^d$ we have $P_{d,\alpha}(h,t)\xrightarrow{t\to 0}0$ and $t^{-\alpha}P_{d,\alpha}(h,t)\xrightarrow{t\to 0}|h|^{-d-\alpha}$. The convergence dominated theorem yields $$U(x,0)= \lim\limits_{t\to0}[P_{d,\alpha}(\cdot, t)*u](x)= u(x)$$ and 
 
\begin{align*}
\lim_{t\to0}	\frac{U(x,t)- U(x,0)}{t^\alpha}
&=\frac{V_{d,\alpha}}{2}  \int_{\R^d}(u(x+h) +u(x-h)-2u(x))\frac{\d h}{|h|^{d+\alpha}}
= -c_{\alpha}(-\Delta)^{\alpha/2}u(x).
\end{align*}
Here $c_{\alpha}=\frac{V_{d,\alpha}}{C_{d,\alpha}}$ is deduced from the expressions of $V_{d,\alpha}$ and  $C_{d,\alpha}$ (see\eqref{eq:pre-value-cds}). Furthermore, since 
$$t^{1-\alpha}\partial_t P_{d,\alpha} (h,t) =  \alpha V_{d,\alpha}\big(|h|^2+t^2\big)^{-\frac{d+\alpha}{2}}-t^2 V_{d,\alpha}(d+\alpha)  \big(|h|^2+t^2\big)^{-\frac{d+\alpha+2}{2}},$$
 for $0<t<1$, we also have the estimate

$$\big|(u(x+h) +v(x-h)-2u(x))t^{1-\alpha} \partial_t P_{d,\alpha} (h,t) \big|\leq 8\|u\|_{C^2_b(\R^d)} (d+\alpha)V_{d,\alpha}(1\land|h|^2)|h|^{-d-\alpha}.$$

\noindent and $t^{1-\alpha} \partial_t P_{d,\alpha} (h,t) \xrightarrow{t\to 0}\alpha|h|^{-d-\alpha}$,  for each $h\in \R^d$. Therefore, by using the differentiation rule under  the integral sign and the convergence dominated theorem  one finds that 
\begin{align*}
	\frac{1}{\alpha}\lim_{t\to 0}t^{1-\alpha}\partial_t U(x,t)
	&= \frac{1}{\alpha} \lim_{t\to 0} t^{1-\alpha}\partial_t \big(U(x,t)-u(x)\big)
=\frac{1}{\alpha} \lim_{t\to 0} [t^{1-\alpha} \partial_t P_{d,\alpha} (\cdot,t)*(u-u(x))](x) \\
	&=  \lim_{t\to 0} \frac{V_{d,\alpha}}{2\alpha}\int_{\R^d}(u(x+h) +u(x-h)-2u(x))\, t^{1-\alpha} \partial_t P_{d,\alpha}(h,t)\d h\\
	&=-c_{\alpha}(-\Delta)^{\alpha/2}u(x). 
\end{align*}
We have shown $(i)$ and $(ii)$, while $(iii)$ is not yet proven. Clearly, for $x\in \R^d$ we have $U^*(x,0)= U(x,0) =u(x)$. Furthermore, for $t>0$, we have  $\Delta_x U^*(x,t)= \Delta_x  U(x, \alpha t^{\frac{1}{\alpha}})$ and 

\begin{align*}
t^\beta \partial^2_{tt} U^*(x,t)= t^\beta\big(\partial_t (t^{\frac{1}{\alpha}-1} \partial_t U(x, \alpha t^{\frac{1}{\alpha}})) \big)= \partial^2_{tt} U(x,\alpha t^{\frac{1}{\alpha}})+ \frac{1-\alpha}{\alpha t^{\frac{1}{\alpha}}}\partial_t U(x, \alpha t^{\frac{1}{\alpha}}). 
\end{align*}

\noindent Recall that $\operatorname{div}(t^{1-\alpha}\nabla )= t^{1-\alpha}\big(\Delta_x +\frac{1-\alpha}{t}\partial_t + \partial^2_{tt}\big)$. Therefore, since $U$ solves \eqref{eq:caffarelli-extension-problem}, letting $z= \alpha t^{\frac{1}{\alpha}}$, we get 
\begin{align*}
\big(\Delta_x+ t^\beta \partial^2_{tt}\big) U^*(x,t)= \big(\Delta_xU(x,z) + \frac{1-\alpha}{z} \partial_t U(x,z)+ \partial^2_{tt} U(x,z\big) = 0. 
\end{align*}
\noindent Moreover, it follows from $(ii)$ that 
\begin{align*}
 \lim_{t\to0}\frac{U^*(x,t)-U^*(x,0)}{t}
 &=  \lim_{t\to0}\frac{U(x,\alpha t^{\frac{1}{\alpha}})-U(x,0)}{t} =\alpha^{\alpha} \lim_{z\to0}\frac{U(x,z)-U(x,0)}{z^\alpha} 
= -\alpha^{\alpha}c_{\alpha}(-\Delta)^{\alpha/2} u(x). 
\end{align*}
\end{proof}
%

%
%
%
%
 
\medskip 
 
%
 \section{Renormalization constant of the fractional Laplacian}\label{sec:cnst-cds}
 
 Here,  we provide a mere alternative method to compute the exact value of the constant $C_{d,s}$ with $0<s<1$  and $d \in\mathbb{N}$ given by
 \begin{align*}
 C_{d,s}:= \Big(\int_{\mathbb{R}^d} \frac{1-\cos(x_1)}{|x|^{d+2s}} \,dx\Big)^{-1}.
 \end{align*}
 Afterwards we  provide its asymptotic behavior when $s\to 1$ and $s\to 0$. 
 Although the constant $C_{d,s} $ has already been computed,  we believe that our approach is simpler and uses
 basic elementary calculus tools.  To the best of our  knowledge the computation of  a similar constant  first appeared in \cite[formula $(1.1.2)$]{Landkof} while studying the Riesz potential. 
 
 \vspace{2mm}
 
\noindent  This constant is nicely deduced  in  \cite[Chapter 1]{bucur-valdinoci-nonlocal}  through Fourier analysis  techniques while establishing a relationship between: (i) the function $U(x,t)$ solution of  the heat equation:  $ \partial_t U(x,t) = \Delta U(x,t) $  for all $(x,t)\in \mathbb{R}^d\times (0,+\infty)$ with initial data $U(x,0)=u(x)$  for some  $u\in \mathcal{S}(\mathbb{R}^d)$, (ii) its  natural semigroup  and (iii) the fractional  power of the associated infinitesimal generator which is nothing but the fractional Laplace operator (see D.11 the Bochner definition for $(-\Delta)^{\alpha/2}$ above ).  Although this  ingenuous idea requires  some further understanding
of the related subjects, the  identification process therein does not provide a direct computation of the integral quantity defining the constant $C_{d,s} $.

\vspace{2mm}
\noindent The art work  in \cite{bucur2016} consists of first showing that the function $u(x) =[\max(1-|x|^2,0)]^s$ is the solution to the Poisson equation $(-\Delta)^s u(x) = \theta_{d,s}$  in $B_1(0)$ and $u(x)= 0$ on $\R^d\setminus B_1(0)$. Here, $\theta_{d,s}$ is a specific constant from \cite{bucur2016}.  Afterwards the Green function $G_s(x,y)$ for the fractional Laplacian in the unit ball is exploited to deduce the value $C_{d,s}$ since one has $1=u(0) = \theta_{d,s}\int_{B_1} G_s(0,y)\d y$. 
 
 \vspace{2mm}
 \noindent In contrast, Mathieu Felsinger  \cite{Fel13} opts for a direct  computational method using polar coordinates. This method  brings into play the notion of Bessel functions of the first order $k\in \mathbb{R}$  the understanding of which requires some advanced  knowledge of integral calculus and the concept of analytic functions. The computational details  therein refer to  several advanced and  pre-established formulae for Bessel functions.  
 
 \vspace{2mm}
 \noindent Instead of these two appointed  techniques, we rather  propose a straightforward computation using basic elementary  methods of integral calculus.  Although the explicit value of $C_{d,s}$ is not of particular interest to
  us and only plays a minor role in our work,  having it at hand  significantly simplified the study of its  asymptotic behavior. Curious readers may look up \cite[Section 4]{Hitchhiker}  for  a lengthy elementary study of the asymptotic behavior of $C_{d,s}$.
 
\noindent To this scope,  we start  writing  $x=(x_1, x')\in \mathbb{R}^d$ where $x'\in\mathbb{R}^{d-1}$ assuming  $d \geq 2$. Thus, performing the change of variables $x'=z|x_1|$ on $\mathbb{R}^{d-1}$ that is $dx'=|x_1|^{d-1}dz$ yields 
 \begin{align}\label{eq:const-expression}
C_{d,s}^{-1}
= \left(\int_{\mathbb{R}} \frac{1-\cos(t)}{|t|^{1+2s}} \,dt\right)\left(  \int_{\mathbb{R}^{d-1}} \frac{1}{(1+|x|^2)^{\frac{d+2s}{2}}} \,dx \right): =\frac{G(s)A(d,s)}{s(1-s)}.
 \end{align}
 %
 With 
 \begin{align}
 G(s)&:= s(1-s)\int_{\mathbb{R}} \frac{1-\cos(t)}{|t|^{1+2s}}\,dt
 \end{align}
 and if $\omega_{d-2}$  denotes the surface measure of the $d-2$-dimensional unit sphere of $\mathbb{R}^{d-1}$  
 \begin{align}
 A(d,s)&:=\int_{\mathbb{R}^{d-1}} \frac{1}{(1+|x|^2)^{\frac{d+2s}{2}}} \,dx = \omega_{d-2}\int_{0}^\infty \frac{r^{d-2}}{(1+r^2)^{\frac{d+2s}{2}}} \,dr.\label{eq:Ans-expression}
 \end{align}
Note that from the expression of  $C(1,s)$ one can postulate the convention that $A(1,s)\equiv 1$ for all $s\in(0,1)$.  Let us recall the well known formula\footnote{This can be obtained by computing the Gaussian  integral $\int_{0}^{\infty} e^{-|x|^2}dx$ in two different 
	fashions:  first using polar coordinates  and second using Fubini 's
theorem observing that $|x|^2= x_1^2+\cdots +x_d^2$.} $\omega_{d-1}:  =  \frac{2\pi^{d/2}}{\Gamma(d/2)}$  where $\Gamma$ is  the  Euler's  Gamma function mapping $\alpha>0$ to 
 %
 %
 %
 %
 \begin{align*}
 \Gamma(\alpha):=\int_{0}^{\infty} e^{-t}t^{\alpha-1}\,dt.
 \end{align*}
 To compute the constant $A(d,s) $  we consider the following general integral
 \begin{align*}
 I(a,b) := \int_{0}^\infty \frac{r^{a-1}}{(1+r^2)^{\frac{b}{2}}} \,dr~~~0\leq a< b.
 \end{align*}
 The change of variables $ u = \frac{1}{1+r^2}$ that is $r=\left(\frac{1}{u}-1\right)^{1/2}$ and  $dr = -\frac{1}{2u^2}\left(\frac{1}{u}-1\right)^{-1/2}du $ yields,
 %
 %
 %
 %
 \begin{align*}
 	I(a,b)&=\frac{1}{2}\int_0^1 \left(1-u\right)^{\frac{a}{2}-1}u^{\frac{b}{2}-\frac{a}{2}-1}\,du = \frac{1}{2} B\left(\frac{a}{2},\frac{b}{2}-\frac{a}{2}\right) = \frac{\Gamma(\frac{a}{2})\Gamma(\frac{b}{2}-\frac{a}{2})}{2\Gamma(\frac{b}{2})}.
 \end{align*}
 Here $B(\cdot,\cdot)$ is the beta function defined for $x,y>0$ and given by the Legendre's duplication formula by
 \begin{align}\label{eq:beta-def}
 B(x,y):= \int_0^1 \left(1-u\right)^{x-1}u^{y-1}\,du = \frac{\Gamma(x)\Gamma(y)}{\Gamma(x+y)}.
 \end{align}
 %
 Then, taking into account that $\omega_{d-2}  = \frac{2\pi^{(d-1)/2}}{\Gamma(\frac{d-1}{2})}$ then from the expression of $I(a,b)$ and \eqref{eq:Ans-expression}, we get 
 \begin{align}
 A(d,s)&= \omega_{d-2}I(d-1, d+2s)=\pi^{(d-1)/2}\frac{\Gamma(s+\frac{1}{2})}{\Gamma(s+\frac{d}{2})}\label{eq:value-Ans}
 \end{align}
 On the other hand, by applying integration by parts, to the expression of $G(s)$ we get,
 
 \begin{align}
 G(s) &= 2s(1-s)\int_{0}^\infty \frac{1-\cos(t)}{t^{1+2s}}\,dt
 = (1-s)\int_{0}^\infty \frac{\sin(t)}{t^{2s}}\,dt\label{eq:value-Gs}.
 %
 %
 \end{align}
 Our main interest is to compute  the last integral term appearing in the expression of $G(s)$. To do this, we need to bring into play the integrals defined  for $\alpha\in\mathbb{R}$ by
 \begin{align}\label{eq:wphi-i}
 \varphi_1(\alpha) =\int_0^\infty \frac{\sin t}{t^\alpha}\,dt \qquad\hbox{and }\qquad \varphi_2(\alpha) =\int_0^\infty \frac{\cos t}{t^\alpha}\,dt.
 \end{align}
 %
 %
 
 \noindent A routine check reveals  that $\varphi_1(\alpha)$ exists if and only if $0<\alpha<2$ whereas, $\varphi_2(\alpha) $ exists if and only if $0<\alpha<1$. Moreover a more  precise  investigation  results in  the below table
 
 \begin{center}
 	\begin{tabular}{|c|c|c|}
 		\hline
 		summing over $(0,\infty)$ & Riemann integrable & Lebesgue integrable  \\ 
 		\hline 
 		$ \frac{\sin t}{t^\alpha} $ & $0<\alpha<2$ & $1<\alpha <2$ \\ 
 		\hline 
 		$\frac{\cos t}{t^\alpha} $ & $0<\alpha<1$ & no value \\ 
 		\hline 
 	\end{tabular}
 \end{center}
 Details of this are left to the interested reader.  Integrals in  \eqref{eq:wphi-i} are somehow linked to the so  called general Fresnel's\footnote{ For $\alpha=2  $ integrals in \eqref{eq:wfresnel-fi} are commonly  known as Fresnel's integrals.} integrals which are defined in their general forms for $\alpha\in\mathbb{R}$ as follows: 
 \begin{align}\label{eq:wfresnel-fi}
 f_1(\alpha ): =  \int_{0}^{\infty} \sin( t^\alpha)\,dt \qquad\hbox{and }\qquad f_2 (\alpha ) :=  \int_{0}^{\infty} \cos (t^\alpha)\,dt.
 \end{align}
 %
 %
 %
 In fact, if we assume $\alpha \neq 1 $ (observe that $f_i(1) =\varphi_i(0)$ does not exists), then enforcing the change of variables $ x= t^\alpha$, i.e $ dt =\frac{1}{\alpha} x^{\frac{1-\alpha}{\alpha}} dx$ we are led to  the following relationships
 %
 %
 %
 %
 \begin{align}\label{eq:f-fi}
 f_i(\alpha)=\frac{\text{sign}(\alpha)}{\alpha}\varphi_i\left(\frac{\alpha-1}{\alpha}\right)\qquad\text{equivalently,}\qquad \varphi_i(\alpha) = \alpha \text{sign}(\alpha) f_i\left(\frac{1}{1-\alpha}\right).
 \end{align}
 In contrast to the functions under the integrals in  \eqref{eq:wphi-i} the functions $t\mapsto \cos(t^\alpha)$ and $t\mapsto \sin(t^\alpha)$ have the advantage of being analytical  on the half complex plan  $\{z\in \mathbb{C}: \operatorname{Re}(z)\geq 0\}$ for $\alpha > 1$. It can  be deduced from the existence conditions on $\varphi_i$  and the above relationships that  $f_1(\alpha)$ exists only for $|\alpha|>1$ and $ f_2(\alpha)$ exists  only for $\alpha >1$.
 %
 %
 %
 Over all, the computation of $\varphi_i(\alpha)$ will be carried out by computing $f_i(\alpha)$ while the
 computation of this latter springs from Cauchy's theorem.
 \begin{theorem}[Cauchy Theorem]
 	Every analytic  function
 	on an open connected set $\Omega\subset \mathbb{C}$ has a null integral over any closed  oriented piecewise smooth simple curve supported in $\Omega$. 
 \end{theorem}
 \begin{proposition}
The values of $\varphi_1(\alpha)$ and  $\varphi_2(\alpha)$ are respectively given by
 	\begin{align}
 	\varphi_1(\alpha) =\int_0^\infty \frac{\sin  t}{t^\alpha}\,dt =\frac{\Gamma(2-\alpha)}{1-\alpha}\cos \left(\frac{\pi\alpha}{2}\right) ~~\alpha\in (0,2)
 	\end{align}
 	and 
 	\begin{align}
 	\varphi_2(\alpha) =\int_0^\infty \frac{\cos  t}{t^\alpha}\,dt=\frac{\Gamma(2-\alpha)}{1-\alpha}\sin \left(\frac{\pi\alpha}{2}\right) ~~\alpha\in (0,1).
 	\end{align}
 	Immediately, one gets 
 	\begin{align}\label{eq:express-Gs}
 	G(s)  = (1-s) \frac{\Gamma(2-2s)}{1-2s}\cos \left(s\pi\right).
 	\end{align}
 \end{proposition}
 
 \begin{proof}
 	First we assume  that $0<\alpha <1$, i.e. $\beta = \frac{1}{1-\alpha}>1$. Next, we introduce the function $z\mapsto e^{\mathrm{i}z^\beta}$ which is analytic on  $\{z\in \mathbb{C}: \operatorname{Re}(z)\geq 0\}$. 
 	Let us also introduce the close  counter-clockwise directed curve denoted by  $\Gamma_R$  and  given by  
 	$\Gamma_R = [O,A] \cup \gamma_R \cup [B, O]$ where $R>1$ is an arbitrarily positive large enough, $
 	\gamma_R = \{ Re^{\mathrm{i} t}:t\in [0,\frac{\pi}{2\beta}]\}$ is the arc of radius $R$ angle $\frac{\pi}{2\beta}
 	$,  $[B, O]=\{xe^{\mathrm{i} \frac{\pi}{2\beta}}, x\in [0,R]\}$ is the inclined segment of angle $\frac{\pi}
 	{2\beta}$ with length $R$  and $[O,A] =[0,R]$ is the segment on the real line of length $R$ (see the figure \ref{fig-curve}).
 	\begin{figure}[htbp!]
 		
 		\centering
 		\begin{tikzpicture}
 		\draw (0,0) node[below] {$O$};
 		\draw (6,0) node[below] {$A$};
 		\draw (6.01,1.25) node[right] {$\gamma_R$};
 		\draw (5.5,2.5) node[above] {$B$};
 		\draw (2.75,1.75) node[above] {$R$};
 		\draw [>=latex, ->] (5.5,2.5) -- (2.75,1.25);
 		\draw (0,0) -- (5.5,2.5) ;
 		\draw[>=latex,->] (0,0) -- (8,0) ;
 		\draw[>=latex,->] (0,0) -- (3,0) ;
 		\draw (3,0) node[below]{$R$} ;
 		\draw[>=latex,->] (0,0) -- (0,6);
 		\draw (6,0) arc (0:25:6) ;
 		\draw[>=latex,->] (6,0) arc (0:13:6) ;
 		\draw (1,0) arc (0:25:1) ;
 		\draw (1.5,0.04) node[above] {$\frac{\pi}{2\beta}$};
 		\end{tikzpicture}
 		\caption{Curve $\Gamma_R = [O,A] \cup \gamma_R \cup [B, O]$ } \label{fig-curve}
 	\end{figure}
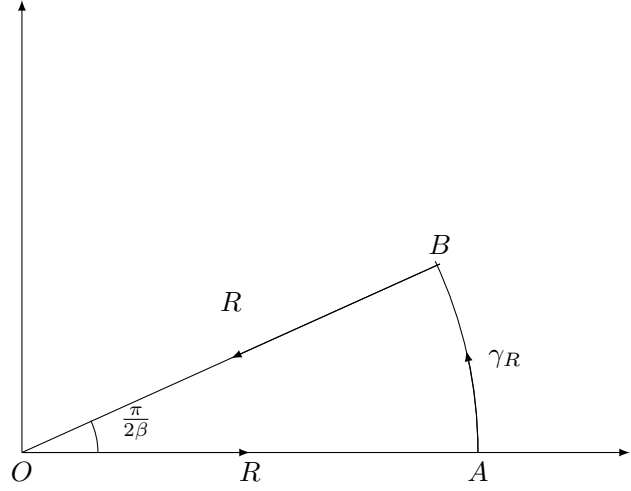

 	\noindent In virtue of  Cauchy's theorem, one has
 	\begin{align}
 	0=\int_{\Gamma_R} f(z) \,dz &= \int_{[O,A]} f(z) \,dz
 	+\int_{\gamma_R} f(z) \,dz
 	+\int_{[B,O]} f(z) \,dz\notag\\
 	&= \int_{0}^R f(x) \,dx 
 	+\int_{0}^{\frac{\pi}{2\beta}} f(Re^{\mathrm{i} t} ) \,d(Re^{\mathrm{i} t})
 	+\int_{R}^0 f( xe^{\mathrm{i} \frac{\pi}{2\beta}}) \,d(xe^{\mathrm{i} \frac{\pi}{2\beta}})\notag\\
 	&= \int_{0}^R e^{\mathrm{i}x^\beta} \,dx 
 	+ \mathrm{i}\int_{0}^{\frac{\pi}{2\beta}} e^{-R^\beta \sin \beta t} e^{\mathrm{i} R^\beta \cos  \beta t} Re^{\mathrm{i} t}\,dt
 	-e^{\mathrm{i} \frac{\pi}{2\beta}}\int_{0}^R e^{-x^\beta} \,dx.\label{eq:Gauss-relation}
 	\end{align}
 	Applying a simple change of variables along with the inequality $\sin t \geq \frac{2}{\pi} t $ for every $t\in [0,\frac{\pi}{2}]$, we have the following
 	\begin{align*}
 	&	\left|  \mathrm{i} \int_{0}^{\frac{\pi}{2\beta}} e^{-R^\beta \sin \beta t} e^{\mathrm{i} R^\beta \cos  \beta t} Re^{\mathrm{i} t}\,dt \right| 
 		%
 		\leq  \frac{R}{\beta}\int_{0}^{\frac{\pi}{2}} e^{-R^\beta \sin t} \,dt\\
 		&\leq  \frac{R}{\beta}\int_{0}^{\frac{\pi}{2}} e^{-R^\beta \frac{2t}{\pi}} \,dt
 		= \frac{\pi}{2\beta R^{\beta-1}}\left(1-e^{-R^\beta}\right)\xrightarrow{R\to\infty} 0 .
 	\end{align*}
 	Hence for $\beta>1 $ the integral of $f(z)$ over $\gamma_R$ vanishes as $R$ goes to infinity. So  that by letting $R\to \infty$ in \eqref{eq:Gauss-relation} after identification, we simultaneously get the formulae 
 	\begin{align*}
 		f_1(\beta ) &= \int_{0}^{\infty} \sin( t^\beta)\,dt
 		= \sin \frac{\pi}{2\beta}\int_0^{\infty} e^{-x^\beta}\,dx 
 		\overset{ t=x^\beta}{=} \frac{1}{\beta}\sin \frac{\pi}{2\beta} \int_0^{\infty}e^{-t} t^{\frac{1}{\beta}-1} \,dt = \frac{1}{\beta}\sin \frac{\pi}{2\beta}\Gamma\left(\frac{1}{\beta}\right) \\
 		f_2 (\beta ) &= \int_{0}^{\infty} \cos (t^\beta)\,dt 
 		=\cos \frac{\pi}{2\beta}\int_0^{\infty} e^{-x^\beta}\,dx  
 		\overset{ t=x^\beta}{=} \frac{1}{\beta} \cos \frac{\pi}{2\beta}\int_0^{\infty}e^{-t} t^{\frac{1}{\beta}-1} =\frac{1}{\beta}\cos \frac{\pi}{2\beta}\Gamma\left(\frac{1}{\beta}\right).
 	\end{align*}
 	This, together with the relation \eqref{eq:f-fi} gives the following for $ 0<\alpha<1,$
 	\begin{align*}
 		\varphi_1(\alpha) &=  \alpha f_1 \left(\frac{1}{1-\alpha} \right) 
 		= \Gamma(1-\alpha)\sin (1-\alpha) \frac{\pi}{2} \\
 		\varphi_2(\alpha)  &= \alpha f_2 \left(\frac{1}{1-\alpha} \right) = \Gamma(1-\alpha)\cos(1-\alpha) \frac{\pi}{2}.
 	\end{align*}

 \noindent For $\varphi_1(\alpha),$ the case $1<\alpha<2$ (i.e. $0<\alpha -1<1$) obviously springs from the previous case and the integration by part as follows: 
 	\begin{align*}
 		\varphi_1(\alpha)= \int_0^\infty \frac{\sin t}{t^\alpha}\,dt &= -\frac{1}{\alpha-1}\underbrace{\left[\frac{\sin t}{t^{\alpha-1}} \right]_0^\infty }_{0} + \frac{1}{\alpha-1} \int_0^\infty \frac{\cos t}{t^{\alpha-1}} = \frac{1}{\alpha-1}\varphi_2(\alpha-1) \\
 		&=  -\frac{\Gamma(2-\alpha)}{1-\alpha}\cos(2-\alpha) \frac{\pi}{2}=  \frac{\Gamma(2-\alpha)}{1-\alpha} \sin(1-\alpha) \frac{\pi}{2}.
 	\end{align*}
 	%
 	%
 	%
 	%
 	%
 	
 \noindent 	Finally, from the relation $\Gamma(x+1) = x\Gamma(x)~~x>0$, it can be deduced from the above steps  that
 	\begin{align}\label{eq:value-phi1}
 	\varphi_1(\alpha) &=\int_0^\infty \frac{\sin  t}{t^\alpha}\,dt 
 	= \frac{\Gamma(2-\alpha)}{1-\alpha}\sin (1-\alpha) \frac{\pi}{2} 
 	=\frac{\Gamma(2-\alpha)}{1-\alpha}\cos \left(\frac{\pi\alpha}{2}\right) ~~\alpha\in (0,2)\\
 	\label{eq:value-phi2}
 	\varphi_2(\alpha) &=\int_0^\infty \frac{\cos  t}{t^\alpha}\,dt = \frac{\Gamma(2-\alpha)}{1-\alpha}\cos (1-\alpha) \frac{\pi}{2}=\frac{\Gamma(2-\alpha)}{1-\alpha}\sin \left(\frac{\pi\alpha}{2}\right) ~~\alpha\in (0,1).
 	\end{align}
 	We deduce from \eqref{eq:value-Gs} and \eqref{eq:value-phi1}  that the expression of $G(s)$ is given by
 	\begin{align*}
 		G(s) &=  (1-s)\int_{0}^\infty \frac{\sin(t)}{t^{2s}}\,dt=(1-s)\varphi_1(2s) = (1-s) \frac{\Gamma(2-2s)}{1-2s}\cos \left(s\pi\right).
 		%
 	\end{align*}
 \end{proof}
 
 %
 
 \noindent We are now in a position to give the explicit value of $C_{d,s}$ and to study its asymptotic behavior as obtained in  \cite{Hitchhiker}. In addition, the Proposition below shows there is a  surprising coincidence between constants  $C_{d,s}$ and  the Brezis-Bourgain-Mironescu constant $K_{d,p}$(see the Section \ref{sec:charact-W1p}),  for $p=2$, which will be useful for studying the convergence of nonlocal structures to local ones.

 \begin{proposition} \label{prop:asymp-cds}
 	The following assertions are true: 
 	\begin{enumerate}[$(i)$]
 		\item Let $0<s< 1$,  we admit that $\Gamma(1-s)= -s\Gamma(-s) $ extending the definition of $\Gamma$ to $-s$. Then,
 		
 		%
 		\begin{align}\label{eq:explicit-consta-Cds}
 		C_{d,s}^{-1} = \pi^{(d-1)/2}\frac{\Gamma(s+\frac{1}{2})\Gamma(2-2s)}{s(1-2s)\Gamma(s+\frac{d}{2})} \cos \left(s\pi\right) =  \frac{\pi^{d/2}|\Gamma(-s)|}{2^{2s}\Gamma(\frac{d+2s}{2})} =  \frac{\pi^{d/2}\Gamma(1-s)}{s2^{2s}\Gamma(\frac{d+2s}{2})} .
 		\end{align}
 		Another simple expression is to set $\alpha= 2s$ so that we get 
 		\begin{align*}
 		 C_{d,s} \equiv C_{d, \alpha}= \frac{2^\alpha \Gamma\big(\frac{d+\alpha}{2}\big)}{\pi^{d/2}\big|\Gamma\big(-\frac{\alpha}{2}\big)\big|}.
 		\end{align*}
 		\item Let $1\leq p<\infty$ and $e$ be any element of  the unit sphere $\mathbb{S}^{d-1}$.  Then, 
 		\begin{align}\label{eq:explicit-const-Kdp}
 		K_{d,p}:= \fint_{\mathbb{S}^{d-1}} |w\cdot e|^p d\sigma_{d-1}(w) = \frac{ \Gamma\left(\frac{d}{2}\right) \Gamma\left(\frac{p+1}{2}\right)}{ \Gamma\left(\frac{1}{2}\right) \Gamma\left(\frac{p+d}{2}\right)}.
 		\end{align}  
 		\item Let $|\mathbb{S}^{d-1}|:=\omega_{d-1}= \frac{2\pi^{d/2}}{\Gamma(\frac{d}{2})}$ be  the surface measure of  unit sphere  $\mathbb{S}^{d-1}$  of $\mathbb{R}^d$. We have,
 		\begin{align}\label{eq:asymptotic-Cds}
 		\lim_{s\to 0^+} \frac{C_{d,s} }{s(1-s)} =\frac{2}{ \omega_{d-1}}\qquad \text{and }\qquad \frac{4}{\omega_{d-1} K_{d,2}} =  \lim_{s\to 1^-} \frac{C_{d,s} }{s(1-s)} =\frac{4d}{ \omega_{d-1}}.
 		\end{align}
 	\end{enumerate}
 \end{proposition}
 
 \bigskip 
 
 \begin{proof}
 		\noindent (i)  The first equality  clearly holds by combining \eqref{eq:const-expression}, \eqref{eq:value-Ans} and \eqref{eq:express-Gs}. For the second equality we need some helping hand formulae related to Gamma function.  Another version of Legendre's duplication formula stipulates  that for all $z>0$, 
 		\begin{align}\label{eq:Legendre}
 		\Gamma(2z) = \frac{2^{2z-1}}{\sqrt{\pi}}\Gamma(z)\Gamma(z+\frac{1}{2}).
 		\end{align} 
 		A justification of this equality comes from the original Legendre duplication formula in \eqref{eq:beta-def} as follows
 		\begin{alignat*}{2}
 			\frac{\Gamma(z)\Gamma(z)}{\Gamma(2z)} &= B(z,z)= \int_0^1 t^{z-1}(1-t)^{z-1}\,dt
 			&&\overset{t=\frac{x+1}{2}}{=} 2^{1-2z} \left[2\int_0^1 (1-x^2)^{z-1}\,dx\right]\\
 			&\overset{u=x^2}{=}  2^{1-2z} \int_0^1 u^{\frac{1}{2}-1}(1-u)^{z-1}\,du
 			&&= 2^{1-2z} B(\frac{1}{2},z) \\
 			&= 2^{1-2z}\frac{\Gamma(\frac{1}{2})\Gamma(z)}{\Gamma(z+\frac{1}{2})} 
 			=2^{1-2z}\frac{\sqrt{\pi}\Gamma(z)}{\Gamma(z+\frac{1}{2})}.
 		\end{alignat*} 
 		Let us also recall that the Schwarz reflexion formula\footnote{The Legendre duplication formula  implies $\Gamma(\varepsilon)\Gamma(1-\varepsilon)= B(\varepsilon, 1-\varepsilon) = \int_{0}^{\infty}\frac{t^\varepsilon dt}{1+t}$ thus the  Schwarz reflexion  formula can be obtained by  applying  the residues theorem using an appropriate domain.}  infers that for every $0<\varepsilon<1$, 
 		\begin{align}\label{eq:schwarz-reflexion}
 		B(\eps, 1-\eps)=\Gamma(\varepsilon)\Gamma(1-\varepsilon) =\frac{\pi}{\sin(\pi \varepsilon)}\qquad\qquad(\textrm{Schwarz's reflexion formula}).
 		\end{align}
 	Wherefrom, by using the parity of the $t\mapsto\sin t$ one obviously gets the relation 
 		\begin{align}\label{eq:xAna-extension}
 		\Gamma(1-\varepsilon)=-\varepsilon \Gamma(-\varepsilon)= \varepsilon |\Gamma(-\varepsilon)|.
 		\end{align}
 		More generally, the analytic extension formula of Gamma function given for all $m\in\mathbb{Z}$ and all $0<|\varepsilon| <1$ can be obtained by  induction:
 		\begin{align}\label{eq:xAna-extension-generalized}
 		\Gamma(\varepsilon-m)= (-1)^{m}\frac{\Gamma(1-\varepsilon)\Gamma(\varepsilon)}{\Gamma(m+1-\varepsilon)}.
 		\end{align}
 		Indeed, considering the extended Gamma function for negative value, Schwarz's reflexion formula above is also valid for $-1<z<0$.  By taking $m=1~ $ and replacing $~\varepsilon$ in \eqref{eq:xAna-extension-generalized}  by $1-\varepsilon$ with $0<\varepsilon <1$,  however one observes that
 		 $$\Gamma(-\varepsilon) = -\frac{\Gamma(1-\varepsilon)}{\varepsilon}\qquad\text{ that is} \qquad |\Gamma(-\varepsilon)| = \frac{\Gamma(1-\varepsilon)}{\varepsilon}.$$ 
 		Within the relation \eqref{eq:xAna-extension}, a routine check shows that for all $0<s<1$ with $s\neq \frac{1}{2}$  one has
 		\begin{align}\label{eq:w-Ana-extension-bis}
 		\Gamma\left(\frac{1}{2}+1-s\right) = \frac{1-2s}{2}\Gamma\left(\frac{1-2s}{2}\right).
 		\end{align} 
 		%
 		Thus, for $0<s<1$ from relations \eqref{eq:Legendre}, \eqref{eq:xAna-extension}, \eqref{eq:w-Ana-extension-bis} and \eqref{eq:schwarz-reflexion} we get

 		\begin{align*}
 			\Gamma(2(1-s)) &= \frac{2^{-2s+1}}{\sqrt{\pi}}\Gamma(1-s)\Gamma\left(\frac{1}{2}+1-s\right)\\
 			&= s(1-2s)\frac{2^{-2s}}{\sqrt{\pi}}|\Gamma(-s)|\Gamma\left(\frac{1-2s}{2}\right)\\
 			&=  s(1-2s)\frac{2^{-2s}}{\sqrt{\pi}}|\Gamma(-s)| \frac{\pi}{\sin(\pi (\frac{1-2s}{2}))\Gamma(s+\frac{1}{2})}\\
 			&= s(1-2s)2^{-2s}\sqrt{\pi} \frac{|\Gamma(-s)|}{\cos(s\pi)\Gamma(s+\frac{1}{2})}.
 		\end{align*}
 		This shows the second equality in  \eqref{eq:explicit-consta-Cds} for $s\neq \frac{1}{2}$  for which the result is also valid for $s=\frac{1}{2}$ simply by taking the limit as $ s\to \frac{1}{2}$ directly in \eqref{eq:explicit-consta-Cds}.

%
 		\vspace{2mm}
 		
 	\noindent (ii) Taking into account the rotation invariance of the Lebesgue measure one glimpses that $ K_{d,p}$ is  independent of the chosen unit vector $e$ and whence it is sufficient to consider $e:=(0, \cdots, 0,1)\in \mathbb{S}^{d-1}$. Now we let $w = (w', t)\in \mathbb{S}^{d-1} $ with $t\in (-1,1)$ so that $w'\in \sqrt{1-t^2}\mathbb{S}^{d-2}$. In virtue of the Jacobian formula for spherical coordinates one has $ \d \sigma_{d-1}(w)=\frac{ \d\sigma_{d-2}(w')dt}{\sqrt{1-t^2}}$ (see \cite[Appendix D.2]{grafakos04}). Therefore, 
 		\begin{alignat*}{2}
 		K_{d,p}&=  \fint_{\mathbb{S}^{d-1}} |w_1|^pd\sigma_{d-1}(w) 
 		&&=\frac{1}{\omega_{d-1}} \int_{-1}^{1}\int_{w'\in \sqrt{1-t^2}\mathbb{S}^{d-2}} |t|^pd
 		\frac{ \sigma_{d-2}(w')dt}{\sqrt{1-t^2}} \\
 		&=\frac{2}{\omega_{d-1}}\int_{0}^{1} t^p\left|\sqrt{1-t^2}\mathbb{S}^{d-2}\right| 
 		\tfrac{ dt}{\sqrt{1-t^2}}
 		&&=\frac{2\omega_{d-2}}{\omega_{d-1}} \int_{0}^{1} (1-t^2)^{\frac{d-3}{2}} t^p dt \\
 		&=\frac{\omega_{d-2}}{\omega_{d-1}} \int_{0}^{1} (1-t)^{\frac{d-1}{2}-1} t^{\frac{p+1}{2}-1} dt 
 		&&= \frac{\omega_{d-2}}{\omega_{d-1}} B\left( \frac{d-1}{2}, \frac{p+1}{2}\right)\\
 		&= \frac{\omega_{d-2}}{\omega_{d-1}}\frac{\Gamma\left(\frac{d-1}{2}\right)\Gamma\left(\frac{p+1}{2}\right)}{\Gamma\left(\frac{d+p}{2}\right)}.
 		\end{alignat*}
 		Here $B(\cdot, \cdot)$ is the beta function and we have used the relation \eqref{eq:beta-def}. The result follows by using the formula $|\mathbb{S}^{d-1}|= \omega_{d-1}= \frac{2\pi^{d/2}}{\Gamma\big(d/2\big)}$ along with the relation $\Gamma(\frac{1}{2})  = \pi^{1/2}$.
 		
 		\vspace{2mm}
 		\noindent (iii) This result plainly follows from expressions of $G(s)$ and $A(d,s)$. Indeed, one can check  that 
 		\begin{align*}
 		G(s) =(1-s) \Gamma(1-2s)\cos \left(s\pi\right)= \frac{1}{2}\frac{\Gamma(3-2s)}{1-2s}\cos \left(s\pi\right) 
 		\end{align*}
 		from where,
 		\begin{align*}
 		\lim_{s\to 0} G(s) = 1 \qquad \text{and}\qquad \lim_{s\to 1} G(s) = \frac{1}{2}.
 		\end{align*}
 		Since $\Gamma(\frac{1}{2})  = \pi^{1/2}$, from  \eqref{eq:Ans-expression} we get the following   after cancellation 
 		\begin{align*}
 		A(d,0) = \frac{\pi^{d/2}}{\Gamma(\frac{d}{2})} =\frac{ \omega_{d-1}}{2}\quad \text{and }\quad A(d,1) = \frac{\pi^{d/2}}{d\Gamma(\frac{d}{2})} =\frac{ \omega_{d-1}}{2d}.
 		\end{align*} 
 		Therefore, as $s(1-s)C_{d,s}^{-1}=G(s)A(d,s)$ and $K_{d,2} =\frac{1}d$,  we obtain the following results:
 		\begin{align*}
 			\lim_{s\to 0} \frac{C_{d,s} }{s(1-s)} =\frac{2}{ \omega_{d-1}}\qquad \text{and }\qquad \lim_{s\to 1} \frac{C_{d,s} }{s(1-s)} =\frac{4d}{ \omega_{d-1}} = \frac{4}{\omega_{d-1} K_{d,2}}.
 		\end{align*}
 		Alternatively, the above limits follows after one easily verifies that  	
 		\begin{align*}
 	\frac{s(1-s)}{C_{d,s}} = \frac{\pi^{d/2}\Gamma(2-s)}{2^{2s}\Gamma(s+\frac{d}{2})} .
 		\end{align*}
 	
 \end{proof}

 \begin{remark}\label{rmk:rieszpotential-cnst}
 	For $0<a<1$, from the foregoing computations, it follows that the constant
 	\begin{align*}
 	\gamma_{d,a}=  2\phi_2(1-a)A(d, -a/2)
 	= \pi^{\frac{d}{2}-a} \frac{\Gamma\left(\frac{d}{2}\right)}{\Gamma\left(\frac{d-a}{2}\right)}
 	\end{align*}
 	is the constant for which the Riesz Potential  (also known as the inverse of the fractional Laplacian)
 	\begin{align*}
 	I_{a}f(x):=\frac{1}{\gamma_{d,a}}\int_{\mathbb{R}^d} \frac{f(y) dy }{|x-y|^{d-a}}
 	\end{align*} 
 	satisfies the relation  $\widehat{I_af}(\xi) = |\xi|^{-a} \widehat{f}(\xi)$ for all $f\in S(\mathbb{R}^d)$ with equality understood in the sense of distributions. Indeed, the invariance of the Lebesgue measure under rotation implies 
 	 \begin{align*}
 	 \widehat{I_af}(\xi) = \frac{1}{\gamma_{d,a}} \widehat{f}(\xi)  \int_{\mathbb{R}^d}\frac{e^{-ix\cdot \xi}}{|x|^{d-a} }dx= \frac{|\xi|^{-a}}{\gamma_{d,a}} \widehat{f}(\xi)  \int_{\mathbb{R}^d}\frac{e^{-i x_1 }}{|x|^{d-a} }dx.
 	 \end{align*}
 Letting $x'= x_1z$  with $z\in \mathbb{R}^{d-1}$ that is $dx' = x_1^{d-1}dz$
 	 
 	 \begin{align*}
 	  \int_{\mathbb{R}^d}\frac{e^{-i x_1 }}{|x|^{d-a} }dx =  \left(\int_{\mathbb{R}}\frac{\cos( t )}{t^{1-a} }dt \right)\left(\int_{\mathbb{R}^{d-1}}\frac{dz}{(1+|z|^2)^{\frac{d-a}{2}} }\right)= 2\varphi_2(1-a)A(d, -a/2).
 	 \end{align*} 
 \end{remark}

 \section{Order of the fractional Laplace operator}
Now, we show that the fractional Laplacian acts on H\"older spaces like an operator of order $2s$. 
   \begin{definition}
 	Let $m\in \mathbb{N}$ and $0<\sigma<1$. The H\"older space $C_b^{m+\sigma}(\mathbb{R}^d)$ also denoted, 
 	 by $C_b^{m,\sigma}(\mathbb{R}^d)$, is the collection of functions in $C_b^{m}(\mathbb{R}^d)$ whose classical derivative of order $|\alpha|=m$ belongs to $ C_b^{\alpha}(\mathbb{R}^d)$.
  \end{definition}
\noindent In the spirit of the estimate \eqref{eq:second-difference}, one can easily establish the following estimates.
 \begin{proposition}\label{prop:estimate-second-order}
 	Let $u \in C_b^{\gamma}(\mathbb{R}^d)$ with $\gamma= 1+\tau$, $m=1$ and $0\leq \tau \leq 1$. 	Define the second order discrete difference of $u$ by  $D^2_hu(x) =u(x+h)+u(x-h)-2u(x).$ Then for all $x,h,z \in \mathbb{R}^d$ we have 
 	\begin{align*}
 	|D^2_h u(x+z)-D^2_h u(x)|\leq C\|u\|_{C_b^{\gamma}(\mathbb{R}^d)}(|z||h|^\tau \land |h|^{1+\tau}).
 	\end{align*}
 	Furthermore, if $u \in C_b^{\gamma}(\mathbb{R}^d)$ with $0\leq \gamma \leq 1$ then
 	\begin{align*}
 	|D^2_h u(x+z)-D^2_h u(x)|\leq C\|u\|_{C_b^{\gamma}(\mathbb{R}^d)}(|z|^\gamma \land |h|^{\gamma}).
 	\end{align*}
 \end{proposition}

 \bigskip 
 
 \noindent The next result shows that the fractional Laplace operator $(-\Delta)^s$ is an integrodifferential of order $2s$.
 \begin{theorem}[Order of fractional Laplacian]
 	Let $\gamma>0$ and $0<s<1$ such that $0<2s<\gamma$. The fractional Laplacian defined $(-\Delta)^s : C_b^{\gamma}(\mathbb{R}^d) \to C_b^{\gamma-2s}(\mathbb{R}^d)$ is a bounded operator.
 	
 \end{theorem}

 \medskip
 
 \begin{proof}
 	First, assume $\gamma -2s= m+ \sigma$ with $m=0,1$ and $0<\sigma\leq 1$. In the case $m= 0$, i.e. $\gamma= 2s+\sigma$ there are subcases: $0<2s<\gamma \leq 1$ or $1<2s<\gamma \leq 2$. Suppose that $0<2s<\gamma \leq 1$ then in view of Proposition \ref{prop:estimate-second-order} we have 
 	\begin{align*}
 	|(-\Delta)^su(x+z)-(-\Delta)^su(x)|&=\frac{C_{d,s}}{2} \left|\int_{\mathbb{R}^d} \frac{D^2_h u(x+z)-D^2_h u(x)}{|h|^{d+2s}}\mathrm{d}h\right|\\
 	&\leq \frac{C_{d,s}}{2} \int_{\mathbb{R}^d} \frac{|z|^\gamma \land |h|^{\gamma}}{|h|^{d+2s}}\mathrm{d}h
 	= \frac{C_{d,s}}{2}|z|^{\gamma -2s} \int_{\mathbb{R}^d} \frac{1\land |h|^{\gamma}}{|h|^{d+2s}}\mathrm{d}h\\
 	&= \frac{C_{d,s}}{2}|z|^{\gamma -2s} c_d\Big(\frac{1}{\gamma-2s}+\frac{1}{2s}\Big)= C_s|z|^{\gamma -2s}.
 	\end{align*}
 	For the second subcase  where $1< 2s<\gamma \leq 2$ we necessarily have  $0<\sigma <1$. Now we put $\gamma= 1+\tau$ with $\tau\in [0,1]$ and recall that $\gamma-2s=\sigma$ so that $\tau= 2s+\sigma-1$. 
 	From Proposition \ref{prop:estimate-second-order} we get
 	\begin{align*}
 	|(-\Delta)^su(x+z)-(-\Delta)^su(x)|
 	&=\frac{C_{d,s}}{2} \left|\int_{\mathbb{R}^d} \frac{D^2_h u(x+z)-D^2_h u(x)}{|h|^{d+2s}}\mathrm{d}h\right|\\
 	&\leq \frac{C_{d,s}}{2} \int_{\mathbb{R}^d} \frac{|z||h|^\tau \land |h|^{1+\tau}}{|h|^{d+2s}}\mathrm{d}h
 	= \frac{C_{d,s}}{2}|z|^{\gamma -2s} \int_{\mathbb{R}^d} |h|^{-d+\sigma-1}(1\land |h|) \mathrm{d}h\\
 	&= \frac{C_{d,s}}{2}|z|^{\gamma -2s} c_d\Big(\frac{1}{\sigma}+\frac{1}{1-\sigma}\Big)= C_s|z|^{\gamma -2s}.
 	\end{align*}
 	In either case, we have $(-\Delta)^su\in C_b^{\gamma-2s}(\mathbb{R}^d)$. 
 	Now, if $m= 1$ that is $\gamma=2s+\sigma+1$ then, $0<2s<1\leq \gamma \leq 2$, 
 	which implies that $0<2s+\sigma<1$. Furthermore,  we have $\nabla u \in C^{2s+\sigma}(\mathbb{R}^d)$. Since $\nabla u $ is bounded, we have $|D^2_h u(x+te_i)-D^2_h u(x)|\leq C|t|$,  and from Proposition \ref{prop:estimate-second-order} one also gets $|D^2_h u(x+te_i)-D^2_h u(x)|\leq C|t||h|^{2s+\sigma}$ so that 
 	\begin{align*}
 	|D^2_h u(x+te_i)-D^2_h u(x) |\leq C |t|(1\land |h|^{2s+\sigma}).
 	\end{align*}
 	On the other hand, $h\mapsto (1\land |h|^{2s+\sigma})|h|^{-d-2s}$ is integrable. Using the dominated convergence theorem one arrives at $\partial_{x_i} (-\Delta)^s u(x)= (-\Delta)^s \partial_{x_i} u(x).$
 	Furthermore, as $\partial_{x_i} u\in C_b^{2s+\sigma}(\mathbb{R}^d),$ 
 	by applying dominated convergence again one gets that the map $x\mapsto (-\Delta)^s \partial_{x_i} u(x) $ is continuous. In sum, $x\mapsto (-\Delta)^su(x)$ is differentiable and $\nabla (-\Delta)^su = (-\Delta)^s \nabla u$. We know from the previous case $\nabla (-\Delta)^su\in C_b^{\sigma}(\mathbb{R}^d)$, that is $(-\Delta)^su\in C_b^{\sigma+1}(\mathbb{R}^d)= C_b^{\gamma-2s}(\mathbb{R}^d)$, since $\gamma-2s= \sigma+1 $. Proceeding by induction with a similar argument the result remains true when $\gamma=m+\sigma> 2$ and $\partial^\alpha (-\Delta)^s u= (-\Delta)^s \partial^\alpha u$ for each multiindex $|\alpha|\leq m$.
 	The boundedness of $(-\Delta)^s$ blatantly follows from the previous estimates. 
 \end{proof}

 \section{Nonlocal elliptic operators}\label{sec:nonlocal-elliptic}
In this section, we define some concepts analogous to those from elliptic partial differential operators of  second order. We then provide some concrete connections between the nonlocal and local notions. Let us recall that the modern theory of elliptic partial differential operators of second order is governed by two influential classes. The class of operators in the divergence form  $\mathscr{A}$ and the class of operators in the non-divergence form $\widetilde{\mathscr{A}}$. To be more precise  these are operators of the following  forms
 
 \begin{align}
 	\mathscr{A}u& = -\operatorname{div}(A(\cdot)\nabla)u+ b\cdot \nabla u+cu= -\sum_{ i,j=1}^d \frac{\partial}{\partial x_i}\big(a_{ij}\frac{\partial}{\partial x_j}\big)u+ \sum_{ i=1}^d b_i \frac{\partial u}{\partial x_i}+ cu\tag{\text{divergence form}}\\
%
 \widetilde{\mathscr{A}} u& = -\operatorname{tr}(A(\cdot) D^2)u+ b\cdot \nabla u+cu= -\sum_{ i,j=1}^d  a_{ij}\frac{\partial^2u}{\partial x_i\partial x_j}+ \sum_{ i=1}^d b_i\frac{\partial u}{\partial x_i}+ cu.\tag{\text{non-divergence form}}
 \end{align}
 
\noindent Here, $c, b_i:\R^d\to \R$ are measurable functions and  $A= (a_{ij})_{ij}: \R^d\to \mathrm{M}(d\times d)$ is a matrix valued measurable function satisfying the ellipticity condition, i.e. there is a constant $\Lambda_0\geq 1$ such that 
\begin{align*}
\Lambda_0^{-1}|\xi|^2\leq \langle A(x)\xi,\xi\rangle = \sum_{ i,j=1}^d a_{ij}(x)\xi_i\xi_j\leq \Lambda_0|\xi|^2\quad\text{for all $x, \xi\in\R^d $}.
\end{align*}

\noindent The operators  $\mathscr{A}$ and $\widetilde{\mathscr{A}}$ are  symmetric if the matrix $A$ is symmetric, i.e. $a_{ij} =a_{ji}$, $1\leq i,j\leq d$ and translation invariant if and only if the coefficients $c, b_i, a_{ij}$  $1\leq i,j\leq d$ are constants. If the coefficients are constants then   $\mathscr{A}=\widetilde{\mathscr{A}}$.  If $A$ is regular enough, we  then have the relation 
\begin{align*}
	\mathscr{A}u= \widetilde{\mathscr{A}} u+ \widetilde{b}\cdot \nabla u\quad\text{with}\quad \widetilde{b}_i=\sum_{j=1}\frac{\partial a_{ij}}{\partial x_j}. 
\end{align*}

\noindent The most studied  elliptic partial differential operator of second order is the Laplacian
\begin{align*}
\Delta u  =\operatorname{div}(\nabla u)= \operatorname{tr}(D^2 u)=\sum_{ i=1}^d\frac{\partial^2 u}{\partial^2 x_i} \quad\text{for $u\in C^2(\R^d)$}.
\end{align*}

\noindent Our purpose here is to introduce the analog notion of an elliptic operator in  divergence and non-divergence form for nonlocal operators especially for integrodifferential operators. 
\begin{definition}[\textbf{Integrodifferential operator in divergence form}]\label{def:nonlocal-divergence-form}
An integrodifferential operator in divergence form is any nonlocal operator $\mathscr{L}$ that can be written in the form
\begin{align}\label{eq:integro-divergence-form}
\mathscr{L} u(x) = \pv\int_{\R^d}(u(x)- u(y))\, \mu(x, \d y), \quad \quad\text{$(x\in \R^d)$}
\end{align}
for a sufficiently smooth function $u:\R^d\to \R$. Here $\big(\mu(x, \d y)\big)_{x\in \R^d}$ is a family of Borel measures satisfying $\mu(x, \{x\})=0$ for every $x\in \R^d$. In practice, for consistency reasons, it is often assumed that the family  $\big(\mu(x, \d y)\big)_{x\in \R^d}$ satisfies the uniform L\'evy integrability condition
\begin{align}\label{eq:divergence-levy}
	\sup_{x\in \R^d}\int_{\R^d}(1\land |x-y|^2)\mu(x,\d y)<\infty. 
\end{align}

We say that $\mathscr{L}$ is symmetric if for all measurable sets $A,B\subset \R^d$ we have 

\begin{align*}
\iil_{AB} \mu(x, \d y) \d x =\iil_{BA} \mu(x, \d y)\d x.
\end{align*}

\end{definition}


\begin{definition}[\textbf{Integrodifferential operator in non-divergence form}]\label{def:nonlocal-non-divergence-form}
	An integrodifferential operator in non-divergence form is any nonlocal operator $\widetilde{\mathscr{L}}$ that can be written in the form
\begin{align}\label{eq:integro-non-divergence-form}
\widetilde{\mathscr{L}}u(x) = \frac12\int_{\R^d}(2u(x)- u(x+h)-u(x-h)\, \widetilde{\mu}(x, \d h), \quad\quad\text{$(x\in \R^d)$}
\end{align}
for a sufficiently smooth function $u:\R^d\to \R$. Here $\big(\widetilde{\mu}(x, \d h)\big)_{x\in \R^d}$ is a family of Borel measures satisfying $\widetilde{\mu}(x, \{0\})=0$ for every $x\in \R^d$. In practice, for consistency reasons, it is often assumed that the family  $\big(\widetilde{\mu}(x, \d y)\big)_{x\in \R^d}$ satisfies the uniform L\'evy integrability condition
\begin{align}\label{eq:non-divergence-levy}
\sup_{x\in \R^d}\int_{\R^d}(1\land |h|^2)\widetilde{\mu} (x,\d h)<\infty. 
\end{align}

\noindent We say that $\widetilde{\mathscr{L}}$ is symmetric if $\widetilde{\mu}(x, A)= \widetilde{\mu}(x, -A) $ for all measurable set $A\subset \R^d$ and all  $x\in \R^d$.
\end{definition}

\medskip

 \begin{remark} If $u\in C_b^2(\R^d)$ then the pointwise evaluation $\widetilde{\mathscr{A}}u(x)$ obviously makes sense and under the condition \eqref{eq:non-divergence-levy}, the pointwise evaluation $\widetilde{\mathscr{L}}u(x)$ is also well defined (it is sufficient  to adapt the Proposition \ref{prop:well-defined}). In general, however, the pointwise evaluation  $\mathscr{L}u(x)$ might not be defined under the condition \eqref{eq:divergence-levy} even  for \emph{bona fide} test functions in $C_c^\infty(\mathbb{R}^d)$. This operator is as good or as bad as the local operator $\mathscr{A}$ can be. In fact, if the coefficients $(a_{ij})_{ij}$ are 
sufficiently rough then the pointwise expression $\mathscr{A}u(x)$ might not make sense.  If we assume  that $\mathscr{L}$ is symmetric and the condition \eqref{eq:divergence-levy} holds  and  that $\mathscr{A}$ is 
symmetric and elliptic, however,  then $\mathscr{A}u$ and $\mathscr{L}u$ can be evaluated in the generalized sense. It other words, for all $u,\varphi\in C_c^\infty(\R^d)$, the expressions $\langle\mathscr{A} u, \varphi \rangle$ and $\langle\mathscr{L} u, \varphi \rangle$ are well defined: 
	\begin{align}
	\langle\mathscr{A} u, \varphi \rangle&:=\int_{\R^d}(A(x)\nabla u(x)\cdot \nabla \varphi(x))\d x\\
	\langle\mathscr{L} u, \varphi \rangle&:=\iil_{\R^d\R^d}(u(x)-u(y))(\varphi(x)-\varphi(y)\mu(x, \d y)\d x.
	\end{align}
	\noindent 	Additional discussions on the operators $\mathscr{L}$  and $\widetilde{\mathscr{L}}$ are included in \cite{Kas07}. 
\end{remark}

\medskip

\begin{remark} Assume $\big(\mu(x, \d y)\big)_{x\in \R^d}$ is symmetric i.e $\mu(x, -A)=\mu(x,A)$ for every $A\subset \R^d$ and  is translation invariant, i.e.  $\mu(x, A+h)=\nu(x,A)$ for every $A\subset \R^d$  and $x,h\in \R^d$ then $\mathscr{L}=\widetilde{\mathscr{L}}$. In  the local setting, this corresponds to the situation where the matrix $A$ is constant, hence $\mathscr{A}=\widetilde{\mathscr{A}}$. In a sense, an  integrodifferential operator  that is symmetric and  translation invariant corresponds in the local setting to the situation where the coefficients $c,b_i, a_{ij}$ are constants. A simple instance is provided if we let $\mu(x,\d y) = \nu(x-y)\d y$ where the function $\nu:\R^d\to [0,\infty]$ is measurable and satisfies $\nu(-h) =\nu(h)$ for all $h\in \R^d$. In this particular case, one recognizes the operator $L= \mathscr{L}=\widetilde{\mathscr{L}}$ discussed in the first section. 
\end{remark}
\medskip

\begin{remark} Assume $\mu(x,\d y) = \nu(x-y)\d y$ and $\widetilde{\mu}(x, \d h) = \nu(h)\d h$ where $\nu$ is radial. Assume $A= I_d$ (identity matrix), then we have $L= \mathscr{L}=\widetilde{\mathscr{L}}$ and $-\Delta = \mathscr{A}=\widetilde{\mathscr{A}}$. Moreover, it is not difficult to show that the operators $L$ and $-\Delta$ are  isotropic, i.e. invariant by under rotations. Therefore,  when $\nu$ is radial,  the operator $L$ appears to be a prototypical example of an elliptic integrodifferential operator, just as the Laplace  operator $-\Delta $ is a prototypical example of an elliptic differential operator of second order.  
\end{remark}

\medskip 

\noindent Next we introduce the ellipticity condition associated with the operators  $\mathscr{L}$  and $\widetilde{\mathscr{L}}$.  We restrict ourselves to the case where for each $x\in \R^d$, the measures $\mu(x, \d y)$  and $\widetilde{\mu}(x, \d y)$ are absolutely continuous with respect to the Lebesgue measure. 
\begin{definition}
A measurable function $\nu: \R^d\to [0,\infty]$  will be called unimodal if $\nu$ is radial and almost decreasing,
i.e. there is $c>0$ such that $\nu(x)\geq c\nu(y)$ whenever $|x|\leq |y|$.
\end{definition}

\begin{definition}
	We say that  the family $\big(\mu(x, \d y)\big)_{x\in \R^d}$ is weakly elliptic if there exist a constant $\Lambda\geq 1$  and a unimodal function $\nu: \R^d\to [0,\infty]$ that is L\'evy integrable, i.e. $\nu\in L^1((1\land|h|^2), \d h)$ such that for every Borel set $A,B\subset \R^d$ we have 
\begin{align*}
\Lambda^{-1}\iil_{AB}(u(x)-u(y))^2\nu(x-y)\d y\d x\leq\iil_{AB}(u(x)-u(y))^2\mu(x, \d y)\d x\leq\Lambda \iil_{AB}(u(x)-u(y))^2\nu(x-y)\d y\d x.
\end{align*}
The operator $\mathscr{L}$ is  weakly elliptic if $\big(\mu(x, \d y)\big)_{x\in \R^d}$ is weakly elliptic. 
\end{definition}

\begin{definition}[\textbf{Elliptic kernel}]\label{def:elliptic-kernel}
We say that a kernel $J:\R^d\times \R^d\setminus\operatorname{diag}\to [0, \infty]$ is elliptic if 
\begin{align}
\sup_{x\in \R^d}\int_{\R^d}(1\land|x-y|^2)J(x,y)\d y<\infty. 
\end{align}

\noindent and there exist a constant $\Lambda\geq 1$  and a unimodal function $\nu: \R^d\to [0,\infty]$ that is L\'evy integrable, i.e. $\nu\in L^1((1\land|h|^2), \d h)$ such that 
\begin{align}
\Lambda^{-1}\nu(x-y)\leq J(x,y)\leq\Lambda \nu(x-y)\quad\text{for all $0<|x-y|\leq 1$}.
\end{align}

\noindent We say that $J$ is globally elliptic if instead we have 
\begin{align}
\Lambda^{-1}\nu(x-y)\leq J(x,y)\leq\Lambda \nu(x-y)\quad\text{for all $x\neq y$}.
\end{align}
\end{definition}
\noindent  Obviously the global ellipticity condition implies the ellipticity  condition and the ellipticity condition implies the weak ellipticity condition. 

\begin{remark}
	It is possible to omit the condition that $\nu$ is almost decreasing  in Definition 
\ref{def:elliptic-kernel} as it does not really influence the concept  of ellipticity  defined here. One should keep in mind, however, that this condition only plays a technical role when dealing with function spaces and Poincar\'e  types inequalities. 
\end{remark}

\begin{definition}[\textbf{Elliptic operator}]
Assume $\mu(x, \d y)= J(x,y)\d y$ and $\widetilde{\mu}(x, \d h)= J(x,x+h)\d h$  for a kernel $J:\R^d\times \R^d\setminus\operatorname{diag}\to [0, \infty]$. We say that the operators $\mathscr{L}$ and $\widetilde{\mathscr{L}}$ are elliptic if the kernel $J$ is elliptic. Recall that here, 
\begin{align*}
\mathscr{L}u(x) = \pv \int_{\R^d}(u(x)- u(y)) J(x,y)\d y\quad\text{and}\quad \widetilde{\mathscr{L}}u(x) = \frac12\int_{\R^d}(2u(x)- u(x+h)-u(x-h)\, J(x,x+h) \d h. 
\end{align*}
\end{definition}

\medskip

\noindent Next, we bridge a transition from elliptic integrodifferential operators of the forms $\mathscr{L}$ and $\widetilde{\mathscr{L}}$ to elliptic partial differential operators of the forms $\mathscr{A}$ and $\widetilde{\mathscr{A}}$. We  intend to convince that the aforementioned notions (symmetry, translation invariance, ellipticity, divergence form and non-divergence form) are correlated. Let us introduce $(\nu_\alpha)_{0<	\alpha<2}$, a family of L\' evy radial functions  approximating the Dirac measure at the origin, i.e. for every $\alpha, \delta > 0$
\begin{align*}
\begin{split}
\nu_\alpha\geq 0,\,\,\text{ is radial}, \quad \int_{\mathbb{R}^d}	(1\land |h|^2)\nu_\alpha (h)\d h=1, \quad \lim_{\alpha\to 2}\int_{|h|>\delta}	\nu_\alpha(h)\d h=0\,.
\end{split}
\end{align*}                      

\noindent Note that several examples of $(\nu_\alpha)_\alpha$ are provided in Section \ref{sec:approx-dirac-mass}. For a family $(J^\alpha)_{0<\alpha<2}$ of  positive symmetric kernels $J^\alpha: \mathbb{R}^d\times \mathbb{R}^d \setminus \operatorname{diag} \to [0, \infty]$ we set-up the following:

\medskip

\begin{itemize}
	\item[(E)] There exists a constant $\Lambda\geq 1$ such that for every $\alpha\in  (0,2)$ and all $x,y \in \mathbb{R}^d$, with $0<|x-y|\leq 1$ 
	\begin{align}\label{eq:elliptic-condition-chap1}\tag{$E$}
	\Lambda^{-1} \nu_\alpha (x-y) \leq J^\alpha(x,y) \leq  \Lambda \nu_\alpha (x-y).
	\end{align}
	\item[(L)] For every $\delta >0$
	\begin{align} \label{eq:integrability-condition-chap1}\tag{$L$}
	\lim_{\alpha \to 2^-}\sup_{x\in \mathbb{R}^d} \int_{|h| > \delta} J^\alpha(x,x+h)dh=0.
	\end{align}
	
	\item[(I)] For each $\alpha \in (0,2)$ the kernel $J^\alpha$ is translation invariant, i.e., for every $h \in \mathbb{R}^d$
	\begin{align}\label{eq:translation-invariance-chap1}\tag{I}
	J^\alpha(x+h, y+h) = J^\alpha(x, y) .
	\end{align}
		\item[(G-E)] There exists a constant $\Lambda\geq 1$ such that for every $\alpha\in  (0,2)$ and all $x,y \in \mathbb{R}^d$, with $x\neq y$, 
	\begin{align}\label{eq:global-elliptic-condition-chap1}\tag{$G$-$E$}
	\Lambda^{-1} \nu_\alpha (x-y) \leq J^\alpha(x,y) \leq  \Lambda \nu_\alpha (x-y).
	\end{align}
\end{itemize}

\noindent It is clear that \eqref{eq:global-elliptic-condition-chap1} implies \eqref{eq:elliptic-condition-chap1} and \eqref{eq:integrability-condition-chap1}.    Let us define the elliptic matrix uniquely determined by the family $(J^\alpha)_\alpha$. Given $x \in \R^d$ and $\delta > 0$, we define the symmetric matrix as $A(x) = (a_{ij}(x))_{1\leq i,j\leq d}$ by
\begin{align}\label{eq:coef-matrix-chap1}
a_{ij}(x) = \lim_{\alpha\to 2^{-}} \int_{B_\delta}  h_ih_j J^\alpha(x,x+h)dh.
\end{align}



\medskip

\begin{remark} 
	$(i)$ Under conditions \eqref{eq:elliptic-condition-chap1} and \eqref{eq:integrability-condition-chap1} the expression $\int_{B_\delta}  h_ih_j J^{\alpha_n}(x,x+h)dx$ converges for a suitable subsequence of $(\alpha_n)$. The existence of the limit in \eqref{eq:coef-matrix-chap1} poses an implicit condition on the family $(J^\alpha)_\alpha$. 
	$(ii)$ \eqref{eq:elliptic-condition-chap1} and \eqref{eq:integrability-condition-chap1} ensure that the quantity $a_{ij}(x)$  does not  depend on the choice of $\delta$ and is bounded as a function in $x$. Indeed for all $\delta, r>0$, 
	\begin{align*}
	a_{ij}(x) = \lim_{\alpha\to 2^{-}} \int_{B_\delta}  h_ih_j J^\alpha(x,x+h)dh=  \lim_{\alpha\to 2^{-}} \int_{B_r}  h_ih_j J^\alpha(x,x+h)dh.
	\end{align*}
	$(iii)$ Under condition \eqref{eq:translation-invariance-chap1} the functions $a_{ij}(x)$ are constant in $x$. 
\end{remark}

\begin{example}\label{Ex:J-guys-singular} 
	The conditions \eqref{eq:elliptic-condition}, \eqref{eq:integrability-condition} and \eqref{eq:translation-invariance} are fulfilled for each of the examples below. 
	\begin{align*}
	J_1^{\alpha}(x,y) &= \nu_\alpha(x-y), \\
	J_2^{\alpha}(x,y) &=  \nu_\alpha(x-y) \mathbbm{1}_{B_1}(x-y)+ (2-\alpha)J(x,y) \mathbbm{1}_{\mathbb{R}^d\setminus B_1}(x-y) \,, 
	\end{align*}
where $J$ is a symmetric and translation invariant kernel  such that  $$\sup_{x\in \mathbb{R}^d} \int_{\mathbb{R}^d\setminus B_\delta} J(x,x+h)dh<\infty~\text{ for every }~~\delta >0.$$ 
	\noindent We can also consider the standard kernels
	\begin{align*}
	J_3^{\alpha}(x,y) &= \frac{C_{d,\alpha}}{2} |x-y|^{-d-\alpha} \,, \\
	J_4^{\alpha}(x,y) &=  \frac{C_{d,\alpha}}{2}  |x-y|^{-d-\alpha}\mathbbm{1}_{B_1}(x-y)+ (2-\alpha)|x-y|^{-d-\beta} \mathbbm{1}_{\mathbb{R}^d\setminus B_1}(x-y)\,, 
	\end{align*}
	Here, $\beta>0$ and $C_{d,\alpha}$ is  the normalization constant of the fractional Laplacian. Another example is given as follows. For $e\in \R^d$ we set 
	\begin{align*}
	J_5^{\alpha}(x,y) &=  \big(2+\cos(e\cdot (x-y)))\nu_\alpha(x-y).
	\end{align*}

\noindent The matrix corresponding  to $J_1$ and $J_2$ above,  is $A(x)= \frac{1}{d}I_d= (\frac{1}{d}\delta_{ij})_{1\leq i,j\leq d}$ and for  $J_3$ and $J_4$ the corresponding matrix is  $A(x)= I_d= (\delta_{ij})_{1\leq i,j\leq d}$, where $I_d$ is the identity matrix; see Proposition \ref{prop:elliptic-matrix-chap1}). 
\end{example}

\bigskip

\begin{proposition}\label{prop:elliptic-matrix-chap1}
Assume \eqref{eq:elliptic-condition-chap1} and \eqref{eq:integrability-condition-chap1}.  Consider the symmetric matrix $A= (a_{ij})_{ij}$ from  \eqref{eq:coef-matrix-chap1}, i.e.
	
	\begin{align*}
	a_{ij}(x) = \lim_{\alpha\to 2^{-}} \int_{B_1}  h_ih_j J^\alpha(x,x+h)dh.
	\end{align*}
	
	\begin{enumerate}[$(i)$]
		\item 
		The matrix $A$ is elliptic and has bounded coefficients.  To be more precise,  we have
		\begin{align*}
		d^{-1}\Lambda^{-1}|\xi|^2\leq \langle A(x) \xi, \xi \rangle \leq d^{-1}\Lambda |\xi|^2, \quad\text{ for every } x, \xi \in \R^d .
		\end{align*}
		
		\item  Under condition \eqref{eq:translation-invariance-chap1} each  $x\mapsto a_{ij}(x)$ is constant function.  In particular, for $J_1^\alpha(x,y)= \nu_\alpha(x-y)$ we have  $A(x) = \frac{1}{d}(\delta_{ij})_{1\leq i,j\leq d}$, i.e., the matrix $A$ equals the identity matrix and for  $J^\alpha(x,y) =\frac{C_{d,\alpha}}{2}|x-y|^{-d-\alpha}$ we have  $A(x) = (\delta_{ij})_{1\leq i,j\leq d}$.

		\item For $u\in C_b^2(\R^d)$ we have 	$$\lim_{\alpha \to 2} \widetilde{\mathscr{L}}_\alpha u(x) =- \frac{1}{2} \operatorname{tr}(A(x) \nabla u)(x)= - \frac{1}{2}\widetilde{\mathscr{A}} u(x), \quad \text{for all $x\in \R^d$}$$ 
		where 
		\begin{align*}
		\widetilde{\mathscr{L}}_\alpha u(x) :=- \frac{1}{2}\int_{\mathbb{R}^d }(u(x+h)+ u(x-h) -2u(x)) J^\alpha(x,x+h)\,\d h.
		\end{align*}
		
		\item  For all $u\in C_b^2(\R^d)$ we have $(-\Delta)^{\alpha/2} u(x) \xrightarrow{\alpha\to 2}-\Delta u(x)$ and $L_{\alpha} u(x) \xrightarrow{\alpha\to 2}-\frac{1}{2d}\Delta u(x)$. Recall that 
		\begin{align*}
			 (-\Delta)^{\alpha/2} u(x)&= \frac{C_{d, \alpha}}{2}\int_{\mathbb{R}^d }(u(x+h)+ u(x-h) -2u(x)) \frac{\d h}{|h|^{d+\alpha}}\\
			L_{\alpha} u(x)&= \frac{1}{2}\int_{\mathbb{R}^d }(u(x+h)+ u(x-h) -2u(x)) \nu_\alpha(h)\,\d h.
		\end{align*}
		\item For $u\in H^1(\R^d)$ and $\varphi\in C_c^\infty(\R^d)$  we have 	
		$$\lim_{\alpha \to 2} \big\langle\mathscr{L}_\alpha u, \varphi\big\rangle = \big\langle-\operatorname{div}(A(\cdot) \nabla u), \varphi\big\rangle= \big\langle\mathscr{A} u, \varphi\big\rangle . $$ 
Here,
		\begin{align*}
		\mathscr{L}_\alpha u(x) &:= \pv\int_{\mathbb{R}^d }(u(x) -u(y)) J^\alpha(x,y)\d y\\
	\big\langle\mathscr{A} u, \varphi\big\rangle &=	\int_{\R^d} \big( A(x)\nabla u(x)\cdot\nabla v(x) \big)  \d x\\
	\langle\mathscr{L}_\alpha u, \varphi \rangle&:=\iil_{\R^d\R^d}(u(x)-u(y))(\varphi(x)-\varphi(y)J^\alpha(x, y)\d y\d x.
		\end{align*}
	
	\end{enumerate}
\end{proposition}

\medskip

\begin{proof}
$(i)$ Let $x, \xi \in \R^d$ and $|h|\leq 1$. The condition \eqref{eq:elliptic-condition} implies that 
		\begin{align*}
\Lambda^{-1} \int_{B_1}[\xi\cdot h]^2\nu_\alpha(h) \d h \leq \int_{B_1} [\xi\cdot h]^2 J^\alpha(x,x+h)\d h\leq \Lambda \int_{B_1} [\xi\cdot h]^2\nu_\alpha(h)\d h. 
		\end{align*}
From the definition of the matrix $A$ we have 
\begin{align*}
\langle A(x) \xi, \xi \rangle= \sum_{ i,j=1}^d a_{ij}(x)\xi_i\xi_j=\lim_{\alpha\to 2^-}\int\limits_{|h|\leq 1}  J^\alpha(x,x+h) \sum_{ i,j=1}^d h_ih_j\xi_i\xi_j\d h=  \lim_{\alpha\to 2^-}\int\limits_{|h|\leq 1} [\xi\cdot h]^2J^\alpha(x,x+h)\d h \,.
\end{align*}
	
		\noindent Since the Lebesgue measure is invariant under rotations,  we have
		
		\begin{align*}
		\lim_{\alpha\to 2^-}\int\limits_{|h|\leq 1}[\xi\cdot h]^2\nu_\alpha(h)\d h 
		&=\lim_{\alpha\to 2^-}\int\limits_{|h|\leq 1}\sum_{1\leq i\neq j\leq d}\xi_i\xi_j h_ih_j\nu_\alpha(h)\d h+  \lim_{\alpha\to 2^-}\int\limits_{|h|\leq 1}\sum_{i=1}^d\xi_i^2 h_i^2\nu_\alpha(h)\d h\\
		&=\lim_{\alpha\to 2^-} \sum_{1\leq i\leq d}\xi_i^2 \int\limits_{|h|\leq 1} h_1^2 \nu_\alpha(h)\d h=\lim_{\alpha\to 2^-} |\xi|^2 \int\limits_{|h|\leq 1} h_1^2 \nu_\alpha(h)\d h\\	
		&=\lim_{\alpha\to 2^-} |\xi|^2 d^{-1}\int\limits_{|h|\leq 1} \sum_{1\leq i\leq d}h_i^2 \nu_\alpha(h)\d h =\lim_{\alpha\to 2^-} |\xi|^2 d^{-1}\int\limits_{|h|\leq 1}  \nu_{\alpha}(h)\d h\\
		&= |\xi|^2 d^{-1}\,.
		\end{align*}
	Indeed, due to the symmetric, the sum over $i\neq j$ vanishes. Altogether this gives 
\begin{align*}
	\Lambda^{-1} d^{-1} |\xi|^2  \le\langle A(x) \xi, \xi \rangle \leq \Lambda d^{-1} |\xi|^2 \quad\text{for all }\quad x,\xi \in \R^d
	\end{align*}
	
	\noindent $(ii)$ Obviously $x\mapsto a_{ij}(x)$ is constant if the condition \eqref{eq:translation-invariance-chap1} holds.  Assume $J^\alpha(x,y) = \nu_\alpha(x-y)$. We show that $A= \frac{1}{d}I_d$.  By symmetry, for $i\neq j$ is easy to show that $a_{ij}=0$. From the fact that the Lebesgue measure is rotationally invariant and Remark \ref{rem:asymp-nu} we find that 
	\begin{align*}
	a_{ii}(x)&= \lim_{\alpha\to 2^-}\int\limits_{|h|\leq 1} h^2_i\nu_\alpha(h)\d h =\lim_{\alpha\to 2^-} \int\limits_{|h|\leq 1} h_1^2 \nu_\alpha(h)\d h\\
	&=\lim_{\alpha\to 2^-} \frac{1}{d}\int\limits_{|h|\leq 1} \sum_{ i=1}^d h_i^2 \nu_\alpha(h)\d h =\lim_{\alpha\to 2^-} \frac{1}{d}\int\limits_{|h|\leq 1}|h|^2  \nu_{\alpha}(h)\d h= \frac{1}{d}\,. 
	\end{align*}
	
	\noindent If $J^{\alpha}(x,y) = \frac{C_{d,\alpha}}{2} |x-y|^{-d-\alpha} $ then  $a_{ii}=1$ and thus we get  $A=I_d$. Indeed, it suffices to proceed as above and accounting Proposition \ref{prop:asymp-cds} which asserts that  $\frac{C_{d,\alpha}}{2d \omega_{d-1} (2-\alpha)} \xrightarrow{\alpha\to 2} 1$.

\noindent $(iii)$	\noindent We know from \eqref{eq:bound-second-difference} that
	\begin{align*}
	\left|u(x+h)+u(x-h)-2u(x)\right|\leq 4\|u\|_{ C^2_b(\mathbb{R}^d)}(1 \land |h|^2).
	\end{align*}
 This combined with the assumption \eqref{eq:integrability-condition-chap1} yields that 
	\begin{align*}
	\lim_{\alpha \to 2}\il_{|h|\geq \delta} \left|u(x+h)+u(x-h)-2 u(x)\right| J^\alpha(x,x+h)\,\d h =0.
	\end{align*}
	%
	%
The fundamental theorem of calculus suggests that
	\begin{align*}
	\widetilde{\mathscr{L}}_\alpha u(x)= & -\frac{1}{2} \il_{|h|< \delta} \left(u(x+h)+u(x-h)-2(x)\right)J^\alpha(x,x+h)\,\d h = -\frac{1}{2}\il_{|h|< \delta} [D^2(x)\cdot h] \cdot h ~J^\alpha(x,x+h)\,\d h 
	\\ &-\frac{1}{2} \int_{0}^{1} \int_{0}^{1} 2t \il_{|h|< \delta} [D^2 u(x-t h + 2 st h)\cdot h-D^2 u(x) \cdot h] \cdot h ~J^\alpha(x,x+h)\,\d h \, \,\d s\,\d t.
	\end{align*}
	Since $D^2 u$ (the Hessian of $u$) is continuous at $x$, for any $\varepsilon>0$  there is  sufficiently small  $\delta>0$ such that $|D^2(x+z) -D^2 u(x)|<\varepsilon $ for all  $|z|<4\delta$. This implies,
	\begin{align*}
	&\lim_{\alpha\to 2}\frac{1}{2} \int_{0}^{1} \int_{0}^{1} 2t \il_{|h|< \delta} |((D^2 u(x-th + 2sth)-D^2u(x) ) \cdot h)\cdot h |~\nu_\alpha(h)\,\d h \, \,\d s\,\d t\\
	&\leq \frac{\varepsilon}{2} \lim_{\alpha\to 2} \il_{|h|< \delta}(1\land |h|^2)~\nu_\alpha(h)\,\d h
	=\frac{\varepsilon}{2}\xrightarrow{\eps\to 0}0.
	\end{align*} 
Hence we get 
\begin{align*}
-2\lim_{\alpha\to 0} \widetilde{\mathscr{L}}_\alpha u(x)=\il_{|h|< \delta} [D^2(x)\cdot h] \cdot h ~J^\alpha(x,x+h)\,\d h = \sum_{ i,j=1}^da_{ij}(x) \frac{\partial^2u(x)}{\partial x_i\partial x_j}= \operatorname{tr}\big(A(x)D^2u(x) \big)= \widetilde{\mathscr{A}}u(x). 
\end{align*}
\noindent Note that $(iv)$ is a consequence of $(ii)$ and $(iii)$ whereas $(v)$ is a particular case of Theorem \ref{thm:quadratic-convergence-BBM}. 
\end{proof}

\section{Mixed L\'evy operators}
Loosely speaking, we define \textit{a mixed L\'evy operator or a L\'evy operator with mixed jumps} or  \textit{ a generalized anisotropic operator} as L\'evy operators that can be viewed as the sum of lower-dimensional L\'evy operators. In recent years, the study of the sum of the one dimensional fractional Laplacian has attracted much attention. The Fourier multiplier of such operators are of the form $\psi(\xi)= |\xi_1|^{\alpha_1}+ |\xi_2|^{\alpha_2}+\cdots+ |\xi_d|^{\alpha_d}$ with $\alpha_i\in (0, 2)$ for all $\xi \in \R^d$. Such operators are known as anisotropic operators, see for instance \cite{Cha17} and several other references therein. We wish to generalize this to lower dimensions greater than one and due to this geometrical consideration, we believe there is more to learn from such operators.

\noindent Let us view $\R^d$ as $\R^{d_1}\times \R^{d_2}\times\cdots \times \R^{d_n}$ where $d=d_1+d_2+\cdots+d_n$ with $1\leq n\leq d$ and $d_j\in \mathbb{N} $. For $x\in \R^d$ we write $x= (x_1^*, x_2^*,\cdots, x_n^*)$ where $x_j^*\in \R^{d_j}$. In addition we define $\widetilde{x}_j= (0,0, \cdots, x_j^*,0,\cdots,0)\in \R^{d_1}\times \R^{d_2}\times\cdots \times \R^{d_n}$ so that $x=\widetilde{x}_1+\widetilde{x}_2+\cdots+\widetilde{x}_n$. We identify $\R^{d_j}$ as a linear  sub-variety of $\R^d$ by the means of the correspondence $\R^{d_j}\ni x_j^*\mapsto \widetilde{x}_j= (0,\cdots,0, x_j^*,0,\cdots,0)\in \R^d$. Let $\nu_j(\d h_j^*)$ be a symmetric L\'evy measure on  $\R^{d_j}$, i.e. 
$\nu_j(A_j)= \nu_j(-A_j)$ for all $A_j\subset \R^{d_j}$, $\nu_j(\{0\})=0$ and 
\begin{align*}
\int_{\R^{d_j}}(1\land |h_j^*|^2)\nu_j(\d h_j^*)<\infty. 
\end{align*}

\begin{definition}
	For each $x\in \R^d$ we define the (mixed) measures
	\begin{align}\label{eq:mixed-levy-measure}
	\mu(x, \d y)&= \sum_{j=1}^{n} \nu_j(x_j^*-\d y_j^*)\prod_{i\neq j}\delta_{x_i^*}(\d y_i^*)\quad \text{and}\quad 
	\widetilde{\mu}(x, \d h) = \sum_{j=1}^{n} \nu_j(\d h_j^*)\prod_{i\neq j}\delta_{0_i^*}(\d h_i^*). 
	\end{align}
\noindent Here $\delta_{x_i^*}(\d y_i^*)$  represents the Dirac measure at $x_i^*\in \R^{d_i}$. 
\end{definition}

It is noteworthy to mention that  for each $x\in \R^d$  the measures  $\widetilde{\mu}(x, \d h) $  is supported on the sub-varieties $ \R^{d_1},\R^{d_2},\cdots ,\R^{d_n}$ whereas the measure $\mu(x, \d h)$ is supported on the sub-varieties $ \widetilde{x}_1+\R^{d_1}, \widetilde{x}_2+\R^{d_2},\cdots, \widetilde{x}_n+\R^{d_n}$.  To be more precise we have the following.

\medskip

	\begin{proposition}(\textbf{Integration rule})\label{prop:rule-integration}
		\begin{enumerate}[$(i)$]
			\item For each $x\in \R^d$ we have $\operatorname{supp}\widetilde{\mu}(x, \d h) \subset \R^{d_1}\cup \R^{d_2}\cup\cdots \cup \R^{d_n}$ and 
			$\operatorname{supp}\mu(x, \d h) \subset \widetilde{x}_1+\R^{d_1}\cup \,\widetilde{x}_2+\R^{d_2}\cup\cdots \cup \, \widetilde{x}_n+\R^{d_n}$. Moreover, for a  Borel set $A\subset \R^d$, if  we identify  $A_j= A\cap \R^{d_j}$ as a subset of $\R^{d_j}$ then, we have 
			\begin{align*}
			\mu(x,A)= \sum_{j=1}^{n} \nu_j(x_j^*+A_j)\quad\text{and}\quad \widetilde{\mu}(x,A)= \sum_{j=1}^{n} \nu_j(A_j).
			\end{align*}
			
			\item 
Let $f:\R^d\to \R$ be measurable then for each $x\in \R^d$ we have 
	\begin{align*}
		\int_{\R^d}f(y)\mu(x,\d y)&= \sum_{j=1}^{n} \int_{\R^{d_j}}f(x+\widetilde{h}_j)\nu_j(\d h_j^*)
			\quad\text{and}\quad \int_{\R^d}f(h)\widetilde{\mu}(x,\d h)= \sum_{j=1}^{n} \int_{\R^{d_j}}f(\widetilde{h}_j)\nu_j(\d h_j^*).
	\end{align*}
\end{enumerate}
	\end{proposition}

\noindent As a consequence, we have the following L\'evy integrability condition:
		\begin{align*}
	\int_{\R^d}(1\land |x-y|^2)\mu(x,\d y)&= \sum_{j=1}^{n} \int_{\R^{d_j}}(1\land |h_j^*|^2)\nu_j(\d h_j^*) = \int_{\R^d}(1\land |h|^2)\widetilde{\mu}(x,\d h)<\infty. 
	\end{align*}

	\begin{definition}
	We shall call a mixed (or anisotropic) L\'evy operator  any L\'evy type integrodifferential operator whose L\'evy measure can be represented in one of the forms in \eqref{eq:mixed-levy-measure}. For instance the following operator $\mathscr{L}$ and $\widetilde{\mathscr{L}}$ are mixed L\'evy operators. 
	
		\begin{align}\label{eq:mixed-example}
		\begin{split}
	\mathscr{L}u(x)&= 	\pv \int_{\R^d}(u(x) -u(y)) \mu(x,\d y)\\
	\widetilde{\mathscr{L}}u(x)&= \frac12\int_{\R^d}(u(x+h)-u(x-h)-2u(x))\widetilde{\mu}(x,\d h).
		\end{split}
		\end{align}
	\end{definition}

\noindent According to Proposition \ref{prop:rule-integration} we get the following proposition.
	\begin{proposition}
Let $u\in C_c^\infty(\R^d)$  and let $\mathscr{L}$ and $\widetilde{\mathscr{L}}$ be as  in \eqref{eq:mixed-example} then for $x\in \R^d$ we have 
		\begin{align*}
\mathscr{L}u(x)&= \sum_{j=1}^{n} \pv  \int_{\R^{d_j}}(u(x)-u(x+\widetilde{h}_j)\nu_j(\d h_j^*)
	\\ 
\widetilde{\mathscr{L}}u(x)&= \sum_{j=1}^{n} \frac12\int_{\R^{d_j}} (u(x+\widetilde{h}_j)+u(x-\widetilde{h}_j)-2u(x))\nu_j(\d h_j^*). 
	\end{align*}
In addition, we have 
\begin{align*}
	\langle \mathscr{L}u, u\rangle = \sum_{j=1}^{n} \iil_{\R^d \R^{d_j}}(u(x)-u(x+\widetilde{h}_j)^2\nu_j(\d h_j^*)\d x. 
\end{align*}
\end{proposition}

\noindent We see that $\mathscr{L}$ and $\widetilde{\mathscr{L}}$ are in some way the sum of lower dimensional L\'evy operators. This is verified within the Fourier symbol. 
\begin{proposition}
Let $\psi$ be the Fourier symbol of $\widetilde{\mathscr{L}}$, i.e. $\widehat{\widetilde{\mathscr{L}} u}(\xi) = \psi(\xi) \widehat{u}(\xi)$ for $\xi \in \R^d$ and $u\in C_c^d(\R^d)$. Then for each $\xi \in \R^d$ we have $\psi(\xi) = \psi_1(\xi_1^*)+ \psi_1(\xi_2^*)+ \cdots+\psi_1(\xi_n^*)$, where $\psi_j(\xi_j^*)$, $j=1,2,\cdots, n$ is the Fourier symbol of the operator $\widetilde{\mathscr{L}}_j$ defined by 
\begin{align*}
\widetilde{\mathscr{L}}_ju(x_j^*)=\frac12\int_{\R^{d_j}}(u(x_j^*+h_j^*)+ u(x_j^*+h_j^*)-2u(x_j^*))\nu_j(\d h_j^*). 
\end{align*}
\end{proposition}
	\begin{proof}
	According to Proposition \ref{prop:rule-integration} we get 
	\begin{align*}
	\psi(\xi) &= \int_{\R^d}(1-\cos(\xi\cdot h)) \widetilde{\mu}(x,\d h)= \sum_{ j=1}^n\int_{\R^{d_j}}(1-\cos(\xi\cdot \widetilde{h}_j))\nu_j(\d h_j^*)\\
	&= \sum_{ j=1}^n\int_{\R^{d_j}}(1-\cos(\xi_j^*\cdot h_j^*))\nu_j(\d h_j^*)
	= \sum_{ j=1}^n\psi_j(\xi_j^*).  
	\end{align*}
	\end{proof}
	
\noindent \textbf{Notation:} To alleviate the notations we write
\begin{align*}
\widetilde{\mathscr{L}}u= \big[\widetilde{\mathscr{L}}_1+\widetilde{\mathscr{L}}_2+ \cdots+ \widetilde{\mathscr{L}}_n\big] u\\
\mathscr{L}u= \big[\mathscr{L}_1+\mathscr{L}_2+\cdots+\mathscr{L}_n\big]u. 
\end{align*}	

\noindent \textbf{Meaning in term of process:} Assume  that $\widetilde{\mathscr{L}}$  is the generator of a L\'evy process $(X_t)_t$, then the  process $(X_t)_t$ jumps according to the following rules:
\begin{itemize}
	\item The process can start at any point $x\in \R^d$. 
	\item  If the process sites at a point $x\in \R^d$ then the process is only allowed to jump to the points of the  form $y=x+\widetilde{h}_j$, $j=1,2,\cdots, n$, where we recall $\widetilde{h}_j=  (0,\cdots,0, h_j^*,0,\cdots,0)\in \R^d$ and $h_j^*\in \R^{d_j}$. 
	\item The rate jump from $x$ to $y=x+\widetilde{h}_j$ is according to the L\'evy measure $\nu_j(\d h_j^*)$ on $\R^{d_j}$. Roughly speaking, the process behaves like a $d_j$-dimensional L\'evy process associated with the generator $\widetilde{\mathscr{L}}_j$. 
\end{itemize}
	
\noindent Let us see some concrete examples. For simplicity one may assume that $\nu_j(\d h_j^*)= \nu_j(h_j^*)\d h_j^*$ where 
$\nu_j: \R^{d_j}\to [0, \infty]$ is a symmetric L\'evy measure on  $\R^{d_j}$, i.e.  $\nu_j(h_j^*)= \nu_j(-h_j^*)$ for all $h_j^*\in \R^{d_j}$ and 
\begin{align*}
\int_{\R^{d_j}}(1\land |h_j^*|^2)\nu_j(h_j^*)\d h_j^*<\infty. 
\end{align*}

\begin{itemize}
\item If $n=1$, (i.e. $d=d_1$),  then $\mu(x, \d y) = \nu(x-y)\d y$ (see the first section).

\item  Let $\nu_j(h_j^*) = C_{d_j, \alpha_j}|h_j^*|^{-d_j-\alpha_j}$ with $\alpha_j\in (0,2)$ then the Fourier symbol of $\widetilde{\mathscr{L}}$ which we denote by $(-\Delta)^{\alpha_1/2,\cdots,\alpha_n/2}$(\textit{mixed fractional Laplacian)} is given by $\psi(\xi) = |\xi_1^*|^{\alpha_1}+ \cdots+ |\xi_n^*|^{\alpha_n}$ and we have 
	\begin{align*}
(-\Delta)^{\alpha_1/2,\cdots,\alpha_n/2}u(x)&= \sum_{j=1}^{n} \frac{C_{d_j, \alpha_j}}{2}\int_{\R^{d_j}} (u(x+\widetilde{h}_j)+u(x-\widetilde{h}_j)-2u(x))\frac{\d h_j^*}{|h_j^*|^{d_j+\alpha_j}}\\
&\equiv \big[(-\Delta_1)^{\alpha_1/2} + (-\Delta_2)^{\alpha_2/2} +\cdots+  (-\Delta_n)^{\alpha_n/2} \big]u(x).
\end{align*}
\end{itemize}
\noindent Here, $(-\Delta_j)^{\alpha_j/2}$ is the fractional Laplacian of order $\alpha_j\in (0,2)$ on $\R^{d_j}$.  It is worth noting that for every  $u\in C_b^2(\R^d)$, we have $(-\Delta)^{\alpha_1/2,\cdots,\alpha_n/2}u(x)\to -\Delta u(x)$ as $\min\{\alpha_1, \cdots, \alpha_n\}\to 2$. Note that the special case $n=d$, i.e. $d_1=d_2=\cdots=d_n=1$,  is considered in \cite{Cha17} where the authors established the Harnack inequality and H\"older regularity for such operators.  Similar work is carried out in  \cite{BS07} for the case $\alpha_1=\cdots=\alpha_n= \alpha$ and $d_1=d_2=\cdots=d_n=1$.

\chapter{Nonlocal Sobolev-like Spaces}\label{chap:nonlocal-sobolev}
In this chapter we introduce some nonlocal Sobolev-like spaces  that are generalizations of Sobolev-Slobodeckij spaces which we will very often encounter. Roughly speaking these are just some refinements of classical Lebesgue $L^p$-spaces (just like  the classical Sobolev spaces $W^{1,p}(\Omega)$ and $W_0^{1,p}(\Omega) $ ) whose additional structures are of importance. In short, one can perceive them as nonlocal versions of classical Sobolev spaces of first order generalizing the usual fractional Sobolev spaces. When needed we will recall some basics on classical Sobolev spaces. Nonetheless to better understand the correlation between the nonlocal spaces and the local spaces, we recommend curious readers to hitch-hike some classical text books like \cite{Adams, Alt16, Bre10, EE87,Dorothee-Triebel}. For a thorough investigation on the theory of Sobolev spaces we recommend, \cite{maz2013sobolev}.  We shall begin this chapter by reviewing some elementary properties of standard Lebesgue spaces and the usual Sobolev spaces on an open set. Non-advanced readers should be aware that some complementary basic notions on Lebesgue spaces are added in Appendix \ref{chap:complemnt-Lebesgue}.  In the next section we first visit the nonlocal Hilbert spaces which  are crucial for the study of complement values problems. After, we introduce nonlocal Sobolev-like spaces in their general form. Since such spaces are less common, we will examine some of their rudimentary properties useful later in our analysis. The usual fractional Sobolev spaces will appear as a particular case of our set-up. After, we show that functions in such spaces can be realized as approximation of smooth functions. Finally under some additional assumptions we derive some compact embedding wherefrom  we prove some Poincar\'{e} type inequalities. The theoretical effort spent in this chapter will be rewarded in the following ones. Throughout this chapter, unless otherwise stated, $\Omega$ is an open subset of $\R^d$ and $1\leq p<\infty$.
If $1\leq p<\infty$ then it is very often assumed that the function $\nu : \R^d\to [0,\infty]$  satisfies the $p$-L\'{e}vy integrability condition
\begin{align}\tag{$J_1$}\label{eq:plevy-integrability}
\nu(-h) = \nu(h) ~~\text{ for all $h\in \mathbb{R}^d$ and} ~~\int_{\mathbb{R}^d}(1\land |h|^p)\nu(h)\d h<\infty\,.
\end{align}
%
 
%

%

\section{Preliminaries}

Our main focus in this section is to present the convolution product in $L^p(\R^d)$ along with  some applications to the approximation by smooth functions and the compactness result like the Riesz Fr\'{e}cht-Kolomogorov theorem. Those are in some ways the cornerstone in the sequel.

\subsection{Convolution product}
Let us recall that the convolution product of two measurable functions $u$ and $v$ on $\mathbb{R}^d$ is given by
\begin{align*}
u*v(x)= \int_{\mathbb{R}^d} v(y)u(x-y) \,\d y
\end{align*}
provided that for almost every $x\in \mathbb{R}^d$, the integral on the right hand side exists and makes sense. Of course, $u*v$ will not exist unless suitable restrictions are imposed upon $u$ and $v$. If it does exist, it is a painless exercise to verify that convolution as product is commutative, associative and distributive over the addition and the multiplication by scalars (at least for integrable functions). 

\medskip

\noindent We commence this section with a preparation result. Recall that for $h\in \mathbb{R}^d $, $\tau_h$ denotes the shift function defined by $\tau_h u(x) = u(x+h)$.

\medskip

\noindent The following result often known as the continuity of the shift operator on $L^p$-spaces, profoundly serves our purposes in many perspectives.
\begin{theorem}[Continuity of shift]
Let $u\in L^p(\mathbb{R}^d)$ with $1\leq p<\infty$ then 
\begin{align*}
 \lim_{|h|\to 0} \|\tau_h u- u\|_{L^p(\mathbb{R}^d)}=0.
\end{align*}
\end{theorem}

\medskip

\begin{proof}
Let $\varepsilon>0$ and let $g\in L^p(\mathbb{R}^d)$ be a simple function such that $\|u-g\|_{L^p(\mathbb{R}^d)}<\varepsilon$. We have 
\begin{align*}
 \|\tau_h u-u\|_{L^p(\mathbb{R}^d)}&\leq \|\tau_h u-\tau_h g\|_{L^p(\mathbb{R}^d)}+ \|\tau_h g-g\|_{L^p(\mathbb{R}^d)}+ \|u-g\|_{L^p(\mathbb{R}^d)}\\
 &\leq 2\varepsilon + \|\tau_h g- g\|_{L^p(\mathbb{R}^d)}.
\end{align*}
\noindent Hence, it suffices to prove the result for the simple function $g$ in particular if $g=\mathds{1}_A$ for a measurable set $A$. Note that in this case $A$ necessarily has finite measure. Thereupon, the result for this case comes from the regularity of Lebesgue measure as follows
\begin{align*}
 \|\tau_h \mathds{1}_A-\mathds{1}_A\|^p_{L^p(\mathbb{R}^d)}= \int_{\mathbb{R}^d}|\mathds{1}_A(x-h)-\mathds{1}_A(x)|\,\d x = |(A-h)\setminus A)|+|A\setminus(A-h)|\xrightarrow{|h|\to 0}0.
\end{align*}
\end{proof}
\noindent Note that the above result holds true for $p=\infty$ if and only if $u$ is uniformly continuous. 

\noindent The most fundamental inequality involving convolutions is Young's inequality which determines some special situations where the convolution of two functions exist. Recall that for $1\leq p\leq \infty$ we shall define the number $p'$ by $\frac{1}{p}+\frac{1}{p'} =1$ with the understanding that $p' = \infty$ if $p = 1$, and $p' = 1$ if $p = \infty$.
 \begin{theorem}[Young's inequality]\label{thm:young-inequality} Let $1\leq p,q,r\leq \infty$ such that $\frac{1}{p}+\frac{1}{q} =1+\frac{1}{r}$. Suppose that $u\in L^p(\mathbb{R}^d)$ and $v\in L^q(\mathbb{R}^d)$ then $u*v\in L^r(\mathbb{R}^d)$ and 
\begin{align*}
 \|u*v\|_{L^r(\mathbb{R}^d)}\leq \|u\|_{L^p(\mathbb{R}^d)}\|v\|_{L^q(\mathbb{R}^d)}\,.
\end{align*}

\noindent Moreover if $r=\infty $, the map $x\mapsto u*v(x)$ is uniformly continuous. 
 \end{theorem}

\begin{proof}
First of all observe that $\frac{q}{p'}+\frac{q}{r} =1$, $\frac{p}{p'}+\frac{p}{r} =1$ and $\frac{1}{r}+\frac{1}{p'}+\frac{1}{q'} =1$. Thus, applying the generalized H\"{o}lder inequality with exponents $p',q' $ and $r$ we get 

\begin{align*}
 |u*v(x)|&\leq \int_{\mathbb{R}^d}|v(y)||u(x-y|)\,\d y\\
 & =\int_{\mathbb{R}^d}|v(y)|^{q/p'}|u(x-y)^{p/q'}[|v(y)|^{q/r}|u(x-y)|^{p/r}]\,\d y\\
 &\leq \|v\|^{q/p'}_{L^q(\mathbb{R}^d)}\|u\|^{p/q'}_{L^p(\mathbb{R}^d)}
 \Big(\int_{\mathbb{R}^d}|v(y)|^{q}|u(x-y)|^{p}\,\d y\Big)^{1/r}. 
\end{align*}
Consequently for almost all $x\in \mathbb{R}^d$,
\begin{align*}
 |u*v(x)|^r&\leq \|v\|^{qr/p'}_{L^q(\mathbb{R}^d)}\|u\|^{pr/q'}_{L^p(\mathbb{R}^d)}
 \int_{\mathbb{R}^d}|v(y)|^{q}|u(x-y)|^{p}\,\d y. 
\end{align*}
Employing  Fubini's theorem and using once more the relations $\frac{q}{p'}+\frac{q}{r} =1$ and $\frac{p}{p'}+\frac{p}{r} =1$ yields 

\begin{align*}
\int_{\mathbb{R}^d} |u*v(x)|^r\,\d x &\leq\|v\|^{qr/p'}_{L^q(\mathbb{R}^d)}\|u\|^{pr/q'}_{L^p(\mathbb{R}^d)}
 \int_{\mathbb{R}^d}\int_{\mathbb{R}^d}|v(y)|^{q}|u(x-y)|^{p}\,\d y\,\d x\\
 &= \|v\|^{qr/p'+q}_{L^q(\mathbb{R}^d)}\|u\|^{pr/q'+p}_{L^p(\mathbb{R}^d)}=\|v\|^{r}_{L^q(\mathbb{R}^d)}\|u\|^{r}_{L^p(\mathbb{R}^d)}
\end{align*}
which proves the desired inequality. Now assume that $r=\infty$ and let $h\in \mathbb{R}^d$. By the previous inequality and the continuity of the shift we get that
 \begin{align*}
 \|\tau_h u*v -u*v\|_{L^\infty(\mathbb{R}^d)} 
 &=\|u*(|\tau_h v-v)\|_{L^\infty(\mathbb{R}^d)}\\ 
 & \leq \||\tau_h v-v\|_{L^q(\mathbb{R}^d)}\|u\|_{L^p(\mathbb{R}^d)}\xrightarrow[]{h\to0}0.
\end{align*}
Thereby providing the uniform continuity and the proof is now complete.

\end{proof}

\bigskip

\noindent It is worthwhile mentioning that in general Young's inequality can be established on every locally compact group furnished with the left invariant Haar measure. 
Some special cases of Young’s inequality are much simpler to establish. For instance, assume $r=\infty$ which means $\frac{1}{p}+\frac{1}{q} =1$ or $q=p'$. Let $u\in L^p(\mathbb{R}^d)$ and $v\in L^{p'}(\mathbb{R}^d)$, it follows from H\"{o}lder inequality that for almost every $x\in \mathbb{R}^d$, 
\begin{align*}
 \Big|\int_{\mathbb{R}^d} v(y)u(x-y)\,\d y\Big| \leq \|v\|_{L^{p'}(\mathbb{R}^d)}\|u\|_{L^p(\mathbb{R}^d)}.
\end{align*}
Which can be rewritten as 
 \begin{align*}
 \|u*v\|_{L^\infty(\mathbb{R}^d)} \leq \|v\|_{L^q(\mathbb{R}^d)}\|u\|_{L^p(\mathbb{R}^d)}.
\end{align*}

\noindent In the case where $q=1$ we have $p=r$ and hence if we let $u\in L^p(\mathbb{R}^d)$ and $v\in L^1(\mathbb{R}^d)$ then applying the H\"{o}lder inequality   again yields the following for almost every $x\in \mathbb{R}^d$.
\begin{align*}
 \Big|\int_{\mathbb{R}^d}v(y)u(x-y)\,\d y\Big|& \leq \int_{\mathbb{R}^d}|v(y)|^{1/p}|u(x-y)||v(y)|^{1/p'}\,\d y
 \\
 &\leq \|v\|^{1/p'}_{L^1(\mathbb{R}^d)} \Big(\int_{\mathbb{R}^d}|v(y)||u(x-y)|^p\,\d y\Big)^{1/p}.
\end{align*}
Fubini's theorem implies 

\begin{align*}
 \int_{\mathbb{R}^d}|u*v(x)|^p \,\d x &\leq 
 \|v\|^{p/p'}_{L^1(\mathbb{R}^d)}\int_{\mathbb{R}^d} \int_{\mathbb{R}^d}|v(y)||u(x-y)|^p\,\d y\,\d x\\ &= \|v\|^{p/p'}_{L^1(\mathbb{R}^d)} \|v\|_{L^1(\mathbb{R}^d)}\|u\|^p_{L^p(\mathbb{R}^d)} = \|v\|^p_{L^1(\mathbb{R}^d)}\|u\|^p_{L^p(\mathbb{R}^d)}.
\end{align*}
which is
\begin{align*}
 \|u*v\|_{L^p(\mathbb{R}^d)} \leq \|v\|_{L^1(\mathbb{R}^d)}\|u\|_{L^p(\mathbb{R}^d)}.
\end{align*}

\noindent This inequality is often referred to as the Minkowski's inequality for convolution. The particular case $q=p=r=1$ gives 
\begin{align*}
 \|u*v\|_{L^1(\mathbb{R}^d)} \leq \|v\|_{L^1(\mathbb{R}^d)}\|u\|_{L^1(\mathbb{R}^d)}.
\end{align*} 
 \noindent In this way $L^1(\mathbb{R}^d )$ is a commutative Banach algebra with the convolution as product.
 
 \subsection{Approximation by smooth functions via convolution}
 
Next we want to approximate a given function $u\in L^p(\mathbb{R}^d)$ by smooth functions. This will be derived as an application of the convolution product and the continuity of the shift. We begin with some basic facts. Let $\varphi\in C_c^\infty(\mathbb{R}^d)$ and $u\in L^1_{\operatorname{loc}}(\mathbb{R}^d)$. Then as $\operatorname{supp}\varphi$ has compact support, it is routine 
 to check that $u*\varphi\in C^\infty(\mathbb{R}^d)$ and $\partial^\alpha (u*\varphi) = u*\partial^\alpha\varphi$ for all multi-indices $\alpha\in \mathbb{N}^d$. Assume $u, v\in C(\mathbb{R}^d)$ are continuous functions. If $u$ and $v$ have compact supports so has $u*v$. Indeed, in general the convolution $u*v$, if it exists, satisfies the inclusion 
 \begin{align*}
  \operatorname{supp} u*v\subset \overline{ \operatorname{supp}u+\operatorname{supp}v}.
\end{align*}
Let $(\eta_\varepsilon)_\varepsilon$ be the standard mollifier family that is $\eta_\varepsilon(x) = \varepsilon^{-d}\eta\big(\frac{x}{\varepsilon}\big)$ with
\begin{align*}
 \eta(x)= c\exp{-\frac{1}{1-|x|^2}}\quad\text{if}~|x|<1\quad\text{and}\quad \eta(x)=0\quad\text{if}~|x|\geq 1,
\end{align*}
where the constant $c>0$ is chosen such that 
\begin{align*}
 \int_{\mathbb{R}^d}\eta(x)\,\d x=1.
\end{align*}
It is not difficult to establish that for each $\varepsilon>0$,
\begin{align*}
 \eta_\varepsilon\in C_c^\infty(\mathbb{R}^d),\quad \operatorname{supp}\eta_\varepsilon\subset B_\varepsilon(0)\quad\text{and}\quad \int_{\mathbb{R}^d}\eta_\varepsilon(x)\,\d x=1.
\end{align*}

\begin{theorem}\label{thm:approx-smooth}
Let $1\leq p<\infty$, for all $g\in L^p(\R^d)$ we have 
\begin{align*}
  \lim_{\varepsilon\to 0}\|g-g*\eta_\varepsilon\|_{L^p(\mathbb{R}^d)}=0.
 \end{align*}
Furthermore, $C_c^\infty(\mathbb{R}^d)$ is a dense subspace of $L^p(\mathbb{R}^d)$. 
\end{theorem}
 
 \begin{proof}
 Let $g\in L^p(\mathbb{R}^d)$, since $ \int_{\mathbb{R}^d}\eta_\varepsilon(x)\,\d x=1$ for every $x\in \mathbb{R}^d$ 
 \begin{align*}
  |g*\eta_\varepsilon(x)-g(x)| &=\Big|\int_{\mathbb{R}^d} (g(x-y)-g(x))\varepsilon^{-d}\eta(\varepsilon^{-1} y)\,\d y\Big|\\
  &\leq \int_{\mathbb{R}^d} \Big|g(x-\varepsilon y)-g(x)\Big|\eta(y)\,\d y.
 \end{align*}
Applying Jensen's inequality with respect to the measure $\eta(y)\,\d y$ in combination with Fubini's theorem leads to 
 \begin{alignat*}{2}
 \|g*\eta_\varepsilon-g\|^p_{L^p(\mathbb{R}^d)} 
&= \int_{\mathbb{R}^d} |g*\eta_\varepsilon(x)-g(x)|^p \,\d x 
&&\leq \int_{\mathbb{R}^d}\Big(\int_{\mathbb{R}^d} \Big|g(x-\varepsilon y)-g(x)\Big|\eta(y)\,\d y\Big)^p\,\d x\\ 
&\leq \iil_{\mathbb{R}^d\mathbb{R}^d} \Big|g(x-\varepsilon y)-g(x)\Big|^p\,\d x\eta(y)\,\d y
&&=\int_{\mathbb{R}^d} \Big\|g(\cdot-\varepsilon y)-g\Big\|^p_{L^p(\mathbb{R}^d)}\eta(y)\,\d y.
 \end{alignat*}
 
\noindent For each $\varepsilon>0$, we have $\eta(y)\|g(\cdot-\varepsilon y)-g(x)\|^p_{L^p(\mathbb{R}^d)}\leq 2\eta(y)\|g\|^p_{L^p(\mathbb{R}^d)}\in L^1(\mathbb{R}^d)$ and by continuity of the shift, $\|g(\cdot-\varepsilon y)-g\|^p_{L^p(\mathbb{R}^d)}\xrightarrow{\varepsilon\to 0}0$. 
In virtue of the foregoing,  the dominated convergence theorem implies that we also have $\|g*\eta_\varepsilon-g\|^p_{L^p(\mathbb{R}^d)}\xrightarrow{\varepsilon\to 0}0$. 

\noindent Let us now prove the density of $C_c^\infty(\mathbb{R}^d)$ in $L^p(\mathbb{R}^d)$. 
Let $u\in L^p(\mathbb{R}^d)$ and fix $\delta>0$. From  the dominated convergence theorem we are able to find $j_0\geq 1$ large enough such that $\|u-g\|_{L^p(\mathbb{R}^d)}<\delta/2$ with $g=u\mathds{1}_{B_{j_0}(0)}$. 
Since $g$ has compact support it turns out that $g*\eta_\varepsilon$ is of compact support too. Furthermore,  we have $g*\eta_\varepsilon\in C_c^\infty(\mathbb{R}^d)$ and as previously shown, there is $\varepsilon>0$ small enough for which $\|g*\eta_\varepsilon-g\|_{L^p(\mathbb{R}^d)}<\delta/2$ so that 
 \begin{align*}
  \|u-g*\eta_\varepsilon\|_{L^p(\mathbb{R}^d)}\leq \|u-g\|_{L^p(\mathbb{R}^d)}+\|g*\eta_\varepsilon-g\|_{L^p(\mathbb{R}^d)}<\delta.
 \end{align*}
 This finishes the proof. 
 \end{proof}

 \begin{corollary}
 Assume $\Omega\subset \mathbb{R}^d $ is an open subset and let $1\leq p<\infty$. Then $C_c^\infty(\Omega)$ is dense in $L^p(\Omega)$.
 \end{corollary}
 
 \begin{proof}
Let $(K_j)_j$ be an exhaustion of compact sets of $\Omega$ with $\operatorname{dist(K_j, \partial \Omega)>\frac{1}{j}}$. For $u\in L^p(\Omega)$ and $\delta>0$ small enough there exists $j\geq1$ sufficiently large such that $\|u-u\mathds{1}_{K_j}\|_{L^p(\Omega)}<\delta/2$. Assume $u\mathds{1}_{K_j}$ is extended by zero to $\mathbb{R}^d$ then by Theorem \ref{thm:approx-smooth} for $\varepsilon<1/2j$ sufficiently small, we have $\|u\mathds{1}_{K_j}-(u\mathds{1}_{K_j})*\eta_\varepsilon\|_{L^p(\mathbb{R}^d)}<\delta/2$ so that 
\begin{align*}
\|u- (u\mathds{1}_{K_j})*\eta_\varepsilon\|_{L^p(\Omega)}\leq \|u-u\mathds{1}_{K_j}\|_{L^p(\Omega)}+ \|u\mathds{1}_{K_j}-(u\mathds{1}_{K_j})*\eta_\varepsilon\|_{L^p(\mathbb{R}^d)} <\delta.
\end{align*}
Moreover $u*\eta_\varepsilon$ belongs to $C_c^\infty(\Omega)$ since $ 
\operatorname{supp} (u\mathds{1}_{K_j})*\eta_\varepsilon\subset \overline{K_j+ B_{1/2j}(0)}\subset K_{2j}$ and $K_{2j}$ is a compact subset of $\Omega$. This achieves the proof.
 \end{proof}

 \subsection{The Riesz-Fr\'{e}chet-Kolmogorov theorem }

 With the help of approximation by means of convolutions, we shall provide a compactness criterion of subsets in $L^p(\mathbb{R}^d)$, which is very effective in applications. The concerned result will be obtained through the Arzel\`{a}-Ascoli theorem which gives a criterion for compactness in spaces of functions. We state the result explicitly for the spaces of interest here and we shall omit the proof which can be found in \cite[chapter 3]{yosida80}.

\begin{theorem}[Arzel\`{a}-Ascoli]\label{thm:ascoli}
Assume that $K$ is a nonempty compact set in $\mathbb{R}^d$. A subset $\mathcal{F}$ of $C(K)$ normed by $\|u\|_{C(K)} =\max\limits_{x\in K} |u(x)|$, is precompact if and only if $\mathcal{F}$ is bounded in $C(K)$ and $\mathcal{F}$ is equicontinuous on $C(K)$, i.e. $$\lim\limits_{|h|\to 0^+} \sup_{u\in \mathcal{F}}\|\tau_hu-u\|_{C(K)}= 0.$$
Explicitly, for every $\varepsilon>0$ there exists $\delta>0$ such that for all $|h|<\delta$, if $x,x+h\in K$ then
\begin{align*}
|u(x+h)-u(x)|<\varepsilon\qquad\text{for all }\quad u\in \mathcal{F}.
\end{align*}
\end{theorem}

\medskip

\noindent In connection to the equicontinuity we also have the following analogous concept in $L^p(\mathbb{R}^d)$.
\begin{definition}
 Let $\mathcal{F}$ be a subset of $L^p(\mathbb{R}^d)$ with $1\leq p<\infty$. A subset $\mathcal{F}$ of $L^p(\mathbb{R}^d)$ is said to be $p$-equicontinuous if 
\begin{align*}
\lim_{|h|\to 0} \sup_{u\in \mathcal{F}} \|\tau_h u-u\|_{L^p(\mathbb{R}^d)}= 0.
\end{align*}
\end{definition}

\medskip

\noindent The following result known as the Riesz-Fr\'{e}chet-Kolmogorov\footnote{
This was originally proved by M. Riesz. A further
characterization, given by Fr\'echet and Kolmogorov, is the approximation of
precompact sets by finite-dimensional ones} Theorem gives conditions analogous to 
the ones in the Arzel\`{a}- Ascoli theorem for a set to be precompact in $L^p(\mathbb{R}^d)$. The proof of the following version follows \cite[Theorem 26]{Bre10}.

\begin{theorem}[Riesz-Fr\'{e}chet-Kolmogorov]\label{thm:riesz-kolomorov-1}
Let $1\leq p<\infty$. Assume $\mathcal{F}$ is a bounded $p$-equicontinuous subset of $L^p(\mathbb{R}^d)$. Then $\mathcal{F}\mid_\Omega$ is precompact in $L^p(\Omega$) for any measurable subset $\Omega\subset \mathbb{R}^d$ with finite measure . 
\end{theorem}

\bigskip

\begin{proof}
Given that $L^p(\Omega)$ is complete it suffices to show that $\mathcal{F}\mid_\Omega$ is totally bounded therein. To this end, we fix $\delta>0$ since $p$-equicontinuous we choose $\varepsilon>0$ arbitrarily small such that
\begin{align}\label{eq:wsup-equicontinuous}
 \sup_{u\in \mathcal{F}}\|u-u*\eta_\varepsilon\|_{L^p(\mathbb{R}^d)} <\delta/2
\end{align}
In truth, with the $p$-equicontinuity at hand, by arguing as for the proof of Theorem \ref{thm:approx-smooth} and applying the dominated convergence theorem, one comes to conclusion that 
\begin{align*}
 \lim_{\varepsilon\to 0} \sup_{u\in \mathcal{F}}\|u-u*\eta_\varepsilon\|^p_{L^p(\mathbb{R}^d)} \leq \lim_{\varepsilon\to0} \int_{\mathbb{R}^d} \sup_{u\in \mathcal{F}} \Big\|u(\cdot-\varepsilon y)-u\Big\|^p_{L^p(\mathbb{R}^d)}\eta(y)\,\d y=0.
\end{align*}

\noindent Meanwhile, for each $\varepsilon>0$ and  each $u\in\mathcal{F}$, Young's inequality yields that we have
\begin{align}\label{eq:young-mollifiers}
 \|u*\eta_\varepsilon\|_{L^\infty(\mathbb{R}^d)}\leq C\|\eta_\varepsilon\|_{L^{p'}(\mathbb{R}^d)} \quad\textrm{and}\quad \|\nabla (u*\eta_\varepsilon) \|_{L^\infty(\mathbb{R}^d)}\leq C\|\nabla \eta_\varepsilon\|_{L^{p'}(\mathbb{R}^d)},
\end{align}
where we have used $\nabla(u*\eta_\varepsilon)= u*\nabla\eta_\varepsilon$ and  set $C=\sup\limits_{u\in \mathcal{F}}\|u\|_{L^p(\mathbb{R}^d)}.$ Let $K$ be a compact subset of $\Omega$ then using the left estimate from the previous display  and \eqref{eq:wsup-equicontinuous} we have that for all $u\in\mathcal{F}$
\begin{align*}
 \|u\|_{L^p(\Omega\setminus K)} \leq \|u-u*\eta_\varepsilon\|_{L^p(\mathbb{R}^d)}+ \|u*\eta_\varepsilon\|_{L^p(\Omega\setminus K)}\leq\delta/2+ C|\Omega \setminus K|^{1/p}\|\eta_\varepsilon\|_{L^{p'}(\mathbb{R}^d)}
\end{align*}
Wherefore, choosing the compact set $K$ large so that $|\Omega \setminus K|$ is small enough, we obtain that 
\begin{align}\label{eq:tight-omega}
 \sup\limits_{u\in \mathcal{F}}\|u\|_{L^p(\Omega\setminus K)}<\delta.
\end{align}
For such a compact set $K$ and fixed $\varepsilon>0$ as above we claim that the family $\mathcal{F}*\eta_{\varepsilon}\hspace{-1ex}\mid_K$ is equicontinuous in $C(K)$. From the second estimate in \eqref{eq:young-mollifiers} it follows that 
\begin{align*}
  \hspace{-1ex}\sup\limits_{u\in \mathcal{F}}\|\tau_h u*\eta_\varepsilon-u*\eta_\varepsilon\|_{C( K)}\leq |h|\sup\limits_{u\in \mathcal{F}}\|\nabla (u*\eta_\varepsilon)\|_{L^\infty(\mathbb{R}^d)}\leq C |h|\|\nabla \eta_\varepsilon\|_{L^{p'}(\mathbb{R}^d)}\xrightarrow{|h|\to 0}0.
\end{align*}
In view of the Arzel\`{a}-Ascoli Theorem \ref{thm:ascoli} the set $\mathcal{F}*\eta_{\varepsilon}\hspace{-1ex}\mid_K$ is precompact in $C(K)$ and hence is totally bounded. Whence there exist $g_1,\cdots,g_N\in \mathcal{F}*\eta_{\varepsilon}\hspace{-1ex}\mid_K$ such that 
\begin{align*}
\mathcal{F}*\eta_{\varepsilon}\hspace{-1ex}\mid_K\subset \bigcup_{i=1}^N B^\infty_{\delta_K}(g_i)\quad\textrm{with}\quad \delta_K= \delta|K|^{-1/p}.
\end{align*}
To conclude that $\mathcal{F}\mid_\Omega$ is totally bounded, we show that 
\begin{align*}
\mathcal{F}\mid_\Omega\subset \bigcup_{i=1}^N B_{3\delta}(\overline{g}_i),
\end{align*}
where $\overline{g}_i$ is the zero extension to $\mathbb{R}^d$ of $g_i$. Let $u\in\mathcal{F}$ then by the previous  
inclusion, for some $i$ we have $ \|u*\eta_\varepsilon-g_i\|_{C( K)}\leq \delta_K$ . This implies that 
\begin{align}\label{eq:wgi-mollifier}
 \|u*\eta_\varepsilon-g_i\|_{L^p(K)}\leq \delta.
\end{align}
Combining \eqref{eq:wsup-equicontinuous}, \eqref{eq:tight-omega} and \eqref{eq:wgi-mollifier} one arrives at
\begin{align*}
 \|u-\overline{g}_i\|^p_{L^p(\Omega)}&= \Big(\int_{\Omega\setminus K} |u(x)|^p\,\d x+ \int_K|u(x)-g_i(x)|^p\,\d x\Big)^{1/p}\\
 &\leq \|u\|_{L^p(\Omega\setminus K)}+ \|u*\eta_\varepsilon-u\|_{L^p(K)}+ \|u*\eta_\varepsilon-g_i\|_{L^p(K)}<3\delta.
\end{align*}
That is $u\in B_{3\delta}(\overline{g}_i)$, and the proof ends here. 
\end{proof}

\bigskip

\begin{theorem}[Riesz-Fr\'{e}chet-Kolmogorov]
Let $1\leq p<\infty$. A subset $\mathcal{F}$ of $L^p(\mathbb{R}^d)$ is precompact if and only if $\mathcal{F}$ is bounded, $p$-tight and $p$-equicontinuous in $L^p(\mathbb{R}^d)$. 
\end{theorem}

\medskip
\begin{proof}
Assume $\mathcal{F}$ is bounded, $p$-tight and $p$-equicontinuous. In light of the $p$-tightness, for $\delta>0$ let $\Omega\subset \mathbb{R}^d$ be of finite measure and such that 
\begin{align*}
  \sup_{u\in \mathcal{F}}\il_{\mathbb{R}^d\setminus \Omega} |u(x)|^p\,\d x<\delta/2.
\end{align*}
Theorem \ref{thm:riesz-kolomorov-1} reveals that $\mathcal{F}\hspace{-0.5ex}|_\Omega$ is precompact in $L^p(\Omega)$ thus it is possible to cover $\mathcal{F}\hspace{-0.5ex}|_\Omega$ by finitely many balls of radii $\delta/2$ centred at $g_1,\cdots, g_n\in L^p(\Omega)$. Let $\overline{g}_i$ be the zero extension to $\mathbb{R}^d$ of $g_i$. For $ u\in \mathcal{F}$ such that $u\in B_{\delta/2}(g_i)$ we have 
\begin{align*}
  \| u-\overline{g}_i\|_{L^p(\mathbb{R}^d)}&\leq \| u-\overline{g}_i\|_{L^p(\mathbb{R}^d\setminus\Omega )}+\| u-\overline{g}_i\|_{L^p(\Omega)}\\
  &= \| u\|_{L^p(\mathbb{R}^d\setminus\Omega )}+\| u-\overline{g}_i\|_{L^p(\Omega)}<\delta. 
\end{align*}
It follows that 
$$ \mathcal{F}\subset \bigcup_{i=1}^{n}B_\delta(\overline{g}_i).$$
Thus $\mathcal{F}$ is totally bounded in $L^p(\mathbb{R}^d)$ and hence precompact.

\noindent Conversely assume $\mathcal{F}$ is precompact. Then it is evidently bounded. Let $\varepsilon>0$ and let there exists $g_1,\cdots, g_n\in L^p(\mathbb{R}^d)$ such that $\mathcal{F}\subset \bigcup_{i=1}^{n}B_\delta(g_i)$. The $p$-tightness and $p$-equicontinuity (by continuity of the shift) of the set $\{g_1, \cdots, g_n\}$ in $L^p(\mathbb{R}^d)$ implies that of $\mathcal{F} $ within the estimates 
\begin{align*}
 &\|\tau_h u-u\|_{L^p(\mathbb{R}^d)} \leq \|\tau_h g_i-g_i\|_{L^p(\mathbb{R}^d)}+2\|u-g_i\|_{L^p(\mathbb{R}^d)}
\intertext{and}
&\|u\|_{L^p(\mathbb{R}^d\setminus B_R(0))} \leq \|u-g_i\|_{L^p(\mathbb{R}^d)}+ \|g_i\|_{L^p(\mathbb{R}^d\setminus B_R(0))}.
\end{align*}
for $u\in \mathcal{F}$, $R>0$, $h\in \mathbb{R}^d$ and $i=1,2,\cdots,n$.
\end{proof}

\vspace{1cm}

 \noindent The following theorem generalizes \cite[Corollary 4.28]{Bre10} or \cite[Theorem 6.23]{ACS14}.
 \begin{theorem}\label{thm:compactness-convolution}Let $1\leq p,q,r\leq \infty$ such that $\frac{1}{p}+\frac{1}{q} =1+\frac{1}{r}$. Suppose that $\mathcal{F}$ is a bounded subset of $L^p(\mathbb{R}^d)$ and let $g\in L^q(\mathbb{R}^d)$ then $\mathcal{F}*g$ is relatively compact in $L^r_{\operatorname{loc}}(\mathbb{R}^d) $.
 \end{theorem}

\medskip

\begin{proof}
Since $ \mathcal{F}$ is bounded in $L^p(\mathbb{R}^d)$, for each $u\in \mathcal{F}$ and $h\in \mathbb{R}^d$ Young's inequality implies 
\begin{align*}
 \|u*(\tau_hg-g)\|_{L^r(\mathbb{R}^d)}\leq \|u\|_{L^p(\mathbb{R}^d)}\|\tau_hg-g\|_{L^q(\mathbb{R}^d)}\leq C\|\tau_hg-g\|_{L^q(\mathbb{R}^d)}.
\end{align*}
Together with the continuity of the shift, we get the $r$-equicontinuity of $\mathcal{F}*g$ as follows:
\begin{align*}
 \lim_{|h|\to 0} \sup_{u\in\mathcal{F}}\|\tau_h(u*g)-u*g\|_{L^r(\mathbb{R}^d)} \leq C\lim_{|h|\to 0}\|\tau_hg-g\|_{L^q(\mathbb{R}^d)} =0.
\end{align*}
In case $r=\infty$ we know from Young's inequality (cf Theorem \ref{thm:young-inequality}) that $\mathcal{F}*g\subset C(\mathbb{R}^d)$ so that by Arzel\'a-Ascoli  Theorem \ref{thm:ascoli} we get that $\mathcal{F}*g$ is relatively compact in $C(K)$ for every compact subset $K$ of $\mathbb{R}^d$ that is to say $\mathcal{F}*g$ is relatively compact in $L^\infty_{\operatorname{loc}}(\mathbb{R}^d)$. If $r<\infty$ then the result readily follows from the Kolmogorov–Riesz–Fr\'{e}chet Theorem \ref{thm:ascoli}. 
\end{proof}

\section{Classical Sobolev spaces}

In this section we go through a rudimentary review of Sobolev spaces. We refer the reader to \cite{Adams, maz2013sobolev} for more discussions on the theory of Sobolev spaces. Assume $\Omega\subset \R^d$ is an open set. Let $\alpha=(\alpha_1,\cdots,\alpha_1) \in \mathbb{N}^d_0$ be a multiindex
and $u\in L^1_{\operatorname{loc}}(\Omega)$(space of locally integrable function on$\Omega$). A function  $g\in L^1_{\operatorname{loc}}(\Omega) $ usually denote by $g=\partial^\alpha  u$ is called weak derivative or distributional derivative of
$u$ of order $\alpha$, if
\begin{align*}
\int_\Omega u(x)\varphi(x)\d x= (-1)^{|\alpha|}\int_\Omega g(x)\varphi(x)\d x\qquad \text{for all}\quad \varphi\in C_c^\infty(\Omega). 
\end{align*} 
\noindent The uniqueness of the weak derivative $g=\partial^\alpha  u$ follows from  the fundamental lemma of calculus of variation which asserts that  in $L^1_{\operatorname{loc}}(\Omega)$, only the null function  $u=0$ a.e. on $\Omega$ satisfies
\begin{align*}
\int_\Omega  u(x) \varphi(x)\d x=0\quad\text{for all}\quad\varphi\in C_c^\infty(\Omega).
\end{align*} 
\noindent Let $m\in \mathbb{N}$ and $1\leq p\leq \infty$. The space $W^{m,p}(\Omega)$ is the equivalence classes of functions $u\in L^{p}(\Omega)$ whose distributional derivatives $D^{\alpha}u$, up to the other  $m$, belong to $L^{p}(\Omega)$. In other words we have  $W^{m,p}(\Omega):= \Big\{ u\in L^{p}(\Omega)~:D^{\alpha}u \in L^{p}(\Omega) ,\, |\alpha|\leq m \Big\}$. 
	The space $W^{m,p}(\Omega)$ is furnished with the norm $	\|\cdot\|_{W^{m,p}(\Omega)}$ defined by 
	\begin{align*}
	\|u\|_{W^{m,p}(\Omega)}&:= \Big( \|u\|^p_{L^{p}(\Omega)} + \sum_{|\alpha|\leq m} \|D^{\alpha}u\|^p_{L^{p}(\Omega)} \Big)^{\frac{1}{p}} \qquad\text{for}\quad1\leq p<\infty,\\
	\|u\|_{W^{m,\infty}(\Omega)}&:=\|u\|_{L^{\infty}(\Omega)} + \sum_{|\alpha|\leq m} \|D^{\alpha}u\|_{L^{\infty}(\Omega)}.
	\end{align*}
	
\noindent The closure of $C_c^\infty(\Omega)$ in $W^{m,p}(\Omega)$ is denoted by $W_0^{m,p}(\Omega)$.
\vspace{2mm}

\noindent\textbf{Notation:} For $p=2$ and $m=1$ we shall write $H^1(\Omega)$(resp. $H_0^1(\Omega)$) in place of $W^{1,2}(\Omega)$ (resp.$W_0^{1,2}(\Omega)$).

\vspace{2mm}
\noindent The space $W^{m,p}(\Omega)$ is a separable Banach space (Hilbert space for $p=2$ ) for $1\leq p<\infty$ and reflexive for $1< p<\infty$ (c.f. \cite{Adams}). The absence of reflexivity for the case $p=1$ gives rise to another type of function space. When $m=1$ this is known to be the space of bounded variation functions.  Roughly speaking it is the space of elements in $L^1(\Omega)$ whose derivatives in the sense of distributions are bounded Radon measures. This is formally defined as follows.
\begin{definition}
	The space of functions with bounded variation on $\Omega$ denoted by $BV(\Omega)$ is defined as the space of functions $u\in L^1(\Omega)$ such that $|u|_{BV(\Omega) }<\infty$ (in which case $u$ is said to has bounded variation on $\Omega$) where
	\begin{align}\label{eq:bounded-variation}
	|u|_{BV(\Omega) }:= \sup\left\lbrace \int_{\Omega} u(x) \operatorname{div} \phi(x) \mathrm{d} x:~ \phi \in C_c^\infty(\Omega, \mathbb{R}^d),~\|\phi\|_{L^\infty(\Omega)}\leq 1\right\rbrace.
	\end{align}
	
\noindent We will still denote the distributional derivative of a function $u\in BV(\Omega)$ by $\nabla u$. Roughly speaking,  $\nabla u =(\Lambda_1, \Lambda_2,\cdots, \Lambda_d)$ can be seen as a vector valued Radon measure\footnote{That is any Borel measure which is inner and outer regular and finite on each Compact subset of $\Omega.$} on $\Omega$ such that 
\begin{align*}
\int_\Omega u(x)\frac{\partial \varphi}{\partial x_i}(x) \d x=
-\int_\Omega \varphi(x)\d \Lambda_i(x), \quad \text{for all }\quad \varphi\in C_c^\infty(\Omega),~~i=1,\cdots, d.
\end{align*}
In particular,  if $u\in W^{1,1}(\Omega)$ then $|u|_{BV(\Omega)}\leq \|\nabla u\|_{L^1(\Omega)}\leq d^2 |u|_{BV(\Omega)}$ that is $u\in BV(\Omega)$ and we have $\partial_{x_i} u(x)\d x= \d\Lambda_i(x)$. Indeed, note that for $\varphi\in C_c^\infty(\Omega)$ and $e\in \mathbb{S}^{d-1}$ then $e\varphi\in  C_c^\infty(\Omega, \R^d)$ and $\operatorname{div}(e\varphi)= \nabla\varphi \cdot  e.$ By duality we have 
\begin{align*}
\|\nabla u\|_{L^1(\Omega)} \leq  \sqrt{d} \sum_{ i=1}^d \|\nabla u\cdot e_i\|_{L^1(\Omega)}
&=  \sqrt{d} \sum_{ i=1}^d \sup_{\|\varphi\|_{L^{\infty}(\Omega)}\leq 1} \Big|\int_{\Omega} u(x) \, \operatorname{div}(e_i\varphi)(x)~\mathrm{d}x\Big|\leq d^2 |u|_{BV(\Omega)}.
\end{align*}
Conversely, since $u\in W^{1,1}(\Omega)$, the integration by part implies the following 
\begin{alignat*}{2}
	|u|_{BV(\Omega)} &= \sup_{\|\phi\|_{L^{\infty}(\Omega)}\leq 1} \Big|\int_{\Omega}\nabla  u(x) \, \cdot \phi(x)~\mathrm{d}x\Big|\leq \sup_{\|\phi\|_{L^{\infty}(\Omega)}\leq 1} \int_{\Omega} |\nabla u(x)| |\phi(x)|~\d x\leq \|\nabla u\|_{L^1(\Omega)}. 
\end{alignat*}
\noindent The quantity $|\nabla u|$  may be regarded as a positive Radon measure whose value on an open set $U\subset \Omega$ is  $ |\nabla u|(U) =|u|_{BV(U)}$. Notationally when no confusion may arise, it is often the case to write $\int_{\Omega} \,\d |\nabla u|$ or $|\nabla u|(\Omega)$ to synonymously denote the semi-norm $|u|_{BV(\Omega)}$ also called the total variation of $|\nabla u|$. The space $BV(\Omega)$ becomes a Banach space under the norm $$\|u\|_{BV(\Omega) }= \|u\|_{L^1(\Omega) } + |u|_{BV(\Omega) }.$$ We recommend the books \cite{EG15,maz2013sobolev,Zi12} for further details on the space of functions with bounded variation. 
\end{definition}

\begin{remark} Given $1<p<\infty$ by the reflexivity of $L^p(\Omega)$ it is possible to show that for 
	\begin{align*}
	|u|_{BV^p(\Omega) }:= \sup\left\lbrace \int_{\Omega} u(x) \operatorname{div} \phi(x) \mathrm{d} x:~ \phi \in C_c^\infty(\Omega, \mathbb{R}^d),~\|\phi\|_{L^{p'}y(\Omega)}\leq 1\right\rbrace,
	\end{align*}
	the space $BV^p(\Omega)$ of $L^p(\Omega)$ such that  $	|u|_{BV^p(\Omega) }<\infty$ coincides with the Sobolev space $W^{1,p}(\Omega)$. We show this fact  implicitly in the proof of Theorem \ref{thm:liminf-nu}. Moreover, as for the $p=1$, one shows that $|u|_{BV^p(\Omega)} =\|\nabla u\|_{L^p(\Omega)}$.
	 However, the inclusion $W^{1,1}(\Omega) \hookrightarrow BV(\Omega)$  is strict and continuous. For instance the weak derivative of the function $u(x)= \mathds{1}_{(0,1)}(x) -\mathds{1}_{(-1,0)}(x)$ is $2\delta_0$ (Dirac mass at the origin) whence it belongs to $BV(-1,1)$ but is not in $W^{1,1}(-1,1)$. Actually, in this specific case the radon measure
	 $|\nabla u|$ equals $2\delta_0$ and  $|\nabla u|(-1,1) =|u|_{BV}(-1,1)= 2$.
\end{remark}


\noindent Later  we shall need the following less common approximation result. 
\begin{theorem}[\hspace*{-1ex}\textnormal{\cite[p.172]{EG15},\cite[Theorem 3.9]{AFP00}}]\label{thm:evans-grapiepy}
	Let $\Omega\subset\mathbb{R}^d$ be open. Assume $u \in BV(\Omega)$. There exist functions $(u_n)_n$ in $BV(\Omega)\cap C^\infty(\Omega)$ such that
	$\|u_n-u\|_{L^1(\Omega)}\xrightarrow{n \to \infty}0$ and $\|\nabla u_n\|_{L^1(\Omega)}\xrightarrow{n \to \infty}|u|_{BV(\Omega)} = |\nabla u|(\Omega)$. 
\end{theorem}
\begin{remark}
	Note that the above approximation theorem does not claim that $|u_n-u|_{BV(\Omega)}\xrightarrow{n\to \infty}0$ but rather implies that $\|u_n\|_{W^{1,1}(\Omega)}\xrightarrow{n \to \infty}\|u\|_{BV(\Omega)}$. Strictly speaking, $BV(\Omega)\cap C^\infty(\Omega)$ is not necessarily dense in $BV(\Omega)$. 
	On the other hand, if a function $u\in L^1(\Omega)$ is regular enough say $u \in W^{1,1}(\Omega)$ then $\|\nabla u\|_{L^1(\Omega)} \asymp |u|_{BV(\Omega)}$. From this we find that $BV(\Omega)\cap C^\infty(\Omega)= W^{1,1}(\Omega)\cap C^\infty(\Omega)$. 
\end{remark}
A separate aim of this chapter is to provide some strong connections between the aforementioned Sobolev spaces and the upcoming nonlocal spaces of higher interest. Most of these properties will be realized under some additional assumption on the geometry of $\Omega$. Especially if $\Omega$ is an extension domain understood in the  following sense.

\begin{definition}[cf. \cite{Adams}] An open set  $\Omega \subset \mathbb{R}^d$ is called an $W^{m,p}$-extension (resp. $BV$-extension) domain if there exists a linear operator $E:W^{m,p}(\Omega)\to W^{m,p}(\mathbb{R}^d)$ (resp. $E: BV(\Omega)\to BV(\mathbb{R}^d)$) and a constant $C: = C(\Omega, d)$ depending only on the domain $\Omega$ and the dimension $d$ such that 
	\begin{align*}
	Eu\mid_{\Omega} &= u \qquad\hbox{and} \qquad \|Eu\|_{W^{m,p}(\mathbb{R}^d)}\leq C \|u\|_{W^{m,p}(\Omega)} \qquad\text{for all}\quad u \in W^{m,p}(\Omega)\\
(\text{resp.}\quad
	Eu\mid_{\Omega}& = u \qquad\hbox{and} \qquad \|Eu\|_{BV(\mathbb{R}^d)}\leq C \|u\|_{BV(\Omega)} \qquad\qquad \text{for all}\quad u \in BV(\Omega)).
	\end{align*}
\end{definition} 

\noindent Extension domains are significant in applications and are necessary in order to extend certain embeddings theorems on function spaces defined on domains. Note that,  bounded Lipschitz domains are  both $W^{1,p}$-extension and $BV$-extension domains. The geometry characterization  of extension domains has been extensively studied in the recent years. The $W^{1,p}$-extension property of an open set $\Omega$ implies certain regularity of the boundary $\partial\Omega$.  For instance, according to \cite[Theorem 2]{HKT08},  a $W^{1,p}$-extension domain $\Omega\subset \R^d$ is necessarily is a $d$-set, i.e. satisfies the volume density condition: there exists 
a constant $c>0$ such that for all $x\in \partial \Omega$ and  $0<r<1$ we have $|\Omega\cap B(x,r)|\geq cr^d$.  Using the Lebesgue differentiation  theorem, it is easy to show that the boundary of  a $d$-set $\Omega$, has a zero Lebesgue measure, i.e. $|\partial\Omega|=0$.  Therefore, given a $W^{1,p}$-extension domain $\Omega$, we have the following 
\begin{align}\label{eq:lp-boundary-extension}
\text{$\int_{\partial \Omega} |\nabla Eu (x)|^p\d x=0$\quad for all $u\in W^{1,p}(\Omega)$.}
\end{align} 

\noindent To the best of our knowledge, we do not know whether the geometric characterization  \eqref{eq:lp-boundary-extension} remains true for a $BV$-extension domain. However, we emphasize that \cite[Lemma 2.4]{HKT08} every $W^{1,1}$-extension is a $BV-$extension domain.   Hence, we  will require a $BV$-extension domain $\Omega$ to satisfy the additional condition 
\begin{align}\label{eq:bv-boundary-extension}
\text{$|\nabla Eu | (\partial \Omega) =\int_{\R^d}\mathds{1}_{\partial \Omega}(x) \d|\nabla Eu | =0$\quad  for all $u\in BV(\Omega)$.}
\end{align}
\noindent Some authors prefer \cite{AFP00} to define a $BV$-extension domain together with the condition \eqref{eq:bv-boundary-extension}.  Discussions on $BV-$extension domains can be found in \cite{KMS10,Lah15}. Several references on extension domains for Sobolev spaces can be found in \cite{Zh15}.

\subsection*{Trace spaces of Sobolev spaces}

Next we recall the  notion of traces spaces that are in a certain sense Sobolev spaces on the boundary of a smooth domain. Loosely speaking trace spaces  are models of Sobolev spaces on lower dimensional smooth manifolds. Some cares are needed in order to properly build such spaces. Indeed for functions of Sobolev spaces the classical restriction to a lower dimensional manifold are meaningless. Because a smooth  manifold of lower manifolds has Lebesgue a vanishing measure and Sobolev functions are solely defined in the almost everywhere sense. Notwithstanding, by means of functional methods one can generalize the concept of\textit{ restriction by introducing the notion of trace}. Let us commence by recalling the trace theorem for a smooth domain $ \Omega$ in the Hilbert setting.

\begin{theorem}[Trace theorem, \cite{Dorothee-Triebel, Ponce16elliptic}] Assume $\Omega\subset \mathbb{R}^d $ is  bounded domain with $C^{2}-$ boundary. Then there exists  a bounded linear map (or trace operator) $\gamma_0:H^{1}(\Omega)\to L^2(\partial \Omega)$ such that  $ \operatorname{Ker}\gamma_0 = H^1_0(\Omega)$,  for every $v\in C(\overline{\Omega})$  we have  $\gamma_0v = v|_{\partial \Omega}$,  and there is  some constant $C>0$  for which 
	$$\| \gamma_0 v \|_{L^{2}(\partial \Omega)}\leq  C \|v\|_{H^1(\Omega)}, \qquad\forall~v\in H^1(\Omega).$$
	
	\noindent Similarly the linear  map $\gamma_1: C^1(\overline{\Omega}) \cap H^{1}(\Omega) \to   C(\partial\Omega) \cap L^2(\partial \Omega)$ defined by
	\begin{equation*}
	\gamma_1 v:= \gamma_0\circ \frac{\partial v}{\partial \nu }= \nabla v\cdot \mathrm{n} |_{\partial \Omega}
	\end{equation*}
	continuously extends to $H^{2}(\Omega) $ and  $ \operatorname{Ker}\gamma_0 \cap \operatorname{Ker}\gamma_1 = H^{2}_0(\Omega).$
\end{theorem}

\vspace{2mm}
\begin{definition}
	Let $\Omega$ be a bounded smooth domain.  The  trace spaces $H^{1/2}(\partial\Omega)$ and $H^{3/2}(\partial\Omega)$ are respectively  the range of the trace operators $\gamma_0$ and  $\gamma_1$  that is   
\begin{equation*}
H^{1/2}(\partial\Omega):= \gamma_0(H^1(\Omega))\qquad\hbox{and}\qquad H^{3/2}(\partial\Omega):= \gamma_1(H^{2}(\Omega))
\end{equation*}
so that the mappings  $ \gamma_{0}: H^{1}(\Omega) \to H^{1/2}(\partial\Omega)$ and $ \gamma_{1}: H^{2}(\Omega)) \to H^{3/2}(\partial\Omega)$  are onto and remain bounded respectively under the natural norms 
\begin{equation*}
\|v\|_{H^{1/2}(\partial \Omega)} := \inf\{ \|w\|_{H^{1}( \Omega)}: w\in H^{1}(\Omega)~\hbox{and}~ v= \gamma_{0} w \}
\end{equation*}
and 
\begin{equation*}
\|v\|_{H^{3/2}(\partial \Omega)} := \inf\{\|w\|_{H^{2}( \Omega)}: w\in H^{1}(\Omega)~\hbox{and}~ v= \gamma_{1} w \}.
\end{equation*}
Usually  $H^{-1/2}(\partial\Omega)$ denote the dual space of $H^{1/2}(\partial\Omega)$. Clearly we have the following continuous embeddings: $ H_0^1(\Omega)\hookrightarrow  H^1(\Omega)\hookrightarrow  H^{1/2}(\partial\Omega)\hookrightarrow L^2(\partial\Omega)$. Note that, if $\Omega$ is not smooth, then an intelligent way to define $H^{1/2}(\partial\Omega)$ is to set $H^{1/2}(\partial\Omega)\equiv H^1(\Omega)/H^1_0(\Omega)$ that is, the quotient space of $H^1(\Omega)$ by $H^1_0(\Omega)$ via the equivalent relation $u\sim v$ if and only if $ u-v\in H^1_0(\Omega)$.
\end{definition}

\vspace{2mm}
\noindent Let us also emphasize  that for a domain with sufficiently smooth boundary, it is possible to define an intrinsic norm \cite{Din96, mikhailov11} on the spaces $H^{s}(\partial \Omega) $ with $s= 1/2 ,3/2$ as follows:
$$\|v\|^2_{H^{s}(\partial \Omega)}:= \int_{\partial \Omega}  |u(x)|^2\mathrm{d}\sigma(x) +   \iint\limits_{\partial \Omega\times \partial \Omega}  \frac{|u(x)-u(y)|^2}{|x-y|^{d-1+2s}}\mathrm{d}\sigma(x)\mathrm{d}\sigma(y).$$
Here $\mathrm{d}\sigma $ represents the Lebesgue surface measure on $\partial \Omega$. This is an avant-gout 
 towards the introduction of the so called fractional Sobolev spaces.  For the general case, where  $s\in (0,1)$ see	\cite[Section 1.3.3]{grisvard11}. 
 \vspace{1mm}
 
\begin{remark}
	Roughly speaking the operator $ \gamma_{0}$ models the restriction of $H^{1}(\Omega)$ functions on the 
	boundary $\partial\Omega$ (indeed the classical restriction on $\partial\Omega$ of such functions apriori does not make sense in general since $\partial\Omega$ has zero Lebesgue measure as $d-1$-dimensional manifold). The operator $ \gamma_{1}$ can be seen as extension of the notion of the normal derivative on $H^{2}(\Omega)$. Meanwhile, one can  prove  that operators $\gamma_0$ and $\gamma_1$ are not defined on $L^{2}(\Omega)$ and $H^{1}( \Omega)$ respectively.  A modern treatise on trace spaces in general is included in \cite{Ponce16elliptic}. 
\end{remark}

\section{Nonlocal Hilbert function spaces}\label{sec:function-spaces}

In this section we  introduce generalized Sobolev-Slobodeckij-like function spaces with respect  to a symmetric L\'{e}vy measure $\nu(h) \d h$ tailor made for $L^2$-theory of nonlocal elliptic complement value problems. We will extend this to $L^p$-spaces later. Our standing object is a function $\nu:\mathbb{R}^d\setminus \{0\} \to [0, \infty]$ which satisfies the L\'{e}vy integrability condition, i.e., $\nu \in L^1(\R^d, (1 \wedge |h|^2) \d h)$ and is symmetric, i.e., $\nu(h)= \nu(-h)$ for all $h\in \mathbb{R}^d$. The function $\nu$ then is the density of a symmetric L\'{e}vy measure. In case $\nu$ is radial we adopt the convention by identifying $\nu$ with its radial profile, i.e., $\nu(h) = \nu(|h|), h\in \R^d.$ Let $\Omega\subset \mathbb{R}^d$ be open.

\vspace{2mm}
\noindent $\bullet$ We define the space $\HnuOm$ by, 
\begin{align*}
 \HnuOm= \Big\{u \in L^2(\Omega): \,\, |u|^2_{\HnuOm}:=\iil_{\Omega\Omega} \big(u(x)-u(y) \big)^2 \, \nu (x-y)\d y\,\d x<\infty \Big \}\,.
\end{align*}
equipped with the norm defined as follows
\begin{align*}
\|u\|^2_{ \HnuOm}= \|u\|^2_{L^{2} (\Omega)}+ \iil_{\Omega\Omega} \big(u(x)-u(y) \big)^2 \, \nu (x-y) \d y\,\d x. 
\end{align*}
Note that we  also have  the following nice representation
\begin{align*}
\iil_{\Omega\Omega} \big(u(x)-u(y) \big)^2 \, \nu (x-y) \d y\,\d x= \iil_{\R^d\R^d} \!\!\big(u(x)-u(y) \big)^2 \, \big[\mathds{1}_\Omega (x)\cdot \mathds{1}_\Omega(y)\big]\nu (x-y) \d y\,\d x.
\end{align*}

\noindent $\bullet$ We also introduce the space $  \VnuOm$ which is of highest interest in this work.
\begin{align*}
\VnuOm = \Big\lbrace u: \R^d \to \R \text{ meas. }: \,\mathcal{E}(u,u) :=\frac{1}{2} \!\!\iil_{(\Omega^c\times \Omega^c)^c} \!\!\big(u(x)-u(y) \big)^2 \, \nu (x-y) \d y\,\d x<\infty \Big\rbrace \,.
\end{align*}

\noindent Note that $(\Omega^c\times \Omega^c)^c=(\mathbb{R}^d\times \mathbb{R}^d)\setminus (\Omega^c\times \Omega^c)= (\Omega\times \Omega) \cup( \Omega^c\times \Omega)\cup( \Omega\times \Omega^c)$. Thus, an equal simple representation of the form $\mathcal{E}(\cdot,\cdot) $ is given by 
\begin{align*}
\mathcal{E}(u,u) :=\frac{1}{2} \!\!\iil_{\R^d\R^d} \!\!\big(u(x)-u(y) \big)^2 \, \big[\mathds{1}_\Omega (x)\lor \mathds{1}_\Omega(y)\big]\nu (x-y) \d y\,\d x.
\end{align*}

\vspace{2mm}

\noindent $\bullet$ Assume  $L$  is the L\'evy operator associated with the measure $\nu$. Define the new space  $V^1_\nu(\Omega|\R^d)$ by
\begin{align*}
V^1_\nu(\Omega|\R^d)= \Big\lbrace u\in \VnuOm:\, \text{ $Lu$ exists weakly and } Lu\in L^2(\Omega)\,\Big\rbrace \,.
\end{align*}
Here the weak integrodifferentiability of $Lu$ is understood in the sense of Definition \ref{def:weak-integro}.

\vspace{2mm}
\noindent $\bullet$ Denote by $ \VnuOmO$ the space of functions that vanish on the complement of $\Omega$ i.e
\begin{align*}
 \VnuOmO= \{ u\in \VnuOm~: ~u=0~~\text{a.e. on } \mathbb{R}^d\setminus \Omega\}\,.
\end{align*} 
\begin{proposition} The space  $\VnuOm$ can be equally defined as follows
	\begin{align*}
\VnuOm = \Big\lbrace u: \R^d \to \R \text{ meas. } | \,|u|^2_{\VnuOm} :=\iil_{\Omega\R^d} \!\!\big(u(x)-u(y) \big)^2 \, \nu (x-y) \d y\,\d x<\infty \Big\rbrace \,.
\end{align*}
Moreover, for all $u\in \VnuOm$ we have $|u|^2_{\VnuOm} \leq \mathcal{E}(u,u) \leq 2|u|^2_{\VnuOm} $.
\end{proposition}

\begin{proof}
A routine check, shows that we have 
	\begin{align*}
	\mathcal{E}(u,u) &=\frac{1}{2} \!\!\iil_{\R^d\R^d} \!\!\big(u(x)-u(y) \big)^2 \, \big[\mathds{1}_\Omega (x)\lor \mathds{1}_\Omega(y)\big]\nu (x-y) \d y\,\d x\\
|u|^2_{\VnuOm} &= \frac12 \iil_{\R^d\mathbb{R}^d} \big(u(x)-u(y) \big)^2\, \big[\mathds{1}_\Omega (x)+\mathds{1}_\Omega(y)\big]\nu (x-y) \d y\,\d x
	\end{align*}
Clearly we have $\frac12 \big[\mathds{1}_\Omega (x)+ \mathds{1}_\Omega(y)\big] \leq \big[\mathds{1}_\Omega (x)\lor \mathds{1}_\Omega(y)\big]\leq 2\cdot \frac12\big[\mathds{1}_\Omega (x) +\mathds{1}_\Omega(y)\big].$
Thus the following comparison holds true
	\begin{alignat*}{2}
	&\iil_{\Omega\mathbb{R}^d} \!\!\big(u(x)-u(y) \big)^2\, \nu (x-y) \d y\,\d x \
	\leq \mathcal{E}(u,u)
	&\leq 2\iil_{\Omega\mathbb{R}^d} \!\!\big(u(x)-u(y) \big)^2\, \nu (x-y) \d y\,\d x\,.
	\end{alignat*}
\end{proof}
	

\noindent Therefore, for some proofs the usage of the quadratic form $\iil_{\Omega\mathbb{R}^d} \!\!\big(u(x)-u(y) \big)^2\, \nu (x-y) \d y\,\d x$ which is simpler to handle in place of $\iil_{(\Omega^c\times \Omega^c)^c} \!\!\big(u(x)-u(y) \big)^2 \, \nu (x-y) \d y\,\d x$ should not alarm the reader. In reality, as we can observe, the notation $\VnuOm$ is to emphasize that the integral of the measurable map $(x,y)\mapsto \big(u(x)-u(y) \big)^2 \, \nu (x-y) $ performed over $\Omega\times \mathbb{R}^d$ is finite. 
 Another possible suitable notation of the space $\VnuOm$ could be $H_\nu(\Omega|\R^d)$. But we shall use only the notation $\VnuOm$.  We shall visit some fundamental results in connection to the aforementioned spaces in wider context. 

\subsection{Natural embeddings of nonlocal energy spaces}
Here we intend to answer the following questions: $(i)$ What is the most suitable way to define a norm on $\VnuOm$? $(ii)$ What are the possible ways to embed $\VnuOm$ into some $L^2$-space of functions defined on $\Omega^c$ or on the whole $\R^d$? The answer to the latter question is of great interest for our framework as it is a cornerstone for the study of the Inhomogeneous Neumann  and Robin complement value problems.

\noindent Before we start with our investigation, let us observe that under certain conditions on $\nu$ and $\Omega$, there exists at least a  natural norm on $\VnuOm$. Let us recall the following. 
\begin{definition}
A function $\nu:\R^d\to [0,\infty]$ is said to be the density of an unimodal L\'{e}vy measure if it is radial, $\nu\in L^1(1\land |h|^2\d h)$ and almost decreasing, i.e., there is $c>0$ such that $|y|\geq |x|$ implies $\nu(y)\leq c \nu(x)$.  We merely say that the function $\nu$ is unimodal. 
\end{definition}
\medskip

\begin{proposition}\label{prop:natural-norm-on-V}
	Let $\nu : \mathbb{R}^d\setminus \{0\}\to \mathbb{R}$ be unimodal and $\Omega\subset \mathbb{R}^d$ be a bounded open set. Assume $\Omega\subset B_{|\xi|/2}(0) $ for some $ \xi \in \mathbb{R}^d$ with $\nu(\xi)\neq 0$. In particular, if $\nu$ is fully supported on $\mathbb{R}^d$. Then $\VnuOm \subset L^2(\Omega)$. Equivalently $\VnuOm \cap L^2(\Omega) = \VnuOm$. 
\end{proposition}

\begin{proof}
	First, if $\Omega\subset B_{|\xi|/2}(0)$, then for all $x,y\in \Omega$ we have $\nu(x-y)\geq c'$ with $c'= c\nu(\xi)>0 $. By Jensen's inequality, we have 
	\begin{align*}
	\iil_{(\Omega^c\times \Omega^c)^c} \big(u(x)-u(y)\big)^2 \nu(x-y) \d x \, \d y&\geq c' \iil_{\Omega\Omega} (|u(x)|-|u(y)|)^2\d x \, \d y\\
	&\geq c'|\Omega| \il_{\Omega} \Big(|u(x)|-\hbox{$\fint_{\Omega}|u|$}\Big)^2\d x.
	\end{align*}
	This shows that the mean value $\fint_{\Omega}|u|$ is finite. We conclude $u \in L^2(\Omega)$ because of 
	\[\int_{\Omega}u^2(x) \d x\leq 2\int_{\Omega} \Big(|u(x)|-\hbox{$\fint_{\Omega}$}|u|\Big)^2\d x+ 2|\Omega| \Big(\hbox{$\fint_{\Omega}|u|$}\Big)^2\,. \]
\end{proof}

\medskip

\noindent Therefore, when $\Omega$ is bounded, under the assumptions of Proposition \ref{prop:natural-norm-on-V} it is natural to endow the space $\VnuOm$ with the norm
\begin{align*}
\|u\|^2_{\VnuOm} &:= \|u\|^2_{L^{2} (\Omega)}+ \mathcal{E}(u,u) \asymp  \|u\|^2_{L^{2} (\Omega)}+ |u|^2_{\VnuOm} .
\end{align*}

\noindent Accordingly on $\VnuOm$ it is convenient to define the corresponding inner product
\begin{align*}
\big(u,v\big)_{\VnuOm} &= \big(u,v\big)_{L^2(\Omega)}+ \mathcal{E}(u,v)~~\text{or} ~~ \big(u,v\big)_{\VnuOm}= \big(u,v\big)_{L^2(\Omega)}+ [u,v]^2\,.
\intertext{where} 
&\mathcal{E}(u,v)=\frac{1}{2} \iil_{(\Omega^c\times \Omega^c)^c} \big(u(x)-u(y) \big) \big(v(x)-v(y) \big) \, \nu (x-y) \d y\,\d x \\
&[u,v]_{\VnuOm}=\iil_{\Omega\R^d} \big(u(x)-u(y) \big) \big(v(x)-v(y) \big) \, \nu (x-y) \d y\,\d x \,.
\end{align*}
When the function $\nu$ is bounded, e.g. in the case $\nu(h) = \mathbbm{1}_{B_1}(h)$, the space $\HnuOm$ 
equals $L^2(\Omega)$. The same holds true if $\nu \in L^1(\mathbb{R}^d)$. Indeed in such situation if $u\in L^2(\Omega)$ then we have
\begin{align*}
&\iil_{\Omega\Omega} \big(u(x)-u(y) \big)^2 \, \nu (x-y)\d y\,\d x
\leq 4 \iil_{\Omega\Omega} u^2(x) \, \nu (x-y)\d y\,\d x\notag\\
&\leq 4 \iil_{\Omega\mathbb{R}^d} u^2(x) \, \nu (h)\mathrm{d}x\,\mathrm{d}h
=4\|\nu\|_{L^1(\mathbb{R}^d)}\|u\|_{L^2(\Omega)}\,.\label{eq:integrable-nu}
\end{align*}

\begin{definition}
From now on, for a general measure $\nu$ with full support, even if $\Omega$ is not bounded, we refer to the space $\VnuOm$ as the space $\VnuOm \cap L^2(\Omega) $ equipped with the norm defined by 
$$\|u\|^2_{\VnuOm}= \|u\|^2_{L^2(\Omega)}+ \mathcal{E}(u,u). $$
The space $ V^1_\nu(\Omega|\R^d)$ shall be  equipped with the norm
\begin{align*}
\|u\|^2_{V^1_\nu(\Omega|\R^d)}=\|Lu\|^2_{L^2(\Omega)}+ \|u\|^2_{L^2(\Omega)}+ \mathcal{E}(u,u)=\|Lu\|^2_{L^2(\Omega)}+  \|u\|^2_{\VnuOm}.
\end{align*}
 \noindent Some authors find it convenient to work with space $\VnuOm \cap L^2(\mathbb{R}^d) $ equipped with the norm
\begin{align*} 
\vertiii{u}^2_{\VnuOm} &:= \|u\|^2_{L^{2} (\mathbb{R}^d)}+  \mathcal{E}(u,u).
\end{align*}
\begin{remark}
The space $V^1_\nu(\Omega|\R^d)$ can be understood as a nonlocal version of $H^2(\Omega)$. Indeed it is well-known from Cald\'{e}r\`{o}n-Zygmund inequality that for a smooth domain $\Omega,$ $H^2(\Omega) = \{u\in H^1(\Omega):\,\, \Delta u\in L^2(\Omega)\}$ where $\Delta u$ is understood in the distributional sense.  Moreover $\|\cdot\|_{H^2(\Omega)} \asymp  \|\cdot\|_{H^1(\Omega)}+ \|\Delta \cdot\|_{L^2(\Omega)} $. 
\end{remark}

\noindent Clearly, $\VnuOm \cap L^2(\mathbb{R}^d)\subset \VnuOm$ and $\|u\|_{\VnuOm} \leq \vertiii{u}_{\VnuOm}$ for all $u\in \VnuOm \cap L^2(\mathbb{R}^d)$ so that the embedding $\VnuOm \cap L^2(\mathbb{R}^d)\hookrightarrow \VnuOm $ is continuous. Note that as in \eqref{eq:integrable-nu}, for an integrable function $\nu$, the space $\VnuOm\cap L^2(\mathbb{R}^d)$ equals $L^2(\mathbb{R}^d)$ whereas, $\VnuOm$ reduces to the subspace of functions $ u\in L^2(\Omega)$ with
\begin{align*}
\iil_{\Omega\Omega^c}(u(x)-u(y))^2\nu(x-y)\d x\d y<\infty.
\end{align*}
This shows that in many situations the space $\VnuOm$ might be a little larger than $\VnuOm\cap L^2(\mathbb{R}^d)$. 
Nevertheless, when $\Omega= \mathbb{R}^d$ the normed spaces $\big(\VnuOm, \|\cdot\|_{\VnuOm}\big)$, $\big(\VnuOm \cap L^2(\mathbb{R}^d), \vertiii{\cdot}_{\VnuOm}\big)$ and $\big(H_\nu (\Omega), \|\cdot\|_{ \HnuOm}\big)$ are all equal and coincide in norms and will be denoted by $(H_\nu(\mathbb{R}^d), \|\cdot \|_{H_\nu(\mathbb{R}^d)})$. 
\end{definition}
\medskip

\begin{theorem}\label{thm:completness-vnup}
	Assume that $\nu:\R^d\setminus \{0\}\to \R$  is a symmetric L\'{e}vy measure. The function spaces  $\big(\VnuOm\cap L^2(\R^d), \vertiii{\cdot}_{\VnuOm}\big)$ and $\big(\HnuOm, \|\cdot\|_{\HnuOm}\big)$ are separable Hilbert spaces.
	
	\noindent  In addition if $\nu$ has full support in $\R^d$ then the same is true for the space $\big(\VnuOm, \|\cdot\|_{\VnuOm} \big)$ and the space $\big(V^1_\nu(\Omega|\R^d), \|\cdot\|_{V^1_\nu(\Omega|\R^d)} \big)$
\end{theorem}  

\vspace{2mm}

\begin{proof}
For the spaces  $\big(\VnuOm\cap L^2(\R^d), \vertiii{\cdot}_{\WnuOmR}\big)$ and $\big(\HnuOm, \|\cdot\|_{\WnuOm}\big)$ and $\big(\VnuOm, \|\cdot\|_{\WnuOmR} \big)$ the statements are contained in Theorem \ref{thm:completness-Wnup}. Now assume $\nu$ is fully supported and let $(u_n)_n\subset V^1_\nu(\Omega|\R^d) $. Then $(u_n)_n$ and $(Lu_n)_n$ are Cauchy sequences in $\big(\VnuOm, \|\cdot\|_{\VnuOm} \big)$ and $L^2(\Omega)$ respectively. Let $u\in \VnuOm $ such that $\|u_n-u\|_{\VnuOm}\xrightarrow{n\to \infty}0$ and $g\in L^2(\Omega)$ such that $\|Lu_n-g\|_{L^2(\Omega)}\xrightarrow{n\to \infty}0$. It remains to show that $Lu=g$ distributionally which is straightforward. According to Definition \ref{def:weak-integro},  since $\|Lu_n-g\|_{L^2(\Omega)}+ \|u_n-u\|_{L^2(\Omega)}\xrightarrow{n\to \infty}0$, we have $Lu= g$. Indeed since each $Lu_n$ is define in the weak sense, for $\varphi\in C_c^\infty(\Omega)$ with $K= \supp \varphi$ we get
\begin{align*}
\int_K g(x)\varphi(x)\d x= \lim_{n\to \infty} \int_K Lu_n(x)\varphi(x)\d x =  \lim_{n\to \infty} \int_K u_n(x) L\varphi(x)\d x =  \int_K u(x) L\varphi(x)\d x.
\end{align*}
\end{proof}

%

\bigskip

\noindent The following lemma shows that it is possible to define certain norms on $\VnuOm $ (with nice properties) which are equivalent to the norm $\|\cdot\|_{\VnuOm}$ when the function $\nu$ is unimodal. 

\begin{lemma}\label{lem:natural-norm-on-V}
	Let $\Omega\subset \mathbb{R}^d $ be open $($not necessarily bounded$)$ such that $\Omega^c$ has a
	 non-empty interior. Assume $\nu : \mathbb{R}^d\setminus \{0\}\to \mathbb{R}$ is unimodal and has full support. Then there exists an almost decreasing radial Radon measure $\widetilde{\nu}: \mathbb{R}^d \to [0,\infty) $ with full support and a constant $C>0$ both depending only on $\nu, d$ and $\Omega$ such that the following assertions hold true. 
	
	\begin{enumerate}[$(i)$]
		\item $\widetilde{\nu}(\mathbb{R}^d)<\infty$.
		\item $0< \widetilde{\nu} \leq C(1\land \nu)$ a.e. Moreover, we have $\widetilde{\nu} \asymp 1\land \nu$ if in addition $\nu $ satisfies the doubling grow condition: there is a constant $\kappa >0 $ such that $\nu(h)\leq \kappa \nu(2h) $ for all $|h|\geq 1$.
		\item We have the continuous embeddings $$\VnuOm \hookrightarrow L^2(\R^d,\widetilde{\nu}) \hookrightarrow L^1(\R^d,\widetilde{\nu}).$$
		\item On $\VnuOm$, the norms $\|\cdot\|^{\#}_{\VnuOm}$ and $\|\cdot\|^{*}_{\VnuOm}$ are equivalent. Where, 
		\begin{align*}
		\|u\|^{*2}_{\VnuOm}&=\int_{\mathbb{R}^d} u^2(x) \widetilde{\nu}(x)\d x+ \iil\limits_{(\Omega^c\times \Omega^c)^c} (u(x)-u(y))^2\nu(x-y)\d x\d y\,,\\
		\|u\|^{\#2}_{\VnuOm}&=\int_{\Omega} u^2(x) \widetilde{\nu}(x)\d x+ \iil\limits_{(\Omega^c\times \Omega^c)^c} (u(x)-u(y))^2\nu(x-y)\d x\d y\,.
		\end{align*} 
	\end{enumerate}
	\noindent	Furthermore, if $\Omega$ is bounded then the norms $\|\cdot\|_{\VnuOm}$ and $\|\cdot\|^{*}_{\VnuOm}$ are also equivalent.
\end{lemma}

\begin{proof}
Assume $\nu $ has full support. We take $R\geq 1$  sufficiently large such that we have $|B_R(0)\cap \Omega|>0$ and $|B_R(0)\cap \Omega^c|>0$. Set $\Omega(R) = B_R(0)\cap \Omega $ (if $\Omega$ is bounded take $R>1$ large so that $\Omega \subset B_R(0)$, i.e. $\Omega(R) = \Omega $). In any case, for all $x\in \Omega(R) $ and all $y \in \mathbb{R}^d$ we have $|x-y|\leq R(1+|y|)$. The monotonicity condition on $\nu$ implies $\nu(R(1+|y|))\leq c\nu(x-y) $. Set $\widetilde{\nu}(h) = \nu(R(1+|h|))$ for $h \in \R^d$, where we abuse the notation and write $\nu(|y|)$ instead of $\nu(y)$ for $y\in \R^d$. Let us show that $\widetilde{\nu}$ satisfies the desired conditions. Firstly, using the scaling and the translation invariance properties of the Lebesgue measure it is easy to show that $\widetilde{\nu}$ also has full support. Note that $|h|\leq R(1+|h|)$ and $R\leq R(1+|h|)$ for all $h \in \R^d. $ Whence, $0<\widetilde{\nu} \leq C(1\land \nu)$ a.e. In addition let us assume that $\nu(h)\leq \kappa \nu(2h) $ for all $|h|\geq 1$. Let $n\geq 1$ be the unique integer such that $2^n\leq R<2^{n+1}$ and consider $\delta = \tfrac{R}{2^{n+1}-R}\geq 1$. For this choice of $\delta$ we have $ R(1+|h|) \leq \dfrac{R|h|(\delta+1)}{\delta} = 2^{n}|h|$ whenever $|h|\geq \delta$ or $|h|/\delta \geq 1$. This entails that $\widetilde{\nu}(h)\geq c\nu(2^n|h|)$ and by assumption we get $\widetilde{\nu}(h)\geq c\kappa^{-n}\nu(|h|) $. Now, if $|h|\leq \delta$ then $ R(1+|h|) \leq R(\delta+1)$ so that $\widetilde{\nu}(h)\geq c\nu(R(\delta+1))=C$. We have proved that $\widetilde{\nu}(h)\geq C(1\land \nu(h)) $ and hence $(ii)$ is verified. Passing through polar coordinates, we have 
	\begin{align*}
	\widetilde{\nu}(\R^d)& = \il_{\R^d} \nu(R(1+|h|))\d h =|\mathbb{S}^{d-1}| \int_{0}^{\infty} \nu(R(1+r)) r^{d-1} \d r\\
	&= |\mathbb{S}^{d-1}|R^{-1} \int_{R}^{\infty} \nu(r) \Big(\frac{r}{R}-1\Big)^{d-1} \d r\leq |\mathbb{S}^{d-1}|R^{-d} \int_{R}^{\infty} \nu(r) r^{d-1} \d r\\
	&= R^{-d} \int_{|h|\geq R} (1\land |h|^2) \nu(h) \d h \leq R^{-d} \int_{\R^d} (1\land |h|^2) \nu(h) \d h< \infty\, .
	\end{align*}
	This proves $(i)$ and hence $L^2(\R^d,\widetilde{\nu}(h)\d h) \subset L^1(\R^d,\widetilde{\nu}(h)\d h)$. We recall that $\VnuOm $ is a subspace of $L^2(\Omega)$ endowed with the norm $\|\cdot\|_{\VnuOm}$. Let $u\in \VnuOm $ since we know that $\widetilde{\nu}(x)\leq C $ and $\widetilde{\nu}(x)\leq c\nu(x-y)$ for all $y\in \Omega(R)$ and all $x\in \mathbb{R}^d $ the following estimates hold
	\begin{align*}
	\int_{\Omega}u^2(x)\d x+\iint\limits_{\Omega\Omega^c} &(u(x)-u(y))^2\nu(x-y)\d y\d x
	\\
	&\geq C^{-1}\int_{\Omega}u^2(x)\widetilde{\nu}(x) \d x+\iint\limits_{\Omega\Omega^c} (u(x)-u(y))^2\nu(x-y)\d y\d x
	\\
	&\geq C^{-1} \int_{\Omega(R)}u^2(x) \widetilde{\nu}(x)\d x+\iint\limits_{\Omega(R)\Omega^c} (u(x)-u(y))^2\nu(x-y)\d y\d x
	\\
	&= C^{-1}\widetilde{\nu}(\Omega^c)^{-1} \iint\limits_{\Omega(R)\Omega^c} u^2(x)\widetilde{\nu}(y)\d y\d x
	+\iint\limits_{\Omega(R)\Omega^c} (u(x)-u(y))^2\nu(x-y)\d y\d x\\
	&\geq ( C^{-1}\widetilde{\nu}(\Omega^c)^{-1}\land c^{-1}) \iint\limits_{\Omega(R)\Omega^c}\Big[ u^2(x)+ (u(x)-u(y))^2 \Big] \widetilde{\nu}(y)\d y\d x\\
	&\geq ( C^{-1}\widetilde{\nu}(\Omega^c)^{-1}\land c^{-1}) \frac{|\Omega(R)|}{2} \int\limits_{\Omega^c}u^2(y) \widetilde{\nu}(y)\d y\, .
	\end{align*}
	The second and the last line imply that $u\in L^2(\R^d, \widetilde{\nu})$ thereby proving that $\VnuOm \subset L^2(\R^d,\widetilde{\nu}(h)\d h)$. Obviously the proof of $(iii)$ is also complete. The first, second and the last line of the above estimates show that there exist two constants $C_1>0$ and $C_2>0$ depending only on $\Omega, \nu, R$ and $d$ such that 
	\begin{align*}
	\|u\|_{\VnuOm}\geq C_1 \| u\|^{\#}_{\VnuOm}\geq C_2 \| u\|^{*}_{\VnuOm}\, .
	\end{align*}
	Together with the trivial inequality $\| u\|^{\#}_{\VnuOm}\leq \| u\|^{*}_{\VnuOm}$, the norms $\| \cdot\|^{\#}_{\VnuOm} $ and $\| \cdot \|^{*}_{\VnuOm}$ turn out to be equivalent. 
	
	\noindent Moreover, if $\Omega$ is bounded, then $R\leq R(1+|h|)\leq R(1+R)$ for all $ h\in \Omega$. The monotonicity of $\nu$ yields $C^{-1}\leq \widetilde{\nu}(h) \leq C$ for all $h\in \Omega$. Hence, $C^{-1}\|u\|_{L^2(\Omega)} \leq\|u\|_{L^2(\Omega, \, \widetilde{\nu}(h)\d h)} \leq C \|u\|_{L^2(\Omega)} $ which implies the equivalent of $\| \cdot \|^{\#}_{\VnuOm} $ and $\| \cdot \|_{\VnuOm}$ thereby proving the equivalence of $\| \cdot\|_{\VnuOm} $ and $\| \cdot \|^{*}_{\VnuOm}$. Part $(iv)$ is proved.
\end{proof}

\vspace{2mm}
\begin{example}
	For $\alpha\in (0,2)$ consider $\nu(h) = |h|^{-d-\alpha}$ and one obtains $\widetilde{\nu}\asymp 1\land \nu $ with $ \widetilde{\nu}(h)= (1+|h|)^{-d-\alpha}$ or $ \widetilde{\nu}(h)= (1+|h|^{d+\alpha})^{-1}$. In the aforementioned case, the space $\HnuOm$ equals the classical Sobolev-Slobodeckij space $H^{\alpha/2}(\Omega)$. For the same choice of $\nu$ we define $V^{\alpha/2}(\Omega|\mathbb{R}^d)$ as the space $\VnuOm$. So that  we have $V^{\alpha/2}(\Omega|\mathbb{R}^d) \hookrightarrow L^2(\R^d, (1+|h|)^{-d-\alpha})$.
\end{example}

\noindent In the special case where $\nu$ is unimodal, the Lemma \ref{lem:natural-norm-on-V} provides several possible ways to define other norms on $\VnuOm$ equivalent to the norm $\|\cdot \|_{\VnuOm}$. In particular we have the continuous embedding $\VnuOm \hookrightarrow L^2(\R^d,\widetilde{\nu}) $. 
Next we obtain similar results for a more general symmetric L\'{e}vy kernel $\nu$.  We start with the following definition 
\begin{definition}\label{def:nuKinfimum}
Let $\nu : \R^d \setminus \{0\} \to [0, \infty]$ be the density of a symmetric L\'{e}vy measure and $K\subset \R^d$ be a measurable  set such that $|K|>0$. We write $\nu_K$ to denote the measurable function  defined by $ \nu_K: \R^d\to [0,\infty]$ such that for $x\in \R^d$,  
$$ \nu_K(x)= \operatorname{essinf}\limits_{y\in K}\nu(x-y).$$
\end{definition}

\begin{remark}
Since $\nu$ is integrable away from the origin it has some decay at infinity. Typically, for  $x\in \R^d$ the infimum of $\nu(x-y)$ is realized when $y$ is far from $x.$
In some sense the function $\nu_K$  destroys the singularity of $\nu$. A simple  way to illustrate this is to consider the example where $\nu$ is unimodal.
\end{remark}

\begin{proposition}\label{prop:unimodal-ball} 
Let $\nu$ be the density of a unimodal L\'evy measure with full support on $\R^d\setminus \{0\}$. Then (i) $\nu_{B_1}\in  L^1(\R^d)$, (ii) $\nu_{B_1}\leq c(1\land \nu)$ and (iii) if in addition $\nu$ satisfies the doubling growth condition:  $\exists\,\,\kappa >0 $ such that $\nu(h)\leq \kappa \nu(2h) $ for all $|h|\geq 1$ then 
\begin{align*}
\nu_{B_1}\asymp 1\land \nu\qquad \text{with} \qquad\nu_{B_1}(x)= \operatorname{essinf}\limits_{y\in B_1}\nu(x-y).
\end{align*}
The same conclusions hold with $B_1$ replaced by any other ball $B$. Furthermore for a ball $B\subset \Omega$ we have  
  \begin{align*}
  \VnuOm\hookrightarrow L^2(\R^d, \nu_{B}) \hookrightarrow L^1(\R^d, \nu_{B}).
  \end{align*}
\end{proposition}

\medskip

\begin{proof}
For  $x \in B_1$ there exists $ x_*\in  B_1$ (diametrically opposite to $x$) such that $2\geq |x-x_*|= 2-|x|\geq 1$. Since $\nu$ is almost decreasing we have $c\nu(2)\leq \nu_{B_1} (x)\leq \nu(x-x_*) \leq c^{-1}\nu(1)$. If $x\in \R^d \setminus B_1$ there exists $ x_*\in B_1$ (diametrically opposite to $x$) such that $|x|\leq  |x|+1=|x-x_*|$ that is $c\nu_{B_1} (x)\leq c \nu(|x|)$ for all $y\in B_1$ $|x-y|\leq  |x|+1\leq 2|x|$ which implies $c\nu(2|x|)\leq \nu_{B_1} (x)$. 
We have $c\nu(|x|)\leq \nu_{B_1} (x)\leq c^{-1}\nu(|x|)$.

\noindent Observe that $\nu_{B_1}$ is bounded on $B_1$ and $\nu_{B_1}\nu$ on $\R^d\setminus B_1$ and  $\nu$ is integrable on $\R^d\setminus B_1$. Therefore, $\nu_{B_1}\in L^1(\R^d)$ and thus $L^2(\R^d, \nu_{B_1}) \hookrightarrow L^1(\R^d, \nu_{B_1})$. 
Further we have shown that $c\nu_{B_1}\leq (1\land \nu)$ and $\nu_{B_1}(x)\geq c (1\land \nu(2|x|))$. If $\nu$ satisfies $\nu(h)\leq \kappa \nu(2|h|)$ for all $|h|\geq1$ then $\kappa^{-1}c(1\land\nu)\leq \nu_{B_1}\leq c^{-1}(1\land \nu)$. One reaches the similar conclusions for any other ball $B$. If $ B\subset \Omega$ we let  $u \in \VnuOm$ and $y\in B$ then 
	\begin{align*}
\int_{\R^d}|u(x)|^2\nu_{B}(x)\d x
&\leq 2|u(y)|^2\|\nu_{B}\|_{L^1(\R^d) }+ 2 \int_{\mathbb{R}^d}(u(x)-u(y))^2\nu(x-y)\d x.
\end{align*}
Integrating both sides with respect to the variable $y$ over $B$ yields
\begin{align*}
\int_{\R^d}|u(x)|^2\nu_{B}(x)\d x &\leq C\int_{B} |u(x)|^2\d y+ C\iint\limits_{K'\mathbb{R}^d}(u(x)-u(y))^2\nu(x-y)\d x\d y\\
&\leq C\int_{\Omega}|u(x)|^2\d y+ C\iint\limits_{\Omega\mathbb{R}^d}(u(x)-u(y))^2\nu(x-y)\d x\d y
\end{align*}
with $C=2|B|^{-1}( \|\nu_{B}\|_{L^1(\R^d) }+1)$. Hence $\VnuOm\hookrightarrow L^2(\R^d, \nu_{B})$ since 
$\|u\|^2_{L^2(\Omega^c, \nu_K)}\leq C\|u\|^2_{\VnuOm}. $ 
\end{proof}

\vspace{1mm}

\begin{remark} $(i)$ The conclusions of Lemma \ref{lem:natural-norm-on-V} remain true with $\widetilde{\nu}$ replaced by $\nu_{B}$ for a ball $B\subset\Omega$.

\vspace{2mm}
\noindent $(ii)$ Under the doubling scaling condition on a unimodal kernel $\nu $,  it is possible to analogously show that for any ball $B$ we have $\nu_B\asymp 1\land \nu$. Hence for any compact set $K$ with a non-empty interior  we have $\nu_K\asymp 1\land \nu$.  This means that almost all $\nu_K$ are comparable. But this is not always the case in general. For instance consider $\nu(h) = |h|^{-d}e^{-|h|^2}$, then we claim that $\nu_{B_1}$ and $\nu_{B_2}$  are not comparable. Indeed assume there is a constant $C>0$ such that $\nu_{B_1}(x)\leq c\nu_{B_1}(x)$  for all $x\in \R^d $.  Let $|x|\geq 4$ then there are two points $x_*\in \partial B_1$ and $x'_*\in \partial B_2$(diametrically opposite to $x$) such that $|x-x_*|= |x|+1$,  $|x-x'_*|= |x|+2$ and 
\begin{align*}
\nu_{B_1}(x)= \nu(x-x_*)  =( |x|+1)^{-d} e^{-(1+|x|)^2} \quad\text{and }\quad \nu_{B_2}(x)= \nu(x-x'_*) =( |x|+2)^{-d}e^{-(2+|x|)^2}.
\end{align*}
That $\nu_{B_1}(x)\leq C\nu_{B_2}(x)$ implies $1\leq C\big(\frac{|x|+2}{ |x|+1}\big)^d e^{-2|x|} \xrightarrow{|x|\to \infty}0$, which is impossible. 

\vspace{2mm}

\noindent $(iii)$ Note that, although the class of  the almost decreasing unimodal L\'evy kernel is fairly large, there exists some radial L\'evy kernels which are not almost decreasing . For example, for $\beta\in [-1, 2)$ define 
\begin{align*}
\nu_\beta(h) = |h|^{-d-\beta}\,\,\rho(h)\text{ with }\,\, \rho(h) :=\Big(\frac{2+\cos |h|}{3}\Big)^{|h|^4}.
\end{align*}
 Note that $\nu_\beta$ is not almost decreasing since  $\rho(2\pi n) =1$ and $\rho(\pi (2n-1)) =3^{-\pi^4 (2n-1)^4}$ for all $n\in \mathbb{N}$. Now if $\beta\in (0,2)$, it is clear that $\nu_\beta$ is L\'evy integrable  since $\rho$ is bounded. If $-1\leq \beta \leq 0$ then $\nu_\beta$  is also L\'evy integrable since the map  $r\mapsto \rho(r)$ is in  $L^1(\R)$.  
\end{remark}

\noindent The following result is more general and provides an alternative to the Lemma \ref{lem:natural-norm-on-V} when $\nu$ is not a unimodal function. 
\medskip

\begin{lemma}\label{lem:natural-weighted-space}
Assume $\nu$ is a symmetric L\'{e}vy kernel and $\Omega\subset \mathbb{R}^d$ is an open set. For any measurable subset $K\subset \Omega$ with positive measure $|K|>0$, consider the map $ \nu_K: \Omega^c\to [0,\infty]$ with $$ \nu_K(x)= \operatorname{essinf}\limits_{y\in K}\nu(x-y).$$ Then $ \nu_K\in L^1(\Omega^c)$ and we have the continuous embeddings $$\VnuOm\hookrightarrow L^2(\Omega^c, \nu_K) \hookrightarrow L^1(\Omega^c, \nu_K) .$$
\end{lemma}

\medskip

\begin{proof}
	
	Given that $|K|>0$, for a suitable $z\in K\subset \Omega $ we have $\delta_z= \operatorname{dist}(z, \partial\Omega)>0$ and $\nu_K(x)\leq \nu(x-z)$ for almost every $x\in \Omega^c$. As $\Omega^c\subset B^c_{\delta_z}(z)$ by L\'{e}vy integrability of $\nu$, we get
	
	\begin{align*}
	\int_{\Omega^c}\nu_K(x)\d x\leq \int_{B^c_{\delta_z}(z)}\nu(x-z)\d x= \int_{B^c_{\delta_z}(0)}\nu(h)\d h<\infty.
	\end{align*}
In particular, $\nu_K\in L^1(\Omega^c)$. This implies that $L^2(\Omega^c, \nu_K)\hookrightarrow L^1(\Omega^c, \nu_K) $. Let $y\in K'\subset K$ be measurable such that $0<|K'|<\infty$ . Proceeding as in the proof of Proposition \ref{prop:unimodal-ball}, one obtains the estimate $\|u\|^2_{L^2(\Omega^c, \nu_K)}\leq C\|u\|^2_{\VnuOm}$ with $C=2|K'|^{-1}( \|\nu_K\|_{L^1(\Omega^c) }+1)$.
	%
%
	%
\end{proof}

\noindent As we have demonstrated above, the main feature of the weight $\nu_K$ consists of removing the singularity of $\nu$. Another way of defining such a weight is as follows.
\begin{definition}\label{def:nuKintegral}
Let $\nu : \R^d \setminus \{0\} \to [0, \infty]$ be the density of a symmetric L\'{e}vy measure and $K\subset\R^d$ be a measurable  set such that $|K|>0$. We write $\mathring{\nu}_K$ to denote the measurable function  defined by $ \mathring{\nu}_K: \R^d\to [0,\infty]$ such that for $x\in \R^d$,  
\begin{align*}
\mathring{\nu}_K(x)= \int_K 1\land \nu(x-y)\d y.
\end{align*}
\end{definition}

\begin{proposition}\label{prop:nuKintegral}
Assume that $\nu$ is the density of a symmetric L\'{e}vy measure and $\Omega\subset \mathbb{R}^d$ is an open set. For any measurable set $S\subset \R^d$ such that $0<|S|<\infty$ we have $ \mathring{\nu}_S \in L^1(\R^d)$. Let $K\subset \Omega$ be measurable with $|K|>0$ then the embedding $\VnuOm\hookrightarrow L^2(\R^d, \mathring{\nu}_K)$ is continuous. Moreover  If $0<|K|<\infty$ then the embeddings $\VnuOm\hookrightarrow  L^2(\R^d, \mathring{\nu}_K) \hookrightarrow L^1(\R^d,\mathring{\nu}_K)$ are continuous. 
\end{proposition}
\medskip

\begin{proof}
The integrability of $ \mathring{\nu}_S$ is readily obtained as follows
\begin{align*}
\int_{\R^d} \mathring{\nu}_S(x)\d x =  \int_{S}   \int_{\R^d} 1\land \nu(x-y)\d x\d y= |S| \int_{\R^d} 1\land \nu(h)\d h\leq |S|\int_{B_1} \d h+ |S|\int_{B_1^c} \nu(h)\d h<\infty.
\end{align*}
\noindent Let $u\in \VnuOm$ and $K\subset \Omega$, we get the estimate
\begin{align*}
\int_{\R^d} |u(x)|^2\mathring{\nu}_K(x)\d x
& = \int_{K}   \int_{\R^d}  |u(x)-u(y)+u(y)|^21\land \nu(x-y)\d x\d y\\
&\leq 2\|1\land\nu\|_{L^1(\R^d)} \int_{\Omega} |u(y)|^2 \d y+ 2\iil_{\Omega\R^d}  (u(x)-u(y))^2 \nu(x-y)\d y\,\d x\\
&\leq C\|u\|^2_{\VnuOm}.
\end{align*}

\end{proof}

\noindent Let us look at the case where $\nu$ is radial. 
\begin{proposition}
	Let $\nu$ be the density of a unimodal L\'evy measure with full support on $\R^d\setminus \{0\}$. In addition $\nu$ satisfies the doubling growth condition, i.e.  $\exists\,\,\kappa >0 $ such that $\nu(h)\leq \kappa \nu(2h) $ for all $|h|\geq 1$ then for any ball $B\subset \R^d$ we have
	\begin{align*}
	\mathring{\nu}_{B}\asymp 1\land \nu\quad \text{with} \quad\mathring{\nu}_{B}(x)= \int_{B}1\land\nu(x-y)\d y. 
	\end{align*}
Furthermore, if  $B\subset \Omega$ then the embeddings $\VnuOm\hookrightarrow L^2(\R^d, \mathring{\nu}_B)\hookrightarrow L^1(\R^d,\mathring{\nu}_B) $ are continuous. 
\end{proposition}
\begin{proof} Since $\nu$ is unimodal, let $c>0$ be a constant such that $c\nu(y)\leq \nu(x)$ if $|y|\geq |x|$.
Let $R\geq 1$ be large enough such that $B\subset B_R$. We know from  Lemma \ref{lem:natural-norm-on-V} $(ii)$ that there exists $\lambda>1$ 
such that $$\lambda^{-1}\widetilde{\nu}(h)\leq 1 \land\nu(h)\leq \lambda \widetilde{\nu}(h)\quad\text{ for all $h\in \R^d$},$$ where $\widetilde{\nu}(h)= \nu(R(1+|h|))$. 

\noindent  Assume $|x|\geq 4R$. For $y\in B$ we  have $|y|\leq R\leq 4R\leq |x|$ and 
thus $R(1+|x-y|) \leq 3R|x|$.  We also have $|x-y|\geq |x|-|y|\geq \frac{|x|}{2}$. Hence, for each $y\in B$, we have  $\frac{|x|}{2}\leq |x-y|\leq R(1+|x-y|)\leq 3R|x|$, which implies $c\nu(3R|x|)\leq \widetilde{\nu}(x-y) \leq c^{-1}\nu(\frac{|x|}{2}).$ Since $\frac{|x|}{2}\geq 2R\geq 1$, the scaling condition implies that $\nu(\frac{|x|}{2})\leq \kappa\nu(|x)$ and there exists a constant $\eta>0$ such that $\eta\,\nu(|x|)\leq \nu(3R|x|)$. 
It follows  that $$c\eta\,\nu(|x|)\leq \widetilde{\nu}(x-y) \leq c^{-1}\kappa\nu(|x|).$$

\noindent  Assume  $|x|\leq 4R$ then for each $y\in B\subset B_R$ we have $R\leq R(1+|x-y|)\leq 6R$. This implies that $c\nu(6R)\leq \widetilde{\nu}(x-y)\leq c^{-1}\nu(R)$ for all $y\in B$.

\noindent Altogether, it follows that for some constant $C>1$ and for all $y\in B$ and all $x\in \R^d$ we have 
$$C^{-1}\,(1\land\nu(x))\leq \widetilde{\nu}(x-y) \leq C(1\land\nu(x)).$$
Finally, for all $y\in B$ and all $x\in \R^d$ we get 
$$\lambda^{-1}C^{-1}\,(1\land\nu(x))\leq 1\land\nu(x-y) \leq \lambda C(1\land\nu(x)).$$
Integrating over $B$ yields, $\mathring{\nu}_B\asymp 1\land\nu$ since
$$\lambda^{-1}C^{-1}|B\,(1\land\nu(x))\leq \mathring{\nu}_B(x) \leq \lambda C(|B| 1\land\nu(x)). $$
If $B\subset \Omega$ then according to  Proposition \eqref{prop:nuKintegral} the embeddings $\VnuOm\hookrightarrow L^2(\R^d, \mathring{\nu}_B)\hookrightarrow L^1(\R^d,\mathring{\nu}_B) $ are continuous.
\end{proof}
\medskip
 
\begin{example}
	Assume $K$ is a compact subset of $\Omega$ then $\delta_K = \operatorname{dist}(K,\partial\Omega)>0$. For the particular choice $\nu(h) = |h|^{-d-\alpha}, ~\alpha\in (0,2)$ one has  that $\nu_K(x) \asymp (1+|x|)^{-d-\alpha} \asymp 1\land \nu(x)$. Indeed, for fixed $x\in \Omega^c$ and choose $y\in K$ such that $\nu_K(x) = \nu(x-y)$. Let $R\geq 1 $ be large enough such that $K \subset B_R(0)$ then $|x-y|\leq R(1+|x|)$ that is $ R^{-d-\alpha} (1+|x|)^{-d-\alpha}\leq \nu_K(x)$. We now show that upper bound in two cases. If $|x|\geq4 R$ then $|x-y|\geq |x|-|y|\geq \frac{|x|}{2}+R\geq \frac{1}{2}(1+|x|) $ so that $\nu_K(x)\leq 2^{d+\alpha}(1+|x|)^{-d-\alpha}$. If $|x|\leq 4R$ since $\delta_K \leq |x-y|$ we get $ \nu_K(x)\leq (5R\delta_K^{-1})^{-d-\alpha} (1+|x|)^{-d-\alpha}$. 
	
	\vspace{1mm}
	\noindent Analogously, we also  have $\mathring{\nu}_K(x) \asymp (1+|x|)^{-d-\alpha} \asymp 1\land \nu(x)$. 
\end{example}

\noindent From this example, the inclusion $\VnuOm \hookrightarrow L^2(\Omega^c, \nu_K)$ turns out to be a variant of the inclusion $V^{\alpha/2}(\Omega|\mathbb{R}^d) \hookrightarrow L^2(\Omega^c, \frac{1}{(1+|x|)^{d+\alpha}})$ from \cite[Proposition 13]{DyKa18}. 

\medskip

\noindent Now we resume the conclusions we have obtained following Definition \ref{def:nuKinfimum} and Definition \ref{def:nuKintegral} .
\begin{theorem}\label{thm:Vnu-in-l2}
Let  $\nu$ be the density of a symmetric L\'{e}vy  measure. Let $\Omega\subset \mathbb{R}^d$ be open. 
The following assertions are true. 
\begin{enumerate}[$(i)$]
	\item For a ball $B\subset \Omega$ we have  $\VnuOm\hookrightarrow L^2(\R^d, \nu_B)  \hookrightarrow L^1(\R^d,\nu_B).$ If $\nu$ is unimodal and satisfies the doubling growth condition then we have  $\nu_B\asymp 1\land \nu$.  The same holds for $\mathring{\nu}_B.$
	
	\item  For a measurable set $K\subset \Omega$ with $|K|>0$  we have $\VnuOm\hookrightarrow L^2(\Omega^c, \nu_K) \hookrightarrow L^1(\Omega^c, \nu_K) .$ If $K$ is compact then $\nu_K\in L^1(\R^d)$ and   $\VnuOm\hookrightarrow L^2(\R^d, \nu_K)  \hookrightarrow L^1(\R^d,\nu_K).$
	\item For a measurable set $K\subset \Omega$ with $|K|>0$ we have  $\VnuOm\hookrightarrow L^2(\R^d, \mathring{\nu}_K)$.  If $0<|K|<\infty$ then $\mathring{\nu}_K)\in L^1(\R^d)$ and $L^2(\R^d, \mathring{\nu}_K)\hookrightarrow L^1(\R^d, \mathring{\nu}_K).$ 
\end{enumerate}
\end{theorem}


\subsection{Nonlocal trace spaces}

The main goal of this part is to design an abstract notion of a trace space of $\VnuOm$ similarly as one does for the space $H^1(\Omega)$. Due the nonlocal feature of the space $\VnuOm$, the trace space assumes functions defined on the complement of $\Omega$. Indeed, one reason is that elements of $\VnuOm$ are defined on the whole $\mathbb{R}^d$. This contrasts with the local situation, where the trace space of $H^1(\Omega)$ (with $\Omega$  sufficiently smooth) consists of  elements defined on the boundary $\partial\Omega$. In this part, we assume that $\nu $ is fully supported on $\mathbb{R}^d$ and $\VnuOm$ is solely endowed with the norm $\|\cdot \|_{\VnuOm} $. We define $\TnuOm$, as  the space of restrictions to $\Omega^c$ of functions of $\VnuOm$. That is, 
\begin{align*}
\TnuOm = \{v: \Omega^c\to \mathbb{R}~\text{meas.} ~~\hbox{ such that }~~ v = u|_{\Omega^c} ~~\hbox{with }~~ u \in \VnuOm\}.
\end{align*}
We endow $\TnuOm $ with its natural norm, 
\begin{align}
\|v\|_{\TnuOm } = \inf\{ \|u\|_{\VnuOm }: ~~ u \in \VnuOm ~~ \hbox{ with }~~ v = u|_{\Omega^c} \}. 
\end{align}
As an immediate consequence of the definition of the space $\big(\TnuOm, \|\cdot\|_{\TnuOm}\big)$, the trace map
$\operatorname{Tr}: \VnuOm \to \TnuOm$ with $u \mapsto \operatorname{Tr}(u) = u\mid_{\Omega^c}$ is continuous, onto, i.e. $\operatorname{Tr}(\VnuOm)= \TnuOm$, $ \ker(\operatorname{Tr} ) =  \VnuOmO$ and satisfies $\|\operatorname{Tr}( u)\|_{\TnuOm}\leq \|u\|_{\VnuOm}$ for all $u\in \VnuOm.$
\medskip

\begin{theorem}
	The space is a separable Hilbert space with the scalar product
	\begin{align*}
	(u,v)_{\TnuOm}= \frac{1}{2}\left( \|u+v\|^{2}_{\TnuOm }-\|u\|^{2}_{\TnuOm }-\|v\|^{2}_{\TnuOm } \right).
	\end{align*}
\end{theorem}
\begin{proof}
	Clearly the norm $\|\cdot \|_{\TnuOm }$ verifies the parallelogram law since the norm $\|\cdot \|_{\VnuOm }$ does. It follows that $\big(\cdot, \cdot\big)_{\TnuOm}$ is a scalar product on $\TnuOm$ with associated norm $\|\cdot \|_{\TnuOm }$ . We want to prove  that $\TnuOm$ is complete under the norm $\|\cdot \|_{\TnuOm }$. Let $(u_n)_n$ be a Cauchy sequence in $\TnuOm$ then up to extraction of a subsequence we may assume that
	\begin{align*}
	\|u_n-u_{n+1}\|_{\TnuOm } <\frac{1}{2^{n+1}}\qquad\hbox{for all }~~~n\geq 1\,. 
	\end{align*}
	Let us fix $\overline{u}_1\in \VnuOm $ such that $u_1= \overline{u}_1|_{\Omega^c}$. By definition of $\|u_1-u_2\|_{\TnuOm}$ there exists $v \in \VnuOm$ such that $u_1-u_2 = v|_{\Omega^c}$ and 
$	\|v\|_{\VnuOm } <\|u_1-u_2\|_{\TnuOm } +4^{-1} \,.$ Letting $\overline{u}_2= v+ \overline{u}_1$, i.e. $v= \overline{u}_1-\overline{u}_2$ yields $\overline{u}_2\in \VnuOm$, $u_2= \overline{u}_2|_{\Omega^c}$ and hence
$\|\overline{u}_1-\overline{u}_2\|_{\VnuOm } <\|u_1-u_2\|_{\TnuOm } +4^{-1}<2^{-1}\,. $ Repeating this process, one constructs a sequence of functions $\overline{u}_n\in \VnuOm $ such that
	\begin{align*}
	\|\overline{u}_n-\overline{u}_{n+1}\|_{\VnuOm } <\frac{1}{2^{n}}\qquad\hbox{for all}~~~ n\geq 1,
	\end{align*}
	which turns out to be a Cauchy sequence in the complete space $\VnuOm$. 
	Let $\overline{u}\in \VnuOm$ be the limit of $(\overline{u}_n)_n$. 
	Clearly, setting $u = \overline{u}\mid_{\Omega^c} $ we have $
	\|u_n-u \|_{\TnuOm } \leq \|\overline{u}_n-\overline{u}\|_{\VnuOm } \xrightarrow{n\to \infty} 0.$
Finally, the original Cauchy sequence $(u_n)_n$ converges up to extraction of subsequence to $u$ and hence converges itself to $u$.
\end{proof}

\vspace{2mm}

\noindent
It is natural to ask the following question: Can the space $\TnuOm$ be defined with an intrinsic scalar product preserving its initial Hilbert structure such that its trivial connection to $\VnuOm$ is less evident? In the local situation, it is possible to define a scalar product on the space $H^{1/2}(\partial \Omega)$ when $\Omega$ is a special Lipschitz domain (see \cite{Din96}) . We give an answer to this question provided that some regularity and growth bound conditions on $\nu$ are assumed. We follow \cite{BGPR17} where the authors enforce the following assumptions.
\begin{align*}
\intertext{\textbf{A1}: $\nu$ is unimodal, twice continuously differentiable and there is a constant $C_1>0$ such that}
&|\nu'(r)|, |\nu''(r)|\leq C_1 \nu(r).
\intertext{\textbf{A2}: There exist constants $\beta\in (0, 2)$ and $C > 0$ such that}
& \nu(\lambda r)\leq C \lambda^{d-\beta}\nu(r),\,\, 0< r,\lambda\leq 1,\quad \text{and}\quad\nu(r)\leq C\nu(r+1),\,\,\, r\geq 1.
\end{align*}
Assume $\Omega^c$ satisfies the volume density condition, i.e. $\exists\, c>0$ such that $|\Omega^c\cap B_r(x)|\geq cr^d$ for all $x\in \partial\Omega$ and all $r>0$. Then under assumptions, \textbf{A1} and \textbf{A2} \cite[Theorem 2.3]{BGPR17} reveals that for any $g\in \TnuOm$, there exists a unique $u_g\in\VnuOm $ such that $u_g|_{\Omega^c}= g$ and
\begin{align}\label{eq:douglas-formula}
\mathcal{H}_\Omega(g,g):= \iil_{\Omega^c\Omega^c} \!\!\big(g(x)-g(y) \big)^2 \, \gamma_\Omega (x, y) \d y\,\d x=\hspace*{-2ex} \iil_{(\Omega^c\times \Omega^c)^c} \!\!\big(u_g(x)-u_g(y) \big)^2 \, \nu (x-y) \d y\,\d x. 
\end{align}
Furthermore, $u_g$ satisfies the weak formulation 
\begin{align}
\iil_{(\Omega^c\times \Omega^c)^c} \!\!\big(u_g(x)-u_g(y) \big)(\phi(x)-\phi(y)) \, \nu (x-y) \d y\,\d x=0,\qquad \text{for all}~~~\phi\in  \VnuOmO. 
\end{align}
The interaction kernel $\gamma_\Omega (x, y)$ is given via the Poisson kernel $P_\Omega(\cdot, \cdot)$ of $\Omega$ by the formula
\begin{align*}
\gamma_\Omega (x, y) = \int_{\Omega}P_\Omega(x,z)\nu(z-y) \d z\,\qquad x,y\in \Omega^c.
\end{align*}
Furthermore, a precise formula for $u_g$ in $\Omega $ is given by the Poisson integral
\begin{align*}
u_g(x)= P_\Omega[g](x) = \int_{\Omega^c}g(y)P_\Omega(x,y)\d y\, \qquad x\in \Omega\,. 
\end{align*}
%
%
\noindent Now let $v\in \TnuOm$ by definition of $\|\cdot\|_{\TnuOm }$ we have 
\begin{align*}
\|v\|^{2}_{\TnuOm } &= \inf\{ \|u\|^{2}_{\VnuOm } :~~ u \in \VnuOm ~~ \hbox{with}~~ v = u|_{\Omega^c} \}\\
&= \inf\Big\{ \int_{\Omega}u^2(x)\d x: ~~ u \in \VnuOm ~~ \hbox{with}~ v = u|_{\Omega^c} \Big\}+ \mathcal{H}_{\Omega}(v,v)\,.
\end{align*}
It is rather challenging to find or to estimate the quantity 
\begin{align*}
\inf\Big\{ \int_{\Omega}u^2(x)\d x: ~~ u \in \VnuOm ~~ \hbox{with}~ v = u|_{\Omega^c} \Big\}.
\end{align*}
Let us remind that our goal here is to explicitly define a norm which is equivalent to $\|\cdot\|_{\TnuOm }$ and has less visible connection to $\VnuOm$. To this end, we bring into play the norm $\|\cdot\|^{*}_{\VnuOm }$ defined as in Lemma \ref{lem:natural-norm-on-V}. 
\vspace{2mm}

\begin{proposition}
	Assume $\Omega$ is open and bounded, such that $\Omega^c$ satisfies the volume density condition. Assume $\nu$ satisfies conditions \textbf{A1} and \textbf{A2} (in particular $\nu$ is unimodal and  has full support on $\mathbb{R}^d$). Let $\widetilde{\nu}$ and $\|\cdot\|^{*}_{\VnuOm }$ be the measure and the norm respectively given in Lemma \ref{lem:natural-norm-on-V}. Then, 
	\begin{align*}
	\TnuOm= \Big\{ v:\Omega^c\to \mathbb{R}~\text{meas.}~~\mathcal{H}_{\Omega}(v,v)= \iil_{\Omega^c\Omega^c} \!\!\big(v(x)-v(y) \big)^2 \, \gamma_\Omega (x, y) \d y\,\d x<\infty \Big\}\,
	\end{align*}
	and the norms $\|\cdot\|_{\TnuOm }$, $\|\cdot\|^{*}_{\TnuOm }$ and $\|\cdot\|^{\dagger}_{\TnuOm }$ are all equivalent. The norms are defined by 
	\begin{align*}
	\|v\|^{*}_{\TnuOm } &= \inf\{ \|u\|^{*}_{\VnuOm } :~~ u \in \VnuOm ~~ \hbox{with}~~ v = u|_{\Omega^c} \}\\
	\|v\|^{\dagger 2}_{\TnuOm } &=\int_{\Omega^c} v^2(x)\widetilde{\nu}(x)\d x+ \iil_{\Omega^c\Omega^c} \!\!\big(v(x)-v(y) \big)^2 \, \gamma_\Omega (x, y) \d y\,\d x\,.
	\end{align*}

\end{proposition}
\begin{proof}
	The equivalence between $\|\cdot\|_{\TnuOm }$ and $\|\cdot\|^{*}_{\TnuOm }$ is an immediate consequence of Lemma \ref{lem:natural-norm-on-V} $(iv)$. By \eqref{eq:douglas-formula} it follows that
	\begin{align*}
	\|v\|^{*2}_{\TnuOm } &= \inf\{ \|u\|^{*2}_{\VnuOm } ~~ u \in \VnuOm ~~ \hbox{with}~~ v = u|_{\Omega^c} \}\\
	&= \inf\Big\{ \int_{\mathbb{R}^d}u^2(x)\widetilde{\nu}(x) \d x ~~ u \in \VnuOm ~~ \hbox{with}~ v = u|_{\Omega^c} \Big\}+ \mathcal{H}_{\Omega}(v,v)\\
	&\geq \int_{\Omega^c} v^2(x)\widetilde{\nu}(x)\d x+ \mathcal{H}_{\Omega}(v,v). 
	\end{align*}
	We have $\|v\|^{\dagger}_{\TnuOm }\leq \|v\|^{*}_{\TnuOm }$. Hence the identity $Id: (\TnuOm, \|\cdot\|^{*}_{\TnuOm }) \to (\TnuOm, \|\cdot\|^{\dagger}_{\TnuOm })$ is continuous. The space $(\TnuOm, \|\cdot\|^{*}_{\TnuOm })$ is a Hilbert space since $\|\cdot\|_{\TnuOm }$ and $\|\cdot\|^{*2}_{\TnuOm }$ are equivalent. Also, using Fatou's lemma one can easily show that $(\TnuOm, \|\cdot\|^{\dagger}_{\TnuOm })$ is a Hilbert space. As a consequence of the open mapping theorem the norms $ \|\cdot\|^{\dagger}_{\TnuOm }$ and $ \|\cdot\|^{*}_{\TnuOm }$ are equivalent. 
	
\end{proof}

	\noindent In the special case $\nu(h)= (2-\alpha)|h|^{-d-\alpha}$, noting $\delta_z = \operatorname{dist}(z, \partial\Omega)$, the authors in \cite[Theorem 3]{DyKa18} claim that if $v\in \VnuOm$ then
	\begin{align}\label{eq:finiteness-ext}
	\iil\limits_{\Omega^c \Omega^c} \frac{\big(v(x)-v(y) \big)^2 }{(|x-y|+\delta_x+\delta_y)^{-d-\alpha}} \d x\d y<\infty\,.
	\end{align}
Conversely if \eqref{eq:finiteness-ext} holds true for $v=g$ on $\Omega^c$ then there exists $u_g\in \VnuOm$ such that $u_g|_{\Omega^c =g}$ and 
	\begin{align}\label{eq:compare-semi-douglas}
	(2-\alpha) \iil\limits_{(\Omega^c \times \Omega^c)^c} \frac{\big(u_g(x)-u_g(y) \big)^2 }{(|x-y|+\delta_x+\delta_y)^{-d-\alpha}} \d y\,\d x\asymp \iil\limits_{\Omega^c \Omega^c} \frac{\big(v(x)-v(y) \big)^2 }{(|x-y|+\delta_x+\delta_y)^{-d-\alpha}} \d x\d y
	\end{align}
	with the constants  independent of $g$ and $u_g$. Therefore, it readily follows that 
	\begin{align*}
	\TnuOm= \Big\{ v:\Omega^c\to \mathbb{R}~\text{meas.}~~ \iil_{\Omega^c\Omega^c} \!\!\frac{\big(v(x)-v(y) \big)^2 }{(|x-y|+\delta_x+\delta_y)^{-d-\alpha}}\, \d y\,\d x<\infty \Big\}\,. 
	\end{align*}
Actually it is not straightforward to conclude on the equivalence of norms in this case. In view of  \eqref{eq:compare-semi-douglas} we only have $\|\cdot\|_{\TnuOm}\leq C \|\cdot\|^{'}_{\TnuOm}$ with the norm 
	\begin{align*}
	\|v\|^{'2}_{\TnuOm} = \int_{\Omega^c} \frac{v^2(x)}{(1+|x|)^{d+\alpha}}\d x+ \iil\limits_{\Omega^c \Omega^c} \frac{\big(v(x)-v(y) \big)^2 }{(|x-y|+\delta_x+\delta_y)^{-d-\alpha}} \d x\d y\,.
	\end{align*}
	

\begin{remark}
	Let us emphasize the nonlocal trace $\operatorname{Tr}$ does not need any special construction via the  
	functional analysis and density argument. Since $\Omega^c$ is still a d-dimensional manifold, it makes sense to consider the restriction of a measurable function on $\Omega^c.$ 
	Moreover no regularity on $\Omega$ is required. Whereas in the local situation, the trace of a Sobolev function 
	$u$ on the boundary $\partial\Omega$ requires the smoothness of both $u$ and $\partial\Omega$.
\end{remark}

\noindent Let us give an important result concerning the notion of nonlocal trace, thereby providing some analogies with the classical notion traces.
\begin{proposition}\label{prop:nonlocal-L2-trace} Let $K\subset \Omega$ be measurable with $|K|>0$. The following properties are true.
	\begin{enumerate}[$(i)$]
		\item The trace map $\operatorname{Tr}: \VnuOm \to L^2(\Omega^c, \nu_K)$ with $u \mapsto \operatorname{Tr}(u) = u\mid_{\Omega^c}$ is linear and continuous.
		\item The trace map $\operatorname{Tr}: \VnuOm \to L^2(\Omega^c, \mathring{\nu}_K)$ with $u \mapsto \operatorname{Tr}(u) = u\mid_{\Omega^c}$ is linear and continuous.
		\item In both cases we have $ \ker(\operatorname{Tr} ) =  \VnuOmO $ and  $\operatorname{Tr}(\VnuOm)= \TnuOm$.
	\end{enumerate}
\end{proposition}
\begin{proof}
	This is a straightforward consequence of Theorem \ref{thm:Vnu-in-l2}. 
	
\end{proof}
\begin{remark}
	 Let us observe that the interplay between , $\VnuOm$, $  \VnuOmO$, $\TnuOm$ and   $L^2(\Omega^c, \nu_K)$ may be view as the nonlocal counterpart of that between $H^1(\Omega)$, $H^1_0(\Omega),$ $H^{1/2}(\partial\Omega)$ and $L^2(\partial\Omega)$. Indeed, it is well known that if the classical trace operator $\gamma_0 :H^1(\Omega)\to L^2(\partial\Omega)$ exists then we have
\begin{itemize}
	\item $\gamma_0 :H^1(\Omega)\to L^2(\partial\Omega)$  is linear and continuous. 
\item $\gamma_0(H^1(\Omega))= H^{1/2}(\partial\Omega)$ and $\ker(\gamma_0)= H^1_0(\Omega)$.
\end{itemize}
 In this regard it is fair to view  $\VnuOm$, $  \VnuOmO$, $\TnuOm$ and   $L^2(\Omega^c, \nu_K)$  respectively as the nonlocal replacement of $H^1(\Omega)$, $H^1_0(\Omega),$ $H^{1/2}(\partial\Omega)$ and $L^2(\partial\Omega)$. 
\end{remark}

\medskip
\begin{proposition}
Assume $\Omega$ is a Lipschitz bounded domain.  Then $T: H^{1/2}(\partial\Omega)\to \TnuOm $ with $Tu = \operatorname{Tr}\circ E\circ \operatorname{Ext} u$ is a linear operator bounded.

\noindent Here $E : H^1(\Omega)\to H^1(\R^d)$ is the Sobolev extension, i.e.  $Eu |_\Omega = u$ and $\|Eu\|_{H^1(\R^d)}\leq C \|u\|_{H^1(\Omega)}$  for all $u\in H^1(\Omega)$.  And $\operatorname{Ext}: H^{1/2}(\partial\Omega)\to H^1(\Omega)$ is the trace-extension operator, i.e. $\gamma_0\circ \operatorname{Ext} u = u $ and $\| \operatorname{Ext} u\|_{H^{1/2}(\partial \Omega)}\leq C \|u\|_{H^1(\Omega)}$ for all $u\in H^{1/2}(\partial \Omega)$.
\end{proposition}
\begin{proof}
The proof is immediate. 
\end{proof}

\begin{proposition}
	Let $C^\infty_c(\overline{\Omega^c} )=C^\infty_c(\mathbb{R}^d)|_{\Omega^c} $ be set of restrictions on $\Omega^c$ of $C^\infty$ functions on $\mathbb{R}^d$ with compact support. Then $C^\infty_c(\overline{\Omega^c} )$ is dense in $\TnuOm$. 
\end{proposition}

\begin{proof}
	For $v \in \TnuOm$ we write $v= u|_{\Omega^c}$ with $u \in \VnuOm$. According to Theorem \ref{thm:density} there exists $u_n\in C^\infty_c(\mathbb{R}^d)$ such that $\|u_n -u\|_{\VnuOm}\to 0.$ Put, $v_n= u_n|_{\Omega^c}$ we get
	\begin{align*}
	\|v_n -v\|_{\TnuOm}\leq \|u_n -u\|_{\VnuOm}\xrightarrow{n\to \infty} 0.
	\end{align*}
	
\end{proof}

\vspace{1mm}
\noindent Note the recent result on the nonlocal trace from \cite{Zoran19} is rather restrictive compared to our existing results. Indeed therein the space $ L^2(\Omega^c, \mu)$ appears to be the trace of 
$\VnuOm \cap L^2(\mathbb{R}^d, m)$ where $m(x) =\mathds{1}_{\Omega}(x)+ \mu(x)\mathds{1}_{\Omega^c}(x) $ and $\mu(x) = \int_{\Omega}\nu(x-y)\d y$ for $x \in \Omega^c\,.$ The corresponding extension operator defined from $ L^2(\Omega^c, \mu)$ to $\VnuOm \cap L^2(\mathbb{R}^d, m)$ is merely the extension operator by zero outside $\Omega$.  The main defect of this approach is that the measure $m$ is too singular across the boundary of $\Omega$. Another point is that our trace space $\TnuOm$ turns out to be larger than $L^2(\Omega^c, \mu)$. In the sense that $L^2(\Omega^c, \mu)$ is continuously embedded in $\TnuOm$. Indeed, for $v\in L^2(\Omega^c, \mu) $ its zero extension $v_0= v\mathds{1}_{\Omega^c} $ belongs to $\VnuOm$ since we have
	\begin{align*}
	\iil_{(\Omega^c\times \Omega^c)^c} \!\!\big(v_0(x)-v_0(y) \big)^2\, \nu (x-y) \d y\,\d x= 2\int_{\Omega^c} v^2(x)\mu(x)\d x.
	\end{align*}
\noindent So that $v \in \TnuOm$ and the continuity is obtained as follows	$\|v\|_{\TnuOm}\leq \|v_0\|_{\VnuOm}\leq \sqrt{2}\|v\|_{L^2(\Omega^c, \mu)}.$

\section{Nonlocal Sobolev-like spaces}

\noindent An important role in our study is played by function spaces. We assume that $\Omega$ an open subset of $\R^d$. Let us introduce generalized Sobolev-Slobodeckij-like spaces with respect to a symmetric $p$-L\'{e}vy measure $\nu(h) \d h$ ($1\leq p<\infty$)  We will show some strong connections with the classical Sobolev spaces. Our standing assumption is that $\nu: \R^d\setminus\{0\}\to [0,\infty]$ satisfies
\begin{align}\tag{$J_1$}\label{eq:plevy-integrability-bis}
\nu(-h) = \nu(h) ~~\text{ for all $h\in \mathbb{R}^d$ and} ~~\int_{\mathbb{R}^d}(1\land |h|^p)\nu(h)\d h<\infty\,.
\end{align} 
In case $\nu$ is radial we adopt the convention by identifying $\nu$ with its radial profile, i.e. $\nu(h) = \nu(|h|), h\in \R^d.$ For several results $\nu $ will be assumed to have full support. Let $\Omega\subset \mathbb{R}^d$ be open.

\noindent $\bullet$ We define the space $\WnuOm$ as 
\begin{align*}
 \WnuOm= \Big\{u \in L^p(\Omega): \,\, |u(x)-u(y)|\nu^{1/p}(x-y)\in L^p(\Omega\times \Omega) \Big \}\,.
\end{align*}
equipped with the norm defined as follows
\begin{align*}
&\|u\|^p_{ \WnuOm}= \|u\|^p_{L^p (\Omega)}+ |u|^p_{\WnuOm}\quad\text{with}\quad 
|u|^p_{\WnuOm}:=  \iil_{\Omega\Omega} \big|u(x)-u(y) \big|^p \, \nu (x-y) \d y\,\d x. 
\end{align*}
Note that norm $\|\cdot\|^p_{ \WnuOm}$ is equivalent to norm defined by 
\begin{align*}
u\mapsto  \|u\|^p_{L^p (\Omega)}+ \iil_{\Omega\Omega\cap \{|x-y|\leq \delta\}} \big|u(x)-u(y) \big|^p \, \nu (x-y) \d y\,\d x\qquad\text{for all}~~~\delta>0. 
\end{align*}
Indeed it suffice to observe that letting $C_\delta = \int_{|h|\geq \delta}\nu(h)\d h>0$
\begin{align*}
 \iil_{\Omega\Omega\cap \{|x-y|\geq \delta\}} \big|u(x)-u(y) \big|^p \, \nu (x-y) \d y\,\d x\leq 2^p\int_\Omega|u(x)|^p\d x\int_{|x-y|\geq \delta}\nu(x-y)\d y=C_\delta\|u\|^p_{L^p (\Omega)}.
\end{align*}
\noindent $\bullet$ We also introduce the space $ \WnuOmR$ defined as follows
\begin{align*}
\WnuOmR = \Big\lbrace u: \R^d \to \R \text{ meas. } : \,\,\, |u(x)-u(y)|\nu^{1/p}(x-y)\in L^p(\Omega\times \R^d) \Big\rbrace \,.
\end{align*}
 
\noindent The space $\WnuOmR $ shall be equipped with the following norm
\begin{align*}
&\|u\|^p_{ \WnuOmR}= \|u\|^p_{L^p (\Omega)}+ |u|^p_{\WnuOmR}\quad\text{with}\quad|u|^p_{\WnuOmR}:=  \iil_{\Omega\R^d} \big|u(x)-u(y) \big|^p \, \nu (x-y) \d y\,\d x. 
\end{align*}
Likewise as above, the norm $\|\cdot\|^p_{ \WnuOmR}$ is equivalent to the norm defined by 

\begin{align*}
u\mapsto  \|u\|^p_{L^p (\Omega)}+ \iil_{\Omega\R^d\cap \{|x-y|\leq \delta\}} \big|u(x)-u(y) \big|^p \, \nu (x-y) \d y\,\d x\qquad\text{for all}~~~\delta>0. 
\end{align*}
\pagebreak[2]
Another possibility is to consider the norm
\begin{align*}
\vertiii{u}^p_{\WnuOmR}= \|u\|^p_{L^p (\R^d)}+ |u|^p_{\WnuOmR}.
\end{align*}

\vspace{1mm}
\noindent $\bullet$  The space $\WnuOmO$ is the space of functions that vanish on the complement of $\Omega$ i.e
\begin{align*}
\WnuOmO= \{ u\in \WnuOmR~: ~u=0~~\text{a.e. on } \mathbb{R}^d\setminus \Omega\}\,.
\end{align*} 
We set $$\|u\|_{ \WnuOmO}= \|u\|_{ \WnuOmR}.$$

\begin{remark}
	\begin{enumerate}[$(i)$]
	\item  For $p=2$ the spaces $\WnuOm, \WnuOmR$ and $\WnuOmO$  become the spaces $\HnuOm,\VnuOm$ and $ \VnuOmO$ of Section \ref{sec:function-spaces} respectively.
	
	\item  It is a  little exercise  to show  that the investigations from Section \ref{sec:function-spaces}  extend to the case $p\neq 2$. 
	\item For $\Omega= \R^d$ the spaces $\WnuOmR$, $\WnuOmO$ and $\WnuOm$ all coincide. We shall denote by $W^p_\nu(\mathbb{R}^d)$ or $H_\nu(\mathbb{R}^d)$
	if $p=2$. 
	
	\item Obviously the norms $\|\cdot \|_{\WnuOmR} $, $ \|\cdot \|_{W^p_\nu(\mathbb{R}^d)}$ and $\vertiii{\cdot}_{\WnuOmR} $ agree on $ \WnuOmO$. Furthermore, $ \WnuOmO$ is a closed subspace of $\WnuOmR$. 
	\item  Note that $h\mapsto |h|^{-d-sp}$ belongs to $L^1(\R^d, 1\land |h|^p)$ if and only if $0<s<1$. In this case taking
	$\nu(h) = |h|^{-d-sp} $ the space $\WnuOm$ turns out to be the fractional Sobolev space $W^{s,p}(\Omega)$. The space $\WnuOmR$ shall be denoted by $W^{s,p}(\Omega|\R^d)$. 
	\item Along the lines of Section \ref{sec:function-spaces}  one can analogously establish that  $u\in L^p(\R^d, (1+|h|)^{-d-sp} \d h)$ for any $u\in \WnuOmR$ provided that $|\Omega|>0$ and $|\R^d \setminus\Omega|>0$. 
	\item  The integrability of the increments $|u(x)-u(y)|^p\nu(x-y)$ over $\Omega\times \R^d$ encodes certain regularity of the function $u$ in $\Omega$ and across the boundary of $\Omega$. We concretely deal with the regularity inside $\Omega$ in Section \ref{sec:charact-W1p} under certain conditions. To picture the regularity across the boundary, for the instance case $\nu(h) = |h|^{-d-sp} $ with $sp\geq 1$, let us assume $u$ is continuous and $u=0$ on $\Omega^c$. If $u\in \WnuOmR$ then
	\begin{align*}
|u|^p_{\WnuOmR}&\geq 2\int_\Omega |u(x)|^p\int_{\Omega^c}|x-y|^{-d-sp} \d y= C \int_\Omega |u(x)|^p \operatorname{dist}(x, \partial \Omega)^{-sp} \d x.
	\end{align*}
	
\noindent For the case $sp\geq 1$, this implies that $u(x)\xrightarrow{\dist(x, \partial \Omega)\to0}0$.
	\end{enumerate}
\end{remark}

\noindent Here we highlight  the reason why it is crucial to assume $\nu$ to be of full support. 

\begin{proposition} Assume $\nu$ has full support then $|u|_{\WnuOmR}=0$ if and only if $u=c$ a.e on $\R^d$
 $(c\in \R)$. In particular, $(\WnuOmR, \|u\|_{\WnuOmR})$ is a normed space. This conclusion may be false if $\nu$ is not of full support. 
\end{proposition}

\vspace{1mm}
\begin{proof} 
In the case $ |u|_{\WnuOm} = 0,$  there exists $N\subset \Omega$ with $|N|=0$ and for all  $\Omega\setminus N$ we have 
\begin{align*}
0=\int_{\R^d}|u(x)-u(y) |^p\nu(x-y)\d y=  \int_{\R^d}|u(x)-u(x+h) |^p\nu(h)\d h. 
\end{align*}
Since $\nu(h)>0$ a.e, for $a\in \Omega\setminus N$ we have $u(a+h)= u(a)$ for almost all $h\in \R^d$.  That $u$ is almost everywhere  constant on $\R^d$. 
If  $ \|u\|_{\WnuOm} = 0,$ then $|u|^p_{\WnuOmR}=0$, $u$ is constant almost everywhere and since  $u=0$ a.e on $\Omega$ it follows that $u=0$ a.e on $\R^d$. This enables $\|\cdot\|_{\WnuOmR}$  to be a  norm on $\WnuOmR. $ 

\vspace{2mm}
\noindent  Next assume that $\Omega$ is bounded and $\nu$ has a compact support. Let $S= \R^d\setminus \big(\Omega\cup \operatorname{supp} \nu+ \Omega\big)$ and consider the function $u(x)= \mathds{1}_S(x)$. A routine verification shows that $\|u\|_{\WnuOmR}=0$ but $u\neq0$. This means that $(\WnuOmR, \|\cdot\|_{\WnuOmR})$ cannot be a normed space.
\end{proof}

\vspace{2mm}
\noindent It is noteworthy to mention that the assumption \eqref{eq:plevy-integrability} is optimal in the sense of the following proposition.
\begin{proposition} Assume that $\nu: \R^d\to [0,\infty]$ is symmetric. The following assertions hold true. 
\begin{enumerate}[$(i)$]
\item If $\nu\in L^1(\R^d)$, then $\WnuOm= L^p(\Omega)$ and $\WnuOmR\cap L^p(\R^d)= L^p(\R^d)$.

\item If $\nu \in L^1(\R^d, 1\land |h|^p)$ and $\Omega$ is bounded, then $\WnuOm$  and $\WnuOmR$ contain all bounded Lipschitz functions. 
\item Assume $\int_{B_\delta} |h|^p\nu(h)\d h= \infty$ for all $\delta>0$
(in particular  $\nu \not\in L^1(\R^d, 1\land |h|^p)$), $\nu$ is radial and $\Omega$ is connected.  Then the only smooth functions contained in $\WnuOm$ are constants. 
\item If $\nu \in L^1(\R^d, 1\land |h|^p)$ and $\nu$ is radial then there exists two constant $C_1, C_2>0$ such that for all $u\in W^{1,p} (\R^d)$ we have $ C_1\|u\|_{W^{1,p} (\R^d)} \leq \|u\|_{W^{p}_\nu (\R^d)}\leq  C_2\|u\|_{W^{1,p} (\R^d)} $.
\end{enumerate}
\end{proposition}

\medspace

\begin{proof}  If $\nu\in L^1(\R^d)$ then $(i)$  is obtained through the following estimate 
\begin{align*}
&\iil_{\Omega\Omega} \big |u(x)-u(y) \big|^p \, \nu (x-y)\d y\,\d x\leq \iil_{\Omega\R^d} \big |u(x)-u(y) \big|^p \, \nu (x-y)\d y\,\d x\\
&\leq 2^p \iil_{\Omega\R^d }|u(x)|^p \, \nu (x-y)\d y\,\d x
\leq 2^p \iil_{\Omega\mathbb{R}^d} |u(x)|^p \, \nu (h)\mathrm{d}x\,\mathrm{d}h
=2^p\|\nu\|_{L^1(\mathbb{R}^d)}\|u\|_{L^p(\Omega)}\,.
\end{align*}

\noindent For a bounded Lipschitz function  $u$, there is a constant $C>0$ such that $|u(x) -u(y)|\leq C(1\land|x-y|)$ for all $x,y\in \R^d$. Hence it clearly follows that $u\in \WnuOm$ and $u\in \WnuOmR$, i.e. $(ii)$ is true.

\noindent Now assume $\int_{B_\delta} |h|^p\nu(h)\d h= \infty$ for all $\delta>0$ and let $u\in C_c^\infty(\Omega)$.  Assume $u\in \WnuOm$, let $ K \subset\Omega $ be a compact and let $\delta>0$ such that $K(\delta)= K+B_\delta(0)\subset \Omega$. Using the fundamental theorem of calculus and passing through polar coordinates yields
\begin{align}\label{eq:optimal-non-levy}
\begin{split}
|u|_{\WnuOm}^p
&= \iil_{ \Omega\Omega}|u(x)-u(y)|^p\nu(x-y)\,\d y\,\d x\\
&\geq \int_{K} \int_{|x-y|\leq \delta} \left| \int_0^1\nabla u(x+t(y-x))\cdot \frac{(y-x)}{|y-x|}\d t\right|^p |x-y|^{p}\nu(x-y) ,\d y\,\d x\\
&= \int_{K(\delta)} \int_{\mathbb{S}^{d-1}} \left|\nabla u(x)\cdot w\right|^p \mathrm{d}\sigma(w) \int_{0}^{\delta} r^{p+d-1} \nu(r)\mathrm{d}r\\
&= |\mathbb{S}^{d-1}|^{-1}\Big( \int_{K(\delta)} \int_{\mathbb{S}^{d-1}} \left|\nabla u(x)\cdot w\right|^p \mathrm{d}\sigma(w) \d x\Big) \Big(\int_{B_\delta} |h|^p\nu(h)\d h\Big).
\end{split}
\end{align}

\noindent However, this is possible only if $\nabla u = 0$ since  $\int_{B_\delta} |h|^p\nu(h)\d h= \infty$. Hence  $u$ must be a constant function. To prove $(iv)$, observe that, taking $K=\R^d$ in \eqref{eq:optimal-non-levy} and noting that the Lebesgue measure is invariant by rotation we have 
\begin{align*}
|u|^p_{W^p_\nu(\R^d)} &\geq|\mathbb{S}^{d-1}|^{-1}\Big( \int_{\R^d} \int_{\mathbb{S}^{d-1}} \left|\nabla u(x)\cdot w\right|^p \mathrm{d}\sigma(w) \d x\Big) \Big(\int_{B_\delta} |h|^p\nu(h)\d h\Big)\\
&=  K_{p,d} \Big(\int_{B_\delta} |h|^p\nu(h)\d h\Big)\int_{\R^d}|\nabla u(x)|^p\d x:= C_\delta\|\nabla u\|^p_{L^p(\R^d)}. 
\end{align*}

\noindent With $K_{p,d} = \fint_{\mathbb{S}^{d-1}} |w\cdot e|^p \d\sigma_{d-1}(w)$ (see  Lemma \ref{lem:BBM-regular}). From this,  the estimate  $C_1\|u\|_{W^{1,p} (\R^d)} \leq \|u\|_{W^{p}_\nu (\R^d)}$ follows. The reverse estimate is a direct consequence of the estimate \eqref{eq:levy-p-estimate} below.
\end{proof}

\vspace{1mm}

\noindent Let us now see some fundamental properties of the spaces under consideration. 
\begin{theorem}\label{thm:completness-Wnup}
Assume that $\nu$ satisfies the condition \eqref{eq:plevy-integrability}. The function spaces  $\big(\WnuOmR, \vertiii{\cdot}_{\WnuOmR}\big)$ and $\big(\WnuOm, \|\cdot\|_{\WnuOm}\big)$ are separable Banach (Hilbert for $p=2$) spaces and reflexive for $1<p<\infty$. 

\noindent  In addition, if $\nu$ has full support in $\R^d$, then the same is true for the space $\big(\WnuOmR, \|\cdot\|_{\WnuOmR} \big)$.
\end{theorem}  

\medskip 
\begin{proof}
	It is not difficult to check that $ \|\cdot\|_{\WnuOmR}$ and $\vertiii{\cdot}_{\WnuOmR} $ are norms on  $\WnuOmR$ and $\WnuOm$ respectively. We know that $\|\cdot\|_{\WnuOmR}$ is a  norm on $\WnuOmR. $ 
	
	\medskip
\noindent Now, let $(u_n)_n$ be a Cauchy sequence in $\big(\WnuOmR, \vertiii{\cdot}_{\WnuOmR}\big)$. It converges to some $u$ in the topology of $L^p(\mathbb{R}^d)$ and pointwise almost everywhere in $\mathbb{R}^d$ up to a subsequence $(u_{n_k})_k$. Fix $k$ large enough, the Fatou's lemma implies 
\begin{align*}
	|u_{n_k}-u|^p_{\WnuOmR} \leq  \liminf_{\ell\to \infty}  \iil_{\Omega\R^d} \big|[u_{n_k}-u_{n_\ell}](x)-([u_{n_k}-u_{n_\ell}](y) \big|^p \, \nu (x-y) \d y\,\d x \,.
\end{align*}
Since $(u_{n_k})_k$ is a Cauchy sequence, the right hand side is finite for any $k$  and tends to $0$ as $k\to \infty$. This implies $u\in \WnuOmR$ and $|u_{n_k}-u|^p_{\WnuOmR} \xrightarrow{k\to \infty} 0$. Finally, $u_n\to u$ in $\WnuOmR$.  Furthermore, the map $\mathcal{I}: \WnuOmR\to L^p(\R^d) \times L^p(\Omega\times \R^d)$ with 
	\begin{align*}
	\mathcal{I}(u) = \Big(u(x), (u(x)-u(y))\nu^{1/p}(x-y)\Big)
	\end{align*}
	is an isometry. From its Banach structure, the space  $\big(\WnuOmR, \|\cdot\|_{\WnuOmR}\big)$, which can be identified with $\mathcal{I}\Big(\WnuOmR\Big)$, is  separable (and reflexive for $1<p<\infty$) as a closed subspace of the separable  (and reflexive for $1<p<\infty$) space  $ L^p(\R^d) \times L^p(\Omega\times \R^d)$.  Analogously, one can show that  $\big(\WnuOm, \|\cdot\|_{\WnuOm}\big)$ is a separable Banach space. 
	
	\medskip
\noindent	It remains to prove that $\big(\WnuOmR, \|\cdot\|_{\WnuOmR}\big)$ is a separable Banach space. Here we assume that $\nu$ has full support on $\R^d$. Without loss of generality we assume $\nu(h)>0$ for every $h\in \R^d$.   
	Let  $(u_n)_n$ be  a Cauchy sequence in $\big(\WnuOmR, \|\cdot\|_{\WnuOmR}\big)$. Then there exist a subsequence $(u_{n_k})_k$ of $(u_{n})_n$, a function $u$ in $L^p(\Omega)$, a function $U \in L^p(\Omega\times \mathbb{R}^d)$, and null sets $N\subset \mathbb{R}^d$  and $\mathcal{R}\subset \Omega\times \mathbb{R}^d$ such that 
	\begin{itemize}
		\item[-] $(u_{n_k})_k$ converges to $u$ in  $L^p(\Omega)$\,,
		\item[-] $(u_{n_k})_k$ converges to $u$ pointwise on $\Omega\setminus N$\,,
		\item[-] $(U_{n_k})_k$ converges to $U$ in  $L^p(\Omega \times \R^d)$\,,
		\item[-] $(U_{n_k})_k$ converges to $U$ pointwise on $(\Omega\times \R^d)\setminus \mathcal{R}$\,,
	\end{itemize}
	where $U_n(x,y) = (u_n(x)-u_n(y))\nu^{1/p}(x-y)$. For$(x,y)\in (\Omega \times \R^d)\setminus \mathcal{R'} $ with $x\neq y$ where $ \mathcal{R'}= \mathcal{R} \cup (N\times \emptyset)$,   
	\begin{align*}
	u_{n_k}(y)= u_{n_k}(x)  - U_{n_k}(x,y)/\nu^{1/p}(x-y)\xrightarrow{k\to \infty}  u(x)  - U(x,y)/\nu^{1/p}(x-y)
	\end{align*}
	
	\medskip
	
\noindent 	Finally, $U(x,y) = (u(x)-u(y))\nu^{1/p}(x-y)\in L^p(\Omega \times \R^d)$ so that $u \in \WnuOmR $. We easily conclude $\|u_n-u\|_{\WnuOmR}\xrightarrow{n\to\infty}  0$  which proves completeness. 

%
\end{proof}


\vspace{2mm}

\noindent Actually the spaces $\WnuOmO, \WnuOmR$ and $\WnuOmO$ (just like the Sobolev space $W^{1,p}(\Omega)$ and $W_0^{1,p}(\Omega) $) are refinement of $L^p(\Omega)$.  Next we highlight some strong connections with the classical Sobolev spaces $W^{1,p}(\Omega)$ and $W_0^{1,p}(\Omega) $) ($1\leq p<\infty$).  To start, let us point out some useful estimates. First of all observe that for $h\in \R^d$ we have $\|\tau_hu-u\|_{L^p(\mathbb{R}^d)}\leq 2\|u\|_{L^p(\mathbb{R}^d)}$. On the other hand since smooth functions of compact support are dense in $W^{1,p}(\R^d)$, using the fundamental theorem of calculus along with Jensen's inequality we find that 
\begin{align*}
\int_{\R^d} |u(x+h)-u(x)|^p\,\d x = \int_{\R^d} \Big|\int_0^1\nabla u(x+th)\cdot h\Big|^p\,\d x\leq |h|^p \|\nabla u\|_{L^{p}(\R^d)}. 
\end{align*}
Therefore, the following estimate holds true for all for $1\leq p<\infty$ 
\begin{align}\label{eq:levy-p-estimate}
\int_{\R^d} |u(x+h)-u(x)|^p\,\d x\leq 2^p(1\land |h|^p)\|u\|_{W^{1,p}(\R^d)},\quad\text{for every}\quad u\in W^{1,p}(\R^d),~\text{and}~h\in \R^d.
\end{align}
\noindent Next, since by Theorem \ref{thm:evans-grapiepy} the $BV$-norm of an element in $BV(\R^d)$ can be approximated by the $W^{1,1}$-norms of elements in $W^{1,1}(\R^d)$, we easily find the following analogous estimate
\begin{align}\label{eq:levy-estimate-Bv}
\int_{\R^d} |u(x+h)-u(x)|\,\d x\leq 2(1\land |h|)\|u\|_{BV(\R^d)},\quad\text{for every}\quad u\in BV(\R^d)~\text{and}~h\in \R^d.
\end{align}

\begin{lemma}\label{lem:boundedness-limsup}
	Let $\nu \in L^1(\R^d, 1\land |h|^p)$ be a nonnegative function with $1\leq p<\infty$. Let $\Omega$ be a $W^{1,p}$-extension (in particular a Lipschitz domain) open subset of $\mathbb{R}^d$. Then there is a constant $ C= C(\Omega, p,d)$ depending only $\Omega$ $p$ and $d$ such that for all $u\in W^{1, p}(\Omega)$
	\begin{align*}
	\iint\limits_{\Omega\Omega}|u(x)-u(y)|^p\nu(x-y) \mathrm{d}y\mathrm{d}x \leq C \|u\|^p_{W^{1,p}(\Omega) } \|\nu\|_{L^{1}(\R^d, 1\land |h|^p) }.
	\end{align*}
\end{lemma}

\medskip

\begin{proof}
	Since $\Omega$ is an extension domain, let $\overline{u}$ be an $ W^{1,p}(\Omega)$-extension of $u$ on $\mathbb{R}^d.$ First of all recall \eqref{eq:levy-p-estimate} that $\|\tau_h\overline{u}-\overline{u}\|^p_{L^p(\mathbb{R}^d)}\leq 2^p(1\land |h|^p) \|\overline{u}\|^p_{W^{1,p}(\mathbb{R}^d)}$ for all $h\in \mathbb{R}^d$. Hence using Fubini's theorem yields
	\begin{alignat*}{2}
	&\iint\limits_{\Omega\Omega}|u(x)-u(y)|^p\nu(x-y) \mathrm{d}y\mathrm{d}x 
	&&\leq \iint\limits_{\mathbb{R}^d\mathbb{R}^d}|\overline{u}(x+h)-\overline{u}(x)|^p \nu(h) \mathrm{d}h\mathrm{d}x\\
	&= \int\limits_{\mathbb{R}^d}\nu(h) \mathrm{d}h \int\limits_{ \mathbb{R}^d} |\overline{u}(x+h)-\overline{u}(x)|^p\mathrm{d}x
	&&\leq \|\overline{u}\|_{W^{1,p}(\mathbb{R}^d)} \int\limits_{\mathbb{R}^d}2^p(1\land |h|^p) \nu(h) \mathrm{d}h \\
	&\leq C\|u\|_{W^{1,p}(\Omega) } \|\nu\|_{L^{1}( \R^d, 1\land |h|^p) }. 
	\end{alignat*}
	
\end{proof}

\medskip

\noindent For the case $p=1$ the above result still holds true with $W^{1,1}(\Omega)$ replaced by the larger space $BV(\Omega)$.

\begin{lemma}\label{lem:boundedness-limsup-1}
	If $\nu\in L^1(\R^d,1\land |h|)$ and $\Omega$ is a $BV$-extension open subset of $\mathbb{R}^d$ there exists a constant $C=C(\Omega, d)>0$ such that for all $u\in BV(\Omega)$,
	\begin{align*}
	\iint\limits_{\Omega\Omega}|u(x)-u(y)|\nu(x-y) \mathrm{d}y\mathrm{d}x \leq C \|u\|_{BV(\Omega) } \|\nu\|_{L^{1}(\R^d, 1\land |h|) }.
	\end{align*}
\end{lemma}

\medskip

\begin{proof}
 Consider $\overline{u}\in BV(\mathbb{R}^d)$ to be the extension on $\mathbb{R}^d$ of a function$u\in BV(\Omega)$, i.e. $\|\overline{u}\|_{BV(\R^d)}\leq C \|u\|_{BV(\Omega)}$ with constant $C>0$ independent of $u$. We know from Theorem \ref{thm:evans-grapiepy} that there is $(u_n)_n$ a sequence of functions of $W^{1,1}(\R^d)$ converging to $\overline{u}$ in $L^1(\R^d)$ and such that $\|\nabla u_n\|_{L^1(\R^d)}\xrightarrow{n \to \infty}|\overline{u}|_{BV(\R^d)}$. This implies $\|u_n\|_{W^{1,1}(\R^d)}\xrightarrow{n \to \infty}\|\overline{u}\|_{BV(\R^d)}$. This together with Lemma \ref{lem:boundedness-limsup} and Fatou's lemma,  yield
	\begin{alignat*}{2}
	&\iint\limits_{\Omega\Omega}|u(x)-u(y)|\nu(x-y) \mathrm{d}y\mathrm{d}x
	\leq \iint\limits_{\R^d\R^d}|\overline{u}(x)-\overline{u}(y)|\nu(x-y) \mathrm{d}y\mathrm{d}x \\
	&\leq \liminf_{n\to\infty} \iint\limits_{\R^d\R^d}|u_n(x)-u_n(y)|\nu(x-y) \mathrm{d}y\mathrm{d}x
	 \leq 2 \liminf_{n\to\infty} \|u_n\|_{W^{1,1}(\R^d) } \|\nu\|_{L^1(\R^d,1\land |h|) }\\
	&= 2 \|\overline{u}\|_{BV(\R^d) } \|\nu\|_{L^1(\R^d,1\land |h|)}
	\leq C\|u\|_{BV(\Omega)} \|\nu\|_{L^{1}(\R^d, 1\land |h|)}.
	\end{alignat*}
	\noindent Willingly, one could proceed as in the Lemma \ref{lem:boundedness-limsup} by means of the inequality \eqref{eq:levy-estimate-Bv}. 
	
\end{proof}

\noindent It is noteworthy to emphasize that Lemma \ref{lem:boundedness-limsup} (respectively Lemma \ref{lem:boundedness-limsup-1} ) does not hold true if $\Omega$ is not an extension domain, see the Counterexample \ref{Ex:counterexample-extension}. As a result we have the following.
\medskip

\begin{theorem}\label{thm:w1p-in nonlocal}
	If $\Omega$ is an $W^{1,p}$-extension (resp. $BV$-extension) domain, then  the following embedding is continuous
\begin{align*}
	W^{1,p}(\Omega)\hookrightarrow \WnuOm\,\quad ( \text{resp. for $p=1$} \quad
	BV(\Omega)\hookrightarrow W^1_\nu(\Omega)\,). 
	\end{align*}
	The  above embeddings fail if $\Omega$ is not an extension domain (see the counterexample \ref{Ex:counterexample-extension}).
\end{theorem}

\vspace{2mm}
\begin{proof}
	Let $\overline{u}\in W^{1,p}(\mathbb{R}^d)$ be an extension of a function $u\in W^{1,p}(\Omega) $ with			 $\|\overline{u}\|_{W^{1,p}(\mathbb{R}^d)} \leq C \|u\|_{W^{1,p}(\Omega)}$ for a constant $C$ depending only on $\Omega$ and $d$. Within the estimate \eqref{eq:levy-p-estimate} we easily get the continuous embedding $W^{1,p}(\Omega)\hookrightarrow \WnuOm\,$. 
	Analogously  using Lemma 
	\ref{lem:boundedness-limsup-1} one can also establish that $BV(\Omega)\hookrightarrow W^1_\nu(\Omega)\, $. If $\Omega$ is not an extension domain  the counterexample \ref{Ex:counterexample-extension} shows that the results do not always hold.
\end{proof}

\bigskip

\noindent Let us collect some trivial embeddings also involving the classical Sobolev spaces. Note the norm $\vertiii{\cdot }_{\WnuOmR}$ is not excluded in the following. 

\medskip 

\noindent $\bullet$ It springs from \eqref{eq:levy-p-estimate} that the following continuous embeddings hold true: 
\begin{align*}
W^{1,p}(\mathbb{R}^d) \hookrightarrow W^p_{\nu}(\mathbb{R}^d) \hookrightarrow \WnuOmR. 
\end{align*}

\noindent $\bullet$ We obviously have the following continuous embeddings 
\begin{align*}
\WnuOmR\hookrightarrow \WnuOm \hookrightarrow L^p(\Omega)\, .
\end{align*} 

\noindent $\bullet$ Let $ W^p_{\nu,0}(\Omega) $ be the closure of $C_c^\infty(\Omega)$ with respect to $\|\cdot \|_{\WnuOm)}$. Note that- the zero extension to $\mathbb{R}^d$ of any function in $ W_0^{1,p}(\Omega)$ belongs to $W^{1,p}(\mathbb{R}^d)$. Hence, using \eqref{eq:levy-p-estimate}  we also have the following continuous embeddings: 
\begin{align*}
&W_0^{1,p}(\Omega)\hookrightarrow \WnuOmO \hookrightarrow W^p_{\nu,0}(\Omega)\hookrightarrow L^p(\Omega)
\quad\text{and}\quad
\WnuOmO \hookrightarrow W_{\nu}(\mathbb{R}^d).
\end{align*} 

\noindent The embeddings $  \WnuOmO \hookrightarrow W^p_{\nu,0}(\Omega)$ and $\WnuOmO \hookrightarrow W_{\nu}(\mathbb{R}^d)$ follow from the fact  that for all $u\in \WnuOmR$
\begin{align*}
\iint\limits_{\mathbb{R}^d\mathbb{R}^d} |u(x)-u(y)|^p\nu(x-y)\d y\,\d x
&= \iint\limits_{\Omega \Omega} |u(x)-u(y)|^p\nu(x-y)\d y\,\d x+ 2\iint\limits_{\Omega \Omega^c} |u(x)|^p\nu(x-y)\d y\,\d x.
\end{align*} 
It is worth noticing  that not every function $u \in H_{\nu,0}(\Omega)$ has its zero extension in $H_{\nu}(\mathbb{R}^d)$. Indeed for this to hold, one would need that
\begin{align*}
\int_{\Omega}|u(x)|^p\d x\int_{\Omega^c} \nu(x-y)\d y<\infty. 
\end{align*}
This is not always true because the measure $\nu$ might be very singular at the origin. This purely nonlocal effect contrasts with the local function space $W_0^{1,p}(\Omega)$ whose elements can be isometrically extended by zero on $\mathbb{R}^d$ as functions of $W^{1,p}(\mathbb{R}^d)$.  This also shows that for some appropriate domain $\Omega$ and for some appropriate measure $\nu$ e.g. $\nu(h) = |h|^{-d-\alpha}$ the spaces $W_{\nu,0}(\Omega)$ and $\WnuOmO$ are strictly different although they both possess $C_c^\infty(\Omega)$ as dense subspace. 

\bigskip

\begin{proposition}\label{prop:nonlocal-embeding} Let $\Omega\subset\R^d $ be open.
\begin{enumerate}[$(i)$]
	\item Let $\nu_1$ and $\nu_2$ satisfying \eqref{eq:plevy-integrability}. Assume that there exist two constants $r>0$ and $k>0$ such that if  $|h|\leq r$  we have
	$$\nu_1(h)\leq k\nu_2(h).$$
	Then the   following embeddings are continuous 
	\begin{align*}
	W^p_{\nu_2} (\Omega)\hookrightarrow  W^p_{\nu_1} (\Omega)\quad \text{and}\quad
	W^p_{\nu_2} (\Omega|\R^d)\hookrightarrow  W^p_{\nu_1} (\Omega|\R^d). 
	\end{align*}
	
	\item  If $\Omega_1\subset \Omega_2$ then we have the continuous embeddings
	\begin{align*}
	W^p_{\nu} (\Omega_2)\hookrightarrow  W^p_{\nu} (\Omega_1)\quad \text{and}\quad
	W^p_{\nu} (\Omega_2|\R^d)\hookrightarrow  W^p_{\nu} (\Omega_1|\R^d). 
	\end{align*}
	\item  For a ball $B\subset \Omega$ we have the continuous embeddings
	\begin{align*}
\WnuOmR \hookrightarrow  L^p(\R^d ,\nu_{B})\quad \text{and}\quad
\WnuOmR \hookrightarrow  L^p(\R^d, \mathring{\nu}_{B}),\
	\end{align*}
	 where for $x\in \R^d$ we let
	 \begin{align*}
	 \nu_{B}(x)= \operatorname{essinf}_{y\in B}\nu(x-y)\quad \text{and}\quad \mathring{\nu}_{B}(x)= \int_B 1\land \nu(x-y)\d y.
	 \end{align*}
\end{enumerate}
\end{proposition}

\medskip

\begin{proof}
Note that $(ii)$ is obvious and for the proof of $(iii)$ we use analogous methods as in Section \ref{sec:function-spaces}. Putting $C_\delta = \int_{|h|\geq \delta} \nu_1(h)\d h$,  then for $u\in W^p_{\nu_2} (\Omega)$ we get
\begin{align*}
\iil_{ \Omega\Omega} |u(x) -u(y)|^p\nu_1(x-y) \d y\d x
&\leq 2^pC_\delta\|u\|^p_{L^p(\Omega)}+ \iil_{ \Omega\Omega \cap\{|x-y|\leq \delta\}} |u(x) -u(y)|^p\nu_1(x-y) \d y\d x\\
&\leq 2^pC_\delta\|u\|^p_{L^p(\Omega)}+ k\iil_{ \Omega\Omega \cap\{|x-y||\leq \delta\}} |u(x) -u(y)|^p\nu_2(x-y) \d y\d x\\
&\leq2^p(k+1)C_\delta\|u\|^p_{W^p_{\nu_2} (\Omega)}
\end{align*}
From this it follows that  $W^p_{\nu_2} (\Omega)\hookrightarrow  W^p_{\nu_1} (\Omega)$ and likewise 
$W^p_{\nu_2} (\Omega|\R^d)\hookrightarrow  W^p_{\nu_1} (\Omega|\R^d). $

\end{proof}

\noindent  In light of Proposition \ref{prop:nonlocal-embeding} it is legitimate to formulate the following definition. 
\begin{definition}
	Given two functions $\nu_1$ and $\nu_2$ satisfying \eqref{eq:plevy-integrability} we shall say that $\nu_2$ is regular than $\nu_1$ there exist two constants $r_0, \kappa>0$ such that for all $h\in\R^d$ with $|h|\leq r_0$  one has $\nu_1(h) \leq \kappa \nu_2(h)$. 

\medskip

\noindent This definition is not fortuitous. In some sense it means that $\nu_2$  is likely to be more singular than  $\nu_1$ near  the vicinity of the original. Should the singularity of $\nu_2$ be higher, more regularity of a function $u$ is required for the integrability  of the increment $(x,y)\mapsto (u(x)-u(y))\nu_2^{1/p}(x-y)$.  A simple example is given as follows.
\end{definition}

\begin{example}
For $0<s<s'<1$ we have $W^{s',p}(\Omega) \hookrightarrow W^{s,p}(\Omega)$ and  $W^{s',p}(\Omega|\R^d) \hookrightarrow W^{s,p}(\Omega|\R^d)$. Indeed, this follows from $(i)$ since for every $|h|\leq 1$ we have $|h|^{-d-sp}\leq |h|^{-d-s'p}$ and hence it suffices to take $\nu_1(h) =|h|^{-d-sp}$, $\nu_2(h) =|h|^{-d-s'p}$.
Moreover we have the continuous embedding $ W^{s,p}(\Omega|\R^d) \hookrightarrow L^p(\R^d, (1+|h|)^{-d-sp} \d h).$ Note that the norm of $W^{s,p}(\Omega)$ can be rescaled so that in the limiting case $s=1$, $W^{s,p}(\Omega)$ is precisely $W^{1,p}(\Omega)$  whose functions are more regular. 
\end{example}


\begin{counterexample}\label{Ex:counterexample-extension}
	For a simple instance, in one dimension, consider $\Omega= (-1,0)\cup (0,1)$ and put $u(x) = -1$ if $x\in (-1,0)$ and $u(x) = 1$ if $x\in (0,1)$. Clearly, we have $u \in W^{1,p}(\Omega)$ for all $1\leq p<\infty$ with $\nabla u=0$. However, $u$ does not belong to any of the fractional Sobolev space $W^{s,p}(\Omega)$ provided that $s\geq 1/p$. Recall that here we have $\nu(h) =|h|^{-1-sp}$. Indeed, since $s\geq 1/p$
	\begin{align*}
	\iil_{\Omega\Omega} \frac{|u(x)- u(y)|^p}{|x-y|^{1+sp}}\,\d x\,\d y= 2^{p+1}\int_{0}^1\int_{-1}^0 \frac{\,\d x\,\d y}{|x-y|^{1+sp}} = 2^{p+1}\int_{0}^1 x^{-sp}+ (1+x)^{-sp} \,\d x= \infty.
	\end{align*}
	Moreover, $\Omega= (-1,0)\cup (0,1)$ is not a $W^{1,p}$-extension domain. Indeed, if $\overline{u}\in W^{1,p}(\R^d)$ is an extension of $u$ defined as above then in particular we would have $\overline{u}\in W^{1,p}(-1,1)$ and $\overline{u}= u$ on $\Omega= (-1,0)\cup (0,1)$. The distributional derivative of $\overline{u}$ on $(-1, 1)$ is $\nabla \overline{u} = 2\delta_0 $ (where $\delta_0$ stands for the Dirac mass at the origin). Hence $\nabla \overline{u}$ is not a function which contradicts the fact that $\overline{u}\in W^{1,p}(\R)$. 
	
	\vspace{2mm}
	\noindent Next we prove that this example persists in the higher dimensional space $d\geq 2$.
	Let $ B^+_1(0)= B_1(0)\cap \{(x',x_d)\in \R^d:~x_d>0\}$ and $ B^-_1(0)= B_1(0)\cap \{(x',x_d)\in \R^d:~x_d<0\}$. Put $\Omega= B_1(0)\setminus \{(x',x_d)\in \R^d:~x_d=0\}= B^+_1(0)\cup B^-_1(0) $ and 
	define the function 
	$u(x) = \mathds{1}_{B^+_1(0)}(x)- \mathds{1}_{B^-_1(0)}(x)$, i.e. $u(x) = 1$ if $x\in B^+_1(0)$ and $u(x) = -1$ if $x\in B^-_1(0)$. Obviously we have $u \in W^{1,p}(\Omega)$ for all $1\leq p<\infty$ with $\nabla u=0$. Further following the discussion above, one can check that $u$ does not have any extension to the whole space. Hence $\Omega$ cannot be an extension domain. On the other hand $u$ does not belong to any of the fractional Sobolev space $W^{s,p}(\Omega)$ provided that $s\geq 1/p$. Recall that here we have $\nu(h) =|h|^{-d-sp}$. Indeed, assume $s\geq 1/p$. Since integrals disregard sets of zero measure, we have
	\begin{align*}
	\iil_{\Omega\Omega} \frac{|u(x)- u(y)|^p}{|x-y|^{d+sp}} \,\d x\,\d y= 2^{p+1}\int_{B^+_1(0)}\int_{B^-_1(0)} \frac{\,\d x\,\d y}{|x-y|^{d+sp}} = 2^{p+1}\int_{B^+_1(0)}\d x \hspace{-3ex}\il_{B_1(x)\cap\{h_d>x_d\}} \hspace{-3ex}|h|^{-d-sp}\,\d h = \infty.
	\end{align*}
	We justify this conclusion as follows. It is easy to show that $D'\subset D$ where 
	\begin{align*}
	D&= \{(x,h)\in \R^d\times \R^d: x\in B^+_1(0),~~h\in B^+_1(x) \cap\{h_d>x_d\}\}\\ 
	D'&=B^+_{1/2}(0)\times \{(h',h_d)\in \R^d:~|h'|<1/4,~~x_d<h_d<1/4\}.
	\end{align*}
\end{counterexample}
Also note that $ \{(x',x_d)\in \R^d:~|x'|<1/4,~~0<x_d<1/4\} \subset B^+_{1/2}(0) $. Using the change of variables $h'=h_dz'$ so that $\d h'=h_d^{d-1}\d z'$ then since $s\geq 1/p$ we get the following 
\begin{align*}
\iil_{\Omega\Omega} \frac{|u(x)- u(y)|^p}{|x-y|^{d+sp}} \,\d x\,\d y
&= 2^{p+1}\iint_{D} |h|^{-d-sp}\,\d x\,\d h\geq 2^{p+1}\iint_{D'} |h|^{-d-sp}\,\d x\,\d h\\
&= 2^{p+1}\il_{B^+_{1/2}(0)}\d x \int_{x_d}^{1/4}\d h_d\il_{|h'|<1/4} \hspace{-2ex}(|h'|^2+h_d^2)^{-(d+sp)/2}\,\d h'\\
&= 2^{p+1}\il_{B^+_{1/2}(0)}\d x \int_{x_d}^{1/4}h_d^{-1-sp}\d h_d
\il_{|h'|<1/4} \hspace{-2ex}(1+|z'|^2)^{-(d+sp)/2}\,\d z'\\
&=C_s\il_{B^+_{1/2}(0)}(x^{-sp}_d-4^{-sp})\d x \qquad\text{with}\quad C_s= 2^{p+1}\il_{|z'|<1/4} \hspace{-2ex}(1+|z'|^2)^{-(d+sp)/2}\,\d z'\\
&\geq C_s \int_{0}^{1/4} (x^{-sp}_d-4^{-sp}) \d x_d\il_{|x'|<1/4}\d x' =\infty.
\end{align*}

\section{Fractional Sobolev spaces}
Let us first observe that for $s\in \R$,  $\nu: h\mapsto |h|^{-d-sp}$ satisfies \eqref{eq:plevy-integrability}, i.e. belongs to $L^p(\R^d, (1\land |h|^p))$ if and only if $s\in (0,1)$. In this case the spaces $\WnuOm, W^p_{\nu,0}(\Omega), \WnuOmO$ and $\WnuOmR$ are respectively denoted by $W^{s,p}(\Omega)$, $W^{s,p}_0(\Omega)$, $W^{s,p}_*(\Omega)$  and $W^{s,p}(\Omega|\R^d)$.
 The spaces $W^{s,p}(\Omega)$ and $W^{s,p}_0(\Omega)$ are recognized as the usual fractional Sobolev spaces.  They are often called \textit{Aronszajn, Gagliardo or Slobodeckij} spaces\cite{Aro55,Gag61,Slo58} after the names of those who introduced them simultaneously. Those spaces  are of particular interest in their own right and have been extensively studied throughout the literature in recent years.  Rigorous treatises on fractional Sobolev spaces can be found in \cite{grisvard11, maz2013sobolev}. 
See \cite{Hitchhiker} for a treatment of Sobolev spaces of fractional order. The terminologies \emph{fractional} is on one hand owed to the fact that for on a  smooth domain the space $W^{s,p}(\Omega)$ can be realized as the $s$-interpolation of $L^p(\Omega)$ and $W^{1,p}(\Omega)$. 
On the other hand the spaces $W^{1-\frac1p,p}(\partial \Omega)$ (with $\partial \Omega$ viewed as an open manifold of $\R^{d-1}$) is the trace space of $W^{1,p}(\Omega)$. Precisely  there exists a bounded operator $T: W^{1,p}(\Omega) \to W^{1-\frac1p,p}(\partial\Omega)$  such that if $u\in C^\infty(\overline{\Omega})\cap W^{1,p}(\Omega)$ we have  $Tu=u\mid_{\partial\Omega}$, i.e. $Tu=u$ on $\partial \Omega$ ). 
The existence of such an operator was proven by Emilio Gagliardo around 1957, a modern treatment of this result can be found in \cite{Ponce16elliptic}. Next we relate the fractional Sobolev spaces to the interpolation of 
$L^p$ and $W^{1,p}$. 
\begin{theorem}\label{thm:log-convex-sobolev} Let $s\in (0,1)$, $1\leq p<\infty$. There is a  constant $\theta(s,d, p)>0$ such that for all $u\in W^{s,p}(\R^d)$ the following log-convex inequality (interpolation estimate) holds
\begin{align}\label{eq:interpolation-inequality}
\|u\|_{W^{s,p}(\R^d)} \leq \theta(s,d, p) \|u\|_{L^p(\R^d )}^{1-s} 
		\|u\|_{W^{1,p}(\R^d)}^s.
\end{align}
With $ 1\leq \theta(s,d, p)< C(s(1-s))^{-1/p}$ where the constant $C$ does not depend on $s$.

\noindent If $\Omega$ is a $W^{1,p}$-extension domain there is $C>0$ also depending on $\Omega$ such that
\begin{align*}
\|u\|_{W^{s,p}(\Omega)} \leq C(d, p, \Omega) (s(1-s))^{-1/p}\|u\|_{L^p(\Omega )}^{1-s} 
\|u\|_{W^{1,p}(\Omega)}^s\quad\text{for all $u\in W^{s,p}(\Omega) $}. 
\end{align*}
\end{theorem}

\medskip

\begin{proof} Assume $u \in W^{1,p}(\R^d)$ and let $r>0$. We have 
	\begin{align*}
	\il_{\R^d}\il_{|h|\geq r} \frac{|u(x) -u(x+h)|^p}{|h|^{d+sp}} \d h\d x\leq 2^{p} \il_{\R^d}|u(x)|^p\d x \il_{|h|\geq r}\hspace{-2ex} |h|^{-d-sp}\d h =2^p |\mathbb{S}^{d-1}|\frac{r^{-sp}}{sp} \|u\|^p_{L^p(\R^d )}.
	\end{align*}
	Furthermore, using Fubini's theorem we get 
	\begin{align*}
	\il_{\R^d}\il_{|h|< r} \frac{|u(x) -u(x+h)|^p}{|h|^{d+sp}} \d h\d x
	&= \il_{|h|< r}  \il_{\R^d} \int_0^1\big|\nabla u(x+th)\cdot h \d t\big|^p|h|^{-d-sp} \d x \d h\\
	&\leq \il_{\R^d} \big|\nabla u(x)\big|^p\d x \il_{|h|< r} |h|^{-d-(1-s)p} \d h= |\mathbb{S}^{d-1}|\frac{r^{p(1-s)}p}{p(1-s)} \|\nabla u\|^p_{L^p(\R^d )}.
	\end{align*}
Adding the two inequalities for\footnote{This value is obtained after solving  for $r,$ $2^p |\mathbb{S}^{d-1}|\frac{r^{-sp}}{sp} \|u\|^p_{L^p(\R^d )}= |\mathbb{S}^{d-1}|\frac{r^{p(1-s)}p}{p(1-s)} \|\nabla u\|^p_{L^p(\R^d )}.$}  
$r=\Big(\tfrac{1-s}{s}\Big)^{1/p} \tfrac{\|\nabla u\|_{L^p(\R^d )} }{\|u\|_{L^p(\R^d )} }$ and letting  $\theta^p(s,d, p) = 2^{sp}\big(1+\tfrac{|\mathbb{S}^{d-1}|}{ps^{1-s}(1-s)^s}\big)$ yields\linebreak[-2]
\begin{align*}
\|u\|^p_{L^p(\R^d )}+ |u|^p_{W^{s,p}(\R^d )}
&\leq  \tfrac{|\mathbb{S}^{d-1}|}{ps^{1-s}(1-s)^s} \|u\|^{p(1-s)}_{L^p(\R^d )} \|\nabla u\|^{sp}_{L^p(\R^d )}\\
&2^{-sp}\theta^p(s,d, p)  \|u\|^{p(1-s)}_{L^p(\R^d )}\big(\|u\|^{sp}_{L^p(\R^d )}+ \|\nabla u\|^{sp}_{L^p(\R^d )}\big)\\
&\leq \Big[ \theta(s,d, p)\|u\|^{1-s}_{L^p(\R^d )}\|u\|^{s}_{W{1,p}(\R^d )} \Big]^p\, . 
\end{align*}
Note that $e^{-1/e} \leq t^t \leq 1$ for $t\in (0,1)$ hence $\theta(s,d, p) \leq C (s(1-s))^{-1/p}$. The claim is proved.
\end{proof}

\noindent We observe that \eqref{eq:interpolation-inequality} is sharp in the sense that tending  $s\to1$(resp. $s\to0^+$) leads to the obvious inequality $ \|u\|^p_{W^{1,p}(\R^d )}\leq C \|u\|^p_{W^{1,p}(\R^d )}$ ( resp. $\|u\|^p_{L^p(\R^d )}\leq C\|u\|^p_{L^p(\R^d )}$).  We discuss this widely in Section \ref{sec:charact-W1p}. 
In addition, the log-convex inequality \eqref{eq:interpolation-inequality} somehow encodes an interpolation between $L^p$ and $W^{1,p}$ whose the resulting interpolation space is $W^{s,p}$.   Truly speaking, the fractional Sobolev space results from the so called \textit{real interpolation.}  A concrete elaboration of this assertion can be found in \cite{triebel1995interpolation, bergh2012interpolation}. 

\vspace{2mm}

\noindent Next we derive the fractional counterpart of the famous Gagliardo–Nirenberg–Sobolev embedding theorem. The spectacularly amazing elementary proof presented here is apparently due to Haim Brezis \cite[Proposition 15.5]{Ponce16elliptic}  from a personal communication. It is important to highlight that earlier proofs of this theorem exist in the literature as well. For instance another lengthy proof also using basic analysis tools is well incorporated in \cite[Theorem 6.5]{Hitchhiker}.  See also \cite{CT04} where the inequality is established with the best constant for $p=2$.

\begin{theorem}[\textbf{Gagliardo–Nirenberg–Sobolev}]
	Let $s \in (0, 1)$ and $1\leq p<\frac{d}s$. For all $u\in L^{p^*_s}(\R^d)$,
\begin{align}
	\|u\|_{L^{p^*_s}(\R^d)}\leq 2^{p^*_s/p} |B_1(0)|^{-1/p-s/d}\Big(\iil_{\R^d\R^d} \frac{|u(x)- u(y)|^p}{|x-y|^{-d-sp}}\d x\d y\Big)^{1/p}. 
\end{align}
Here $\dfrac{1}{p^*_s} = \dfrac1p- \dfrac {s}d>0$ is the so-called \textit{fractional critical exponent or fractional Sobolev conjugate}. In particular, the embedding $ W^{s,p}(\R^d) \hookrightarrow L^{p^*_s}(\R^d)$ is continuous. 
\end{theorem}

\vspace*{1mm}

\begin{proof}First we fix $ x\in \R^d$ and $r>0$. 
Integrating the inequality $|u(x)|\leq |u(y)|+ |u(x)-u(y)|$ over $y\in B_r(x)$  and using Jensen's inequality gives 

\begin{align*}
|u(x)|&\leq \fint_{B_r(x)}|u(y)|\d y+ \fint_{B_r(x)}|u(x)-u(y)|\d y\\
&\leq \Big(\fint_{B_r(x)}|u(y)|^{p^*_s}\d y\Big)^{1/p^*_s}+ \Big(\fint_{B_r(x)}|u(x)-u(y)|^p\d y\Big)^{1/p}\\
&\leq  \Big(\fint_{B_r(x)}|u(y)|^{p^*_s}\d y\Big)^{1/p^*_s}+ \Big(r^{d+sp}\fint_{B_r(x)}\frac{|u(x)-u(y)|^p}{|x-y|^{d+sp}}\d y\Big)^{1/p}\\
&= r^{-d/p^*_s}|B_1(0)|^{-1/p^*_s}\Big(\il_{B_r(x)}|u(y)|^{p^*_s}\d y\Big)^{1/p^*_s}+ r^s|B_1(0)|^{-1/p}\Big(\il_{B_r(x)}\frac{|u(x)-u(y)|^p}{|x-y|^{d+sp}}\d y\Big)^{1/p}\\
&\leq r^{-d/p^*_s}|B_1(0)|^{-1/p^*_s}\Big(\il_{\R^d}|u(y)|^{p^*_s}\d y\Big)^{1/p^*_s}+ r^s|B_1(0)|^{-1/p}\Big(\il_{\R^d}\frac{|u(x)-u(y)|^p}{|x-y|^{d+sp}}\d y\Big)^{1/p}.
\end{align*}

\noindent Next solving the following equation for  $r$,
 $$ r^{-d/p^*_s}|B_1(0)|^{-1/p^*_s}\Big(\il_{\R^d}|u(y)|^{p^*_s}\d y\Big)^{1/p^*_s}= r^s|B_1(0)|^{-1/p}\Big(\il_{\R^d}\frac{|u(x)-u(y)|^p}{|x-y|^{d+sp}}\d y\Big)^{1/p}$$
gives 
\begin{align*}
r(x)= r=|B_1(0)|^{1/d-p/dp^*_s}\Big(\il_{\R^d}|u(y)|^{p^*_s}\d y\Big)^{p/dp^*_s}\Big(\il_{\R^d}\frac{|u(x)-u(y)|^p}{|x-y|^{d+sp}}\d y\Big)^{-1/d}.
\end{align*} 
Substituting this specific $r(x)$ in the preceding estimate leads to 
\begin{align*}
|u(x)|^{p^*_s}
&\leq 2^{p^*_s} r^{-d}(x)|B_1(0)|^{-1}\Big(\il_{\R^d}|u(y)|^{p^*_s}\d y\Big)\\
&= 2^{p^*_s} |B_1(0)|^{-2+p/p^*_s}\Big(\il_{\R^d}|u(y)|^{p^*_s}\d y\Big)^{1-p/p^*_s}\Big(\il_{\R^d} \frac{|u(x)-u(y)|^p}{|x-y|^{d+sp}}\d y\Big). 
\end{align*}
Integrating both sides with respect to the variable $x$ gives 

\begin{align*}
\il_{\R^d}|u(x)|^{p^*_s}\d x
&\leq  2^{p^*_s} |B_1(0)|^{-2+p/p^*_s}\Big(\il_{\R^d}|u(y)|^{p^*_s}\d y\Big)^{1-p/p^*_s}\Big(\iil_{\R^d\R^d} \frac{|u(x)-u(y)|^p}{|x-y|^{d+sp}}\d y\d x\Big). 
\end{align*}
Cancellation on both sides provides the desired estimate. 

\end{proof}

\noindent It is interesting to know what could be the suitable definition of fractional Sobolev space $W^{s,p}$ for $s\geq 1$. It turns out that the definition of $|\cdot|_{W^{s,p}(\Omega)}$ cannot be extended to the situation where $s\geq 1$. This is justified by the following result. 
\begin{proposition}\label{prop:fracasymp}
	Let $\Omega \subset \R^d$ be open and let $u \in C_c^\infty(\Omega)$ with $u\neq 0$. Then for $1\leq p<\infty$
	\begin{align*}
	|u|_{W^{s,p}(\Omega)}= \infty\quad\text{for all $s\geq 1$ and }\quad \lim_{s\to 1^-}|u|_{W^{s,p}(\Omega)}= \infty.
	\end{align*}
\end{proposition} 

\begin{proof} 
	Let $ K \subset\Omega $ be a compact and let $\delta>0$ such that $K(\delta)= K+B_\delta(0)\subset \Omega$. Since $u\in C^2_c(\Omega)$, by mimicking \ref{eq:optimal-non-levy} we find that 
	\begin{align*}
	|u|_{W^{s,p}(\Omega)}^p
	&\geq \int_{K(\delta)} \mathrm{d}x \int_{\mathbb{S}^{d-1}} \left|\nabla u(x)\cdot w\right|^p \mathrm{d}\sigma(w) \int_{0}^{\delta} r^{(1-s)p-1} \mathrm{d}r.
	\end{align*}
	The claim follows since for all $s\geq 1$, $\int_{0}^{\delta} r^{(1-s)p-1} \mathrm{d}r= \infty $ and we also have $$\int_{0}^{\delta} r^{(1-s)p-1} = \frac{\delta^{(1-s)p}}{ (1-s)p}\xrightarrow{s\to 1^-}\infty.$$
\end{proof}

\noindent On a smooth domain $\Omega$ the spaces $W^{s,p}(\Omega)$ with $0<s<1$ can be realized as the interpolation space between $L^{p}(\Omega)$ and $W^{1,p}(\Omega)$. Therefore it is interesting to know how close is the space $\|\cdot\|_{W^{s,p}(\Omega)}$ to $\|\cdot\|_{L^p(\Omega)}$ as $s\to 0$ and to $\|\cdot\|_{W^{1,p}(\Omega)}$ as $ s\to 1^-$. Unpleasantly, for the asymptotic $\to 1^-$ we know from Proposition \ref{prop:fracasymp} that $\lim\limits_{s\to 1^-}|u|_{W^{s,p}(\Omega)}= \infty$. This divergence was initially observed by Haim Brezis, Jean Bourgain and Pierre Mironescu in 2001 \cite{BBM01}. As one can foresee from the above computation, to correct this anomaly, putting the factor $(1-s)$ in front of the term $|u|^p_{W^{s,p}(\Omega)}$ annihilates the singularity. In fact, they prove that for all $u\in W^{1,p}(\Omega)$
\begin{align}\label{eq:frac-BMM}
\lim_{s\to 1^-}(1-s)|u|^p_{W^{s,p}(\Omega)} = K_{d,p}\|\nabla u\|_{L^{p}(\Omega)}^p:= 
K_{d,p}|u|^p_{W^{1,p}(\Omega)}.
\end{align}
for some appropriate universal constant $ K_{d,p}$ depending only on $d$ and $p\geq 1$. We shall see this in a more general context later on. The case $s\to 0$ was solved in \cite{MS02} claiming that for all $u\in \bigcup\limits_{0<s<1}W^{s,p}(\R^d)$, 
\begin{align*}
\lim_{s\to 0^+} s|u|^p_{W^{s,p}(R^d)} =\frac{2}{p}|\mathbb{S}^{d-1}|\|u\|_{L^{p}(\R^d)}^p.
\end{align*}

\noindent Seemingly, Proposition \ref{prop:fracasymp} suggests that the space $W^{s,p}(\Omega)$ for $s> 1$ must be defined in a different way. Assuming for now that the relation \eqref{eq:frac-BMM} is true the following definition plainly makes sense. 
\medskip

\begin{definition} Let $\Omega\subset \R^d$ be an open set, $1\leq p<\infty$ and $s>0$. We put $s=m+\sigma$ if $s$ is not an integer, where $m\in \mathbb{N}$ and $\sigma\in (0,1)$. The space $W^{s,p}(\Omega)$ assumes the equivalence classes of functions $u\in W^{m,p}(\Omega)$ whose distributional derivatives $D^{\alpha}u$ of order $|\alpha|=m$ belong to $W^{\sigma,p}(\Omega)$, namely
	\begin{align*}
	W^{s,p}(\Omega):= \Big\{ u\in W^{m,p}(\Omega)~:\; D^{\alpha}u \in W^{\sigma,p}(\Omega) \ \text{for any}\ \alpha\ \text{s.t.} \ |\alpha|=m \Big\}
	\end{align*}
	The space $W^{s,p}(\Omega)$ is a Banach space with respect to the norm defined by 
	\begin{align*}
	\|u\|_{W^{s,p}(\Omega)}:= \Big( \|u\|^p_{W^{m,p}(\Omega)} + \sum_{|\alpha|=m} \|D^{\alpha}u\|^p_{W^{\sigma,p}(\Omega)} \Big)^{\!\frac{1}{p}}.
	\end{align*}
	\vspace{-2ex}
	
		\noindent Purposely we emphasize that 
	$ W^{\sigma,p}(\Omega) = \Big\{ u\in L^p(\Omega) : \dfrac{|u(x)-u(y)|}{|x-y|^{\frac{d}{p} + \sigma}} \in L^p(\Omega\times\Omega) \Big\}$ and 
	$$ \|u\|^p_{W^{\sigma,p}(\Omega)} =  \int_\Omega |u(x)|^p \; dx + \iil_{\Omega\Omega} \frac{|u(x)-u(y)|^p}{|x-y|^{d+ \sigma p}} \; dx dy. $$
	\vspace{-2ex}
	
\noindent Furthermore, $W^{-s,p}(\Omega)$ is the topological dual of $W_0^{s,p'}(\Omega) = \overline{C_c(\Omega)}^{ W^{s,p'}(\Omega)}$ where $1/p'=1-1/p$.  The space $W^{s,\infty}(\Omega)$ boils down to the usual H\"{o}lder space $C^{m, \sigma}(\Omega)$. 

\noindent Clearly, if $s=m$ is an integer, the space $W^{s,p}(\Omega)$ coincides with the Sobolev space $W^{m,p}(\Omega)$.  
\end{definition}

  \begin{remark}
  The fractional Sobolev spaces $W^{s,p}(\Omega)$ fill the gaps between the classical Sobolev of integers orders. Actually there are two main different schools when it comes to define the notions of fractional Sobolev spaces. 
  The spaces  $W^{s,p}(\Omega)$, $s\in \R$ belong to the schools of Besov-Gagliardo-Nirenberg. There are also the spaces $H^{s,p}(\Omega)$ (which we define below) belonging to the school of Triebel-Linzorkin. The latter spaces can be realized as the \textit{complex interpolation} between $L^p(\R^d)$ and $W^{1,p}(\R^d)$ whereas $W^{s,p}(\Omega)$  is the \textit{real interpolation} between $L^p(\R^d)$ and $W^{1,p}(\R^d)$. 
  \end{remark}

\begin{definition}[Via Fourier transform]
	Let $\Omega\subset \R^d$ be open, $1\leq p<\infty$ and $s\in\mathbb{R}$. According to the definition due to Triebel-Linzorking the fractional Sobolev space denoted $H^{s,p}(\mathbb{R}^d )$ is the completion of the space
	$$H^{s,p}(\mathbb{R}^d):=\Big\{u\in\mathcal{S}(\mathbb{R}^{n}) : \|(\langle{\xi}\rangle^{s}\widehat{u})^{\vee}\|_{L^{p}(\R^d)}<\infty\Big\},$$
	equipped with the norm
	$$\|u\|_{H^{s,p}(\R^d)}=\|\big[(1+|\xi|^2)^s\widehat{u}\big]^{\vee}\|_{L^{p}(\R^d)}.$$
	Furthermore, we have the space $H^{s,p}(\Omega):=\Big\{u\mid_{\Omega}:\, u\in H^{s,p}(\mathbb{R}^d)\Big\},$
	equipped with the norm
	$$\|u\|_{H^{s,p}(\Omega)}=\inf\{\|v\|_{H^{s,p}(\R^d)} :\, u=v\mid_{\Omega}\, , v\in H^{s,p}(\R^d)\}.$$
\end{definition}

\vspace{2mm}

\noindent Let us now illuminate some correlations between the spaces $W^{s,p}(\R^d)$ and $H^{s,p}(\R^d)$.

\begin{itemize}
\item By using the Fourier transform (see the presentation in Chapter \ref{chap:basics-nonlocal-operator}) one gets that for all $s\in \R$,  $W^{s,2}(\R^d)= H^{s,2}(\R^d)$. 

\item According to  Mikhlin's  Theorem \cite[Theorem 2.6.1]{Zi12} we have $W^{m,p}(\R^d)= H^{m,p}(\R^d)$ for all $m\in \mathbb{Z}$ and $1<p<\infty$. Moreover the equality holds with equivalent norms. This fails for $p=1$ or $p=\infty$.

\item Furthermore, according to \cite{grisvard11},  for all $s\in \R^d$ we have
\begin{align*}
&W^{s,p}(\mathbb{R}^d) \subset H^{s,p}(\mathbb{R}^d)\quad\text{for $1<p\leq 2$}\quad \text{and}\quad H^{s,p}(\mathbb{R}^d) \subset W^{s,p}(\mathbb{R}^d)\quad\text{for $2\leq p<\infty$}.
\end{align*}

\item  In addition Lions and Peetre (1964) \cite{grisvard11} have proved that if $s_1<s_2<s_3$ then 
$$W^{s_3,p}(\mathbb{R}^d)\subset H^{s_2,p}(\mathbb{R}^d) \subset W^{s_1,p}(\mathbb{R}^d) .$$ 

\item  If $\Omega$ is a Lipschitz domain then combining Proposition \ref{prop:nonlocal-embeding} and Theorem \ref{thm:w1p-in nonlocal} it is easy to establish the following continuous embedding  for all $1\leq s_1\leq s_2$
$$W^{s_2,p}(\Omega)\hookrightarrow W^{s_1,p}(\Omega).$$  
Similarly, for all $1\leq s_1\leq s_2$, one can show that 
$H^{s_2,p}(\Omega)\hookrightarrow H^{s_1,p}(\Omega).$  
\end{itemize}

\noindent Let us now quote some interesting results in connection to the spaces $W^{s,p}(\Omega)$. 

\begin{theorem}[\cite{Dorothee-Triebel, grisvard11}] Let $\Omega\subset \mathbb{R}^d$ be open.
\begin{enumerate}[$(i)$]	
\item $C_c^\infty(\mathbb{R}^d)$ is dense in $W^{s,p}(\mathbb{R}^d)$ for all $s\in \R$. 
\item  If $\Omega$ has a continuous boundary, then the space  $W_{*}^{s,p}(\Omega)= \{u\in W^{s,p}(\mathbb{R}^d):\, u=0~\text{a.e on}~\Omega^c \}$ contains $C_c^\infty(\Omega)$ as a dense subspace.
\item  	Assume $\Omega$ has a continuous boundary, then $C_c^\infty(\overline{\Omega})$ is dense in $W^{s,p}(\Omega)$ for all $s>0$, where  we recall that $C_c^\infty(\overline{\Omega})$ is the space of all functions which are restriction of $C^\infty$ functions with compact support in $\mathbb{R}^d$ to $\Omega$. 
\item  	Assume $\Omega$ has  a Lipschitz boundary, then $C_c^\infty(\Omega)$ is dense in $W^{s,p}(\Omega)$ for $0\leq s\leq 1/p$, i.e.,  $W_0^{s,p}(\Omega)= W^{s,p}(\Omega)$ for $s\leq 1/p$. 
		
		
\item 	Let $\Omega$ be a bounded open subset of $\mathbb{R}^d$ with a Lipschitz boundary, then $ W_{*}^{s,p}(\Omega)= W_0^{s,p}(\Omega)$ except if $s-1/p$ is an integer. 
		Furthermore, if $0<s<1/p$ , $ W_{*}^{s,p}(\Omega)= W_0^{s,p}(\Omega) = W^{s,p}(\Omega)$.
\end{enumerate}
\end{theorem}
	
\noindent The proofs of $(i) -(iv)$ are included in Section \ref{sec:density} in a more general context.
 \begin{remark}
	To some authors  e.g \cite{FKV15, Dyda-Kasmmann}, the space $W_0^{s,p}(\Omega)$  is the closure in  $W^{s,p}(\mathbb{R}^d)$ of $C^\infty_0(\Omega)$
	in this  case we will denote  simply by $\overline{ C_0^{\infty}(\Omega)}^{W^{s,p}(\mathbb{R}^d)}$.    Whereas others e.g \cite{Hitchhiker, Dorothee-Triebel}
	define $W_0^{s,p}(\Omega)$  to be closure of $C^\infty_0(\Omega)$ in  $W^{s,p}(\Omega)$. One  should be careful since  (see $(v)$) in general  both  definitions do not coincide especially on non-smooth domains.
\end{remark}

\vspace{2mm}

\section{Approximations by smooth functions}\label{sec:density}
In practice it is rather laborious and demanding to work in a function space not containing smooth functions. In this section we show from different perspectives  that functions from the nonlocal spaces of interest for us can be  approximated by smooth functions. This part is stimulated by \cite{EE87}. Throughout this section, the function $\nu: \R^d \setminus\{0\} \to [0,\infty]$ is assumed to satisfy the condition \eqref{eq:plevy-integrability} see (page \pageref{eq:plevy-integrability}). In addition we fix a function $\phi \in C_c^{\infty}(\mathbb{R}^d)$ supported in $B_1(0)$ such that $\phi \geq 0$ and $ \int_{B_1} \phi = 1$ we denote the corresponding mollifier by $\phi_\delta(x)= \frac{1}{\delta^d}\phi\left(\frac{x}{\delta}\right)$. It is not difficult to establish that for each $\delta>0$,
\begin{align*}
\phi_\delta\in C_c^\infty(\mathbb{R}^d),\quad \phi_\delta\geq 0,\quad \operatorname{supp}\phi_\delta\subset B_\delta(0)\quad\text{and}\quad \int_{\mathbb{R}^d}\phi_\delta(x)\,\d x=1.
\end{align*}
A standard example is given by taking
 $\phi(x)= c\exp{-\tfrac{1}{1-|x|^2}}$ if $|x|<1$ and $\phi(x)=0$ if $|x|\geq 1$ where the constant $c>0$ is chosen such that $\int_{\mathbb{R}^d}\phi(x)\,\d x=1.$

\vspace{7mm}

\noindent We start by recalling some basic remarks and some useful Lemmata. 
\begin{remark} Note that if we notationally set $U(x,y)= (u(x)- u(y))\nu^{1/p}(x-y)$ then 
	%
	\begin{itemize}
	\item  $u\in \WnuOm$  if and only if $(u,U) \in L^p(\Omega)\times L^p(\Omega \times \Omega) $
	\item $u\in \WnuOmR$  if and only if $(u,U) \in L^p(\Omega)\times L^p(\Omega \times \R^d) $.
	\end{itemize}
\noindent Employing  Jensen's inequality we get
\begin{align}\label{eq:estimate-conv}
\begin{split}
|\phi_\delta*u -u|^p_{\WnuOm} 
&= 
\iil_{\Omega\Omega} \big|[\phi_\delta*u-u](x)- [\phi_\delta*u-u](y)\big|^p \nu(x-y)\d y \d x\\
&= \iil_{\Omega\Omega} \Big|\il_{ \R^d} \phi_\delta(z)\big\{ U(x-z, y-z) -U(x,y) \big\}\d z\Big|^p \d y \d x\\
&=  \iil_{\Omega\Omega} \Big|\il_{ \R^d} \phi(z) \big\{  U(x-\delta z, y-\delta z) -U(x,y) \big\}\d z\Big|^p \d y \d x\\
&\leq  \il_{ B_1} \phi(z)\d z \iil_{\Omega\Omega} \big| U(x-\delta z, y-\delta z) -U(x,y)\big|^p \d y \d x
\end{split}
\end{align}
\end{remark}

%
%
%
%

\begin{lemma}\label{lem:local-regularization} Assume $\Omega\subset \R^d$ is an open set.  Let $\Omega' \subset  \Omega$  be open and such that $\operatorname{dist}(\Omega', \partial \Omega)>0$ (if $\Omega\neq \R^d$).  For $u\in \WnuOm$ then $\phi_\delta*u \in C^\infty(\Omega') \cap W^p_\nu(\Omega')$ for all $0<\delta< \dist(\Omega', \partial \Omega)$ and 
\begin{align*}
\|\phi_\delta*u-u\|_{W^p_\nu(\Omega')} \xrightarrow{\delta\to 0}0.
\end{align*}
The same is true with $W^p_\nu(\Omega')$ replaced with $W^p_\nu(\Omega'|\R^d)$. Moreover if $\Omega=\R^d$ then one can take $\Omega'= \R^d$. 
\end{lemma}

\vspace{2mm}

\begin{proof} We will only prove the first statement and the others statements will follow analogously. Fix $0<\delta< \dist(\Omega', \partial \Omega)$. Note that for $x,y\in \Omega'$ and $z\in B_1$, then $x-\delta z, y-\delta z\in \Omega'+ B_\delta\subset \Omega$. Hence, by mimicking the estimate  \eqref{eq:estimate-conv} and using a trivial change of variables we find that 
\begin{align*}
\begin{split}
|\phi_\delta*u|^p_{W^p_\nu(\Omega')} 
&\leq  \il_{ B_1} \phi(z)\d z \iil_{\Omega'\Omega'} \big| U(x-\delta z, y-\delta z)\big|^p \d y \d x \\
&\leq   \il_{ B_1} \phi(z)\d z \iil_{\Omega'+ B_\delta\, \Omega'+ B_\delta} \big| U(x,y)\big|^p \d y \d x
\leq |u|^p_{\WnuOm}<\infty.
\end{split}
\end{align*}
Thus, $\phi_\delta*u \in W^p_\nu(\Omega')$.  By continuity of the shift in $L^p(\R^d\times \R^d),$
\[ \iil_{\Omega'\Omega'}\Big| U(x-\delta z, y-\delta z) -U(x,y)\Big|^p \d y \d x \xrightarrow{\delta\to 0}0. \]
Further, we have 
\[
\phi (z)  \iil_{\Omega'\Omega'}\Big| U(x-\delta z, y-\delta z) -U(x,y)\Big|^p \d y \d x  \leq 2^p \phi(z) \iil_{\Omega\Omega}\big|U(x,y)\big|^p \d y \d x\leq 2^p \phi(z)\|u\|^p_{\WnuOm}\in L^1(\R^d)
\]
By dominated convergence we conclude that 
$|\phi_\delta*u -u|^p_{W^p_\nu(\Omega')} \xrightarrow{\delta\to 0}0$ and according  to Theorem \ref{thm:approx-smooth} we know that  $\|\phi_\delta*u -u\|_{L^p(\Omega')} \xrightarrow{\delta\to 0}0$ and so we get $\|\phi_\delta*u-u\|_{W^p_\nu(\Omega'|\R^d)} \xrightarrow{\delta\to 0}0.$
%
\end{proof}

\vspace{2mm}
\noindent The next lemma plays a determinant  role in approximation functions. 
\begin{lemma}\label{lem:cutt-off} Assume $\Omega \subset \R^d$ is open. The following assertions are true 
\begin{enumerate}[$(i)$]
	\item Let $\Omega' \subset  \Omega$  be open and such that $\operatorname{dist}(\Omega', \partial \Omega)>0$. For $u\in \WnuOm$ with $\operatorname{supp} u\subset \Omega'$ the zero extension $\overline{u}$ of $u$ outside $\Omega$ belongs to $W^p_\nu(\R^d)$. Moreover, there is $C>0$  independent of $u$ with
	\begin{align*}
	\|\overline{u}\|_{W^p_\nu(\R^d)}  \leq C \|u\|_{\WnuOm}.
	\end{align*}
	 
	\item  Let $\varphi \in C_c^\infty(\R^d)$  then for all $u\in L^p(\Omega)$ and all $x,y\in \R^d,$ we have the estimate
	\begin{align}\label{eq:cutt-off-estimate}
	\begin{split}
	|\varphi(x)u(x) -\varphi(y)u(y)|^p &\leq 2^p\|\varphi\|^p_{W^{1,\infty}}\Big(\mathds{1}_{\operatorname{supp} \varphi}(y)|u(x) -u(y)|^p+ |u(x)|^p (1\land |x-y|^p)\Big).
	\end{split}
	\end{align}
\item Let $\varphi \in C_c^\infty(\R^d)$.  If $u\in \WnuOm$ then  $\varphi u\in \WnuOm$ and for some constant  independent of $u$, 
 \begin{align*}
\|\varphi u\|_{\WnuOm}  \leq C \|u\|_{\WnuOm}.
\end{align*}
The same holds with $\WnuOm$ replaced by $\WnuOmR$. 

\item Let $\varphi \in C_c^\infty(\Omega)$.  If $u\in \WnuOm$ then  $\varphi u\in W^p_\nu(\R^d)$ and for some constant  independent of $u$
\begin{align*}
\|\varphi u\|_{W^p_\nu(\R^d)}  \leq C \|u\|_{\WnuOm}.
\end{align*}
The same holds true with $\WnuOm$ replaced by $\WnuOmR$. 
\end{enumerate}
\end{lemma}

\vspace{2mm}

\begin{proof} Observe that $(iv)$ follows by combining $(i)$ and $(iii)$. Next, let us prove $(i)$.  Set $\delta= \operatorname{dist}(\Omega', \partial \Omega)>0$ and $C_\delta=\int_{|h|\geq \delta} \nu(h)\d h>0$. It follows that
	\begin{align*}
	\iil_{\R^d\R^d} |\overline{u}(x) -\overline{u}(y)|^p\nu(x-y)\d y \d x
	&= \iil_{\Omega\Omega}  |u(x) -u(y)|^p\nu(x-y)\d y \d x + 2\int_{\Omega'}|u(x)|^p \d x\int_{\R^d\setminus \Omega} \nu(x-y)\d y\\
	&\leq |u|_{\WnuOm} + 2\int_{\Omega'}|u(x)|^p \d x\int_{|x-y|\geq \delta} \nu(x-y)\d y\\
	& =  |u|^p_{\WnuOm} + 2C_\delta\|u\|^p_{L^p(\Omega)}.
\end{align*}	
Whence, $\|\overline{u}\|^p_{W^p_\nu(\R^d)}  \leq (2C_\delta+ 1) \|u\|^p_{\WnuOm}$.  Now we prove $(ii)$. 
Observe that for all $x,y \in \R^d $ we have 
\begin{align*}
\begin{split}
|\varphi(x)u(x) -\varphi(y)u(y)| 
&\leq  |\varphi(y)||u(x) -u(y)| + |u(x)|  |\varphi(x)-\varphi(y)|\\
&= |\varphi(y)|\mathds{1}_{\operatorname{supp} \varphi}(y)|u(x) -u(y)| + |u(x)|  \Big|\int_0^1\nabla \varphi(x+t(y-x))\cdot (y-x)\Big|\\
&\leq \|\varphi\|_{L^\infty(\R^d)}\mathds{1}_{\operatorname{supp} \varphi}(y)|u(x) -u(y)| + \|\varphi\|_{W^{1,\infty}}|u(x)|  (1\land |x-y|). 
\end{split}
\end{align*}
Thus implies $(ii)$ follows.
On the other hand,  $(iii)$ follows by integrating both sides of the estimate  \eqref{eq:cutt-off-estimate}, over $\R^d\times \R^d$,  with respect to the measure $\nu(x-y)\d y\d x$. Indeed, letting $\delta= \dist(\supp \varphi, \Omega^c)>0$ we have 
	\begin{align}\label{eq:estimate-extending-to-rd}
\begin{split}
\iil_{ \R^d \R^d} &|\varphi(x)u(x) -\varphi(y)u(y)|^p\nu(x-y)\d x \d y\\
&\leq 2^p \|\varphi\|^p_{W^{1,\infty}}\Big(\iil_{\Omega\Omega}|u(x) -u(y)|^p\nu(x-y)\d x \d y+ \il_{ \Omega}|u(x)|^p\d x\il_{\R^d}(1\land|h|^p)\nu(h) \d h\Big)\\
&+2 \il_{ \operatorname{supp}\varphi}|u(x)|^p \d x\il_{|h|\geq \delta}\nu(h) \d h.
\end{split}
\end{align}

\end{proof}

\vspace{2mm}

\begin{lemma}\label{lem:truncate-approx}
Assume that $\Omega \subset \R^d$ is open. Let $\chi \in C_c^\infty(\R^d)$ such that $ 0\leq \chi \leq 1$ and $\chi(x)=1$ for $|x|<1$ and  $\chi(x)=0$ for $|x|\geq 2$. If $u \in \WnuOm$ then $\chi(\cdot/j)u\in \WnuOm$ and
 $$\|\chi(\cdot/j)u-u \|_{\WnuOm} \xrightarrow{j\to \infty}0.$$ 
The same holds with $\WnuOm$ replaced by $\WnuOmR$. 
\end{lemma}

\medskip

\begin{proof} Pointwise we have $\chi(x/j) \xrightarrow{j\to \infty} 1$. So that for $u \in\WnuOm$,  $\chi(x/j)u(x)\xrightarrow{j\to \infty} u(x)$.  We know that $\chi(\cdot/j) u\in \WnuOm$. On the other hand, $|\chi(\cdot/j)u|\leq |u|$ then $|\chi(\cdot/j) u-u|_{L^p(\Omega)} \xrightarrow{n\to \infty}0$ and using the estimate \eqref{eq:cutt-off-estimate}, for all $x,y \in \R^d $ and $j\geq 1$ we find that
	\begin{align*}
	|\chi(x/j) u(x) -\chi(y/j)u(y)|^p\nu(x-y) &\leq|u(x) -u(y)|^p\nu(x-y) \\
	&+ 2^p\|\chi\|^p_{W^{1,\infty}}|u(x)|^p  (1\land |x-y|^p)\nu(x-y) \in L^1(\Omega \times \Omega).
	\end{align*}
The dominated convergence  theorem implies that $|\chi(\cdot/j) u-u|_{\WnuOm} \xrightarrow{j\to \infty}0$ and thus, it follows that $\|\chi(\cdot/j)u-u\|_{\WnuOm} \xrightarrow{j\to \infty}0.$ The case $u\in \WnuOmR$ follows likewise. 
\end{proof}


\noindent Let us immediately start with a simple case. 
\begin{theorem} \label{thm:smooth-Rd} Let $1\leq p<\infty$ and $\nu: \R^d\setminus\{0\} \to [0,\infty]$ satisfying \eqref{eq:plevy-integrability}. Then $C_c^\infty(\R^d)$ is dense in $W^p_\nu(\R^d)$, i.e. for a function $u\in W^p_\nu(\R^d)$, there exist functions $u_n\in  C_c^\infty(\R^d)$ such that 
	$$\|u_n-u\|_{W^p_\nu(\R^d)} \xrightarrow{n\to \infty}0.$$
	
\end{theorem}


\begin{proof}
Let $\chi \in C_c^\infty(\R^d)$ such that $ 0\leq \chi \leq 1$ and $\chi(x)=1$ for $|x|<1$ and  $\chi(x)=0$ for $|x|\geq 2$.  Letting $\overline{u}_j= \chi(\cdot/j)u$, by Lemma \ref{lem:truncate-approx} we can find  $j\geq 1$ such that $\|\overline{u}_j-u\|_{W^p_\nu(\R^d)} <\eps/2$ with $\eps>0$. Note that $\phi_\delta*\overline{u}_j\in  C_c^\infty(\R^d)$.   Setting $U_j(x,y)= (\overline{u}_j(x)- \overline{u}_j(y))\nu^{1/p}(x-y) \in L^p(\R^d\times \R^d).$ By mimicking the estimate \eqref{eq:estimate-conv} for $\Omega= \R^d$  we have
\begin{align*}
|\phi_\delta*\overline{u}_j -\overline{u}_j|^p_{W^p_\nu(\R^d)} 
\leq  \il_{ B_1} \phi(z) \d z \iil_{\R^d\R^d} \big| U_j(x-\delta z, y-\delta z) -U_j(x,y)\big|^p \d y \d x\,. 
\end{align*}
By the continuity of the shift in $L^p(\R^d\times \R^d)$
\[ \iil_{\R^d\,\R^d}\Big| U_j(x-\delta z, y-\delta z) -U_j(x,y)\Big|^p \d y \d x \xrightarrow{\delta\to 0}0. \]
Further, we have 
\[
\phi (z)  \iil_{\R^d\,\R^d}\Big| U_j(x-\delta z, y-\delta z) -U_j(x,y)\Big|^p \d y \d x  \leq 2^p \phi(z) \iil_{\R^d\,\R^d}\big|U_j(x,y)\big|^p \d y \d x = 2^p \phi(z)|u|^p_{W^p_\nu(\R^d)} \in L^1(\R^d).
\]
By dominated convergence we conclude that $|\phi_\delta*\overline{u}_j -\overline{u}_j|^p_{W^p_\nu(\R^d)} \xrightarrow{\delta\to 0}0$ and according  to Theorem \ref{thm:approx-smooth} we know that  $\|\phi_\delta*\overline{u}_j -\overline{u}_j\|_{L^p(\R^d)} \xrightarrow{\delta\to 0}0$. Accordingly we can find $\delta >0$ sufficiently small such that $\|\phi_\delta*\overline{u}_j -\overline{u}_j\|_{W^p_\nu(\R^d)}<\eps/2$.
Finally we get $\phi_\delta*\overline{u}_j\in C_c^\infty(\R^d)$ and  it is not difficult to show that $\phi_\delta*\overline{u}_j\in W^p_\nu(\R^d)$ and we have 
\begin{align*}
\|\phi_\delta*\overline{u}_j -u\|_{W^p_\nu(\R^d)}\leq \|\overline{u}_j -u\|_{W^p_\nu(\R^d)}+ \|\phi_\delta*\overline{u}_j -\overline{u}_j\|_{W^p_\nu(\R^d)} <\eps.
\end{align*}
\end{proof}


\noindent We now deal with what can be seen as the Meyers-Serrin density type result for nonlocal spaces. 

\begin{theorem}\label{thm:meyer-serrin} Let $1\leq p<\infty$ and $\nu: \R^d\setminus\{0\} \to [0,\infty]$  satisfying \eqref{eq:plevy-integrability}.  Then $ C^\infty(\Omega) \cap \WnuOm$ is dense in $\in \WnuOm$.
\end{theorem}

\begin{proof}
We shall only assume $\Omega\neq \R^d$ because the case $\Omega= \R^d$ is greatly covered by Theorem  \ref{thm:smooth-Rd}.  We follow the standard arguments as for the classical Sobolev spaces. 
Set $$O_j = \big\{x\in \Omega: \dist(x, \partial \Omega)>2^{-j},\,\, |x|<2^j\big\} $$
 and $V_1= O_3$ and $V_j =O_{j+2} \setminus \overline{O}_{j}$ for $j\geq 2$. We have $\Omega= \bigcup\limits_{j=1}^{\infty} V_j$. Consider $(\xi_j)_j$, a partition of the unity subordinate to the family $(V_j)_j$. That is,
\begin{align*}
\xi_j\in C_c^\infty(\R^d),\quad\supp(\xi_j)\subset V_j, \quad 0\leq \xi_j \leq 1, \quad \sum_{j=1}^\infty \xi_j =1. 
\end{align*}
According to Lemma \ref{lem:cutt-off}  for $u\in \WnuOm$, each $u_j=\xi_j u \in \WnuOm$. We also have $\supp u_j \subset V_j\subset  \Omega_{j+2}$.  
Chose $\delta>0$ small enough so that $\supp u_j\cup \supp(\varphi* u_j ) \subset O_{j+3} \setminus O_{j-1}$. Note that $\dist(O_{j+3}, \partial \Omega)>2^{-j-3}$.  In accordance to Lemma \ref{lem:local-regularization} we have
\begin{align*}
\|\phi_\delta*u_j-u_j\|_{W^p_\nu(O_{j+4})} \xrightarrow{\delta\to 0}0.
\end{align*}
On the other hand, since $\dist(O_{j+3},  \Omega \setminus O_{j+4})\geq 2^{-1}$, letting $C= \int_{|h|\geq 2^{-1}} \nu(h)\d h>0 $ we get
\begin{align*}
\int_{O_{j+3}} |\phi_\delta*u_j(x)-u_j(x)|^p\d x\int_{\Omega\setminus O_{j+4}}\hspace*{-4ex}\nu(x-y)\d y\leq C\int_{O_{j+3}} |\phi_\delta*u_j(x)-u_j(x)|^p\d x \xrightarrow{\delta\to 0}0.
\end{align*}
Altogether we have $\|\phi_\delta*u_j-u_j\|^p_{\WnuOm}\xrightarrow{\delta\to 0}0$ because
\begin{align*}
\|\phi_\delta*u_j-u_j\|^p_{\WnuOm} = \|\phi_\delta*u_j-u_j\|^p_{W^p_\nu(O_{j+4})} +  
2\int_{O_{j+3}} |\phi_\delta*u_j(x)-u_j(x)|^p\d x
\int_{\Omega\setminus O_{j+4}} \hspace*{-4ex}\nu(x-y)\d y.
\end{align*}
Accordingly, given $\varepsilon>0$ we can find $\delta_j>0$ such that 
\begin{align*}
\|\phi_{\delta_j}*u_j-u_j\|^p_{\WnuOm} \leq \varepsilon2^{-j}.
\end{align*}
Given that $O_{j+3} \setminus O_{j-1} { }^{'s}$ can only overlap at most five times and $ \phi_\delta*u_j\in C_c^\infty(O_{j+3} \setminus O_{j-1})$, the function $v= \sum\limits_{j=1}^\infty \phi_{\delta_j}*u_j$ is well defined and belongs to  $C^\infty(\R^d)$. Noticing that $u = \sum\limits_{j=1}^\infty \xi_j u= \sum\limits_{j=1}^\infty u_j$, from the above we get on the one hand that $\|v-u\|_{\WnuOm} \leq \varepsilon$ since
\begin{align*}
\|v-u\|_{\WnuOm} &= \Big\| \sum_{j=1}^\infty   \phi_{\delta_j}*u_j-u_j\Big\|_{\WnuOm}
\leq  \sum_{j=1}^\infty \| \phi_{\delta_j}*u_j-u_j\|_{\WnuOm}\|
\leq \sum_{j=1}^\infty  \varepsilon2^{-j}= \varepsilon. 
\end{align*} 
 And on the other hand, $v\in C^\infty(\R^d)\cap \WnuOm$ since 
 \begin{align*}
\|v\|_{\WnuOm} \leq  \|v-u\|_{\WnuOm} + \|u\|_{\WnuOm} \leq \|u\|_{\WnuOm} + \varepsilon.
\end{align*}
This completes the proof. 
\end{proof}
\vspace{2mm}

\noindent Note that the main point in the above proof is that if $u\in \WnuOm$ then by truncation it is possible to shrink the support of $u$ inside $\Omega$ using a suitable partition of unity, so that  convoluting makes sense afterwards. On the other hand, the question of approximating a function of $\WnuOmR$ by smooth functions is more delicate, but can be reduced to the question whether  $\WnuOmR$ is closed under the shift operator, i.e. for $u\in \WnuOmR$, do we have $u(\cdot-h)=\tau_h u\in \WnuOmR$?  This is equivalent to say that for $\delta>0$ sufficiently small, $u\in W^p_\nu(\Omega+B_\delta|\R^d)$ when $u\in \WnuOmR$. Whereas this is not totally obvious because in this situation the function $u$ is already defined on the whole $\R^d$.  Next we give a partial answer to this question when $\Omega$ has a compact Lipschitz boundary.

\vspace{2mm}

\begin{definition}\label{def:domain-regular}  An open set $\Omega\subset \R^d$ shall be called a local graph domain if for each $x\in \partial \Omega$ there exists a $r>0$ and  a function $\gamma:\R^{d-1}\to \R $ such that 
	\begin{align*}
	\Omega\cap B_r(x) &= \big\{ x\in B_r(x):\, x_d>\gamma(x')\big\}\\
	\partial\Omega\cap B_r(x) &= \big\{ x\in B_r(x):\, x_d=\gamma(x')\big\}\\
	\Omega^c\cap B_r(x) &= \big\{ x\in B_r(x):\, x_d\leq \gamma(x')\big\}.
	\end{align*}
	
	\noindent In addition, if $\gamma$ is of class $C^{m,\sigma}$ then $\Omega$ is called to be a $C^{m,\sigma}$-domain. For $m=0,\sigma=1$, $\Omega$ is called a Lipschitz domain. 
\end{definition}


\noindent We will need the following result. 
\begin{lemma}\label{lem:locp-integrable}
	Let $1\leq p<\infty$ and assume $\nu: \R^d\setminus\{0\} \to [0,\infty]$  satisfying \eqref{eq:plevy-integrability} is radially almost decreasing, i.e. $\nu(x) \leq c\nu(y)$ ($c>0$ ) whenever $|y|\geq |x|$.  Then for $\Omega\subset \R^d$ open, $\WnuOmR\subset L^p_{\operatorname{loc}}(\R^d). $ Moreover for each compact $K\subset \R^d$ there is a constant $C= (K, \Omega, \nu)$ such that for all $u\in \WnuOm$,
	\begin{align*}
	\int_{K} |u(x)|^p\d x \leq C\|u\|^p_{\WnuOm}.
	\end{align*}
\end{lemma}

\vspace{2mm}

\begin{proof}
If $u\in \WnuOmR$ then $u\in L^p(\Omega)$. For a compact set $K\subset \R^d$ write $K=K_1\cup K_2$ with $K_1= \Omega\cap K$ and $K_2= \Omega^c\cap K$. Now choose $K'\subset \Omega$ to be any compact set and $R>0$ sufficiently large such that $K'\cup K_2\subset B_R(0)$.
	Clearly $u\in L^p(K_1)$ and  $u\in L^p(K')$. It remains to show that $u \in L^p(K_2).$  For every  $x\in K'$ and $y\in K_2$ we have $|x-y|\leq R$ so that $\nu(x-y)\geq c_R$ for some constant $c_R >0$. Applying this and Jensen's inequality we get the following 
	\begin{align*}
	\infty>\iil_{ \Omega\R^d}|u(x) -u(y)|^p\nu(x-y)\d y\d x
	&\geq c_R \int_{K'} \int_{K_2} |u(x) -u(y)|^p\d y\d x \\
	&\geq c_R |K'|\int_{K_2} |u(x) -\mbox{$\fint_{K'}$}u|^p \d x. 
	\end{align*}
	The conclusion is reached because
	\begin{align*}
	\int_{K_2} |u(x)|^p\d x &\leq 2^p |K_2| \fint_{K'} |u(x)|^p\d x + 2^p \int_{K_2} |u(x) -\mbox{$\fint_{K'}$}u|^p \d x\\
	&\leq 2^p (c_R^{-1}+1)|K_2|  |K'|^{-1} \|u\|^p_{\WnuOm}<\infty.
	\end{align*}
\end{proof}

\bigskip

\noindent The following density result ameliorates the analogous one from \cite{Voi17} whose proof is more likely not to be fully satisfactory. 
\begin{theorem}\label{thm:density}
	Assume $\Omega\subset \mathbb{R}^d$ is open with a compact Lipschitz boundary $\partial \Omega$. Let $\nu$ satisfy \eqref{eq:plevy-integrability} and in addition, assume $\nu$ is radially almost decreasing. Then $C_c^\infty(\mathbb{R}^d)$ is dense in $\WnuOmR $ with respect to the norm $\|\cdot\|_{\WnuOmR}$, i.e. for $u \in \WnuOmR$ there exists a sequence $(u_n) \subset C_c^\infty(\mathbb{R}^d)$ with 
	\begin{align*}
	\|u_n-u\|_{\WnuOmR} \longrightarrow 0 \text{ as } n\to\infty \,.
	\end{align*} 
\noindent Moreover for $p=2$, then $C_c^\infty(\mathbb{R}^d)$ is dense in $V^1_\nu(\Omega|\R^d) $ with respect to the norm $$\|u\|^2_{V^1_\nu(\Omega|\R^d)}= \|Lu\|^2_{L^2(\Omega)}+ \|u\|^2_{\VnuOm}. $$
\end{theorem}

\medskip

\begin{proof}
Let $u \in \WnuOmR$. We prove that  $|v_\varepsilon-u|_{\WnuOm}\xrightarrow{\eps\to0}0$ where$v_\varepsilon\in C^\infty_c(\R^d)$. This implies
$
\|v_\varepsilon-u\|_{ \WnuOmR}  \xrightarrow{\eps\to0}0
$
because convergence of the $L^p$-norms follows by standard arguments. Note that the sequence $(v_\varepsilon)$ is constructed by translation and convolution of the function $u$ with a mollifier
Thus for $u\in V^1_\nu(\Omega|\R^d)$ since $\partial \Omega$ is Lipschitz, taking into account Proposition  \ref{prop: properties-Levy-op}, by a subsequence shift arguments one can show that $\|Lu-Lv_\eps\|_{L^2(\Omega)}\xrightarrow{\eps\to 0}0$. 
\medskip

\noindent \textbf{Step 1:} Let $x_0\in\partial\Omega$. Since is $\partial\Omega$ Lipschitz, there exists $r>0$ and a Lipschitz  function $\gamma:\R^{d-1}\to\R$ with Lipschitz constant $k>0$, such that (upon relabelling the coordinates) 

\begin{align*}
\Omega\cap B_r(x_0) &= \{ x=(x',x_d) \in B_r(x_0)| x_d > \gamma(x')\}.
\end{align*}
For the sake of convenience, we choose $r>0$ so small such that $|\Omega\cap B^c_r(x_0)|>0$. For $x\in B_{r/2}(x_0)$, $\tau >1+k$  and $0<\eps <\frac{r}{2(1+\tau)}$  we define the shifted point $x_\eps = x +\tau \eps e_d\,.$ We define 
\begin{align*}
u_\eps(x) = u(x_\eps)= u(x+\tau \eps e_d)\quad\text{and}\qquad v_\eps = \phi_\eps \ast u_\eps,
\end{align*} 
where $\phi_\eps$ is a smooth mollifier having support in $B_\eps(0)$. It is noteworthy to  emphasize that the smoothness of the function $\phi_\eps* u_\eps$ is warranted by Lemma \ref{lem:locp-integrable}, which assures that $u\in L^p_{\operatorname{loc}}(\R^d)$ when $u\in \WnuOmR$ and $\nu$ is radially almost decreasing. Indeed we have $C_c^\infty(\R^d)* L^1_{\operatorname{loc}}(\R^d)\subset C^\infty(\R^d)$.  
 
\medskip

\noindent \textbf{Step 2:} Let us assume 
$\supp u\Subset B_{r/4}(x_0)$. In this case $v_\eps \in C^\infty_c (B_r(x_0))$. In this step, we aim to prove 
\[ 
|v_\eps - u|_{\WnuOmR} \longrightarrow 0 \quad \text{ as } \eps \to 0 \,. 
\]
Due to the nonlocal nature of the semi-norm, this step turns out to be rather challenging. We begin with a geometric observation.

\begin{lemma}\label{lem:guy-geo}
	Let $z\in B_1(0)$. Let $\Omega^z_\eps = \Omega+ \eps(\tau e_d-z)$. Then $\Omega^z_\epsilon\cap B_{r/2}(x_0) \subset \Omega\cap B_{r}(x_0)$. 
\end{lemma}

\begin{proof}
	For $h \in \Omega^z_\epsilon\cap B_{r/2}(x_0)$, let us write $ h= t+\eps \tau e_d-\eps z $ with $t \in \Omega$. Note that  since $\eps<\frac{r}{2(\tau+1)}$ we get $|t-x_0|\leq|t-h|+|h-x_0| < \eps(\tau+1) +r/2 <r$. So that  $t\in B_{r}(x_0) $, $h' = t'-\eps z'$ and $h_d= t_d+ \eps(\tau -z_d)$. Since $\gamma$ is Lipschitz with Lipschitz constant $k<\tau-1$ and $t\in \Omega\cap B_{r}(x_0) = \{ x\in B_{r}(x_0)| x_d > \gamma(x')\}$, i.e. $t_d>\gamma(t')$,  we obtain
	\begin{align*}
	\gamma(h')&\leq \gamma(t')+ |\gamma(h')-\gamma(t')|  <t_d+ \eps k|z'|\\
	& <t_d+ \eps k< t_d+ \eps(\tau -z_d)=  h_d.
	\end{align*}
	Thus, $h \in B_r(x_0)$ and $h_d>\gamma(h')$. We have shown $h\in \Omega \cap B_r(x_0)$ as desired. 
\end{proof}

The main technical tool of the argument below is the Vitali convergence theorem, (see Theorem \ref{thm:vitali} in Appendix). Since $u$ belongs to the space $\WnuOmR$, for every $\delta > 0$ there is $\eta > 0$ such that for all sets $E \subset \Omega$, $F \subset \R^d$ with $|E \times F| < \eta$ we know
\begin{align}\label{eq:equi-int-u}
\iil_{E F} \big|u(x)-u(y)\big|^p \nu(x-y) \d y \d x &< \delta \text{ and } \int_{EF}  |u(y)|^p \d y < \delta \,.
\end{align}
The second estimate uses the fact that $u$ has compact support and belongs to $L^p_{\operatorname{loc}}(\R^d)$.

\begin{lemma}\label{lem:equi-int-u-eps}
	For every $\delta > 0$ there is $\eta > 0$ such that for all sets $E \subset \Omega$, $F \subset \R^d$ with $|E \times F| < \eta$ 
	\begin{align}\label{eq:equi-int-u_eps-z}
	\sup\limits_{z \in B_1(0)} \sup\limits_{\eps > 0} \iil_{E F} \big|u^z_\eps(x)-u^z_\eps(y)\big|^p \nu(x-y) \d y \d x < \delta \,,
	\end{align}
	where $ u^z_\eps(\xi) = u_\eps(\xi-\eps z) = u(\xi+\eps \tau e_d-\eps z)$.
\end{lemma}

\begin{proof}
	Let $\delta > 0$. Choose $\eta > 0$ as in \eqref{eq:equi-int-u}. Let $\eps > 0$, $z \in B_1(0)$. Let $E \subset \Omega$, $F \subset \R^d$ be sets with $|E \times F| < \eta$. Then
	\begin{align}
	\iil_{E F} \big|u^z_\eps(x)-u^z_\eps(y)\big|^p \nu(x-y) \d y \d x = \iil_{E^z_\eps F^z_\eps} \big|u(x)-u(y)\big|^p \nu(x-y) \d y \d x \,,
	\end{align}
	where $E^z_\eps = E + \eps(\tau e_d-z)$ and $F^z_\eps$ defined analogously. We decompose $E^z_\eps$ as follows $E^z_\eps = E^z_\eps \cap B_{r/2}(x_0) \cup E^z_\eps \cap B^c_{r/2}(x_0)$. Note 
	\[ E^z_\eps \cap B_{r/2}(x_0) \subset \Omega^z_\eps \cap B_{r/2}(x_0) \subset \Omega \cap B_{r/2}(x_0) \,,\]
	where we apply Lemma \ref{lem:guy-geo}. We directly conclude
	\begin{align}\label{eq:equi-one}
	\iil_{E^z_\eps   F^z_\eps} \mathbbm{1}_{B_{r/2}(x_0)}(x) \big|u(y)-u(x)\big|^p \nu(x-y) \d y \d x \leq \delta \,.
	\end{align}
	With regard to the remaining term since $\operatorname{supp} u \subset  B_{r/4}(x_0)$, note
	\begin{align}\label{eq:equi-two}
	\begin{split}
	\iil_{E^z_\eps   F^z_\eps}& \mathbbm{1}_{B^c_{r/2}(x_0)}(x) \big|u(x)-u(y)\big|^p \nu(x-y) \d y \d x \\ 
	&= 
	\iil_{E^z_\eps   F^z_\eps} \mathbbm{1}_{B^c_{r/2}(x_0)}(x) \mathbbm{1}_{B_{r/4}(x_0)}(y) |u(y)|^p \nu(x-y) \d y \d x\\
	&\leq c(r, \nu) \iil_{E^z_\eps   F^z_\eps} \mathbbm{1}_{B^c_{r/2}(x_0)}(x) \mathbbm{1}_{B_{r/4}(x_0)}(y) |u(y)|^p \d y \d x \leq c \iil_{E^z_\eps   F^z_\eps}  |u(y)|^p \d y \d x\\
	&= c \iil_{E F^z_\eps}  |u(y)|^p \d y \d x\leq c \delta \,.
	\end{split}
	\end{align}
	The positive constant $c(r,\nu)$ depends on $r$ and on the shape of $\nu $. Summation over \eqref{eq:equi-one} and \eqref{eq:equi-two} completes the proof after redefining $\delta$ accordingly. 
\end{proof}

\noindent The next lemma shows the tightness of $u^z_\eps(x)-u^z_\eps(y)$ uniformly for $z\in B_1(0)$ and $\eps >0$.

\begin{lemma}\label{lem:tigh-u-eps-z}
	For every $\delta>0$ there exists $E(\delta) \subset \Omega$ and $F_\delta \subset \mathbb{R}^d$ such that $|E(\delta) \times F_\delta |<\infty$  and 
	\begin{align}\label{eq:tight-u_eps}
	\sup\limits_{z \in B_1(0)} \sup\limits_{\eps > 0} \iil_{(\Omega \times \mathbb{R}^d) \setminus ( E(\delta) \times F_\delta)} \big|u^z_\eps(x)-u^z_\eps(y)\big|^p \nu(x-y) \d y \d x < \delta.
	\end{align}
\end{lemma}

\begin{proof}
	Fix $\eps>0$ and $z\in B_1(0)$. Let $ \bar{R}= \sup\limits_{\xi \in \Omega} |\xi-x_0|$ which is finite since $\Omega $ is bounded.  Note that $\supp u^z_\eps \subset B_{r/2}(x_0)$. Choose $R>0$ so large such that $[B^c_{R}(x_0)]_{\eps}^z= B^c_{R}(x_0) +\eps (\tau e_d+z)\subset B^c_{R/2}(x_0) $  and  $|x-y|\geq R/2-\bar{R}$  for $x\in B^c_{R/2}(x_0) $ and $y \in \Omega$. Thus,
	\begin{align*}
	\iil_{(\Omega \times \mathbb{R}^d) \setminus ( \Omega \times B_R(x_0))} & \big|u^z_\eps(x)-u^z_\eps(y)\big|^p  \nu(x-y) \d y \d x
	= \iil_{\Omega B^c_R(x_0)} \big|u^z_\eps(x)\big|^p \nu(x-y) \d y \d x\\
	&= \il_{\Omega^z_\eps \cap B_{r/2}(x_0) } |u(x)|^p\d x \il_{ [B^c_{R}(x_0)]_{\eps}^z}  \nu(x-y) \d y \leq \il_{\Omega} |u(x)|^p\d x  \il_{ B^c_{R/2-\bar{R}}(x)}  \nu(x-y) \d y\\ 
	&= \|u\|^p_{L^p(\Omega)}  \il_{ B^c_{R/2-\bar{R}}(0)} \nu(h) \d h. 
	\end{align*}
	The  desired result follows by taking   $E(\delta) = \Omega$ and $F_\delta =B_R(x_0)$ with  $R>0$ large enough such that $ \il_{ B^c_{R/2-\bar{R}}(0)} \nu(h) \d h<\delta \|u\|^{-2}_{L^p(\Omega)}$ . 
\end{proof}
\begin{lemma} There exists a constant $C(\Omega,r, \nu)$ depending on $\Omega,r$ and $\nu$ such  that
	\begin{align}\label{eq:estimate-seminorm}
	|u^z_\eps|^p_{\WnuOmR}\leq C(\Omega,r, \nu) | u|^p_{\WnuOmR}\quad\text{for all $z\in B_1(0)$ and all $\eps>0$}. \quad
	\end{align} 
\end{lemma}

\begin{proof}
	Note that $|x-y|\geq r/4$ for $x\in B^c_{r/2}(x_0) $ and $y \in B_{r/4}(x_0)$.  Note that is $\nu$ radially almost decreasing, i.e. $\nu(x-y)<c(r,\nu)=c\nu(r/4) $. Let us choose
	\[ C= 1+ \hspace{-3ex}\sup_{y\in B_{r/4}(x_0)}\Big(\hspace{-1ex}
	\il_{\Omega\cap B_r^c(x_0)} \hspace{-2ex}\nu(x-y)\d x\Big)^{-1} \hspace{-2ex}\int\limits_{ B^c_{r/4}(0)} \hspace{-2ex}\nu(h)   \d h\leq 1+ (c(r,\nu)|\Omega\cap B_r^c(x_0)|)^{-1} \hspace{-3ex}\int\limits_{ B^c_{r/4}(0)} \hspace{-2ex}\nu(h)<\infty.\]
 
\noindent 	Therefore, for each $z\in B_1(0)$ and each $\eps>0 $ we have 
	\begin{alignat*}{2}
	&\iint\limits_{ \Omega^z_\epsilon\cap B^c_{r/2}(x_0)\times \R^d} \left|u(x)|^p - u(y) \right|^p \nu(x-y)  \d y \d x
	&&= \int\limits_{ B_{r/4}(x_0) } |u(y)|^p\d y \int\limits_{ \Omega^z_\epsilon\cap B^c_{r/2}(x_0)}  \nu(x-y)   \d x
	\\& \leq \int\limits_{ B_{r/4}(x_0) } |u(y)|^p\d y \int\limits_{ B^c_{r/4}(y)}  \nu(h)   \d h
	&&\leq  C\hspace{-3ex}  \int\limits_{ B_{r/4}(x_0) } |u(y)|^p\d y \int\limits_{ \Omega \cap B^c_{r}(x_0)}  \nu(x-y)   \d x\\
	&= C \hspace{-3ex}\iint\limits_{ \Omega \cap B^c_{r}(x_0)\times\R^d }|u(x)-u(y)|^p  \nu(x-y)  \d y \d x.\
	\end{alignat*}
	By applying change of variables, this  and Lemma \ref{lem:guy-geo}, we have 
	\begin{alignat*}{2}
	|u^z_\eps|^p_{\WnuOmR} &= \iint\limits_{\Omega\R^d} \big|u^z_\eps(x) - u^z_\eps(y) \big|^p \nu(x-y)  \d y \d x
	&&= \iint\limits_{ \Omega^z_\epsilon\R^d} \big|u(x) - u(y) \big|^p\nu (x-y)  \d y \d x\\
	&= \iint\limits_{ \Omega^z_\epsilon\cap B_{r/2}(x_0)\times \R^d}\hspace*{-2ex} \big|u(x)- u(y) \big|^p \nu(x-y)  \d y \d x
	&&+ \hspace{-3ex} \iint\limits_{ \Omega^z_\epsilon\cap B^c_{r/2}(x_0)\times \R^d} \big|u(x)- u(y) \big|^p \nu(x-y)  \d y \d x\\
	&\leq C\hspace{-3ex} \iint\limits_{ \Omega\cap B_{r}(x_0)\times \R^d}\hspace*{-2ex} \big|u(x) - u(y) \big|^p \nu(x-y)  \d y \d x
	&&+C\hspace{-3ex} \iint\limits_{ \Omega \cap B^c_{r}(x_0)\times\R^d }|u(x)-u(y)|^p  \nu(x-y)  \d y \d x\\
	&\leq  C |u|^p_{\WnuOmR}.
	\end{alignat*}
\end{proof}	

\noindent We are now in a position to prove the main result of this step. By mimicking  the estimate \eqref{eq:estimate-conv}, we get
\begin{align*}
&|v_\eps - u|^p_{\WnuOmR} = \iil_{\Omega\,\R^d} ((v_\eps(x)-v_\eps(y)) -(u(x)-u(y))|^p
\nu(x-y)\d y \d x\\
%
&\leq \int\limits_{B_1(0)}  \phi(z)  \d z \iint\limits_{\Omega \R^d} \big| (u_\eps(x-\eps z) - u_\eps(y-\eps z) ) 
-(u(x)-u(y))\big|^p\nu(x-y)  \d y \d x\, \,\\
&= \int\limits_{B_1(0)}  |u^z_\eps-u|^p_{\WnuOmR}\phi(z) \d z\,.
\end{align*}
\noindent For each fixed $z \in B_1(0)$ the family of functions 
$(x,y)\mapsto \big|(u^z_\eps(x) - u^z_\eps(y) ) -(u(x)-u(y)) \big|^p\nu(x-y) $ with $(x,y) \in \Omega\times\mathbb{R}^d$, $\eps>0$  is equiintegrable (by Lemma \ref{lem:equi-int-u-eps}),   is tight (by Lemma \ref{lem:tigh-u-eps-z}) and converges to $0$ a.e on $\Omega\times\mathbb{R}^d$. Also note that according to the estimate \eqref{eq:estimate-seminorm}, each member of this family is integrable (this  follows from the equiintegrability). Thus for fixed $z\in B_1(0)$ the Vitali's convergence theorem gives 
\begin{align*}
\iil_{\Omega\,\R^d} \big| (u_\eps(x-\eps z) - u_\eps(y-\eps z) ) 
-(u(x)-u(y))\big|^p\nu(x-y)  \ \d y \ \d x\overset{\eps \to 0}{\longrightarrow} 0\,.
\end{align*}	
That is, $|u^z_\eps-u|^p_{\WnuOmR} \to 0, $ as $\eps \to 0$ for each $ z\in B_1(0)$.  Further, from estimate \eqref{eq:estimate-seminorm} the function $ z \mapsto  \phi(z)[u^z_\eps-u]^p_{\WnuOmR} $ is bounded by $2C [u]_{\WnuOmR}$ for all $\eps>0$ and a.e. $z\in B_1(0)$. Thus, by Lebesgue's 
dominated convergence theorem  
\[ \int\limits_{B_1(0)}  |u^z_\eps-u|^p_{\WnuOmR}\phi(z) \d z \overset{\eps \to 0}{ \longrightarrow 0}.  
\]
This implies  $ [v_\eps-u]_{\WnuOmR}\to 0$ as $\eps \to 0$.


\medskip

\noindent \textbf{Step 3:} 	Let $u\in \WnuOmR$ be arbitrary. Let $R>0$ such that $\Omega\subset B_{R}(0)$. Let $f_R\in C_c^\infty(B_{3R}(0))$ with $f_R\leq1$ and 
$f_R(x)=1$ for all $x\in B_{2 R}(0)$. 
Define $u_R = f_R u$. Then according to Lemma \ref{lem:truncate-approx}, we have $\supp(u_R)\subset B_{3R}(0)$, $u_R\in \WnuOm$ and $[u-u_R]_{\WnuOmR} \to 0$ as $R \to \infty$. 

\medskip

\noindent \textbf{Step 4:} Let $x_i\in\partial\Omega$, $r_i>0$, $i=1,..,N$, such that
\[
\partial\Omega \subset \bigcup_{i=1}^N B_{r_i/2}(x_i), 
\]
where the $r_i$ are chosen small enough, such that up to relabelling the coordinates, we can assume 
\begin{align*}
\Omega\cap B_{4r_i}(x_i) &= \{ x\in B_{4r_i}(x_i)| x_d > \gamma_i(x')\}\\
\end{align*}
for some smooth $\gamma_i:\R^{d-1}\to\R$ as in Step 1. Let $\Omega^*= \{x\in \R^d| \dist(x,\Omega)>\frac12 \min_{i=\{1,..,N\}} r_i\}$ 
and $\Omega_0 = \{x\in \Omega| \dist(x,\Omega^c)>\frac12 \min_{i=\{1,..,N\}} r_i\}$. Then 
\[
\bigcup_{i=1}^N B_{r_i}(x_i) \cup \Omega^*\cup\Omega_0 = \R^d .
\]
Let $\{\xi_i\}_{i=0}^{N+1}$ be a smooth partition of unity subordinated to the above constructed sets. 
We define 
\[ u_i = \xi_i\cdot u_R \text{ for all } i\in\{0,..,N+1\}, \] 
and thus 
\begin{align*}
&\supp u_i \subset  B_{r_i}(x_i) \text{ for }i\in\{1,..N\}, \\
&\supp u_0 \subset \Omega_0, \\
&\supp u_{N+1} \subset \Omega^*. 
\end{align*}

\medskip

\noindent \textbf{Step 5:} Let $\delta>0$ and  $i\in\{1,..,N\}$. By Step 2 there exists a sequence $v^i_\eps\in C_c^\infty(B_{4r_i}(x_i))$ such that 
\[
[u_i-v^i_\eps]_{\WnuOmR} \longrightarrow 0 
\]
for $\eps\to 0$. Thus we can choose $\eps_0>0$ such that $[u_i-v^i_\eps]_{\WnuOmR}< \frac{\delta}{N+2}$
for all $i\in\{1,..,N\}$.

For $i=N+1$ define $v^{N+1}_\eps = \phi_\eps \ast u_{N+1}$ and set $r=\frac14 \min_{i\in\{1,..,N\}} r_i$. 
Choosing $\eps< r$ and since  
$\supp u_{N+1} \subset \Omega^*$ for all $x\in\Omega$, $y\in\R^d$ and $z\in B_\eps(0)$ 
\[
U_{N+1}(x,y) = U_{N+1}(x-z, y-z)= 0  \quad \text{ or } \quad \bet{x-y}>r.
\]
where we set $U_{N+1} (x, y) = (u(x)-u(y)\nu^{1/p}(x-y)$. 
Thus, following the estimate \eqref{eq:estimate-conv} , we have
\begin{align*}
|v^{N+1}_\eps-u_{N+1}|^p_{\WnuOmR} &=  |\phi_\eps*u^{N+1}_\eps-u_{N+1}|^p_{\WnuOmR}\\
&\leq\int_{B_1(0)} \phi(z)\d z  \iint\limits_{\Omega\R^d} | U_{N+1}(x-\eps z,y-\eps z)-U_{N+1}(x,y) |^p \d y \d x.  
\end{align*}
By the continuity of the shift in $L^p(\R^d\times \R^d)$,
\[ \iil_{\Omega\,\R^d} |U_{N+1}(x-\eps z, y-\eps z)-U_{N+1}(x,y)|^p \d y \d x \longrightarrow 0. \]
Further,  since $U_{N+1}(x,y) = U_{N+1}(x-z, y-z)= 0 $ or $ \bet{x-y}>r,$  for any $z\in B_1(0)$, then the map
\[
z \mapsto \phi(z)  \iint\limits_{\Omega\,\R^d}  U_{N+1}(x-\eps z,y-\eps z)-U_{N+1}(x,y) |^p 
\nu(x-y) \d y \d x 
\]
is bounded. Thus $[v^{N+1}_\eps-u_{N+1}]_{\WnuOmR}\to 0$ by dominated convergence and we find
$\eps_0>0$, such that  $|v^{N+1}_\eps-u_{N+1}|_{\WnuOmR} < \frac{\delta}{N+2}$ for all $\eps<\eps_0$.
We define $v^{0}_\eps = \phi_\eps \ast u_{0}$. Thus for $\eps<r$ 
\[\supp v^0_\eps\Subset \Omega.\]
The convergence $v^0_\eps \to u_0$ follows by the same arguments as above and we find $\eps_0>0$ such 
that $[v_\eps^0-u_0]_{\WnuOmR} < \frac{\delta}{N+2}$ for all $\eps<\eps_0$. 

\medskip

\noindent \textbf{Step 6:} Define $v_\eps = \sum\limits_{i=0}^{N+1} v^i_\eps\in C^\infty_c(\R^d)$. Since $u_R(x) = \sum\limits_{i=0}^{N+1} u_i(x)$, we have 
\begin{align*}
|u_R-v_\eps|_{\WnuOmR} &\leq \Big|\sum_{i=0}^{N+1} \big(v^i_\eps - u_i\big)\Big|_{\WnuOmR} \\
& \leq \sum_{i=0}^{N+1} |v^i_\eps -u_i|_{\WnuOmR} \leq (N+2) \frac{\delta}{N+2} .
\end{align*}
Choosing $\frac{1}{\eps}<R$ in Step 3 such that $|u-u_{R}|_{\WnuOmR}<\delta$, concludes 
\[|u-v_\eps|_{\WnuOmR}\leq |u-u_{R}|_{\WnuOmR}+ |u_R-v_\eps|_{\WnuOmR}\leq 2\delta. \]
The convergence in $L^p(\Omega)$ follows from the continuity of the shift in $L^p(\R^d)$.  
\end{proof}

\bigskip

\noindent We now derive another variant of Theorem \ref{thm:density} for $\WnuOm\cap  L^p(\R^d)$ only assuming \eqref{eq:plevy-integrability}. Recall $$\vertiii{u}^p_{\WnuOmR} =\|u\|^p_{L^p(\R^d)} + \|u\|^p_{\WnuOmR}.$$
\begin{theorem}\label{thm:density-bis}
	Assume $\Omega\subset \mathbb{R}^d$ is open with a compact Lipschitz boundary $\partial \Omega$. Assume that $\nu$ only satisfies \eqref{eq:plevy-integrability}. Then, $C_c^\infty(\mathbb{R}^d)$ is dense in $\WnuOmR\cap  L^p(\R^d)$ with respect to the norm $\vertiii{\cdot}_{\WnuOmR}$, i.e. for $u \in \WnuOmR$,  there exists a sequence $(u_n) \subset C_c^\infty(\mathbb{R}^d)$ with 
	\begin{align*}
 \vertiii{u_n-u}_{\WnuOmR} \longrightarrow 0 \text{ as } n\to\infty \,.
	\end{align*} 
\end{theorem}

\medskip

\begin{proof} The proof of the assertion here solely follows the scheme of the proof of Theorem \ref{thm:density}  so we only point out the crucial Step 2 where the radially of $\nu$ enters into play. The estimate \eqref{eq:estimate-seminorm}, i.e.
\begin{align*}\label{eq:estimate-seminorm-repeated}
		|u^z_\eps|^p_{\WnuOmR}\leq C(\Omega,r, \nu) | u|^p_{\WnuOmR}
		\end{align*} 
remains true with the constant 
\[ C= 1+ \hspace{-3ex}\sup_{y\in B_{r/4}(x_0)}\Big(\hspace{-1ex}
\il_{\Omega\cap B_r^c(x_0)} \hspace{-2ex}\nu(x-y)\d x\Big)^{-1} \hspace{-2ex}\int\limits_{ B^c_{r/4}(0)} \hspace{-2ex}\nu(h)\d h<\infty.\]
We recall, $|x-y|\geq r/4$ for $x\in B^c_{r/2}(x_0) $ and $y \in B_{r/4}(x_0)$.  

\noindent However since in this case we have $u\in L^p(\R^d)$ it is easy to obtain the following analog estimate
\begin{align}\label{eq:estimate-seminorm-bis}
|u^z_\eps|^p_{\WnuOmR}\leq C(r, \nu) \| u\|^p_{\WnuOmR},
\end{align} 
with the mere constant $ C= 1+\int\limits_{ B^c_{r/4}(0)} \hspace{-2ex}\nu(h)\d h<\infty.$
Indeed, a close look at the proof of the estimate \eqref{eq:estimate-seminorm-bis} shows that one gets the estimate \eqref{eq:estimate-seminorm-bis} by applying the following
\begin{alignat*}{1}
&\iint\limits_{ \Omega^z_\epsilon\cap B^c_{r/2}(x_0)\times \R^d} \left|u(x)|^p - u(y) \right|^p \nu(x-y)  \d y \d x\\
&= \int\limits_{ B_{r/4}(x_0) } |u(y)|^p\d y \int\limits_{ \Omega^z_\epsilon\cap B^c_{r/2}(x_0)}  \nu(x-y)   \d x\leq  C \int_{\R^d} |u(y)|^p\d y.
\end{alignat*}
Next by mimicking the estimate \eqref{eq:estimate-conv} we have 
\begin{align*}
|v_\eps -u|^p_{\WnuOm} 
&\leq  \il_{ B_1} \phi(z)\d z  \iil_{\Omega\R^d} \big| U(x-\eps(\tau e_d-z),y- \eps(\tau e_d-z)) -U(x,y)\big|^p \d y \d x
\end{align*}
with $U(x,y)= (u(x)-u(y))\nu^{1/p}(x-y)$. The estimate \eqref{eq:estimate-seminorm-bis} shows that $u, u^z_\eps\in \WnuOmR$ we recall $u^z_\eps(x) =u(x+\eps(\tau e_d-z))$ equivalently  $ U(\cdot-\eps(\tau e_d-z), \cdot-\eps(\tau e_d-z)), U\in L^p(\Omega\times \R^d)$. So that by continuity of the shift, 
\begin{align*}
 \iil_{\Omega\R^d} \big| U(x-\eps(\tau e_d-z),y- \eps(\tau e_d-z)) -U(x,y)\big|^p \d y \d x\xrightarrow{\varepsilon\to 0}0.
\end{align*}
 Analogously by invoking dominated convergence we get  $|v_\eps -u|^p_{\WnuOm} \xrightarrow{\varepsilon\to 0}0$ and step 2 follows. The remaining details from the proof of Theorem \ref{thm:density} are unchanged. However we emphasize that the smoothness of $v_\eps= \phi_\eps* u_\eps$  is provided since $u\in L^p(\R^d)\subset L^p_{\operatorname{loc}}(\R^d)$. 
\end{proof}

\medskip

\begin{theorem}\label{thm:density-zero-outside}
Let $\Omega$ be a continuous domain such that the boundary $\partial \Omega$ is compact. Assume $\nu$ satisfies 
\eqref{eq:plevy-integrability} (see page  \pageref{eq:plevy-integrability}).  The following assertions hold.
\begin{enumerate}[$(i)$]
\item $C_c^\infty(\Omega)$ is dense in $\WnuOmO$. 
\item  $C_c^\infty(\overline{\Omega})$ is dense in $\WnuOm$, where  $C_c^\infty(\overline{\Omega})$ is the space of all functions which are restriction of $C^\infty$ functions with compact support in $\mathbb{R}^d$  to $\Omega$. 
\end{enumerate}
\end{theorem}
\bigskip

\begin{remark}
Note that an incomplete proof of Theorem \ref{thm:density-zero-outside} $(i)$ is provided in \cite[Theorem A.4]{BGPR17} under the additional condition that $\nu$ radially almost decreasing. However for our setting, this is just reminiscent of the main result in \cite{FSV15} for fractional Sobolev spaces from which we borrow some ideas. \linebreak[-2]
\end{remark}

\begin{proof} $(i)$ Since  $\partial\Omega$ is compact, let $x_i\in\partial\Omega$, $i=1,..,N$ and $r>0$ such that
\[
\partial\Omega \subset \bigcup_{i=1}^N B_{r/2}(x_i), 
\]
where the $r>0$ is chosen small enough, such  that up to relabelling the coordinates, we can assume 
\begin{align*}
\Omega\cap B_{r}(x_i) &= \{ x=(x',x_d)\in B_{r}(x_i)| x_d > \gamma_i(x')\}\\
\Omega^c\cap B_{r}(x_i) &= \{ x=(x',x_d)\in B_{r}(x_i)| x_d \leq \gamma_i(x')\}.
\end{align*}
for some continuous functions $\gamma_i:\R^{d-1}\to\R$. Let $\Omega_{r/2}= \{x\in \Omega: \dist(x,\partial\Omega)>r/2\}$. Then, 
\[
\Omega\subset \bigcup_{i=1}^N B_{r}(x_i) \cup \Omega_r .
\]
Let $\{\xi_i\}_{i=0}^{N}$ be a smooth partition of unity subordinated to the above constructed sets. That is we have $\xi_i\geq 0$, $\sum\limits_{i=0}^N\xi_i=1$ and $\xi_0\in C_c^\infty(\Omega_{r/2})$ and $ \xi_i\in C_c^\infty(B_{r}(x_i)).$ Let  $u\in \WnuOmO$ and define 
\[ u_i = \xi_i\cdot u\quad \text{ for all } i\in\{0,..,N\}. \] 
%
For $0<\eta<r/2$ sufficiently small, we define
\begin{align*}
	u^i_\eta (x) = u^i(x',x_d-\eta) \quad\text{ for all } i\in\{0,..,N\}.
\end{align*}
Recall that $u=0 $ on $\Omega^c$, so  it is clear that $u\in W^p_\nu(\R^d)$ and according to Lemma \ref{lem:cutt-off} $u^i=\xi_i u \in W^p_\nu(\R^d)$. Therefore, regarding the finiteness of the $\|u^i\|_{W^p_\nu(\R^d)}$, no problem should raise while shifting the function $u^i$.  Thus $u^i, u^i_\eta  \in W^p_\nu(\R^d)$. Invoking the continuity of the shift in both $L^p(\R^d)$ and $L^p(\R^d\times \R^d)$, for instance like in the proof of Theorem \ref{thm:approx-smooth}, we find that $\|u^i-u^i_\eta\|_{W^p_\nu(\R^d)}\xrightarrow{\eta\to 0}0$. Given $\eps>0$, from now  we fix $\eta>0$ sufficiently small that 
\begin{align*}
\|u^i-u^i_\eta\|_{W^p_\nu(\R^d)}\leq \frac{\eps}{2(N+1)}. 
\end{align*}
Meanwhile we claim that $u^i_\eta$ is compactly supported in $\Omega$. We even have $\supp u^i_\eta\subset \overline{V^i_\eta} $ with  $V^i_\eta= \{x\in B_r(x_i):\, x_d-\eta>\gamma_i(x')\}.$ Indeed assume $$x= (x',x_d)\not \in V^i_{\eta}= \{ z= (z',z_d)\in B_r(x_i):\, _d-\eta>\gamma_i(z')\}$$ 
equivalently 
$$(x',x_d-\eta)\in \{z= (z',z_d)\in B_r(x_i):\, z_d\leq \gamma_i(z')\} = \Omega^c\cap B_r(x_i).$$
Which implies $u^i_\eta(x) = [\xi_iu](x',x_d-\eta) =0$ because  $(x',x_d-\eta) \in \Omega^c$ and $u^i= \xi_i u=0$ on $\Omega^c$. We now get
$$\supp u^i_\eta\subset \overline{V^i_\eta} \subset\Omega\cap B_r(x_i).$$
 Indeed, by the continuity of the $\gamma_i's$ we get $$ \overline{V^i_\eta} = \{x\in B_r(x_i):\, x_d-\eta\geq \gamma_i(x') \}\subset  \{x\in B_r(x_i):\, x_d>\gamma_i(x') \}= \Omega \cap B_r(x_i).$$ 
 Meanwhile the same arguments used in the proof of  Theorem \ref{thm:approx-smooth} yields that 
$\|\phi_\delta *u^i_\eta- u^i_\eta\|_{W^p_\nu(\R^d)}\xrightarrow{\delta\to 0}0$. 
For $0<\delta <\dist(\supp u^i_\eta, \partial \Omega)/2$, we have 
\begin{align*}
	\supp \phi_\delta*u^i_\eta\subset \overline{B_\delta+  \supp u^i_\eta} \subset \Omega. 
\end{align*}
We have $\phi_\delta* u^i_\eta \in C_c^\infty(\Omega)$.  We can find $0<\delta <\dist(\supp u^i_\eta, \partial \Omega)/2$ such that $\phi_\delta* u^i_\eta \in C_c^\infty(\Omega)$ and 
\begin{align*}
\|\phi_\delta *u^i_\eta- u^i_\eta\|_{W^p_\nu(\R^d)}\leq \frac{\eps}{2(N+1)}.
\end{align*}
For $u_0$ the assertion $\phi_\delta*u_0\in C_c^\infty(\Omega)$, is much easier since there is no need to shift it. Of course $v= \sum\limits_{i=0}^N \phi_\delta* u^i_\eta \in C_c^\infty(\Omega)$ and $u= \sum\limits_{i=0}^N u^i $ since $\sum\limits_{i=0}^N \xi_i=1$. 
Altogether we get 
\begin{align*}
\|v-u\|_{W^p_\nu(\R^d)}& =\Big\|\sum_{i=0}^N\phi_\delta *u^i_\eta- u^i\Big\|_{W^p_\nu(\R^d)} \\
&\leq\sum_{i=0}^N \big[\|\phi_\delta *u^i_\eta- u^i_\eta\|_{W^p_\nu(\R^d)}+ \|u^i_\eta- u^i\|_{W^p_\nu(\R^d)} \big]\\
&\leq  \sum_{i=0}^N\frac{\eps}{N+1}= \eps.
\end{align*}

\noindent To show $(ii)$, in view of Theorem \ref{thm:meyer-serrin} it suffices to prove  it for $u\in C^\infty(\Omega)\cap \WnuOm$. With the previous notations we consider $u^i_\eta(x) =[u\xi_i](x', x_d+\eta)$, 
$i=1, \cdots, N$. Note that $u^i_\eta$ is $C^\infty$ and supported in $B_r(x_i)$ because $\supp \xi_i\subset B_r(x_i)$. Furthermore,  for $0<\eta <r$ small enough, it is possible to show that $x\in \Omega^c\cap B_{\eta/2}(x_i) $ implies $(x', x_d+ \eta)\in \Omega\cap B_r(x_i)$ and hence $u^i_\eta(x', x_d+\eta) $ is well define. 
In other words, $u^i_\eta\in C_c^\infty(\overline{\Omega})$. Also note that for $\eta>0$ small enough we also  
have $u^i_\eta\in \WnuOm$. With analogous arguments as previously, for $v= \sum\limits_{i=0}^N \phi_\delta* u^i_\eta \in C_c^\infty(\overline{\Omega})$ one gets $\|v-u\|<\eps$.
 
%
\end{proof}



%
%
	\noindent Now, we state an extension result for the space $\WnuOm$.  Up to a rigorous modification of the extension result of \cite[Theorem 5.4]{Hitchhiker} one
	 is able to obtain the following extension theorem. 
\begin{theorem}
Let $1\leq p<\infty$ and $\Omega$ be open with a bounded Lipschitz boundary. Let $\nu: \R^d\setminus\{0\}\to [0,\infty]$ be radial, almost decreasing and satisfies the $p$-L\'evy integrability \eqref{eq:plevy-integrability}. Then there exists linear bounded operator $E:\WnuOm\to W_{\nu}^{p}(\R^d)$  such that $Eu|_\Omega = u$ for all $u\in \WnuOm$. 
\end{theorem}
 
\noindent If the conclusion of the above theorem  holds then one says that $\Omega$ is an extension domain for $\WnuOm$.  Let us recall, in the special case where $\nu(h)= |h|^{-d-sp}$,  $s\in (0,1)$, the  main result of  \cite{Zh15} states that $\Omega$ is an extension domain  for $W^{s,p}(\Omega)$ if and only if $\Omega$ is a $d$-set, i.e. there exists $c>0$ such that for every $x\in \partial\Omega$ and $r\in (0,1)$ we have $|B(x,r)\cap \Omega|\geq cr^d$. 

\section{Compact embeddings and Poincar\'e type inequalities}\label{sec:compactness}

\noindent In this section we are concerned with the compact embeddings of the spaces $\WnuOm$, $\WnuOmR$ and $\WnuOmO$ into $L^p(\Omega)$. Let us start with some basic observations and formulate some sufficient assumptions on $\nu$ and $\Omega$. We shall temporarily modify our general assumption on the function $\nu: \mathbb{R}^d\setminus\{0\}\to [0, \infty]$. Like previously for $1\leq p<\infty$, we assume that $\nu$ satisfies
\begin{align}\tag{$I_1$}\label{eq:integrability-condition+even}
\nu(-h) = \nu(h) ~\text{ for all $h\in \mathbb{R}^d$ and} ~\int_{\mathbb{R}^d}(1\land |h|^p)\nu(h)\d h<\infty\,.
\end{align}

\noindent It is an obvious fact that if $\nu\in L^1(\mathbb{R}^d)$, then the space $\WnuOm$ coincides with $L^p(\Omega)$ and thus cannot be compactly embedded into $L^p(\Omega)$. Likewise, if $\nu\in L^1(\mathbb{R}^d)$, then the spaces $\WnuOmR \cap L^p(\mathbb{R}^d)$ and $W^p_{\nu}(\mathbb{R}^d)$ coincide with $L^p(\mathbb{R}^d)$ which is not even locally compactly embedded in $L^p(\Omega)$. In other words, the least necessary condition for compact embeddings to hold is that $\nu$ must not be integrable. Therefore, it is necessary to consider the following non-integrability condition 
\begin{align}\tag{$I_2$}\label{eq:non-integrability-condition}
\int_{\mathbb{R}^d}\nu(h)\d h= \infty.
\end{align}
If the condition \eqref{eq:integrability-condition+even} holds true, it is possible to strengthen the condition \eqref{eq:non-integrability-condition} by merely assuming  that 
\begin{align}\tag{$I'_2$}\label{eq:limit-at-0-explode}
\lim_{|h|\to 0}|h|^d \nu(h)= \infty. 
\end{align}
\noindent The conditions \eqref{eq:integrability-condition+even} and \eqref{eq:limit-at-0-explode} clearly imply \eqref{eq:integrability-condition+even} and \eqref{eq:non-integrability-condition}. The latter conditions are sufficient to obtain locally compact embeddings of $\WnuOm$ and $\WnuOmR$ into $L^p(\Omega)$ 
(see Corollary \ref{cor:local-compatcness}). The global compactness needs some extra regularity assumptions on $\Omega$ compatible with $\nu$ that we will state later on. First and foremost, we the following well known from functional analysis. 
%
\begin{theorem}[Chapter X, p.278, \cite{yosida80}]  \label{thm:closeness-compact-oper}Given $X$ and $Y$ two Banach spaces denote by $\mathcal{L}(X,Y)$ (resp. $\mathcal{K}(X,Y)$) the space of linear bounded operators( resp. linear compact operators) from $X$ into $Y$. Then $\mathcal{K}(X,Y)$ is closed in $\mathcal{L}(X,Y)) $ with respect top the topology associated to the norm $	\|\cdot\|_{\mathcal{L}(X,Y)) }$.
	\begin{align*}
	\|T\|_{\mathcal{L}(X,Y)) }:= \sup_{\|x\|_X=1}\|Tx\|_Y	.
	\end{align*} 

\end{theorem}

\noindent For a measurable subset $D\subset\mathbb{R}^d$ we adopt the notation $R_D$ to denote the restriction operator assigning $u|_D$ to a function $u$. The following lemma is a consequence of Theorem \ref{thm:compactness-convolution} applied with $p=r$ and $q=1$. 

\begin{lemma} \label{lem:compactness-convolution} 
	Let $w \in L^1(\mathbb{R}^d)$. Then the convolution operator $T_w : L^p(\mathbb{R}^d )\to L^p(\mathbb{R}^d)$ with $T_wu=w*u$ is continuous, locally compact and its norm is not greater than $\|w\|_{L^1(\mathbb{R}^d)}$. 
\end{lemma}

\medskip

\noindent In what follows, we denote by $\nu_\delta$ the kernel $\nu$ truncated away from the ball $B_\delta(0)$, i.e.,  for every $h\in \mathbb{R}^d $ and every $\delta>0$, we have 
$$\nu_\delta = \mathbbm{1}_{\mathbb{R}^d\setminus B_\delta(0)}(h) \nu(h).$$ 
The following result is reminiscent of \cite[Theorem 1.2]{JW19} under a slightly weaker assumption. Recall that $u\in \WnuOmO$  if and only if $u\in \WnuOmR$ and $u= 0$ on $\Omega^c$. 
\begin{theorem}\label{thm:local-compactness}
	Let $\nu:\mathbb{R}^d\setminus\{0\}\to [0, \infty]$ be a measurable function for which the conditions \eqref{eq:non-integrability-condition} and \eqref{eq:integrability-condition-near-zero} hold (in particular if \eqref{eq:integrability-condition+even} holds true ) with 
	\begin{align}\tag{$I'_1$}\label{eq:integrability-condition-near-zero}
	\nu(-h) = \nu(h) ~\text{ for all $h\in \mathbb{R}^d$ and} ~\int_{|h|\geq \delta}\nu(h)\d h<\infty~\text{ for all $\delta>0$} .
	\end{align}
	\noindent Then the embedding $W^p_\nu(\mathbb{R}^d) \hookrightarrow L^p(\mathbb{R}^d)$ is locally compact. Moreover, for $\Omega\subset \mathbb{R}^d$ open and bounded, the embedding $\WnuOmO \hookrightarrow L^p(\Omega)$ is compact. 
\end{theorem}

\begin{proof} This proof is in the spirit of \cite{JW19}. Let $\delta>0$ small enough, the assumptions \eqref{eq:integrability-condition-near-zero} and \eqref{eq:non-integrability-condition} imply $0<\|\nu_\delta\|_{L^1(\mathbb{R}^d)}<\infty$; we set $w_\delta =\frac{\nu_\delta}{\|\nu_\delta\|_{L^1(\mathbb{R}^d)}}$. Thus, $\|w_\delta\|_{L^1(\mathbb{R}^d)} =1$ and for fixed $u \in L^p(\mathbb{R}^d)$, by evenness of $\nu $ for all $x\in \mathbb{R}^d$ we have 
	\begin{align*}
	T_{w_\delta} u(x)= \int_{\mathbb{R}^d} w_\delta(y) u(x-y)\d y = \int_{\mathbb{R}^d} w_\delta(y) u(x+y)\d y, \qquad \,. 
	\end{align*} 
	Thus, by Jensen's inequality
	\begin{align*}
	\|u- T_{w_\delta} u\|^p_{L^p(\mathbb{R}^d)} &= \int_{\mathbb{R}^d} \Big| \int_{\mathbb{R}^d} [u(x)-u(x+h) ]w_\delta (h)\d\, h\Big|^p\d x\\
	&\leq \iint\limits_{\mathbb{R}^d\mathbb{R}^d} |u(x)-u(x+h) |^p w_\delta (h)\d\, h\d x\\
	& \leq \|\nu_\delta\|^{-1}_{L^1(\mathbb{R}^d)} \iint\limits_{\mathbb{R}^d\mathbb{R}^d} |u(x)-u(x+h) |^p \nu (h)\d\, h\d x\\
	& \leq \|\nu_\delta\|^{-1}_{L^1(\mathbb{R}^d)} \|u\|^p_{W^p_\nu(\mathbb{R}^d)} \,.
	\end{align*}
	
	\noindent Accordingly, for a compact subset $K$ of $\mathbb{R}^d$, taking into account the assumption \eqref{eq:non-integrability-condition} leads to
	\begin{align*}
	\|R_K -R_K T_{w_\delta} \|_{\mathcal{L}\big(W^p_\nu(\mathbb{R}^d),\, L^p(K)\big)}\leq \|\nu_\delta\|^{-1}_{L^1(\mathbb{R}^d)}\xrightarrow[]{\delta\to 0}0\,.
	\end{align*}
	
	\noindent Thus the embedding $R_K:~ W^p_\nu(\mathbb{R}^d)\to L^p(K)$ is compact since by Lemma \ref{lem:compactness-convolution}, the operator $R_K T_{w_\delta}$ is also compact for every $\delta>0$. In view of Theorem \ref{thm:closeness-compact-oper} $\mathcal{K}\big(W^p_\nu(\mathbb{R}^d),\, L^p(K)\big)$ is closed in $\mathcal{L}\big(W^p_\nu(\mathbb{R}^d),\, L^p(K)\big)$. This prove the locally compactness of the embedding $W^p_{\nu}(\mathbb{R}^d) \hookrightarrow L^p(\mathbb{R}^d)$. Furthermore, it springs directly from the continuous embeddings, $\WnuOmO \hookrightarrow W^p_{\nu}(\mathbb{R}^d) \hookrightarrow L^p(\Omega)$ and the ideal property of compactor operators that the embedding $\WnuOmO\hookrightarrow L^p(\Omega)$ is also compact. 
\end{proof}

\begin{remark}
	Note that \eqref{eq:integrability-condition-near-zero} and \eqref{eq:non-integrability-condition} do not only capture the class of L\'{e}vy integrable functions which are non-integrable but also functions with strong singularity at the origin. For instance, consider $\nu(h)= |h|^{-d-\beta}$ for $h \neq0$ with $\beta>0$. In this case, the  $p$-L\'{e}vy integrability in \eqref{eq:integrability-condition+even} fails for $\beta\geq p$ while \eqref{eq:integrability-condition-near-zero} and \eqref{eq:non-integrability-condition} remain true and hence Theorem \ref{thm:local-compactness} applies.
\end{remark}

\medskip

	\noindent In the sequel we will frequently use the family of cut-off functions introduced ion the following lemma. 
\begin{lemma}\label{lem:cut-off-existence}
	For $\delta>0$ small enough,  recall that $\Omega_\delta= \{x\in\Omega:\, \dist(x,\partial\Omega)>\delta\}$.
	Let $\phi \in C_c^\infty(\mathbb{R}^d)$ supported in the unit ball $B_1(0)$, $\phi\geq 0$ and $\int_{\mathbb{R}^d} \phi(x)\d x =1$. 
	Then the  function $\varphi_\delta(x) =\phi_{\delta/8}*\mathds{1}_{\Omega_{5\delta/8}}(x) $ with $ \phi_{\delta}(x)=\frac{1}{ \delta^d}\phi(\frac{x}{\delta})$ satisfies: $\varphi_\delta\in C_c^\infty(\Omega),\, \supp \varphi_\delta\subset \Omega_{\delta/2},\, $ $0\leq \varphi_\delta\leq 1, \, \varphi_\delta =1 \,\text{ on $\Omega_\delta$},\, $ $\varphi_\delta\xrightarrow{\delta\to 0}1$ and $|\nabla \varphi_{\delta} |\leq c/\delta, (\mbox{with $c =8\int|\nabla \phi|$}).$
	%
	
\end{lemma}

\medskip

\begin{proof}
	Indeed, since $\phi_{\delta/8}$ is supported in $B_{\delta/8}(0)$,
	\begin{align*}
	\varphi_{\delta}(x)= \int\limits_{B_{\delta/8}(0)} \phi_{\delta/8} (y)\mathds{1}_{\Omega_{5\delta/8}}(x-y)\,\d y= \int\limits_{\Omega_{5\delta/8}} \phi_{\delta/8} (x-y)\,\d y.
	\end{align*}
	
	\noindent Fix $x\in \Omega_\delta$. Let $z\in \partial\Omega$ and $y \in B_{\delta/8}(0)$ then $|x-z-y|\geq |x-z|-|y|\geq \frac{3\delta}{4}$. Since $z\in \partial\Omega$ is arbitrarily chosen, $\dist(x-y, \partial \Omega)\geq \frac{3\delta}{4}>\frac{5\delta}{8}$  we get $x-y\in \Omega_{5\delta/8}$ i.e $ \mathds{1}_{\Omega_{5\delta/8}}(x-y)=1$ which implies 
	$$\varphi_{\delta}(x)= \int\limits_{B_{\delta/8}(0)} \phi_{\delta/8} (y)\,\d y=1.$$ 
	That is $\varphi_\delta =1,$ on $\Omega_\delta$. Hence we also get $\varphi_\delta\xrightarrow{\delta\to 0}1$. Now if $x\in \Omega\setminus \Omega_{\delta/2}$ then for $y \in \Omega_{5\delta/8}$ we have 
	\vspace{-2ex}
	
	\begin{align*}
	|x-y|\geq |\operatorname{dist}(y,\partial\Omega)-\operatorname{dist}(x,\partial\Omega)| \geq \delta/8\,
	\end{align*}
	so that $x-y\not\in B_{\delta/8}(0)$, i.e. $\phi_{\delta/8}(x-y) =0$ for all $y \in \Omega_{5\delta/8}$ which means that $\varphi_\delta(x)=0$ that is $\operatorname{supp }\varphi_\delta\subset \Omega_{\delta/2}$. Last we note that $\nabla \varphi_\delta(x) =\frac{8}{\delta}[\nabla\phi]_{\delta/8}*\mathds{1}_{\Omega_{5\delta/8}}(x) $ so that $|\nabla \varphi_{\delta} |\leq c/\delta$, (\text{with $c= 8\int |\nabla \phi|$}). \\
\end{proof}

\medskip

\noindent As an immediate consequence of Theorem \ref{thm:local-compactness} we get  the local compactness of $\WnuOm$ in $L^p(\Omega)$. 

\begin{corollary}\label{cor:local-compatcness}
	Let $\Omega\subset\mathbb{R}^d$ be open bounded. Assume that $\nu:\mathbb{R}^d\setminus\{0\} \to \mathbb{R}$ fulfills conditions \eqref{eq:integrability-condition+even} and \eqref{eq:non-integrability-condition}. The embedding $\WnuOm\hookrightarrow L^p(\Omega) $ is locally compact. Furthermore, for every bounded sequence $(u_n)_n$ there exits $u\in \WnuOm$ and subsequence $(u_{n_j})_j$ converging to $u$ in $L^p_{\operatorname{loc}}(\Omega)$. Moreover, 
		\begin{align*}
	\|u\|_{\WnuOm} \leq \liminf_{n\to \infty}\|u_n\|_{\WnuOm}.
	\end{align*}
\end{corollary}

\vspace{2mm}

\begin{proof}
	In view of Lemma \ref{lem:cutt-off} we claim that for $\varphi \in C_c^\infty(\Omega)$, the mapping $J_\varphi: \WnuOm\to \WnuOmO$, with $J_\varphi u = u\varphi $ is continuous and hence by Theorem \ref{thm:local-compactness}  and the ideal property of compact operator, the linear map $ J_\varphi: \WnuOm\to  \WnuOmO\to L^p(\Omega)$ is compact which therefore implies the locally compactness of the embedding $\WnuOm\hookrightarrow L^p(\Omega)$. 
	
	%
	%
	
	\vspace{2mm}
	
	\noindent Let us prove the second statement. 
	For $\delta>0$ small enough, we let $\varphi_\delta \in C^\infty(\Omega)$ be  such that $\varphi_\delta =1$ on $\Omega_\delta$ see in Lemma \ref{lem:cut-off-existence}. 
	
	\vspace{2mm}
	
	\noindent Next we employ Cantor's diagonalization procedure to show that every bounded sequence of $\WnuOm$ has a converging subsequence in $L_{\operatorname{loc}}^2(\Omega)$ to some function $u\in \WnuOm$. 
	To this end, for $k\in \mathbb{N}$, we let $\delta_k = \frac{1}{2^{k+k_0}}$ and merely set the shorthand notations $\Omega_k = \Omega_{\delta_k}$ and $\varphi_k = \varphi_{\delta_k}$ where $k_0$ is large enough so that $ \Omega_{\delta_0}$ is non-empty. For every $k\geq 1$, $\overline{\Omega}_k \subset \Omega_{k+1}$ and $\varphi_k=1$ on $\overline{\Omega}_k$. 
	
	\noindent Assume $(u_n)_n$ is a bounded sequence in $\WnuOm$. By the above remark, for each $k\geq 1$ the sequence $(\varphi_k u_n)_n$ is also bounded in $W^p_\nu(\mathbb{R}^d)$ and hence relatively compact in $L^p_{\operatorname{loc}}(\mathbb{R}^d)$ by Theorem \ref{thm:local-compactness}. In particular for fixed $k,$ $(u_n)_n$ is relatively compact in $L^p(\Omega_k)$ because $\varphi_k u_n= u_n$ on $\Omega_k$. Thus there is a subsequence $(u_{\theta_k(n)})_n$ of $(u_n)_n$ converging to some function $u_k$ in $L^p(\Omega_k)$ and almost everywhere in $\Omega_k$. Applying this argument again on $(\varphi_{k+1}u_{\theta_k(n)})_n$ it turns out that the subsequence $(u_{\theta_k(n)})_n$ possesses a further subsequence $(u_{\theta_{k+1}(n)})_n$ converging to some function $u_{k+1}$ in $L^p(\Omega_{k+1})$ and a.e. on $\Omega_{k+1}$. 
	
	\vspace{2mm}
	
	\noindent By this procedure, assume that for each $k$ we have constructed a subsequence $(u_{\theta_k(n)})_n$ having the property that the $(k+1)^{th}$ subsequence $(u_{\theta_{k+1}(n)})_n$ is a subsequence of the preceding $k^{th}$ subsequence $(u_{\theta_{k}(n)})_n$ and additionally the sequence $(u_{\theta_{k}(n)})_n$ converges in $L^p(\Omega_k)$ and almost everywhere on $\Omega_k$ to some function $u_k$. 
	Clearly, for each $k,$ the restriction of $u_{k+1}$ on $\Omega_k$ coincides with $u_k$. 
	Therefore the function $ u:\Omega\to \mathbb{R}$ coinciding with $u_k$ on $\Omega_k$ for every $k\in \mathbb{N}$ is well defined and measurable since it is easy to verify that $u = u_0+ \sum\limits_{k=1}^{\infty} u_k\mathds{1}_{\Omega_{k}\setminus\Omega_{k-1}}~~a.e$. 
	Now we consider the diagonal sequence $(u_{\theta_{n}(n)})_n$ which is a bounded subsequence in $L^p(\Omega)$ of the original sequence $(u_n)_n$ converging almost everywhere to $u$ on $\Omega$. Indeed this follows immediately since $u_{\theta_{k}(n)} \to u_k~a.e $ on $\Omega_k$ and $(\Omega_k)_k$ is an exhaustion of $\Omega$. We conclude that $u \in \WnuOm$ since $(u_{\theta_{n}(n)})_n$ is bounded in $\WnuOm$, and by Fatou's lemma we have 
	\begin{align*}
	\|u\|_{\WnuOm} \leq \liminf_{n\to \infty}\|u_{\theta_{n}(n)}\|_{\WnuOm}<\infty.
	\end{align*}
\end{proof}

\medskip

\noindent  Meanwhile, the global compactness requires some extra  compatibility assumptions 
between  $\Omega$ and $\nu$.  We establish the global compactness, by exploiting the recent results from \cite{JW19} and \cite{DMT18}. We intend to provide an alternative approach to the compactness result in \cite[Theorem 2.2]{CP18}. The technique therein is adapted from \cite[Theorem 7.1]{Hitchhiker} for fractional Sobolev spaces which uses the Sobolev extension property of the corresponding domain. However, the proof provided in \cite{CP18} seems to be valid only for domains which can be written as a finite union of cubes; unless 
the corresponding nonlocal function space possesses the extension property. Our approach is rather standard and 
follows the idea  used to prove the classical  Rellich–Kondrachov theorem, i.e. the compactness of the embedding $H^1(\Omega)\hookrightarrow L^2(\Omega)$ for $\Omega$  sufficiently smooth. It consists of applying the local local compactness and using an approximation argument near the boundary of $\Omega$. 

\medskip

\noindent Let us  introduce some regimes relating $\Omega$ and $\nu$ under which the global compactness holds true. We will enumerate these assumptions on the couple $(\nu, \Omega)$ into different classes.  We say that  the couple $(\nu, \Omega)$ is in the class $\mathscr{A}_i,~~i=1,2,3$ if $\Omega\subset \mathbb{R}^d$ is an open bounded set, $\nu: \mathbb{R}^d\setminus\{0\}\to [0, \infty]$ satisfies the conditions \eqref{eq:integrability-condition+even-intro} and \eqref{eq:non-integrability-condition-intro} and  additionally $\nu$ and $\Omega$ satisfy:  

\vspace{1mm}

\noindent $\bullet$ The class $\mathbf{\mathscr{A}_1}$: there exists an $\WnuOm$-extension operator $E: \WnuOm\to W^p_{\nu}(\mathbb{R}^d)$. That is there is a constant $C: C(\nu, \Omega, d)>0$ such that  $\|u\|_{W^p_\nu(\mathbb{R}^d)}\leq C\|u\|_{\WnuOm}$ and $Eu|_\Omega =u$ for every $u\in \WnuOm$.

\vspace{1mm}

\noindent $\bullet$ The class $\mathbf{\mathscr{A}_2}$: $\Omega$ has Lipschitz boundary, $\nu$ is radial and 
\begin{align}\label{eq:class-lipschitz}
q(\delta):= \frac{1}{\delta^p}\int_{B_\delta(0)} |h|^p\nu(h)\d h \xrightarrow{\delta \to 0}\infty. 
\end{align}

\vspace{1mm}

\noindent $\bullet$ The class $\mathbf{\mathscr{A}_3}$: letting $\Omega_\delta= \{x\in \Omega: \operatorname{dist}(x,\partial\Omega)>\delta\}$ for $\delta>0$,  the following condition holds true

\begin{align}\label{eq:class-sing-boundary}
\widetilde{q}(\delta): = \inf_{a\in \partial\Omega}\int_{\Omega_\delta}\nu(h-a)\d h \xrightarrow[]{\delta \to 0}\infty.
\end{align}

\medskip
 \noindent Let us also introduce a fourth class $\mathbf{\mathscr{A}_4}$ of interest. 
 
\noindent $\bullet$ The class $\mathbf{\mathscr{A}_4}:$ we say that the couple $(\nu, \Omega)$ is in the class $\mathbf{\mathscr{A}_4}$ if  $\Omega$ is any open bounded subset of $\mathbb{R}^d$ and $\nu:\mathbb{R}^d\setminus\{0\}\to [0, \infty]$ is a unimodal L\'{e}vy measure that is, $\nu$ is radial, almost decreasing and $\nu \in L^1(\R^d, 1\land |h|^2\d h)$. Note that in the class $\mathscr{A}_4,$ $\nu$ is not necessarily singular near $0$.

\vspace*{2mm}


%
%
\noindent Let us make some comments about the newly introduced classes. We  note that the monotonicity of $\nu$ is not required here. The condition \eqref{eq:class-sing-boundary} reflects a certain correlation between $\Omega$ and the singularity of $\nu$ near the origin when shifted on the boundary $\partial \Omega$ of $\Omega$. 
In a sense, the singularity of $\nu$ is compatible with the boundary $\partial \Omega$. On the other hand, $\widetilde{q}(\delta)<\infty$ since for each $a\in \partial\Omega$ and each $\delta>0$, $ \Omega_\delta\subset B^c_\delta(a)$ and hence by \eqref{eq:integrability-condition+even} we get $$ \widetilde{q}(\delta)\leq \int_{ B^c_\delta(0) }\nu(h)\d h<\infty.$$

\noindent Straight away, we would like to show some examples of elements of the classes $\mathscr{A}_i$, $i=1,2,3$. To this end, let us recall some concepts about the regularity of a domain. 

\smallskip

\noindent Recall that $\Omega$ is of class $C^{1,1} $ if for every $a\in \partial \Omega$ there is $r>0$ for which $B_r(a)\cap \partial \Omega= \{x=(x',x_d)\in B_r(a)~: x_d=\gamma(x')\}$ represents the graph of a $C^{1,1}$ function $\gamma: \mathbb{R}^{d-1}\to \mathbb{R} $. That is to say $\gamma$ is a $C^1$ function whose gradient is Lipschitz. 
The main result in \cite{BaS09} shows that an open set $\Omega$ is $C^{1,1} $ if and only if $\Omega$ satisfies the interior and exterior sphere condition. We say that $\Omega$ satisfies the interior and exterior sphere condition at some scale $r>0$ if for every $ a\in \partial \Omega$ one can find $a'\in \Omega$ and $a''\in \overline{\Omega}^c$ for which $B_r(a')\subset \Omega $, $B_r(a'')\subset \overline{\Omega}^c$ and $ \overline{B_r(a')}\cap \overline{B_r(a'')}= \{a\}$. 
The interior and exterior sphere condition holds for every scale $r\in (0,r_0)$ once it holds for $r_0$. This characterization entails that a $C^{1,1}$ set $\Omega$ is a $d$-set (or \emph{volume density condition} according to some authors): that is, there exist two positive constants $ c>0$ and $ r_0>0$ such that for every $r\in (0,r_0)$ and every $a \in \partial \Omega$
\begin{align*}
|\Omega \cap B_r(a)|\geq cr^d.
\end{align*}

\noindent Now we show that the classes $\mathscr{A}_i,~i=1,2,3$ are not empty. Let us simply assume that \eqref{eq:limit-at-0-explode} holds true, take for instance $\nu(h)= |h|^{-d-sp},$($s\in (0,1)$), which together with \eqref{eq:integrability-condition+even} implies \eqref{eq:non-integrability-condition}. It is easy to see that 
$$q(\delta)= \frac{1}{\delta^p} \int_{ B_\delta(0) }|h|^p\nu(h) \d h= \frac{|\mathbb{S}^{d-1}|}{p(1-s)}\delta^{-sp}\xrightarrow{\delta\to 0}\infty.$$
This shows that \eqref{eq:class-lipschitz} is verified. Wherefore, considering any Lipschitz domain $\Omega$, $(\nu, \Omega)$ in the class $\mathscr{A}_2$.

\smallskip

\noindent If $\Omega$ is of class $C^{1,1}$ we would like to show that $(\nu, \Omega)$ is in the class $\mathscr{A}_3$. Consider $R>0$, $\delta_0$ and $r_0>0$ as above. Fix $a\in \partial\Omega$, by the interior sphere condition, consider $\delta\in (0,\delta_0/4)$ small enough and $x\in \Omega$ such that $B_{2\delta}(x) \subset \Omega$, $\operatorname{dsit}(x, \partial\Omega)= |x-a|=2\delta$ and $\overline{B_{2\delta}(x)} \cap \partial\Omega= \{a\}$ then obviously, $B_\delta(x) \subset \Omega_\delta \cap B_{2\delta}(x) \subset \Omega_\delta \cap B_{4\delta}(a) $. This yields that
\begin{align}\label{eq:strong-density}
|\Omega_\delta \cap B_{4\delta}(a)|\geq d\omega_d\delta^d, \quad\text{with}~~\omega_d=|\mathbb{S}^{d-1}|.
\end{align}
Therefore, recalling that $\nu(h-a)\geq M|h-a|^{-d}\geq \frac{M}{4^d\delta^d}$ when $h\in B_{4\delta}(a)$ we have 

\begin{align*}
\int_{\Omega_\delta}\nu(h-a)\d h \geq \frac{R}{4^d\delta^d} \int_{\Omega_\delta\cap B_{4\delta}(a) } \hspace{-3ex}\d h= \, \frac{R}{4^d\delta^d}|\Omega_\delta \cap B_{4\delta}(a)|\geq \frac{ d\omega_d}{4^d}R.
\end{align*}
Finally, we get $\widetilde{q}(\delta) \geq \frac{ d\omega_d}{4^d}R$ which means that the condition \eqref{eq:class-sing-boundary} is verified since $M$ can be arbitrarily large. Thus if $\Omega$ is $C^{1,1} $ and $\nu$ satisfying \eqref{eq:limit-at-0-explode}, $(\nu,\Omega)$ belongs to $\mathscr{A}_3$. 

\vspace{2mm}

\noindent On the other hand, it is well known from \cite{Zh15} that $\Omega$ is an extension domain for $W^{s,p}(\Omega),~s\in (0,1)$ if and only if $\Omega$ is a $d$-set and thus, $ (|\cdot|^{-d-sp}, \Omega)$ is an element of the class $\mathscr{A}_1$. 

\vspace{1mm}

\begin{figure}
	\centering
	\includegraphics[scale=0.35]{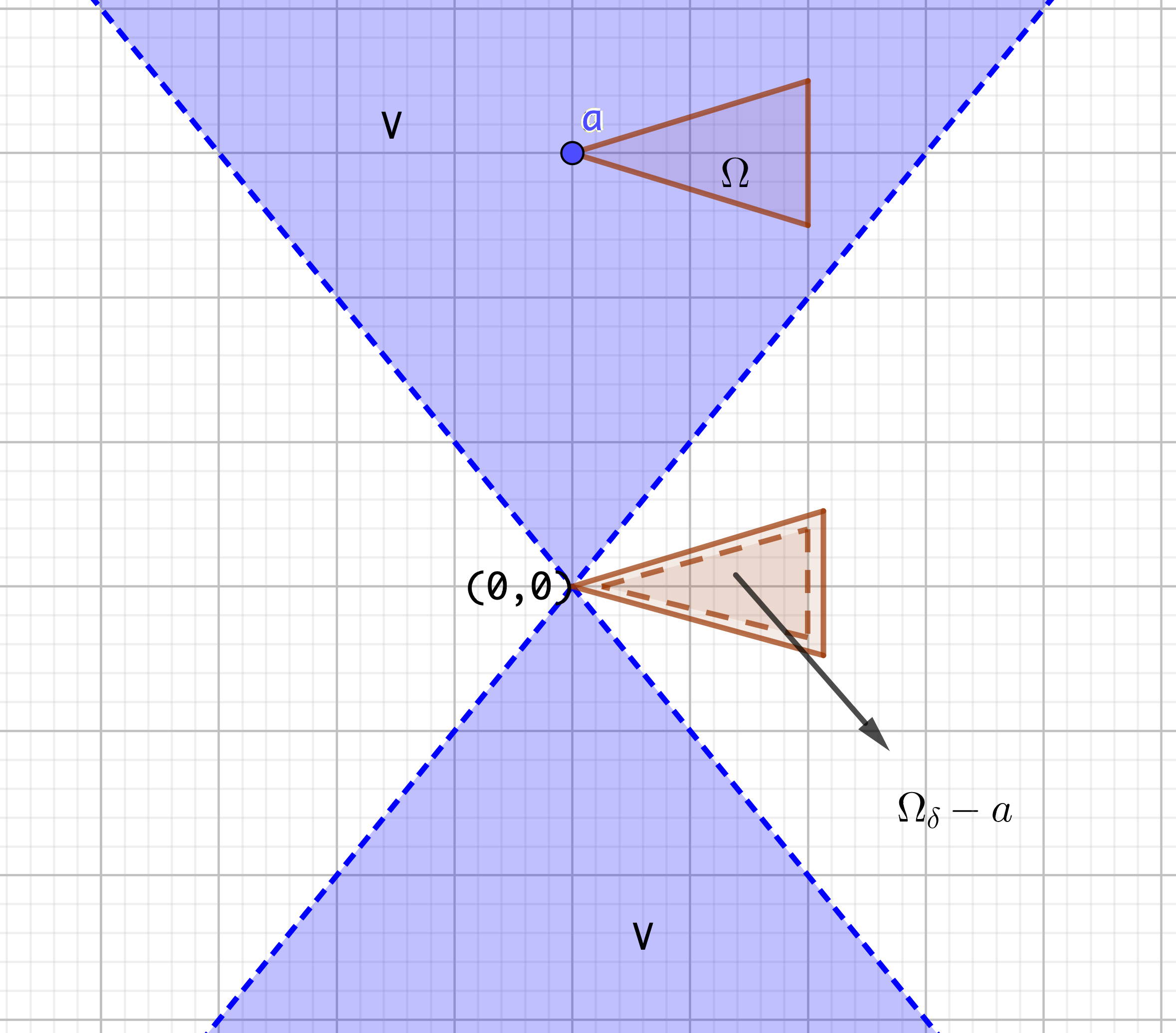} 
	\caption{ }\label{fig:cone-contra-example}
\end{figure}

\noindent Assuming \eqref{eq:integrability-condition+even} and \eqref{eq:non-integrability-condition} it is an interesting question to know under which additional minimal conditions on boundary $\partial\Omega$ and $\nu$ the condition \eqref{eq:class-sing-boundary} holds true. We illustrate this interest, with a simple counter example. 
In the Euclidean plane consider $\nu(h)= |h|^{-2-\alpha}\mathds{1}_V(h)\,\,(d=2, p=2, \alpha\in (0,2))$ with $V=\{(x_1,x_2)\in\mathbb{R}^{2}:~ |x_1|<|x_2|\}$ and $\Omega= \{(x_1,x_2)\in\mathbb{R}^{2}:~ 4|x_2-6|<x_1,~0<x_1<4\}$ whose boundary is continuous. Considering $a=(0,6)\in\partial \Omega$ one has $V\cap (\Omega_\delta-a)= \emptyset$ for every $\delta>0$ (see figure \ref{fig:cone-contra-example}). Therefore we have 
\begin{align*}
\widetilde{q}(\delta)\leq \int_{\Omega_\delta}\nu(h-a)\d h=0. 
\end{align*}
\noindent In conclusion, the condition \eqref{eq:class-sing-boundary} fails although conditions \eqref{eq:integrability-condition+even} and \eqref{eq:non-integrability-condition} are satisfied.

\vspace{2mm}

\newpage

\noindent It is interesting to know whether for small $\delta>0$, $\Omega_\delta$ inherits the regularity of $\Omega$. As proven in \cite[Section 6.14]{GT15} if $\Omega$ is of class $C^k$ with $k\geq 2$ then so is $\Omega_\delta$. Regarding the relevance of the relation \eqref{eq:strong-density}, we introduce the notion of \textit{strong volume density condition} which, together with \eqref{eq:limit-at-0-explode} will imply \eqref{eq:class-sing-boundary}. 

\begin{definition}
	We say that $\Omega$ satisfies \emph{the strong volume density condition} if there exist positive constants $r_0>0$ and $c>0$ such that for all $\delta, r\in (0, r_0)$ and 
	$a\in \partial\Omega$, one has  $	|\Omega_\delta \cap B_r(a)|\geq cr^d.$
\end{definition}

\bigskip

\noindent 
Let us resume with our quest toward the global compactness. In the spirit of \cite{DMT18} and \cite{Ponce2004} we will need some estimates near the boundary $\partial\Omega$. We start with an inequality involving cut-off functions.\linebreak[-2]

\begin{lemma}\label{lem:estimate-cut-off}
	Let $\Omega \subset \mathbb{R}^d$ be open and bounded . Assume $\nu:\mathbb{R}^d\setminus\{0\} \to [0, \infty]$ be an even measurable function. For $\delta>0$ small enough, let $\varphi \in C^\infty(\Omega)$ be such that $ \varphi=0$ on $\Omega_\delta$, $ \varphi=1$ on $\Omega\setminus \Omega_{\delta/2}$, $0\leq \varphi\leq 1$ and $|\nabla \varphi|\leq c/\delta$ (cf. Lemma \ref{lem:cut-off-existence}  by taking $\varphi = 1-\varphi_\delta$). Then for every $u\in L^p(\Omega)$, the following estimate holds true
	\begin{align}\label{eq:estimate-cut-off}
	\iint\limits_{\Omega\Omega}\big|[u\varphi](x)-[u\varphi](y)\big|^p\nu(x-y)\d x\d y
	\leq \frac{C}{\delta^p}\int\limits_{\Omega_{\delta/2}}|u(x)|^p\d x+ 2^{p}\iint\limits_{\Omega\Omega}|u(x)-u(y)|^p\nu(x-y)\d x\d y
	\end{align}
	where, $C=2^{p}c^p \int\limits_{B_R(0)}|h|^p\nu(h)\d h$ and $R=\operatorname{diam}(\Omega)$.
	
\end{lemma}
\begin{proof}
	Firstly, since $\varphi=1$ on $\Omega\setminus \Omega_{\delta/2}$ we have 
	
	\begin{align*}
	\iint\limits_{\Omega\setminus \Omega_{\delta/2}~ \Omega\setminus \Omega_{\delta/2}}\hspace*{-4ex} \big|[u\varphi](x)-[u\varphi](y)\big|^p\nu(x-y)\d x\d y 
	&=\hspace*{-3ex} \iint\limits_{\Omega\setminus \Omega_{\delta/2}~ \Omega\setminus \Omega_{\delta/2}} \hspace*{-4ex}|u(x)-u(y)|^p\nu(x-y)\d x\d y\\
	&\leq \iint\limits_{\Omega\Omega}|u(x)-u(y)|^p\nu(x-y)\d x\d y\,.
	\end{align*}
	
	\noindent In view of the fact that $0\leq \varphi\leq 1$ and $|\varphi(x)-\varphi(y)|\leq c/\delta|x-y|$ for every $x,y \in \Omega$, we have
	\begin{align}
	\big|[u\varphi](x)-[u\varphi](y)\big|^p &= \big( \varphi(y)(u(x)-u(y)) + u(x)(\varphi(x)-\varphi(y)) \big)^p \notag \\
	&\leq 2^{p-1}|u(x)-u(y)|^p+ \frac{2^{p-1}c^p}{\delta^p} |u(x)|^p|x-y|^p\,.\label{eq:split}
	\end{align} 
	Secondly, noticing that $\Omega\subset B_R(x)$ for all $ x\in \Omega$ where $R= \operatorname{diam}(\Omega)$ and integrating both sides of \eqref{eq:split} over $\Omega_{\delta/2}\times \Omega_{\delta/2}$ we obtain the following estimate

	\begin{align*}
	&\iint\limits_{\Omega_{\delta/2} \Omega_{\delta/2}}\big|[u\varphi](x)-[u\varphi](y)\big|^p\nu(x-y)\d x\d y\\
	&\leq 2^{p-1}\iint\limits_{\Omega\Omega}|u(x)-u(y)|^p\nu(x-y)\d x\d y+ \frac{2^{p-1}c^p}{\delta^p} \int\limits_{ \Omega_{\delta/2}}|u(x)|^p\d x \int\limits_{B_R(x)} |x-y|^p\nu(x-y)\d y \\
	&= 2^{p-1}\iint\limits_{\Omega\Omega}|u(x)-u(y)|^p\nu(x-y)\d x\d y+ \frac{2^{p-1}c^p}{\delta^p} \Big(\int\limits_{B_R(0)} |h|^p\nu(h)\d h\Big) \int\limits_{ \Omega_{\delta/2}}|u(x)|^p\d x \,.
	\end{align*}
	Similar to the previous estimate, using \eqref{eq:split} we get 
	\begin{align*}
	& \iint\limits_{\Omega_{\delta/2}\times \Omega\setminus \Omega_{\delta/2}} \big|[u\varphi](x)-[u\varphi](y)\big|^p\nu(x-y)\d x\d y\\
	%
	%
	&\leq 2^{p-1}\iint\limits_{\Omega\Omega}|u(x)-u(y)|^p\nu(x-y)\d x\d y+ \frac{2^{p-1}c^p}{\delta^p} \Big(\int\limits_{B_R(0)} |h|^p\nu(h)\d h\Big) \int\limits_{ \Omega_{\delta/2}}|u(x)|^p\d x \,.
	\end{align*}
	Altogether, the desired estimate follows as claimed since by symmetry we can use the split
	
	\begin{align*}
	\iint\limits_{\Omega\times \Omega}= \iint\limits_{\Omega_{\delta/2}\times \Omega_{\delta/2}} +\, 2
	\hspace*{-2ex} \iint\limits_{\Omega_{\delta/2}\times \Omega\setminus \Omega_{\delta/2}}+ \iint\limits_{\Omega\setminus\Omega_{\delta/2}\times\Omega\setminus\Omega_{\delta/2}}\,.
	\end{align*}

\end{proof}

\medskip

\noindent Since Lipschitz domains play an important role in this section let us recall another handy characterization of such domains according to  \cite[Theorem 1.2.2.2]{grisvard11}.
\begin{theorem} A bounded open subset $\Omega$ of $\mathbb{R}^d$ has the uniform cone properties if and only if its boundary is Lipschitz. In particular any open bounded convex set is Lipschitz.

\vspace{2mm}
\noindent Let us recall the following: $\Omega$ has the uniform segment property (resp. cone property ) if for every $x\in \partial \Omega$, there exists a neighborhood of $x$ in $\mathbb{R}^d$ and new coordinate system $\{y_1,y_2, \cdots, y_d\}$ such that 
	\begin{enumerate}[$(i)$]
		\item V is a hypercube, $V= \{(y_1,y_2, \cdots, y_d) ~: -a_i\leq y_i\leq a_i ~~i=1,2,\cdots, d\}$.
		\item $y-z\in \Omega$ whenever $y\in \overline{\Omega }\cap V$ and $z\in C$ where $C $ is the open segment $ \{(0,0, \cdots, z_d) ~: 0<z_d<h\}$ (resp. the open cone $\{(z', z_d)~: \cot \theta |z'|< z_d<h\}$ for some $\theta\in (0,\pi /2] $) for some $h>0$. 
	\end{enumerate} 
Note that $\tan \theta < 1/k$ where is $k$ the uniform Lipschitz constant or $\Omega$. 
\end{theorem}

\noindent We  need  the following lemma in the sequel. 
\begin{lemma}\label{lem:estimate-near-boundary}
	Assume $ \Omega\subset \mathbb{R}^d$ is open bounded and $\nu:\mathbb{R}^d\setminus\{0\}\to [0, \infty]$ is an even function. Then for every $u \in L^p(\Omega)$ and every $\delta>0$ small enough there exists a positive constant $C>0$ independent of $\delta$ such that
	\begin{align}\label{eq:estimate-boundary}
	\int_{\Omega} |u(x)|^p\d x\leq \frac{C}{\delta^p \widetilde{q}(2\delta)}\int_{\Omega_{\delta/2}}|u(x)|^p \d x+ \frac{2^{p+1}}{\widetilde{q}(2\delta)}\iint\limits_{\Omega\Omega}|u(x)-u(y)|^p\nu(x-y)\d x\d y.
	\end{align}
	Moreover if $\Omega$ has Lipschitz boundary and $\nu$ is radial then 
	\begin{align}\label{eq:estimate-boundary-lip}
	\int_{\Omega} |u(x)|^p\d x\leq \frac{C}{\delta^p q(2\delta)}\int_{\Omega_{\delta/2}}|u(x)|^p \d x+ \frac{2^{p+1}}{q(2\delta)}\iint\limits_{\Omega\Omega}|u(x)-u(y)|^p\nu(x-y)\d x\d y.
	\end{align}
\end{lemma}

\begin{proof}
	Let $\varphi$ be as in Lemma \ref{lem:estimate-cut-off} and fix $a\in \partial\Omega$. A routine check reveals that $\Omega_{2\delta}-a\subset \Omega_{\delta}-x$ for every $x\in \Omega \cap B_\delta(a)$ which yields,
	\begin{align*}
	\int_{\Omega \cap B_\delta(a)} [\varphi u]^p (x) \d x \int_{\Omega_\delta} \nu(x-y)\d y &\geq \int_{\Omega \cap B_\delta(a)} [\varphi u]^p (x) \d x \int_{\Omega_\delta-x} \nu(h)\d h\\
	&\geq \int_{\Omega \cap B_\delta(a)} [\varphi u]^p (x) \d x \int_{\Omega_{2\delta}-a } \nu(h)\d h\\
	&\geq \widetilde{q}(2\delta) \int_{\Omega \cap B_\delta(a)} [\varphi u]^p (x) \d x.
	\end{align*}
	%
	
	\noindent By a compactness argument there exist $a^1,a^2,\cdots a^n\in \partial\Omega$ such that $\partial\Omega\subset \bigcup\limits_{i=1}^n B_{\delta/2}(a^i)$. Thus, $\Omega\setminus \Omega_{\delta/2}\subset \bigcup\limits_{i=1}^n \Omega\cap B_{\delta}(a^i) \subset \Omega\setminus \Omega_{\delta}$.
	Since $\varphi u=0$ on $\Omega_{\delta/2}$ so on $\Omega_\delta$. Therefore with the aid of the above estimate we get the following estimate
	\begin{align*}
	&\iint\limits_{\Omega\Omega} \big|[\varphi u] (x)- [\varphi u] (y)\big|^p\nu(x-y)\d x\d y
	\geq 2 \iint\limits_{\Omega\setminus \Omega_{\delta} ~ \Omega_\delta} [\varphi u]^p (x) \nu(x-y)\d x\d y\\
	&\geq \int\limits_{ \bigcup\limits_{i=1}^n \Omega\cap B_{\delta}(a^i)} [\varphi u]^p (x) \d x \int_{ \Omega_\delta} \nu(x-y)\d y
	\geq 2\widetilde{q}(2\delta)\int\limits_{ \bigcup\limits_{i=1}^n \Omega\cap B_{\delta}(a^i)} [\varphi u]^p (x) \d x \\
	&\geq 2\widetilde{q}(2\delta)\int\limits_{ \Omega\setminus \Omega_{\delta/2}} [\varphi u]^p (x) \d x= 2\widetilde{q}(2\delta)\int\limits_{ \Omega\setminus \Omega_{\delta/2}} |u(x)|^p (x) \d x \,.
	\end{align*}
	This combined with \eqref{eq:estimate-cut-off} gives \eqref{eq:estimate-boundary}. 
	
	\vspace{2mm}
	\noindent Now we assume  that $\nu$ is radial and  $\Omega$ is Lipschitz.  Let $a\in \partial \Omega$ up to translating we may $a=0$. Since $ [u\varphi] = 0$ on $\Omega_\delta$, according to \cite[Eq. 22]{Ponce2004} there is a constant $C= C(\Omega,d,p)$ not depending on $u\varphi$ such that for all $w\in \Lambda\cap \mathbb{S}^{d-1}$
	\begin{align}\label{eq:ponce-on-boundary}
	\il_{\Omega \cap B_{\delta/2}(0)} | [u\varphi](x)|^p\d x\leq C\delta^p\int_{\Omega\cap B_{3\delta}(0)} \frac{|[u\varphi](x+rw) -[u\varphi]((x)|^p}{r^p}\d x 
	\end{align}
	Here $\Lambda$ is a half cone locally related to $\Omega$ and defined by $\Lambda=\{x=(x', x_d) :\, |x'|\leq x_d\}$. 
	
\noindent	In connection to the polar coordinates, since $\nu $ is radial, integrating the above inequality with respect to the measure $ |h|^p\nu(h)\d h:= r^{p+d-1}\nu(r)\d \sigma_{d-1}(w) \d r$ over $\Lambda\cap \mathbb{S}^{d-1} \times (0, \delta) $ and letting $c_1= |\Lambda\cap \mathbb{S}^{d-1}|/|\mathbb{S}^{d-1}|$ gives
	\begin{align*}
	\begin{split}
	\Big(c_1\il_{B_{\delta}(0)} |h|^p\nu(h) \d h\Big)\il_{\Omega \cap B_{\delta/2}(0)} |[u\varphi](x)|^p\d x
	&\leq C\delta^p\iil_{ \Omega \cap B_{3\delta}(0)\,B_{\delta}(0)}  | [u\varphi](x+h) - [u\varphi](x)|^p\nu(h)\d h\d x \\
	&\leq C\delta^p\iil_{ \Omega \cap B_{4\delta}(0)\,\Omega\cap B_{4\delta}(0)}  | [u\varphi](y) - [u\varphi](x)|^p\nu(y-x)\d y\d x,
	\end{split}
	\end{align*} 
	where the last inequality holds since $x+h \in \Omega \cap B_{4\delta}(0)$ for $x \in \Omega \cap B_{3\delta}(0)$ and $h \in B_{\delta}(0)$. Up to doubling the parameter $\delta>0$ relabelling the constant $C>0$ that is 
	\begin{align*}
	\Big(\il_{B_{2\delta}(0)} |h|^p\nu(h) \d h\Big)\il_{\Omega \cap B_{\delta}(0)} | [u\varphi](x)|^p\d x
	&\leq C\delta^p\iil_{ \Omega \cap B_{8\delta}(0)\,\Omega\cap B_{8\delta}(0)}  | [u\varphi](y) - [u\varphi](x)|^p\nu(y-x)\d y\d x 
	\end{align*} 
	and using once more  a recovering argument as previously (since $\Omega\setminus \Omega_{\delta/2}\subset \bigcup\limits_{i=1}^n \Omega\cap B_{\delta}(a^i) \subset\bigcup\limits_{i=1}^n \Omega\cap B_{8\delta}(a^i)$)  one reaches the following estimate 
	\begin{align}\label{eq:estimate-near-Lipschitz-boundary}
	\begin{split}
	\Big(\il_{B_{2\delta}(0)} |h|^p\nu(h) \d h\Big)\il_{\Omega \setminus \Omega_{\delta/2}} | [u\varphi](x)|^p\d x
	&\leq C\delta^p\iil_{ \Omega \Omega}  | [u\varphi](y) - [u\varphi](x)|^p\nu(y-x)\d y\d x. 
	\end{split}
	\end{align}
	Hence,
	\begin{align*}
	2q(2\delta) \int\limits_{ \Omega\setminus \Omega_{\delta/2}} [\varphi u]^p (x) \d x\leq C\iint\limits_{\Omega\Omega} \big|[\varphi u] (x)- [\varphi u] (y)\big|^p\nu(x-y)\d x\d y.
	\end{align*}
	which combined with \eqref{eq:estimate-cut-off} implies \eqref{eq:estimate-boundary-lip}.
		
	%
\end{proof}

\medskip

\noindent Here is the global compactness result.
\begin{theorem}\label{thm:embd-compactness}
	Let $\Omega\subset \mathbb{R}^d$ be open, bounded and $\nu:\mathbb{R}^d\setminus\{0\} \to [0,\infty]$ be a measurable function. If the couple $(\nu, \Omega)$ belongs to  one of the class $\mathscr{A}_i,~ i=1,2,3$ then the embedding $\WnuOm \hookrightarrow L^p(\Omega)$ is compact. In particular, the embedding $\WnuOmR \hookrightarrow L^p(\Omega)$ is compact.
\end{theorem}

\medskip 
\begin{proof} Given the continuous embedding $\WnuOmR\hookrightarrow\WnuOm$, it will be sufficient to prove that the embedding $\WnuOm\hookrightarrow L^p(\Omega)$ is compact. For $(\nu, \Omega)$ belonging to the class $\mathscr{A}_1$ the result is a direct consequence of Theorem \ref{thm:local-compactness}. Now assume $(\nu, \Omega)$ belongs to o the class $\mathscr{A}_2$ (resp. $\mathscr{A}_3$) then for $\varepsilon>0$ there is $\delta>0$ small enough such that $2^{p+1}q^{-1}(2\delta)<\varepsilon$ (resp. $2^{p+1} \widetilde{q}^{-1}(2\delta)<\varepsilon$) If $(u_n)_n$ is a bounded sequence of $\WnuOm$ then Corollary \ref{cor:local-compatcness} infers the existence of a subsequence $(u_{n_j})_j$ of $(u_n)_n$ converging to some $u \in \WnuOm$ in $L^p(\Omega_{\delta/2})$ i.e $\|u_{n_j}-u\|_{L^p(\Omega_{\delta/2})}\to0$ as $j\to \infty$. In any case, in view of Lemma \ref{lem:estimate-near-boundary}, passing to the limsup in \eqref{eq:estimate-boundary} or in \eqref{eq:estimate-boundary-lip} applied to $u_{n_j}-u$ we get 
	\begin{align*}
	\limsup_{j\to \infty}\int_{\Omega}|u_{n_j}(x)-u(x)|^p\d x\leq M\varepsilon
	\end{align*}
	where $M= 2^p\|u\|^p_{\WnuOm}+ 2^p\sup\limits_{n\geq 1}\|u_n\|^p_{\WnuOm}<\infty.$ Finally, letting $\eps\to 0$ gives  $\limsup\limits_{j\to \infty}\|u_{n_j}-u\|_{L^p(\Omega)}=0$.  Which ends the proof. 
\end{proof}

\medskip 

\noindent A noteworthy consequence of Theorem \ref{thm:embd-compactness} is the well known and established  Rellich-Kondrachov's\footnote{The classical Rellich-Kondrachov theorem is often known as the Kondrachov compactness theorem after V.Kondrachov who generalized Rellich's result for compact map $W^{1,p}_{0}(\Omega)$ into $L^q(\Omega)$ whenever $1\leq q\leq d p/(d-p)$.}. 

\begin{corollary}[Rellich-Kondrachov's Theorem]
	Assume $\Omega\subset \R^d$ is open and bounded. The embedding  $W_0^{1,p}(\Omega)\hookrightarrow L^p(\Omega) $ is compact. Further if $\Omega$ is Lipschitz, then  the embedding $W^{1,p}(\Omega)\hookrightarrow L^p(\Omega)$ is compact. 
\end{corollary}

\begin{proof}
 The embedding $W_0^{1,p}(\Omega)\hookrightarrow L^p(\Omega)$ is compact because, with the choice $\nu(h) = |h|^{-d-sp},$ since in view of Theorem \ref{thm:local-compactness} the embedding $\WnuOmO \hookrightarrow L^p(\Omega)$  is compact and the embedding $W_0^{1,p}(\Omega)\hookrightarrow \WnuOmO $ is continuous.

 \noindent Now if $\Omega$ is  Lipschitz then by Theorem \ref{thm:w1p-in nonlocal} we know that the embedding $W^{1,p}(\Omega)\hookrightarrow \WnuOm $ is continuous whereas from Theorem \ref{thm:embd-compactness} it comes that the embedding $\WnuOm \hookrightarrow L^p(\Omega)$  is compact. It thus turns out that  $W^{1,p}(\Omega)\hookrightarrow L^p(\Omega)$ is compact too. 
\end{proof}

\medskip

\noindent The efforts made to establish the compactness  result of Theorem \ref{thm:embd-compactness} will be rewarded for the elaboration of the Poincar\'e type inequality which will be useful in the forthcoming section.



\begin{theorem}[Poincar\'e inequality]\label{thm:poincare-inequality}
	Let $\Omega$ be an open bounded subset of $\mathbb{R}^d$ and $\nu:\mathbb{R}^d\setminus\{0\} \to [0,\infty]$ be a measurable function with full support. Assume the couple $(\nu, \Omega)$ belongs to one of the class $\mathscr{A}_i,~ i=1,2,3,4$. Then there exists a positive constant $C= C(d,\Omega, \nu )$ depending only on $d,~ \Omega$ and $\nu$ such that
	
	\begin{align}\label{eq:poincare-inequalityH}
	\big\|u-\mbox{$\fint_{\Omega}$}u\big\|^p_{L^p(\Omega)} \leq C\iint\limits_{\Omega\Omega}|u(x)-u(y)|^p\nu(x-y)\d x\d y\qquad\text{for every }~ u \in L^p(\Omega)
	\end{align}
	and immediately, 
	\begin{align}\label{eq:poincare-inequalityV}
	\big\|u-\mbox{$\fint_{\Omega}$}u\big\|^p_{L^p(\Omega)} \leq C\iint\limits_{\Omega\R^d}|u(x)-u(y)|^p\nu(x-y)\d x\d y\,\qquad\text{for every }~ u \in L^p(\Omega)\, .
	\end{align}
\end{theorem}

\begin{proof}
	Assume such constant does not exist then we can find a sequence $(u_n)_n$ elements of $\WnuOm$ such that for every $n$, $\fint_{\Omega}u_n =0$, $\|u_n\|_{L^p(\Omega)}=1$ and 
	\begin{align*}
	\iint\limits_{\Omega\Omega}|u_n(x)-u_n(y)|^p\nu(x-y)\d x\d y\leq \frac{1}{2^n}\,.
	\end{align*}

	\noindent The sequence $(u_n)_n$ is thus bounded in $\WnuOm$ which by Theorem \ref{thm:embd-compactness} is compactly embedded in $L^p(\Omega)$ whenever $(\nu,\Omega)$ belongs to one of  the class $\mathscr{A}_i,~ i=1,2,3$. Therefore, if it is the case, passing through a subsequence, $(u_n)_n$ converges in $L^p(\Omega)$ to some function $u$. Clearly it follows that $\fint_{\Omega}u =0$ and $\|u\|_{L^p(\Omega)}=1$. Moreover, by Fatou's lemma we have 
	\begin{align*}
	\iint\limits_{\Omega\Omega}|u(x)-u(y)|^p\nu(x-y)\d x\d y\leq \liminf_{n \to \infty} \iint\limits_{\Omega\Omega}|u_n(x)-u_n(y)|^p\nu(x-y)\d x\d y=0
	\end{align*}
	which implies that $u$ equals the constant function $x\mapsto \fint_{\Omega}u =0$ almost everywhere on $\Omega $. This goes against the fact that $\|u\|_{L^p(\Omega)}=1$ hereby showing that our initial assumption was wrong. 
	
	\smallskip
	
	\noindent Next assume $(\nu,\Omega)$ belongs to the class $\mathscr{A}_4$ then, as $\nu$ has a full support, is radial and is almost decreasing and $\Omega$ is bounded, there is a constant $c'>0$ such that $\nu(x-y)\geq c'$ for all $x,y \in \Omega$. Using this and Jensen's inequality we get the desired inequality as follows
	\begin{align*}
	\iint\limits_{\Omega\Omega}|u(x)-u(y)|^p\nu(x-y)\d x\d y&\geq c' |\Omega| \int_{\Omega} \fint_{\Omega}|u(x)-u(y)|^p\d x\d y\\
	&\geq c' |\Omega| \|u-\mbox{$ \fint_{\Omega} u $} \|^p_{L^p(\Omega)}\,.
	\end{align*}
	Which ends the proof since \eqref{eq:poincare-inequalityV} is clearly a consequence of \eqref{eq:poincare-inequalityH}.
\end{proof}

\medskip

\noindent The above Poincar\'e inequality \eqref{eq:poincare-inequalityH}-\eqref{eq:poincare-inequalityV} can be seen has the nonlocal counterpart of the classical Poincar\'e inequality which states that for a connected bounded Lipschitz $\Omega$, there is $C>0$ for which 
\begin{align*}
\big\|u-\mbox{$\fint_{\Omega}$}u\big\|_{L^p(\Omega)} \leq C\|\nabla u\|_{L^p(\Omega)}\,,\qquad\text{for all }~ u \in L^p(\Omega) 
\end{align*}
where by convention we assume $\|\nabla u\|_{L^p(\Omega)} =\infty$ if $|\nabla u|$ is not in $L^p(\Omega)$. Alongside this, we also recall the classical Poincar\'e-Friedrichs inequality: there is $C>0$ such that
\begin{align*}
\|u\|_{L^p(\Omega)} \leq C\|\nabla u\|_{L^p(\Omega)}\,\qquad\text{for all }~ u \in W_0^1(\Omega)\,. 
\end{align*}
\vspace{2mm}

\noindent With a close look at the proof of Theorem \ref{thm:completness-Wnup}  one readily establishes  the following generalization.

 \begin{theorem}[Poincar\'e inequality]\label{thm:poincare-inequality-general}
 Let $\Omega$ be an open bounded subset of $\mathbb{R}^d$ and $\nu:\mathbb{R}^d\setminus\{0\} \to [0,\infty]$ be a measurable function with full support. Assume the couple $(\nu, \Omega)$ belongs to the $\mathscr{A}_i,~ i=1,2,3$.  Assume that there $G\subset L^p(\Omega)$ is a nonempty closed subset of $L^p(\Omega)$ not containing non-zero constant functions.   Then there exists a positive constant $C= C(d,\Omega, \nu, G )$ such that
 	\begin{align*}
 	\big\|u\big\|^p_{L^p(\Omega)} \leq C\iint\limits_{\Omega\Omega}|u(x)-u(y)|^p\nu(x-y)\d x\d y\qquad\text{for every }~ u \in G. 
 	\end{align*}
 \end{theorem}
\noindent By a subsequent analogy,  one can also establish the following result as well. 
\begin{theorem}[Poincar\'e inequality]\label{thm:poincare-inequality-general-bis}
	Let $\Omega$ be an open bounded subset of $\mathbb{R}^d$ and $\nu:\mathbb{R}^d\setminus\{0\} \to [0,\infty]$ be a measurable function with full support. Assume the couple $(\nu, \Omega)$ belongs to the $\mathscr{A}_i,~ i=1,2,3$.  Let $G\subset\WnuOmR$ be a nonempty closed subset of $\WnuOmR$ not containing non-zero constant functions.    Then there exists a positive constant $C= C(d,\Omega, \nu, G )$ such that
	\vspace{-1ex}
	\begin{align*}
	\big\|u\big\|^p_{L^p(\Omega)} \leq C\iint\limits_{\Omega\R^d}|u(x)-u(y)|^p\nu(x-y)\d x\d y\qquad\text{for every }~ u \in G.
	\end{align*}
\end{theorem}

\medskip

\begin{remark}
	$(i)$ Let $E\subset \Omega$ be a measurable set with a positive measure and $\delta>0$. Some special closed subsets of $L^p(\Omega)$ to which Theorem \ref{thm:poincare-inequality-general} applies are given by, $G_E=\{u\in L^p(\Omega): \,\mbox{$\fint_E u$}=0\,\}$, $
	G^0_E=\{u\in L^p(\Omega): \, u=0\,\text{ a.e. on $E$}\}$ and 
	$G_\delta=\{u\in L^p(\Omega): \, |\{u=0\}|\geq \delta\}.$
	

	\vspace{2mm}
	\noindent $(ii)$ Let $E\subset \R^d$ be a measurable set with a positive measure and $\delta>0$. Some special closed subsets of $\WnuOmR$ to which Theorem \ref{thm:poincare-inequality-general-bis} applies are given by, 
		$G_E=\{u\in \WnuOmR: \,\mbox{$\fint_E u$}=0\,\}$, $
	G^0_E=\{u\in \WnuOmR: \, u=0\,\text{ a.e. on $E$}\}$ and 
	$G_\delta=\{u\in \WnuOmR: \, |\{u=0\}|\geq \delta\}.$
	
\end{remark}

%
%

\vspace{2mm}
\noindent In the same spirit, as we will see below the corresponding nonlocal Poincar\'e-Friedrichs inequality $\WnuOmO$ (which we recall is the closure of the $C_c^\infty(\Omega)$ in $\WnuOmR$) is much easier to obtain and no compactness argument is required. This provides an effortless alternative  way to proof  the Poinca\'e-Friedrichs  inequality from \cite[Lemma 2.7]{FKV15} which uses iterated convolutions. Nevertheless, we point out that under the condition that the embedding is $\WnuOmO \hookrightarrow L^p(\Omega)$ is compact, a similar inequality is derived in \cite{JW19} when $\Omega$ is only  bounded in one direction. 

\medskip  

\begin{theorem}[Poincar\'e-Friedrichs inequality for $\WnuOmO$]\label{thm:poincare-friedrichs-ext}
	Let $\Omega$ be a open bounded subset of $\mathbb{R}^d$ and $\nu:\mathbb{R}^d\setminus\{0\}\to [0, \infty]$ be a function satisfying \eqref{eq:integrability-condition-near-zero}. There is a constant $C=C(d,\Omega, \nu)$ such that \linebreak[-2]
	
	\vspace{-4ex}

	\begin{align}\label{eq:poincare-friedrichs-ext}
	\|u\|^p_{L^p(\Omega)}\leq C |u|^p_{\WnuOmR}\, \quad \text{for every $u\in \WnuOmO$}.
	\end{align}
	
\end{theorem}

\medskip

\begin{proof}
	Let $R>0$ be the diameter of $\Omega$ then for all $x\in \Omega$ we have $B^c_R(x)\subset \Omega^c$. For $u\in \WnuOmO$ we recall that $u=0~a.e$ on $\Omega^c$. It suffices to take C= $(2\|\nu_R\|_{L^1(\mathbb{R}^d)})^{-1}$ with $\nu_R=\nu \mathds{1}_{\mathbb{R}^d\setminus B_R(0)}$ indeed
	\begin{align*}
	|u|^p_{\WnuOmR}&=\iint\limits_{\Omega\Omega}|u(x)-u(y)|^p\nu(x-y)\d x \d y+ 2\int_{\Omega}|u(x)|^p\d x\int_{\Omega^c} \nu(x-y)\d y\\
	&\geq 2\int_{\Omega}|u(x)|^p\d x\int_{B^c_R(x)} \nu(x-y)\d y= 2\|\nu_R\|_{L^1(\mathbb{R}^d)}\|u\|^p_{L^p(\Omega)}.
	\end{align*}
	
\end{proof}

\vspace{5mm}

\noindent In view of the compactness of $W^p_{\nu,0}(\Omega)$ (the closure of $C_c^\infty(\Omega)$ in $\WnuOm$) into $L^p(\Omega)$ it is very tempting to say that the nonlocal counterpart of the Poincar\'e-Friedrichs inequality also holds on $W^p_{\nu,0}(\Omega)$ under the only assumption that $(\nu,\Omega)$ is in the class $\mathscr{A}_i,~i=1,2,3$. But this is not warranted especially if one considers the simple case $\nu(h)= |h|^{-d-sp}$ with $0<s <1/p$. Indeed, is well known that $C_c^\infty(\Omega)$ is dense in $W^{s,p}(\Omega)$ when $0<s <1/p$ (c.f. \cite{grisvard11,Tri83}) that is to say in this situation we have $\WnuOm =W^p_{\nu,0}(\Omega) $. In fact, since $\WnuOm$ contains constant functions, it is impossible to find a $C>0$ for which the following Poincar\'e-Friedrichs inequality holds
\begin{align}\label{eq:poincare-friedrichs}
\|u\|^p_{L^p(\Omega)}\leq C\iint\limits_{\Omega\Omega}|u(x)-u(y)|^p\nu(x-y)\d x\d y\, \qquad\text{for all}~~ u\in C_c^\infty(\Omega)\,.
\end{align}

\noindent More generally, following the scheme of the proof of Theorem \ref{thm:poincare-inequality} one is able to establish the following.

\begin{theorem}[Poincar\'e-Friedrichs inequality for $W^p_{\nu,0}(\Omega)$] Let $\Omega\subset \R^d$ be an open bounded set  and a measurable function $\nu:\mathbb{R}^d\to [0,\infty]$ satisfying \eqref{eq:integrability-condition-near-zero} such that the embedding $\WnuOm \hookrightarrow L^p(\Omega)$ is compact. Then the Poincar\'e-Friedrichs inequality \eqref{eq:poincare-friedrichs} holds true on the space $ W^p_{\nu,0}(\Omega)$ if and only if the constant $u=1$ is not an element of $W^p_{\nu,0}(\Omega)$. 
\end{theorem}

\vspace{2mm}
\noindent We are also in the mood to establish the so called nonlocal Friedrichs type inequality. To do this let us start with the following well known result. 

\begin{theorem}\label{thm:banach-compact}
	Let $X,Y$ and $Z$ be three Banach spaces such that $X$ is reflexive. Let $K\in \mathcal{L}(X,Y)$ be a compact operator and $S\in \mathcal{L}(X,Z)$ be one-to-one. Then for every $\varepsilon>0$ there exists a constant $C_\varepsilon>0$ such that for all $x\in X$
	\begin{align*}
	\|Kx\|_Y\leq \varepsilon \|x\|_X+C_\varepsilon \|Sx\|_Z,
	\end{align*}
\end{theorem} 

\begin{proof}
	Assume the claim fails for some $\varepsilon>0$. Then there exists a sequence $(x_n)_n$ of X such  that for all $n\in \mathbb{N}$ we have $\|Kx_n\|_Y>\varepsilon \|x_n\|_X+2^n\|Sx_n\|_Z.$
Preferably we assume $\|x_n\|_X =1$. Hence, since $X$  is reflexive, a  subsequence $(x_n)$ weakly converges to some $x\in X$. So that $Sx_n$ also converges weakly in $Z$ to $Sx$  which implies $\|Sx\|_Z \leq \liminf\limits_{n\to \infty} \|Sx_n\|_Z$ and  by compactness of $K$, $Kx_n$ converges to $Kx$ in $Y$. Simultaneously, the above relation forces $\|Sx_n\|_Z\to 0$ that is $Sx=0$ which means $x=0$ since $S$ is one-to-one and  hence $\|Kx\|_Y>\varepsilon$ which is a contradiction.
\end{proof}
\noindent As a consequence of this we have the following nonlocal version of Friedrichs' inequality. 
\begin{theorem}[Nonlocal Friedrichs inequality]
	Let $\Omega$ be an open bounded subset of $\mathbb{R}^d$, $1<p<\infty$ and $\nu:\mathbb{R}^d\setminus\{0\} \to [0,\infty]$ be a measurable function with full support. Assume the couple $(\nu, \Omega)$ belongs to the $\mathscr{A}_i,~ i=1,2,3$.  Assume $K\subset \Omega$ (eventually $K=\Omega$ ) and  $O\subset \Omega^c$(eventually $O=\Omega^c$)  be a measurable sets such that $|K|>0$ and $|O|>0$. There exists a constant $C>0$ such that for all $u \in \WnuOmR$, 
	\begin{align*}
	\|u\|^p_{L^p(\Omega)}&\leq C|u|^p_{\WnuOmR}+ C \|u\|^p_{L^p(O, \, \nu_K)},\\
		\|u\|^p_{L^p(\Omega)}&\leq C|u|^p_{\WnuOmR}+ C \|u\|^p_{L^p(O, \,\mathring{\nu}_K)}
	\end{align*}
 where we recall that 
\begin{align*}
\nu_K(x) = \operatorname{essinf}_{y\in K}\nu(x-y)\qquad \text{and} \qquad\mathring{\nu}_K(x) = \int_K 1\land \nu(x-y)\d y.
\end{align*}
In particular the following norms are mutually equivalent to $\|\cdot\|_{\WnuOmR}$,
\begin{align*}
	&u \mapsto \Big( |u|^p_{\WnuOmR} +\|u\|^p_{L^p(O, \nu_K)}\Big)^{1/p}
	\quad\text{and} \quad u \mapsto \Big( |u|^p_{\WnuOmR} +\|u\|^p_{L^p(O, \mathring{\nu}_K)}\Big)^{1/p}.
\end{align*}
\end{theorem}
\medskip

\begin{proof}
	The embedding $\WnuOmR\hookrightarrow L^p(\Omega)$ is compact. A mere adaption of Proposition \ref{prop:nonlocal-L2-trace} shows that $\operatorname{Tr}: \WnuOmR \to L^p(\Omega^c, \nu_K)$ and $\operatorname{Tr}: \WnuOmR \to L^p(\Omega^c, \mathring{\nu}_K)$ with $u \mapsto \operatorname{Tr}(u) = u\mid_{\Omega^c}$ are linear and continuous. On the other hand we trivially have that the embeddings $ L^p(\Omega^c, \nu_K) \hookrightarrow L^p(O, \nu_K)$ and  $ L^p(\Omega^c, \mathring{\nu}_K) \hookrightarrow L^p(O, \mathring{\nu}_K)$ are continuous.  Whence the mappings $S: \WnuOmR\to L^p(\Omega\times \mathbb{R}^d)\times L^p(O, \nu_K)$ and $\mathring{S}: \WnuOmR\to L^p(\Omega\times \mathbb{R}^d)\times L^p(O, \nu_K)$  with 
	$$\mathring{S}u= Su = \big(u(x)-u(y))\nu^{1/p}(x-y), \operatorname{Tr}u\big) $$ 
	are linear bounded and one-to-one.  In virtue of Theorem \ref{thm:banach-compact}, for $\eps>0$ we have 
	\begin{align*}
	\|u\|^p_{L^p(\Omega)}& \leq \varepsilon \|u\|^p_{\WnuOmR}+  C_\varepsilon \|Su\|^2_{ L^2(O, \nu_K)}\\
	&=\varepsilon\|u\|^p_{L^p(\Omega)}+  (\varepsilon+ C_\varepsilon)|u|^p_{\WnuOmR} + C_\varepsilon\|u\|^p_{L^p(O, \nu_K)}.
	\end{align*}
Likewise, we have 
	\begin{align*}
\|u\|^p_{L^p(\Omega)}& \leq \varepsilon\|u\|^p_{L^p(\Omega)}+  (\varepsilon+ C_\varepsilon)|u|^p_{\WnuOmR} + C_\varepsilon\|u\|^p_{L^p(O, \mathring{\nu}_K)}.
\end{align*}
	Taking $\varepsilon=1/2$ provides the required inequalities. That the norms are equivalent blatantly follows.  
\end{proof}

\noindent Let us recall the classical Friedrichs inequality whose proof can be derived analogously.

\begin{theorem}[c.f. \cite{maz2013sobolev}]
Assume $\Omega\subset \R^d$ has a Lipschitz boundary. Let $\Gamma_0\subset \partial\Omega$ be a surface with a positive area, i.e. $|\Gamma_0|>0$. Let $1<p<\infty$. Then there exists a constant $C>0$ such that 
\begin{align*}
\|u\|_{L^p(\Omega)}\leq C\|u\|_{W^{1,p}(\Omega)} + C \|u\|_{L^p(\Gamma_0)}
\end{align*}
Moreover, the norm $\|\cdot\|_{W^{1,p}(\Omega)}$ is equivalent to the norm
$$u \mapsto \Big( |u|^p_{W^{1,p}(\Omega)} +\|u\|^p_{L^p(\Gamma_0)}\Big)^{1/p}.$$
\end{theorem}

\chapter[Complement Value Problems]{ \relsize{-1}{Complement Value Problems For Nonlocal Operators}}\label{chap:IDEs}

The overreaching goal of this chapter is to investigate weak solutions to some specific integrodifferential equations (IDEs) involving nonlocal operators of L\'{e}vy type. In many cases, this will be equivalent to proving the existence of minimizers via the direct method of calculus of variations\footnote{Calculus of variations is a branch of mathematical analysis dealing with optimization problems to find the extremum for a functional.}. At first, we shall be concerned with the well-posedness of Dirichlet (first), Neumann (second), Robin (third) and mixed complement value type problems for integrodifferential operators (IDEs) of  L\'{e}vy type.  Afterwards, we study the spectral decomposition of L\'evy type operators that are subject of the aforementioned complement conditions. The latter will
allow us to study evolution IDEs problems on bounded domains, Dirichlet-to-Neumann map  and essentially self-adjointness for integrodifferential operators. Our approaches substantially consist of developing the aforementioned notions using gadgets from $L^2$-theory. Analogous approaches treating standard elliptic  PDEs  of
the second order are referenced in \cite{AA15,Ev10,Hunter14,Dorothee-Triebel,Jost-Jurgen,LeDret16,Giovanni13,Mikhailov78,Taylor}. We begin by reviewing some important results from the theory of calculus of variations.

\section{Review of variational principles}
 
\begin{definition}
Let $ (V,\|\cdot\|_V)$ be a normed space. Let $\ell: V\to \mathbb{R}$ be a linear form and $a(\cdot, \cdot): V\times V\to \mathbb{R}$ be a bilinear form.

\begin{itemize}
 \item $\ell $ is bounded if there exists $M>1$ such that 
 \begin{align*}
 |\ell(v)|\leq M\|v\|_V \quad \text{for all} ~~ v\in V.
 \end{align*}
\item $a(\cdot,\cdot) $ is bounded if there exists $M>0$ such that 
 \begin{align*}
 |a(u,v)|\leq M\|u\|_V\|v\|_V \quad \text{for all} ~~u, v\in V.
 \end{align*}
 \item $a(\cdot,\cdot) $ is called to be coercive or $V-$elliptic if there exists $0<\theta<1$ such that 
 \begin{align*}
 a(v,v)\geq \theta\|v\|^2_{V} \quad \text{for all} ~~v \in V.
 \end{align*}
 In short, we will simply say  that $a(\cdot, \cdot )$ is $\theta$-coercive. 
\end{itemize}
\end{definition}
 
\noindent The dual space $V'$ is the collection of linear continuous forms on $V$ and is a Banach space under the norm
\begin{align*}
 \|\ell \|_{V'} =\sup_{v\in V,~\|v\|_V=1} |\ell(v)|\,.
\end{align*}

 \medskip 
 \noindent We omit the proof of the next theorem. 
 \begin{theorem}\label{thm:minimizer-var}
 Let $ (V,\|\cdot\|_V)$ be a normed space and $K\subset V$ be a convex subset. Let 
 $\ell: V\to \mathbb{R}$ be a linear form and let $a(\cdot, \cdot): V\times V\to \mathbb{R}$ be a symmetric 
 and positive definite bilinear form. Then the functional $J: V\to \mathbb{R}$ with 
 $$ J(v) =\frac{1}{2}a(v,v)-\ell(v) $$
 is strictly convex. Moreover, there is at most one $u \in K$ such that 
 $$ J(u)=\min\limits_{v\in K}J(v).$$ 
Moreover, this minimization is equivalent to the variational inequality 
 \begin{align*}
 a(u, v-u)\geq \ell(v-u)\quad\text{for all}~v\in K\,.
 \end{align*}
 Furthermore, in the special case where $K$ is an affine subspace of $V$, i.e. $K = v_0+U$ with $v_0\in V$ and $U$ is a closed subspace of $V$, the above variational inequality becomes
 \begin{align*}
 a(u, v)= \ell(v)\quad\text{for all}~v\in K\,.
 \end{align*}
 \end{theorem}
 
 \medskip
 
 \noindent The next theorem shows that in the setting of Banach spaces, the above variational 
 inequality is well-posed in the sense of Hadamard. In other words, it possesses a unique 
 solution which continuously depends upon the data. This is useful to show the well-posedness of many variational equations. 
 
 \begin{theorem}[Stampachia\footnote{This theorem was established by Guido Stampachia in the 
 setting of Hilbert spaces in 1964, he extended it later in 1967 to Banach spaces in a joint work with Jacques-Louis Lions. See for example \cite{ET09} and other references therein}]\label{thm:stampachia}
 Let $V$ be a Banach space and $K\subset V$ be a nonempty, closed, convex set. Let $a(\cdot, \cdot): V\times V\to \mathbb{R}$ be a continuous and  $\theta$-coercive bilinear form. Then for every continuous 
 linear form $\ell:V\to \mathbb{R}$ there is a unique $u\in K$ such that 
 \begin{align}\label{eq:var-inequality}
 a(u, v-u)\geq \ell(v-u)\quad \text{for all}~ v\in K\,.
 \end{align}
 Moreover, if $\overline{u}\in K$ corresponds to  another continuous linear form $\overline{\ell}:V\to \mathbb{R}$ then,
 \begin{align}\label{eq:well-posed}
 \|u-\overline{u}\|_V\leq \frac{1}{\theta}\|\ell-\overline{\ell}\|_{V'}
 \end{align}
 \end{theorem}

 \medskip
 
 \begin{proof}
Assume $u, \overline{u}\in K$ be solutions corresponding to $\ell, \overline{\ell}$. Since $K$ is convex, testing both $u$ and $\overline{u}$ with $v= \frac{u+\overline{u}}{2}\in K$ and adding both inequalities one easily arrives at
\begin{align*}
 -a(u-\overline{u},u-\overline{u})\geq (\ell-\overline{\ell})(u-\overline{u})\,.
\end{align*}
The coercivity yields 
\begin{align*}
 \theta \|u-\overline{u}\|^2_{V}\leq a(u-\overline{u}, u-\overline{u})\leq \|\ell-\overline{\ell}\|_{V'}\|u-\overline{u}\|_V\,.
\end{align*}
 This entails the estimate \eqref{eq:well-posed} from which the uniqueness follows subsequently. 
 Now we prove the existence which in virtue of Theorem \ref{thm:minimizer-var} corresponds to show that the functional
 \begin{align*}
 J(v) = \frac{1}{2}a(v,v)-\ell(v)\,.
 \end{align*}
 has a minimizer on $K$. First of all, the quantity $d=\inf\limits_{v\in K} J(v)$ exists,  since $J$ is bounded below on $V$. Indeed, for $v\in V,$
 \begin{align*}
J(v)\geq \frac{\theta}{2}\|v\|^2_V -\|\ell\|_{V'}\|v\|_V= \Big(\sqrt{\frac{\theta}{2}}\|u\|_V^2-\frac{1}{\sqrt{2\theta}} \|\ell\|_{V'}\Big)^2 - \frac{1}{2\theta }\|\ell\|^2_{V'} \geq -\frac{1}{2\theta }\|\ell\|^2_{V'}\,.
 \end{align*}
 
%
\noindent First, we assume that $a(\cdot, \cdot)$ is symmetric. For each $n\in \mathbb{N}$ we consider $u_n \in K$ such that $d\leq J(u_n)\leq d+\frac{1}{2^n}$
Whence, since $J$ is convex, from the relation $2ab = a^2+b^2- 4\big(\frac{a+b}{2}\big)^2$ we get 
\begin{align*}
 \theta \|u_n-u_m\|^2_{V}&\leq a(u_n-u_m, u_n-u_m)\\
 &= 4J(u_n)+4J(u_m)-8J\big(\frac{u_n+u_m}{2}\big)\\ 
 &\leq 4\big(d+ \frac{1}{2^n}\big)+ 4\big(d+ \frac{1}{2^m}\big)-8d= 4\big( \frac{1}{2^n}+ \frac{1}{2^m}\big)\,.
\end{align*}
Therefore, $(u_n)_n$ is a Cauchy sequence in the Banach space $V$ and hence converges to some $u\in K$ since $K$ is closed. By continuity of $J$, we have $d=J(u)$. 

\noindent If $a(\cdot, \cdot)$ is not symmetric, then for fixed $t\in [0,1]$ we write
\begin{align*}
 a_t(u,v) = a_0(u,v)+ tb(u,v)
\end{align*}
with 
\begin{align*}
 a_0(u,v) = \frac{1}{2}\big( a(u,v)+ a(v,u)\big)\quad\text{and}\quad b(u,v) = \frac{1}{2}\big( a(u,v)- a(v,u)\big)\,.
\end{align*}
Clearly, the bilinear forms $b(\cdot,\cdot), a_t(\cdot,\cdot)$ are bounded and $a_t(\cdot,\cdot)$ is $\theta$-coercive since for every $v\in V$ we have $a_t(v,v) = a_0(v,v) = a(v,v)$. For fixed $w\in V$, define
the bounded linear form $\ell_w(v) = \ell(v) -tb(w,v)$. Given that $a_0(\cdot,\cdot)$ is symmetric, 
from the previous case, there is a unique $u=u (w)\in K$ satisfying the variational inequality 
\begin{align*}
 a_0(u,v-u)\geq \ell_w(v-u)\qquad\text{for all}~~ v\in K.
\end{align*}
Accordingly, the map $T: V\to K$ with $u(w)=T w $ is well defined. Choosing $t$ such that $ 0\leq t\leq \frac{\theta}{2M}$ then utilizing the estimate \eqref{eq:well-posed} leads to
\begin{align*}
 &\|Tw- Tw'\|_{V}\leq\frac{1}{\theta} \|\ell_w-\ell_{w'}\|_{V'}\\
 &= \frac{t}{\theta} \|b(w-w',\cdot)\|_{V'}\leq \frac{M}{\theta}t \|w-w'\|_V\leq \frac{1}{2} \|w-w'\|_V\,.
\end{align*}
This shows that $T$ is a contraction on $V$ and thus has a unique fixed point $u_t\in K$. We have $u_t=Tu_t$ which by definition implies that for every $v\in K$,
\begin{align*}
a_t(u_t, v-u_t)&= a_0(u_t, v-u_t)+ tb(u_t, v-u_t)\\
&\geq \ell_{u_t}(v-u_t)+ tb(u_t, v-u_t)= \ell(v-u_t)\,.
\end{align*}
We have shown that for $ 0\leq t\leq \frac{\theta}{2M}$, there is a unique $u_t\in K$ such that
\begin{align*}
a_t(u_t, v-u_t)\geq \ell(v-u_t)\quad\text{for all}~~\in v\in K\,.
\end{align*}
\noindent A recursive argument shows that this remains true when $\frac{\theta n}{2M}\leq t\leq \frac{\theta(n+1)}{2M} $ for all $n\in \mathbb{N}$ and thus for all $t\geq 0$. This terminates the proof since for $t=1$ we have  $a_1(u,v) = a(u,v)$ for all $u,v \in V$.
 \end{proof}

 \bigskip
 
\noindent The celebrated Lax-Milgram lemma turns out to be a consequence of the Theorem \ref{thm:stampachia}. 
 
 \begin{corollary}[Lax-Milgram lemma]
 Let $V$ be a Banach space and $K\subset V$ be a nonempty closed, convex set. Let $a(\cdot, \cdot): V\times V\to \mathbb{R}$ be a continuous $\theta$-coercive bilinear form. Then for every continuous linear form $\ell:V\to \mathbb{R}$ there is a unique $u\in V$ such that 
 \begin{align*}
 a(u, v)= \ell(v)\quad \text{for all}~ v\in V\,.
 \end{align*}
 Moreover, the mapping $\ell \mapsto u$ is linear and continuous from $V'$ to $V$ with 
 \begin{align*}
 \|u\|_V\leq \frac{1}{\theta}\|\ell\|_{V'}\,.
 \end{align*}
 \end{corollary}
 
 \medskip
 
 \section{Lagrange multipliers}
 
 The Lagrange multiplier method is one of the most powerful tools used to solve certain types
 of constrained minimization problems in the setting of Banach spaces. For a fair exposition we need some basics on differential calculus. We recommend \cite{cheney2013} for supplementary notions on differential calculus. In this section, 
 $V,W$ are two Banach spaces and $U\subset V$ is an open set of $V$. We recall the notion of Fr\'echet derivative.
 %
 A function $f: U\to W $ is said to be Fr\'{e}chet differentiable at a point $a\in U$ if there is a linear bounded operator $L_a:V\to W$ i.e $L_a\in \mathcal{L}(V,W)$ such that 
 \begin{align*}
 \lim_{\|h\|_{V}\to 0}\frac{\|f(a+h)-f(a) -L_a(h)\|_{W}}{\|h\|_{V}}=0. 
 \end{align*}
 It is often common that the operator $L_a$ is synonymously denoted by d$f(a)$ or $f'(a)$. We will adopt the latter notation. 
 \noindent We say that $f$ is $C^1$ or continuously differentiable if $f'(a)$ exists for every $a\in U$ and the function $f':U\to \mathcal{L}(V,W)$ with $a\mapsto f'(a)$ is continuous. 
 \noindent Note that $ \mathcal{L}(V,W)$ is equipped with the norm
 \begin{align*}
 \|T\|_{ \mathcal{L}(V,W)}:=\sup\{\|Tv\|_{W}:~\|v\|_{V}=1\}. 
 \end{align*}
 
 \noindent If the Fr\'{e}chet derivative of $f$ at $a$ exists then
 \begin{align*}
 \lim_{t\to0}\frac{f(a+t v)-f(a)}{t}= f'(a)v\qquad\text{for all}~v\in V. 
 \end{align*}
In general, if the above expression exists then $f$ is said to be g\^{a}teau differentiable at $a$ in the direction $v$. However, the g\^{a}teau differentiability does not always imply the Fr\'{e}chet differentiability. The Fr\'{e}chet derivative is related to the open mapping theorem. 

\medskip

 \begin{theorem}[Open mapping Theorem]\cite[Section 3.4]{cheney2013}\label{thm:open-map}
 Let $f:U\to W$ be a continuously differentiable function where $U\subset V$ is open. 
 If the Fr\'{e}chet derivative $f'(a) $ of $f$ at point $a\in U$ is surjective, then $f(U)$ is a neighborhood of $f(a)$ in $W$.
 \end{theorem}

\vspace{2mm}

\begin{theorem}[Lagrange Multiplier]\cite[Section 3.5]{cheney2013}
Let $U\subset V$ be open and let $f: U\to \mathbb{R}$ and $g:V\to W$ be continuously differentiable functions. Let $K=\{x\in U:~g(x)=0\}$. Assume $a \in U\cap K $ is the a local minimum of $f$ on $U\cap K$ such that $g'(a)\in \mathcal{L}(V,W)$ is surjective. Then there exists a bounded linear form $\ell\in W'$ so that $f'(a)=\ell\circ g'(a) $. To be more precise, $f'(a)(v)=\ell\circ g'(a)(v)$ for all $v\in V. $
%
\end{theorem}

\noindent This theorem can be proved by applying the implicit function theorem, see \cite{cheney2013}. However, the Lagrange multiplier theorem above can be reformulated in the special case where $W= \mathbb{R}$ as follows.

\vspace{2mm}

\begin{theorem}[Lagrange Multiplier]\label{thm:lagrange-multiplier}
Let $U\subset V$ be open and let $f: U\to \mathbb{R}$ and $g:U\to W$ be continuously differentiable 
functions. Let $K=\{x\in U:~g(x)=0\}$. Assume $a \in U\cap K $ is the a local minimum of $f$ on 
$U\cap K$ such that $g'(a)\neq 0 $. Then there is $\lambda \in \mathbb{R}$ so that $f'(a)=\lambda g'(a)$. To be more precise, $f'(a)(v)=\lambda g'(a)(v)$ for all $v\in V. $
%
\end{theorem}

\vspace{2mm}

\begin{proof}
Without loss of generality, assume that $f(a) = \min\limits_{x\in U\cap K} f(x)$. Define $F:U\to \mathbb{R}^2$ with $F(x) = (f(x), g(x))$. Then $F(a)= (f(a) ,0)$ and $F$ is also $C^1$ and for 
$x\in U$, $F'(x)(v)=(f'(x)(v), g'(x)(v))$ for all $v\in V$. As $f(a)$ is the local minimum of $f$, it appears that for $\varepsilon>0$ small enough, $(f(a)-\varepsilon,0))$ is not an element of $F(U)$. This means that $F(U)$ is not a neighborhood of $F(a)$. 
Whence the function $F$ cannot be open which in virtue of the open mapping Theorem \ref{thm:open-map} means that $F'(a)$ is not surjective.
Accordingly, the range $\operatorname{Im} (F'(a))$ of $F'(a)$ is at most of dimension one. Hence there is a bounded linear form $ \ell \in V'$ and an element $\gamma=(\gamma_1, \gamma_2)\in \mathbb{R}^2 $ such that $F'(a)(v) =\gamma\ell(v) $ for all $v\in V$. We have $f'(a)(v) = \gamma_1\ell(v)$ and $g'(a)(v) = \gamma_2\ell(v)$. The fact that $g'(a)\neq 0$ implies that $\gamma_2\neq 0$. Letting $\lambda= \gamma_1/ \gamma_2$ we obtain $f'(a) = \lambda g'(a)$.
\end{proof}

\bigskip

\section{Integrodifferential equations (IDEs)}
Here we study the well-posedness of various types of integrodifferential equations (IDEs) associated with  symmetric integrodifferential operators of L\'evy type. These are operators of the form
\begin{align*}
Lu(x) = \pv \int_{\R^d} (u(x) - u(y)) \nu(x-y)\d y, \qquad\text{($x\in \R^d$)}
\end{align*}
where $\nu: \R^d\setminus\{0\}\to [0,\infty]$ is even and L\'evy integrable, i.e. 
\begin{align*}
\nu(h) = \nu(-h)\,\, \text{for $h\neq 0$}\quad\text{and}\quad\int_{\R^d} (1\land |h|^2)\nu(h)\d h<\infty.
\end{align*}
\noindent We are interested in showing the well-posedness of IDEs on a domain $\Omega\subset \R^d$ with Neumann, Dirichlet, Robin and mixed complement condition. 
It is worth emphasizing that due to the nonlocal feature of $L$, the condition for each of the aforementioned problems is imposed on the complement of $\Omega$.
In each case, we briefly recall the local analog. From now on we use several notations from the previous chapters. 
In particular, we associate with  the L\'evy density $\nu$,  the bilinear form $\mathcal{E}(\cdot, \cdot)$ defined for $u,v\in\VnuOm$ by 
\begin{align*}
\mathcal{E}(u,v)=\frac{1}{2} \iil_{(\Omega^c\times \Omega^c)^c} \big(u(x)-u(y) \big) \big(v(x)-v(y) \big) \, \nu (x-y) \mathrm{d}x \, \mathrm{d}y.
 	\end{align*}
In addition, the function spaces introduced in Chapter \ref{chap:nonlocal-sobolev} play an important role in this section. It is worth recalling that the space $\VnuOm$ is always endowed with the norm $\|v\|^2_{\VnuOm}= \|v\|^2_{L^2(\Omega)} + \mathcal{E}(u,u)$. Also recall that $|\cdot|^2_{\VnuOm}\leq \mathcal{E}(\cdot, \cdot )\leq 2|\cdot|^2_{\VnuOm}$, where 
\begin{align*}
|v|^2_{\VnuOm}=\iil_{\Omega\R^d} \big(v(x)-v(y) \big)^2 \, \nu (x-y) \mathrm{d}x \, \mathrm{d}y.
\end{align*}

\subsection{Integrodifferential equations (IDEs) with Neumann condition}

Assume $\Omega\subset \mathbb{R}^d$ is an open set. Let $f:\Omega\to \mathbb{R}$ and $g: \mathbb{R}^d\setminus \Omega\to \mathbb{R}$ be measurable functions. The Neumann problem for the operator $L$ is to find a measurable  function $u:\mathbb{R}^d\to \mathbb{R}$ such that 
\begin{align}\label{eq:nonlocal-Neumann}\tag{$N$}
L u = f \quad\text{in}~~~ \Omega \quad\quad\text{ and } \quad\quad \mathcal{N} u= g ~~~ \text{on}~~~ \mathbb{R}^d\setminus\Omega, 
\end{align} 
where $\mathcal{N}$, also called for an obvious reason the \emph{nonlocal normal derivative operator} (see \cite{DROV17, DGLZ12}) on $\Omega$ with respect to $\nu$ is the integrodifferential operator defined by 
\begin{align}\label{eq:nonlocal-derivative}
\mathcal{N}u(y) = \int_{\Omega}(u(x)-u(y))\nu(x-y)\,\mathrm{d}x\qquad\qquad (y\in \mathbb{R}^d\setminus\Omega). 
\end{align}

\vspace{2mm}

\noindent Let us derive the so called \textit{nonlocal Green-Gauss formula} cf. \eqref{eq:green-gauss-nonlocal}  which provides a nonlocal version to the classical Green-Gauss formula  for $ u\in H^{2}(\Omega) $ and $v\in H^{1}(\Omega)$ (See \cite[Theorem 2.20]{FSU19}),
\begin{align}\label{eq:green-Gauss}
\int_{\Omega} (-\Delta) u(x) v(x) \, \mathrm{d}x = \int_{\Omega} \nabla u(x) \cdot \nabla v(x) \, \mathrm{d}x- \int_{\partial \Omega} \gamma_{1} u(x) \gamma_{0}v (x)\, \mathrm{d}\sigma(x), 
\end{align}
 where we recall that  $\gamma_1= \gamma_0 \circ\nabla $, and $\gamma_0$ stands for the trace operator on $H^1(\Omega)$. 
 
\medskip

\begin{proposition}[\textbf{Green-Gauss formula}]
Let $\Omega\subset \R^d$ be open and bounded. For every $u\in C_b^2(\R^d)$  and $v\in C_b^1(\R^d)$  the following Green-Gauss formula holds true
	%
	\begin{align}\label{eq:green-gauss-nonlocal}
	\int_{\Omega} [Lu(x)]v(x)\mathrm{d}x= \mathcal{E}(u,v) -\int_{\Omega^c}\mathcal{N}u(y)v(y)\mathrm{d}y.
	\end{align}
	In particular, letting $v=1 $ one gets the integration by part formula, 
	\begin{align}\label{eq:integration-by-part-nonlocal}
	\int_{\Omega} Lu(x)\mathrm{d}x=-\int_{\Omega^c}\mathcal{N}u(y)\mathrm{d}y.
	\end{align}
Furthermore, if $u, v\in C_b^2(\R^d)$ then we have the following second Green-Gauss formula, 
	\begin{align}\label{eq:green-gauss-nonlocal-second}
\int_{\Omega} [Lu(x)]v(x)- [L v(x)]u(x)\mathrm{d}x=  \int_{\Omega^c}[\mathcal{N}v(y)]u(y)-[\mathcal{N}u(y)]v(y)\mathrm{d}y.
\end{align}
\end{proposition}

\medskip

\begin{proof} It is sufficient to prove \eqref{eq:green-gauss-nonlocal} since it implies \eqref{eq:green-gauss-nonlocal-second}.  Note that for $\varphi\in C_b^1(\R^d)$, we have 
	\begin{align}\label{eq:first-order-diff}
	|\varphi(x)- \varphi(x+h)|\leq 2\|\varphi\|_{C_b^1(\R^d)}(1\land |h|)\quad \text{for all $x,h\in \R^d$}. 
	\end{align}
	Let $u \in C_b^{2}(\mathbb{R}^d)$ and $v \in C_b^{1}(\mathbb{R}^d)$. With the aid of Proposition \ref{prop:uniform-cont} we can write 
	\begin{align*}
	&\int_{\Omega} [Lu(x)] u(x)v(x)\mathrm{d}x =\lim_{\varepsilon\to 0} \int\limits_{\Omega} v(x) \mathrm{d}x 
	\int\limits_{\mathbb{R}^d\setminus B_\varepsilon(x)} ((u(x)-u(y))\nu(x-y)\,\mathrm{d}y\\
	&= \lim_{\varepsilon\to 0} \int\limits_{\Omega} \int\limits_{\Omega \setminus B_\varepsilon(x)}(u(x)-u(y))v(x)\nu(x-y)\,\mathrm{d}y \mathrm{d}x + \int\limits_{\Omega} \int\limits_{\Omega^c }(u(x)-u(y))v(x)\nu(x-y)\,\mathrm{d}y \mathrm{d}x 
	\end{align*}
On one side, by a symmetry argument we have 
	\begin{align*}
	&\qquad \lim_{\varepsilon\to 0} \int\limits_{\Omega} \int\limits_{\Omega \setminus B_\varepsilon(x)}(u(x)-u(y))v(x)\nu(x-y)\,\mathrm{d}y \mathrm{d}x
	= \lim_{\varepsilon\to 0} \hspace{-2ex}\iint\limits_{\Omega \times \Omega\cap \{|x-y|>\varepsilon\}}\hspace*{-3ex}(u(x)-u(y))v(x)\nu(x-y)\,\mathrm{d}y \mathrm{d}x \\
	&=\hspace{-0.2ex} \lim_{\varepsilon\to 0} \frac{1}{2}\hspace{-4ex} \iint\limits_{\Omega \times \Omega\cap \{|x-y|>\varepsilon\}} \hspace{-4ex} (u(x)-u(y))(v(x)-v(y))\nu(x-y)\,\mathrm{d}y \mathrm{d}x =\hspace{-0.5ex}
	\frac{1}{2} \iint\limits_{\Omega \Omega } (u(x)-u(y))(v(x)-v(y))\nu(x-y)\,\mathrm{d}y \mathrm{d}x
	\end{align*}
	where one gets rid of the principal value using the estimate \eqref{eq:first-order-diff} applied to $u$ and $v$. 
	On the other side, with the help of Fubini's theorem we have
	\begin{align*}
	& \iint\limits_{\Omega \Omega^c}(u(x)-u(y))v(x)\nu(x-y)\,\mathrm{d}y \mathrm{d}x \\
	&= \iint\limits_{\Omega \Omega^c} (u(x)-u(y))(v(x)-v(y))\nu(x-y)\,\mathrm{d}y \mathrm{d}x
	+ \int\limits_{\Omega^c} v(y)\mathrm{d}y \int\limits_{\Omega }(u(x)-u(y))\nu(x-y)\, \mathrm{d}x\\
	&= \frac{1}{2}\iint\limits_{\Omega \Omega^c} (u(x)-u(y))(v(x)-v(y))\nu(x-y)\,\mathrm{d}y\, \mathrm{d}x
	+ \frac{1}{2}\iint\limits_{\Omega^c \Omega}(u(x)-u(y))(v(x)-v(y))\nu(x-y)\,\mathrm{d}y \mathrm{d}x\\
	& \qquad\qquad- \int_{\Omega^c}\mathcal{N} u(y) v(y)\mathrm{d}y \, .
	\end{align*}
	Altogether inserted in the initial relation provide the desired relation.
\end{proof}
\vspace{2mm}

\noindent Regarding the density of $C_c^\infty(\R^d)$ in  $V^1_{\nu}(\Omega|\R^d)$ and in $ \VnuOm$  (see Theorem \ref{thm:density}),  the nonlocal version of Green-Gauss formula \eqref{eq:green-Gauss} is given as follows.
\begin{theorem}
Assume $\Omega\subset \R^d$ is open and bounded with Lipschitz boundary.  Recall the definition of  the space $V^1_{\nu}(\Omega|\R^d)$ see Section \ref{sec:function-spaces}, for every $u \in V^1_{\nu}(\Omega|\R^d)$ and every $v \in \VnuOm$ we have 
\begin{align*}
\int_{\Omega} [Lu(x)]v(x)\mathrm{d}x= \mathcal{E}(u,v) -\int_{\Omega^c}\mathcal{N}u(y)v(y)\mathrm{d}y
\end{align*}
where, the operators $L$ and $\mathcal{N}$ are understood in the generalized sense. 
\end{theorem}

\medskip

\noindent Let us look at a connection between the trace space $\TnuOm$ and  the nonlocal normal derivative $\mathcal{N}$. 

\begin{theorem}\label{thm:linear-form-characto}
	Assume $\TnuOm$ is endowed with the norm $\|\cdot\|_{ \TnuOm}$. 
	Then for any linear continuous form $\ell :\TnuOm\to \mathbb{R}$ there exists $w\in \VnuOm$ such that for every $v\in C^\infty_c(\overline{\Omega}^c)$ 
	\begin{align}
	\ell(v) = \int_{\Omega^c} \mathcal{N} w(y)v(y)\mathrm{d}y
	\end{align}
	In particular, given a measurable function $g:\Omega^c \to \mathbb{R}$, if the linear mapping $ \ell_g: v\mapsto \int_{\Omega^c} g(y)v(y)\mathrm{d}y$ is continuous on $\TnuOm $, then there exists $w\in \VnuOm $ such that $g= \mathcal{N} w$ almost everywhere on $\Omega^c$.
\end{theorem}

\medskip

\begin{proof}
Let $\ell \in (\TnuOm)'$. Because of the continuity of the trace operator $\operatorname{Tr}: \VnuOm\to \TnuOm$, the linear form $\ell \circ \operatorname{Tr}$ is also continuous on $ \VnuOm $. By Riesz's representation theorem there exists $w \in \VnuOm$ such that $\ell \circ \operatorname{Tr}(v) = \left(v, w\right)_{( \VnuOm } $ for each $v \in \VnuOm$. In particular, for $v\in C^\infty_c(\overline{\Omega}^c)$ identified with its zero extension on $\Omega$ so that $\operatorname{Tr}(v) =v$, we obtain the following 
	
	\begin{align*}
	\ell(v) &= \int_{\Omega} w(x)v(x) \mathrm{d}x + \iint\limits_{(\Omega^c\times\Omega^c)^c} (w(x)-w(y))(v(x)-v(y)) \nu(x-y)\,\mathrm{d}x\,\mathrm{d}y\\ 
	%
	%
	&= \int\limits_{\Omega^c}v(y)\mathrm{d}y \int\limits_{\Omega} (w(y)-w(x))\nu(x-y)\,\mathrm{d}x= \int_{\Omega^c} \mathcal{N} w(y)v(y)\mathrm{d}y.
	\end{align*} 
	Furthermore, if $g: \Omega^c\to \mathbb{R}$ is such that $\ell_g$ is continuous on $\TnuOm$ then by the above computation, it follows that $g=\mathcal{N} w$ almost everywhere on $\Omega^c$ since 
\begin{align*}
	\int_{\Omega^c} g(y)v(y)\d y= \int_{\Omega^c} \mathcal{N} w(y)v(y)\mathrm{d}y,\quad \hbox{for all }~~v\in C_c^\infty(\overline{\Omega}^c).
	\end{align*}
\end{proof}

\vspace{2mm}

\begin{remark}
	The second statement of Theorem \ref{thm:linear-form-characto} particularly suggests that the space of all measurable functions $g:\Omega^c\to \mathbb{R}$ for which linear the form $v\mapsto\int_{\Omega^c}g(y)v(y)\d y$ is continuous on $\TnuOm$ is contained in $\mathcal{N}(\VnuOm)$ (the range of $\mathcal{N}$). 
\end{remark}

\medskip

\noindent In light of the relation \eqref{eq:green-gauss-nonlocal} it is reasonable to define weak solutions of the Neumann problem as follows. 
\begin{definition}
	A measurable function $u:\mathbb{R}^d\to \mathbb{R}$ is a weak solution or a variational solution of the inhomogeneous Neumann problem \eqref{eq:nonlocal-Neumann} if $u \in \VnuOm$ and satisfies the relation
	\begin{align}\label{eq:var-nonlocal-Neumann}\tag{$V$}
	\mathcal{E}(u,v) = \int_{\Omega} f(x)v(x)\mathrm{d}x +\int_{\Omega^c} g(y)v(y)\mathrm{d}y,\quad \mbox{for all}~~v \in \VnuOm\,.
	\end{align}
	
	\noindent In particular, if $\Omega$ is bounded then taking $v=1$, \eqref{eq:var-nonlocal-Neumann} turns to the so called compatibility condition
	\begin{align}\label{eq:compatible-nonlocal}\tag{$C$}
	\int_{\Omega} f(x)\mathrm{d}x +\int_{\Omega^c} g(y)\mathrm{d}y=0.
	\end{align}
\end{definition} 

\medskip

\begin{remark}
	The compatibility condition \eqref{eq:compatible-nonlocal} is an implicit necessary requirement that the data $f$ and $g$ must fulfill before any attempt at solving the problems \eqref{eq:var-nonlocal-Neumann} and \eqref{eq:nonlocal-Neumann}. The local counterpart of this compatibility condition, where $g$ is defined on $\partial\Omega$, is given by 
	\begin{align}
	\int_{\Omega} f(x)\mathrm{d}x +\int_{\partial\Omega} g(y)\mathrm{d}\sigma(y)=0.
	\end{align} 
	Let us recall that the variational formulation of the classical Neumann problem $-\Delta u=f$ in $\Omega$ and $\frac{\partial u}{\partial n} =g$ on $\partial \Omega$ is to find $u\in H^1(\Omega)$ such that
	\begin{align}
	\int_{\Omega} \nabla u(x)\cdot \nabla v(x) \mathrm{d}x = \int_{\Omega} f(x)v(x)\mathrm{d}x +\int_{\partial\Omega} g(y)v(y)\mathrm{d}\sigma(y),\quad \mbox{for all}~~v \in H^1(\Omega)\,.
	\end{align} 
\end{remark}

%

\noindent Both integrodifferential operators $L$ and $\mathcal{N}$ annihilate additive constants. 
Therefore, as long as $u$ is a solution to the system \eqref{eq:nonlocal-Neumann} or to the variational problem \eqref{eq:var-nonlocal-Neumann} so does the function $\widetilde{u}= u+c$ for any $c\in \mathbb{R}$. Accordingly, both problems are ill-posed in the sense of Hadamard. The situation is likewise in the local setting with the operators $L$ and $\mathcal{N}$ respectively replaced by the operators $-\Delta$ and $\frac{\partial}{\partial n}$. In order to overcome this anomaly we introduce an appropriate functional space $ \VnuOm^{\perp}$ consisting of functions in $ \VnuOm$ with zero mean over $\Omega$. To be more precise, 
\begin{align*}
\VnuOm^{\perp}:= \Big\{ u\in \VnuOm: \int_{\Omega}u(x)\mathrm{d}x=0\Big\}.
\end{align*}
Assuming that $\Omega$ is bounded, the space $\VnuOm^{\perp}$ endowed with the scalar product of $\VnuOm$ is also a Hilbert as well. With this at hand, we rewrite the variational problem  \eqref{eq:var-nonlocal-Neumann} as  follows
\begin{align}\label{eq:var-nonlocal-Neumann-bis}\tag{$V'$}
\mathcal{E}(u,v) = \int_{\Omega} f(x)v(x)\mathrm{d}x +\int_{\Omega^c} g(y)v(y)\mathrm{d}y,\quad \mbox{for all}~~v \in \VnuOm^{\perp}\,.
\end{align}
It is noteworthy emphasizing that in contrast to \eqref{eq:var-nonlocal-Neumann}, the variational problem \eqref{eq:var-nonlocal-Neumann-bis} possesses at most one solution since $\mathcal{E}(\cdot, \cdot)$ defines a scalar product on $\VnuOm^\perp$. Analogous observations can be carried out in the local setting by introducing the space $H^1(\Omega)^\perp=\big\{ u\in H^1(\Omega): \int_{\Omega}u(x)\mathrm{d}x=0\big\}.$

\medskip 

\noindent Under additional regularity assumptions, we show that both problems \eqref{eq:nonlocal-Neumann} and \eqref{eq:var-nonlocal-Neumann} are equivalent.
\begin{proposition}
	Let $\Omega$ be an open bounded subset of $\mathbb{R}^d$ with Lipschitz boundary. Let $u\in C^2_b(\mathbb{R}^d)$, $ f\in L^{2}(\Omega)$ and $g\in L^2(\Omega^c, \nu^{-1}_K)$ where $\nu_K (x)=\operatorname{essinf}_{y\in K} \nu(x-y)$ for some measurable  set $K\subset \Omega$ with $|K|>0$. Then $u$ is a solution of \eqref{eq:nonlocal-Neumann} if and only if $f$ and $g$
	are compatible in the sense of \eqref{eq:compatible-nonlocal} and $u$ is a solution of \eqref{eq:var-nonlocal-Neumann}. 
	
\end{proposition}
\medskip

\begin{proof} If $u$ solves \eqref{eq:nonlocal-Neumann}, i.e. $Lu=f$ in $\Omega$ and $\mathcal{N}u =g$ on $\Omega^c$, then by the Green-Gauss formula \eqref{eq:green-gauss-nonlocal} we obtain the following 
	\begin{align}\label{eq:var-Neumann-regular}
	\mathcal{E}(u,v)= \int_{\Omega} f(x)v(x)\mathrm{d}x + \int_{\Omega^c} g(y)v(y)\mathrm{d}y, \quad \mbox{for all}~~v \in C^1_b(\mathbb{R}^d).
	\end{align}
	As shown in \eqref{eq:linearform-f}-\eqref{eq:linearform-g} below, all terms involved in \eqref{eq:var-Neumann-regular} are linear and continuous on $\VnuOm$ with respect to the variable $v$. Moreover, smooth functions of compact support are dense in $\VnuOm$ hence the relation in \eqref{eq:var-Neumann-regular} remains true for functions $v$ in $\VnuOm$ so \eqref{eq:var-nonlocal-Neumann} is satisfied. In particular, taking $v=1$ yields the condition \eqref{eq:compatible-nonlocal}. Conversely, assume $u$ solves \eqref{eq:var-nonlocal-Neumann} then inserting the Green-Gauss formula \eqref{eq:green-gauss-nonlocal} with $v\in C^1_b(\mathbb{R}^d)\subset \VnuOm$ in \eqref{eq:var-Neumann-regular} yields 
	%
	\begin{align*}
		\int_{\Omega} L u(x) v(x) \mathrm{d}x - \int_{\Omega} f(x)v(x)\mathrm{d}x = \int_{\Omega^c} g(y)v(y)\mathrm{d}y-\int_{\Omega^c}\mathcal{N} u(y)v(y)\mathrm{d}y, \quad \mbox{for all}~~v \in C^1_b(\mathbb{R}^d).
		\end{align*}
	%
	Specializing this relation for $ v \in C^\infty_c(\Omega)$ and $v \in C^\infty_c(\mathbb{R}^d\setminus \overline{\Omega})$ respectively, we end up with 
	\begin{align*}
	&\int_{\Omega} L u(x) v(x) \mathrm{d}x - \int_{\Omega} f(x)v(x)\mathrm{d}x = 0 \qquad \mbox{for all}~~v \in C^\infty_c(\Omega),
	\\
	&\int_{\Omega^c} g(y)v(y)\mathrm{d}y-\int_{\Omega^c}\mathcal{N} u(y)v(y)\mathrm{d}y= 0 \qquad \mbox{for all}~~v \in C^\infty_c(\mathbb{R}^d\setminus \overline{\Omega}). 
	\end{align*}
According to Proposition \ref{prop:uniform-cont}, $Lu$ is well defined and bounded hence belongs to $L^2(\Omega)$. Similarly $\mathcal{N}u$ is well defined and bounded, i.e. belongs to $L^\infty(\Omega^c)$. Thus, up to null sets, the above equations lead to \eqref{eq:nonlocal-Neumann}. Precisely, we have $L u = f \,\,\text{in}~ \Omega ~\text{and}~\mathcal{N} u= g ~\text{on}~ \mathbb{R}^d\setminus\Omega.$
\end{proof}

\vspace{2mm}

\noindent By standard procedure, a solution of the variational problem \eqref{eq:var-nonlocal-Neumann} is characterized as a critical point (a minimizer) of the functional 
\begin{align}
\begin{split}
\mathcal{J}(v) &= \frac{1}{2} \mathcal{E}(v,v) - \int_\Omega f v\, \mathrm{d}x- \int_{\Omega^c} g v \mathrm{d}x\\
&= \frac{1}{4}\iint\limits_{(\Omega^c\times\Omega^c)^c} (v(x)-v(y))^2\nu(x-y)\,\mathrm{d}x\,\mathrm{d}y - \int_\Omega f v\, \mathrm{d}x- \int_{\Omega^c} g v \mathrm{d}x.
\end{split}
\end{align}
\begin{proposition}
	A function $ u\in\VnuOm^{\perp} $ is solution to \eqref{eq:var-nonlocal-Neumann-bis} if and only if $u$ is also a solution of the minimization problem
	\begin{align}\label{eq:nonlocal-Neumann-min}\tag{$M'$}
	\mathcal{J}(u) =\min_{v\in \VnuOm^{\perp}}\mathcal{J}(v)
	\end{align}
	Moreover, $ u\in\VnuOm^{\perp} $ solves \eqref{eq:var-nonlocal-Neumann-bis} if and only if for any $c\in \mathbb{R}$, $u+c$ solves the variational problem \eqref{eq:var-nonlocal-Neumann} and the latter problem is equivalent to the minimization problem 
	\begin{align}\label{eq:nonlocal-Neumann-min-const}\tag{$M$}
	\mathcal{J}(u) =\min_{v\in \VnuOm}\mathcal{J}(v).
	\end{align}
\end{proposition} 

\vspace{1mm}

\begin{proof}
	Let $ u,v \in\VnuOm^{\perp}$.  Employing, Cauchy-Schwartz inequality yields
	\begin{align*}
	\mathcal{E}(u,v) 
	%
	&\leq \frac{1}{2} \mathcal{E}(u,u)+ \frac{1}{2} \mathcal{E}(v,v)= \mathcal{E}(u,u)- \frac{1}{2} \mathcal{E}(u,u)+ \frac{1}{2} \mathcal{E}(v,v).
	\end{align*} 
Thus, if \eqref{eq:var-nonlocal-Neumann-bis} holds true for all $ v\in\VnuOm^{\perp}$ then we get $\mathcal{J}(u)\leq \mathcal{J}(v)$ and thus $u $ solves \eqref{eq:nonlocal-Neumann-min}. 
Conversely, assume that $ u\in \VnuOm^\perp$ satisfies \eqref{eq:nonlocal-Neumann-min}.  For all $v\in \VnuOm^\perp$ and all $t\in \R$ \eqref{eq:nonlocal-Neumann-min} implies  $\mathcal{J}(u)\leq \mathcal{J}(u+tv)$ since  $u+tv \in \VnuOm^\perp$. Therefore, for fixed $ v\in\VnuOm^\perp$ the polynomial of second order $\mathcal{J}(u+\centerdot v): \mathbb{R} \to \mathbb{R}$, 
	\begin{align*}
	t\mapsto \mathcal{J}(u+tv)= \mathcal{J}(u)+ t\Big[\mathcal{E}(u,v)- \int_{\Omega} f(x)v(x)\mathrm{d}x\Big] + \frac{t^2}{2} \mathcal{E}(v,v) 
	\end{align*}
	has a critical point at $t=0$. From this we get that \eqref{eq:var-nonlocal-Neumann-bis} is verified since 
	\begin{align*}
	0= \lim_{t\to 0}\frac{\mathcal{J}(u+tv)-\mathcal{J}(u)}{t} = \mathcal{E}(u,v)- \int_{\Omega} f(x)v(x)\mathrm{d}x-\int_{\Omega^c} g(y)v(y)\mathrm{d}y .
	\end{align*}
\noindent Meanwhile, if the compatibility condition \eqref{eq:compatible-nonlocal} holds, then it is easy to observe that the 
relation in \eqref{eq:var-nonlocal-Neumann-bis} remains unchanged under additive constant, i.e. $\mathcal{J}(v+c)=\mathcal{J}(v)$ for all $v \in\VnuOm$ and all $c\in \mathbb{R}.$ Accordingly, if $u\in \VnuOm^\perp$ solves \eqref{eq:var-nonlocal-Neumann-bis}, then we have $\mathcal{J}(u+c) =\min\limits_{v\in \VnuOm}\mathcal{J}(v)$ which by a similar arguments as above is equivalent to \eqref{eq:var-nonlocal-Neumann}. 
\end{proof}

\medskip

\noindent We are now in a position to state the well-posedness  of \eqref{eq:var-nonlocal-Neumann-bis} and hence of \eqref{eq:var-nonlocal-Neumann} up to additive constant. 

\medskip

\begin{theorem} \label{thm:nonlocal-Neumann-var}
	Assume that $\Omega\subset \mathbb{R}^d$ is a bounded open set and the function $\nu: \mathbb{R}^d\to [0,\infty]$ is the density of a symmetric L\'{e}vy measure with full support such that the couple $(\nu, \Omega)$ belongs to one of the classes $\mathscr{A}_i,~i=1,2,3,4$. For 	a set   $K\subset \Omega$ with positive measure, assume $ \nu_K(x)= \operatorname{essinf}_{y\in K}~\nu(x-y)>0~a.e$ for almost $x\in \Omega^c$. 
	Then, given $ f \in L^2 (\Omega)$ and $g\in L^2(\Omega^c, \nu_K^{-1})$, there exists a unique solution $u\in\VnuOm^{\perp} $ to the variational problem \eqref{eq:var-nonlocal-Neumann-bis}. Further, the solutions to \eqref{eq:var-nonlocal-Neumann} are of the form $u+c$ with $c\in \mathbb{R}$ provided that $f$ and $g$ are compatible in the sense that \eqref{eq:compatible-nonlocal} is verified. 
	
	\noindent Moreover, there exists a constant $C: = C(d,\Omega, \nu)>0$ independent of $f$ and $g$ such that any solution $w$ of $\eqref{eq:var-nonlocal-Neumann}$ satisfies the following weak regularity estimate
	\begin{align}\label{eq:weak-regular}
	\|w-\hbox{$\fint_{\Omega}w$} \|_{\VnuOm}\leq C \left(\|f\|_{L^2(\Omega)}+\|g\|_{L^{2}(\Omega^c, \nu_K^{-1})}\right)\,.
	\end{align}
	
	\noindent In particular, the operator $\Phi : L^{2}(\Omega) \times L^{2}(\Omega^c, \nu_K^{-1})\to\VnuOm^\perp $ mapping the Neumann data $(f,g) \in L^{2}(\Omega)\times L^{2}(\Omega^c, \nu_K^{-1})$ 
	to the unique solution $ u\in\VnuOm^\perp $ of the variational problem
	\eqref{eq:var-nonlocal-Neumann-bis} is linear, one-to-one, continuous and we have
	\begin{align*}
	\|\Phi(f,g)\|_{\VnuOm}\leq C \|(f,g)\|_{L^{2}(\Omega)\times L^{2}(\Omega^c, \nu_K^{-1})}.
	\end{align*} 
\end{theorem}

\medskip

\begin{remark}
	In some sense, a solution  to the variational problem \eqref{eq:var-nonlocal-Neumann}  exists only if the data $f$ and $g$ satisfy the compatibility condition \eqref{eq:compatible-nonlocal}. 
	This constraint corresponds to the situation arising in finite dimension when solving linear equations $Ax=b$ with $b\in \mathbb{R}^d$ and $A\in \mathbb{R}^{d\times d}$ where a unique solution exists if and only if $ \d im(\ker A)=0$. Recall that $\d im(\ker A)+ \d im(\operatorname{Im} A)=d\,.$ This generalizes in infinite dimensional spaces via the so called Fredholm alternative (see Theorem \ref{thm:fredhoml-alternative}). 
\end{remark}

\noindent The next theorem offers an alternative formulation of Theorem \ref{thm:nonlocal-Neumann-var} with  a relaxed condition on $\nu_K$ and a different (possibly larger) function space for $g$. 

\medskip

\begin{theorem} \label{thm:nonlocal-Neumann-var-weighted}
	Assume that $\Omega\subset \mathbb{R}^d$ is a bounded open set and the function $\nu: \mathbb{R}^d\to [0,\infty]$ is the density of a symmetric L\'{e}vy measure with full support such that the couple $(\nu, \Omega)$ belongs to one of the classes $\mathscr{A}_i,~i=1,2,3,4$ (see page  \pageref{eq:class-lipschitz}). For some $K\subset \Omega$ with positive measure, consider $ \nu_K(x)= \operatorname{essinf}_{y\in K}~\nu(x-y)$ for almost all $x\in \Omega^c$. 
	
	\smallskip
	
	\noindent Given $ f \in L^2 (\Omega)$ and $g\in L^2(\Omega^c, \nu_K)$, there exists a unique solution $u_{*}\in\VnuOm^{\perp} $ to the following variational problem \eqref{eq:var-nonlocal-Neumann-bis-weigthed} 
	\begin{align}\label{eq:var-nonlocal-Neumann-bis-weigthed}\tag{$V_{*}'$}
	\mathcal{E}(u_{*}, v) = \int_{\Omega} f(x)v(x)\mathrm{d}x +\int_{\Omega^c} g(y)v(y)\nu_K(y)\mathrm{d}y,\quad \mbox{for all}~~v \in \VnuOm^{\perp}\,.
	\end{align}
	Additionally, if $f$ and $g$ verify the condition \eqref{eq:compatible-nonlocal-weighted},  then all solutions to the problem \eqref{eq:var-nonlocal-Neumann-weigthed} are of the form $u_{*}+c$ with $c\in \mathbb{R}$, where we let
	\begin{align}\tag{$C_{*}$}\label{eq:compatible-nonlocal-weighted}
	\int_{\Omega} f(x)\mathrm{d}x +\int_{\Omega^c} g(y)\nu_K(y)\mathrm{d}y=0
	\end{align}
	\begin{align}\label{eq:var-nonlocal-Neumann-weigthed}\tag{$V_{*}$}
	\mathcal{E}(u,v) = \int_{\Omega} f(x)v(x)\mathrm{d}x +\int_{\Omega^c} g(y)v(y)\nu_K(y)\mathrm{d}y,\quad \mbox{for all}~~v \in \VnuOm\,.
	\end{align}
\noindent Moreover, there exists a constant $C: = C(d,\Omega, K, \nu)>0$ independent of $f$ and $g$ such that any solution $w$ of \eqref{eq:var-nonlocal-Neumann-weigthed} satisfies the following weak regularity estimate
	
	\begin{align}\label{eq:weak-regular-weighted}
	\|w-\hbox{$\fint_{\Omega}w$} \|_{\VnuOm}\leq C \left(\|f\|_{L^2(\Omega)}+\|g\|_{L^{2}(\Omega^c, \nu_K)}\right)\,.
	\end{align}
	In particular, the operator $\Phi_{*}: L^{2}(\Omega) \times L^{2}(\Omega^c, \nu_K)\to\VnuOm^\perp $ mapping the Neumann data $(f,g) \in L^{2}(\Omega)\times L^{2}(\Omega^c, \nu_K)$ 
	to the unique solution $ u\in\VnuOm^\perp $ of the variational problem
	\eqref{eq:var-nonlocal-Neumann-weigthed} is linear, one-to-one, continuous and we have
	\begin{align*}
	\|\Phi_{*}(f,g)\|_{\VnuOm}\leq C \|(f,g)\|_{L^{2}(\Omega)\times L^{2}(\Omega^c, \nu_K)}\,.
	\end{align*} 
\end{theorem}

\vspace{1mm}

\begin{remark}
	$(i)$ It is worthwhile noticing that Theorem \ref{thm:nonlocal-Neumann-var} and  Theorem \ref{thm:nonlocal-Neumann-var-weighted} 
	remain true with the weight $\nu_K$ replaced by $\mathring{\nu}_K$, where we recall that 
	\begin{align*}
	\mathring{\nu}_K(x) = \int_K 1\land \nu(x-y)\d y.
	\end{align*}
	The reason is that,  according to  Theorem \ref{thm:Vnu-in-l2}, the embedding $\VnuOm \hookrightarrow L^2(\Omega^c, \mathring{\nu}_K)$ is also continuous. 
	
	\vspace{1mm}
	
	\noindent $(ii)$ Let $f\in L^2(\Omega)$, for $g=0$ it is worthwhile to see that the variational problem \eqref{eq:var-nonlocal-Neumann} coincides with \eqref{eq:var-nonlocal-Neumann-weigthed} and both correspond to the variational (weak) formulation of the homogeneous Neumann problem $Lu=f$ in $\Omega$ and $\mathcal{N}u=0$ on $\Omega^c$. 
\end{remark}

\vspace{1mm}

\begin{proof}[Proof of Theorem \ref{thm:nonlocal-Neumann-var} and 
	Theorem \ref{thm:nonlocal-Neumann-var-weighted}.] 
	The existence and the uniqueness of solutions of \eqref{eq:var-nonlocal-Neumann-bis} and 
	\eqref{eq:var-nonlocal-Neumann-bis-weigthed} on $\VnuOm^{\perp}$ plainly spring from the Lax-Milgram lemma using the arguments below. The bilinear form $\mathcal{E}(\cdot, \cdot)$ is continuous on $\VnuOm^\perp$. From the Poincar\'{e} inequality \eqref{eq:poincare-inequalityV} (Theorem \ref{thm:poincare-inequality}), for some constant $C>0$ we have 
	\begin{align*}
	\|v\|_{L^2(\Omega)}^2\leq C \mathcal{E}(v,v) \quad\quad \text{for all }~~~v\in \VnuOm^{\perp}\,.
	\end{align*}
	Equivalently, this inequality amounts to the coercivity of $\mathcal{E}(\cdot, \cdot)$ on $\VnuOm^\perp$ and we obtain
	\begin{align}\label{eq:coercivity-quadratic}
	\mathcal{E}(v,v)\geq \big(1+ C\big)^{-1}\|v\|^2_{\VnuOm}.
	\end{align}

	Let us show the continuity of the linear forms involved. For every  $v \in\VnuOm^{\perp}$ we have
	\begin{align}\label{eq:linearform-f}
	\Big|\int_{\Omega} f(x)v(x)\mathrm{d}x \Big|\leq \|f\|_{L^{2}(\Omega)}\|v\|_{L^{2}(\Omega)} \leq \|f\|_{L^{2}(\Omega)}\|v\|_{\VnuOm}\,.
	\end{align}
If $g\in L^2(\Omega^c, \nu_K^{-1})$ along with the continuity of $ \operatorname{Tr}: \VnuOm \hookrightarrow L^2(\Omega^c, \nu_K)$  (see Theorem \ref{thm:Vnu-in-l2}) then,
	\begin{align}\label{eq:linearform-g}
	\Big|\int_{\Omega^c} g(x) v(x)\mathrm{d}x \Big|\leq \|g\|_{ L^2(\Omega^c, \nu_K^{-1})}\|v\|_{L^2(\Omega^c, \nu_K)} \leq \|g\|_{ L^2(\Omega^c, \nu_K^{-1})}\|v\|_{\VnuOm}.
	\end{align}
	Similarly, if $g\in L^2(\Omega^c, \nu_K)$ we get
	\begin{align}\label{eq:linearform-weighted-g}
	\Big|\int_{\Omega^c} g(x) v(x)\nu_K(x)\mathrm{d}x \Big| \leq \|g\|_{L^2(\Omega^c, \nu_K)}\|v\|_{\VnuOm}.
	\end{align}	
	This conclude the existence and the uniqueness of a solution  $u\in \VnuOm^\perp$ (resp. $u_*\in \VnuOm^\perp$)  to the problem \eqref{eq:var-nonlocal-Neumann-bis} (resp. \eqref{eq:var-nonlocal-Neumann-bis-weigthed}).  
	
\vspace{1mm}	
	\noindent For $v\in\VnuOm$ set $v= \widetilde{v}+c'$ with $c'=\mbox{$ \fint_{\Omega} v\mathrm{d}x$}$ so that $ \widetilde{v} \in\VnuOm^\perp$. 
	In addition, a constant function $w=c$ belongs to $\VnuOm$ for every $c\in \mathbb{R}^d$ because $\Omega$ is bounded. This means that $\VnuOm= \VnuOm^{\perp}\oplus \mathbb{R}$. From this observation along with the identity 
	$\mathcal{E}(u+c,v+c')=\mathcal{E}(u,v) $ for all $c,c'\in \mathbb{R}$ and the uniqueness of $u\in\VnuOm^{\perp}$
	solving \eqref{eq:var-nonlocal-Neumann-bis}, it becomes easy to check under the compatibility condition
	\eqref{eq:compatible-nonlocal} that all solutions of \eqref{eq:var-nonlocal-Neumann} are of the form $u+c$. 
Analogously, if $u_*\in\VnuOm^{\perp}$ solves \eqref{eq:var-nonlocal-Neumann-bis-weigthed} then under condition \eqref{eq:compatible-nonlocal-weighted} all solutions of \eqref{eq:var-nonlocal-Neumann-weigthed} are of the form $u_*+c$. 

\vspace{1mm}
\noindent Conversely, if the problem \eqref{eq:var-nonlocal-Neumann} (resp. \eqref{eq:var-nonlocal-Neumann-weigthed})
 has a solution, then testing with $v=1$ provides the compatibility condition \eqref{eq:compatible-nonlocal} (resp. \eqref{eq:compatible-nonlocal-weighted}). 
 
 \vspace{2mm}

	\noindent Next observe that, if $w\in \VnuOm$ solves \eqref{eq:var-nonlocal-Neumann} then $w-\mbox{$ \fint_{\Omega} w\mathrm{d}x$}\in \VnuOm^\perp$ solves \eqref{eq:var-nonlocal-Neumann-bis}, and thus by uniqueness we have $u= w-\mbox{$ \fint_{\Omega} w\d x$}$.  Plunging $u \in \VnuOm^{\perp}$ into \eqref{eq:coercivity-quadratic}-\eqref{eq:linearform-g} and \eqref{eq:var-nonlocal-Neumann-bis} one easily arrives at \eqref{eq:weak-regular}. The estimate \eqref{eq:weak-regular-weighted} is obtained  analogously by taking into account the estimate \eqref{eq:linearform-weighted-g} instead of \eqref{eq:linearform-weighted-g}.  
The continuity of the linear mappings $\Phi$ and $\Phi_*$ follow immediately. 

\end{proof}


\subsection{Integrodifferential equations (IDEs) with Dirichlet condition}

Assume that $\Omega\subset \mathbb{R}^d$ is an open set. Let $f:\Omega\to \mathbb{R}$ and $g: \mathbb{R}^d\setminus \Omega\to \mathbb{R}$ be measurable functions. The Dirichlet problem for the operator $L$ associated with the data $f$ and $g$ is to find a measurable  function $u:\mathbb{R}^d\to \mathbb{R}$ such that 
\begin{align}\label{eq:nonlocal-Dirichlet}\tag{$D$}
L u = f \quad\text{in}~~ \Omega \quad\quad\text{ and } \quad\quad  u= g ~~ \text{on}~~ \mathbb{R}^d\setminus\Omega.
\end{align} 

\begin{definition}
 A measurable function $u:\mathbb{R}^d\to \mathbb{R}$ is a weak solution or a variational solution of the inhomogeneous Dirichlet problem \eqref{eq:nonlocal-Dirichlet} if $u\in \VnuOm$ and satisfies
	\begin{align}\label{eq:var-nonlocal-Dirichlet}\tag{$V_0$}
	u-g \in \VnuOmO\quad\text{and}\quad \mathcal{E}(u,v) = \int_{\Omega} f(x)v(x)\mathrm{d}x, \quad \mbox{for all}~~v \in \VnuOmO\,.
	\end{align}
\end{definition}

\vspace{2mm}

\noindent The weak formulation \eqref{eq:var-nonlocal-Dirichlet} makes sense regarding the Green-Gauss formula \eqref{eq:green-gauss-nonlocal}. Let us now see the corresponding minimization formulation. 
\begin{proposition}
A function $ u\in\VnuOm$ is solution to \eqref{eq:var-nonlocal-Dirichlet} if and only if 
	\begin{align}\label{eq:nonlocal-Dirichlet-min}\tag{$M_0$}
	\mathcal{J}_0(u) =\min_{v\in \VnuOmO+g}\mathcal{J}_0(v)
	\end{align}
where
\begin{align*}
\mathcal{J}_0(v) &= \frac{1}{2} \mathcal{E}(v,v) - \int_\Omega f(x) (v(x)-g(x))\, \d x \quad\text{ $ v\in g+\VnuOmO$}.
\end{align*}
\end{proposition} 

\begin{proof}
	Let $ u\in\VnuOmO+g$ satisfying \eqref{eq:var-nonlocal-Dirichlet}. Note that for $ v\in\VnuOm+g$ we have $u-v\in \VnuOmO$ thus from the relation \eqref{eq:var-nonlocal-Dirichlet} and Cauchy-Schwartz inequality we get
	\begin{align*}
\mathcal{E}(u,u)&= \mathcal{E}(u,u-v) + \mathcal{E}(u,v) = \int_\Omega f(x)(u(x)-v(x))\d x+ \mathcal{E}(u,v) \\
	&\leq \frac{1}{2} \mathcal{E}(u,u)+\int_\Omega f(x)(u(x)-g(x))\d x+  \frac{1}{2} \mathcal{E}(v,v)- \int_\Omega f(x)(v(x)-g(x))      \d x
	\end{align*} 
	Equivalently, we have $\mathcal{J}_0(u)\leq \mathcal{J}_0(v)$ and hence $u $ verifies \eqref{eq:nonlocal-Dirichlet-min}. Conversely, assume that $ u\in \VnuOmO+g$ satisfies \eqref{eq:nonlocal-Dirichlet-min}.  For all $v\in \VnuOmO$ and all $t\in \R$ \eqref{eq:nonlocal-Dirichlet-min}  we have $u+tv \in \VnuOmO+g$ which implies  $\mathcal{J}_0(u)\leq \mathcal{J}_0(u+tv)$. Therefore, for fixed $ v\in\VnuOmO$ the polynomial of second order $\mathcal{J}_0(u+\centerdot v): \mathbb{R} \to \mathbb{R}$, 
	\begin{align*}
	t\mapsto \mathcal{J}_0(u+tv)= \mathcal{J}_0(u)+ t\left[\mathcal{E}(u,v)- \int_{\Omega} f(x)v(x)\mathrm{d}x\right] + \frac{t^2}{2} \mathcal{E}(v,v) 
	\end{align*}
has a critical point at $t=0$. Hence \eqref{eq:var-nonlocal-Dirichlet} is verified since we have 
	\begin{align*}
	0= \lim_{t\to 0}\frac{\mathcal{J}_0(u+tv)-\mathcal{J}_0(u)}{t} = \mathcal{E}(u,v)- \int_{\Omega} f(x)v(x)\mathrm{d}x. 
	\end{align*}
\end{proof}

\begin{theorem} \label{thm:nonlocal-Dirichlet-var}
Let  $\Omega\subset \mathbb{R}^d$ be bounded and open. Assume that the function $\nu: \mathbb{R}^d\to [0,\infty]$ satisfies the mild condition \eqref{eq:integrability-condition-near-zero}.
	Given $ f \in L^2 (\Omega)$ and $g\in \TnuOm$, there exists a unique solution $u\in\VnuOm$  to the variational problem \eqref{eq:var-nonlocal-Dirichlet}.  Moreover, there exists a constant $C: = C(d,\Omega, \nu)>0$ independent of $f$ and $g$ such that for any  $\overline{g}\in \VnuOm$ extending $g$, we have 
	\begin{align}\label{eq:energy-estimate}
	\mathcal{E}(u,u)\leq C \big(\|f\|^2_{L^2(\Omega)}+\mathcal{E}(\overline{g}, \overline{g})\big).
	\end{align}
In addition, the following weak regularity estimate holds
	\begin{align}\label{eq:weak-regular-Dirichlet}
	\|u\|_{\VnuOm}\leq C \big(\|f\|_{L^2(\Omega)}+\|g\|_{\TnuOm}\big)\,.
	\end{align}
	
	\noindent In particular, the operator $\Phi' : L^{2}(\Omega) \times \TnuOm \to\VnuOm$ mapping the Dirichlet  data $(f,g) \in L^{2}(\Omega)\times\TnuOm$ 
	to the unique solution $ u\in\VnuOm$ of the variational problem
	\eqref{eq:var-nonlocal-Dirichlet} is linear, one-to-one, continuous and we have
	\begin{align*}
	\|\Phi'(f,g)\|_{\VnuOm}\leq C \|(f,g)\|_{L^{2}(\Omega)\times\TnuOm}.
	\end{align*} 
\end{theorem}

\vspace{2mm}

\begin{proof}
 From the Poincar\'e-Friedrichs inequality \eqref{eq:poincare-friedrichs-ext} (see Theorem \ref{thm:poincare-friedrichs-ext}) one finds a constant $C= C(\Omega, \nu, d)>0$ such that  the following coercivity holds
\begin{align*}
\|v\|^2_{\VnuOm}\leq C \mathcal{E}(v,v),\quad \text{for all }~~~v\in \VnuOmO\,.
\end{align*}
%
Let $\overline{g}\in \VnuOm$ be an extension of $g$. The linear form $v\mapsto -\mathcal{E}(\overline{g}, v)+ \int_\Omega f(x)v(x)\d x$ is continuous on $\VnuOmO$. It is clear that, the symmetric bilinear form $\mathcal{E}(\cdot, \cdot)$ is also continuous. By Lax-Milgram, there exists a unique $u_0\in \VnuOmO$ such that 
\begin{align*}
\mathcal{E}(u_0,v)= -\mathcal{E}(\overline{g}, v)+\int_\Omega f(x)v(x)\d x \quad\text{for all $v\in \VnuOmO$}.
\end{align*} 
\noindent Hence $u= u_0+\overline{g}$  solves \eqref{eq:var-nonlocal-Dirichlet}  since $u-g\in \VnuOmO$ and 
\begin{align*}
\mathcal{E}(u,v)= \int_\Omega f(x)v(x)\d x, \quad\text{for all $v\in \VnuOmO$}. 
\end{align*}
With the aid of the coercivity one can easily show that $u$ is unique and does not depend on the choice of $\overline{g}$.
Meanwhile, since $u-\overline{g}\in \VnuOmO$, using the definition of $u$  and the Poincar\'e-Friedrichs inequality we get 
\begin{align*}
\mathcal{E}(u,u)= \mathcal{E}(u,u-\overline{g})+ \mathcal{E}(u,\overline{g}) 
&= \int_\Omega f(x)( u(x) -\overline{g} (x))\d x + \mathcal{E}(u,\overline{g})\\
&\leq \|f\|_{L^2(\Omega)} \|u-\overline{g}\|_{L^2(\Omega)}+ \mathcal{E}(u,u)^{1/2}\mathcal{E}(\overline{g},\overline{g})^{1/2}\\
&\leq C \|f\|_{L^2(\Omega)} \mathcal{E}(u-\overline{g},u-\overline{g})^{1/2}+ \mathcal{E}(u,u)^{1/2}\mathcal{E}(\overline{g},\overline{g})^{1/2}\\
&\leq\mathcal{E}(u,u)^{1/2} ( C \|f\|_{L^2(\Omega)} +\mathcal{E}(\overline{g},\overline{g})^{1/2})+ C \|f\|_{L^2(\Omega)} \mathcal{E}(\overline{g},\overline{g})^{1/2} \\
&\leq \frac12\mathcal{E}(u,u)+\frac12( C \|f\|_{L^2(\Omega)} +\mathcal{E}(\overline{g},\overline{g})^{1/2})^2+ C \|f\|_{L^2(\Omega)} \mathcal{E}(\overline{g},\overline{g})^{1/2} \\
&=\frac12\mathcal{E}(u,u)+\frac12(C^2 \|f\|^2_{L^2(\Omega)} +\mathcal{E}(\overline{g},\overline{g}))+ 2C \|f\|_{L^2(\Omega)} \mathcal{E}(\overline{g},\overline{g})^{1/2} \\
&\leq \frac12\mathcal{E}(u,u)+\frac12(C^2+1)( \|f\|^2_{L^2(\Omega)} +\mathcal{E}(\overline{g},\overline{g}))+ C ( \|f\|^2_{L^2(\Omega)} +\mathcal{E}(\overline{g},\overline{g}) \\
&= \frac12\mathcal{E}(u,u)+\frac12(C+1)^2( \|f\|^2_{L^2(\Omega)} +\mathcal{E}(\overline{g},\overline{g})), 
\end{align*}
where all along we used the young inequality $|ab|\leq \frac12(a^2+b^2)$. After cancellation we get 
\begin{align*}
\mathcal{E}(u,u)\leq (C+1)^2( \|f\|^2_{L^2(\Omega)} +\mathcal{E}(\overline{g},\overline{g}).
\end{align*}
Hence, \eqref{eq:energy-estimate} follows. Next, since $u-\overline{g}\in \VnuOmO$, the Poincar\'e-Friedrichs inequality \eqref{eq:poincare-friedrichs-ext} implies
\begin{align*}
\|u\|^2_{L^2(\Omega)}&\leq 2\|u-\overline{g}\|^2_{L^2(\Omega)}+2\|\overline{g}\|^2_{L^2(\Omega)}\\
&\leq 2C \mathcal{E}(u-\overline{g},u-\overline{g})+2\|\overline{g}\|^2_{L^2(\Omega)}\\
&\leq 4C \big(\mathcal{E}(u,u)+ \|\overline{g}\|^2_{\VnuOm}\big).
\end{align*}
 Recalling that $\|g\|_{\TnuOm}=\inf\big\{\|\overline{g}\|_{\VnuOm}: \operatorname{Tr}(\overline{g}) =g\big\}$, the estimate \eqref{eq:weak-regular-Dirichlet}  follows as well by combining \eqref{eq:energy-estimate} with the above estimate. 
\end{proof}

\vspace{2mm}

\noindent Next, we establish a nonlocal version of Weyl's decomposition lemma which the local case asserts that $H^1(\Omega) = H_0^1(\Omega)\oplus H_\gamma(\Omega)$ where $H_\gamma(\Omega)=  \big\{u\in H^1(\Omega):  -\Delta u= \gamma u, \text{ in $\Omega$} \text{ and }\, u=0\text{ on $\partial \Omega$}\big\}.$

\begin{proposition}[\textbf{Weyl's lemma}]
	For $\gamma>0$ define the scalar product $\mathcal{E}_\gamma(u,u)= \mathcal{E}(u,u) + \gamma (u,u)_{L^2(\Omega)}$ and the space of $\gamma-$Harmonic functions by 
	$$ V_{\nu,\gamma}(\Omega|\R^d) =\big\{u\in \VnuOm: \, Lu= \gamma u \text{ in $\Omega$}\, \text{ and }\, u=0\text{ on $\Omega^c$}\big\}.$$
	The following $\mathcal{E}_\gamma(\cdot, \cdot) $-orthogonal decomposition is true $$\VnuOm= \VnuOmO\oplus V_{\nu,\gamma}(\Omega|\R^d).$$
\end{proposition}
%


\begin{proof}
	Let $u\in\VnuOm$. We know there is a unique $u'\in \VnuOm$ such that $u''= u-u'\in \VnuOmO$ and $\mathcal{E}_\gamma(u',v) =0$ for all $v\in \VnuOmO$ that is $u'\in V_{\nu,\gamma}(\Omega|\R^d)$. In particular we have $\mathcal{E}_\gamma(u',u'') =0$, i.e. $u'$ and $u''$ are $\mathcal{E}_\gamma(\cdot, \cdot) $-orthogonal. Thus, $u=u'+u''\in  \VnuOmO\oplus V_{\nu,\gamma}(\Omega|\R^d).$
\end{proof}

\subsection{Integrodifferential equations (IDEs) with mixed condition}

Here we treat a special IDE whose complement condition is a mixture of the Dirichlet and the Neumann complement conditions. Assume $\Omega\subset \mathbb{R}^d$ is an open set. Let $D\subset \Omega^c$ and  $N\subset \Omega^c$ be measurable sets such that $|D|, |N|>0$, $\Omega^c= D\cup N$ and $|D\cap N|=0$.
Let $f:\Omega\to \mathbb{R}$, $g_D: D\to \mathbb{R}$  and  $g_N: N\to \mathbb{R}$ be measurable functions. The mixed  problem for the operator $L$ associated with the data $f$, $g_D$ and $g_N$ is to find a measurable  function $u:\mathbb{R}^d\to \mathbb{R}$ such that 
\begin{align}\label{eq:nonlocal-mixed}
L u = f \quad\text{in}~~~ \Omega, \quad u=g_D\quad \text{on}~~D\quad \text{and} \quad\mathcal{N} u= g_N~~~ \text{on}~~N.
\end{align} 
\noindent The local counterpart of the problem \eqref{eq:nonlocal-mixed} for the Laplace operator is to find a measurable function $u: \Omega\to \R$ such that 
\begin{align}\label{eq:local-mixed}
-\Delta u = f \quad\text{in}~~~ \Omega, \quad u=g_D\quad \text{on}~~D\quad \text{and} \quad\frac{\partial u }{\partial n} = g_N~~~ \text{on}~~N.
\end{align} 
where the functions $g_D$ and $g_N$ are respectively defined on $D\subset \partial \Omega$ and $N\subset \partial \Omega$ so that $D\cup N= \partial \Omega$  with $|D\cap N|=0$. As for the  Neumann and the  Dirichlet complement conditions, a  weak solution of \eqref{eq:nonlocal-mixed} is defined as follows. 
\begin{definition}
A measurable function $u:\mathbb{R}^d\to \mathbb{R}$ is a weak solution or a variational solution of the inhomogeneous mixed problem \eqref{eq:nonlocal-mixed} if $u\in \VnuOm$ and satisfies
	\begin{align}\label{eq:var-nonlocal-mixed}
	u-g_D \in V_D\quad\text{and}\quad \mathcal{E}(u,v) = \int_{\Omega} f(x)v(x)\mathrm{d}x+ \int_{N} g_N(y)v(y)\d y \quad \mbox{for all}~~v \in V_D\,
	\end{align}
	where $V_D=\big\{v\in \VnuOm: v=0 ~~\text{a.e. on D}\big\}$.
\end{definition}

\medskip

\begin{theorem} \label{thm:nonlocal-mixed-var}
	Assume $\Omega\subset \mathbb{R}^d$ is a bounded open set and the function $\nu: \mathbb{R}^d\to [0,\infty]$ is the density of a symmetric L\'{e}vy measure with full support such that the couple $(\nu, \Omega)$ belongs to one of the class $\mathscr{A}_i,~i=1,2,3,4$ (see page  \pageref{eq:class-lipschitz}). For a set  $K\subset \Omega$ with positive measure, assume $ \nu_K(x)= \operatorname{essinf}_{y\in K}~\nu(x-y)>0~a.e$ for almost $x\in \Omega^c$. 
	Then given $ f \in L^2 (\Omega)$, $g_D\in \TnuOm$ and $g_N \in L^2(\Omega^c, \nu_K^{-1})$, there exists a unique solution $u\in\VnuOm $ to the variational problem \eqref{eq:var-nonlocal-mixed}. 
	
	\noindent Moreover, there exists a constant $C: = C(d,\Omega, \nu)>0$ independent of $f$ and $g$ such that $u$ satisfies the following weak regularity estimate
	\begin{align}\label{eq:weak-regular-m}
	\|u\|_{\VnuOm}\leq C \left(\|f\|_{L^2(\Omega)}+ \|g_D\|_{\TnuOm}+\|g_N\|_{L^{2}(\Omega^c, \nu_K^{-1})}\right)\,.
	\end{align}
\end{theorem}

\begin{proof}
Since $|D|>0$,  by Theorem \ref{thm:poincare-inequality-general-bis} the Poincar\'e inequality holds on $V_D$, i.e. there exists a constant $C=C(d,\Omega, D,\nu,)$ such that
	\begin{align*}
	\|v\|_{L^2(\Omega)}\leq C\mathcal{E}(v,v)\qquad\text{for all $v\in V_D$}.
	\end{align*}
Therefore, the proof follows by adapting the proofs of Theorem \ref{thm:nonlocal-Neumann-var}  and Theorem \ref{thm:nonlocal-Dirichlet-var}.
\end{proof}

\subsection{Integrodifferential equations (IDEs) with Robin condition}
We now treat the Robin type problem with respect to the nonlocal operator $L$ on $\Omega$. In the classical setting for the Laplace operator, the Robin boundary problem also known as Fourier boundary problem or third boundary problem is a combination of the Dirichlet and Neumann boundary problem  and is often given as follows.
\begin{align*}
-\Delta u = f \quad\text{in}\quad \Omega\quad\text{and}\quad \frac{\partial u}{\partial n} + b u= g \quad\text{on}\quad \partial\Omega.
\end{align*}
Here $f:\Omega\to \R$ and  $b, g: \partial\Omega\to \mathbb{R}$ are given measurable functions. For a more extensive study on the local Robin boundary value problems see \cite{Dan00, Ro14}. Analogously, in the nonlocal set up, we assume that $b, h: \Omega^c\to \mathbb{R}$ are measurable functions. We now formulate the Robin problem with respect to $L$ with data $f,h$ and $b$. It consists of  finding a measurable function $u:\mathbb{R}^d\to \mathbb{R}$ such that 
\begin{align}\label{eq:robin-problem}
L u = f \quad\text{in}\quad \Omega\quad\text{and}\quad \mathcal{N}u + b u= h \quad\text{on}\quad \Omega^c.
\end{align}
\noindent
For the sake of simplicity we assume $b(x)= \beta(x)\nu_K(x)$ and $h(x) = g(x)\nu_K(x)$, $x\in \Omega^c$, for a measurable set $K\subset \Omega$ with a positive measure. We recall that $\nu_K(x)= \operatorname{essinf}_{y\in K}~\nu(x-y)$.  If $\beta=0$ one recovers a Neumann problem and if $\beta\to \infty $ then it leads to a Dirichlet problem.
As for the weak formulation of the Neumann problem, we define a weak solution to $\eqref{eq:robin-problem}$ as follows. 
\begin{definition}
 A function $u:\mathbb{R}^d\to \mathbb{R}$ is a weak solution of the Robin problem \eqref{eq:robin-problem} if 
\begin{align}\label{eq:weak-robin-problem}
Q_\beta(u,v)=\int_{\Omega} f(x)v(x)\d x+ \int_{\Omega^c} g(y)v(y)\nu_K(y)\d y ~~\text{for all }~v\in \VnuOm,
\end{align}
 where we introduce the bilinear form
\begin{align*}
Q_\beta(u,v)= \mathcal{E}( u,v)+ \int_{\Omega^c} u(y)v(y)\beta(y)\nu_K(y)\d y. 
\end{align*} 
\end{definition}

\medskip

\begin{theorem} \label{thm:nonlocal-Robin-var}
	Assume $\Omega\subset \mathbb{R}^d$ is a bounded open set and the function $\nu: \mathbb{R}^d\to [0,\infty]$ is the density of a symmetric L\'{e}vy measure with full support such that the couple $(\nu, \Omega)$ belongs to one of  the class $\mathscr{A}_i,~i=1,2,3$ (see page \pageref{eq:class-lipschitz}). 
	For a  set  $K\subset \Omega$ with a positive measure, let $\nu_K(x)= \operatorname{essinf}_{y\in K}~\nu(x-y)$ for almost every $x\in \Omega^c$. Let $ \beta: \Omega^c\to [0, \infty) $ be essentially bounded such that $\beta \nu_K>0~ a.e$ on a subset of positive of $\Omega^c$. 
	
\noindent Then given $ f \in L^2 (\Omega)$ and $g\in L^2(\Omega^c, \nu_K)$, there exists a unique solution $u\in\VnuOm $ to the variational problem \eqref{eq:weak-robin-problem}. Moreover, there exists a constant $C: = C(d,\Omega, \nu, K, \beta)$ independent of $f$ and $g$ such that
	\begin{align}\label{eq:weak-regular-robin}
	\|u \|_{\VnuOm}\leq C \left(\|f\|_{L^2(\Omega)}+\|g\|_{L^{2}(\Omega^c, \nu_K)}\right).
	\end{align}
	
	\noindent In particular, the operator $\Psi : L^{2}(\Omega) \times L^{2}(\Omega^c, \nu_K)\to\VnuOm$ mapping the Robin data $(f,g) \in L^{2}(\Omega)\times L^{2}(\Omega^c, \nu_K)$ to the unique solution $ u\in\VnuOm$ of the variational problem \eqref{eq:weak-robin-problem} is linear, one-to-one, continuous and we have
	\begin{align*}
	\|\Psi(f,g)\|_{\VnuOm}\leq C \|(f,g)\|_{L^{2}(\Omega)\times L^{2}(\Omega^c, \nu_K)}.
	\end{align*} 
\end{theorem}

\medskip

\begin{proof}
	First of all, we claim that the form $Q_\beta(\cdot, \cdot)$ is coercive on $\VnuOm$. Assume it is not true then for each $n\geq 1$ there exists $u_n \in\VnuOm $ preferably, $\|u_n\|_{\VnuOm}=1$ such that 
	\begin{align*}
	Q_\beta(u_n , u_n )<\frac{1}{2^n}\,. 
	\end{align*}
	\noindent In virtue of compactness Theorem \ref{thm:embd-compactness}, $(u_n)_n$ converges up to a subsequence in $L^2(\Omega)$ to some $u$. It turns out that $\|u\|_{L^2(\Omega)}= 1$, since $\mathcal{E}(u_n , u_n ) \xrightarrow[]{n \to \infty}0$ and for all $n\geq 1$, $\|u_n\|_{\VnuOm}=1$. That $\mathcal{E}(u_n , u_n )\xrightarrow[]{n \to \infty} 0$ and $\|u_n-u\|_{L^2(\Omega)} \xrightarrow[]{n \to \infty}0$ imply that $u_n$ converges to $u$ in $\VnuOm$ with $\mathcal{E}(u,u)=0$ and thus $u $ is constant almost everywhere in $\mathbb{R}^d$. On the other hand, since $\beta$ is bounded and the embedding $\VnuOm\hookrightarrow L^2(\Omega^c, \nu_K)$ is continuous, we have 
	\begin{align*}
	\int_{\Omega^c} u^2(y)\beta(y)\nu_K(y)\d y&\leq 2 \int_{\Omega^c} u_n^2(y)\beta(y)\nu_K(y)\d y+2\|\beta\|_\infty \int_{\Omega^c} (u_n(y)-u(y))^2\nu_K(y)\d y\\
	&\leq 2Q_\beta(u_n, u_n)+ C\|u_n-u\|^2_{\VnuOm}\xrightarrow[]{n \to \infty}0\,.
	\end{align*}
	
	\noindent From this, we have $u=0$ since we know that $u$ is a constant function and $\beta\nu_K>0$ almost everywhere on a set of positive measure $U\subset \Omega^c$ on which $u $ vanishes. This negates the fact that $\|u\|_{L^2(\Omega)}=1$ and hence our initial assumption was wrong. Therefore, there exists a constant $C=C(d, \Omega, \nu, K, \beta)>0$ such that
	\begin{align}\label{eq:coercive-robin}
	Q_\beta(u,u)\geq C\|u\|^2_{\VnuOm}\quad\text{for all }~u\in\VnuOm.
	\end{align}
	
	\medskip
	
	\noindent Beside this, it is easy to prove that the remaining assumptions of the Lax-Milgram lemma are met in order to guarantee the existence and uniqueness of a solution to \eqref{eq:weak-robin-problem}. The estimate
	\eqref{eq:weak-regular-robin} is a direct consequence of \eqref{eq:coercive-robin}.
\end{proof}

\subsection*{Comments}
A separate objective was to show that  the function spaces studied in Chapter \ref{chap:nonlocal-sobolev} are adapted for the study of IntegroDifferential Equations (IDEs) on domains with complement conditions similar to those present in the study of  boundary values problems associated with elliptic Partial Differential Equations (PDEs)  of second order.  It is noteworthy to highlight that our context,  however easily extends to the setting where the interacting kernel $\nu(x-y)$ is substituted by \emph{an elliptic symmetric kernel $k:\R^d\times \R^d \to [0,\infty]$}, i.e. satisfying the  following weak elliptic condition: 

\vspace{1mm}
\noindent $(E_k)$: There exists a constant $\Lambda>1$ such that for every measurable function $u:\R^d\to \R$ we have 
\begin{align}\label{eq:weak-elliptic}\tag{$E_k$}
\Lambda^{-1}\mathcal{E}(u,u) \leq \mathcal{E}^k(u,u)\leq \Lambda\mathcal{E}(u,u),
\end{align}

\noindent where
\begin{align*}
\mathcal{E}^k(u,u)= \frac{1}{2}\iil_{ (\Omega^c \times\Omega^c)^c} (u(x)-u(y) )^2 k(x,y)\d x \d y. 
\end{align*}

\noindent It is worth noting that choosing a kernel $k:\R^d\times \R^d \to [0,\infty]$ in the nonlocal setting is the same as choosing a matrix $A:\R^d\to \R^{d^2}$, $A(x) = (a_{ij} (x))_{1\leq i,j\leq d}$ in the local setting. A concrete explanation is encoded in Chapter \ref{chap:basics-nonlocal-operator} and we elaborate more on this fact in Chapter \ref{chap:Nonlocal-To- Local}. 

\medskip

\noindent  Let us mention that the Dirichlet problem for nonlocal operators has been extensively studied throughout the literature. The Dirichlet problem for nonsymmetric kernels is thoroughly carried out in \cite{FKV15}. See also \cite{Ros15} and further references therein. The Neumann problem for nonsymmetric kernel requires some additional knowledge on a way to control the set of functions $u$ on which the linear form $\mathcal{E}^k(u,\cdot)$ is degenerated. This is not  straightforward since if $k$ is nonsymmetric, one can easily check that $\mathcal{E}^k(u,u)= \mathcal{E}^{k_s}(u,u)$ where $k_s(x,y) = \frac{1}{2}\big(k(x,y)+k(y,x\big)$ is the symmetric part of $k$.

\noindent If $\nu$ is fully supported on $\R^d$  and $k$ is symmetric and satisfying \eqref{eq:weak-elliptic},  then  $\mathcal{E}^k(u,\cdot)$ is degenerated if and only if $u$ is a constant. From this the solvability of the nonlocal Neumann for a symmetric kernel $k$ satisfying \eqref{eq:weak-elliptic} easily extends since in this case the compatibility condition \eqref{eq:compatible-nonlocal}  persists and it thus suffices to work on the space $\VnuOm^\perp$. 
 \vspace{2mm}
 
\noindent Furthermore, when working on the space $V_k(\Omega|\R^d)$ of measurable functions $u:\R^d\to \R$ for which $\mathcal{E}^k(u,u)<\infty$ in replacement of the space $\VnuOm$, the lower elliptic bound condition \eqref {eq:lower-weak-elliptic} below would be sufficient as it guaranties the coercivity of the form $\mathcal{E}^k(\cdot,\cdot)$ via Poincar\'e type inequalities.

\vspace{2mm}

\noindent ($L$-$E_k$): There exists a constant $\Lambda>1$ such that for every measurable function $u:\R^d\to \R$ we have 

\begin{align}\label{eq:lower-weak-elliptic}\tag{$L$-$E_k$}
\Lambda^{-1}\mathcal{E}(u,u) \leq \mathcal{E}^k(u,u). 
\end{align}

\section{Spectral decomposition of nonlocal operators}

The main subject of interest in this section is the spectral decomposition of a singular nonlocal operator of the 
form $L$ on a bounded domain $\Omega$ subject to the Neumann, Dirichlet or Robin boundary condition. 
Precisely, we are concerned about the pairs $(\lambda,u)$ consisting of a real number $\lambda$ called an eigenvalue and a function $u\in \VnuOm$ called eigenfunction so that (in the strong sense) we have $Lu = \lambda u$ on $\Omega$ plus some additional conditions which we formally consider in the next definitions. These spectral decomposition have some properties of important interest which we will 
explore here and later. Let us first define the weak version eigenvalues and eigenfunctions.


\begin{definition}[Neumann eigenvalue of $L$]\label{def:neumann-eigenvalue}
A non zero function $u \in \VnuOm$ will be called a Neumann eigenfunction of the operator $L$ on $\Omega$ if there exists a real number $\mu $ which is the eigenvalue associated with $u$ such that the following holds 
\begin{align}\label{eq:neumann-eigenvalue}
 \mathcal{E}( u ,v) = \mu \int_{\Omega} u(x)v(x)\mathrm{d}x \quad \mbox{for all}~~v \in \VnuOm\,.
 \end{align}
 We will formally write $Lu=\mu u $ in $\Omega$ and $\mathcal{N}u=0$ on $\Omega^c$ which in 
 fact corresponds to the weak formulation in \eqref{eq:neumann-eigenvalue} provided that $u $ is regular enough.
\end{definition} 
 
 \noindent It is worth noticing that if $u$ is a Neumann eigenfunction of $L$ with associated eigenvalue $\mu$ then automatically, either $u \in \VnuOm^\perp$ when $\mu\neq 0$ or else $\mu=0$ and the constant functions $u=c, ~c\in \mathbb{R}\setminus\{0\}$ as the related eigenfunctions. 

 \medskip

\begin{definition}[Dirichlet eigenvalue of $L$]
A non zero function $u \in \VnuOmO$ will be called a Dirichlet eigenfunction of the operator $L$ on $\Omega$ if there exists a real number $\lambda $ which is the eigenvalue associated with $u$ such that the following holds 
\begin{align}\label{eq:def-dirichlet-eigenvalue}
 \mathcal{E}( u ,v) = \lambda \int_{\Omega} u(x)v(x)\mathrm{d}x \quad \mbox{for all}~~v \in \VnuOmO\,.
 \end{align}
 We will formally write $Lu=\lambda u $ in $\Omega$ and $u=0$ on $\Omega^c$ which in fact corresponds to the weak formulation in \eqref{eq:def-dirichlet-eigenvalue} provided that $u $ is regular enough.
\end{definition}

\medskip

\begin{definition}[Robin eigenvalue of $L$]
Let $\beta: \Omega^c\to [0, \infty]$ be a measurable function and $K\subset \Omega $ be a measurable set such that $|K|>0$. A non zero function $u \in \VnuOm$ will be called a Robin eigenfunction of the operator $L$ on $\Omega$ if there exists a real number $\gamma(\beta)$ which is the eigenvalue associated with  $u$ such that the following holds 
\begin{align}\label{eq:def-robin-eigenvalue}
 Q_\beta(u,v) = \gamma(\beta) \int_{\Omega} u(x)v(x)\mathrm{d}x \quad \mbox{for all}~~v \in \VnuOm\,,
 \end{align}
 where we recall that
 \begin{align*}
 Q_\beta(u,v)= \mathcal{E}( u, v)+ \int_{\Omega^c}u(y)v(y)\beta(y)\nu_K(y)\,\d y\,.
 \end{align*}
 We will formally write $Lu=\gamma(\beta) u$ in $\Omega$ and $\mathcal{N}u+\beta\nu_K u=0$ on $\Omega^c$ which in fact corresponds to the weak formulation in \eqref{eq:def-robin-eigenvalue} provided that $u $ is regular enough.
\end{definition} 

\medskip

\begin{theorem} \label{thm:existence-of-eigenvalue-Neumann}
Assume $\Omega\subset \mathbb{R}^d$ is a bounded open set and the function $\nu: \mathbb{R}^d\to [0,\infty]$ is the density of a symmetric L\'{e}vy measure with full support such that the couple $(\nu, \Omega)$ belongs to one of the class $\mathscr{A}_i,~i=1,2,3$. Then there exist a sequence $(\phi_n)_{\in \mathbb{N}_0}$ in $\VnuOm$, orthonormal basis of $L^2(\Omega)$ and an increasing sequence of real numbers $ 0=\mu_0<\mu_1\leq \cdots\leq \mu_n\leq \cdots. $ such that $\mu_n \to \infty$ as $n \to \infty$ and each $\phi_n$ is a Neumann eigenfunction of $L$ whose corresponding eigenvalue is $\mu_n$. Each eigenvalue is listed with its geometrical multiplicity. 
\end{theorem}


\begin{proof}
 For $f_1, f_2\in L^2(\Omega)$ let us denote $ u_{f_k}=\Phi_0(f_k) =\Phi(f_k,0)\in \VnuOm^{\perp}, ~k=1,2$ the unique solution of \eqref{eq:var-nonlocal-Neumann-bis} with Neumann data $f=f_k$ and $g=0$. Precisely, 
\begin{align}\label{eq:var-nonlocal-Neumann-test}
 \mathcal{E}( \Phi_0(f_k),v) = \int_{\Omega} f_k(x)v(x)\mathrm{d}x \quad \mbox{for all}~~v \in \VnuOm^{\perp}\,.
 \end{align}

\noindent Testing \eqref{eq:var-nonlocal-Neumann-test} against $v= \Phi_0(f_2)$ and $v=\Phi_0(f_1)$ successively when $k=1$ and $k=2$ yields 
\begin{align*}
\big(f_1, \Phi_0(f_2)\big)_{L^2(\Omega)} = \mathcal{E}(\Phi_0(f_1), \Phi_0(f_2)) = \mathcal{E}(\Phi_0(f_2), \Phi_0(f_1)) = \big(f_2,\Phi_0(f_1)\big)_{L^2(\Omega)} .
\end{align*}
Therefore, the operator $R_\Omega\circ \Phi_0: L^2(\Omega)\xrightarrow[]{\Phi_0} \VnuOm^{\perp}\xrightarrow[]{R_\Omega} L^2(\Omega)^{\perp}$ is compact (by Theorem \ref{thm:embd-compactness}) and symmetric hence self-adjoint. Recall that $R_\Omega$ denotes the restriction operator to $\Omega$.  It is a fact from the spectral theory of compact self-adjoint operators that $L^2(\Omega)^\perp$ possesses an orthonormal basis $(e_n)_n$ whose elements are eigenfunctions of $R_\Omega\circ \Phi_0$ and the sequence of the corresponding 
eigenvalues are non-negative real numbers $(r_n)_n$ which we assume are ordered in the decreasing order, $r_1\geq r_2\geq \cdots\geq r_n\geq \cdots 0$ such that $r_n\to 0$ as $n\to \infty$. 
Precisely, for each $n\geq 1$, $R_\Omega\circ \Phi_0(e_n) = r_n e_n$ or simply write $\Phi_0(e_n) = r_n e_n$ a.e in $\Omega$. Combining the latter relation with  the definition of $ \Phi_0(e_n)$ we get 
\begin{align*}
 \mathcal{E}( \Phi_0(e_n),v) &= \int_{\Omega} e_n(x)v(x)\mathrm{d}x= r^{-1}_n \int_{\Omega}\Phi_0( e_n) (x)v(x)
 \mathrm{d}x\, \quad \mbox{for all}~~v \in \VnuOm^{\perp}\,. 
 \end{align*}
Equivalently, letting $\mu_n= r^{-1}_n$ and $\phi_n = \Phi_0(e_n)/\|\Phi_0(e_n)\|_{L^2(\Omega)}= r_n^{-1} \Phi_0(e_n) $ which is clearly an element of $\in \VnuOm^\perp$ yields 
\begin{align*}
 \mathcal{E}( \phi_n,v) &= \mu_n \int_{\Omega} \phi_n (x)v(x)\mathrm{d}x\, \quad \mbox{for all}~~v \in \VnuOm^{\perp}\, .
 \end{align*}
 Hereby, along with $\mu_0=0$ and $\phi_0= |\Omega|^{-1}$ provides the sequences sought for. 
 Now if we assume $\mu_1=0$ then we have $\phi_1\in \VnuOm^\perp$ and $\mathcal{E}(\phi_1, v) = 0$ for all $v \in \VnuOm^\perp$. In particular $\mathcal{E}(\phi_1, \phi_1) = 0$, i.e $\phi_1$ is a constant function in $\VnuOm^\perp$ necessarily $\phi_1=0$ since $\phi_1$ has zero mean over $\Omega$. We have therefore reached a contradiction as $\phi_1 $ is supposed to be an eigenfunction, i.e. $\phi_1\neq 0$. Thus, $\mu_1>0$ and this  completes the proof. 
\end{proof}

\medskip

\noindent Next, by exploiting \cite{Nic11}, we characterize the eigenpairs of the operator $L$ using the Rayleigh quotient.

\begin{theorem}\label{thm:critical-rayleigh}
Assume $\Omega\subset \mathbb{R}^d$ is a bounded open set and the function $\nu: \mathbb{R}^d\to [0,\infty]$ is the density of a symmetric L\'{e}vy measure with full support such that the 
couple $(\nu, \Omega)$ belongs to one of  the class $\mathscr{A}_i,~i=1,2,3$ (see page \pageref{eq:class-lipschitz}) . The following assertions are true.
\begin{enumerate}[$(i)$]
 \item There exists $u\in \VnuOm^\perp$ which is a global minimum of $ v\mapsto \mathcal{E}(v,v)$ in $\VnuOm^\perp $ subject to the constraint $\|v\|_{L^2(\Omega)}=1$. In short, 
 \begin{align*}
 \mathcal{E}(u,u) = \min\{\mathcal{E}(v,v):~ v\in \VnuOm^\perp,~ \|v\|_{L^2(\Omega)}=1\}.
 \end{align*}
 \item Let $u$ be a local extremum (minimum or maximum)  of $ v\mapsto \mathcal{E}(v,v)$ in $\VnuOm^\perp $  subject to the constraint $\|v\|_{L^2(\Omega)}=1$, then $u$ is an eigenfunction of $L$ subject to Neumann condition. Moreover, the corresponding eigenvalue is a strictly positive real number given by $\mu = \mathcal{E}(u,u)$. 
\end{enumerate}
\end{theorem}

\medskip

\begin{proof}
To prove $(i)$, note that the quantity $m= \inf\{\mathcal{E}(v,v): v\in \VnuOm^\perp,\|v\|_{L^2(\Omega)}=1\}$ exists according to the Poincar\'{e} inequality \eqref{eq:poincare-inequalityV}. Next, let  $(u_n)_n\subset \VnuOm^\perp $ be a minimizing sequence such that for all $n\geq 1$, $\|u_n\|_{L^2(\Omega)}=1$ and 
$$ m\leq \mathcal{E}(u_n,u_n)\leq m+\frac{1}{2^n}.$$ 
The sequence $(u_n)_n$ is clearly bounded in $\VnuOm$. Therefore, according to Theorem \ref{thm:embd-compactness}, there exists $u\in \VnuOm$ so that a subsequence of $(u_n)_n$  weakly converges in $\VnuOm$ and strongly converges in $L^2(\Omega)$ to $u$. Recalling that $ \|\cdot\|^2_{\VnuOm)}=\|\cdot \|^2_{L^2(\Omega)}+ \mathcal{E}(\cdot, \cdot)$, the weak convergence combined with $\|u_n\|_{L^2(\Omega)}=\|u\|_{L^2(\Omega)}=1$ yields, 
\begin{align*}
\mathcal{E}(u,u)\leq \liminf_{n\to\infty}\mathcal{E}(u_n,u_n)\leq\liminf_{n\to\infty} (m+ \frac{1}{2^n})= m.
\end{align*}
That is, $ \mathcal{E}(u,u)\leq m$. In addition, the strong convergence implies that $\fint_\Omega u =0$ so that $u\in \VnuOm^\perp$ and $\|u\|_{L^2(\Omega)}=1$. It follows that
 $$m=\mathcal{E}(u,u)= \min\{\mathcal{E}(v,v):~ v\in \VnuOm^\perp,~\|v\|_{L^2(\Omega)}=1\}.$$

\noindent Let us now tackle $(ii)$. On the space $\VnuOm^\perp$, let us define $\Lambda: v \mapsto \mathcal{E}(v,v)$ and $ G: v\mapsto \|v\|^2_{L^2(\Omega)}-1$. So that $u$ is a local extremum of $\Lambda$ subject to the constraint $G(u) =0$. It is  routine to show that $\Lambda $ and $G$ are continuously Fr\'{e}chet differentiable. Furthermore, the first variation of $\Lambda$ and $G$ reveals that 
\begin{align*}
\Lambda'(u)(v)= 2 \mathcal{E}(u,v) \quad\text{and}\quad G'(u)(v)= 2(u,v)_{L^2(\Omega)}\qquad\text{for all } ~v\in \VnuOm^\perp. 
\end{align*}
Moreover, $G'(u)\neq0$ since $G'(u)(u) = 2$, given that $G(u) = 0$. Applying the Lagrange multiplier Theorem \ref{thm:lagrange-multiplier} with the space $\VnuOm^\perp$ serving the role of the Banach space, 
one finds a real number $\mu$ such that $\Lambda'(u) =\mu G'(u) $. Equivalently, we have 
\begin{align*}
\mathcal{E}(u,v) = \mu \int_{\Omega}u(x)v(x)\,\d x \qquad\text{for all } ~v\in \VnuOm^\perp. 
\end{align*}
 We have reached the relation \ref{eq:neumann-eigenvalue} claiming $u$ as a Neumann eigenfunction of $L$ and we have $\mu= \mathcal{E}(u,u)$ since $G(u) =0$. 
 If $\mu= \mathcal{E}(u,u)= 0$ then $u$ is constant with a vanishing mean over $\Omega$ since $\VnuOm^\perp$. 
 Necessarily, $u=0$ which goes against the fact that $\|u\|_{L^2(\Omega)}=1$. Thus $\mu>0$ which achieves the proof. 

\end{proof}

\medskip

\begin{remark}
The assertions $(i)$ and $(ii)$ of Theorem \ref{thm:critical-rayleigh} reveal that a global minimum $u\in \VnuOm^\perp$ of $v\mapsto\mathcal{E}(v,v)$ subject to the constraint $ \|v\|_{L^2(\Omega )} =1$ exists. Moreover, this minimum is a Neumann eigenfunction of $L$ whose corresponding eigenvalue is $\mu= \mathcal{E}(u,u)>0$. In fact $\mu$ is the smallest positive Neumann eigenvalue of $L$. For this reason, we shall call this eigenvalue $\mu_1$ and one of the associated eigenfunction $\varphi_1$. Note that the smallest eigenvalue is $\mu_0 =0$ and the associated eigenfunctions are constant functions $\varphi_0=c$ with $c\neq 0$. 

\noindent Another important remark is that, minimizing $v\mapsto\mathcal{E}(v,v)$ on $\VnuOm^\perp$ subject to 
$\|v\|_{L^2(\Omega)} =1$ is equivalent as minimizing the quotient $v\mapsto\tfrac{\mathcal{E}(v,v)}{ \|v\|^2_{L^2(\Omega )}}$ with $v$ running in $\VnuOm^\perp$. This 
is the \emph{Rayleigh quotient}. Wherefore, we have the following explicit representation 
\begin{align*}
\mu_1= \inf_{v\in \VnuOm^\perp}\Big\{ \frac{\mathcal{E}(v,v)}{ \|v\|^2_{L^2(\Omega )}}~ \Big\} = \min_{v\in \VnuOm^\perp}\{\mathcal{E}(v,v):~ \|v\|_{L^2(\Omega)}=1\}. 
\end{align*}
\end{remark}

\bigskip

\noindent Next, we explain how to find the other eigenvalues by employing a similar technique.
The following theorem is a refinement Theorem \ref{thm:existence-of-eigenvalue-Neumann} providing a more constructive approach to the eigenvalues $(\mu_n)_n$. 

\bigskip

\begin{theorem}\label{thm:rayleigh-construction}
Let $\Omega\subset \mathbb{R}^d$ be a bounded open set and the function $\nu: \mathbb{R}^d\to [0,\infty]$ is the density of a symmetric L\'{e}vy measure with full support. Assume that the 
couple $(\nu, \Omega)$ belongs to one of the class $\mathscr{A}_i,~i=1,2,3$ (see page \pageref{eq:class-lipschitz}). Then there exist a sequence $(\varphi_n)_{n\in\mathbb{N}_0}$ of positive real numbers and a sequence of functions $(\varphi_n)_{n\in\mathbb{N}_0}$ in $\VnuOm$ such that for each $n\geq0$ the pair $(\mu_n,\varphi_n)$ satisfies
\begin{align}\label{eq:Rayleigh-algorithm-Neumann}
 \mu_n =\mathcal{E}(\varphi_n,\varphi_n) =\min_{v\in V_n}\{\mathcal{E}(v,v):~ \|v\|_{L^2(\Omega)}=1\}= \inf_{v\in V_n}\Big\{ \frac{\mathcal{E}(v,v)}{ \|v\|^2_{L^2(\Omega )}}~ \Big\}. 
\end{align}
\noindent Here we denote $V_0 = \VnuOm$, $\varphi_0 = |\Omega|^{-1}$ and for $n\geq 1$, 
\begin{align*}
 V_n = \{\varphi_0,\varphi_1,\cdots,\varphi_{n-1}\}^\perp = \{v\in \VnuOm : (v, \varphi_i)_{L^2(\Omega)} =0,~~i=0,1,\cdots,n-1\}.
\end{align*}
Moreover, the following assertions hold true. 

\begin{enumerate}[$(i)$]
 \item For each $n\geq 0,$ $\mu_n$ is a Neumann eigenvalue of $L$ the operator $L$ constrained with the Neumann boundary condition $\mathcal{N} u =0$ on $\Omega^c$ and the associated eigenfunction is $\varphi_n.$
 \item $0=\mu_0<\mu_1\leq \cdots\leq \mu_n\leq \cdots $ and $\lim\limits_{n\to \infty}\mu_n =\infty$. 
 \item The family $(\varphi_n)_{n\in \mathbb{N}_{0}}$ is an orthonormal basis of $L^2(\Omega)$. 
\end{enumerate}
\end{theorem}

\medskip

\begin{proof}
Obviously, we have $V_0= \VnuOm $, $\mu_0=0$ and $\varphi_0= |\Omega|^{-1}$. We also have $V_1= \VnuOm^\perp = \{\varphi_0\}^\perp = \{v\in \VnuOm : \int_\Omega v =0\}$. 
The existence of $(\mu_1, \varphi_1)$ clearly springs from Theorem 
\ref{thm:existence-of-eigenvalue-Neumann}.
Repeating this procedure inductively by modifying appropriately with $V_1$ replaced by $V_2, V_3,\cdots $ one obtains the existence of $\mu_n$ and $\varphi_n$ for each $n\geq 1$. 
The argument is based on the grounds that for each $n\geq 1$ the linear forms $v \mapsto (v, \varphi_i)_{L^2(\Omega)}~i=0,1,\cdots, n-1$ are continuous on $\VnuOm$ which implies that $V_n$ is a closed subspace of $\VnuOm$ and hence can be thought as a Hilbert space in its own right. 
Thus \eqref{eq:Rayleigh-algorithm-Neumann} is verified and $(i)$ follows by identical argument as in Theorem \ref{thm:existence-of-eigenvalue-Neumann}. 
It is worth noticing that by construction all elements of the sequence $(\varphi_n)_n$ are mutually orthogonal in $L^2(\Omega)$. 
Besides, $V_0\supset V_1\supset V_2\supset \cdots\supset V_n\supset \cdots$ wherefrom it follows that $0=\mu_0\leq \mu_1\leq \cdots \leq \mu_n\leq \cdots$ thus we have 
\begin{align*}
 \lim_{n\to \infty} \mu_n = B\quad \text{with }\quad B=\sup_{n\in \mathbb{N}} \mu_n. 
\end{align*}
Assume $B$ is a real number, i.e. $\sup_{n\in \mathbb{N}} \mu_n= \sup_{n\in \mathbb{N}} \mathcal{E}(\varphi_n, \varphi_n) =B<\infty$.  For all $n\geq 1$ we also have $\|\varphi_n\|_{L^2(\Omega)} =1$.
In other words, the sequence $(\varphi_n)_n$ is bounded in $\VnuOm$. In virtue of the compactness Theorem \ref{thm:embd-compactness} there exists a subsequence $(\varphi_{n_k})_k$ converging in $L^2(\Omega)$. 
Which is impossible since by orthogonality 
\begin{align*}
 \| \varphi_{n_{k+1}}- \varphi_{n_k} \|^2_{L^2(\Omega)} = \| \varphi_{n_{k+1}} \|^2_{L^2(\Omega)}+\|\varphi_{n_k} \|^2_{L^2(\Omega)}=2\quad\text{for all}\quad k\geq0.
\end{align*}
Consequently the assumption that $B<\infty$ is impossible hence 
$\lim\limits_{n\to \infty}\mu_n= B=\infty$. We have proved $(ii)$. Now to show $(iii)$ we solely need to show that the span of $(\varphi_n)_{n\in \mathbb{N}_{0}}$ is dense in $L^2(\Omega)$. Let $f\in L^2(\Omega)$. For $\varepsilon>0$ small enough consider $g\in C_c^\infty(\mathbb{R}^d)\subset \VnuOm$ such that $\|f-g\|_{L^2(\Omega)}< \varepsilon$. Since the elements $(\varphi_n)_{n\in \mathbb{N}_0}$ are mutually orthogonal, the Bessel inequality tells us that 
 \begin{align*}
 \Big|\sum_{k=0}^{n} (f-g, \varphi_k)_{L^2(\Omega)}\varphi_k\Big|\leq \|f-g\|_{L^2(\Omega)}<\varepsilon.
 \end{align*}
 Therefore, we have the following estimate for all $n\geq 1$
 \begin{align}
 \Big\|f -\sum_{k=0}^{n} (f, \varphi_k)_{L^2(\Omega)}\varphi_k\Big\|_{L^2(\Omega)}
 &\leq \|f-g\|_{L^2(\Omega)} + \Big\|g -\sum_{k=0}^{n} (g, \varphi_k)_{L^2(\Omega)}\varphi_k\Big\|_{L^2(\Omega)}\notag\\
 &+ \Big\|\sum_{k=0}^{n} (f-g, \varphi_k)_{L^2(\Omega)}\varphi_k\Big\|_{L^2(\Omega)}\notag \\
 &<2\varepsilon +\Big\|g -\sum_{k=0}^{n} (g, \varphi_k)_{L^2(\Omega)}\varphi_k\Big\|_{L^2(\Omega)}\label{eq:wapprox}
 \end{align}
%
 
 \noindent For  $n\geq 1$, consider $g_n = g-\sum\limits_{k=0}^{n} (g, \varphi_k)_{L^2(\Omega)}\varphi_k.$
Since $g\in \VnuOm$, by  definition of $(\mu_j, \varphi_j)$ we have 
 $$\mathcal{E}(g, \varphi_j) = \mu_j (g,\varphi_j)_{L^2(\Omega)} \qquad\text{for all}~~~ 0\leq j\leq n.$$
 Also by orthogonality, we have
 $$\mathcal{E}(\varphi_k, \varphi_j) = \mu_k (\varphi_k,\varphi_j)_{L^2(\Omega)} = \mu_k \delta_{k,j} \qquad\text{for all} ~~~ 0\leq j,k\leq n.$$
 Whence, for $0\leq j\leq n$, combining the above relations we get 
 \begin{align*}
 \mathcal{E}(g_n, \varphi_j) &= \mathcal{E}(g, \varphi_j)- \sum_{k=0}^{n} (g, \varphi_k)_{L^2(\Omega)}\mathcal{E}(\varphi_k, \varphi_j) \\
 &= \mu_j(g,\varphi_j)_{L^2(\Omega)}- \sum_{k=0}^{n}\mu_k 
 (g, \varphi_k)_{L^2(\Omega)}\delta_{k,j}\\
 &=\mu_j(g,\varphi_j)_{L^2(\Omega)}-\mu_j(g,\varphi_j)_{L^2(\Omega)} =0
 \end{align*}

 \noindent which yields $\mathcal{E}(g_n, g-g_n)= \mathcal{E}\big(g_n, \sum\limits_{j=0}^{n} (g, \varphi_j\big)_{L^2(\Omega)}\varphi_j) =0.$
 %
 %
Therefore we get 
 \begin{align}\label{eq:bessel-inequality}
 \mathcal{E}(g,g) = \mathcal{E}(g_n,g_n)+ \mathcal{E}(g-g_n,g-g_n)\geq \mathcal{E}(g_n,g_n).
 \end{align}
 Again by orthogonality we have 
 \begin{align*}
 (g_n, \varphi_j)_{L^2(\Omega)} &= (g, \varphi_j)_{L^2(\Omega)}- \sum_{k=0}^{n} (g, \varphi_k)_{L^2(\Omega)}(\varphi_k, \varphi_j)_{L^2(\Omega)} \\
 &= (g,\varphi_j)_{L^2(\Omega)}- \sum_{k=0}^{n} 
 (g, \varphi_k)_{L^2(\Omega)}\delta_{k,j}\\
 &= (g,\varphi_j)_{L^2(\Omega)}-(g,\varphi_j)_{L^2(\Omega)} =0.
 \end{align*}
 So that $g_n\in V_{n+1} = \{\varphi_0,\varphi_1, \cdots, \varphi_n\}^\perp$. By the expression of $\mu_{n+1}$ together with \eqref{eq:bessel-inequality} we get
 \begin{align*}
 \mu_{n+1} \|g_n\|^2_{L^2(\Omega)} \leq \mathcal{E}(g_n,g_n) \leq \mathcal{E}(g,g).
 \end{align*}
 Since this holds for arbitrarily chosen $n\geq 1$ the fact $\mu_{n+1}\xrightarrow{n\to \infty}0$ forces $\|g_n\|^2_{L^2(\Omega)}\xrightarrow{n\to \infty}0$ that is
 \begin{align*}
 \lim_{n\to \infty} \big\|g-\sum_{k=0}^{n} (g, \varphi_k)_{L^2(\Omega)}\varphi_k\big\|^2_{L^2(\Omega)}=0.
 \end{align*}

\noindent Inserting this in our initial estimate \eqref{eq:wapprox} yields 
\begin{align*}
 \lim_{n\to \infty} \big\|f-\sum_{k=0}^{n} (f, \varphi_k)_{L^2(\Omega)}\varphi_k\big\|^2_{L^2(\Omega)}\leq 2\varepsilon.
 \end{align*}
 Wherefore letting $\varepsilon\to 0$ leads to the relation
\begin{align*}
 f = \sum_{n=0}^{\infty} (f, \varphi_n)_{L^2(\Omega)} \varphi_n \qquad\textrm{in the $L^2(\Omega)$ sense}. 
 \end{align*}
 This gives precisely the density result we are striving for and the proof is completed. 
\end{proof}

\medskip

\noindent By analogous arguments as in Theorem \ref{thm:rayleigh-construction} one is able to prove the following similar result when the operator $L$ is subject to Dirichlet or Robin boundary condition. 

\begin{theorem}\label{thm:rayleigh-construction-Dirichlet}
Let $\Omega\subset \mathbb{R}^d$ be a bounded open set and let the function $\nu: \mathbb{R}^d\to [0,\infty]$  be the density of a symmetric L\'{e}vy measure. Assume the conditions \eqref{eq:non-integrability-condition} and \eqref{eq:integrability-condition-near-zero} (see page  \pageref{eq:non-integrability-condition}) hold. Then there exist a sequence of positive real numbers $(\lambda_n)_{n\in\mathbb{N}}$ and a sequence of functions $(\psi_n)_{n\in\mathbb{N}}$ in $\VnuOmO$ such that for each $n\geq1$ the pair $(\lambda_n, \psi_n)$ satisfies the relations
\begin{align}\label{eq:Rayleigh-algorithm-Dirichlet}
 \lambda_n =\mathcal{E}(\psi_n,\psi_n) =\min_{v\in V_n}\{\mathcal{E}(v,v):~ \|v\|_{L^2(\Omega)}=1\}= \inf_{v\in V_n}\Big\{ \frac{\mathcal{E}(v,v)}{ \|v\|^2_{L^2(\Omega )}}~ \Big\}. 
\end{align}
\noindent Here we denote $V_1 = \VnuOmO$ and for $n\geq 2$, 
\begin{align*}
 V_n = \{\psi_1,\cdots,\psi_{n-1}\}^\perp = \{v\in \VnuOmO : (v, \psi_i)_{L^2(\Omega)} =0,~~i=1,\cdots,n-1\}.
\end{align*}
Moreover, the following assertions hold true. 

\begin{enumerate}[$(i)$]
 \item 
Each $\lambda_n$ is a Dirichlet eigenvalue of the operator $L$ constrained with the Dirichlet boundary condition $u =0$ on $\Omega^c$ and the associated eigenfunction $\psi_n.$
 \item $0<\lambda_1\leq \cdots\leq \lambda_n\leq \cdots $ and $\lim\limits_{n\to \infty}\lambda_n =\infty$. 
 \item The family $(\psi_n)_{n\in \mathbb{N}}$ is an orthonormal basis of $L^2(\Omega)$. 
\end{enumerate}
\end{theorem}

\medskip

\begin{theorem}\label{thm:rayleigh-construction-Robin}
Let $\Omega\subset \mathbb{R}^d$ be a bounded open set and let the function $\nu: \mathbb{R}^d\to [0,\infty]$
be the density of a symmetric L\'{e}vy measure with full support. Assume that the 
couple $(\nu, \Omega)$ belongs to one of the class $\mathscr{A}_i,~i=1,2,3$ (see page  \pageref{eq:class-lipschitz}). Let $K\subset \Omega $
be a measurable set such that $|K|>0$ and let $\beta: \Omega^c\to [0, \infty]$ 
be a measurable function such that $\beta\nu_K>0~~a.e$ on a set of positive measure.
Then there exist a sequence of positive real numbers $(\gamma(\beta)_n)_{n\in\mathbb{N}}$ and a sequence of functions $(\chi_n)_{n\in\mathbb{N}}$ in $\VnuOm$ such that for each $n\geq1$ the pair $(\chi_n, \gamma(\beta)_n)$ satisfies the relations 
\begin{align}\label{eq:Rayleigh-algorithm-Robin}
 \gamma(\beta)_n &=Q_\beta(\chi_n,\chi_n) =\min_{v\in V_n}\{Q_\beta(v,v):~ \|v\|_{L^2(\Omega)}=1\}= \inf_{v\in V_n}\Big\{ \frac{Q_\beta(v,v)}{ \|v\|^2_{L^2(\Omega )}}~ \Big\}
 \intertext{with}
&Q_\beta(u,v)= \mathcal{E}( u,v)+ \int_{\Omega^c} u(y)v(y)\beta(y)\nu_K(y)\d y.\notag 
\end{align}

\noindent Here we denote $V_1 = \VnuOm$ and for $n\geq 2$, 
\begin{align*}
 V_n = \{\chi_1,\chi_2,\cdots,\chi_{n-1}\}^\perp = \{v\in \VnuOm : (v, \chi_i)_{L^2(\Omega)} =0,~~i=1,\cdots,n-1\}.
\end{align*}
Moreover the following assertions hold true. 

\begin{enumerate}[$(i)$]
 \item 
Each $\gamma(\beta)_n$ is a Robin eigenvalue of the operator $L$ constrained with the Robin boundary condition $\mathcal{N} u+ \beta\nu_K u =0$ on $\Omega^c$ and the associated eigenfunction is $\chi_n$. 
 \item $0< \gamma(\beta)_1\leq \cdots\leq \gamma(\beta)_n\leq \cdots $ and $\lim\limits_{n\to \infty}\gamma(\beta)_n =\infty$. 
 \item The family $(\chi_n)_{n\in \mathbb{N}}$ is an orthonormal basis of $L^2(\Omega)$. 
\end{enumerate}
\end{theorem}

\begin{remark}
\noindent Recall that it is possible to replace $\nu_K$ by $\mathring{\nu}_K$ in the foregoing. It is worthwhile mentioning that the constants $\mu_1>0$ and $\lambda_1>0$ respectively satisfy the following Poincar\'e type inequalities. 
\begin{align*}
\mathcal{E}(u,u)\geq \mu_1\|u\|^2_{L^2(\Omega)}, \,\,\text{for all}~ u \in \VnuOm^\perp
\quad\text{and}\quad
\mathcal{E}(u,u)\geq \lambda_1\|u\|^2_{L^2(\Omega)}, \,\,\text{for all}~ u \in V_\nu^\Omega(\Omega|\mathbb{R}^d).
\end{align*}

\noindent In addition,  the following comparison holds
\begin{align*}
\mu_{n-1} \leq \gamma_n(\beta)\leq \lambda_n, \quad \text{for all $n\geq 1$}. 
\end{align*}
In short, one merely writes Neumann $\leq$ Robin $\leq $ Dirichlet. 
\end{remark}

\section{Helmholtz equation for nonlocal operators}
There are several applications of spectral decomposition. We shall illuminate this through the Helmholtz type equations subject with  the Neumann complement condition. The other types of complement  conditions can be derived likewise. The Helmholtz equation with the Neumann  condition  is of  the form 
\begin{align}\label{eq:nonlocal-Neumann-Helmotz}
L u -\lambda u= f \quad\text{in}~~~ \Omega \quad\quad\text{and} \quad\quad \mathcal{N} u= 0 ~~ \text{on}~~ \mathbb{R}^d\setminus\Omega.
\end{align} 

\noindent Here $\lambda$ is a potential and for brevity we assume that $\lambda\in \R$. Define $\mathcal{E}_{-\lambda}(u,v) = \mathcal{E}(u,v)-\lambda (u,v)_{L^2(\Omega)}.$ Obviously, the weak formulation of \eqref{eq:nonlocal-Neumann-Helmotz}  consists of finding $u\in\VnuOm$ such that 
\begin{align}\label{eq:var-nonlocal-Neumann-Helmotz}
\mathcal{E}_{-\lambda}(u,v) = \int_\Omega f(x) v(x)\d x\quad\text{for all $v\in\VnuOm$}. 
\end{align} 

\noindent To study \eqref{eq:var-nonlocal-Neumann-Helmotz},  we will need the so called Fredholm alternative. 

\begin{theorem}[Appendix D, Theorem 5, \cite{Ev10}]\label{thm:fredhoml-alternative}
	Let $K:H\to H$ be a compact operator on a Hilbert space $(H , (\cdot, \cdot)_H)$. The following assertions are true. 
	\begin{enumerate}[$(i)$]
		\item $\ker(I-K)$ has finite dimension and $\operatorname{Im}(I-K)$ is closed. 
		\item If $K^*$ denotes the adjoint of $K$ then $\ker(I-K)= \operatorname{Im}(I-K^*)^\perp $ and $\dim \operatorname{Im}(I-K)= \dim \operatorname{Im}(I-K^*)$. 
		\item $\ker(I-K)=\{0\}$ if and only if $\operatorname{Im}(I-K)=H$. 
	\end{enumerate}
\end{theorem}

\medskip 

\begin{theorem}
	Assume $\Omega\subset \mathbb{R}^d$ is a bounded open set and the function $\nu: \mathbb{R}^d\to [0,\infty]$ is the density of a symmetric L\'{e}vy measure with full support such that the couple $(\nu, \Omega)$ belongs to one of  the class $\mathscr{A}_i,~i=1,2,3 $ (see page  \pageref{eq:class-lipschitz}).  Let $0=\mu_0<\mu_1\leq \mu_2\leq \cdots\leq \mu_n\leq \cdots$ be the eigenvalues of $L$ subject to the Neumann condition (c.f Theorem \ref{thm:existence-of-eigenvalue-Neumann}).  Let $ f \in L^2 (\Omega)$. The following assertions are true. 
	\begin{enumerate}[$(i)$]
	\item If $\lambda  \in \R^d\setminus\{\mu_n: n\in \mathbb{N}_0\}$ there exists a unique function $u\in\VnuOm^\perp$ such that 
\begin{align}\label{eq:var-nonlocal-Neumann-Helmotz-bis}
\mathcal{E}_{-\lambda}(u,v) =\int_\Omega f(x) v(x)\d x\quad\text{for all $v\in \VnuOm^\perp$}. 
\end{align}

\noindent Moreover, there is a constant $C: = C(d,\Omega, \nu, \lambda)>0$ such that
	\begin{align}\label{eq:weak-regular-Helmotz}
	\|u\|_{\VnuOm}\leq C \|f\|_{L^2(\Omega)}\,.
	\end{align}
	\item Let  $j\geq 1$. Assume $\lambda= \mu_j$ and let  $\varphi_{j}, \varphi_{j+1},\cdots, \varphi_{j+r_j}\in \VnuOm^\perp$ be the linear independent eigenfunctions corresponding to the eigenvalue $\mu_j$. Consider the direct decomposition $\VnuOm^\perp= W_j\oplus W_j^\perp$ where $	W_j= \operatorname{span}\{\varphi_j,\varphi_{j+1}, \cdots, \varphi_{j+r_j}\}$ and $W_j^\perp=\{v\in \VnuOm^\perp: (v,\varphi_j)_{L^2(\Omega)}=\cdots=(v,\varphi_{j+r_j})_{L^2(\Omega)}=0\}.$
Then there exists a unique $u_j\in W_j^\perp$ such that 
	\begin{align}
	\mathcal{E}_{-\mu_j}(u_j,v) = \int_\Omega f(x) v(x)\d x\quad\text{for all $v\in W_j^\perp$}. 
	\end{align} 
	\end{enumerate}

\noindent Moreover, there exists $C: = C(d,\Omega, \nu, \mu_j)>0$  independent of $f$ such that the following estimate holds 
\begin{align}\label{eq:weak-regular-Helmotz-bis}
\|u_j\|_{\VnuOm}\leq C \|f\|_{L^2(\Omega)}\,.
\end{align}
Furthermore, the problem \eqref{eq:var-nonlocal-Neumann-Helmotz-bis} has a solution if and only if $ (f,\varphi_j)_{L^2(\Omega)}=\cdots=(f,\varphi_{j+r_j})_{L^2(\Omega)}=0$. In this case, all solutions to the problem \eqref{eq:var-nonlocal-Neumann-Helmotz-bis}  are of the form $u_j+\varphi$ for any $\varphi\in W_j$. 
\end{theorem}

\medskip 
\begin{proof}
$(i)$ Consider the linear operator $\Phi_0 : L^{2}(\Omega) \to\VnuOm^\perp $ mapping $f \in L^{2}(\Omega)$ to $ u=  \Phi_0  f\in\VnuOm^\perp $ the unique solution of the problem \eqref{eq:var-nonlocal-Neumann-bis} with $g=0$, i.e. 
\begin{align*}
\mathcal{E}(u,v) = \int_\Omega f(x) v(x)\d x\quad\text{for all $v\in \VnuOm^\perp$}. 
\end{align*}
\noindent As in the proof of Theorem \ref{thm:existence-of-eigenvalue-Neumann}, we know that the operator $S= R_\Omega\circ \Phi_0: L^{2}(\Omega)  \to L^2(\Omega)^\perp$ is compact.  Let $w\in \ker(I- \lambda S)$  that is $\lambda Sw=  w$. Equivalently, we have
\begin{align*}
\mathcal{E}(\Phi_0w,v) = \lambda\int_\Omega \Phi_0w(x) v(x)\d x\quad\text{for all $v\in \VnuOm^\perp$}. 
\end{align*}

\noindent Necessarily, $w=0$ since by assumption, $\lambda$ is not a Neumann eigenvalue. Hence $\ker(I- \lambda S)=\{0\}$. According to Theorem \ref{thm:fredhoml-alternative} (iii), we have 
$\operatorname{Im}(I-\lambda S) = L^2(\Omega)$. Hence for $f\in L^2(\Omega)$, we can write $f= w-\lambda Sw= w-\lambda R_\Omega\circ \Phi_0 w$ for some $ w\in L^2(\Omega)$.  We have $w = f+ \lambda R_\Omega\circ \Phi_0 w $. Let $v\in \VnuOm^\perp$, by definition of $ \Phi_0 w$ we have
\begin{align*}
\mathcal{E}(\Phi_0 w,v) = \int_\Omega w(x) v(x)\d x=  \int_\Omega f(x) v(x)\d x+ \lambda \int_\Omega \Phi_0w(x) v(x)\d x.
\end{align*}
Therefore the function $u= \Phi_0 w$, satisfies \eqref{eq:var-nonlocal-Neumann-Helmotz-bis}. Furthermore, the uniqueness of $u$ is due to the fact that $\lambda$ is not a Neumann eigenvalue. 
Next, consider the  operator $T_\lambda : L^{2}(\Omega) \to\VnuOm^\perp $ such that for $f\in L^2(\Omega)$, $u=T_\lambda f$ solves \eqref{eq:var-nonlocal-Neumann-Helmotz-bis}. Let $(f_n)_n$ be a sequence in $L^2(\Omega)$ tending to zero such that $u_n= \Phi f_n$ converges to $u$ in $\VnuOm$.  Let $v\in \VnuOm^\perp$ then for $n\geq 1,$ we have 
\begin{align*}
\mathcal{E}(u_n,v)-\lambda\int_\Omega u_n(x)v(x)\d x =\int_\Omega f_n(x) v(x)\d x. 
\end{align*}
Letting $n\to \infty$, this amounts to $\mathcal{E}_{-\lambda}(u,v)=0$. 
%
 Necessarily $u=0$ since by assumption $\lambda$ is 
 not  a Neumann eigenvalue. Therefore, in virtue of the closed graph Theorem $T_\lambda$ is a bounded operator.
 There exists a constant $C>0$ such that 
 \begin{align*}
 	\|T_\lambda h\|_{\VnuOm}\leq C\|h\|_{L^2(\Omega)}\quad\text{for all $h\in L^2(\Omega)$.}
 \end{align*}
  \noindent Hence the estimate \eqref{eq:weak-regular-Helmotz} follows. 
 
\vspace{2mm}

\noindent $(ii)$ Assume $\lambda= \mu_j$.  The procedure is analogous to the proof of $(i)$.   Consider the linear operator $\Phi_j: L^{2}(\Omega) \to W_j^\perp$ mapping $f \in L^2(\Omega)$ to $ u_j=  \Phi_j  f $ the unique function in $W_j^\perp$ satisfying
 \begin{align*}
 \mathcal{E}(u_j,v) = \int_\Omega f(x) v(x)\d x\quad\text{for all $v\in W_j^\perp$}. 
 	\end{align*}
\noindent This makes sense since $W_j^\perp$ is a closed subspace of $\VnuOm^\perp$ and hence Lax-Milgram's lemma can be applied on $W^\perp_j$. Similarly, we can show that $S= R_\Omega\circ \Phi_j: L^2(\Omega)  \to L^2(\Omega)^\perp$ is compact. Moreover, $\ker(I-\mu_j S)= \{0\}$. Indeed, if $w\in \ker(I-\mu_j S)$ then since $\mu_j$ is a Neumann eigenvalue we can easily show that 
$w=\lambda S w\in W_j\cap W_j^\perp= \{0\}$. 
Again by Fredholm alternative (Theorem \ref{thm:fredhoml-alternative} $(iii)$) for $f\in L^2(\Omega)$ there exists $u_j\in W_j$ such that $f= u_j-\mu_j S u_j$ and $u_j$ is unique since $W_j\cap W_j^\perp= \{0\}$. We get that
\begin{align*}
\mathcal{E}_{-\mu_j}(u_j,v) =\int_\Omega f(x) v(x)\d x \quad\text{for all $v\in W_j$}. 
\end{align*}
 One can also show the estimate \eqref{eq:weak-regular-Helmotz-bis} by applying the closed graph theorem as above. Let us recall that $W_j$ is the eigenspace associated with the eigenvalue $\mu_j$. If $u\in \VnuOm^\perp$ solves \eqref{eq:var-nonlocal-Neumann-Helmotz-bis} then since $\mathcal{E}_{-\mu_j}(u,\varphi)= 0$ for all $\varphi\in W_j$. In particular this implies that $\big(f, \varphi_j \big)_{L^2(\Omega)}=,\cdots,= \big(f, \varphi_{j+r_j} \big)_{L^2(\Omega)} =0$. The converse remains true. Indeed, every  $\varphi\in W_j$  is an eigenfunction associated with the eigenvalue $\mu_j$, i.e. $\mathcal{E}_{-\mu_j}(\varphi,v) =0 \quad\text{for all $v\in \VnuOm^\perp$}.$ Accounting the fact that $\VnuOm^\perp= W_j\oplus W_j^\perp $ one easily show that 
 \begin{align*}
\mathcal{E}_{-\mu_j}(u_j+\varphi, v) =0 \quad\text{for all $v\in \VnuOm^\perp$} .
 \end{align*}
Conversely, if $u$ solves \eqref{eq:var-nonlocal-Neumann-Helmotz-bis} then  one can easily show that $u-u_j\in W_j$ and $f\in W_j^\perp$. Thus all solutions of \eqref{eq:var-nonlocal-Neumann-Helmotz-bis} are of the form $u_j+\varphi$ with $\varphi\in W_j$ if and only if $f\in W_j^\perp$. 
\end{proof}

\vspace{2mm}

\begin{remark}
It is easy to deduce solutions of the problem \eqref{eq:var-nonlocal-Neumann-Helmotz} form those of the problem \eqref{eq:var-nonlocal-Neumann-Helmotz-bis}. Note that the situation where $\lambda=0$ is already encoded in the variational problem \eqref{eq:var-nonlocal-Neumann}. Furthermore, in a similar manner one is able to treat the Helmholtz type equation associated with other types of complement value conditions.
\end{remark}

%
\section{Profiling solutions of evolution of IDEs }
In this section the spectral decomposition also plays an essential  role. We do some accustomed heuristic on certain IDEs with time evolution. 
Our exposition here is mostly for illustrative purposes, so we will not give precise arguments.
 We shall only focus on Neumann complement type condition.  Assume $\Omega\subset \R^d$ is an open set and let $0< T\leq\infty$. Let us consider the heat equation with complement Neumann condition, 
\begin{align}\label{eq:parabolic-Neumann-first}
\partial_t u + Lu = f ~~\text{in $\Omega\times [0, T)$},
\quad \mathcal{N} u =g ~~\text{on $\Omega^c\times [0, T)$},
\quad u=u_0~~ \text{on  $\Omega\times \{0\}$}. 
\end{align}
 Here $f,g: \Omega\times [0,T)\to \R$ and $u_0: \R^d\to \R$ are given data. In practice $f$ is called the forcing term and $u_0$ is the initial condition. For  the sake of simplicity, from now on we assume that $g=0$, i.e. 
 \begin{align}\label{eq:parabolic-Neumann}
 \partial_t u + Lu = f ~~\text{in $\Omega\times [0, T)$},
 \quad \mathcal{N} u =0 ~~\text{on $\Omega^c\times [0, T)$},
 \quad u=u_0~~ \text{on  $\Omega\times \{0\}$}.
 \end{align}
 According to \cite{DROV17} the problem \eqref{eq:parabolic-Neumann} encodes the evolution of a L\'evy process
  whose jumps from $x$ to $y$ are driven by $\nu(x-y)$ with the following features:
 \begin{enumerate}
 	
 \item $u(x,t)$ represents the probability distribution of the position of a L\'evy process moving randomly inside $\Omega$. In other words the generator of this process is $-L$ (See Chapter \ref{chap:basics-nonlocal-operator}). Thus $u(x,t)$ satisfies 
 \begin{align*}
 	 \partial_t u + Lu = 0 ~~\text{in $\Omega\times [0, T)$}. 
 \end{align*}
 
\item When the process exits $\Omega$, it instantaneously comes back to $\Omega$ almost surely.  If the process sites at the position $x\in \Omega^c$ at the time $t>0$, it will jump back at the same time to some position $y\in \Omega$ with probability $1$ according to the jump intensity $\nu(x-y)$. Thus at the position $x\in \Omega^c$ and at the time $t$ the probability $u(x,t)$ is the sum of all possible positions that process may occupy in $\Omega$ at time $t$. This results with the following relation
\begin{align*}
u(x,t) = \kappa(x)\int_\Omega u(y,t)\nu(x-y)\d y.
\end{align*}

\noindent  The constant $\kappa$ ensures that the function $y\mapsto \kappa(x)\nu(x-y)$ is a probability density on $\Omega$. That is $\kappa^{-1}(x) =\int_\Omega \nu(x-y)\d y.$ 
Therefore the above relation becomes 
 	\begin{align*}
 	\int_\Omega u(x,t)\nu(x-y)\d y = \int_\Omega u(y,t)\nu(x-y)\d y~~\text{equivalently,}~~ \mathcal{N}u(x,t)=0. 
 	\end{align*}
 	\item The process starts in $\Omega$ at time $t=0$ with some  probability distribution $u_0$. 
 \end{enumerate}
In sum the probability distribution solves the equation \eqref{eq:parabolic-Neumann} with $f=0$, i.e.
\begin{align*}
\partial_t u + Lu = 0 ~~\text{in $\Omega\times [0, T)$}
\quad \mathcal{N} u =0 ~~\text{on $\Omega^c\times [0, T)$}
\quad u=u_0~~ \text{on  $\Omega\times \{0\}$}.
\end{align*}
\noindent In practice, it would makes sense to call such a process \textit{gravity process} since $\Omega$ it acts on the process like an attractor. Another possible name could be  \textit{reflected jumps process} or \textit{pull-back process}.

\medskip

 \begin{remark}
 	Assume the above process is modified in such a way that it terminates in $\Omega^c$ upon exiting 
 	$\Omega$ instead of coming back then the corresponding complement value condition is given by $u(x,t) = 0$ 
 	on $[0,T)\times\Omega^c$ (Dirichlet condition). This process is well-known as \textit{a 
killed process} and  its probability described by the following equation
\begin{align*}
\partial_t u + Lu = 0 ~~\text{in $\Omega\times [0, T)$},
\quad u =0 ~~\text{on $\Omega^c\times [0, T)$},
\quad u=u_0~~ \text{on  $\Omega\times \{0\}$}.
\end{align*}
 \end{remark}
 \medskip
 
\noindent Now we want to analyze the equation \eqref{eq:parabolic-Neumann} heuristically. We adopt the notation $v(t)=v(\cdot,t)$. 
 \begin{definition} A function $u: \R^d\times  [0,T)\to \R$ is called to be a weak solution of the evolution problem
 	 \eqref{eq:parabolic-Neumann} if $u\in L^2(0,T; \VnuOm)$, $\partial_t u \in L^2(0,T; \VnuOm')$ and 
 	\begin{align}\label{eq:weak-parabolic-Neumann}
 		\big(\partial_t u(t) , v\big)_{L^2(\Omega)}+ \mathcal{E}(u(t), v)= 	\big(f(t) , v\big)_{L^2(\Omega)}\quad\text{for all $v\in \VnuOm$ and $t\in [0,T)$}.
 	\end{align}
 	Henceforward, the spaces $ L^2(0,T; L^2(\Omega))$, $ L^2(0,T; \VnuOm)$ and $L^2(0,T; (\VnuOm)')$ are understood in the Bochner sense. 
 \end{definition}

\medskip

\noindent Next, we are interested in establishing some basic properties of a weak solution of \eqref{eq:parabolic-Neumann} whenever it exists. The spectral decomposition of $L$ is one of the key points to do this.  Before, let us see what could be the profile of a weak solution.  Under the assumption of Theorem \ref{thm:existence-of-eigenvalue-Neumann}, assume $(\phi_n, \mu_n)_{n\in \mathbb{N}_0}$ constitutes the Neumann eigenpairs of the operator $L$. We emphasize that $(\phi_n)_n$ is an orthonormal basis of $L^2(\Omega)$ and $\phi_n\in \VnuOm$ for all $n\geq 0$. Recall that if $h\in L^2(\Omega)$ and $h_k= \big(h, \phi_k\big)_{L^2(\Omega)}$ stands for the $k^{\text{th}}$ Fourier coefficient of $h$ then the Fourier expansion and the Parseval identity of $h$ respectively give
\begin{align*}
h= \sum_{k=0}^\infty h_k\phi_k \quad\text{and}\quad \|h\|^2_{L^2(\Omega)}= \sum_{k=0}^\infty |h_k|^2.
\end{align*}

\noindent Assume $u_0\in \VnuOm$ and $f\in L^2(0,T; L^2(\Omega))$. Assume  $u\in L^2(0,T; \VnuOm)$ solves \eqref{eq:weak-parabolic-Neumann} such that $\partial_t u(t)\in L^2(\Omega)$ for all $t\in [0,T)$. By definition of $\phi_k$ we have $\mathcal{E}(u(t), \phi_k)= \mu_k \big(u(t),  \phi_k\big)_{L^2(\Omega)}= \mu_k u_k(t)$, thus testing \eqref{eq:weak-parabolic-Neumann} against $v=\phi_k$ amounts to the ordinary differential equation 
\begin{align*}
 u'_k(t) + \mu_k u_k(t) = f_k(t)\quad\text{equivalently}\quad  \big(e^{\mu_k t} u_k(t)\big)'= e^{\mu_k t} f_k(t),\qquad t\in [0,T).
\end{align*}

\noindent Integrating gives
\begin{align*}
	u_k(t) = e^{-\mu_k t}u_k(0)+ \int_0^t e^{\mu_k (s-t)} f_k(s)\d s.
\end{align*}
Therefore, since $u(0) =u_0$, $u_k(0)= \big(u_0(t), \phi_k\big)_{L^2(\Omega)}$ and $f_k(t)= \big(f(t), \phi_k\big)_{L^2(\Omega)}$ we have
\begin{align}\label{eq:parabolic-expression-Neumann}
u(x,t) = 
\sum_{k=0}^\infty \big(u_0(t), \phi_k\big)_{L^2(\Omega)} e^{-\mu_k t}\phi_k(x)+ \sum_{k=0}^\infty \phi_k(x)\int_0^t  e^{-\mu_k (t-s)}\big(f(s), \phi_k\big)_{L^2(\Omega)} \d s. 
\end{align}

\noindent  The expression \eqref{eq:parabolic-expression-Neumann} shows that $u(x,t)$ is uniquely defined whenever it exists. Further, it is worth emphasizing that the equality in \eqref{eq:parabolic-expression-Neumann} is understood in $L^2(\Omega)$ sense. On the other hand, assume $T<\infty$ then  for  $k\geq 1$ given that $0<\mu_1\leq \mu_k$ we find that 
\begin{align*}
|u_k(t)|^2\leq 2e^{-2\mu_1 t}|u_k(0)|^2+ 2T\int_0^T |f_k(t)|^2\d t. 
\end{align*}
Recalling that $\phi_0 = |\Omega|^{-1}$ and $\mu_0=0$, the Parseval identity together with the preceding estimate yield
\begin{align*}
\begin{split}
\big\|u(t) -\hbox{$\fint_\Omega u_0- \int_0^t\fint_\Omega f(s)\d s$}\big\|^2_{L^2(\Omega)}
&\leq 2e^{-2\mu_1 t}\|u_0-\hbox{$\fint_\Omega u_0$}\|^2_{L^2(\Omega)}+ 2T\int_0^T\|f(t)-\hbox{$\fint_\Omega f(t) $} \|^2_{L^2(\Omega)}\d t .
\end{split}
\end{align*} 
Meanwhile, using $v=1$ in \eqref{eq:weak-parabolic-Neumann} implies that 
\begin{align*}
\int_\Omega u(x,t) \d x= \int_\Omega u_0(x) \d x+ \int_0^t\int_\Omega f(x,s) \d x\d s\qquad\text{for all $0\leq t\leq T$}.
\end{align*}
Altogether we get
\begin{align}\label{eq:parabolic-estimate}
\begin{split}
\big\|u(t) -\hbox{$\fint_\Omega u(t)$}\big\|^2_{L^2(\Omega)}
&\leq 2e^{-2\mu_1 t}\|u_0-\hbox{$\fint_\Omega u_0$}\|^2_{L^2(\Omega)}+ 2T\int_0^T\|f(t)-\hbox{$\fint_\Omega f(t) $} \|^2_{L^2(\Omega)}\d t .
\end{split}
\end{align} 
It is worth emphasizing that the estimate \eqref{eq:parabolic-estimate} can also be established by using Poincar\'e inequality and Gronwall's lemma.  In particular, we have 
\begin{align}\label{eq:parabolic-estimate-uniform}
\begin{split}
\sup_{0\leq t\leq T}\big\|u(t) -\hbox{$\fint_\Omega u(t)$}\big\|^2_{L^2(\Omega)}
&\leq 2\|u_0-\hbox{$\fint_\Omega u_0$}\|^2_{L^2(\Omega)}+ 2T\int_0^T\|f(t)-\hbox{$\fint_\Omega f(t) $} \|^2_{L^2(\Omega)}\d t .
\end{split}
\end{align}

\noindent Exploiting the relation \eqref{eq:parabolic-Neumann} together with the Poincar\'e inequality whose sharp constant is given by $\mu_1^{-1}$,  leads to the following  estimate 
\begin{align}\label{eq:parabolic-estimate-V}
\begin{split}
\big\|u(t) -\hbox{$\fint_\Omega u(t)$}\big\|^2_{\VnuOm}
&\leq Ce^{-2\mu_1 t}\|u_0-\hbox{$\fint_\Omega u_0$}\|^2_{L^2(\Omega)}+ C\int_0^T\|f(t)-\hbox{$\fint_\Omega f(t) $} \|^2_{L^2(\Omega)}\d t .
\end{split}
\end{align} 
 \noindent Here  $C$ is a generic constant only depending on $T$ and $\mu_1$. Finally, \eqref{eq:parabolic-estimate-V} ,\eqref{eq:parabolic-estimate-uniform}  and \eqref{eq:weak-parabolic-Neumann} yield
 \begin{align}\label{eq:parabolic-estimate-global}
 \begin{split}
 \big\|u-\hbox{$\fint_\Omega u$}\big\|^2_{L^\infty(0,T;L^2(\Omega))}&+ \big\|u-\hbox{$\fint_\Omega u$}\big\|^2_{L^2(0,T;\VnuOm)}+ \big\|u-\hbox{$\fint_\Omega u$}\big\|^2_{L^2(0,T;(\VnuOm)')}\\
 &\leq C\|u_0-\hbox{$\fint_\Omega u_0$}\|^2_{L^2(\Omega)}+ C\|f-\hbox{$\fint_\Omega f $} \|^2_{L^2(0,T;L^2(\Omega))}.
 \end{split}
 \end{align} 

\noindent A rigorous analysis using some basics on Bochner integrals allows us to establish the existence and uniqueness of a weak solution solving \eqref{eq:weak-parabolic-Neumann} under appropriate choices of the initial and complement condition. In practice this analysis is based on the so called Galerkin approximation method which consists of projecting the problem on the finite dimensional Hilbert space generated by the $n^{\text{th}}$ first vectors $\phi_0\cdots, \phi_{n-1}$. The Galerkin approximation method also allows us to establish the estimate \eqref{eq:parabolic-estimate-global} more consistently. A curious readers may consult \cite{Ev10,Hunter14,LeDret16} for more about this theoretical approach. Roughly, we can summarize the well-posedness of \eqref{eq:weak-parabolic-Neumann} as follows. 

\medskip

\begin{theorem}
	Under the assumption of Theorem \ref{thm:existence-of-eigenvalue-Neumann}, i.e $(\nu, \Omega)$ belongs 
	to one of the class $\mathscr{A}_i, ~i=1,2,3$ (see page  \pageref{eq:class-lipschitz}). Then for $f\in L^2(0,T; L^2(\Omega))$ and $u_0\in \VnuOm$ there exists a unique solution  $u\in C(0,T; L^2(\Omega))\cap L^2(0,T;\VnuOm)$ of \eqref{eq:weak-parabolic-Neumann} such that $\partial_t u \in L^2(0,T; (\VnuOm)')$. Here we assume $T<\infty$. Moreover, the estimate \eqref{eq:parabolic-estimate-global} holds true. 
\end{theorem}

\medskip 

\noindent One is able to establish an analogous result replacing the  Neumann complement condition with the Dirichlet complement condition.  Next, we  would like to investigate the long time behavior of $u(x,t)$ upon assuming $f=0$ and $T=\infty$.  That is $\partial_t u+ Lu =0$ in $\Omega\times [0,\infty)$, $\mathcal{N} u =0$ on $\Omega^c\times [0,\infty)$ and $u=u_0$ on $\R^d\times \{0\}$. Thus \eqref{eq:weak-parabolic-Neumann} reduces to 
	\begin{align}\label{eq:weak-parabolic-Neumann-0}
\big(\partial_t u(t) , v\big)_{L^2(\Omega)}+ \mathcal{E}(u(t), v)=0 \quad\text{for all $v\in \VnuOm$ and $t\in [0,T)$}.
\end{align}
To avoid technicality, we assume $(x,t)\mapsto u(x,t)$ is sufficiently regular so that the interchanges between integral sign and derivative are possible.

\noindent \textbf{Mass conservation:} As a consequence of \eqref{eq:weak-parabolic-Neumann-0} with $v=1$ it follows that 
\begin{align*}
\int_\Omega u(x,t)\d x=\int_\Omega u_0(x)\d x\qquad\text{for all $t\geq 0$}.
\end{align*}
Now taking $v=u(t)$ we get that $\frac{d}{dt} \|u(t)\|^2_{L^2(\Omega)}=-\mathcal{E}(u(t), u(t))\leq 0$ for all $t\geq 0$ which implies
%

\begin{align*}
\|u(t)\|^2_{L^2(\Omega)}\leq \|u_0\|^2_{L^2(\Omega)}\qquad\text{for all $t\geq 0$}.
\end{align*}
\noindent \textbf{Energy conservation:} Assume $\partial_t u(t)\in \VnuOm$ and let the energy $E(t)= \mathcal{E}(u(t), u(t))$ then testing \eqref{eq:weak-parabolic-Neumann-0} with $v=\partial_t u(t)$ yields 
\begin{align*}
E'(t)= \frac{d}{dt}\mathcal{E}(u(t), u(t)) = 2\mathcal{E}(u(t), \partial_t u(t)) =-2\|\partial_t u(t)\|^2_{L^2(\Omega)}\leq 0. 
\end{align*}
Therefore, we have $E(t)\leq E(0)$, i.e. 
\begin{align*}
\mathcal{E}(u(t), u(t))\leq \mathcal{E}(u_0, u_0)\qquad\text{for all $t\geq 0$}.
\end{align*}
\noindent \textbf{Dissipation effect:} A close look at the proof the estimate \eqref{eq:parabolic-estimate} reveals that 
\begin{align}\label{eq:disp-effect}
\big\|u(t) -\hbox{$\fint_\Omega u_0$}\big\|^2_{L^2(\Omega)}
&\leq e^{-2\mu_1 t}\|u_0-\hbox{$\fint_\Omega u_0$}\|^2_{L^2(\Omega)}\xrightarrow{t\to \infty}0. 
\end{align}
In other words $u(\cdot,t)$ exponentially tends to $\fint_\Omega u_0$ in $L^2(\Omega)$ as $t\to \infty$.

\vspace{2mm}

\noindent \textbf{Stability at the equilibrium:} Now assume that in \eqref{eq:parabolic-Neumann-first}, $f$ and $g$ are independent of the time variable.  An equilibrium solution of \eqref{eq:parabolic-Neumann} is any solution which we denote $u_e$ such that $\partial_t u_e=0$ that  is $u_e$ satisfies the Neumann problem
\begin{align*}
Lu_e = f ~~\text{in $\Omega$}
\quad \text{and}\quad \mathcal{N} u_e=g ~~\text{on $\Omega^c$}. 
\end{align*}
Next assume $u$ is a solution to\eqref{eq:parabolic-Neumann-first} then $w(x,t)= u(x,t)-u_e(x)$ with $w_0(x)= u_0(x)-u_e(x)$ satisfies 

\begin{align}\label{eq:parabolic-Neumannt}
\partial_t w+ Lw= 0~~\text{in $\Omega\times [0, T)$},
\quad \mathcal{N} w =0 ~~\text{on $\Omega^c\times [0, T)$},
\quad w=w_0= u_0-u_e~~ \text{on  $\Omega\times \{0\}$}. 
\end{align}

\noindent From the estimate in \eqref{eq:disp-effect} we have

\begin{align}
\big\|u(t)-u_e -\hbox{$\fint_\Omega [u_0-u_e]$}\big\|^2_{L^2(\Omega)}
&\leq e^{-2\mu_1 t}\|u_0-u_e-\hbox{$\fint_\Omega [u_0-u_e]$}\|^2_{L^2(\Omega)}\xrightarrow{t\to \infty}0. 
\end{align}

\noindent Note that if $f\in L^2(\Omega)$ and $g\in L^2(\Omega^c, \nu_K^{-1})$ are compatible then, Theorem \ref{thm:nonlocal-Neumann-var} implies that for all $t>0$,
\begin{align}
\big\|u(t)-u_e -\hbox{$\fint_\Omega [u_0-u_e]$}\big\|_{L^2(\Omega)}
&\leq Ce^{-\mu_1 t} \big( \|u_0-\hbox{$\fint_\Omega u_0$}\|_{L^2(\Omega)}+ \|f\|_{L^2(\Omega)}+\|g\|_{L^2(\Omega^c,\nu_K^{-1})}\big).
\end{align}

\vspace{2mm}
\noindent Let us highlight other types of evolution IDEs.  The spectral decomposition can be further applied with a parallel analysis on nonlocal Schr\"odinger equation of the form 
\begin{align}\label{eq:schrodinger-Neumann}
i\partial_t u + Lu = f ~~\text{in $\Omega\times [0, T)$},
\quad \mathcal{N} u =0 ~~\text{on $\Omega^c\times [0, T)$},
\quad u=u_0~~ \text{on  $\Omega\times \{0\}$}.
\end{align}
If $u_0\in \VnuOm$ and $f\in L^2(0,T;L^2(\Omega))$ then a profile of solution of \eqref{eq:schrodinger-Neumann} is of the form 
\begin{align*}
u(x,t) = 
\sum_{k=0}^\infty \big(u_0(t), \phi_k\big)_{L^2(\Omega)} e^{i\mu_k t}\phi_k(x)-i\sum_{k=0}^\infty \phi_k(x)\int_0^t  e^{i\mu_k (t-s)}  \big(f(s), \phi_k\big)_{L^2(\Omega)}\d s. 
\end{align*}
 
 \noindent Note that if $f=0$ then in this case we have mass conservation $\|u(t)\|_{L^2(\Omega)}= \|u_0\|_{L^2(\Omega)}$ for all $t\geq 0$.  We also point out the nonlocal hyperbolic equation of the form
 \begin{align}\label{eq:wave-Neumann}
 \partial^2_{tt} u + Lu = f ~~\text{in $\Omega\times [0, T)$},
 \quad \mathcal{N} u =0 ~~\text{on $\Omega^c\times [0, T)$},
 \quad u=u_0\,\,\text{and}\,\, \partial_t u=u_1~~ \text{on  $\Omega\times \{0\}$}.
 \end{align}
 If $u_0,u_1\in \VnuOm$ and $f\in L^2(0,T;L^2(\Omega))$ then a profile of solution of \eqref{eq:wave-Neumann} is of the form 
 \begin{align*}
 u(x,t) &= \sum_{k=0}^\infty \Big[\cos(\sqrt{\mu_k} t)\big(u_0(t), \phi_k\big)_{L^2(\Omega)}+ \frac{\sin(\sqrt{\mu_k} t)}{\sqrt{\mu_k}}\big(u_1(t), \phi_k\big)_{L^2(\Omega)} \Big] \phi_k(x)\\
 &-\sum_{k=0}^\infty \phi_k(x)\int_0^t \frac{\sin(\sqrt{\mu_k}(t-s))}{\sqrt{\mu_k}}\big(f(s), \phi_k\big)_{L^2(\Omega)} \d s. 
 \end{align*}

\section{Dirichlet-to-Neumann map for nonlocal operators}

In this section we wish to define the Dirichlet-to-Neumann map related to  the nonlocal operator $L$ under consideration.  After we show its eigenvalues are strongly connected to the Robin eigenvalues of the operator $L$. Our exposition here is largely influenced by \cite{WoRa07, WoRa12} where the Dirichlet-to-Neumann map is treated in the local setting for the Laplacian. However, we point out that an attempt to define the Dirichlet-to-Neumann map is proved in \cite{Zoran19}. For the case for the fractional Laplacian $L=(-\Delta)^{\alpha/2})$,  a sightly different Dirichlet-to-Neumann map to ours is derived in \cite{GSU16}. Let us start this section by recalling that if $f\in L^2(\Omega)$ and $g\in \TnuOm$, there exists a unique weak solution $u\in \VnuOm$ to the nonlocal Dirichlet problem $Lu=f$ in $\Omega$ and $u=g$ on $\Omega^c$, i.e.  we have $u=g$ a.e. on $\Omega^c$ and 
\begin{align}
\mathcal{E}(u,v) = \int_{\Omega} f(x) v(x)\d x\quad \text{for all }~~ v\in V_\nu^\Omega(\Omega|\mathbb{R}^d).
\end{align}
\noindent Moreover, there exists a constant $C>0$ such that 
\begin{align}\label{eq:dirichlet-regular}
\|u\|_{\VnuOm}\leq C(\|f\|_{L^2(\Omega)} + \|g\|_{\TnuOm}). 
\end{align}

\noindent  It is noteworthy recalling that under the conditions \eqref{eq:integrability-condition-near-zero} and \eqref{eq:non-integrability-condition} (see page  
\pageref{eq:non-integrability-condition}), by Theorem \ref{thm:rayleigh-construction-Dirichlet} there exist of a family of element $(\psi_n)_n$ elements of $V_\nu^\Omega(\Omega|\mathbb{R}^d)$, 
orthonormal basis of $L^2(\Omega)$ and an increasing sequence of real number $ 0<\lambda_1\leq \cdots\leq \lambda_n\leq \cdots $ such that $\lambda_n \to \infty$ as $n \to \infty$ and each $\psi_n$ is a Dirichlet eigenfunction of $L$ whose corresponding eigenvalue is $\lambda_n$ namely
\begin{align*}
\mathcal{E}(\psi_n,v) = \lambda_n \int_{\Omega} \psi_n(x) v(x)\d x\quad \text{for all }~~ v\in V_\nu^\Omega(\Omega|\mathbb{R}^d).
\end{align*}
%
\noindent Before we formally define the Dirichlet-to-Neumann map, some prerequisites are required. Let $f\in L^2(\Omega)$ and $g\in \TnuOm$. Assume $\lambda< \lambda_1$, then the bilinear form 

$$\mathcal{E}_{-\lambda}(u, u) = \mathcal{E}(u, u) -\lambda \|u\|^2_{L^2(\Omega)}$$ 
is coercive on $V_\nu^\Omega(\Omega|\mathbb{R}^d)$. Thus there exists a function $u \in V_\nu^\Omega(\Omega|\mathbb{R}^d)$ which is a unique weak solution to the Dirichlet problem $Lu - \lambda u = f$ in $\Omega$ and $u=g$ on $\Omega^c$. Explicitly, $u=g$ on $\Omega^c$ and 
\begin{align}\label{eq:weak-dirichlet-lambda}
\mathcal{E}(u,v)- \lambda \int_{\Omega} u(x) v(x)\d x = \int_{\Omega} f(x) v(x)\d x \quad \text{for all }~~ v\in V_\nu^\Omega(\Omega|\mathbb{R}^d).
\end{align}
Moreover, the estimate \eqref{eq:dirichlet-regular} (with the estimating constant depending on $\lambda$) remains true. More generally, by the mean of the Fredholm alternative and the closed graph Theorem, the preceding facts \eqref{eq:weak-dirichlet-lambda} and \eqref{eq:dirichlet-regular} respectively remain true for the operator $L -\lambda$, whenever $\lambda\in \mathbb{R}\setminus\{\lambda_n:n\geq1\}$. 

\medskip

\noindent From now on we suppose $f=0$ and $\lambda\in \mathbb{R}\setminus\{\lambda_n:n\geq1\}$ and label the solution of \eqref{eq:weak-dirichlet-lambda} by $u =u_g$. Then the mapping $g\mapsto u_g$ is linear and continuous from $\TnuOm$ to $\VnuOm$ since by \eqref{eq:dirichlet-regular} we have 
\begin{align*}
\|u_g\|_{\VnuOm}\leq C \|g\|_{\TnuOm}. 
\end{align*} 
Given $v \in \TnuOm$, put $\widetilde{v}= \operatorname{ext}(v)\in \VnuOm$ be an extension of $v$. 

\medskip

\begin{definition}
	Let $\lambda\in \mathbb{R}\setminus\{\lambda_n:n\geq1\}$. We call the Dirichlet-to-Neumann map with respect to the operator $L-\lambda$, the mapping $\mathscr{D}_{\lambda}: \TnuOm \to (\TnuOm)'$ defined by $g \mapsto \mathscr{D}_\lambda g= \mathcal{E}_{-\lambda}(u_g, \widetilde{\cdot})$ such that $\langle\mathscr{D}_\lambda g, v\rangle = \mathcal{E}_{-\lambda}(u_g, \widetilde{v})$.  Here $\langle\cdot, \cdot \rangle$ stands for the dual pairing between $\TnuOm$ and  $(\TnuOm)'$. 
\end{definition}

\medskip

\begin{theorem}\label{thm:DN-map}
	The Dirichlet-to-Neumann operator $\mathscr{D}_{\lambda}: \TnuOm \to (\TnuOm)'$ with $g \mapsto \mathscr{D}_\lambda g= \mathcal{E}_{-\lambda}(u_g, \widetilde{\cdot})$ is well defined, linearly bounded and self-adjoint. Moreover there exists $c\in \mathbb{R}$ such that for all $g \in \TnuOm$
	\begin{align*}
	\langle \mathscr{D}_\lambda g , g\rangle\geq c \|u_g\|^2_{\VnuOm}.
	\end{align*}
\end{theorem}

\medskip

\begin{proof}
Let $v' \in \VnuOm$ be another extension of $v$ thus $v'-\widetilde{v}\in V_\nu^\Omega(\Omega|\mathbb{R}^d)$. By definition of $u_g$ we get'
	\begin{align*}
	\mathcal{E}_{-\lambda}(u_g,v'- \widetilde{v})=0\quad \text{or equally}\quad \mathcal{E}_{-\lambda}(u_g, \widetilde{v})=\mathcal{E}_{-\lambda}(u_g,v').
	\end{align*}
	Therefore the mapping $v\mapsto\mathcal{E}_{-\lambda}(u_g, \widetilde{v})$ is well defined, linear and bounded on $\TnuOm$. Indeed, 
	\begin{align*}
	|\mathcal{E}_{-\lambda}(u_g, \widetilde{v})|\leq (|\lambda|+1)\|u_g\|_{\VnuOm} \|\widetilde{v}\|_{\VnuOm}.
	\end{align*}
	Since the extension $\widetilde{v}$ of $v$ is arbitrarily chosen, upon the estimate \eqref{eq:weak-dirichlet-lambda} we obtain
	\begin{align*}
	|\mathcal{E}_{-\lambda}(u_g, \widetilde{v})|\leq C\|g\|_{\TnuOm} \|v\|_{\TnuOm}.
	\end{align*}
	This shows that $\mathcal{E}_{-\lambda}(u_g, \widetilde{\cdot})$ belongs $(\TnuOm)'$. Subsequently it also 
	follows from this estimate that the mapping $\mathscr{D}_{\lambda}: \TnuOm \to (\TnuOm)'$ with $g \mapsto \mathscr{D}_\lambda g= \mathcal{E}_{-\lambda}(u_g, \widetilde{\cdot})$ is linearly bounded. Now let $g,h \in \TnuOm$, specializing the definition of $\mathscr{D}_\lambda$ with $\widetilde{g} = u_g$ and $\widetilde{h} =u_h$. The self-adjointness is obtained as follows 
	\begin{align*}
	\langle \mathscr{D}_\lambda g , h\rangle= \mathcal{E}_{-\lambda}(u_g, u_h) = \mathcal{E}_{-\lambda}(u_h, u_g)= \langle \mathscr{D}_\lambda h, g\rangle. 
	\end{align*}
	We end this proof by taking $c=\min(1,-\lambda)$, since $\langle \mathscr{D}_\lambda g , g\rangle= \mathcal{E}_{-\lambda}(u_g, u_g)\geq \min(1,-\lambda) \|u_g\|^2_{\VnuOm}.$
	
\end{proof}

\medskip

\begin{remark}
	The above definition is plainly motivated by the following observation. Assume $u_g$ be as before and let $ \varphi\in C_c^\infty(\mathbb{R}^d)$, the Green-Gauss formula \eqref{eq:green-gauss-nonlocal} gives 
	\begin{align}\label{eq:DN-weak-formula}
	\langle \mathscr{D}_\lambda g ,\varphi\rangle &= \mathcal{E}_{-\lambda}(u_g, \varphi) = \int_{\Omega^c} \mathcal{N}u_g(y)\varphi(y)\d y = \int_{\Omega^c} \nu_K^{-1}(y)\mathcal{N}u_g(y)\varphi(y)\nu_K(y)\d y\, . 
	\end{align}
	From the second  equality we can identify $\mathscr{D}_\lambda g= \mathcal{N}u_g\in L^2(\Omega^c, \nu^{-1}_K) \subset (\TnuOm)'$. Hence $\mathscr{D}_\lambda g\mapsto \mathcal{N}u_g$ which, regarding the definition of the operator $\mathcal{N}$, agrees with conceptual idea behind the Dirichlet-to-Neumann map. For the sake of consistency, in the next result, we will rather consider the following equivalent alternative identification (up to a multiplicative weight) of the Dirichlet-to-Neumann operator $\mathscr{D}_\lambda : \TnuOm \to L^2(\Omega^c, \nu_K) $ with $\mathscr{D}_\lambda g=\nu_K^{-1} \mathcal{N}u_g $. 
\end{remark}

\medskip 

\begin{theorem}
	Let  the assumptions of Theorem \ref{thm:nonlocal-Robin-var} be in force. For $\beta>0$, denote $L_\beta$ the operator $L$ subject to the Robin boundary condition $\mathcal{N} u+ \beta \nu_K u=0$. Then the point spectrum $\sigma_p(L_\beta)=(\gamma_n(\beta))_n $ of $L_\beta$ is infinitely countable say $0<\gamma_1(\beta)\leq \gamma_2(\beta) \leq \cdots\leq \gamma_n(\beta)\leq \cdots$ and the corresponding eigenfunctions are elements of $\VnuOm$ and form an orthonormal basis of $L^2(\Omega)$. 
\end{theorem}

\medskip

\begin{proof}
	It suffices to proceed as the in proof of Theorem \ref{thm:existence-of-eigenvalue-Neumann}. 
\end{proof}

\medskip

\noindent Next we observe the relation between the spectrum of the operator $L$ subject to the Robin boundary condition and that of Dirichlet-to-Neumann operator.

\begin{theorem} Let $\lambda\in \mathbb{R}\setminus\{\lambda_n:n\geq1\}$ and $\beta \in \mathbb{R}$. Under the previous notations, consider the Dirichlet-to-Neumann map $\mathscr{D}_\lambda g = \nu_K^{-1}\mathcal{N}u_g$. Let $\sigma_p(\mathscr{D}_\lambda)$ and  $\sigma_p(L_\beta)$ respectively denote pure point spectrum of $\mathscr{D}_\lambda$ and $L_\beta$.
	Then, $-\beta\in \sigma_p(\mathscr{D}_\lambda)$ if and only if $\lambda\in \sigma_p(L_\beta)$. In addition, $\dim \ker ( L_\beta -\lambda) = \dim \ker (\mathscr{D}_\lambda + \beta) $.
\end{theorem}

\medskip

\begin{proof}
	Let $u\in \ker ( L_\beta -\lambda)$ then for all $v\in \VnuOm$, 
	\begin{align*}
	Q_\beta(u,v) = \lambda\int_{\Omega} u(x)v(x)\d x\quad
	\text{equivalently} \quad
	\mathcal{E}_{-\lambda}(u,v) = -\beta\int_{\Omega^c}u(y)v(y)\nu_K(y)\d y.
	\end{align*}
	
	\noindent Set $g= \operatorname{Tr}(u)= u|_{\Omega^c}$, with the aid of \eqref{eq:DN-weak-formula} the above relation reduces to 
	\begin{align*}
	\int_{\Omega^c}\nu_K^{-1}(y)\mathcal{N} u_g(y)v(y)\nu_K(y)\d y= -\beta\int_{\Omega^c}g(y)v(y)\nu_K(y)\d y.
	\end{align*}
	Thus $ g\in \ker ( \mathscr{D}_\lambda + \beta)$. We have shown that the mapping $T: \ker ( L_\beta -\lambda)\to \ker ( \mathscr{D}_\lambda + \beta)$ with $u\mapsto
	\operatorname{Tr}(u)$ is well defined and is onto. Both assertions will follow once we show that $T$ defines a bijection. 
	In other words we only have to show that $T$ is one-to-one. For $u\in \ker ( L_\beta -\lambda)$ if $\operatorname{Tr}(u)=0$ then from the first relation above, we have $\mathcal{E}(u,v) = \lambda\int_{\Omega} u(x)v(x)\d x$ for all $v\in V_\nu^\Omega(\Omega|\mathbb{R}^d)$. Necessarily, $u=0$ otherwise $\lambda$ is a Dirichlet eigenvalue which is not the case by assumption.
\end{proof}

\section{Essentially self-adjointness for nonlocal operators}
We now investigate the unboundedness aspect of the operator $L$ subject with the Dirichlet and the  Neumann complement condition. Beforehand, let us recall some basics on unbounded operators. We refer to \cite{Dav96,Kow09,ReBa80} for a more extensive expositions on unbounded operators. Let  $(H, (\cdot, \cdot)_H)$  be a Hilbert space.  Throughout this section, we write $(T, \mathcal{D}(T))$ to denote a densely defined linear operator, i.e.  $T:\mathcal{D}(T)\to H$ is a linear operator whose domain $\mathcal{D}(T)$ is a dense subspace of  $H$. The adjoint of the operator $(T, \mathcal{D}(T))$ will be denoted by $(T^*, \mathcal{D}(T^*))$ where 
\begin{align*}
\mathcal{D}(T^*)&= \big\{w\in H:~ \text{ such that } v\mapsto (Tv, w)_H~ \text{is continunous on $\mathcal{D}(T)$}\big\}\\
& =\big\{w\in H: \exists! \, w^*\in H \text{ such that } (Tv, w)_H= (v,w^*)_H \text{ for all $v\in \mathcal{D}(T)$}\big\}
\end{align*}
and $T^*w= w^*$ for all $w\in \mathcal{D}(T)$. Note that $T^*w$ is well defined since the existence and uniqueness of $w^*$ is due to the Riesz representation theorem. 

\noindent 
Let  $(T_1, \mathcal{D}(T_1))$ and $(T_2, \mathcal{D}(T_2))$ be two densely defined operators on $H$ such that $\mathcal{D}(T_1)\subset\mathcal{D}(T_2)$ and $T_2\mid_{ \mathcal{D}(T_1)}$ then we say that $(T_1, \mathcal{D}(T_1))$  is the restriction  on $ \mathcal{D}(T_1)$ of  $(T_2, \mathcal{D}(T_2))$  and $(T_2, \mathcal{D}(T_2))$ is an extension of  $(T_1, \mathcal{D}(T_1))$ and write $T_1\subset T_2$.  In particular, $T_1=T_2$ if and only if $T_1\subset T_2$ and $T_2\subset T_1$. 

\noindent A densely defined operator $(T, \mathcal{D}(T))$ is said to be self-adjoint if $T^*=T$. Obviously, we have $\mathcal{D}(T)\subset \mathcal{D}(T^*)$ if the operator $(T, \mathcal{D}(T))$ is symmetric i.e  $(Tu, v)_H= (u, Tv )_H$ for all $u,v\in \mathcal{D}(T)$. 

\noindent The operator $(T, \mathcal{D}(T))$ is called to be closed if 
$(\mathcal{D} (T), \|\cdot\|_T )$ is Hilbert space under the graph norm defined be $\|v\|^2_T = \|v\|^2_H+ \|Tv\|^2_H $ for $v\in \mathcal{D}(T)$. In other words,  $(T, \mathcal{D}(T))$ is closed if and only its graph 
$\Gamma(T)$ is a closed subspace of $H\times H$.  Recall that $\Gamma(T):= \{(v, Tv)\in H\times H: ~ v\in  \mathcal{D}(T)\}$. 

\noindent Note that the adjoint $(T^*, \mathcal{D} (T^*))$  of $(T, \mathcal{D} (T))$ is always closed.

\noindent Instead, $(T, \mathcal{D} (T))$ is called closable if it possesses a closed extension. Equivalently, $(T, \mathcal{D} (T))$ is closable if for every sequence $(v_n)_n\subset\mathcal{D} (T))$ such that $v_n\to0$ and $T v_n\to y$ for $y\in H$ we have that $y=0$.  

\medskip

\noindent We write  $(\overline{T},\mathcal{D} (\overline{T}))$ to denote the closure of  a closable operator $(T, \mathcal{D} (T))$, i.e. the smallest closed extension of $(T, \mathcal{D} (T))$. Note that $(\overline{T},\mathcal{D} (\overline{T}))$ is the closure of $(T, \mathcal{D} (T))$ is equivalent to say that $\overline{ \Gamma(T)}= \Gamma(\overline{T})$. 

\begin{definition}
A symmetric operator $(T, \mathcal{D} (T))$ is said to be essentially self-adjoint if it admits a unique self-adjoint extension. This is merely equivalent to saying that the closure  $(\overline{T},\mathcal{D} (\overline{T}))$ of $(T, \mathcal{D} (T))$ is self-adjoint.
\end{definition}

\noindent The following result is adapted from \cite[Lemma 5.10]{Kow09}. One should note that the author intentionally ignores the identification of the space $H$ with $\ell^2(\mathbb{N})$ although it might be natural to do so.
\begin{lemma}[Lemma 5.10, \cite{Kow09}]\label{lem:essentailly-self-adjoint}
	Let H be a separable Hilbert space and $(T,\mathcal{D}(T))$ a positive symmetric unbounded operator on H. Assume  there exists an orthonormal basis $(e_j)_j$ of H such that $e_j\in \mathcal{D}(T) $ for all  $j\geq 1$ and which are eigenfunctions of $T$, i.e. $ Te_j = \lambda_j e_j$  for some $\lambda_j > 0$ for all $j \geq  1$. Then $T$ is essentially self-adjoint and its closure is unitarily equivalent with the multiplication operator $(M,\mathcal{D}(M))$  on $\ell^2(\mathbb{N})$, i.e. $M:\mathcal{D}(M)\to \ell^2(\mathbb{N})$ given by
	\begin{align*}
\mathcal{D}(M)= \big\{ (u_j)_j\in \ell^2(\mathbb{N}) ~: \sum_{j=1}^{\infty} \lambda^2_j|u_j|^2  <\infty\big\}, \quad \text{and}\quad  M((u_j)_j)= \big(\lambda_j u_j\big)_j. 
	\end{align*}
To be more precise we have $\overline{T}= UMU^{-1} $ where $U: \ell^2(\mathbb{N})\to H$ is the unitary operator defined for $(u_j)_j\in \ell^2(\mathbb{N})$ by $U ((u_j)_j) = \sum\limits_{j=1}^\infty u_je_j$ and $U^{-1}$ 
is defined for $u\in H$ by $U^{-1}(u) =((u,e_j)_H)_j$. 
\end{lemma}

\noindent One can formally write $Mu(x)= \lambda(x) u(x)$  for $x\in \mathbb{N}$ where $u,\lambda: \mathbb{N}\to \mathbb{R}$ with $\lambda(j) = \lambda_j$ and $u \in  H$ with $u(j)= u_j$. So that if $m$ denotes the counting measure on $\mathbb{N},$ one has
\begin{align*}
\sum_{j=1}^{\infty} \lambda^2_j|u_j|^2= \int_{\mathbb{N}} \lambda^2(x)u^2(x)m(\d x).
\end{align*}
Note that if the sequence $(\lambda_j)$ has an accumulation point, the spectrum is not the same as the set of eigenvalues (this already occurred for compact operators, where $0$ may belong to the spectrum even if the kernel is trivial).

\begin{proof} First of all for $u\in H$ consider the identification $u_j = (u,e_j)_H$. From the  Bessel and the  Parseval identities we know that 
\begin{align*}
	u= \sum_{ j=1}^\infty u_j e_j\quad\text{and}\quad	\|u\|_H= \sum_{ j=1}^\infty |u_j|^2. 
\end{align*}
\noindent This shows that the operators $U$ and $U^{-1}$ are well defined. Moreover, $U$ is a unitary operator. Indeed for $(u_j)_j \in \ell^2(\mathbb{N})$ and $v\in H$, set $u = U((u_j)_j)\in H$ with $u_j = (u,e_j)_H$ since $U^{-1}(v)= ((v,e_j)_H)_j$ we get 
\begin{align*}
( U((u_j)_j), v)_H= (u, v)_H= \sum_{ j=1}^\infty (u,e_j)_H (v,e_j)_H=\sum_{ j=1}^\infty u_j (v,e_j)_H= \big((u_j)_j, U^{-1}(v)\big)_{\ell^2(\mathbb{N})}. 
\end{align*} 
We will show  $\overline{T}= UMU^{-1}$is the closure of $T$ and is self adjoint. We proceed in several steps. 

\vspace{2mm}
 \noindent \textbf{Step 1:} We show that $M$ is self-adjoint, this would imply that  $\overline{T}= UMU^{-1}$ is self-adjoint too since $U^*= U^{-1}$. It is  sufficient to show that $\mathcal{D}(M^*)\subset \mathcal{D}(M)$. In fact, since $M$ is symmetric, we get 
$\mathcal{D}(M)\subset \mathcal{D}(M^*) $. Let $w\in \mathcal{D}(M^*)$ then by definition, the linear form $v\mapsto (Mv,w)_{\ell^2(\mathbb{N}}$ is bounded on $\mathcal{D}(M)\subset \ell^2(\mathbb{N})$. In other words, there is a constant $C>0$ such that 
$$|(Mv,w)_{\ell^2(\mathbb{N})}|\leq C\|v\|_{\ell^2(\mathbb{N})}, \quad\text{ for all $v\in \mathcal{D}(M)$}.$$
 
 \noindent Specializing with $v=v_n\in \mathcal{D}(M)$ where $v_{n,j}= \lambda_j(w, e_j)e_j$  for $j\leq n$ and $v_{n,j}=0$ for $j\geq n+1$ for fixed  $n\geq 1,$ we have  $v_n = ((v_{n,j}))_j \in \mathcal{D}(M)$ and $\|v_n\|^2_H= \sum_{ j=1}^n\lambda_j^2|(w, e_j)|^2= (Mv_n,w)_{\ell^2(\mathbb{N})}.$ This amounts the above inequality to 
 $$\sum_{ j=1}^n\lambda_j^2|(w, e_j)|^2 \leq C\quad\text{ for all $n\geq1$}.$$
 Hence  we conclude that $Mw = (\lambda_j w_j)_j\in \ell^2(\mathbb{N})$, i.e. $w\in \mathcal{D}(M)$. Finally $\mathcal{D}(M)= \mathcal{D}(M^*)$ altogether with symmetry we get that $M$ is self-adjoint.  

\noindent \textbf{Step 2:} $\overline{T} = UMU^{-1}$ is closed since it is self-adjoint.
	
\noindent \textbf{Step 3:} Let us check that  $\overline{T}= UMU^{-1}$ is an extension of  $T$.  Consider  $u\in \mathcal{D}(T)$ and set $u_j= \big( u, e_j\big)_H$. Since each $e_j\in \mathcal{D}(T)$ is an eigenfunction,  the symmetry of $T$ gives $(Tu, e_j)_H= (u, Te_j)_H = \lambda_j (u, e_j)_H$. Thus accounting  that $(e_j)_j$ is an orthonormal basis of $H$ we find that
	\begin{align*}
	Tu=  \sum_{j=1}^{\infty} \big( Tu, e_j\big)_He_j=   \sum_{j=1}^{\infty} \lambda_j\big( u, e_j\big)_He_j= U\big((\lambda_ ju_j)_j\big) = UM(( u_j)_j)= UM U^{-1} u.
	\end{align*}
 In other words, $T\subset  UMU^{-1} $, i.e. $UMU^{-1}$ is an extension of $T$ since $Tu = UMU^{-1} u $ for all $u \in \mathcal{D}(T)$ and
	\begin{align*}
	\sum_{j=1}^{\infty} \lambda^2_j|\big(u, e_j\big)_H|^2=  \|Tu\|^2_H<\infty.
	\end{align*}
	\noindent  
	
\noindent \textbf{Step 4:} It remains to show that $\overline{\Gamma(T)}=\Gamma(\overline{T})$ which will imply that $\overline{T}$ is the closure of $T$. From Step 3 we know that $\Gamma(T)\subset \Gamma(\overline{T})$. Let $(u, \overline{T}u )\in \Gamma(\overline{T}) $ with $u\in \mathcal{D}(\overline{T})$. For all $n\geq 1$ $u_n = \sum_{j=1}^{n} (u,e_j)_H e_j\in \mathcal{D}(T)$ since $e_j\in \mathcal{D}(T)$ and $Tu_n = \sum_{j=1}^{n} \lambda_j (u,e_j)_H e_j$ since $T$ is symmetric and each  $e_j$ 
is an eigenvalue of $T$. Thus, $(u_n, Tu_n)\in \Gamma(T)$. Meanwhile, we know that

\begin{align*}
&\|u\|^2_H=\sum_{j=1}^{\infty} |(u,e_j)_H|^2<\infty \quad\text{and}\quad \|\overline{T} u\|^2_H=\sum_{j=1}^{\infty} \lambda_j^2|(u,e_j)_H|^2<\infty \\
\intertext{since}
&u= \sum_{j=1}^{\infty} (u,e_j)_H e_j\in H\quad\text{and}\quad \overline{T}u= UMU^{-1}u= \sum_{j=1}^{\infty} \lambda_j (u,e_j)_H e_j\in H.
	\end{align*} 
	
\noindent Therefore, it follows that $\|u_n-u\|_H\xrightarrow[]{n\to \infty}0$ and $\|Tu_n- Tu\|_H\xrightarrow[]{n\to \infty}0$. Hence, $\overline{\Gamma(T)}=\Gamma(\overline{T})$. 
From the foregoing,  it follows that $\overline{T}= UMU^{-1}$ is self-adjoint and is the closure of $T$ which means that $T$ is essentially self-adjoint. 

%
\end{proof}

\noindent Armed with these prerequisites, let us turn our attention to some concrete examples to which the above notions apply. 
Assume $\Omega\subset \R^d$ is an open bounded set. 
%
%
On $L^2(\Omega)$, consider the operator $(L_D, V_D)$ where  $V_D= \VnuOmO$ and $L_D$ is the integrodifferential operator $L$ subject to the Dirichlet complement condition.  In other words, $L_Du= f$ if and only $f\in L^2(\Omega)$ and 
\begin{align*}
\mathcal{E}(u,v)= \int_\Omega f(x) v(x)\d x\quad\text{for all $v\in \VnuOmO$}. 
\end{align*}

\noindent Likewise, on $L^2(\Omega)$, consider the operator $(L_N, V_N)$  where  $V_N=\big\{ v\in \VnuOm: ~\mathcal{N} v=0~ \text{on $\Omega^c$}\big\} $ and  $L_N$ is the integrodifferential operator $L$ subject to the Neumann complement condition.  In other words, $L_Nu= f$ if and only $f\in L^2(\Omega)$ and 
\begin{align}\label{eq:essential-neumman-L}
\mathcal{E}(u,v)= \int_\Omega f(x) v(x)\d x\quad\text{for all $v\in \VnuOm$}. 
\end{align}
\noindent We adopt the convention that $\mathcal{N} u=0$ on $\Omega^c$ for function $u\in\VnuOm$ if there is $f\in L^2(\Omega)$ such that \eqref{eq:essential-neumman-L} holds.  The next theorem is the core result of this section showing that under appropriated conditions on $\Omega$ and the jump kernel $\nu$ the operators $L_D$ and $L_N$ are essentially self-adjoint. 
\begin{theorem}\label{thm:essentially-self-adjoint}
The operators $(L_D, V_D)$ and $(L_N, V_N)$ are symmetric and positive on $L^2(\Omega)$. Furthermore, the following assertions are true. 
\begin{enumerate}[$(i)$]
	\item Assume  $(\nu, \Omega)$ belongs to one of the classes $\mathscr{A}_i,~ i=1,2,3$ (see page  \pageref{eq:class-lipschitz}), then $(L_N, V_N)$ is a densely defined, unbounded and essentially self-adjoint operator on $L^2(\Omega)$.
\item Assume $\nu$ satisfies the conditions \eqref{eq:non-integrability-condition} and \eqref{eq:integrability-condition-near-zero} (see page  \pageref{eq:non-integrability-condition}), then $(L_D, V_D)$ is a densely defined, unbounded and essentially self-adjoint operator on $L^2(\Omega)$. 

\end{enumerate}
\end{theorem}
\begin{proof} Proceeding as in the proof of Theorem \ref{thm:existence-of-eigenvalue-Neumann} one finds that $(L_D, V_D)$ and $(L_N, V_N)$ are symmetric and positive on $L^2(\Omega)$. Assume $(\nu, \Omega)$ belongs to one of the classes $\mathscr{A}_i,~ i=1,2,3$ then by Theorem \ref{thm:rayleigh-construction} there exists a sequence of eigenpairs $(\mu_n, \varphi_n)_{n\in\mathbb{N}_0}$ such that for each $n\geq 0$,  $\phi_n\in V_N$, $L_N\phi_n= \mu_n \phi_n$ and the family $(\varphi_n)_{n\in \mathbb{N}_{0}}$ is an orthonormal basis of $L^2(\Omega)$.  It turns out that $V_N$ is dense in $L^2(\Omega)$ and thus $(L_N, V_N)$ is essentially self-adjoint on $L^2(\Omega)$ according to Lemma \ref{lem:essentailly-self-adjoint}. Furthermore, due to the orthonormal basis $(\phi_n)_{\mathbb{N}_0}$, it becomes clear that $(L_N, V_N)$ is unbounded on $L^2(\Omega)$. Likewise, from Theorem \ref{thm:rayleigh-construction-Dirichlet} and Lemma \ref{lem:essentailly-self-adjoint} it follows that 
$(L_D, V_D)$ is densely defined and  essentially self-adjoint on $L^2(\Omega)$.  
\end{proof}

\begin{remark} Under the  assumption that  $(\nu, \Omega)$ belongs to one of the classes $\mathscr{A}_i,~ i=1,2,3$, the corresponding formulation of the above theorem  is true if the operator $L$ is subject to the mixed or Robin complement condition. 
\end{remark}

\begin{corollary} Let $\Omega\subset \R^d$ be open and bounded. Let $\nu(h) = C_{d, \alpha}|h|^{-d-\alpha}$
$h\neq 0$ for some $\alpha\in (0, 2)$. Let $V_D^{\alpha/2}(\Omega|\R^d)= \big\{v\in V^{\alpha/2}(\Omega|\R^d): \, v=0\, \text{on $\Omega^c$}\big\}$.
Then $((-\Delta)^{\alpha/2}, V_D^{\alpha/2})$ is essentially self-adjoint and also $((-\Delta)^{\alpha/2}, V_N^{\alpha/2})$ if $\Omega$ is Lipschitz  where $V_N^{\alpha/2}(\Omega|\R^d)= \big\{v\in V^{\alpha/2}(\Omega|\R^d): \, \mathcal{N}v=0\, \text{on $\Omega^c$}\big\}$.
\end{corollary}
\noindent  It is relevant to mention that the  characterization  of essentially self-adjointness of the fractional Laplacian $(-\Delta)^{\alpha/2}$ subject to the Dirichlet complement condition on an unbounded open set $\Omega\subset \R^d $  is the central discussion in \cite{HKM17}  and several other references therein. Next, in spirit of Theorem \ref{thm:essentially-self-adjoint} we have the following result in the local setting %
for the Laplace operator. 
\begin{theorem} Assume $\Omega\subset\R^d$ is open and bounded. 
	Let $(\Delta_D, H_D)$ be the Dirichlet Laplacian $\Delta_D$ with $H_D= H_0^1(\Omega)$.  Let $(\Delta_N, H_N)$ be the Neumann Laplacian with  $H_N=\big\{ v\in H^1(\Omega): ~\frac{\partial u}{\partial n}=0~ \text{on $\partial \Omega$}\big\} $ where by convention we say that $\frac{\partial u}{\partial n} =0$ on $\partial \Omega$  if there exists $f\in L^2(\Omega)$ such that 
	\begin{align*}
	\int_\Omega\nabla u(x)\cdot \nabla v(x)\d x= \int_\Omega f(x) v(x)\d x\quad\text{for all $v\in H^1(\Omega)$}. 
	\end{align*}
Then $(\Delta_D, H_D)$  is densely defined, unbounded and essentially self-adjoint. If $\Omega$ is Lipschitz, the same holds for $(\Delta_N, H_N)$. 
\end{theorem}

\chapter{From Nonlocal To Local }\label{chap:Nonlocal-To- Local}

This Chapter is essentially  devoted to investigating  local objects as limits of nonlocal ones. We especially focus on the following concepts: Sobolev spaces, elliptic operators, energy forms, Poincar\'e type inequalities, elliptic   PDEs  of second order with Dirichlet or Neumann boundary condition,  Dirichlet and Neumann eigenvalues.  Roughly speaking, this chapter serves as a needed bridge between basic elliptic partial differential equations of second order and elliptic integrodifferential equations that require a significant background in functional analysis. To begin this journey let us first introduce our central tool. 


\section{Approximation  of the Dirac mass by L\'evy measures}\label{sec:approx-dirac-mass}

\noindent In this section we introduce the tool that will help us to move from nonlocal objects to their corresponding local notions. We begin by coining the concept of $p$-L\'evy measure for $1\leq p<\infty$.
\begin{definition}
A nonnegative Borel measure $\nu(\,\d h)$ on $\mathbb{R}^d$ will be called  a $p$-L\'evy measure for $1\leq p<\infty$ if $\nu(\{0\})=0$ and it satisfies the $p$-L\'evy integrability condition. That is to say
\begin{align*}
\int_{\mathbb{R}^d} (1\land |h|^p)\nu(\,\d h)<\infty. 
\end{align*}
Patently, one recovers the usual definition a L\'evy measure when $p=2$. Notationally, we shall intentionally ignore the dependence of $\nu$ on $p$. In addition, for a family $(\nu_\varepsilon)_{0<\varepsilon <\varepsilon_0}$ for a fixed real number $\varepsilon_0>0$ we shall merely write $(\nu_\varepsilon)_\varepsilon$ or $(\nu_\varepsilon)_{\varepsilon>0}$. It is important to keep in mind that $\varepsilon $ is a quantity near $0$ ($\varepsilon\to 0^+$).
\end{definition}

\medskip

\noindent We start with the following motivating rescaling result. 
\begin{proposition}\label{prop:rescalled}
Assume $\nu : \mathbb{R}^d \to [0, \infty]$ is a positive measurable function satisfying the $p$-L\'evy integrability condition, i.e. $\nu \in L^1((1\land |h|^p)dh$. Define the rescaled family $(\nu_\varepsilon)_\varepsilon$ as follows

\begin{align}\label{eq:rescalled-levy-measure}
\begin{split}
\nu_\varepsilon(h) = 
\begin{cases}
\varepsilon^{-d-p}\nu\big(h/\varepsilon\big)& \text{if}~~|h|\leq \varepsilon\\
\varepsilon^{-d}|h|^{-p}\nu\big(h/\varepsilon\big)& \text{if}~~\varepsilon<|h|\leq 1\\
\varepsilon^{-d}\nu\big(h/\varepsilon\big)& \text{if}~~|h|>1.
\end{cases}
\end{split}
\end{align}

\noindent Then 
\begin{align*}
&\int_{\mathbb{R}^d}(1\land |h|^p)\nu_\varepsilon(h)d h= \int_{\mathbb{R}^d}(1\land |h|^p)\nu(h)d h\,
\quad\text{for each $\delta>0$ and}\quad
\lim_{\varepsilon\to 0 } \il_{|h|\geq \delta }(1\land |h|^p)\nu_\varepsilon(h)d h=0.
\end{align*}

\end{proposition}

\medskip

\begin{proof}
Observe that since $\nu \in L^1((1\land |h|^p)dh)$ by dominated convergence theorem we get
\begin{align*}
\lim_{\varepsilon\to 0 } \il_{|h|\geq \delta }(1\land |h|^p)\nu_\varepsilon(h)d h=\lim_{\varepsilon\to 0 } \il_{|h|\geq \delta/\varepsilon }(1\land |h|^p)\nu(h)d h=0.
\end{align*}
We omit the remaining details as it solely involves straightforward computations.
\end{proof}

\medskip

\noindent There are two keys observations that govern the rescaled family $(\nu_\varepsilon)_\varepsilon$ for $p=2$. The first is that it gives raise to a family of measures with a concentration property. Secondly, from the probabilistic point of view one obtains a family of pure jumps L\'evy processes $(X_\varepsilon)_\varepsilon$ each associated with the measure $\nu_\varepsilon$ from a L\'evy process $X$ associated with $\nu$. Viewed from this scope we aim to show that $(X_\varepsilon)_\varepsilon$ converges in some sense (to be clarified later see Theorem \ref{thm:mosco-characterisation} and Theorem \ref{thm:Mosco-convergence}) to a Brownian motion provided that one in addition assumes that $\nu$ is radial. This, in certain sense, could be one more argument to justify the ubiquity of the Brownian motion.

\medskip

\noindent From now on, for $1\leq p<\infty$, $(\nu_\varepsilon)_{\varepsilon>0}$ denotes a family of nonnegative functions such that 
\begin{align}
& \nu_\varepsilon~~\text{is radial}\qquad \text{and} \qquad\int_{\mathbb{R}^d} (1\land|h|^p) \nu_\varepsilon(h)\,\d h=1\label{eq:normalized-inetrals},\\
& \lim_{\varepsilon\to 0}\il_{|h|\geq \delta} (1\land|h|^p) \nu_\varepsilon(h)\,\d h=0 \quad\text{for all }\quad \delta >0 .\label{eq:concentration-property}
\end{align}
Note that the relation \eqref{eq:concentration-property} is often known as the concentration property and because of \eqref{eq:normalized-inetrals} it is merely equivalent to 
\begin{align*}
\lim_{\varepsilon\to 0}\il_{|h|\geq \delta }\nu_\varepsilon(h)\,\d h=0\quad\text{for all }\quad \delta >0.
\end{align*}
 Consequently, we also have 
 \begin{align*}
 \lim_{\varepsilon\to 0}\il_{|h|\leq \delta }(1\land |h|^p)\nu_\varepsilon(h)\,\d h= \lim_{\varepsilon\to 0}\il_{|h|\leq \delta }|h|^p\nu_\varepsilon(h)\,\d h=1\quad\text{for all }\quad \delta >0.
 \end{align*}
 \begin{definition}\label{def:plevy-approx}
Under the conditions \eqref{eq:normalized-inetrals} and \eqref{eq:concentration-property} the family $(\nu_\varepsilon)_\varepsilon$ will be called \emph{a radial Dirac approximation of $p$-L\'{e}vy measures.}
 \end{definition}
 \bigskip
 \begin{remark}\label{rem:asymp-nu}
 	Under the conditions \eqref{eq:normalized-inetrals} and \eqref{eq:concentration-property} we have that for all $ \beta, R>0$ 
 	
 	\begin{align*}
 	\lim_{\varepsilon\to 0}\il_{|h|\leq R} (1\land |h|^\beta)\nu_\varepsilon(h)\mathrm{d}h =
 	\begin{cases}
 	0\quad \text{if}\quad \beta >p\\
 	1\quad \text{if}\quad \beta =p.
 	\end{cases}. 
 	\end{align*}
Indeed, for fixed $\delta>0$
 	\begin{align*}
 	&\lim_{\varepsilon \to 0} \il_{ \delta<|h|\leq R} (1\land |h|^p) \nu_\varepsilon(h) \mathrm{d}h
 	\leq \lim_{\varepsilon \to 0} \il_{ |h|>\delta} (1\land |h|^p) \nu_\varepsilon(h) \mathrm{d}h=0,
 	\\
 	&\lim_{\varepsilon \to 0}\il_{ |h|<\delta} (1\land |h|^p) \nu_\varepsilon(h) \mathrm{d}h
 	=1- \lim_{\varepsilon \to 0} \il_{ |h|\geq \delta} (1\land |h|^p) \nu_\varepsilon(h) \mathrm{d}h=1. 
 	\end{align*}
 	Thus, for $\beta>p$, we have 
 	\begin{align*}
 	\lim_{\varepsilon \to 0}\il_{|h|\leq R} (1\land |h|^\beta)\nu_\varepsilon(h)\mathrm{d}h
 	&\leq  \lim_{\varepsilon \to 0} \Big( R^{\beta-p} \il_{\delta< |h|\leq R} (1\land |h|^p)\nu_\varepsilon(h) \mathrm{d}h + \delta^{\beta-p} \il_{|h|<\delta} (1\land |h|^p)\nu_\varepsilon(h) \mathrm{d}h\Big) \\
 	& = \delta^{\beta-p}.
 	\end{align*}
 	Letting $\delta \to 0$ provides the claim. 
 	
 \end{remark}

\bigskip

\noindent The following result  infers the weak convergence of the family $(\nu_\varepsilon)_\varepsilon$ to the Dirac mass at origin and thus withstands Definition \ref{def:plevy-approx}.
\begin{proposition}\label{prop:dirac-approx} Assume $(\nu_\varepsilon)_\varepsilon$ satisfies conditions \eqref{eq:normalized-inetrals} and \eqref{eq:concentration-property} for $p=1$. Then $\nu_\varepsilon\rightharpoonup \delta_0$(weakly) in the sense of distributions, i.e. for every $\varphi\in C_c^\infty(\R^d)$, $\langle\nu_\varepsilon, \varphi\rangle\xrightarrow{\varepsilon\to 0} \langle\delta_0, \varphi\rangle = \varphi(0).$ Here $\delta_0$ stands for the Dirac mass at the origin.
\end{proposition}

\medskip

\begin{proof}
Take $\varphi\in C_c^\infty(\R^d)$ and in virtue of \eqref{eq:normalized-inetrals} write 
\begin{align*}
\langle\nu_\varepsilon, \varphi\rangle &=\int_{\R^d}\varphi(h)\nu_\varepsilon(h)\d h \\
&= \varphi(0) + \int_{|h|< 1}(\varphi(h)-\varphi(0)\, \nu_\varepsilon(h) \d h + \int_{|h|\geq 1}(\varphi(h)- \varphi(0))\nu_\varepsilon(h)\d h\\
&= \varphi(0) + \int_{|h|< 1}(\varphi(h)-\varphi(0)-\nabla \varphi(0)\cdot h)\, \nu_\varepsilon(h) \d h + \int_{|h|\geq 1}(\varphi(h)- \varphi(0))\nu_\varepsilon(h)\d h\\
&= \varphi(0) + \int_{|h|< 1}\int_0^1(\nabla \varphi(th)-\nabla \varphi(0)\cdot h)\, \nu_\varepsilon(h) \d t\d h + \int_{|h|\geq 1}(\varphi(h)- \varphi(0))\nu_\varepsilon(h)\d h \\
&=\varphi(0) + \int_{|h|< 1}\int_0^1\int_0^1 s((D^2 \varphi(tsh)\cdot h)\cdot h)\, \nu_\varepsilon(h) \d s\d t\d h + \int_{|h|\geq 1}(\varphi(h)- \varphi(0))\nu_\varepsilon(h)\d h .
\end{align*}
The conclusion clearly follows since 
\begin{align*}
\Big| \int_{|h|\geq 1}(\varphi(h)- \varphi(0))\nu_\varepsilon(h)\d h\Big|\leq 2\|\varphi\|_\infty \int_{|h|\geq 1}\nu_\varepsilon(h)\d h\xrightarrow{\varepsilon\to 0}0
\end{align*}
and by Remark \ref{rem:asymp-nu} we have 
\begin{align*}
\Big|\int_{|h|< 1}\int_0^1\int_0^1 s((D^2 \varphi(tsh)\cdot h)\cdot h)\, \nu_\varepsilon(h) \d s\d t\d h \nu_\varepsilon(h)\d h\Big|\leq \|D^2\varphi\|_\infty \int_{|h|< 1}|h|^2\nu_\varepsilon(h)\d h\xrightarrow{\varepsilon\to 0}0.
\end{align*}
\end{proof}

\bigskip

\noindent Another motivating reason to consider the family $(\nu_\varepsilon)_\varepsilon$ is given by the following.
\begin{proposition}
	For $p=2$ and $ \alpha= 2-\varepsilon\in (0,2)$ under the conditions \eqref{eq:normalized-inetrals} and \eqref{eq:concentration-property} suppose the function $u: \mathbb{R}^d\to \mathbb{R}$ is bounded and $C^2$ on a neighborhood of $x$ then 
	$$\lim_{\alpha \to 2} L_\alpha u(x) =-\frac{1}{2d}\Delta u(x)$$
	where $L_\alpha$ is the nonlocal operator
	\begin{align*}
	L_\alpha u(x) := -\frac{1}{2}\int_{\mathbb{R}^d }(u(x+h)+ u(x-h) -2u(x)) \nu_{2-\alpha}(h)\,\d h.
	\end{align*}
\end{proposition}

\bigskip

\begin{proof}
	It suffices to adapt the proof of Proposition \ref{prop:elliptic-matrix-chap1} $(iv)$.
\end{proof}

\bigskip

 \vspace{0.5cm}
\noindent Let us mention some examples of particular interest. The first class relates to the fractional Laplacian.
\begin{example}\label{Ex: stable-class-2}
For $p=2$ consider $(\nu_\alpha)_\alpha$ the family of $\alpha$-stable kernel defined for $\alpha= 2-\varepsilon\in (0, 2)$ and $h \neq 0$ by $\nu_\alpha(h) = a_{d,\alpha}  |h|^{-d-\alpha}$ with $ a_{d,\alpha}  = \tfrac{\alpha(2-\alpha)}{2 |\mathbb{S}^{d-1}|} .$
Let us show that \eqref{eq:normalized-inetrals} and \eqref{eq:concentration-property} are fulfilled for $p=2$. Passing through polar coordinates yields 
	\begin{align*}
		& \int_{\mathbb{R}^d} (1\land |h|^2) |h|^{-d-\alpha}\,\d h= 	 |\mathbb{S}^{d-1}| \Big(\int_{0}^{1} r^{1-\alpha}\,\d r+ \int_{1}^{\infty}r^{-1-\alpha}\,\d r\Big) \\
		&= |\mathbb{S}^{d-1}| \Big(\frac{1} {2-\alpha}+ \frac{1} {\alpha}\Big)= |\mathbb{S}^{d-1}| \frac{2} {\alpha(2-\alpha)}= a^{-1}_{\alpha,d}. 
		\end{align*}
For $\delta>0$ a similar computation gives
	\begin{align*}
&a_{d,\alpha}
\il_{|h|\geq \delta} (1\land |h|^2) |h|^{-d-\alpha}\,\d h\leq	\frac{\alpha(2-\alpha)}{2} \int_{\delta}^{\infty}r^{-1-\alpha}\,\d r
= \frac{(2-\alpha)\delta^{-\alpha}} {2}\xrightarrow{\alpha \to 2}0.
\end{align*}
Note that $ \lim\limits_{\alpha \to 2} \frac{C_{d,\alpha} }{ a_{d,\alpha} }=1$ (see Proposition \ref{prop:asymp-cds}) where $C_{d,\alpha}$ is the norming constant of the fractional Laplacian. 
 \end{example}

\noindent More generally, in connection with the fractional Sobolev spaces we have the following example. 
\begin{example}\label{Ex: stable-class}
The family $(\nu_\varepsilon)_\varepsilon$ of kernels defined for $h \neq 0$ by
 $$\nu_\varepsilon(h) = a_{\varepsilon, d,p} |h|^{-d-p+\varepsilon}\quad\text{with }\quad a_{\varepsilon, d,p} = \frac{\varepsilon(p-\varepsilon)}{p|\mathbb{S}^{d-1}|}$$ fulfills the conditions \eqref{eq:normalized-inetrals} and \eqref{eq:concentration-property}. 
	
\end{example}

\noindent The next class is that of Proposition \ref{prop:rescalled}.
\begin{example}
Assume $\nu: \mathbb{R}^d \to [0, \infty]$ is radial and consider the family $(\nu_\varepsilon)_\varepsilon$ such that each $\nu_\varepsilon$ is the rescaling of $\nu$ defined as in \eqref{eq:rescalled-levy-measure} provided that 

	\begin{align*}
	\int_{\mathbb{R}^d}(1\land |h|^p)\nu(h)\,\d h=1. 
\end{align*}
A subclass is obtained if one considers an integrable radial function $\rho:\mathbb{R}^d \to [0,\infty]$ and defines $\nu(h) = c|h|^{-p} \rho(h)$ for a suitable normalizing constant $c>0$. 
\end{example}

\medskip

\noindent The following class  is of important interest and  is related to the so called approximation of the identity. 
\begin{example}\label{Ex:rho-var}
	Assume $(\rho_\varepsilon)_\varepsilon$ is an approximation of the identity, i.e. $\varepsilon>0$, $\rho_\varepsilon: \mathbb{R}^d \to [0, \infty]$ is radial,
\begin{align*}
	\int_{\mathbb{R}^d}\rho_\varepsilon (h)\,\d h=1\quad\text{and}\quad \lim_{\varepsilon\to 0 }\il_{|h|\geq \delta} \rho_\varepsilon(h) \,\d h=0\quad\text{for all }\quad \delta >0.
\end{align*}
For instance,  define $\rho_\varepsilon (h) =\varepsilon^{-d} \rho(h/\varepsilon)$ where $\rho:\mathbb{R}^d \to [0,\infty]$ is radial and such that $\int_{\R^d}\rho(h)\d h=1$. 
Define the family $(\nu_\varepsilon)_\varepsilon$ by 
$$\nu_\varepsilon(h) = c_\eps |h|^{-p} \rho_\varepsilon(h)$$
where $c_\eps>0$ is a normalizing constant for which \eqref{eq:normalized-inetrals} holds true. 
\medskip

\noindent The resulting class obtained here is somewhat restrictive and does not capture good classes 
such that as the kernels from Example \ref{Ex: stable-class}. To be more precise, there is no sequence $(\rho_\eps)_\eps$ satisfying the conditions above, for which $\nu_\eps(h) = c_\eps |h|^{-p} \rho_\varepsilon(h)= |h|^{-d-p+\varepsilon}$ for all $h$. 

\vspace{1mm}
\noindent To remedy to this, the truncation $\widetilde{\rho}_\varepsilon(h) = c_\varepsilon \rho_\varepsilon(h)\mathds{1}_{B_1(0)}$ is 
actually sufficient where $c_\varepsilon>0$ is such that $\int \widetilde{\rho}_\varepsilon(h) \d h=1$. Afterwards define $\widetilde{\nu}_\varepsilon(h) = \widetilde{\rho}_\varepsilon(h)|h|^{-p} $. 

\noindent Note that the class $(\widetilde{\nu}_\varepsilon )_\varepsilon$ (where each $\widetilde{\rho}_\varepsilon(h)  $ is supported in a ball centered at the origin) is precisely the one considered in \cite{BBM01,Ponce2004} and used as main tool on related topics. 
\end{example}

\bigskip

\noindent Now we collect some concrete examples of sequence $(\nu_\varepsilon)_\varepsilon$ satisfying \eqref{eq:normalized-inetrals} and \eqref{eq:concentration-property}. 
\begin{example}  

 \label{Ex:example-poincre1} Let $0<\varepsilon < 1$ and $\beta >-d$. Set 
\begin{align*}
\nu_\varepsilon(h) = \frac{d+\beta}{ |\mathbb{S}^{d-1}|\varepsilon^{d+\beta}} |h|^{\beta-p}\mathds{1}_{B_{\varepsilon}}(h). 
\end{align*} 
 For the limiting case $\beta=-d$ consider $0<\varepsilon <\varepsilon_0<1$ and put
\begin{align*}
\nu_\varepsilon(h) = \frac{1}{|\mathbb{S}^{d-1}|\log(\varepsilon_0/\varepsilon )}|h|^{-d-p}\mathds{1}_{ B_{\varepsilon_0}\setminus B_\varepsilon}(h). 
\end{align*} 
It is important to notice that $ \frac{\log(\varepsilon_0/\varepsilon )}{ |\log(\varepsilon )|}\to 1.$ Some special cases are obtained with $ \beta = 0$, $\beta =p$ and $\beta = (1-s)p-d$ for $s\in (0,1)$. 

\end{example}

\begin{example} \label{item-example-poincre2}
		
	  Let $0<\varepsilon < \varepsilon_0<1$ and $\beta >-d$. For $x \in \mathbb{R}^d$ consider
		\begin{align*}
		\nu_\varepsilon(h) = \frac{(|h|+\varepsilon)^{\beta}|h|^{-p}}{|\mathbb{S}^{d-1}| b_\varepsilon}\mathds{1}_{B_{\varepsilon_0}}(h)\qquad \hbox{with }\quad  b_\varepsilon = \varepsilon^{d+\beta}\int^{1}_{\frac{\varepsilon}{\varepsilon+\varepsilon_0}} t^{-d-\beta-1}(1-t)^{d-1}dt. 
		\end{align*} 
For the limiting case $\beta=-d$ consider
		
		\begin{align*}
		\nu_\varepsilon(h) = \frac{(|h|+\varepsilon)^{-d}|h|^{-p}}{|\mathbb{S}^{d-1}| b_\varepsilon}\mathds{1}_{B_{\varepsilon_0}(h)}\qquad \hbox{with }\quad  b_\varepsilon = \int^{1}_{\frac{\varepsilon}{\varepsilon+\varepsilon_0}} t^{-1}(1-t)^{d-1}dt. 
			\end{align*}
In either case the constant $b_\varepsilon $ is such that $\int_{\mathbb{R}^d}(1\land |h|^p) \nu_\varepsilon(h) dh =1$. Note that one can check that  $$ a_\varepsilon:= \frac{|\log(\varepsilon)|}{b_\varepsilon}\to 1\quad \hbox{as }~~~\varepsilon\to 0^+. $$ Another example familiar to the case $\beta=-d$ is 
	\begin{align*}
\nu_\varepsilon(h) = \frac{(|h|+\varepsilon)^{-d-p}}{|\mathbb{S}^{d-1}| b_\varepsilon}\mathds{1}_{B_\varepsilon(h)}\qquad \hbox{with }\quad  b_\varepsilon = \int^{1}_{\frac{\varepsilon}{\varepsilon+\varepsilon_0}} t^{-1}(1-t)^{d+p-1}dt. 
\end{align*} 
\end{example}

\section{Characterization of classical Sobolev spaces}\label{sec:charact-W1p}

In this section we show the convergence of the nonlocal Sobolev like spaces $(W^p_{\nu_\varepsilon}(\Omega))_\varepsilon$ to the usual Sobolev space $W^{1,p}(\Omega)$. 
 Further we will see that the $W^{1,p}(\Omega)$ with $1<p<\infty$ can be characterized by the mean of $(W^p_{\nu_\varepsilon}(\Omega))_\varepsilon$. The case $p=1$ will be obtained separately as it involves the space of bounded variation functions $BV(\Omega)$. Our exposition here mainly relies upon \cite{BBM01,Brezis-const-function}. Different approaches can be found in \cite{Pon06,PS17,LS11}.
\begin{lemma}\label{lem:BBM-regular}
Let $u\in C^1_c(\mathbb{R}^d)$ and let $\Omega\subset \mathbb{R}^d$ be open (not necessarily bounded). Assume that the family $(\nu_\varepsilon)_\varepsilon$ fulfills the conditions \eqref{eq:normalized-inetrals} and \eqref{eq:concentration-property} then for $1\leq p<\infty$ the following convergence occurs in both pointwise and $L^1(\Omega)$ sense 
\begin{align*}
	\lim_{\varepsilon\to 0}\int_{\Omega}|u(x)-u(y)|^p\nu_\varepsilon(x-y) \mathrm{d}y = K_{d,p} | \nabla u(x)|^p
	\end{align*}
here $	K_{d,p}$ is a universal constant\footnote{The geometrical  constant $K_{d,p}$ is the same as initially established for  in \cite{BBM01}. A similar constant also appears in \cite[Section 7]{JN10}.} independent of the geometry of $\Omega$ and is defined for any unit vector $e\in \mathbb{S}^{d-1}$ by 
	\begin{align*}
	K_{d,p} = \fint_{\mathbb{S}^{d-1}} |w\cdot e|^p\mathrm{d}\sigma_{d-1}(w) = \frac{\Gamma(d/2) \Gamma((p+1)/2)}{\Gamma((p+d)/2)\Gamma(1/2)}.
	\end{align*} 
\end{lemma}

\medskip 

\begin{proof}	Let $\sigma>0$ be sufficiently small. By assumption $\nabla u$ is uniformly continuous and hence one can find  $0<\eta<1$  such that if  $|x-y|<\eta$ then 
	\begin{align}\label{eq:continuity-estimate}
	\frac{1}{2} |\nabla u(x)|\leq 	|\nabla u(y)| \leq \frac{3}{2} |\nabla u(x)|\qquad\hbox{and}\qquad |\nabla u(y)-\nabla u(x)| \leq \sigma.
	\end{align}
Let $\eta_x = \min (\eta, \delta_x)$ with $\delta_x= \operatorname{dist}(x, \partial \Omega)$ so that $B(x, \eta_x)\subset \Omega$. Consider the mapping $F:\Omega\times (0,1)\to \mathbb{R}$ with
\begin{align*}
F(x,\varepsilon):= \int\limits_{\Omega\cap \{|x-y|\leq \eta_x\}} \hspace{-2ex}|u(x)-u(y)|^p \nu_\varepsilon(x-y)\mathrm{d}y= \int\limits_{|h|\leq \eta_x} |u(x)-u(x+h)|^p \nu_\varepsilon(h)\,\d h. 
\end{align*} 
In virtue of the fundamental theorem of calculus, we have 
\begin{align*}
F(x,\varepsilon)&=\int\limits_{|h|\leq \eta_x} \Big|\int_{0}^{1}\nabla u(x+th)\cdot h\mathrm{d}t\Big|^p \nu_\varepsilon(h)\,\d h \\
&=\int\limits_{|h|\leq \eta_x} |\nabla u(x)\cdot h|^p\, \nu_\varepsilon(h)\,\d h + R(x,\varepsilon)
\end{align*} 
with
 \begin{align*}
R(x,\varepsilon) = \int\limits_{|h|\leq \eta_x} \left(\Big|\int_{0}^{1}\nabla u(x+th)\cdot h\mathrm{d}t\Big|^p- \Big|\int_{0}^{1}\nabla u(x)\cdot h\mathrm{d}t\Big|^p \right)\nu_\varepsilon(h)\,\d h.
\end{align*}
Consider the function $s \mapsto G_p(s)=|s|^p$ which belongs to $C^1(\mathbb{R}^d\setminus \{0\})$ and $G'_p(s) = pG_p(s) s^{-1}$. Further, we are able to write
\begin{align*}
G_p(b)-G_p(a)	= (b-a) \int_{0}^{1} G'_p(a+s(b-a))\,\d s. 
\end{align*}
Applying this with $a= \nabla u (x)\cdot h$ and $b= \int_0^1\nabla u (x+th)\cdot h \,\d t$ by taking into account \eqref{eq:continuity-estimate} leads to the following estimates
\begin{align*}
|G_p(b)-G_p(a)|	&\leq |b-a| \int_{0}^{1} |a+s(b-a)|^{p-1}\,\d s\\
&\leq c_p |\nabla u(x)|^{p-1} |h|^{p-1}\int_0^1 |\nabla u(x+th)-\nabla u(x)||h|\,\d t\\
&\leq c_p \sigma|\nabla u(x)|^{p-1} |h|^{p}. 
\end{align*}
From this we get 
	\begin{align*}
		|R(x,\varepsilon)|&:= \Big| \int\limits_{|h|\leq \eta_x} \left(\Big|\int_{0}^{1}\nabla u(x+th)\cdot h\mathrm{d}t\Big|^p- \Big|\int_{0}^{1}\nabla u(x)\cdot h\mathrm{d}t\Big|^p \right)\nu_\varepsilon(h)\,\d h\Big|\\
		&\leq c_p \sigma|\nabla u(x)|^{p-1} \int\limits_{|h|\leq \eta_x} |h|^p\nu_\varepsilon(h)\,\d h\\
		&= c_p \sigma |\nabla u(x)|^{p-1} \int\limits_{|h|\leq \eta_x} (1\land |h|^p)\nu_\varepsilon(h)\,\d h\xrightarrow{\varepsilon, \sigma \to 0} 0.
	\end{align*}
%
On the other side, if we consider $O\in \mathcal{O}(d)$ to  be a rotation\footnote{An affine mapping $O$ on $\mathbb{R}^d$ is called a rotation if $O^\top O= I_d.$} assigning $\tfrac{\nabla u(x)}{|\nabla u(x)|}$ to the unit vector $e$ then passing through polar coordinates, the rotation invariance of the $d-1$-dimensional Lebesgue surface measure together with Remark \ref{rem:asymp-nu} we find that
\begin{align*}
\int\limits_{|h|\leq \eta_x}|\nabla u(x)\cdot h|^p \, \nu_\varepsilon(h)\,\d h
	&= \int_{0}^{\eta_x} r^{d+p-1} \mathrm{d}r \int\limits_{\mathbb{S}^{d-1}}|\nabla u(x)\cdot w|^p\mathrm{d}\sigma_{d-1}(w)\\
	&= |\nabla u(x)|^p \int\limits_{|h|\leq \eta_x}|h|^p\nu_\varepsilon(h)\,\d h \fint\limits_{\mathbb{S}^{d-1}}\big| w\cdot \tfrac{\nabla u(x)}{|\nabla u(x)|}\big|^p\mathrm{d}\sigma_{d-1}(w)\\
	&= |\nabla u(x)|^p \int\limits_{|h|\leq \eta_x}|h|^p\nu_\varepsilon(h)\,\d h \fint\limits_{\mathbb{S}^{d-1}}| w\cdot e|^p\mathrm{d}\sigma_{d-1}(w)\\
	&=K_{d,p} |\nabla u(x)|^p \int\limits_{|h|\leq \eta_x}|h|^p\nu_\varepsilon(h)\,\d h \xrightarrow{\varepsilon \to 0} K_{d,p} |\nabla u(x)|^p.
\end{align*}
Therefore, we have shown that 
\begin{align*}
\lim_{\varepsilon\to 0}F(x, \varepsilon) = K_{d,p} |\nabla u(x)|^p.
\end{align*}
Further, a close look to our reasoning reveals that we have subsequently shown that for all $\delta>0$
%
%
\begin{align}\label{eq:pointwise-close-limit}
\lim_{\varepsilon\to 0}\il_{\Omega\cap \{|x-y|\leq \delta\}}\hspace{-2ex} |u(x)-u(y)|^p\nu_\varepsilon(x-y)\,\d y = K_{d,p} |\nabla u(x)|^p.
\end{align}
In fact since for all $\delta >0$
\begin{align}\label{eq:limit-point-far}
\int\limits_{ \Omega\cap \{|x-y|\geq\delta \}} \hspace{-3ex}|u(x)-u(y)|^p \nu_\varepsilon(x-y) \, \mathrm{d}y 
\leq 2^p\|u\|^p_{L^{\infty}(\mathbb{R}^d)} \il_{|h|\geq \delta} \nu_\varepsilon(h) \mathrm{d}h
\xrightarrow{\varepsilon \to 0} 0.
\end{align}
We thus have the  pointwise convergence as claimed 
\begin{align*}
\lim_{\varepsilon\to 0}\int\limits_{\Omega} |u(x)-u(y)|^p \nu_\varepsilon(x-y) \d y\,\d x = K_{d,p} |\nabla u(x)|^p. 
\end{align*}
In order to proceed with the convergence in $L^1(\Omega)$, for $0<\eta <1$ as above we write
\begin{align*}
\iil_{\Omega\Omega}|u(x)-u(y)|^p\nu_\varepsilon(x-y)\,\d y \,\d x
 &= \iil_{\Omega\times \Omega\cap 
	\{|x-y|\leq \eta\}}|u(x)-u(y)|^p\nu_\varepsilon(x-y)\,\d y \,\d x\\
&+ \iil_{\Omega\times \Omega\cap 
\{|x-y|>\eta\}} \hspace{-2ex} |u(x)-u(y)|^p\nu_\varepsilon(x-y)\,\d y \,\d x.
\end{align*}
First and foremost, observe that for this specific choice of $0<\eta<1$, using \eqref{eq:continuity-estimate} one gets
\begin{align*}
 \il_{ \Omega\cap 
	\{|x-y|\leq \eta\}} \hspace{-2ex}|u(x)-u(y)|^p\nu_\varepsilon(x-y)\,\d y 
&\leq \il_{|h|\leq \eta} |u(x)-u(x+h)|^p\nu_\varepsilon(h)\,\d h \\
&=\il_{ |h|\leq \eta} \Big|\int_0^1 \nabla u(x+th)\cdot h \,\d t \Big|^p\nu_\varepsilon(h)\,\d h\\
&\leq 2^p|\nabla u(x)|^p \il_{ |h|\leq \eta} |h|^p\nu_\varepsilon(h)\,\d h \leq 2^p |\nabla u(x)|^p. 
\end{align*}
Since $x\mapsto |\nabla u(x)|^p$ belongs to $L^1(\Omega)$, the above estimate in combination with the pointwise limit in \eqref{eq:pointwise-close-limit} and the dominated convergence theorem yield that 
\begin{align*}
\lim_{\varepsilon \to 0}\iil_{ \Omega\times \Omega \cap 
	\{|x-y|\leq \eta\}} \hspace{-2ex}|u(x)-u(y)|^p\nu_\varepsilon(x-y)\,\d y \,\d x=K_{d,p} \int_\Omega |\nabla u(x)|^p \,\d x.
\end{align*}
Thus it follows that 
\begin{align*}
\lim_{\varepsilon \to 0}\iil_{ \Omega\times \Omega} \hspace{-1ex}|u(x)-u(y)|^p\nu_\varepsilon(x-y)\,\d y \,\d x= K_{d,p}\int_\Omega |\nabla u(x)|^p \,\d x
\end{align*}
 since one has
 \begin{align*}
 \iint\limits_{\Omega \times \Omega\cap \{|x-y|>\eta \}} \hspace{-3ex}|u(x)-u(y)|^p \nu_\varepsilon(x-y) \d y\,\d x 
 \leq 2^p\|u\|^p_{L^{p}(\Omega)} \il_{|h|>\eta} \nu_\varepsilon(h) \mathrm{d}h
 \xrightarrow{\varepsilon \to 0} 0. 
\end{align*}	
The claimed convergence in $L^1(\Omega)$ is a straightforward consequence of Corollary \ref{cor:Scheffe-lemma} (Sch\'effé lemma). Moreover, the explicit value of $K_{d,p}$ is already established in Proposition \ref{prop:asymp-cds}.

\end{proof}

\medskip

\begin{proof}[\textbf{Alternative proof of Lemma \ref{lem:BBM-regular} under a weaker condition}]
	We assume in addition that $u \in C_c^2(\mathbb{R}^d)$. Note that $G_p\in C^2(\mathbb{R}^d\setminus \{0\})$ with $G_p(s)=|s|^p$. The Taylor formula implies 
	\begin{align*}
	&u(y)-u(x)= \nabla u(x)\cdot (y-x) + O(|x-y|^2)\quad x,y\in \R^d, \\\\
	&	G_p(b)- G_p(a) = G'_p(a)(b-a) +O(b-a)^2\quad a,b\in \R\setminus\{0\}.
	\end{align*}
	Hence for almost all $x, y\in \R^d$, we have 
	\begin{align*}
	|u(y)-u(x)|^p =G_p( \nabla u(x)\cdot (y-x) + O(|y-x|^2))= |\nabla u(x)\cdot (y-x)|^p + O(|y-x|^{p+1}).
	\end{align*}
	Set $\delta_x= \operatorname{dist}(x,\partial\Omega)$. Passing through polar coordinates and  the rotation invariance of the Lebesgue measure yields
	\begin{align*}
	\int\limits_{ B(x,\delta_x)}\hspace{-2ex} |u(x)-u(y)|^p\nu_\varepsilon(x-y) \mathrm{d}y 
	&= \int\limits_{|h|\leq \delta_x} \Big|\nabla u(x)\cdot h\Big |^p \nu_\varepsilon(h) \mathrm{d}h+ O\Big( \int\limits_{|h|\leq \delta_x} |h|^{p+1}\nu_\varepsilon(h) \mathrm{d}h \Big)\\
	&= \int\limits_{\mathbb{S}^{d-1}} \Big|\nabla u(x)\cdot w \Big|^p\mathrm{d}\sigma_{d-1} (w) \int_{0}^{\delta_x} r^{d-1}\nu_\varepsilon(r) \mathrm{d}r+ O\Big( \hspace{-1ex}\int\limits_{|h|\leq \delta_x} |h|^{p+1} \nu_\varepsilon(h) \mathrm{d}h \Big)\\
	&= K_{d,p}\left|\nabla u(x)\right|^p \int\limits_{|h|\leq \delta_x} \nu_\varepsilon(h) \mathrm{d}h + O\Big( \int\limits_{|h|\leq \delta_x} |h|^{p+1} \nu_\varepsilon(h) \mathrm{d}h \Big).
	\end{align*}
	
	
	\noindent Therefore, letting $\varepsilon\to0 $ in the latter expression, by taking into account Remark \ref{rem:asymp-nu} gives 
	\begin{align*}
	&\lim_{\varepsilon\to 0} \int\limits_{ B(x,\delta_x)} |u(x)-u(y)|^p\nu_\varepsilon(x-y) \mathrm{d}y 
	= K_{d,p}|\nabla u(x)|^p . 
	\end{align*}
	
	\noindent On the other side, we have 
	\begin{align*}
	\il_{\Omega\setminus B(x,\delta_x)} |u(x)-u(y)|^p\nu_\varepsilon(x-y) \mathrm{d}y \leq 2^p\| u\|^p_{L^\infty(\Omega)} \il_{|x-y|\geq \delta_x} \nu_\varepsilon(x-y) \mathrm{d}y\xrightarrow{\eps\to0}0.
	\end{align*}
	We have proved that 
	\begin{align*}
	\lim_{\varepsilon\to 0} \int\limits_{ \Omega} |u(x)-u(y)|^p\nu_\varepsilon(x-y) \mathrm{d}y = K_{d,p}|\nabla u(x)|^p . 
	\end{align*}
	For the remaining details one can proceed as in the previous proof. 
\end{proof}

\noindent Following \cite[Theorem 2]{BBM01}, we derive the next theorem which is a  direct consequence of Lemma \ref{lem:BBM-regular}.
\begin{theorem}\label{thm:BBM-result}
Let $\Omega\subset \mathbb{R}^d$ be a $W^{1,p}$-extension set (not necessarily bounded), $1\leq p<\infty$. Assume that the family $(\nu_\varepsilon)_\varepsilon$ fulfills the conditions \eqref{eq:normalized-inetrals} and \eqref{eq:concentration-property} then for all $u\in W^{1,p}(\Omega)$ 
\begin{align*}
\lim_{\varepsilon\to 0}\iil_{\Omega\Omega}|u(x)-u(y)|^p\nu_\varepsilon(x-y) \mathrm{d}y \,\d x= K_{d,p} \int_\Omega | \nabla u(x)|^p\,\d x. 
\end{align*}

\end{theorem}
 
\medskip 

\begin{proof}
By Lemma \ref{lem:boundedness-limsup} there is $C>0$ independent of $\varepsilon$ such that for $u, v\in W^{1,p}(\Omega)$,

\begin{align*}
\big|\|U_\varepsilon\|_{L^p(\Omega\times \Omega)}-\|V_\varepsilon\|_{L^p(\Omega\times \Omega)}\big| \leq\|U_\varepsilon-V_\varepsilon\|_{L^p(\Omega \times \Omega)}\leq C \|u-v\|_{W^{1,p}(\Omega)}. 
\intertext{where we set}
U_\varepsilon(x,y) = |u(x)-u(y)|\nu_\varepsilon^{1/p}(x-y)\quad \text{and}\quad V_\varepsilon(x,y) = |v(x)-v(y)|\nu_\varepsilon^{1/p}(x-y)\,.
\end{align*}
\noindent 	Therefore, it suffices to establish the result for $u$ in a dense subset of $W^{1,p}(\Omega)$. Note that,  $C_c^\infty(\mathbb{R}^d )$ is dense in $W^{1,p}(\Omega)$ since $\Omega$ is an $W^{1,p}$-extension domain. We conclude  by using Lemma \ref{lem:BBM-regular}.  Alternatively, the result follows by combining the forthcoming Theorem \ref{thm:liminf-BBM} and Theorem \ref{thm:limsup_BBM}. 
\end{proof}

\begin{remark}
The counterexample \ref{Ex:counterexample-extension} shows that the conclusion of Theorem \ref{thm:BBM-result} might be erroneous whenever $\Omega$ is not an extension domain. However, even if $\Omega$ is not an extension domain,  the  next result, very close to Theorem \ref{thm:BBM-result}, shows that one can still get the following inequality which also appears in \cite{Ponce2004} and was first established in \cite{Brezis-const-function} for the case $\Omega=\R^d$.
\end{remark}

\begin{theorem}\label{thm:liminf-BBM} Let $\Omega\subset \mathbb{R}^d$ be open (not necessarily bounded). Assume that the family $(\nu_\varepsilon)_\varepsilon$ fulfills the conditions \eqref{eq:normalized-inetrals} and \eqref{eq:concentration-property}. Let $u\in L^p(\Omega)$ with $1 <p<\infty$ or $u\in W^{1,1}(\Omega)$. Adopting the convention that $\|\nabla u\|_{L^p(\Omega)} =\infty$ if $|\nabla u|$ does not belong to $L^p(\Omega)$ we have 
	\begin{align*}
 K_{d,p} \int_\Omega | \nabla u(x)|^p\,\d x\leq 	\liminf_{\varepsilon\to 0}\iil_{\Omega\Omega}|u(x)-u(y)|^p\nu_\varepsilon(x-y) \mathrm{d}y \,\d x .
	\end{align*}
	
\end{theorem}

\medskip

\begin{proof} For  $\delta >0$ sufficiently small, $\Omega_\delta =\{x\in \Omega: \operatorname{dist}(x,\partial\Omega)>\delta\}.$ Define the mollifier $\phi_\delta(x)= \frac{1}{\delta^d}\phi\left(\frac{x}{\delta}\right)$ with support in $B_\delta(0)$ where $\phi \in C_c^{\infty}(\mathbb{R}^d)$ is supported in $B_1(0)$, $\phi \geq 0$ and $ \int_{} \phi = 1$. 
For the sake of the simplicity we will assume $u$ is extended by zero off $\Omega$ and let $u^\delta = u*\phi_\delta$ denote the convolution product of $u $ and $\phi_\delta$. 

\noindent Assume $z\in\Omega_\delta $ and $|h|\le \delta$ then $z-h\in \Omega_\delta-h \subset \Omega$ so that by a simple change of variables, we have
	\begin{align*}
	\iil_{\Omega\Omega}|u(x) -u(y)|^p \nu_\varepsilon(x-y)\d y\,\d x
	&\geq\hspace{-2ex} \iil_{\Omega_\delta-h\Omega_\delta-h}|u(x) -u(y)|^p \nu_\varepsilon(x-y)\d y\,\d x\\
	&=\iil_{\Omega_\delta\Omega_\delta} |u(x-h) -u(y-h)|^p \nu_\varepsilon(x-y)\d x\d y.
	\end{align*}
	Thus given that $\int_{} \phi_\delta = 1$ integrating both sides over the ball $B_\delta(0)$ with respect to 
	$\phi_\delta(h)dh$ and employing Jensen's inequality afterwards, yields 
	\begin{align*}
	\iil_{\Omega\Omega}|u(x) -u(y)|^p \nu_\varepsilon(x-y)\d y\,\d x
&\geq \int_{\R^d} \phi_\delta(h)\d h \iil_{\Omega_\delta\Omega_\delta} |u(x-h) -u(y-h)|^p \nu_\varepsilon(x-y)\d y\,\d x\\ 
	&= \iil_{\Omega_\delta\Omega_\delta} \int_{\R^d} |u(x-h) -u(y-h)|^p \phi_\delta(h) \d h \,\nu_\varepsilon(x-y)\d x\d y\\ & \geq \iil_{\Omega_\delta\Omega_\delta} \Big| \int_{\R^d} u(x-h) -u(y-h) \phi_\delta(h) \d h\Big|^p \nu_\varepsilon(x-y)\d x\d y\\
	&= \iil_{\Omega_\delta\Omega_\delta} |u*\phi_\delta (x) -u*\phi_\delta (y)|^p \nu_\varepsilon(x-y)\d x\d y.
\end{align*}
	Thus we have 
	\begin{align}\label{eq:molification-convex-Jessen}
 \iil_{\Omega_\delta\Omega_\delta} |u^\delta (x) -u^\delta (y)|^p \nu_\varepsilon(x-y)\d x\d y\leq \iil_{\Omega\Omega}|u(x) -u(y)|^p \nu_\varepsilon(x-y)\d x\d y.
	\end{align}
Note that $u^\delta\in C^\infty(\R^d)$ and $\Omega_{\delta,j}= \Omega_\delta\cap B_j(0)$ has a compact closure for each $j\geq 1$. Then for each $j\geq 1$ the Lemma \ref{lem:BBM-regular} implies
	\begin{align*}
K_{d,p}\int_{\Omega_{\delta,j} } |\nabla u^\delta (x)|^p \d x
&=	\lim_{\varepsilon\to 0} \iil_{\Omega_{\delta,j}\Omega_{\delta,j}}|u(x) -u(y)|^p \nu_\varepsilon(x-y)\d x\d y\\
&\leq 	\liminf_{\varepsilon\to 0} \iil_{\Omega\Omega}|u(x) -u(y)|^p \nu_\varepsilon(x-y)\d x\d y.
\end{align*}
Tending $j\to\infty$ in the latter we get
		\begin{align}\label{eq:limit-approx}
K_{d,p}\int_{\Omega_\delta} |\nabla u^\delta (x)|^p \d x\leq 	\liminf_{\varepsilon\to 0} \iil_{\Omega\Omega}|u(x) -u(y)|^p \nu_\varepsilon(x-y)\d x\d y.
	\end{align}
The only interesting scenario occurs if the right hand side is finite. If it is the case, then the upcoming result Theorem \ref{thm:charact-w1p} ensures that $u\in W^{1,p}(\Omega)$. Thus clearly $\nabla u^\delta= \nabla( u*\phi_\delta)
	= \nabla u*\phi_\delta $ and we know that $\| \phi_\delta* \nabla u - \nabla u \|_{L^p(\Omega)}\to 0$ as $\delta \to 0$.
\noindent The desired inequality follows by letting $\delta\to 0$ since 
	\begin{align*}
\Big| \| \nabla u \|_{L^p(\Omega)} -\|\nabla u *\phi_\delta \|_{L^p(\Omega_\delta)} \Big|\leq \| \nabla u *\phi_\delta- \nabla u \|_{L^p(\Omega)} + \Big(\int_{\Omega\setminus \Omega_\delta} |\nabla u (x)|^p\d x\Big)^{1/p}\xrightarrow{\delta\to 0} 0.
\end{align*}
\end{proof}
\bigskip

\noindent The next theorem is a the counterpart of Theorem \ref{thm:liminf-BBM} and is a refinement version of Theorem \ref{thm:limpsup}. 
\begin{theorem}\label{thm:limsup_BBM} Let $\Omega\subset \mathbb{R}^d$ be a $W^{1,p}$-extension set (not necessarily bounded). Assume that the family $(\nu_\varepsilon)_\varepsilon$ fulfills the conditions \eqref{eq:normalized-inetrals} and \eqref{eq:concentration-property} with $1 \leq p<\infty$. Adopting the convention that $\|\nabla u\|_{L^p(\Omega)} =\infty$ if $|\nabla u|$ does not belong to $L^p(\Omega)$ then for all $u\in L^p(\Omega)$ we have 
	\begin{align*}
	\limsup_{\varepsilon\to 0}\iil_{\Omega\Omega}|u(x)-u(y)|^p\nu_\varepsilon(x-y) \mathrm{d}y \d x\leq	K_{d,p} \int_\Omega | \nabla u(x)|^p\,\d x\,.
	\end{align*}
\end{theorem}

\medskip

\begin{proof}
	The case where $\|\nabla u\|_{L^p(\Omega)} =\infty$ is trivial. Now for $u\in W^{1,p}(\Omega)$, we let $\overline{u}\in W^{1,p}(\R^d)$ be its extension to $\R^d$. Consider 
	$\Omega(\delta) = \Omega+ B_\delta(0) = \{x\in\R^d~:\operatorname{dist}(x, \Omega)<\delta \}$ be a neighborhood of $\Omega$ with $0<\delta<1$ sufficiently small. Let us start by proving the following estimate	which holds for each $\varepsilon$,
	\begin{align}\label{eq:split-estimate}
	\iil_{\Omega\Omega}|u(x)-u(y)|^p\nu_\varepsilon(x-y) \mathrm{d}y \d x\leq K_{d,p}\int_{\Omega(\delta) }|\nabla \overline{u}(x)|^p\d x+ 2^p\|u\|^p_{L^p(\Omega)}\int_{|h|\geq \delta}\nu_\varepsilon(h)\,\d h.
	\end{align}
	We know that there is $(u_n)_n$ a sequence of functions in $C_c^\infty(\R^d)$ converging to $\overline{u}$ in $W^{1,p}(\R^d)$. For each $n\geq 1$, exploiting the computations from the proof of Lemma \ref{lem:BBM-regular} (see the step before \eqref{eq:pointwise-close-limit}) we find that
	\begin{align*}
	\int_{\Omega} \il_{ |x-y|\leq \delta}|u_n(x)-u_n(y)|^p\nu_\varepsilon(x-y)\,\d y \,\d x
	&\leq \il_{|h|\leq \delta}\int_0^1\int_\Omega |\nabla u_n(x+th)\cdot h|^p\,\d x ~\,\d t~\nu_\varepsilon(h)\,\d h\\
	&\leq \il_{|h|\leq \delta}\int_{\Omega(\delta)}|\nabla u_n(z)\cdot h|^p\,\d z ~\nu_\varepsilon(h)\,\d h\\
	&= K_{d,p} \Big(\int_{\Omega(\delta)}|\nabla u_n(z)|^p\,\d z \Big)\Big(\il_{|h|\leq \delta}(1\land |h|^p)\nu_\varepsilon(h)\,\d h\Big)\\
	&\leq K_{d,p} \int_{\Omega(\delta)}|\nabla u_n(z)|^p\,\d z. 
	\end{align*}
	Therefrom, applying Fatou's lemma we get 
	\begin{align*}
	\int_{\Omega} \il_{ |x-y|\leq \delta}|u(x)-u(y)|^p\nu_\varepsilon(x-y)\,\d y \,\d x
	&\leq \liminf_{n\to\infty}\int_{\Omega} \il_{|x-y|\leq \delta}|u_n(x)-u_n(y)|^p\nu_\varepsilon(x-y)\,\d y \,\d x\\
	&\leq \lim_{n\to\infty} K_{d,p} \int_{\Omega(\delta)}|\nabla u_n(z)|^p\,\d z\\
	&= K_{d,p} \int_{\Omega(\delta)}\, |\nabla \overline{u}(x)|^p\d x.
	\end{align*}
	On the other hand,  we obviously have
	\begin{align*}
	\int_{\Omega} \int_{\Omega\cap \{|x-y|\geq \delta\}} \hspace{-3ex} |u(x)-u(y)|^p\nu_\varepsilon(x-y)\,\d y \,\d x \leq 2^p \|u\|^p_{L^p(\Omega)}\int_{|h|\geq \delta}\nu_\varepsilon(h)\,\d h. 
	\end{align*}
	Altogether, \eqref{eq:split-estimate} clearly follows. Therefore, letting $\varepsilon \to0$ and $\delta \to0$ successively the relation \eqref{eq:split-estimate} becomes
\begin{align*}
\limsup_{\varepsilon\to 0} \iil_{\Omega\Omega}|u(x)-u(y)|^p\nu_\varepsilon(x-y) \mathrm{d}y \d x\leq K_{d,p}\int_{\Omega}|\nabla u(x)|^p\d x
\end{align*}
since recalling that $u= \overline{u}\mid_\Omega$, we have 
$$\int_{|h|\geq \delta}\nu_\varepsilon(h)\,\d h\xrightarrow{\varepsilon\to 0}0\qquad\text{and}\qquad \int_{\Omega(\delta) }|\nabla \overline{u}(x)|^p\d x \xrightarrow{\delta\to 0} \int_{\Omega}|\nabla u(x)|^p\d x.$$
 The desired estimate is proved. 
 
\end{proof}

\begin{remark}
	The counterexample \ref{Ex:counterexample-extension} clearly shows that the conclusion of Theorem \ref{thm:limsup_BBM} might be false whenever $\Omega$ is not an extension domain. It is also highly remarkable that Theorem \ref{thm:liminf-BBM} and Theorem \ref{Ex:counterexample-extension} together imply Theorem \ref{thm:BBM-result}. 
\end{remark}

\noindent The following  lemma is relevant for the sequel and is somewhat a revisited version of \cite[Lemma 1]{BBM01}.

\begin{lemma}\label{lem:boun-integration-bypart}	Let $ u\in L^p(\mathbb{R}^d)$, $\nu \in L^1(\mathbb{R}^d, 1\land |h|^p )$ with $1\leq p<\infty$ and $\varphi \in C_c^\infty(\mathbb{R}^d)$ then for any unit vector $e\in \mathbb{S}^{d-1}$ we have 
\begin{relsize}{-0}
		\begin{align*}
	\Big | \iint\limits_{(y-x)\cdot e\geq 0} \hspace{-3ex} u(x) \frac{\varphi(y) - \varphi(x)}{|x-y|}(1\land |x-y|^p) \nu(x-y) \mathrm{d}y \mathrm{d}x\Big |&+ \Big | \iint\limits_{(y-x)\cdot e\leq 0} \hspace{-3ex} u(x) \frac{\varphi(y) - \varphi(x)}{|x-y|}(1\land |x-y|^p) \nu(x-y) \mathrm{d}y \mathrm{d}x\Big |\\
	&\leq \iint\limits_{\mathbb{R}^d\mathbb{R}^d} \frac{|u(x)-u(y)|}{|x-y|}|\varphi(x)|(1\land |x-y|^p) \nu(x-y) \mathrm{d}y \mathrm{d}x. 
	\end{align*}
\end{relsize}
	
\end{lemma}

\medskip
\begin{proof}
	We begin by introducing $\widetilde{\nu_\delta}(h)= |h|^{-1}(1\land |h|^p) \nu(h)~\mathds{1}_{\mathbb{R}^d\setminus B_\delta }(h)$ for $\delta>0$ which enables us to rule out an eventual singularity of $\nu$ at the origin. Moreover, note that $\nu_\delta\in L^1(\mathbb{R}^d)$. It tuns that the mappings $(x,y)\mapsto u(x) \varphi(y) \widetilde{\nu_\delta}(x-y)$ and $(x,y)\mapsto u(x) \varphi(x) \widetilde{\nu_\delta}(x-y)$ are integrable. Indeed, using H\"older inequality combined with Fubini's theorem lead to the following
	%
	\begin{align*}
	&\iil_{\mathbb{R}^d\mathbb{R}^d} |u(x) \varphi(x)| \widetilde{\nu_\delta}(x-y)\,\d y \,\d x 
	= \iil_{{ |x-y|\geq \delta}} |u(x) \varphi(x)| |x-y|^{-1}(1\land |x-y|^p) \nu(x-y) \,\d y \,\d x \\
	&\leq \delta^{-1} \Big(\hspace{-2ex}\iil_{{ |x-y|\geq \delta}} \hspace{-2ex}| u(x)|^p (1\land |x-y|^{p})\nu(x-y) \,\d y \,\d x \Big)^{1/p} 
	\Big(\hspace{-2ex}\iil_{{ |x-y|\geq \delta}} \hspace{-2ex} |\varphi(x)|^{p'}(1\land |x-y|^p) \nu(x-y) \,\d y \,\d x\Big)^{1/p'} \\
	&\leq \delta^{-1}\|\varphi\|_{L^{p'}(\mathbb{R}^d)}\|u\|_{L^p(\mathbb{R}^d)}\int_{\mathbb{R}^d}(1\land |h|^p)\nu(h)\,\d h<\infty.
	\end{align*}
	%
	%
	Likewise we also have
	\begin{align*}
	\iil_{\mathbb{R}^d \mathbb{R}^d} |u(x) \varphi(y)| \widetilde{\nu_\delta}(x-y)\,\d y \,\d x 
	&\leq \delta^{-1}\|\varphi\|_{L^{p'}(\mathbb{R}^d)}\|u\|_{L^p(\mathbb{R}^d)}\int_{\mathbb{R}^d}(1\land |h|^p)\nu(h)\,\d h<\infty.
	\end{align*}
Consequently, by interchanging $x$ and $y$, using Fubini's theorem and the symmetry of $\nu$ we obtain 
		\begin{alignat*}{2}
	& \iil_{(y-x)\cdot e\geq 0} u(x) \varphi(x) \widetilde{\nu_\delta}(x-y)\,\d y \,\d x 
	& &= \iil_{(x-y)\cdot e\geq 0} u(y) \varphi(y)\widetilde{\nu_\delta}(x-y)\,\d y \,\d x \\
	&= \int_{\mathbb{R}^d} u(y)\varphi(y)\,\d y \il_{h\cdot e\geq 0} \widetilde{\nu_\delta}(h) \,\d h
	& &= \int_{\mathbb{R}^d} u(y)\varphi(y)\,\d y \il_{h\cdot e\leq 0} \widetilde{\nu_\delta}(h) \,\d h\\
	&= \int_{\mathbb{R}^d} u(y)\varphi(y)\,\d y \il_{h\cdot e\geq 0} \widetilde{\nu_\delta}(h) \,\d h
	& &= \int_{\mathbb{R}^d} u(y)\varphi(y)\,\d y \il_{(y-x)\cdot e\geq 0} \widetilde{\nu_\delta}(y-x) \,\d x. 
	\end{alignat*}
	Therefore, we have 
\begin{relsize}{-0}
	\begin{align*}
	\hspace{-2ex} \Big|\int_{\mathbb{R}^d} u(x)\,\d x \hspace{-3ex}\il_{(y-x)\cdot e\geq 0} \hspace*{-3ex}(\varphi(y)-\varphi(x)) \widetilde{\nu_\delta}(x-y)\,\d y\Big| 
	&=\Big |\iil_{(y-x)\cdot e\geq 0} \hspace*{-3ex} u(x) \varphi(y) \widetilde{\nu_\delta}(x-y)\,\d y \,\d x -\hspace{-3ex} \iil_{(y-x)\cdot e\geq 0} \hspace*{-3ex} u(x) \varphi(x) \widetilde{\nu_\delta}(x-y)\,\d y \,\d x \Big| \\
	&=\Big |\hspace{-3ex}\iil_{(y-x)\cdot e\geq 0} \hspace{-3ex}u(x) \varphi(y) \widetilde{\nu_\delta}(x-y)\,\d y \,\d x - \hspace{-3ex}\iil_{(y-x)\cdot e\geq 0} u(y) \varphi(y)| \widetilde{\nu_\delta}(x-y)\,\d y \,\d x \Big| \\
	&=\Big | \int_{\mathbb{R}^d} \varphi(y) \,\d y \il_{(y-x)\cdot e\geq 0} \hspace*{-3ex} (u(x)-u(y)) \widetilde{\nu_\delta}(x-y) \,\d x\Big|
	\\& \leq \int_{\mathbb{R}^d} |\varphi(y)|\,\d y \il_{(y-x)\cdot e\geq 0} \hspace*{-3ex} \frac{|u(x)-u(y)|}{|x-y|} (1\land |x-y|^p)\nu(x-y) \,\d x\\
	&= \int_{\mathbb{R}^d} |\varphi(x)|\,\d x \il_{(y-x)\cdot e\leq 0} \hspace*{-3ex} \frac{|u(y)-u(x)|}{|x-y|} (1\land |x-y|^p)\nu(x-y) \,\d y.
	\end{align*}
	\end{relsize}
Thus, letting $\delta \to 0$ implies 
\begin{relsize}{-0.5}
	\begin{align}\label{eq:delta-estimate1}
	\hspace{-2ex}\Big| \hspace{-3ex}\iil_{(y-x)\cdot e\geq 0} \hspace{-3ex} u(x)(\varphi(y)-\varphi(x)) (1\land |x-y|^p)\nu(x-y) \,\d y \,\d x\Big| 
	\leq \hspace{-3ex}\iil_{(y-x)\cdot e\leq 0} \hspace*{-3ex}|\varphi(x)| \frac{|u(y)-u(x)|}{|x-y|} (1\land |x-y|^p)\nu(x-y) \,\d y\,\d x.
	\end{align}
\end{relsize}
Likewise one has
\begin{relsize}{-0.5}
	\begin{align}\label{eq:delta-estimate2}
		\Big| \hspace{-3ex}\iil_{(y-x)\cdot e\leq 0} \hspace{-3ex} u(x)(\varphi(y)-\varphi(x))(1\land |x-y|^p)\nu(x-y)\,\d y \,\d x\Big| 
	\leq \hspace{-2ex} \iil_{(y-x)\cdot e\geq 0} \hspace*{-3ex}|\varphi(x)| \frac{|u(y)-u(x)|}{|x-y|} (1\land |x-y|^p)\nu(x-y) \,\d y\,\d x.
	\end{align}
\end{relsize}
	Adding \eqref{eq:delta-estimate1} and \eqref{eq:delta-estimate2} gives the desired inequality
	\begin{relsize}{-.4}
	\begin{align*}
	&\Big|\int_{\mathbb{R}^d} u(x)\,\d x \hspace{-3ex}\il_{(y-x)\cdot e\geq 0} \hspace*{-3ex} \frac{\varphi(y)-\varphi(x)) }{|x-y|} (1\land |x-y|^p)\nu(x-y) \,\d y\Big| + \Big|\int_{\mathbb{R}^d} u(x)\,\d x \hspace{-3ex}\il_{(y-x)\cdot e\leq 0} \hspace*{-3ex} \frac{\varphi(y)-\varphi(x)) }{|x-y|} (1\land |x-y|^p)\nu(x-y) \,\d y\Big| \\
	&\leq \int_{\mathbb{R}^d} |\varphi(x)|\,\d x\hspace{-3ex} \il_{(y-x)\cdot e\geq 0} \hspace*{-3ex} \frac{|u(y)-u(x)|}{|x-y|} (1\land |x-y|^p)\nu(x-y) \,\d y+ \int_{\mathbb{R}^d} \hspace{-1ex}|\varphi(x)|\,\d x \hspace{-3ex}\il_{(y-x)\cdot e\leq 0} \hspace*{-3ex} \frac{|u(y)-u(x)|}{|x-y|} (1\land |x-y|^p)\nu(x-y) \,\d y\\
	&= \int_{\mathbb{R}^d} |\varphi(x)|\,\d x \int_{\mathbb{R}^d} \frac{|u(y)-u(x)|}{|x-y|} (1\land |x-y|^p)\nu(x-y) \,\d y. 
	\end{align*}
	\end{relsize}
\end{proof}

\medskip

\begin{proposition}\label{prop:bounded-liminf-estimate} Assume $\Omega\subset \R^d$ is an open set. Let $u \in L^p(\Omega)$ with $1\leq p<\infty $. Then for every $\varphi\in C_c^\infty(\mathbb{R}^d)$ with support in $\Omega$ and for every unit vector $e\in \mathbb{S}^{d-1}$ the following estimate holds true
\begin{align}\label{eq:cont-weak-derivatve}
\left|\int_{\Omega} u(x) \nabla \varphi (x)\cdot e ~\mathrm{d}x \right| \leq \frac{A_p^{1/p}}{ K_{d,1}}\|\varphi\|_{L^{p'}(\Omega)}.
\end{align} 

With
\begin{align*}
A_p:= \liminf_{\varepsilon \to 0} \iil_{\Omega \Omega}|u(y)-u(x)|^p\nu_\varepsilon(x-y) \,\d x\,\d y 
\end{align*} 
and the family $(\nu_\varepsilon)_\varepsilon$ satisfies the conditions \eqref{eq:normalized-inetrals} and \eqref{eq:concentration-property}. 

\end{proposition}

\medskip 

\begin{proof}
 Throughout, to alleviate the notation we denote $\pi_\varepsilon(x-y)= (1\land |x-y|^p)\nu_\varepsilon(x-y)$. Let $\overline{u}\in L^p(\mathbb{R}^d)$ be the zero extension of $u$ outside $\Omega$. Let $\varphi\in C_c^\infty(\mathbb{R}^d)$ with support in $\Omega$. First of all we have the identity
\begin{align*}
\int_{\mathbb{R}^d} |\varphi(x)|\,\d x \int_{\mathbb{R}^d} \frac{|\overline{u}(y)-\overline{u}(x)|}{|x-y|} \pi_\varepsilon(x-y) \,\d y
&= \iil_{\Omega \Omega}\frac{|u(y)-u(x)|}{|x-y|} |\varphi(x)| \pi_\varepsilon(x-y) \,\d x\,\d y \\
&+ \il_{\supp(\varphi)} |\varphi(x)|\,\d x \il_{\mathbb{R}^d\setminus \Omega } \frac{|u(x)|}{|x-y|} \pi_\varepsilon(x-y) \,\d y
\end{align*}
There are two keys observations. First, since $\delta= \operatorname{dist}(\supp(\varphi), \partial\Omega)>0$, H\"older inequality implies 
\begin{align*}
\hspace{-2ex}\il_{\supp(\varphi)} \hspace{-1ex}|\varphi(x)|\,\d x\hspace{-1ex} \il_{\mathbb{R}^d\setminus \Omega }\hspace{-1ex} \frac{|u(x)|}{|x-y|} \pi_\varepsilon(x-y) \,\d y\leq 
\delta^{-1} \|u\|_{L^{p}(\Omega)} \|\varphi\|_{L^{p'}(\Omega)} \hspace{-1ex}\il_{|h|\geq \delta } \hspace{-1ex}(1\land |h|^p)\nu_\varepsilon(h) \,\d h
\xrightarrow{\varepsilon\to 0}0. 
\end{align*}
Second, from H\"older inequality and $|h|^{-p}(1\land |h|^p)\leq 1$ we find that
\begin{align*}
\hspace{-1ex}\iil_{\Omega \Omega}\frac{|u(y)-u(x)|}{|x-y|} |\varphi(x)| &\pi_\varepsilon(x-y) \,\d x\,\d y\\
& \leq \Big( \iil_{\Omega \Omega}\frac{|u(y)-u(x)|^p}{|x-y|^p} \pi_\varepsilon(x-y) \,\d x\,\d y \Big)^{1/p} \Big( \iil_{\Omega \Omega}|\varphi(x)|^{p'} \pi_\varepsilon(x-y) \,\d x\,\d y \Big)^{1/p'}\\
&\leq  \|\varphi\|_{L^{p'}(\Omega)}\Big( \iil_{\Omega \Omega}|u(y)-u(x)|^p\nu_\varepsilon(x-y) \,\d x\,\d y \Big)^{1/p} \Big( \int_{\mathbb{R}^d }(1\land |h|^p)\nu_\varepsilon(h) \,\d h\Big)^{1/p'}\\
&\leq  \|\varphi\|_{L^{p'}(\Omega)}\Big( \iil_{\Omega \Omega}|u(y)-u(x)|^p\nu_\varepsilon(x-y) \,\d x\,\d y \Big)^{1/p} . 
\end{align*}
\noindent Therefore inserting these two observations in the previous identity and combining the resulting estimate with that of Lemma \ref{lem:boun-integration-bypart} and that of the assumptions imply
%
	\begin{align}\label{eq:nonlocal-liminf}
\begin{split}
&\liminf_{\varepsilon \to 0}\Big|\int_{\Omega} u(x)\,\d x \hspace{-3ex} \il_{(y-x)\cdot e\geq 0} \hspace{-3ex} \frac{(\varphi(y)-\varphi(x)) }{|x-y|} (1\land |x-y|^p)\nu_\varepsilon(x-y)\,\d y \Big|~+\\
&\liminf_{\varepsilon \to 0} \Big|\int_{\Omega} u(x)\,\d x \hspace{-3ex} \il_{(y-x)\cdot e\leq 0} \hspace{-3ex}\frac{(\varphi(y)-\varphi(x)) }{|x-y|}(1\land |x-y|^p)\nu_\varepsilon(x-y) \,\d y \Big|\leq A_p^{1/p} \|\varphi\|_{L^{p'}(\Omega)}.
\end{split}
\end{align}
%
\noindent The next step of the proof will consist into computing the limits appearing on the left hand side \eqref{eq:nonlocal-liminf}. We have
\begin{align*}
\lim_{\varepsilon \to 0} \il_{(y-x)\cdot e\geq 0} \hspace{-3ex} &\frac{(\varphi(y)-\varphi(x)) }{|x-y|} (1\land |x-y|^p)\nu_\varepsilon(x-y) \,\d y\\ 
&= \lim_{\varepsilon \to 0} \il_{h\cdot e\geq 0} \int_0^1\nabla \varphi(x+th)\cdot \frac{h}{|h|} \,\d t (1\land |h|^p)\nu_\varepsilon(h) \,\d h\\
&= \lim_{\varepsilon \to 0} \il_{h\cdot e\geq 0} \int_0^1\big[\nabla \varphi(x+th)-\nabla\varphi(x)\big]\cdot \frac{h}{|h|}\,\d t (1\land |h|^p)\nu_\varepsilon(h) \,\d h \\
&+ \lim_{\varepsilon \to 0} \il_{h\cdot e\geq 0} \nabla\varphi(x)\cdot \frac{h}{|h|} (1\land |h|^p)\nu_\varepsilon(h) \,\d h . 
\end{align*}
\noindent Recall that $|\nabla \varphi (x+z)- \nabla \varphi(x)|\leq C (1\land |z|)$ for all $x, z\in \R^d$ so that using Remark \ref{rem:asymp-nu} we obtain
 \begin{align*}
 	\lim_{\varepsilon \to 0} \il_{h\cdot e\geq 0} \int_0^1\Big|\big[\nabla \varphi(x+th)-\nabla\varphi(x)\big]\cdot \frac{h}{|h|}\Big|\,\d t (1\land |h|^p)\nu_\varepsilon(h) \,\d h\leq C\lim_{\varepsilon \to 0} \il_{\mathbb{R}^d} (1\land |h|^{p+1})\nu_\varepsilon(h) \,\d h=0. 
 \end{align*}
Thus, passing through polar coordinates in the other term of the above expression and taking into account \eqref{eq:normalized-inetrals} implies 
\begin{align*}
\lim_{\varepsilon \to 0} \il_{(y-x)\cdot e\geq 0} \hspace{-3ex} \frac{(\varphi(y)-\varphi(x)) }{|x-y|} &(1\land |x-y|^p)\nu_\varepsilon(x-y) \,\d y\\
&= \lim_{\varepsilon \to 0} \il_{h\cdot e\geq 0} \nabla\varphi(x)\cdot \frac{h}{|h|} (1\land |h|^p)\nu_\varepsilon(h) \,\d h\\
&= \lim_{\varepsilon \to 0} \il_{\mathbb{S}^{d-1} \cap \{w\cdot e\geq 0\}} \hspace{-4ex}\nabla\varphi(x)\cdot w \,\d \sigma_{d-1}(w)\int_0^\infty (1\land |r|^p)\nu_\varepsilon(r) \,\d r \\
&= \Big(\il_{\mathbb{S}^{d-1} \cap \{w\cdot e\geq 0\}} \hspace{-4ex}\nabla\varphi(x)\cdot w \,\d \sigma_{d-1}(w) \Big)\times |\mathbb{S}^{d-1}|^{-1}\lim_{\varepsilon \to 0} \int_{\mathbb{R}^d} (1\land |h|^p)\nu_\varepsilon(h) \,\d h\\
&= |\mathbb{S}^{d-1}|^{-1}\il_{\mathbb{S}^{d-1} \cap \{w\cdot e\geq 0\}} \hspace{-3ex}\nabla\varphi(x)\cdot w \,\d \sigma_{d-1}(w). 
\end{align*}

\noindent Let $(e, v_2, \cdots v_d)$ be an orthonormal basis of $\mathbb{R}^d$ in which we write the coordinates $w= (w_1,w_2, \cdots, w_d)= (w_1,w')$ that is $w_1 = w\cdot e$ and $ w_i= w\cdot v_i$. Similarly, in this basis one has $\nabla \varphi (x) =(\nabla \varphi (x)\cdot e, (\nabla \varphi (x))' )$. Observe that $\nabla \varphi (x) \cdot w = [\nabla \phi (x)]_1 w_1+ \cdots +[\nabla \phi (x)]_d w_d= \big(\nabla \varphi (x)\cdot e\big) (w\cdot e)+ [\nabla \phi (x)]'\cdot w' $. From this we find that
\begin{align*}
\il_{\mathbb{S}^{d-1} \cap \{w\cdot e\geq 0\}} \hspace{-3ex} \nabla \varphi (x) \cdot w d \sigma_{d-1}(w) = \il_{\mathbb{S}^{d-1} \cap \{w\cdot e\geq 0\}}\hspace{-3ex} ( \nabla \varphi (x)\cdot e )( w\cdot e) d \sigma_{d-1}(w)+ \il_{\mathbb{S}^{d-1} \cap \{w\cdot e\geq 0\}} \hspace{-3ex} (\nabla \varphi (x))' \cdot w' d \sigma_{d-1}(w). 
\end{align*}
\noindent Consider the rotation $O(w) = (w_1, -w') = (w\cdot e, -w') $ then the rotation invariance of the Lebesgue measure entails that $\d \sigma_{d-1}(w) = \d \sigma(O(w)) $ and we have 
\begin{align*}
\il_{\mathbb{S}^{d-1} \cap \{w\cdot e\geq 0\}} \hspace{-3ex} (\nabla \varphi (x))' \cdot w' d \sigma_{d-1}(w) = -\il_{\mathbb{S}^{d-1} \cap \{w\cdot e\geq 0\}} \hspace{-3ex} (\nabla \varphi (x))' \cdot w' d \sigma_{d-1}(w) =0. 
\end{align*}
Whereas, by symmetry we have 
\begin{align*}
\il_{\mathbb{S}^{d-1} \cap \{w\cdot e\geq 0\}} \hspace{-3ex} w\cdot e \, \,\d \sigma_{d-1}(w) = -\il_{\mathbb{S}^{d-1} \cap \{w\cdot e\leq 0\}} \hspace{-3ex} w\cdot e d \sigma_{d-1}(w) = \frac{ 1}{2} \il_{\mathbb{S}^{d-1}} |w\cdot e |d \sigma_{d-1}(w). 
\end{align*}

\noindent Altogether yields that 
\begin{align*}
|\mathbb{S}^{d-1}|^{-1}\il_{\mathbb{S}^{d-1} \cap \{w\cdot e\geq 0\}}\hspace{-4ex} \nabla\varphi(x)\cdot w \,\d \sigma_{d-1}(w)= \frac{\nabla\varphi(x)\cdot e}{2}\fint_{\mathbb{S}^{d-1}} | w \cdot e|\,\d \sigma_{d-1}(w) = \frac{ 1}{2}K_{d,1}\nabla \varphi (x)\cdot e. 
\end{align*}
In conclusion,
\begin{align}\label{eq:K1d}
&\lim_{\varepsilon \to 0} \il_{(y-x)\cdot e\geq 0} \hspace{-3ex} \frac{(\varphi(y)-\varphi(x)) }{|x-y|} (1\land |x-y|^p)\nu_\varepsilon(x-y) \,\d y= \frac{1}{2}K_{d,1}\nabla \varphi(x)\cdot e.
\intertext{Analogously one is able to show that}\label{eq:K1d-p}
&\lim_{\varepsilon \to 0} \il_{(y-x)\cdot e\leq 0} \hspace{-3ex} \frac{(\varphi(y)-\varphi(x)) }{|x-y|} (1\land |x-y|^p)\nu_\varepsilon(x-y) \,\d y = \frac{1}{2}K_{d,1}\nabla \varphi(x)\cdot e.
\end{align}
\noindent By substituting the two relations \eqref{eq:K1d} and \eqref{eq:K1d-p} in \eqref{eq:nonlocal-liminf}, using the dominate convergence theorem one readily ends up with the desired estimate
\begin{align*}
\left|\int_{\Omega} u(x) \nabla \varphi (x)\cdot e ~\mathrm{d}x \right| \leq \frac{A_p^{1/p}}{ K_{d,1}}\|\varphi\|_{L^{p'}(\Omega)}.
\end{align*} 
\end{proof}

\bigskip

 \noindent Next we proceed with the following little observation which can be easily derived from the foregoing steps. 
\begin{theorem}\label{thm:limpsup}
	Let $\Omega$ be a $W^{1,p}$-extension set with $1\leq p <\infty$. Assume the family $(\nu_\varepsilon)_\varepsilon$ satisfies the conditions \eqref{eq:normalized-inetrals} and \eqref{eq:concentration-property}. There is a constant $C>0$ only depending on $\Omega,~ p$ and $d$ such that for all $u\in W^{1,p}(\Omega)$,
	\begin{align*}
	\limsup_{\varepsilon\to 0} \iil_{\Omega \Omega } |u(x)-u(y)|^p\nu_\varepsilon(x-y)\,\d x\,\d y \leq C \|u\|_{W^{1,p}(\Omega)}.
	\end{align*}
If $p=1$ and $\Omega$ is a $BV$-extension set then for all $u\in BV(\Omega)$,
	\begin{align*}
	\limsup_{\varepsilon\to 0} \iil_{\Omega \Omega } |u(x)-u(y)|\nu_\varepsilon(x-y)\,\d x\,\d y \leq C\|u\|_{BV(\Omega)}.
	\end{align*}
\end{theorem}

\bigskip

\begin{proof}
	The claims blatantly follow from Lemma \ref{lem:boundedness-limsup} and Lemma \ref{lem:boundedness-limsup-1}.. 
\end{proof}
\bigskip

\noindent Conversely to Theorem \ref{thm:limpsup}, the following result helps to characterize functions in $W^{1,p}(\Omega)$.

\begin{theorem}\label{thm:liminf-nu} Assume that the family $(\nu_\varepsilon)_\varepsilon$ fulfills the conditions \eqref{eq:normalized-inetrals} and \eqref{eq:concentration-property}. Assume $\Omega$ is an open set of $\mathbb{R}^d$ and let $u\in L^p(\Omega)$ with $1<p<\infty$ such that 
	\begin{align*}
	A_p:=	\liminf_{\varepsilon\to 0} \iil_{\Omega \Omega } |u(x)-u(y)|^p\nu_\varepsilon(x-y)\,\d x\,\d y <\infty. 
	\end{align*}
	
	\noindent Then $u\in W^{1,p}(\Omega)$ and the following estimate holds
	\begin{align}\label{eq:estim-weak-derivatve}
	\int_\Omega |\nabla u(x)|^p\,\d x\leq d^2 \frac{A_p^{1/p}}{K_{d,1}}.
	\end{align}
\end{theorem}

\medskip

\begin{proof}
An obvious observation is that \eqref{eq:cont-weak-derivatve} holds true for all $\varphi\in C_c^\infty (\Omega)$, all $1\leq p<\infty$ and all $e\in \mathbb{S}^{d-1}$. Now we assume $e=e_i,~~i=1,\cdots,d$ so that $\nabla \varphi (x)\cdot e_i = \partial_{x_i}\varphi (x)$. In virtue of the density of $C_c^\infty (\Omega)$ in $L^{p'}(\Omega)$, it readily follows from \eqref{eq:cont-weak-derivatve} that for each $i=1,\cdots,d$ the linear mapping
\begin{align*}
\varphi\mapsto \int_{\Omega} u(x) \partial_{x_i}\varphi (x)~\mathrm{d}x 
\end{align*}
uniquely extends as a linear and continuous functional on $L^{p'}(\Omega)$. Note that  $1<p'<\infty$ and hence referring to \cite[Theorem 4.11 \& 4.14]{Bre10}, the Riesz representation\footnote{ The result infers that if $1\leq p<\infty$ then $L^{p'}(\Omega)$ is isomorphic to the dual of $L^p(\Omega)$ and any linear continuous functional on $L^p(\Omega)$ is of the form $u\mapsto\int_\Omega g u$ with $g\in L^{p'}\Omega)$. } for Lebesgue spaces reveals that there exists a unique $g_i\in L^p(\Omega)$ such that
%
%

	\begin{align*}
\int_{\Omega} u(x) \partial_{x_i}\varphi (x)~\mathrm{d}x = \int_{\Omega} g_i(x) \varphi (x)~\mathrm{d}x=- \int_{\Omega} \partial_{x_i} u(x) \varphi (x)~\mathrm{d}x \quad\text{for all} \quad \varphi\in C_c^\infty(\Omega), 
\end{align*}

\noindent where we  let $\partial_{x_i} u= -g_i$. In order words, we have $u\in W^{1,p}(\Omega)$. Furthermore, with the aid of \eqref{eq:cont-weak-derivatve} and the  $L^p$-duality,  we obtain the estimate \eqref{eq:estim-weak-derivatve} as follows
\begin{align*}
\|\nabla u\|_{L^p(\Omega)} \leq \sqrt{d}\sum_{i=1}^d \|\partial_{x_i} u\|_{L^p(\Omega)}
&=\sqrt{d}\sum_{i=1}^d  \sup_{\|\varphi\|_{L^{p'}(\Omega)}=1} \Big|\int_{\Omega} u(x) \, \nabla\varphi (x)\cdot e_i\d x\Big|\leq d^{2}\frac{A_p^{1/p}}{ K_{d,1}}.
\end{align*}
\end{proof}

\bigskip

\noindent Note that the estimate \eqref{eq:estim-weak-derivatve} is not  better than the estimate   $\|\nabla u\|^p_{L^p(\Omega)}\leq \frac{A_p}{K_{d,p}}$ of Theorem \ref{thm:liminf-BBM}. Indeed,  Jensen's inequality implies $K_{d,1}^p d^{-2p}\leq K_{d,1}^p\leq K_{d,p}$ because 
$$\Big(\fint_{\mathbb{S}^{d-1}} |w\cdot e|\d \sigma_{d-1}(w)\Big)^p\leq \fint_{\mathbb{S}^{d-1}} |w\cdot e|^p\d \sigma_{d-1}(w).$$ 
The counterpart of Theorem \ref{thm:liminf-nu} for the case $p=1$  rather falls into the class of functions with bounded variations $BV(\Omega)$(see Theorem \ref{thm:liminf-nu-1}). As we shall see later in the proof, because of the lack of reflexivity of $L^1(\Omega)$, assuming that $A_1<\infty$ is not enough to conclude that $u\in W^{1,1}(\Omega)$.
\bigskip

\noindent Next we resume the  characterization  of $W^{1,p}(\Omega)$ with $1<p<\infty$ as follows. 
\begin{theorem}\label{thm:charact-w1p}
	Let $\Omega\subset \mathbb{R}^d$ be a $W^{1,p}$-extension set with $1<p<\infty$. 
	Assume that the family $(\nu_\varepsilon)_\varepsilon$ fulfills the conditions \eqref{eq:normalized-inetrals} and \eqref{eq:concentration-property}. Let $u\in L^p(\Omega)$. 
	Then $u\in W^{1,p}(\Omega)$ if and only if 
	\begin{align*}
	\liminf_{\varepsilon\to 0}\iint\limits_{\Omega\Omega}|u(x)-u(y)|^p\nu_\varepsilon(x-y) \mathrm{d}y\mathrm{d}x <\infty. 
	\end{align*} 
	\noindent Moreover, with the convention that $\int_{\Omega}|\nabla u(x)|^p\mathrm{d}x = \infty$ when $|\nabla u|$ is not in $L^p(\Omega)$ we have 
	\begin{align*}
	\lim_{\varepsilon\to 0}\iint\limits_{\Omega\Omega}|u(x)-u(y)|^p\nu_\varepsilon(x-y) \mathrm{d}y\mathrm{d}x = K_{d,p} \int_{\Omega}|\nabla u(x)|^p\mathrm{d}x. 
	\end{align*}
\end{theorem}
\medskip

\begin{proof}	
	If $u\in W^{1,p}(\Omega)$ then Theorem \ref{thm:BBM-result} asserts that 
	\begin{align*}
	\limsup_{\varepsilon\to 0}\iint\limits_{\Omega\Omega}|u(x)-u(y)|^p\nu_\varepsilon(x-y) \mathrm{d}y\mathrm{d}x = K_{d,p} \int_{\Omega}|\nabla u(x)|^p\mathrm{d}x<\infty. 
	\end{align*}
	%
	The converse is a direct consequence of Theorem \ref{thm:liminf-nu}. Now assume $|\nabla u|$ does not belong to $L^p(\Omega)$ then once more Theorem \ref{thm:liminf-nu}. implies that 
	\begin{align*}
	\liminf_{\varepsilon\to 0}\iint\limits_{\Omega\Omega}|u(x)-u(y)|^p\nu_\varepsilon(x-y) \mathrm{d}y\mathrm{d}x =\infty,
	\end{align*}
	in other words we have 
	\begin{align*}
	\lim_{\varepsilon\to 0}\iint\limits_{\Omega\Omega}|u(x)-u(y)|^p\nu_\varepsilon(x-y) \mathrm{d}y\mathrm{d}x =\int_{\Omega}|\nabla u(x)|^p\mathrm{d}x = \infty. 
	\end{align*}
\end{proof}

\noindent Let us now mention some concrete examples. 
\begin{corollary}\label{cor:BBM-fractional}
	Assume $\Omega\subset \R^d$ is an extension domain. Then the following 
\begin{align*}
&\lim_{s\to 1} (1-s)\iint\limits_{\Omega\Omega}\frac{|u(x)-u(y)|^p}{|x-y|^{d+sp}}\mathrm{d}y\mathrm{d}x =\lim_{\varepsilon \to 0} \varepsilon\iint\limits_{\Omega\Omega}\frac{|u(x)-u(y)|^p}{|x-y|^{d+p(1-\varepsilon)}} \mathrm{d}y\mathrm{d}x =\frac{|\mathbb{S}^{d-1}|}{p}K_{d,p} \int_{\Omega}|\nabla u(x)|^p\mathrm{d}x \\
\\
&\lim_{\varepsilon \to 0} \varepsilon^{-d}\iint\limits_{ \Omega\times\Omega\cap \{|x-y|<\varepsilon\} }|u(x)-u(y)|^p \mathrm{d}y\mathrm{d}x = \frac{|\mathbb{S}^{d-1}|}{d}K_{d,p}\int_{\Omega}|\nabla u(x)|^p\mathrm{d}x \\\\
& 	
\lim_{\varepsilon \to 0} \frac{1}{|\log \varepsilon|}\iint\limits_{\Omega\times\Omega\cap \{|x-y|>\varepsilon\}} \frac{|u(x)-u(y)|^p}{|x-y|^{d+p}}\mathrm{d}y\mathrm{d}x = |\mathbb{S}^{d-1}|K_{d,p} \int_{\Omega}|\nabla u(x)|^p\mathrm{d}x .
\end{align*}
\end{corollary}

\begin{proof}
For the first relation, take $\nu_\varepsilon(h) = a_{\varepsilon, d,p} |h|^{-d-p(1-\varepsilon)}$ with $ a_{\varepsilon, d,p} = \tfrac{p\varepsilon(1-\varepsilon)}{|\mathbb{S}^{d-1}|}$. For the second take $\nu_\varepsilon(h) = \tfrac{d}{|\mathbb{S}^{d-1}|} |h|^{-p}\mathds{1}_{B_\varepsilon}(h)$. Last, take $\nu_\varepsilon(h) = \tfrac{1}{|\mathbb{S}^{d-1}||\log\varepsilon|}|h|^{-d-p} \mathds{1}_{B_1\setminus B_\varepsilon}(h)$ or more generally, for fixed $\varepsilon_0\geq 1$, take \\ $\nu_\varepsilon(h) = \tfrac{ b_\varepsilon}{|\mathbb{S}^{d-1}||\log\varepsilon|}\mathds{1}_{B_{\varepsilon_0}\setminus B_\varepsilon}(h)$ with $b_\varepsilon = \tfrac{|\log\varepsilon|}{ (1-\varepsilon_0^{-p})/p+ |\log\varepsilon|}\to 1$.
\end{proof}

\vspace{2mm}

\noindent A noteworthy application of this section is given by the following result. It infers that the spaces $(W_{\nu_\eps}^{p}(\Omega|\R^d))_\eps$ collapse to the space $W^{1,p}(\Omega)$ as $\eps\to0$.
 
\begin{theorem}[Collapsing convergence across the boundary ]\label{thm:collapsing-accross-boundary}
	Let  $\Omega\subset \R^d$ be open with a compact Lipschitz boundary.  The following are true for every $u\in W^{1,p}(\R^d)$
	\begin{align*}
	\lim_{\eps\to 0}\iil_{ \Omega \Omega^c} |u(x)-u(y)|^p\nu_\eps(x-y)\d x\d y= 0. 
	\end{align*}
	Consequently, 
	\begin{align*}
	\lim_{\eps\to 0}\iil_{ \Omega \R^d} |u(x)-u(y)|^p\nu_\eps(x-y)\d x\d y= \lim_{\eps\to 0}\iil_{ (\Omega^c \times\Omega^c)^c} |u(x)-u(y)|^p\nu_\eps(x-y)\d x\d y= K_{d,p}\int_\Omega|\nabla u(x)|^p\ dx. 
	\end{align*}
\end{theorem}

\begin{proof}
Note that both $\Omega$ and $\Omega^c$ have the same compact Lipschitz boundary and are thus $W^{1,p}$-extension domains. Whence  the claims follow since for $u\in W^{1,p}(\R^d)$ we get
\begin{align*}
&\lim_{\eps\to 0}\iil_{\R^d \R^d} |u(x)-u(y)|^p\nu_\eps(x-y)\d x\d y= K_{d,p}\int_{ \R^d} |\nabla u(x)|^p\d x,\\ 
&\lim_{\eps\to 0}\iil_{ \Omega \Omega} |u(x)-u(y)|^p\nu_\eps(x-y)\d x\d y = K_{d,p}\int_{ \Omega} |\nabla u(x)|^p\d x,\\
&\lim_{\eps\to 0}\iil_{ \Omega^c \Omega^c} |u(x)-u(y)|^p\nu_\eps(x-y)\d x\d y= K_{d,p}\int_{ \Omega^c} |\nabla u(x)|^p\d x. 
\end{align*} 
Thus it suffices to observe that $\Omega\times \Omega^c\cup \Omega^c\times \Omega=  (\R^d\times \R^d)\setminus[\Omega\times \Omega \cup  \Omega^c\times \Omega^c]$,   $\Omega\times \R^d= \Omega\times \Omega\cup \Omega\times \Omega^c$ and $ (\Omega^c\times \Omega^c)^c = (\R^d\times \R^d)\setminus(\Omega^c\times \Omega^c)$.


\end{proof}

\section{Characterization of the spaces of bounded variation}

\noindent A natural question raised by the authors in \cite{BBM01} was to know if Theorem \ref{thm:BBM-result} persists for functions in $BV(\Omega)$. They were able to obtain a positive answer in one dimension with $\Omega= (0,1)$. To be more precise, recalling Example \ref{Ex:rho-var} they showed that 
$$\lim_{\varepsilon\to 0 }\int_0^1\int_0^1\frac{|u(x)-u(y)|}{|x-y|}\rho_\varepsilon(x-y)\,\d y \,\d x=\ |u|_{BV(\Omega)}.$$
The case $d\geq 2$ was thoroughly completed later for a bounded Lipschitz domain by \cite{Dav02}. However, our proof of this result as stated in the below (see Theorem \ref{thm:limit-BV}) is much simpler. Let us start with a variant of Theorem \ref{thm:liminf-BBM} for the case $p=1$.

\medskip

\begin{theorem}\label{thm:liminf-BV} Let $\Omega\subset \mathbb{R}^d$ be any open set (not necessarily bounded). Assume that for $p=1$ the family $(\nu_\varepsilon)_\varepsilon$ fulfills the conditions \eqref{eq:normalized-inetrals} and \eqref{eq:concentration-property}. Adopting the convention that $|\nabla u|_{BV(\Omega)}=\infty$ if the Radon measure $|\nabla u|$ does not have finite total variation, then for all $u\in L^1(\Omega)$ we have 
	\begin{align*}
	K_{d,1} | \nabla u|_{BV(\Omega)}\leq 	\liminf_{\varepsilon\to 0}\iil_{\Omega\Omega}|u(x)-u(y)|\nu_\varepsilon(x-y) \mathrm{d}y \,\d x .
	\end{align*}
	
\end{theorem}

\medskip

\begin{proof}
We extend the proof of Theorem \ref{thm:liminf-BBM} wherein we recall $u^\delta = u*\phi_\delta \in C^\infty(\R^d)$, $\phi_\delta)(x)=\delta^{-d} \phi(\frac{x}{\delta})$ is a mollification sequence with $\phi\in C_c^\infty(\R^d)$ and $u$ is extended off $\Omega$ by $0$. The relation \eqref{eq:limit-approx} asserts that for all $\delta>0,$
	\begin{align*}
K_{d,1}\int_{\Omega_\delta} |\nabla u^\delta (x)| \d x\leq \liminf_{\varepsilon\to 0} \iil_{\Omega\Omega}|u(x) -u(y)| \nu_\varepsilon(x-y)\d x\d y.
\end{align*}
Hence the estimate we are seeking  occurs once we show that 
\begin{align*}
K_{d,1} |u|_{BV(\Omega)}\leq 
\liminf_{\delta\to 0}\int_{\Omega_\delta} |\nabla u^\delta (x)| \d x.
\end{align*}
Let $\chi\in C_c^\infty(\Omega, \R^d)$ such that $\|\chi\|_\infty\leq 1$,  we find that 
\begin{align*}
\Big| \int_{\Omega} u(x) \operatorname{div} \chi(x)\d x &- \int_{\Omega_\delta} u^\delta (x) \operatorname{div} \chi(x)\d x\Big|\\
&=\Big| \int_{\Omega_\delta} (u (x) -u*\phi_\delta (x)) \operatorname{div} \chi(x)\d x + \int_{\Omega\setminus \Omega_\delta} \hspace{-3ex} u*\phi_\delta(x) \operatorname{div} \chi(x)\d x \Big|\\
&\leq \|u*\phi_\delta-u\|_{L^1(\Omega)} +\|\operatorname{div} \chi\|_\infty \|\phi\|_\infty\int_{\Omega\setminus \Omega_\delta} |u (x) |\d x\xrightarrow{\delta\to 0}0. 
\end{align*}
This implies the following by using the fact that $u$ is a distribution on $\Omega$. 
\begin{align*}
\int_{\Omega} u(x) \operatorname{div} \chi(x)\d x&=\lim_{\delta\to 0} \int_{\Omega_\delta} u^\delta (x) \operatorname{div} \chi(x)\d x\\
&= \lim_{\delta\to 0} -\int_{\Omega} \nabla u^\delta (x) \cdot \chi(x)\d x- \lim_{\delta\to 0} \int_{\Omega\setminus \Omega_\delta} \hspace{-3ex} u*\phi_\delta(x) \operatorname{div} \chi(x) \d x\\
&= \lim_{\delta\to 0} -\int_{\Omega_\delta} \nabla u^\delta (x) \cdot \chi(x)\d x - \lim_{\delta\to 0} \int_{\Omega\setminus \Omega_\delta} \hspace{-3ex} u*\nabla \phi_\delta (x) \cdot \chi(x)+ u*\phi_\delta(x) \operatorname{div} \chi(x) \d x\\
&\leq \liminf_{\delta\to 0}\int_{\Omega_\delta} |\nabla u^\delta (x)| \d x +C \lim_{\delta\to 0}\int_{\Omega\setminus \Omega_\delta} |u (x) |\d x\\
&= \liminf_{\delta\to 0}\int_{\Omega_\delta} |\nabla u^\delta (x)| \d x 
\end{align*}
with $C=\big(\|\chi\|_\infty \|\nabla \phi\|_\infty+\|\operatorname{div} \chi\|_\infty \|\phi\|_\infty \big)$. 
Therefore, since the above holds for arbitrarily chosen $\chi\in C_c^\infty(\Omega,\R^d)$ such that $\|\chi\|_\infty\leq 1$, by definition of $BV(\Omega)$ (see the relation \eqref{eq:bounded-variation}) we get
\begin{align*}
 K_{d,1} | u|_{BV(\Omega)}\leq
 \liminf_{\delta\to 0}\int_{\Omega_\delta} |\nabla u^\delta (x)| \d x 
\end{align*}
which completes the proof.
\end{proof}

\medskip

\noindent Now we establish the limsup inequality which is a counterpart of Theorem \ref{thm:liminf-BV} and extends Theorem \ref{thm:limsup_BBM}.
\begin{theorem}\label{thm:limitsup-BV}
	Assume $\Omega\subset \mathbb{R}^d$ is an open $BV$-extension set (not necessarily bounded). Assume that for $p=1$ the family $(\nu_\varepsilon)_\varepsilon$ fulfills the conditions \eqref{eq:normalized-inetrals} and \eqref{eq:concentration-property} then for all $u\in BV(\Omega)$ 
	\begin{align*}
	\limsup_{\varepsilon\to 0}\iil_{\Omega\Omega}|u(x)-u(y)|\nu_\varepsilon(x-y) \mathrm{d}y \,\d x\leq K_{d,1} |u|_{BV(\Omega)} . 
	\end{align*}
\end{theorem}

\medskip

\begin{proof}
The case where $|\nabla u|_{BV(\Omega)} =\infty$ is trivial. Now for $u\in BV(\Omega)$, we let $\overline{u}\in BV(\R^d)$ be its extension to $\R^d$. Consider 
$\Omega(\delta) = \Omega+ B_\delta(0) = \{x\in\R^d~:\operatorname{dist}(x, \Omega)<\delta \}$ with $0<\delta<1$ sufficiently small. We know from Theorem \ref{thm:evans-grapiepy} that there is $(u_n)_n$ a sequence of functions of $C^\infty(\R^d)\cap W^{1,1}(\R^d)$ converging to $\overline{u}$ in $L^1(\R^d)$ and such that $\|\nabla u_n\|_{L^1(\R^d)}\xrightarrow{n \to \infty}|\overline{u}|_{BV(\R^d)}$. For each $n\geq 1$, by the estimate \eqref{eq:split-estimate} we have 
\begin{align*}
\iil_{\Omega\Omega}|u_n(x)-u_n(y)|\nu_\varepsilon(x-y) \mathrm{d}y \d x\leq K_{d,1}\int_{\Omega(\delta) }|\nabla u_n(x)|\d x+ 2\|u_n\|_{L^1(\Omega)}\int_{|h|\geq \delta}\nu_\varepsilon(h)\,\d h.
\end{align*}
The Fatou lemma implies
	\begin{align*}
\iil_{\Omega\Omega}|u(x)-u(y)|\nu_\varepsilon(x-y) \mathrm{d}y \d x
	&\leq \liminf_{n\to\infty}\iil_{\Omega\Omega}|u_n(x)-u_n(y)|\nu_\varepsilon(x-y) \mathrm{d}y \d x \\
	&\leq \lim_{n\to\infty} K_{d,1} \int_{\Omega(\delta)}|\nabla u_n(x)|\,\d x + 2\|u_n\|_{L^1(\Omega)}\int_{|h|\geq \delta}\nu_\varepsilon(h)\,\d h.\\
	&= K_{d,1} |\overline{u}|_{BV(\Omega(\delta))} + 2\|u\|_{L^1(\Omega)}\int_{|h|\geq \delta}\nu_\varepsilon(h)\,\d h.
	\end{align*}
Correspondingly, we have 
	\begin{align}\label{eq:split-estimate-BV}
\iil_{\Omega\Omega}|u(x)-u(y)|\nu_\varepsilon(x-y) \mathrm{d}y \d x
\leq K_{d,1} | \overline{u}|_{BV(\Omega(\delta))} + 2\|u\|_{L^1(\Omega)}\int_{|h|\geq \delta}\nu_\varepsilon(h)\,\d h.
\end{align}
The claimed estimate follows by letting $\varepsilon\to 0$ and $\delta\to 0$ successively. Indeed, since $\overline{u}\in BV(\R^d) $ and $u= \overline{u}\mid_\Omega$ we have 
$$\int_{|h|\geq \delta}\nu_\varepsilon(h)\,\d h\xrightarrow{\varepsilon\to 0}0\qquad\text{and}\qquad |\overline{u}|_{BV(\Omega(\delta))}\xrightarrow{\delta\to 0} |u|_{BV(\Omega)}.$$

\end{proof}

\bigskip

\begin{theorem}[c.f \cite{Dav02}]\label{thm:limit-BV}
	Assume $\Omega\subset \mathbb{R}^d$ is an open $BV$-extension set(not necessarily bounded). Assume that for $p=1$ the family $(\nu_\varepsilon)_\varepsilon$ fulfills the conditions \eqref{eq:normalized-inetrals} and \eqref{eq:concentration-property} then for all $u\in BV(\Omega)$ 
	\begin{align*}
	\lim_{\varepsilon\to 0}\iil_{\Omega\Omega}|u(x)-u(y)|\nu_\varepsilon(x-y) \mathrm{d}y \,\d x = K_{d,1} |u|_{BV(\Omega)} . 
	\end{align*}
	
\end{theorem}

\medskip

\begin{proof}
The result blatantly follows by combining Theorem \ref{thm:liminf-BV} and Theorem \ref{thm:limitsup-BV} since 
\begin{align*}	
\limsup_{\varepsilon\to 0}\iil_{\Omega\Omega}|u(x)-u(y)|\nu_\varepsilon(x-y) \mathrm{d}y \,\d x\leq K_{d,1} |u|_{BV(\Omega)} \leq \liminf_{\varepsilon\to 0}\iil_{\Omega\Omega}|u(x)-u(y)|\nu_\varepsilon(x-y) \mathrm{d}y \,\d x.
\end{align*}
\end{proof}

\bigskip

\noindent The following theorem is a revisited version of \cite[Lemma 2]{Dav02} which provides an alternative to Theorem \ref{thm:BBM-result} and \ref{thm:limit-BV} if $\Omega$ is not an extension domain..
\begin{theorem}\label{thm:limit-weak-BV}
	Let $\Omega\subset \mathbb{R}^d$ be an open set. Assume that for $p=1$ the family $(\nu_\varepsilon)_\varepsilon$ fulfills the conditions \eqref{eq:normalized-inetrals} and \eqref{eq:concentration-property}. For every $u\in W^{1,p}(\Omega)$ the family of  measures $(\mu_{\varepsilon})_\eps$ with 
	\begin{align*}
	\d \mu_{\varepsilon}(x) = \int_{ \Omega}|u(x)-u(y)|^p\nu_\varepsilon(x-y) \mathrm{d}y\d x. 
	\end{align*}
	\noindent converges weakly on $\Omega$ (in the sense of Radon measures) to the Radon measure  $\d \mu(x)  =K_{d,p}|\nabla u(x)|^p\d x$, 
	i.e. $ \mu_\varepsilon (E)\xrightarrow{\eps\to 0}\mu(E)$ for every compact set $E\subset \Omega$. Moreover, if $p=1$ and $u\in BV(\Omega)$ then $\d \mu(x)  =K_{d,1}\d |\nabla u|(x)$. 
\end{theorem}

\medskip

\begin{proof}
Let $E\subset \Omega$ be compact with a nonempty interior.  Consider the open set $E(\delta) = E+ B_\delta(0)$ where $0<\delta< \operatorname{dist}(\partial\Omega, E)$.  The family of functions $(\mu_\varepsilon)_\varepsilon$ is bounded in $L^1(E)$. Indeed, if we denote $ \d |\nabla u|^p(x)= |\nabla u(x)|^p\d x$ for  $u\in W^{1,p}(\Omega)$, then the estimates \eqref{eq:split-estimate} and \eqref{eq:split-estimate-BV} with $\Omega$ replaced by $E$ imply 
\begin{align}\label{eq:split-estimate1}
\int_E \mu_\varepsilon(x)\,\d x
&\leq K_{d,p}\int_{E(\delta)} \hspace{-2ex}\d |\nabla u|^p(x)+  2^p \|u\|^p_{L^p(\Omega)}\int_{|h|>\delta} \nu_\varepsilon(h)\,\d h,\\
\text{with} \qquad\qquad& \qquad\qquad\int_{|h|>\delta} \nu_\varepsilon(h)\,\d h \leq (1\land\delta^p)^{-1}, \notag
\end{align}

\noindent   In virtue of the  weak compactness of $L^1(E)$, (see  \cite[p.116]{Bre10}) we may assume that $(\mu_\varepsilon)_\varepsilon$ converges in the weak-* sense to a Radon measure $\mu_E$ otherwise, one may pick a converging subsequence, i.e  $\langle\mu_\eps-\mu_E, \varphi\rangle\xrightarrow{\eps\to 0}0$ for all $\varphi\in C(E)$. For a suitable $(\Omega_j)_{j\in \mathbb{N}}$ exhaustion of $\Omega$, i.e. $\Omega_j's$ are open, each $K_j=\overline{\Omega}_j$ is compact, $K_j=\overline{\Omega}_j\subset \Omega_{j+1}$ and $\Omega=\bigcup_{j\in \mathbb{N}}\Omega_j$, it would suffice to let  $\mu= \mu_{K_j}= K_{d,p}|\nabla u|^p$ on $K_j$. We aim to show that $\mu= K_{d,p}|\nabla u|^p$. Noticing $\mu$  and $K_{d,p}|\nabla u|^p$ are Radon measures it sufficient to show that both measures coincide open compact sets, i.e. we have to show that  
\begin{align*}
\mu_E(E)  = K_{d,p}\int_{E} \,\d |\nabla u|^p(x). 
\end{align*} 

\noindent On the one hand,  given that  $\mu_\varepsilon(E)\to \mu(E)$ and $\int_{|h|>\delta}\nu_\varepsilon(h)\,\d h\to0$ as $\eps \to 0, $ the fact that $u\in W^{1,p}(\Omega)$ or $u\in BV(\Omega)$ enables us to successively let $\varepsilon\to 0$ and $\delta\to 0$ in \eqref{eq:split-estimate1} which becomes
\begin{align}\label{eq:absolute-cont}
\int_E \d\mu(x)\,\leq K_{d,p}\int_{E} \,\d |\nabla u|^p(x). 
\end{align}

\noindent On other hand, since $E$ has an nonempty interior,  from Theorem \ref{thm:liminf-BBM} we get 
\begin{align*}
K_{d,p}\int_{E} \,\d |\nabla u|^p(x) &\leq	\liminf_{\varepsilon\to 0} \iil_{EE}|u(x) -u(y)|^p \nu_\varepsilon(x-y)\d x\d y
\leq  \lim_{\varepsilon\to 0} \int_{E}\mu_\eps(x)\d x = \int_{E}\d\mu_E(x). 
\end{align*}
\noindent Therefore, we get $\d \mu =K_{d,p} \d |\nabla u|^p$ as claimed.

	\noindent Let us provide an alternative proof to the case $p=1$ and $u\in BV(\Omega)$.  Exploiting Lemma \ref{lem:boun-integration-bypart} together with the relations \eqref{eq:K1d} and \eqref{eq:K1d-p} it clearly appears that for $\varphi\in C_c^\infty(\Omega)$ with $\varphi\geq 0$ and $e\in \mathbb{S}^{d-1}$ the following occurs
	\begin{align*}
	&K_{d,1}\Big|\int_\Omega \varphi(x) \,\d (\nabla u(x)\cdot e) \Big| 
	\leq \liminf_{\varepsilon \to 0} \iil_{\R^d \R^d} |u(x)-u(y)|\varphi(x)\nu_\varepsilon(x-y)\,\d y\,\d x\\
	&= \liminf_{\varepsilon \to 0} \int_{\Omega} \varphi(x)\,\d x \int_{\Omega} |u(x)-u(y)|\nu_\varepsilon(x-y)\,\d y+\liminf_{\varepsilon \to 0} \iil_{\Omega\mathbb{R}^d\setminus \Omega} |u(x)-u(y)| \varphi(x) \nu_\varepsilon(x-y)\,\d y\,\d x\\
	&\leq \liminf_{\varepsilon \to 0} \int_{\Omega} \mu_\varepsilon(x)\varphi(x)\,\d x + 2\|\varphi\|_{L^\infty(\R^d)}\|u\|_{L^1(\Omega)}\lim_{\varepsilon \to 0} \il_{|h|\geq \delta}\nu_\varepsilon(h)\,\d h =\int_{ \Omega}\varphi(x)\,\d \mu(x).
	\end{align*}
	Provided that $ \mu_\varepsilon \rightharpoonup^*\mu $ and the vector measure $\nabla u$ is the weak derivative of $u$ and where we put $\delta= \operatorname{dist}(\supp \varphi, \partial\Omega)>0$. This proves that for every $ \varphi\in C_c^\infty(\Omega),$ $\varphi\geq 0$ and $e\in \mathbb{S}^{d-1}$.
	\begin{align*}
	K_{d,1}\Big|\int_\Omega \varphi(x) \,\d (\nabla u(x)\cdot e) \Big| 
	\leq \int_{ \Omega}\varphi(x)\,\d \mu(x).
	\end{align*}
	%
	In a sense, this implies\footnote{Technically since we are dealing with Radon measures, it is sufficient to show the due inequality for all compacts subset of $\Omega$. The latter statement can be accomplished by choosing appropriate cut-off functions $\varphi\in C_c^\infty(\Omega)$.} that for all Borel set $B\subset \Omega$ we have 
	\begin{align*}
	K_{d,1}(\nabla u\cdot e) (B))\leq \mu(B), \quad\text{for every }\quad e\in \mathbb{S}^{d-1}.
	\end{align*}
	\noindent Choosing especially the unit vector $e= \nabla u(B)/|\nabla u(B)|$ (vector measure of $B$ ) so that 
	$(\nabla u\cdot e) (B) = |\nabla u| (B)$ we get that 
	\begin{align}\label{eq:absolute-continous-bis}
	K_{d,1}|\nabla u|(B)\leq \mu(B).
	\end{align}
	\bigskip Recalling \eqref{eq:absolute-cont} we end up showing that $ K_{d,1}|\nabla u|= \mu $ as claimed. 
	%

\end{proof}

\bigskip

\noindent The following theorem is the appropriate variant of Theorem \ref{thm:liminf-nu} in the case $p = 1$. 
\begin{theorem}\label{thm:liminf-nu-1} Assume that the family $(\nu_\varepsilon)_\varepsilon$ fulfills the conditions \eqref{eq:normalized-inetrals} and \eqref{eq:concentration-property}. Assume $\Omega$ is an open set of $\mathbb{R}^d$ and let $u\in L^1(\Omega)$ such that 
	\begin{align*}
	A_1:=\liminf_{\varepsilon\to 0} \iil_{\Omega \Omega } |u(x)-u(y)|\nu_\varepsilon(x-y)\,\d x\,\d y <\infty. 
	\end{align*}
	Then $u\in BV(\Omega)$ and the following estimate holds
	\begin{align}\label{eq:estim-weak-derivatve-1}
	|u|_{BV(\Omega)}\leq \frac{A_1}{ K_{d,1}}.
	\end{align}
\end{theorem} 

\medskip

\begin{proof}
Let $\chi =(\chi_1,\chi_2,\cdots,\chi_d)\in C_c^\infty(\Omega, \R^d)$ such that $\|\chi\|_{L^\infty(\Omega)}\leq 1$ and $e=e_i,\,\,i=1,2\cdots, d$. Since $\nabla \chi_i \cdot e_i= \partial_{x_i}\chi_i$ with $\chi_i\in C_c^\infty(\Omega)$,  the estimate \eqref{eq:cont-weak-derivatve} implies  
\begin{align*}
\Big|\int_{\Omega} u(x) \, \operatorname{div}\chi \d x\Big| 
&=\Big|\sum_{i=1}^d\int_{\Omega} u(x) \, \nabla\chi_i(x)\cdot e_i\d x\Big|\leq d\frac{A_1}{K_{d,1}}. 
\end{align*} 
\noindent Hence $u\in BV(\Omega)$ and we have $|u|_{BV(\Omega)}\leq d\frac{A_1}{K_{d,1}}$

\end{proof}
\bigskip

 \noindent In the same spirit, a variant of Theorem \ref{thm:charact-w1p} with $p=1$, characterizing the space $BV(\Omega)$ is given as follows. 
\begin{theorem}\label{thm:charact-BV}
	Let $\Omega\subset \mathbb{R}^d$ be a $BV$-extension domain. Assume that for $p=1$ the family $(\nu_\varepsilon)_\varepsilon$ fulfills the conditions \eqref{eq:normalized-inetrals} and \eqref{eq:concentration-property}. Let $u\in L^1(\Omega)$. Then $u\in BV(\Omega)$ if and only if 
	\begin{align*}
	\liminf_{\varepsilon\to 0}\iint\limits_{\Omega\Omega}|u(x)-u(y)|\nu_\varepsilon(x-y) \mathrm{d}y\mathrm{d}x <\infty. 
	\end{align*} 
	\noindent Moreover, with the convention that $|u|_{B(\Omega)}= \infty$ when $|\nabla u|$ is not of finite total variation we have 
	\begin{align*}
	\lim_{\varepsilon\to 0}\iint\limits_{\Omega\Omega}|u(x)-u(y)|\nu_\varepsilon(x-y) \mathrm{d}y\mathrm{d}x = K_{d,1} |u|_{BV(\Omega)}. 
	\end{align*}
	\noindent 
\end{theorem}
\bigskip

\begin{proof}	
	If $u\in BV(\Omega)$ then Theorem \ref{thm:limit-BV} asserts that 
	\begin{align*}
	\liminf_{\varepsilon\to 0}\iint\limits_{\Omega\Omega}|u(x)-u(y)|\nu_\varepsilon(x-y) \mathrm{d}y\mathrm{d}x = K_{d,1} \int_{\Omega}|\nabla u(x)|\mathrm{d}x<\infty. 
	\end{align*}
	%
	The converse is a direct consequence of Theorem \ref{thm:liminf-nu-1}. Now assume $|\nabla u|$ does not have finite total variation then once more Theorem \ref{thm:liminf-nu-1} implies that 
	\begin{align*}
	\liminf_{\varepsilon\to 0}\iint\limits_{\Omega\Omega}|u(x)-u(y)|\nu_\varepsilon(x-y) \mathrm{d}y\mathrm{d}x =\infty
	\end{align*}
	in other words we have 
	\begin{align*}
	\lim_{\varepsilon\to 0}\iint\limits_{\Omega\Omega}|u(x)-u(y)|\nu_\varepsilon(x-y) \mathrm{d}y\mathrm{d}x =|u|_{BV(\Omega)}= \infty. 
	\end{align*}
\end{proof}


\section{Asymptotically compactness}

In the foregoing sections we dealt with convergence of seminorms when the associated functions are  independent of the parameter $\eps>0$. Throughout this section we examine the case where the functions also depend on $\eps>0$. Therefore it is expected that we will encounter a situation where sequential compactness arguments are needed. In order to alleviate the notations we shall replace $\nu_\eps$ with $\nu_{\eps_n}$ where $(\eps_n)_n$ is any sequence of positive real numbers tending to 0 as $n \to \infty$.  We recall that our standing assumptions on $\nu_{\eps_n}'s$ is as follows: $1\leq p<\infty$, 
\begin{align}
& \nu_{\eps_n}~~\text{is radial}\qquad \text{and} \qquad\int_{\mathbb{R}^d} (1\land|h|^p) \nu_{\eps_n}(h)\,\d h=1\label{eq:normalized-inetrals-n}\\
& \lim_{n\to \infty}\il_{|h|\geq \delta} (1\land|h|^p) \nu_{\eps_n} (h)\,\d h=0 \quad\text{for all }\quad \delta >0 .\label{eq:concentration-property-n}
\end{align}
The authors in \cite{Ponce2004,BBM01} solely deal with the restrictive case where $\nu_{\eps_n}(h) = \mathds{1}_{B_1}(h)|h|^{-p}\rho_{\eps_n}(h)$ where $\rho_{\eps_n}$'s satisfy the following condition
\begin{align*}
\rho_{\eps_n}~~\text{is radial}\quad \text{and} \quad\int_{\mathbb{R}^d} \rho_{\eps_n}(h)\,\d h=1\quad \text{and} \quad \lim_{n\to \infty}\il_{|h|\geq \delta} \rho_{\eps_n} (h)\,\d h=0 \quad\text{for all }\quad \delta >0.
\end{align*}
In the present section we upgrade this within a fairly  large class of family $(\nu_{\eps_n})_n$ under consideration.

\medskip

\noindent Next,  we start by borrowing some results from \cite{Ponce2004}. 
\begin{lemma}\label{lem:G-estimate} For $\, 0<s<t$, write $t= ks+ \theta s$ with $k = \lfloor\frac{t}{s}\rfloor\in \mathbb{N}$ and $\theta\in [0,1)$. For $g\in L^p(\R)$ we have 
\begin{align*}
\frac{G(t)}{t^p} \leq 2^{p-1}\Big( \frac{G(s)}{s^p}+ \frac{G(\theta s)}{t^p} \Big),\quad\quad\text{with}\quad
	G(t) =\int_{\R}|g(\tau+t) -g(\tau)|^p\d \tau. 
	\end{align*}
\end{lemma}

\vspace{2mm}

\begin{proof} Applying Jensen's inequality, we obtain 
	\begin{align*}
|g(\tau+t) -g(\tau)|^p &\leq 2^{p-1} |g(\tau+ks) -g(\tau)|^p+ 2^{p-1}  |g(\tau+ks) -g(\tau+ks+\theta s)|^p\\
&= 2^{p-1} \Big|\sum_{j=0}^{k-1}g(\tau+js+s) -g(\tau+js)\Big|^p+2^{p-1}  |g(\tau+ks) -g(\tau+ks+\theta s)|^p\\
&\leq 2^{p-1}k^{p-1} \sum_{j=0}^{k-1}|g(\tau+js+s) -g(\tau+js)|^p+2^{p-1}  |g(\tau+ks) -g(\tau+ks+\theta s)|^p.
	\end{align*}
\noindent Integrating both sides and noticing that $k\leq t/s$ the result follows since
	  \begin{align*}
	  G(t) =2^{p-1}k^{p} G(s)  + 2^{p-1}G(\theta s) \leq  2^{p-1}(t/s)^{p} G(s)  + 2^{p-1}G(\theta s). 
	  \end{align*}
\end{proof}

\noindent Let us derive the following integration formula on the unit sphere. 
\begin{theorem}\label{thm:integration-on-sphere}
Assume $ d\geq 2$. Let $h: \mathbb{S}^{d-1}\to [0, \infty]$ be a Borel measurable function. For $v\in  \mathbb{S}^{d-1}$ we consider that $v$-section of $\mathbb{S}^{d-1}$ denoted by $\mathbb{S}^{d-2}_v= \mathbb{S}^{d-1}\cap (\R^d_v)^\perp$ with $\R^d_v= \{sv:\, s\in \R\}$. Then 
\begin{align}
\int_{\mathbb{S}^{d-1}}h(w) \d \sigma_{d-1}(w) =  \int_{\mathbb{S}^{d-1}} \Big( \int_{\mathbb{S}^{d-2}_v} h(w) \d \sigma_{d-2}(w) \Big)\frac{ \d \sigma_{d-1}(v)}{|\mathbb{S}^{d-2}|}. 
\end{align}
\end{theorem}
\begin{proof}
By standards approximation of positive measurable functions, it is sufficient to prove the claim solely for indicator functions. Let $A\subset\mathbb{S}^{d-1} $ be a Borel measurable set. Define the measure 
\begin{align*}
\mu(A)=  \int_{\mathbb{S}^{d-1}} \sigma_{d-2}(A\cap  \mathbb{S}^{d-2}_v ) \frac{ \d \sigma_{d-1}(v)}{|\mathbb{S}^{d-2}|}.  
\end{align*}
Given that the Borel(Hausdorff) measures $\sigma_{d-1}$ and $\sigma_{d-2}$ are respectively rotation invariant on $\mathbb{S}^{d-1}$ and $\mathbb{S}^{d-2}$  so does the measure $\mu$ on $\mathbb{S}^{d-1}$. Indeed, for each $v\in \mathbb{S}^{d-1}$ and $O\in \mathcal{O}(d)$, a rotation on $\mathbb{R}^d$, one has that $O(A) \cap  \mathbb{S}^{d-2}_v ) = O(A \cap  \mathbb{S}^{d-2}_{O^\top(v)})$ since $O^\top O= I_d$. The rotation invariance of $\sigma_{d-2}$ implies 
\begin{align*}
\sigma_{d-2}(O(A)\cap  \mathbb{S}^{d-2}_v ) = \sigma_{d-2}(O(A\cap  \mathbb{S}^{d-2}_{O^\top(v)})) = \sigma_{d-2}(A\cap  \mathbb{S}^{d-2}_{O^\top(v)}).  
\end{align*}
Hence, by the rotation invariance of $\sigma_{d-1}$ we find that
\begin{align*}
\mu(O(A))&=  \int_{\mathbb{S}^{d-1}} \sigma_{d-2}(O(A)\cap  \mathbb{S}^{d-2}_v ) \frac{ \d \sigma_{d-1}(v)}{|\mathbb{S}^{d-2}|} =  \int_{\mathbb{S}^{d-1}} \sigma_{d-2}(A\cap  \mathbb{S}^{d-2}_{O^\top(v)}) \frac{ \d \sigma_{d-1}(v)}{|\mathbb{S}^{d-2}|}\\
&=     \int_{\mathbb{S}^{d-1}} \sigma_{d-2}(A\cap  \mathbb{S}^{d-2}_{v}) \frac{ \d \sigma_{d-1}(O(v))}{|\mathbb{S}^{d-2}|}= \mu(A). 
\end{align*}
That is, $\mu(O(A)) =  \mu(A) $ and the claim is proved. Recalling that rotation invariant Borel measures on $\mathbb{S}^{d-1} $ are unique up to a multiplicative constant, we get $\mu= c\sigma_{d-1}$.  On the other hand, once again the rotation invariance of $\sigma_{d-2}$ implies that $\sigma_{d-2}(\mathbb{S}^{d-1}\cap \mathbb{S}^{d-2}_v)  =  \sigma_{d-2}(\mathbb{S}^{d-2}_v)$ which implies that $\mu(\mathbb{S}^{d-1})= \sigma_{d-1} (\mathbb{S}^{d-1}) =|\mathbb{S}^{d-1}|$. In conclusion, $\mu= \sigma_{d-1}$. This ends the proof. 

\end{proof} 

\vspace{2mm}

\begin{lemma}\label{lem:radial-monotone}  Given $u\in L^p(\R^d)$ define
	\begin{align*}
U(h) =\int_{\R^d}|u(x+h) -u(x)|^p\d x \quad\text{for all $h\in \R^d$}. 
	\end{align*}
Then for every $0<s<t$ we have 
\begin{align*}
\int_{\mathbb{S}^{d-1}}\frac{U(tw)}{t^p} \d \sigma_{d-1}(w) \leq 2^{2p-1}\int_{\mathbb{S}^{d-1}}\frac{U(sw)}{s^p}\d \sigma_{d-1}(w).
\end{align*}
\end{lemma}
\vspace{2mm}

\begin{proof}
	Applying Lemma \ref{lem:G-estimate} to the function $g(t) = u(x+tw)$ yields 
	\begin{align*}
	\int_{\mathbb{S}^{d-1}}\frac{U(tw)}{t^p} \d \sigma_{d-1}(w) =2^{p-1} \Big(\int_{\mathbb{S}^{d-1}}\frac{U(sw)}{s^p}\d \sigma_{d-1}(w)+  \int_{\mathbb{S}^{d-1}}\frac{U(\theta sw)}{t^p}\d \sigma_{d-1}(w)\Big).
	\end{align*}
	Thus it would suffice to show that 
	\begin{align*}
	\int_{\mathbb{S}^{d-1}}\frac{U(\theta sw)}{t^p}\d \sigma_{d-1}(w)\leq 2^{p} \int_{\mathbb{S}^{d-1}}\frac{U(sw)}{t^p}\d \sigma_{d-1}(w). 
	\end{align*}
We distinguish two cases.

\noindent If $d$ is even then there exists $O\in \mathcal{O}(d)$, a rotation on $\R^d$,  such that $Ow \cdot w= 0$ for all $w\in \R^d$. Such a rotation can be constructed by using the block matrix
$\big(\begin{smallmatrix}
0&-1\\1&0
\end{smallmatrix}\big)$. 
Thereupon, consider
\begin{align*}
O_1 w = \frac{\theta}{2} w+ \frac{\theta'}{2} Ow\quad\text{and}\quad O_2 w = \frac{\theta}{2} w- \frac{\theta'}{2}Ow\;\;\text{ with  $\theta' = \sqrt{4-\theta^2}$}. 
\end{align*}
A routine check shows that $O_1, O_2\in \mathcal{O}(d)$ are rotations too and $\theta w= O_1w+O_2w$. Given that $\d \sigma_{d-1}$ is invariant under rotation, we obtain
\begin{align*}
\int_{\mathbb{S}^{d-1}}\frac{U(\theta sw)}{t^p}\d \sigma_{d-1}(w)\leq 2^{p-1} \int_{\mathbb{S}^{d-1}}\frac{U(sO_1w)+ U(sO_2w)}{t^p}\d \sigma_{d-1}(w) = 2^{p} \int_{\mathbb{S}^{d-1}}\frac{U	(sw)}{t^p}\d \sigma_{d-1}(w). 
\end{align*}
 
 \noindent Now if $d$ is odd then $d-1$ is even thus the previous case reveals that 
 \begin{align*}
 \int_{\mathbb{S}^{d-2}_v}\frac{U(\theta sw)}{t^p}\d \sigma_{d-2}(w)\leq 2^{p} \int_{\mathbb{S}^{d-2}_v}\frac{U(sw)}{t^p}\d \sigma_{d-2}(w)\quad\text{for all $v\in \mathbb{S}^{d-1}$}. 
 \end{align*}
Combining this with the integration formula of Theorem \ref{thm:integration-on-sphere}  finishes the proof since
\begin{align*}
\int_{\mathbb{S}^{d-1}}\frac{U(\theta sw)}{t^p}\d \sigma_{d-1}(w)\leq 2^{p}\int_{\mathbb{S}^{d-1}}\frac{U(sw)}{t^p}\d \sigma_{d-1}(w).
\end{align*}
\end{proof}

\vspace{1mm}

\noindent The next result is a variant of the Riesz-Fr\'echet-Kolmogorov theorem. 
\begin{theorem}\label{thm:asympto-local-compact}
Assume $d\geq 2$. Let $(u_n)_n$ be a bounded sequence in $L^p(\R^d)$ such that 
\begin{align*}
	A_p:=\sup_{n\geq 1}\iil_{\R^d\R^d}|u_n(x) -u_n(y)|^p\nu_{\eps_n}(x-y)\d y\d x<\infty. 
\end{align*}
Then under the assumptions \eqref{eq:concentration-property-n} the sequence  $(u_n)_n$ is relatively compact in $L^p_{\operatorname{loc}}(\R^d)$. 
\end{theorem}

\vspace{2mm}
\begin{proof} For fixed $\delta>0$, from the assumption \eqref{eq:concentration-property-n} we can choose $n_\delta\geq 1$ such that 
\begin{align*}
	\int_{B_\delta(0)} (1\land |h|^p) \nu_{\eps_n}(h) \d h \geq \frac{1}{2}\qquad\text{for all}\quad n\geq n_\delta.
\end{align*}
In virtue of Lemma \ref{lem:radial-monotone}, for $0<s<\delta < t$ we have 
\begin{align*}
\int_{\mathbb{S}^{d-1}}\frac{U_n(tw)}{t^p} \d \sigma_{d-1}(w) \leq 2^{2p-1}\int_{\mathbb{S}^{d-1}}\frac{U_n(sw)}{s^p}\d \sigma_{d-1}(w).
\end{align*}
Thus,
\begin{align*}
\frac{1}{2|\mathbb{S}^{d-1}|}\int_{\mathbb{S}^{d-1}}\frac{U_n(tw)}{t^p} \d \sigma_{d-1}(w)
&\leq\frac{1}{|\mathbb{S}^{d-1}|} 	\int_{B_\delta(0)} (1\land |h|^p) \nu_{\eps_n}(h) \d h \int_{\mathbb{S}^{d-1}}\frac{U_n(tw)}{t^p} \d \sigma_{d-1}(w)\\
&= 	\int_0^\delta \int_{\mathbb{S}^{d-1}}\frac{U_n(tw)}{t^p} \d \sigma_{d-1}(w)\, s^{d-1}(1\land s^p) \nu_{\eps_n}(s) \d s\\
 &\leq 2^{2p-1}\int_0^\delta \int_{\mathbb{S}^{d-1}}\frac{U_n(sw)}{s^p} \d \sigma_{d-1}(w) \,  s^{d-1}(1\land s^p) \nu_{\eps_n}(s) \d s\\
 &\leq 2^{2p-1}\int_0^\infty\int_{\mathbb{S}^{d-1}} U_n(sw)\d \sigma_{d-1}(w) \,  s^{d-1} \nu_{\eps_n}(s) \d s\\
 &= 2^{2p-1}\iil_{\R^d\R^d} |u_n(x+h) -u_n(x)|^p\nu_{\eps_n}(h) \d h\ dx\leq  2^{2p-1}A_p.  
\end{align*}
 We have shown that
 \begin{align*}
 \int_{\mathbb{S}^{d-1}}\frac{U_n(tw)}{t^p} \d \sigma_{d-1}(w)
 \leq  2^{2p}|\mathbb{S}^{d-1}|A_p,\quad\text{ for all $t\geq \delta$ and all $n \geq n_\delta$}.
 \end{align*} 
In particular letting $C_p= 2^{2p}|\mathbb{S}^{d-1}|A_p$ we have 
 \begin{align*}
\int_{\mathbb{S}^{d-1}} U_n(\delta w)\d \sigma_{d-1}(w)
\leq  C_p\delta^p\quad\text{ for  all $n \geq n_\delta$}.
\end{align*} 

Now for $0<s<\delta$, writing $s=-\delta+ (s+ \delta)$,  we also have 
\begin{align*}
\int_{\mathbb{S}^{d-1}} U_n(s w)\d \sigma_{d-1}(w)
&\leq 2^{p-1}\int_{\mathbb{S}^{d-1}} U_n(\delta w)\d \sigma_{d-1}(w)+2^{p-1} \int_{\mathbb{S}^{d-1}} U_n((s+\delta) w)\d \sigma_{d-1}(w)\\
&\leq  2^{p-1}C_p(\delta^p+ (\delta+s)^p)\leq 2^{2p}C_p\delta^p.
\end{align*} 
 Therefore we get 
 \begin{align*}
 \int_{\mathbb{S}^{d-1}} U_n(sw)\d \sigma_{d-1}(w)
 &\leq 2^{4p}A_p|\mathbb{S}^{d-1}|\delta^p,\quad\text{ for all $0<s<\delta $ and  for  all $n \geq n_\delta$}.
 \end{align*} 
 Now assume $\phi \in C_c^\infty(\mathbb{R}^d)$ is radial, supported in the unit ball $B_1(0)$, $\phi\geq 0$, $\int_{\mathbb{R}^d} \phi(x)\d x =1$ and $\phi_{\delta}(x)=\frac{1}{ \delta^d}\phi(\frac{x}{\delta})$. Using this, the  
Jensen inequality and the last estimate it follows that, for $n\geq n_\delta$, we have 
\begin{align*}
\int_{\R^d}|\phi_\delta *u_n(x)- u_n(x)|^p\d x
&\leq \int_{B_\delta} \phi_\delta(h)\d h\int_{\R^d}| u_n(x+h) -u_n(x)|^p\d x\\
&= \int_0^\delta  \phi_\delta(s) s^{d-1} \d s\int_{\mathbb{S}^{d-1}} \int_{\R^d}| u_n(x+sw) -u_n(x)|^p\d x \d \sigma_{d-1}(w)\\
&\leq 2^{4p}A_p|\mathbb{S}^{d-1}|\delta^p \int_0^\delta  \phi_\delta(s) s^{d-1} \d s=2^{4p}A_p\delta^p. 
\end{align*}
That is we have 
\begin{align*}
\sup_{n\geq n_\delta}\int_{\R^d}|\phi_\delta *u_n(x)- u_n(x)|^p\d x\leq 2^{4p}A_p\delta^p. 
\end{align*}
For $\eps>0$, we fix $\delta_\eps:=\delta>0$ such that 
\begin{align}\label{eq:wwphi}
\sup_{n\geq n_\delta} \int_{\R^d}|\phi_\delta *u_n(x)- u_n(x)|^p\d x<\varepsilon/2.
\end{align}
Let $K\subset \R^d$ be a compact with positive measure, to conclude we need to show that $\mathcal{F}=\{u_n:\, n\geq1\}$ is totally bounded in $L^p(K)$. Note that $(\phi_\delta*u_n)$ is a bounded and equicontinuous sequence in $C(K)$. 
Indeed, if we let $C=\sup\limits_{n\geq 1}\|u_n\|_{L^p(\mathbb{R}^d)}$ then by Young's inequality we find that 
\begin{align*}
&\sup_{n\geq 1}\|u_n*\phi_\delta\|_{L^\infty(\mathbb{R}^d)}\leq C\|\phi_\delta\|_{L^{p'}(\mathbb{R}^d)},\\ 
&\hspace{-1ex}\sup\limits_{n\geq1 }\|\tau_h u_n*\phi_\delta-u_n*\phi_\delta\|_{C( K)}
\leq C |h|\|\nabla \phi_\delta\|_{L^{p'}(\mathbb{R}^d)}\xrightarrow{|h|\to 0}0.
\end{align*}
According to the Arzel\`{a}-Ascoli Theorem \ref{thm:ascoli}, the set $\mathcal{F}*\phi_{\delta}\mid_K$ is precompact in $C(K)$ and hence is totally bounded. Whence there exist $g_1,\cdots,g_N\in \mathcal{F}*\phi_{\delta}\mid_K$ such that 
$$\mathcal{F}*\phi_\delta\mid_K\subset \bigcup_{i=1}^N B^\infty_{\eps_K}(g_i),\quad\text{with $\eps_K=|K|^{-1/p} \eps/2$ }.$$
From this we have $$\mathcal{F}\mid_K\subset\{u_1,\cdots, u_{n_\delta}\}\cup \bigcup_{i=1}^N B_{\eps}(g_i). $$ Indeed if $\phi_\delta*u_n \in B^\infty_{\eps_K}(g_i)$, with $n\geq n_\delta$ then   $u_n \in B_{\eps}(g_i)$ since
\begin{align*}
\|u_n-g_i\|_{L^p(K)}\leq \|u_n-\phi_\delta*u_n\|_{L^p(K)}+ \|\phi_\delta*u_n-g_i\|_{L^p(K)}< \eps/2+ \eps_K|K|^{1/p}=\eps. 
\end{align*}
We conclude that $\mathcal{F} =(u_n)_n$ is precompact in $L^p(K)$. This achieves the proof. 
\end{proof}

\noindent Let us now  look to the one dimensional case. The proof in one dimension  can be obtained by modifying the higher dimension case analogously.   
\begin{theorem}\cite[Lemma 7.1]{Ponce2004}
	Assume $d=1$ and  without loss of generality that $\nu_{\eps_n}$ is defined for every $x\in \R$. Assume that  there exist $\theta_0\in (0,1)$ and $c_0>0$ such that 
	\begin{align}\label{eq:one-dimension-condition}
	\nu_{n, \theta_0}(x) := \inf_{\theta_0\leq \theta \leq 1}\nu_{\eps_n}(\theta x)>\frac{c_0}{4}.
	\end{align}
 Let $(u_n)_n$ be a bounded sequence in $L^p(\R)$ such that 
\begin{align*}
A_p:=\sup_{n\geq 1}\iil_{\R\R}|u_n(x) -u_n(y)|^p\nu_{\eps_n}(x-y)\d y\d x<\infty. 
\end{align*}
Then the sequence  $(u_n)_n$ is relatively compact in $L^p_{\operatorname{loc}}(\R)$. 
\end{theorem}
\medskip

\noindent We are now in a position to establish asymptotic compactness on a bounded domains $\Omega$.  Let us proceed with the following result.

%
\begin{theorem}\label{thm:ponce-inquality}
 Let $\Omega\subset \R^d$ be open.  Assume $(u_n)_n$ is a sequence converging to $u$  in  $L^p(\Omega)$ and 
\begin{align*}
A_p:=\liminf_{n\to \infty}\iil_{\Omega\Omega}|u_n(x) -u_n(y)|^p\nu_{\eps_n}(x-y)\d y\d x<\infty. 
\end{align*}
\begin{enumerate}[$(i)$]
\item If $1<p<\infty$, then $u\in W^{1,p}(\Omega)$ and $\|\nabla u\|^p_{L^p(\Omega)}\leq A_p K_{d,p}^{-1}$.
\item If $p=1$, then   $u\in BV(\Omega)$ and $|u|_{BV(\Omega)} \leq A_1 K_{d,1}^{-1}$.
\end{enumerate}
\end{theorem}

\medskip

\begin{proof}
	Let us fix $\delta >0$ sufficiently small and recall, $\Omega_\delta =\{x\in \Omega: \operatorname{dist}(x,\partial\Omega)>\delta\}.$ Define the mollifier $\phi_\delta(x)= \frac{1}{\delta^d}\phi\left(\frac{x}{\delta}\right)$ with support in $B_\delta(0)$ where $\phi \in C_c^{\infty}(\mathbb{R}^d)$ is supported in $B_1(0)$, $\phi \geq 0$ and $ \int_{} \phi = 1$. 
	Assume that $u_n$ and $u$ are extended by zero off $\Omega$. Let $u^\delta_n = u_n*\phi_\delta$ and $u^\delta = u*\phi_\delta$. Regarding \eqref{eq:molification-convex-Jessen} we get 
\begin{align*}
	\iil_{\Omega_\delta\Omega_\delta}|u^\delta_n(x) -u^\delta_n(y)|^p\nu_{\eps_n}(x-y)\d y\d x\leq \iil_{\Omega\Omega} |u_n(x)-u_n(y)|^p\nu_{\eps_n}(x-y)\d y\d x.
	\end{align*}
By Minkowski's inequality  we have 
\begin{align*}
\Big(\iil_{\Omega_\delta\Omega_\delta}|u^\delta(x) -u^\delta(y)|^p\nu_{\eps_n}(x-y)\d y\d x\Big)^{1/p}
&\leq 
\Big(\iil_{\Omega_\delta\Omega_\delta}|u^\delta_n(x) -u^\delta_n(y)|^p\nu_{\eps_n}(x-y)\d y\d x\Big)^{1/p}\\
&+ \Big(\iil_{\Omega_\delta\Omega_\delta}|[u^\delta-u^\delta_n](x) -[u^\delta-u^\delta_n](y)|^p\nu_{\eps_n}(x-y)\d y\d x\Big)^{1/p}\\
&\leq 
\Big(\iil_{\Omega\Omega}|u_n(x) -u_n(y)|^p\nu_{\eps_n}(x-y)\d y\d x\Big)^{1/p} \\
&+ \Big(\iil_{\R^d \R^d}|[u^\delta-u^\delta_n](x) -[u^\delta-u^\delta_n](y)|^p\nu_{\eps_n}(x-y)\d y\d x\Big)^{1/p}
\end{align*}
Passing to the liminf, according to  the estimate \eqref{eq:levy-p-estimate}, the assumption \eqref{eq:normalized-inetrals-n}  and Young's inequality we obtain
\begin{align}
\begin{split}
\liminf_{n\to \infty}\Big(\iil_{\Omega_\delta\Omega_\delta}|u^\delta(x) -u^\delta(y)|^p
&\nu_{\eps_n}(x-y)\d y\d x\Big)^{1/p}\\
&\leq  A_p^{1/p}+2 \| u^\delta-u^\delta_n\|_{W^{1,p}(\R^d)} \int_{\R^d}(1\land |h|^p)\nu_{\eps_n}(h)\d h\\
&= A_p^{1/p}+2 \liminf_{n\to \infty}\| \phi_\delta*(u-u_n)\|_{W^{1,p}(\R^d)} \\
&\leq  A_p^{1/p}+2\|\phi_\delta\|_{W^{1,p}(\R^d)} \liminf_{n\to \infty}\|u-u_n\|_{L^p(\R^d)}\\
&= A_p^{1/p}.
\end{split}
\end{align}

\noindent In view of Theorem \ref{thm:liminf-BBM} we get 

\begin{align}
\begin{split}
K_{d,p}\int_{\Omega_\delta}|\nabla u^\delta(x)|^p\d x
\leq \liminf_{n\to \infty}\iil_{\Omega_\delta\Omega_\delta}|u^\delta(x) -u^\delta(y)|^p\nu_{\eps_n}(x-y)\d x\d y \leq A_p\quad\text{ for all $\delta>0$ .}
\end{split}
\end{align}
Let $\varphi\in C_c^\infty(\Omega)$ and  $i=1,\cdots, d$ then it is not difficult to check that 
\begin{align*}
\Big|\int_\Omega \partial_{x_i} \varphi(x) u(x)\d x\Big|
&= \lim_{\delta\to 0} \Big|\int_{\Omega_\delta} \varphi(x) \partial_{x_i} u^\delta(x)\d x\Big|\\
&\leq \limsup_{\delta\to 0} \|\varphi\|_{L^{p'}{\Omega}} \Big(\int_{\Omega_\delta} |\partial_{x_i} u^\delta(x)|^p \d x\Big)^{1/p}\\ 
&\leq \|\varphi\|_{L^{p'}{\Omega}} A_p^{1/p}K_{d,p}^{-1}. 
\end{align*}
Hence if $1<p<\infty$ then the linear form $\varphi\mapsto \int_\Omega \partial_{x_i} \varphi(x) u(x)\d x$
uniquely extends to a continuous linear form on $L^{p'}(\Omega)$. Wherefore there exists $g_i\in L^p(\Omega)$ such that
\begin{align*}
\int_\Omega \partial_{x_i} \varphi(x) u(x)\d x= \int_\Omega g_i u(x)\d x,\quad\text{ for all $\varphi\in C_c^\infty(\Omega)$.}
\end{align*}
It turns out that $u\in W^{1,p}(\Omega)$ and $\partial_{x_i}  u= -g_i$. Therefore $\nabla u^\delta = \phi_\delta*\nabla u$ so that letting $\delta\to 0$  in
\begin{align*}
	K_{d,p}\int_{\Omega}|\nabla u(x)|^p\d x = \lim_{\delta\to 0}  K_{d,p}\int_{\Omega_\delta}|\phi_\delta*\nabla u(x)|^p\d x\leq A_p.
\end{align*}

\noindent If $p=1$ then  $u\in BV(\Omega)$  because an analogous procedure  yields that (see \eqref{eq:bounded-variation}) 
\begin{align*}
|u|_{BV(\Omega)}= \sup\Big\{\Big|\int_\Omega \operatorname{div} \varphi(x) u(x)\d x\Big|: \, \|\varphi\|_{L^{\infty}(\Omega)} \leq 1\Big\} \leq  A_1 K_{d,1}^{-1}. 
\end{align*}
%

\end{proof}

\medskip

\noindent Next we proceed with some estimates near the boundary of a Lipschitz domain. 

\begin{lemma}\label{lem:estimate-cut-off-sharp}
	Let $\Omega \subset \mathbb{R}^d$ be open bounded. $\nu:\mathbb{R}^d\setminus\{0\} \to [0, \infty]$ satisfies \eqref{eq:plevy-integrability} (see page  \pageref{eq:plevy-integrability}). For $\delta>0$ small enough, let $\varphi \in C^\infty(\Omega)$ be such that $ \varphi=0$ on $\Omega_\delta$, $ \varphi=1$ on $\Omega\setminus \Omega_{\delta/2}$, $0\leq \varphi\leq 1$ and $|\nabla \varphi|\leq c/\delta$ (cf. Lemma \ref{lem:cut-off-existence}  by taking $\varphi = 1-\varphi_\delta$). Then for every $u\in L^p(\Omega)$, the following estimate holds true
	%
	\begin{align}\label{eq:estimate-cut-off-sharp}
	\begin{split}
	\iint\limits_{\Omega\Omega}\big|[u\varphi](x)-[u\varphi](y)\big|^p&\nu(x-y)\d x\d y
	\leq   2^{p-1}\iint\limits_{\Omega\Omega}|u(x)-u(y)|^p\nu(x-y)\d x\d y\\&+ \frac{2^{p-1}  c^p}{\delta^p} \Big( \int\limits_{ \Omega_{\delta/4}}|u(x)|^p\d x \il_{B_R(0)} \hspace{-2ex}|h|^p\nu(h)\d h
	+ \int\limits_{ \Omega}|u(x)|^p\d x 
	\hspace{-2ex}\il_{B^c_{\delta/4}(0)}\hspace{-2ex} (1\land |h|^p)\nu(h)\d  h\Big). 
	\end{split}
	\end{align}
\end{lemma}
\begin{proof}  Since  $\big|[u\varphi](x)-[u\varphi](y)\big|^p \leq 2^{p-1}|u(x)-u(y)|^p+ \frac{2^{p-1}c^p}{\delta^p} |u(x)|^p|\varphi(x)-\varphi(y)|^p$ we get 
	\begin{align*}
	&\iint\limits_{\Omega\Omega}\big|[u\varphi](x)-[u\varphi](y)\big|^p\nu(x-y)\d x\d y\\
	&\leq 2^{p-1}\iint\limits_{\Omega\Omega}|u(x)-u(y)|^p\nu(x-y)\d x\d y+ 2^{p-1} \int\limits_{ \Omega}|u(x)|^p\d x \int_{\Omega}\big|\varphi(x)-\varphi(y)\big|^p\nu(x-y)\d y.
	\end{align*}
Note that  $\Omega\subset B_R(x)$ for all $ x\in \Omega$ with $R= \operatorname{diam}(\Omega)$ and for all
 $x\in \Omega\setminus\Omega_{\delta/4}$ we have  $|x-y|\geq \big|\dist(x,\partial\Omega)- \dist(y,\partial\Omega)\big|\geq \delta/4$  and $ \varphi(x)= \varphi(y)=1$. Moreover, $|\varphi(x)-\varphi(y)|\leq \frac{c}{\delta}(1\land |x-y|)$ for every $x,y \in \Omega$, since $1\leq \varphi\leq 1$ and  $|\varphi(x)-\varphi(y)|\leq \frac{c}{\delta}|x-y|$.  Hence,  we get the following estimates 
	%
	\begin{align*}
	&\int\limits_{ \Omega}|u(x)|^p\d x \int_{\Omega}\big|\varphi(x)-\varphi(y)\big|^p\nu(x-y)\d y \\
	&=\int\limits_{ \Omega_{\delta/4}}|u(x)|^p\d x \int_{\Omega}\big|\varphi(x)-\varphi(y)\big|^p\nu(x-y)\d y
	+ \int\limits_{ \Omega\setminus \Omega_{\delta/4}}|u(x)|^p\d x 
	\il_{\Omega_{\delta/2}}\big|\varphi(x)-\varphi(y)\big|^p\nu(x-y)\d y\\
	&\leq \frac{c^p}{\delta^p}  \int\limits_{ \Omega_{\delta/4}}|u(x)|^p\d x \int_{\Omega}|x-y|^p\nu(x-y)\d y
	+ \frac{c^p}{\delta^p}  \int\limits_{ \Omega\setminus \Omega_{\delta/4}}|u(x)|^p\d x 
	\il_{\Omega_{\delta/2}} ( 1\land |x-y|^p)\nu(x-y)\d y\\
	&\leq   \frac{c^p}{\delta^p}  \int\limits_{ \Omega_{\delta/4}}|u(x)|^p\d x \int_{B_R(x)}|x-y|^p\nu(x-y)\d y
	+ \frac{c^p}{\delta^p}  \int\limits_{ \Omega}|u(x)|^p\d x 
	\il_{B^c_{\delta/4}(x)} ( 1\land |x-y|^p)\nu(x-y)\d y\\
	&= \frac{c^p}{\delta^p}  \int\limits_{ \Omega_{\delta/4}}|u(x)|^p\d x \int_{B_R(0)}|h|^p\nu(h)\d h
	+ \frac{c^p}{\delta^p}  \int\limits_{ \Omega}|u(x)|^p\d x 
	\il_{B^c_{\delta/4}(0)} ( 1\land |h|^p)\nu(h)\d h,
	\end{align*}
	\noindent  Altogether, gives \eqref{eq:estimate-cut-off-sharp}.  
	%
\end{proof}

\vspace{2mm}
\begin{lemma}\label{lem:estimate-near-boundary-sharp}
	Assume  $ \Omega\subset \mathbb{R}^d$ is open bounded with a Lipschitz boundary.  Assume $\nu:\mathbb{R}^d\setminus\{0\}\to [0, \infty]$ is a radial function and satisfies \eqref{eq:plevy-integrability} (see page  \pageref{eq:plevy-integrability}). Then there exists a positive constant $C=C(\Omega, d, p)$  independent of $\nu$ such that for every $u \in L^p(\Omega)$  and every $\delta>0$ small enough, 
	\begin{align}\label{eq:estimate-boundary-sharp}
	\begin{split}
\int_{\Omega} |u(x)|^p\d x
&\leq  CB(\delta,\nu) \int_{ \Omega_{\delta/4}} |u(x)|^p\d x +  CD(\delta,\nu)\int\limits_{ \Omega}|u(x)|^p\d x\\ 
&+C\delta^p A(\delta, \nu)  \iint\limits_{\Omega\Omega}|u(x)-u(y)|^p\nu(x-y)\d x\d y
	\end{split}
	\intertext{where}
	\begin{split}
		A(\delta, \nu)&= \Big(\il_{B_{2\delta}(0)} |h|^p\nu(h) \d h\Big)^{-1},~~\\
	B(\delta, \nu)&=  A(\delta, \nu)\il_{B_R(0)} \hspace{-2ex}|h|^p\nu(h)\d h,\,\\
	D(\delta, \nu)&= A(\delta, \nu) \il_{B^c_{\delta/4}(0)}\hspace{-2ex} (1\land|h|^p)\nu(h)\d  h. 
		\end{split}
\end{align}

\end{lemma}

\begin{proof}
\noindent Recall that by  the relation \eqref{eq:estimate-near-Lipschitz-boundary} we have 
\begin{align*}
\begin{split}
\il_{\Omega \setminus \Omega_{\delta/2}} | [u\varphi](x)|^p\d x
&\leq C\delta^p A(\delta, \nu)\iil_{ \Omega \Omega}  | [u\varphi](y) - [u\varphi](x)|^p\nu(y-x)\d y\d x. 
\end{split}
\end{align*} 
Combining this with Lemma \ref{lem:estimate-cut-off-sharp} leads to the following 
\begin{align*}
\int_{ \Omega} |u(x)|^p\d x
&\leq  \int_{ \Omega_{\delta/4}} |u(x)|^p\d x+ \int_{ \Omega\setminus\Omega_{\delta/2}} |u(x)|^p\d x\\
&\leq \int_{ \Omega_{\delta/4}} |u(x)|^p\d x+ C\delta^pA(\delta, \nu) \Big\{ 2^{p-1}\iint\limits_{\Omega\Omega}|u(x)-u(y)|^p\nu(x-y)\d x\d y\\
&+ \frac{2^{p-1}  c^p}{\delta^p} \Big( \int\limits_{ \Omega_{\delta/4}}|u(x)|^p\d x \il_{B_R(0)} \hspace{-1ex}|h|^p\nu(h)\d h
+ \int\limits_{ \Omega}|u(x)|^p\d x 
\hspace{-2ex}\il_{B^c_{\delta/4}(0)}\hspace{-2ex} (1\land|h|^p) \nu(h)\d  h\Big)\Big\}\\
&\leq  CB(\delta,\nu) \int_{ \Omega_{\delta/4}} |u(x)|^p\d x +  CD(\delta,\nu)\int\limits_{ \Omega}|u(x)|^p\d x\\ &\qquad+C\delta^p A(\delta, \nu)  \iint\limits_{\Omega\Omega}|u(x)-u(y)|^p\nu(x-y)\d x\d y. 
\end{align*}
\noindent Note that here $C$ is a generic constant which neither  depends on $\nu$ nor on $\delta$. 
\end{proof}

\bigskip

\begin{theorem}\label{thm:asymp-compactness}
	Assume $d\geq 2, 1\leq p<\infty$.  Let $\Omega\subset \R^d$ be open with a Lipschitz boundary.
Assume that the sequence $(\nu_{\eps_n})_n$  satisfies the conditions \eqref{eq:normalized-inetrals-n} and \eqref{eq:concentration-property-n}. Assume $(u_n)_n$ is a  bounded sequence in $L^p(\Omega)$ such that
	\begin{align*}
	A_p:=\liminf_{n\to \infty}\iil_{\Omega\Omega}|u_n(x) -u_n(y)|^p\nu_{\eps_n}(x-y)\d y\d x<\infty. 
	\end{align*}
	Then  $(u_n)_n$ has a converging subsequence in $L^p(\Omega)$. Moreover, assume $u\in L^p(\Omega)$
	is the limit of a subsequence of $(u_n)_n$ then the following hold true.
	\begin{enumerate}[$(i)$]
		\item If $1<p<\infty$, then $u\in W^{1,p}(\Omega)$ and $\|\nabla u\|^p_{L^p(\Omega)}\leq A_p K_{d,p}^{-1}$.
		\item If $p=1$, then   $u\in BV(\Omega)$ and $|u|_{BV(\Omega)} \leq A_1 K_{d,1}^{-1}$.
	\end{enumerate}
The same holds true in dimension $d=1$, provided that the condition  \eqref{eq:one-dimension-condition} holds. 
%
\end{theorem}

\vspace{2mm}

\begin{proof} 
	For $\delta>0$ small enough, we let $\varphi_\delta \in C_c^\infty(\Omega)$ be  such that $\varphi_\delta =1$ on $\Omega_\delta$  and  $0\leq \varphi_\delta\leq 1$ as in Lemma \ref{lem:cut-off-existence}. From  estimate \eqref{eq:estimate-extending-to-rd} we get
	%
	%
	\begin{align*}
&\iil_{ \R^d \R^d}|\varphi_\delta(x)u_n(x) -\varphi_\delta(y)u_n(y)|^p\nu_{\eps_n}(x-y)\d x \d y\\
	&\leq 2^p\iil_{\Omega\Omega}|u_n(x) -u_n(y)|^p\nu_{\eps_n}(x-y)\d x \d y+
	2^p \|\varphi_\delta\|^p_{W^{1,\infty}} \il_{ \Omega}|u_n(x)|^p \d x \int_{\R^d}(1\land|x-y|^p)\nu_{\eps_n}(h)\d h\\
	&+2 \il_{ \Omega_{\delta/2}}|u_n(x)|^p \d x \il_{|h|\geq \delta/2}\nu_{\eps_n}(h)\d h\leq 2^pA_p + C_\delta\il_{ \Omega}|u_n(x)|^p\d x. 
	\end{align*}
Therefore, for each $\delta>0$, $(\varphi_\delta u_n)_n$ is bounded in $L^p(\R^d)$  and we have
	\begin{align*}
	\liminf_{n \to \infty}\iil_{ \R^d \R^d}|\varphi_\delta(x)u_n(x) -\varphi_\delta(y)u_n(y)|^p\nu_{\eps_n}(x-y)\d x \d y
	<\infty.
	\end{align*}
According to Theorem \ref{thm:asympto-local-compact}, the sequence $(\varphi_\delta u_n)_n$ is relatively compact in $L^p(\Omega)$.   Employing Cantor's diagonalization procedure as in the proof of Corollary \ref{cor:local-compatcness} one is able to find that  $(u_n)_n$  has a subsequence that we still denote by  $(u_n)_n$ converging to some measurable function $u$ in $L^p_{\operatorname{loc}}(\Omega) $ and  a.e in $\Omega$.  Necessarily, $u\in L^p(\Omega)$. Indeed by Fatou's lemma,  $\|u\|_{L^p(\Omega)}\leq \liminf\limits_{n \to \infty}\|u_n\|_{L^p(\Omega)}<\infty$ since $(u_n)_n$ is bounded in $L^p(\Omega)$. 

\vspace{1mm}
\noindent Next, we show that $\|u_n-u\|_{L^p(\Omega)}\xrightarrow{n\to \infty}0.$ In reference to Lemma 
\ref{lem:estimate-near-boundary-sharp} for $0<\delta<1$ sufficiently small, we have the following
\begin{align*}
\begin{split}
\int_{\Omega} |u_n(x) -u(x)|^p\d x
&\leq  CB(\delta,\nu_{\eps_n}) \int_{ \Omega_{\delta/4}} |u_n(x) -u(x)|^p\d x +  CD(\delta,\nu_{\eps_n})\int\limits_{ \Omega}|u(x)|^p\d x,\\ 
&+C\delta^p A(\delta, \nu_{\eps_n})  \iint\limits_{\Omega\Omega}|u_n(x) -u(x)-[u_n(y) -u(y)]|^p\nu_{\eps_n}(x-y)\d x\d y,\\
&\leq  CB(\delta,\nu_{\eps_n}) \int_{ \Omega_{\delta/4}} |u_n(x) -u(x)|^p\d x +  CM\big\{ D(\delta,\nu_{\eps_n}) +\delta^p A(\delta, \nu_{\eps_n}) \big\}.
\end{split}
\end{align*}
\noindent Here, the constant $C=C(d,p,\Omega)$ is independent of $\delta$ and $n$, and 

$$M= 2^p(\|u\|^p_{\WnuOm}+\sup_{n}\|u_n\|^p_{L^p(\Omega)}+A_p)<\infty.$$ 
Further, in view of the assumptions \eqref{eq:normalized-inetrals-n} and \eqref{eq:concentration-property-n}, by taking into account Remark \ref{rem:asymp-nu} we have 
\begin{align*}
&A(\delta, \nu_{\eps_n})= \Big(\il_{B_{2\delta}(0)} |h|^p\nu_{\eps_n}(h) \d h\Big)^{-1}\xrightarrow{n\to \infty}1\\
&B(\delta, \nu_{\eps_n})=  A(\delta, \nu_{\eps_n})\il_{B_R(0)} \hspace{-2ex}|h|^p\nu_{\eps_n}(h)\d h \xrightarrow{n\to \infty}1\\
&D(\delta, \nu_{\eps_n})= A(\delta, \nu_{\eps_n}) \il_{B^c_{\delta/4}(0)}\hspace{-2ex} (1\land|h|^p)\nu_{\eps_n}(h)\d  h\xrightarrow{n\to \infty}0.
\end{align*}
By the convergence in $L^p_{\operatorname{loc}}(\Omega)$, it follows that $\|u_n-u\|_{L^p(\Omega_{\delta/4})}\xrightarrow{n\to \infty}0.$ Altogether, passing the previous estimate to the $\limsup$,  we remain with the following 
\begin{align*}
\limsup_{n \to \infty}
\int_{\Omega} |u_n(x) -u(x)|^p\d x
&\leq CM\delta^p. 
\end{align*}
Letting $\delta\to 0$, it follows that $\|u_n-u\|_{L^p(\Omega)}\xrightarrow{n\to \infty}0$  since we have $\limsup\limits_{n \to \infty}\|u_n-u\|_{L^p(\Omega)}=0.$
 The proof is complete since the statements $(i)$ and $(ii)$ follow immediately from Theorem \ref{thm:ponce-inquality}.
\end{proof}

\bigskip

\section{Robust Poincar\'e inequalities}

Poincar\'e type inequalities are crucial tools to show the well-posedness for certain classes of variational problems. 
Moreover, the sharp constants for such inequalities are strongly connected to the analysis of eigenvalues of 
certain  operators. In this section, we aim to establish\emph{ robust Poincar\'e  type inequalities} for the class of 
\textit{nonlocal Sobolev-like spaces} introduced in Chapter \ref{chap:nonlocal-sobolev}. The Robustness should be understood in the sense that within such inequalities, one is able to recover the corresponding  classical Poincar\'e inequalities for classical  Sobolev spaces. See the upcoming Corollary \ref{cor:clasical-poincare}.

\vspace{2mm}

\noindent Recall  that for an open bounded connected and Lipschitz set $\Omega\subset \mathbb{R}^d$ ($d\geq 1$) and $1\leq p<\infty$, in the local setting we say that a non-empty subset $G$ of  $L^p(\Omega)$ is admissible for the Poincar\'e inequality if there exists a constant $A=A(\Omega, d, p, G)>0$  such that
\begin{align}\label{eq:poincare-local}
\int_{\Omega}|u(x)|^pdx \leq A\int_{\Omega} |\nabla u(x)|^p\, dx \quad\quad\hbox{for all} ~ ~u\in G. 
\end{align} 
Here $\nabla u$ is the distributional gradient of $u$. By convention, we let $\il_{\Omega} |\nabla u(x)|^p\, dx =\infty $ if $\nabla u \not \in L^p(\Omega)$.   

%

\vspace{2mm}
\begin{remark}
  It is worth noticing  that an admissible set $G$ for inequalities of types \eqref{eq:poincare-local} does not contain non-zero constant functions. This is obvious since these inequalities cannot not  hold true for such  functions.
Let us mention some standard examples of admissible  set $G$ on which the  inequality $\eqref{eq:poincare-local}$ holds true. 
\end{remark}

\begin{itemize}
\item $G_0=W_0^{1,p}(\Omega)$. In this case the  inequality $\eqref{eq:poincare-local}$ is known as  Poincar\'e-Friedrichs inequality. 
	
\item $G_1=  L^{p}(\Omega)^{\perp}=  \big\{u\in L^p(\Omega):\, \mbox{$\fint_{\Omega}u :=  \frac{1}{|\Omega|}\int_{\Omega } u(x)dx =0$}\big\}.$ In this case the  inequality $\eqref{eq:poincare-local}$ corresponds to the Poincar\'e inequality.  For all $u\in L^p(\Omega)$ then,  $u- \fint_{\Omega } u dx \in L^p(\Omega)^\perp$. Thus the inequality \eqref{eq:poincare-local}  can be rewritten as follows 
\begin{align*}
\int_{\Omega}|u(x) - \mbox{$\fint_{\Omega } u$}|^pdx \leq A\int_{\Omega} |\nabla u(x)|^p\, dx, \quad\quad\hbox{for all} ~ ~u\in L^p(\Omega). 
\end{align*} 

\item $G_2=\big\{u\in L^p(\Omega):\, \big| \{u=0 \}\big| \geq \delta \big\}$ for some $0<\delta<|\Omega|$.
\vspace{1mm}

\item $G_3=\big\{u\in L^p(\Omega ): \, \mbox{$\fint_E u =0$ } \big\}$  for a measurable set $E\subset \Omega$ with a positive measure.   
\item  $G_4=\big\{u\in L^p(\Omega ): u= 0 \quad \hbox{a.e on }~E \big\}$  for a measurable set $E\subset \Omega$ with a positive measure.  
\end{itemize}

\vspace{2mm}
\noindent It is important to remark that the dominant  link between $G_i^{'s}$  is that the null function $u=0$ is the only constant belonging to $G_i$.

\begin{theorem}[Robust Poincaré inequality]\label{thm:robust-poincare}
Let $\Omega\subset \R^d$ be a bounded connected Lipschitz open set.  Assume that the family $(\nu_\eps)_\eps$  satisfies the conditions \eqref{eq:normalized-inetrals} and \eqref{eq:concentration-property} with $1\leq p<\infty$. Let $G$ be a nonempty closed  subset of $L^p(\Omega)$ such that $c\not\in G$ for all $c\in \R\setminus \{0\}$. Assume $d\geq 2$ then there exist $\eps_0>0$ and a constant $B=B(\eps_0, d, p, \Omega, G)$  such  that for every  $0<\eps<\eps_0$ and every  $u\in G$ 
\begin{align}\label{eq:robust-poincare}
\|u\|^p_{L^p(\Omega)}\leq B\iil_{ \Omega \Omega}|u(x) -u(y)|^p\nu_\eps(x-y)\d y\d x.
\end{align}
The same is true in dimension $d=1$ provided that the condition  \eqref{eq:one-dimension-condition} holds. 
%
\end{theorem}

\medskip 

\begin{proof} 
Assume the statement is false. Then  for all $n\geq 1$ we can find $0<\eps_n<\frac{1}{2^n}$  and $u_n\in G$ such that $\|u_n\|^p_{L^p(\Omega)}=1$ and 
\begin{align*}
\iil_{ \Omega \Omega}|u_n(x) -u_n(y)|^p\nu_{\eps_n}(x-y)\d y\d x< \frac{1}{2^n}. 
\end{align*}
%
%
According to Theorem \ref{thm:asymp-compactness} the sequence $(u_n)_n$ has a subsequence (which we again denote by $u_n$) converging  in $L^p(\Omega)$ to some $u$ with $u\in W^{1,p}(\Omega)$ for $1<p<\infty$ and $u\in BV(\Omega)$ if $p=1$. Moreover, in either case we have that 
\begin{align*}
K_{d,p}\int_\Omega|\nabla u|^p\leq \liminf_{n \to \infty}\iil_{ \Omega \Omega}|u_n(x) -u_n(y)|^p\nu_{\eps_n}(x-y)\d y\d x=0.
\end{align*}
This implies that $\nabla u=0$ almost everywhere  on $\Omega$ which is a connected set. Necessarily, $u=c$ is a constant function. In addition, since $\|u_n\|^p_{L^p(\Omega)}=1$, it follows that $\|u\|^p_{L^p(\Omega)}=1$ which implies that $u\neq 0$.  Hence, $u=c\neq 0$, and by assumption we have $u=c\not \in G$. On the other hand, $u\in G$ because $\|u_n-u\|_{L^p(\Omega)}\xrightarrow{n\to \infty}0$ and by assumption $G$ is a closed subspace of $L^p(\Omega)$. This contradicts the fact that $u\not \in G$. We have reached a contradiction which means that our initial assumption was wrong. 
\end{proof}

\vspace{2mm}

\begin{corollary}\label{cor:robust-poincare}
	Under the assumption of Theorem \ref{thm:robust-poincare}, there exist $\eps_0>0$ and a constant $B= B(\eps_0, d,p, \Omega)$ such that the following inequalities hold
	\begin{align*}
\|u-\mbox{$\fint_\Omega$}\|^p_{L^p(\Omega)}&\leq B \iil_{ \Omega\Omega}|u(x)-u(y)|^p\nu_\eps(x-y)\d x\d y\quad\text{  for all $\eps\in (0,\eps_0)$ and $u\in L^p(\Omega)$}\\
\|u\|^p_{L^p(\Omega)}&\leq B \iil_{ \Omega\Omega}|u(x)-u(y)|^p\nu_\eps(x-y)\d x\d y\quad\text{ for all $\eps\in (0,\eps_0)$ and $u\in G_i,\,\, i=2,3,4$}.
	\end{align*}
\end{corollary}


\vspace{2mm}

\noindent Note that the constant  $B=B(\eps_0, d,p, \Omega, G_i)$ is  independent of $\eps>0$. Therefore, a noteworthy consequence of Theorem \ref{thm:robust-poincare} is obtained letting $\eps\to 0$ by means of Theorem \ref{thm:charact-w1p} and Theorem \ref{thm:charact-BV} thereby recovering the classical Poincar\'e type inequality.

\begin{corollary}\label{cor:clasical-poincare}
	Assume $d\geq 2, 1\leq p<\infty$.  Let $\Omega\subset \R^d$ be a connected Lipschitz open set. Assume $G$ is a closed subset of $L^p(\Omega)$ which does not contains non-zero constant functions. Then we have 
	\begin{align}\label{eq:poincare-local-bis} 
	\|u\|^p_{L^p(\Omega)}\leq B K_{d,p}\il_{ \Omega}|\nabla u(x)|^p\d x\quad \text{for every $u\in G$.}
	\end{align}

\end{corollary}

\bigskip

\noindent The analog robust Poincar\'e-Friedrichs inequality, is delicate and deserves a different treatment. 
\begin{theorem}[Robust Poincar\'e-Friedrichs inequality]\label{thm:robust-poincare-friedrichs}
	Let $\Omega\subset \R^d$ be a bounded connected open set.  Assume that the family $(\nu_\eps)_\eps$  satisfies the conditions \eqref{eq:normalized-inetrals} and \eqref{eq:concentration-property} with $1\leq p<\infty$.  Assume $d\geq 2$ then there exist $\eps_0>0$ and a constant $B=B(\eps_0, d, p, \Omega)$  such  that for every  $0<\eps<\eps_0$ and every  $u\in C_c^\infty(\R^d)$ with $\operatorname{supp}u \subset \Omega$ it holds
	\begin{align}\label{eq:robust-poincare-friedrichs}
	\|u\|^p_{L^p(\Omega)}\leq B\iil_{ \R^d\R^d}|u(x) -u(y)|^p\nu_\eps(x-y)\d y\d x
	\end{align}
The same is true in dimension $d=1$ provided that the condition  \eqref{eq:one-dimension-condition} holds. 
	
\end{theorem}

\bigskip

\begin{proof} 
	Assume the statement is false. Then  for each $n\geq 1$ we can find $0<\eps_n<\frac{1}{2^n}$  and $u_n\in C_c^\infty(\Omega)$ such that $\|u_n\|^p_{L^p(\Omega)}=1$ and 
	\begin{align*}
\iil_{ \Omega \Omega}|u_n(x) -u_n(y)|^p\nu_{\eps_n}(x-y)\d y\d x\leq 	\iil_{ \R^d \R^d}|u_n(x) -u_n(y)|^p\nu_{\eps_n}(x-y)\d y\d x< \frac{1}{2^n}. 
	\end{align*}
	%
	%
	According to Theorem \ref{thm:asymp-compactness} the sequence $(u_n)_n$ has a subsequence (which we again denote by $u_n$) converging  in $L^p(\Omega)$ to some $u$ with $u\in W^{1,p}(\Omega)$ for $1<p<\infty$ and $u\in BV(\Omega)$ if $p=1$. Moreover, in either case we have that 
	\begin{align*}
	K_{d,p}\int_\Omega|\nabla u|^p\leq \liminf_{n \to \infty}\iil_{ \Omega \Omega}|u_n(x) -u_n(y)|^p\nu_{\eps_n}(x-y)\d y\d x=0. 
	\end{align*}
	This implies that $\nabla u=0$ almost everywhere  on $\Omega$ which is a connected set. Necessarily, $u=c$ is a constant function. In addition, since $\|u_n\|^p_{L^p(\Omega)}=1$, it follows that $\|u\|^p_{L^p(\Omega)}=1$ which implies that $u\neq 0$.  Hence, $u=c\neq 0$. Next we need to show that $u\in W_0^{1,p}(\Omega)$. For $\delta >0$ small enough since each $\nu_\eps $ is radial, by mimicking the relation \eqref{eq:optimal-non-levy} we get the following 
	\begin{align*}
\frac{1}{2^n}&\geq 	\iint_{ \R^d \R^d}|u_n(x) -u_n(y)|^p\nu_{\eps_n}(x-y)\d y\d x\\
&\geq  \int_{ \R^d }\int_{|h|\leq \delta}|u_n(x) -u_n(x+h)|^p\nu_{\eps_n}(h)\d h\d x\\
& =|\mathbb{S}^{d-1}|^{-1}\Big( \int_{\R^d} \int_{\mathbb{S}^{d-1}} \left|\nabla u_n(x)\cdot w\right|^p \mathrm{d}\sigma_{d-1}(w) \d x\Big) \Big(\int_{B_\delta} |h|^p\nu_{\eps_n}(h)\d h\Big)\\
&= K_{d,p} \Big( \int_{\R^d} |\nabla u_n(x)|^p \d x\Big) \Big(\int_{B_\delta} |h|^p\nu_{\eps_n}(h)\d h\Big).
	\end{align*}
Therefore, from Remark \ref{rem:asymp-nu} we have 
\begin{align*}
K_{d,p} \int_{\R^d} |\nabla u_n(x)|^p \d x &\leq  \frac{1}{2^n}\Big(\int_{B_\delta} |h|^p\nu_{\eps_n}(h)\d h\Big)^{-1}\xrightarrow{n\to \infty}0. 
\end{align*}
Finally, we have $\|u_n-u\|_{L^p(\Omega)}\xrightarrow{n\to \infty}0$ and $\|\nabla u_n\|_{L^p(\Omega)}\xrightarrow{n\to \infty}0$. Thus $\|u_n-u\|_{W^{1,p}(\Omega)}\xrightarrow{n\to \infty}0$ 
and we get that $u\in W_0^{1,p}(\Omega)$ since $u_n\in C_c^\infty(\Omega)$. 
This contradicts the fact that $u=c\neq0$ since the null function is the only constant function belonging to $W_0^{1,p}(\Omega)$. 
\end{proof}

\vspace{2mm}
\begin{remark}
Note that statement for  the analog inequality
	\begin{align}\label{eq:fails-poincare-frie}
	\|u\|_{L^p(\Omega)}\leq B\iil_{ \Omega\Omega}|u(x) -u(y)|^p\nu_{\eps}(x-y)\d y\d x
	\end{align}
	holding for all $u\in C_c^\infty(\Omega)$ and all $\eps\in (0,\eps_0)$, is not fully satisfactory since the kernels $(\nu_\eps)_\eps$ are also allowed to be integrable. If this is the 
	case, then  the spaces $L^p(\Omega)$ and $W^p_{\nu_\eps}(\Omega)$ coincide. This would imply that the estimate \eqref{eq:fails-poincare-frie} holds true for all functions in $L^p(\Omega)$, which is impossible since it fails for the constant function $u=1$. 
	Example \ref{Ex:example-poincre1} (with $\beta= p$), i.e. $ \nu_\varepsilon(h) = \frac{d+p}{ |\mathbb{S}^{d-1}|\varepsilon^{d+p}} \mathds{1}_{B_{\varepsilon}}(h) $ is typical illustration  withstanding the present argument applies. It might be challenging to classify families $(\nu_\varepsilon)_\varepsilon$ for which the inequality \eqref{eq:fails-poincare-frie} for all $u\in C_c^\infty(\Omega)$ and all $\eps\in (0,\eps_0)$. 
For instance, we believe it \eqref{eq:fails-poincare-frie} is true for $ \nu_\varepsilon (h) =\frac{\varepsilon(p-\varepsilon)}{p|\mathbb{S}^{d-1}|} |h|^{-d-p+\varepsilon}$ with $ 0<\varepsilon<1/p$ and all $u\in C_c^\infty(\Omega)$. 
\end{remark}

\medskip

\noindent As a consequence of Theorem \ref{thm:robust-poincare-friedrichs}, by letting $\eps\to 0$, one recovers the following classical Poincar\'e-Friedrichs inequality.

\begin{corollary}\label{cor:clasical-poincare-friedrichs}  Assume $\Omega\subset \R^d$ is open, bounded and connected. Then we have 
	\begin{align}\label{eq:poincare-local-friedrichs}
	\|u\|^p_{L^p(\Omega)}\leq B K_{d,p}\il_{ \Omega}|\nabla u(x)|^p\d x\quad \text{for every $u\in W_0^{1,p}(\Omega)$ }.
	\end{align}
\end{corollary}

\vspace{2mm}

\noindent Next, we would like to study the asymptotic behavior of the sharp constant of estimate \eqref{eq:robust-poincare}. It is not so surprising that this asymptotic behavior is related to the sharp  constant of estimate \eqref{eq:poincare-local}. 
\begin{theorem}[\textbf{convergence of sharp constants I}]\label{thm:convergence-sharp-const-poinc} Let $\Omega\subset \R^d$ be open, bounded and connected Lipschitz. Let $1\leq p<\infty$.  Assume that the family $(\nu_\eps)_\eps$  satisfies the conditions \eqref{eq:normalized-inetrals} and \eqref{eq:concentration-property}. Let $G$ be a nonempty closed  subset of $L^p(\Omega)$ such that $c\not\in G$ for all $c\in \R\setminus \{0\}$. Let $\Lambda_{1,p} $ be the sharp constant of the inequality \eqref{eq:poincare-local-bis} and $\Lambda_{1-\varepsilon,p} $ be the sharp constant of the inequality \eqref{eq:robust-poincare} that is, 
	\begin{align*}
	\Lambda_{1,p} := \sup \left\{\|u\|_{L^p(\Omega)}^p:  u \in G, ~~\|\nabla u\|_{L^p(\Omega)}=1 \right\}
\intertext{and for each $\varepsilon\in (0,\varepsilon_0)$,}
	\Lambda_{1-\varepsilon,p} := \sup \left\{\|u\|_{L^p(\Omega)}^p:  u \in G ~~\mathcal{E}^{\varepsilon,p}(u,u) =1 \right\}.
	\end{align*}
Here we  adopt the notation 
	\begin{align*}
	\mathcal{E}_\Omega^{\varepsilon,p}(u,u) = \iint\limits_{\Omega \Omega} |u(x)-u(y)|^p\nu_\eps(x-y)\, dx\, dy.
	\end{align*}
Then the  map $\varepsilon\mapsto \Lambda_{1-\varepsilon,p} $ with $\eps\in (0, \eps_0)$,  is  bounded and we have 
	\begin{align*}
	\lim_{\varepsilon\to 0^+} \Lambda_{1-\varepsilon,p}= \Lambda_{1,p}K^{-1}_{d,p}.
	\end{align*}
\end{theorem} 

\vspace{2mm}
\begin{proof} The boundedness comes from Theorem \ref{thm:robust-poincare}. For every $\varepsilon\in (0,\eps_0)$ and every $u\in G$ we have 
	\begin{align*} 
	\int_{\Omega}|u(x)|^pdx\leq \Lambda_{1-\varepsilon,p} \iint\limits_{\Omega\Omega} |u(x)-u(y)|^p\nu_\eps(x-y)\, dx\, dy.
	\end{align*}
	According to  Theorem \ref{thm:charact-w1p} and Theorem \ref{thm:charact-BV}, letting   $\varepsilon\to 0^+$ implies
\begin{align*} 
\int_{\Omega}|u(x)|^pdx\leq K_{d,p}\liminf_{\varepsilon \to 0}\Lambda_{1-\varepsilon,p} \int\limits_{\Omega} |\nabla u(x)|^p\, \d x.
\end{align*}
	
	\noindent It follows that 
	\begin{align*}
	\Lambda_{1,p} \leq  K_{d,p}\liminf_{\varepsilon\to 0^+} \Lambda_{1-\varepsilon,p}.
	\end{align*}
	On the other hand, if we  fix $\eta\in (0, ~\Lambda_{1,p} K^{-1}_{d,p} )$ then, for each   $\varepsilon>0$  there exists $u_\varepsilon \in G$ such that
	\begin{align}\label{eq:sup-characto}
	\Lambda_{1-\varepsilon,p} -\eta< \|u_\varepsilon\|^p_{L^p(\Omega)} \qquad\hbox{and} \qquad \mathcal{E}_\Omega^{\varepsilon,p}(u_\varepsilon,u_\varepsilon) =1.
	\end{align}
	In virtue of Theorem \ref{thm:robust-poincare}  the family $(u_{\varepsilon})_\varepsilon$ is bounded in $L^p(\Omega)$ since $\|u_\varepsilon\|^p_{L^p(\Omega)}\leq B \mathcal{E}^{\varepsilon,p}(u_\varepsilon,u_\varepsilon)= B$ with the constant $B$ is  independent of $\varepsilon$.  Passing through a subsequence we may assume  as for the proof of Theorem \ref{thm:robust-poincare},  that it  converges in  $L^p(\Omega)$ to some $u$ with $u\in W^{1,p}(\Omega)$ when  $p>1$ and  $u\in BV(\Omega)$ when $p=1$. Further, applying Theorem \ref{thm:liminf-BBM} and  Theorem \ref{thm:liminf-BV}  yields
	\begin{align*}
	\|\nabla u\|^p_{L^p(\Omega)}\leq K^{-1}_{d,p}\liminf \mathcal{E}_\Omega^{\varepsilon,p}(u_\varepsilon,u_\varepsilon) = K^{-1}_{d,p}.
	\end{align*}
	Subsequently, we have $u\in G$ since $G$ is closed in $L^p(\Omega)$. Thus, by applying  the Poincar\'e inequality \eqref{eq:poincare-local} to $u$ and  taking into account \eqref{eq:sup-characto}  we obtain the following 
	\begin{align*}
	\limsup_{\varepsilon\to 0^+} \Lambda_{1-\varepsilon,p} -\eta\leq \|u\|^p_{L^p(\Omega)}\leq  
	\Lambda_{1,p} \|\nabla u\|^p_{L^p(\Omega)}
	\leq  \Lambda_{1,p}  K^{-1}_{d,p}.
	\end{align*}
	Finally, letting $\eta\to0$, completes the proof  since  we reach the following inequality
	\begin{align*}
	\limsup_{\varepsilon\to 0^+} \Lambda_{1-\varepsilon,p} \leq  \Lambda_{1,p}  K^{-1}_{d,p}. 
	\end{align*}
\end{proof}

\vspace{1mm}

\noindent Analogously we also have the following  convergence of sharps constants. 
\begin{theorem}[\textbf{convergence of sharp constants II}]\label{thm:convergence-sharp-const-poinc-friedr} Let $\Omega\subset \R^d$ be open, bounded and connected Lipschitz. Let $1\leq p<\infty$.   Assume that the family $(\nu_\eps)_\eps$  satisfies the conditions \eqref{eq:normalized-inetrals} and \eqref{eq:concentration-property}. Let $\Lambda'_{1,p} $ be the sharp constant of the inequality \eqref{eq:poincare-local-friedrichs} and $\Lambda'_{1-\varepsilon,p} $ be the sharp constant of the inequality \eqref{eq:robust-poincare-friedrichs}, i.e., 
	\begin{align*}
	\Lambda'_{1,p} := \sup \left\{\|u\|_{L^p(\Omega)}^p:  u \in C_c^\infty(\Omega), ~~\|\nabla u\|_{L^p(\Omega)}=1 \right\}
	\intertext{and for each $\varepsilon\in (0,\varepsilon_0)$,}
	\Lambda'_{1-\varepsilon,p} := \sup \left\{\|u\|_{L^p(\Omega)}^p:  u \in C_c^\infty(\Omega), ~~\mathcal{E}^{\varepsilon,p}(u,u) =1 \right\}.
	\end{align*}
	Here we  adopt the notation 
	\begin{align*}
	\mathcal{E}^{\varepsilon,p}(u,u) = \iint\limits_{\R^d\R^d} |u(x)-u(y)|^p\nu_\eps(x-y)\, dx\, dy.
	\end{align*}
	Then the  map $\varepsilon\mapsto \Lambda'_{1-\varepsilon,p} $ with $\eps\in (0, \eps_0)$ is  bounded and we have 
	\begin{align*}
	\lim_{\varepsilon\to 0^+} \Lambda'_{1-\varepsilon,p}= \Lambda'_{1,p}K^{-1}_{d,p}.
	\end{align*}
\end{theorem} 

\vspace{2mm}

\noindent Next, we deal with a situation where the robustness holds globally in terms of the parameter $\eps>0$.

\begin{theorem}\label{thm:robust-poincare-bis}
	Assume $d\geq 2, 1\leq p<\infty$.  Let $\Omega\subset \R^d$ be a bounded connected Lipschitz open set.  Assume that the family $(\nu_\eps)_{\eps\in (0, \eps_*)}$  satisfies the conditions \eqref{eq:normalized-inetrals} and \eqref{eq:concentration-property}. In addition assume that for $R>0$ and  $\tau\in (0, \eps_*)$ there exists $\theta>0$ such that 
	\begin{align}\label{eq:lower-bound-general}
		\nu_\eps(h)\geq \theta,\quad\quad\text{for all $|h|\leq R$ and all $\eps\in (\tau, \eps_*)$}.
	\end{align}
	Then there exists $C= C( d, p, \Omega)$ independent of $\eps>0$ such  that for all  $\eps\in (0,\eps_*)$ and $u\in L^p(\Omega)$ 
	\begin{align*}
	\big\|u-\hbox{ $\fint_{\Omega }$}u\big\|^p_{L^p(\Omega)}\leq C\iil_{ \Omega \Omega}|u(x) -u(y)|^p\nu_\eps(x-y)\d y\d x.
	\end{align*}
	The same is true in dimension $d=1$ provided that the condition  \eqref{eq:one-dimension-condition} holds. 
\end{theorem}
\begin{proof}
	According to Theorem \ref{thm:robust-poincare}, there exist $\eps_0\in (0,\eps_*)$ and a constant $B= B(d,p,\Omega)$ such that  for all $u\in L^p(\Omega)$ and $\eps\in (0, \eps_0)$
		\begin{align*}
	\big\|u-\hbox{ $\fint_{\Omega }$}u\big\|^p_{L^p(\Omega)}\leq B\iil_{ \Omega \Omega}|u(x) -u(y)|^p\nu_\eps(x-y)\d y\d x.
	\end{align*}
		
\noindent The condition \eqref{eq:lower-bound-general} with $\tau = \varepsilon_0$ implies that $\nu_{\varepsilon}(x-y) \geq \theta$ for all $\eps\in (\eps_0, \eps_*)$ and all $|x-y|\leq R$,  where  $R=\operatorname{diam}(\Omega)$ is  the diameter of $\Omega$. This, together  with Jensen's inequality  yield
\begin{align*}
 \iint\limits_{\Omega \Omega} |u(x)-u(y)|^p\nu_{\varepsilon}(x-y)\, dx\, dy  
	&\geq  \theta   \iint\limits_{\Omega \Omega} |u(x)-u(y)|^p \, dx\, dy  \\
	&\geq  \theta|\Omega| \int_{\Omega } \big |u(x)-\hbox{ $\fint_{\Omega }$}u\big|^p  \, dx.  
\end{align*}
The proof is now complete. 
\end{proof}

\vspace{2mm}
\noindent As a consequence of Theorem \ref{thm:robust-poincare-bis} we have following. 
\begin{corollary}\label{cor:robust-poincare-fractional} Assume $d\geq 1,~1<p<\infty $. Let $\Omega\subset \R^d$ be open, bounded and connected with Lipschitz boundary. Then there exists a constant, $C=C(p,d,\Omega)>0$  depending only on $d,p$ and $ \Omega$  such that  for every $s\in(0 ,1)$ and  every $u\in L^p(\Omega)$ we have 
	\begin{align}\label{eq:robust-poincare-fractional}
\big\|u-\hbox{ $\fint_{\Omega }$}u\big\|^p_{L^p(\Omega)}\leq C(1-s)\iint\limits_{\Omega\Omega} \frac{|u(x)-u(y)|^p}{|x-y|^{d+sp}}\, dx\, dy. 
	\end{align}
\end{corollary}

\begin{remark}
It is worth pointing out that Corollary \ref{cor:robust-poincare-fractional} provides a  global  estimate in terms of the fractional order $s\in (0,1)$ . 
	In other words, the estimate holds true for all $s\in (0,1)$ and the estimating constant $C$ is  independent of $s$. In \cite{Ponce2004}, the author  established the inequality \eqref{eq:robust-poincare-fractional}, provided that $s>s_0$ for some $s_0$ depending on  $d,p$ and $ \Omega$. Our 
	present result fills the gap $0<s\leq s_0$. See also \cite{limiting-embedding-BBM, Brasco-Lindgren} for further applications. 
\end{remark}

\vspace{2mm}
\noindent As a consequence of Theorem \ref{thm:robust-poincare-friedrichs} we have following. 
\begin{corollary}\label{cor:robust-poincare-fractional-frid} Assume $d\geq 1,~1<p<\infty $. Let $\Omega\subset \R^d$ be a bounded connected Lipschitz open set. Then there exist $s_0\in (0,1)$ and a constant $C=C(p,d,\Omega)>0$  depending only on $d,p$ and $ \Omega$  such that  for every $s\in(s_0 ,1)$ and  every $u\in C_c^\infty(\R^d)$ with $\operatorname{supp}u \subset \Omega$, we have
	\begin{align}\label{eq:robust-poincare-fractional-fri}
\big\|u\big\|^p_{L^p(\Omega)}\leq C(1-s)\iint\limits_{\R^d\R^d} \frac{|u(x)-u(y)|^p}{|x-y|^{d+sp}}\, dx\, dy. 
	\end{align}
\end{corollary}

\section{Convergence of Hilbert spaces}

In order to establish the main goal of this chapter, we are forced to deal with the family of Hilbert spaces that vary according to limiting variables. Therefore, we need to understand certain aspects and notions of convergence in varying spaces. Here we have in mind concepts like the Gamma convergence, the Mosco convergence and the convergence in varying Hilbert spaces. Our exposition in this section closely follows \cite{KS03}. Throughout this section $(H,\big(\cdot, \cdot \big)_H)$ and $\big(H_n,\big(\cdot, \cdot \big)_{H_n}\big)_n$ is a family of separable Hilbert spaces over $\R$. 
\medskip

\begin{definition}\label{def:convergence-hilbert-space}
	We say that a sequence of Hilbert spaces $(H_n)_n$ converges to a Hilbert space $H$ if $H$ has dense a subspace $\mathscr{C}$ and for each $n\geq 1$ there is a linear map $\Phi_n: \mathscr{C}\to H_n$ such that for all $u\in \mathscr{C}$
	\begin{align}\label{eq:hilbertspace-conv}
	\lim_{n\to \infty}\|\Phi_n u\|_{H_n}= \|u\|_{H}.
	\end{align}
This is equivalent to say that for all $u,v\in \mathscr{C}$,
	\begin{align*}
	\lim_{n\to \infty} (\Phi_nu, \Phi_nv )_{H_n} =(u,v)_{H}.
	\end{align*}
	Indeed, it suffices to write
	\begin{align*}
	(\Phi_nu, \Phi_nv )_{H_n} =\frac{1}{4} \Big( \|\Phi_n(u+v\|^2_{H_n}
	- \|\Phi_n(u-v)\|^2_{H_n}\Big) \quad\text{and}\quad
	(u,v)_{H}= \frac{1}{4} \Big( \|u+v\|^2_{H}- \|u-v\|^2_{H}\Big).
	\end{align*}
	In practice it is common to take $\Phi_n$ be the inclusion operator when $H\subset H_{n+1}\subset H_n$ for all $n\geq 1$.
\end{definition}
\begin{example}
	We let  $H=H^1(\Omega)$ be the standard Sobolev space and $H_n= H^{\alpha_n/2}(\Omega)$ be the fractional Sobolev space where $\Omega\subset \R^d$ is a Lipschitz open set and $(\alpha_n)_n$ is any sequence, $0<\alpha_n<2$ tending to $2$. It results as a consequence of Corollary \ref{cor:BBM-fractional} that $(H_n)_n$ converges to $H$. 
\end{example}

\noindent In what follows, it is assumed in this section that $H$ has dense a subspace $\mathscr{C}$ such that \eqref{eq:hilbertspace-conv} holds true. 

\begin{definition}\label{def:varying-convergence}
	Let $(u_n)_n$ be a sequence such that $u_n\in H_n$ and $u\in H$. 
	
	\begin{enumerate}[$(i)$]
		\item 	 We say that $(u_n)_n$ strongly converges (or simply converges ) to $u$ if there exists another sequence $(\widetilde{u}_m)_m\subset \mathscr{C}$ such that 
		\begin{align*}
		\lim_{m\to \infty}\|\widetilde{u}_m-u\|_{H}= 0\qquad\text{and}\qquad \lim_{m\to \infty}\limsup_{n\to \infty}\|\Phi_n\widetilde{u}_m-u_n\|_{H_n}= 0.
		\end{align*}
		Note that the first limit says that $(\widetilde{u}_m)_m$ tends to $u$ in the topology of $H$.
		\item 	 We say that $(u_n)_n$ converges in the weak sense (\emph{weakly}) to $u$ if 	for every sequence $(v_n)_n$ with  $v_n\in H_n$, strongly converging to $v\in H$ we have
		\begin{align*}
		\lim_{n\to \infty}\big(u_n,v_n\big)_{H_n}= \big(u,v\big)_H.
		\end{align*}
	 
	\end{enumerate}
	
\end{definition}

\noindent Let us now visit some related properties. 
\begin{proposition}\label{prop:varying-hilbert-property}
	Assume that the sequence of Hilbert spaces $(H_n)_n$	converges to a Hilbert space H. 
	Let $(u_n)_n$ and $(v_n)_n$ be sequence such that $u_n, v_n\in H_n$ and $u,v\in H$. 
	
	\begin{enumerate}[$(i)$]
		\item 	$(u_n)_n$ converges to $0\in H$ if and only if $\|u_n\|_{H_n} \xrightarrow{n\to \infty}0$.
		\item If $(u_n)_n$ converges to $u$ then $\|u_n\|_{H_n}\xrightarrow{n\to \infty}\|u\|_{H}$.
		\item If $(u_n)_n$ converges to $u$ and $(v_n)_n$ converges to $v$ then $(u_n+\lambda v_n)_n$ converges to $u+\lambda v$ for $\lambda\in \R$.
		\item If $(u_n)_n$ converges to $u$ and $(v_n)_n$ converges to $v$ then $(u_n,v_n)_{H_n}\xrightarrow{n\to \infty}(u,v)_{H}$.
		\item $\|u_n-v_n\|_{H_n}\xrightarrow{n\to \infty}0$ and $(u_n)_n$ converges to $u$ then $(v_n)_n$ converges to $u$.
		\item For any $w\in H$ there is a sequence $(w_n)_n$ with $w_n\in H_n$ such that $(w_n)_n$ converges to $w$. 
		\item If $(u_n)_n$ converges to $u$ and $(v_n)_n$ converges to $u$ then $\|u_n-v_n\|_{H_n}\xrightarrow{n\to \infty}0$.
	\end{enumerate}
\end{proposition}

\medspace

\begin{proof}
	Throughout, we assume that the sequence $(\widetilde{u}_m)_m\subset \mathscr{C}$ is such that $\|\widetilde{u}_m-u\|_{H}\xrightarrow{m\to \infty}0$ and
	\begin{align*}
	\lim_{m\to \infty}\limsup_{n\to \infty}\|\Phi_n\widetilde{u}_m-u_n\|_{H_n}= 0.
	\end{align*}
	
	\noindent $(i)$ If $u_n\to 0$ then $\|\widetilde{u}_m\|_{H}\xrightarrow{m\to \infty}0$ and since $\|\Phi_n\widetilde{u}_m\|_{H_n} \xrightarrow{n\to \infty} \|\widetilde{u}_m\|_{H}$ we find that 
	\begin{align*}
	\lim_{n\to \infty}\|u_n\|_{H_n} &\leq \lim_{m\to \infty}\limsup_{n\to \infty}\|\Phi_n\widetilde{u}_m-u_n\|_{H_n}+ \lim_{m\to \infty}\limsup_{n\to \infty}\|\Phi_n\widetilde{u}_m\|_{H_n}\\
	&=  \lim_{m\to \infty}\|\widetilde{u}_m\|_{H}=0.
	\end{align*}
	Conversely, assume $\|u_n\|_{H_n}\to 0$, by density we can approximate $0\in H$ i.e there exists a sequence $(u'_m)_m\subset \mathscr{C}$ such that $\|u'_m\|_H\to0$. We know that $\|\Phi_n u'_m\|_{H_n} \xrightarrow{n\to \infty} \|u'_m\|_{H}$ for each $m$. Hence the conclusion follows since
	\begin{align*}
	\lim_{m\to \infty}\limsup_{n\to \infty}\|\Phi_n u'_m-u_n\|_{H_n}
	&\leq 
	\lim_{m\to \infty}\limsup_{n\to \infty}\|\Phi_n u'_m\|_{H_n}+ \lim_{n\to \infty}\|u_n\|_{H_n}\\
	&= \lim_{m\to \infty}\|u'_m\|_{H}+ \lim_{n\to \infty}\|u_n\|_{H_n}=0.
	\end{align*}
	
	\noindent $(ii)$ For all $n,m$ by triangle inequality it is effortless to see that
	\begin{align*}
	\|\Phi_n\widetilde{u}_m\|_{H_n}- \|\Phi_n\widetilde{u}_m-u_n\|_{H_n} \leq \|u_n\|_{H_n}\leq \|\Phi_n\widetilde{u}_m\|_{H_n}+ \|\Phi_n\widetilde{u}_m-u_n\|_{H_n}
	\end{align*}
	from which we get 
	\begin{align*}
	\Big|\|u_n\|_{H_n}-\|u\|_{H}\Big|\leq\|\Phi_n\widetilde{u}_m-u_n\|_{H_n}+ \Big|\|\Phi_n\widetilde{u}_m\|_{H_n}-\|u\|_{H}\Big|.
	\end{align*}
	Hence the conclusion holds since $\|\Phi_n\widetilde{u}_m\|_{H_n} \xrightarrow{n\to \infty} \|\widetilde{u}_m\|_{H}$, $\|\widetilde{u}_m-u\|_{H}\xrightarrow{m\to \infty}0$ and 
	\begin{align*}
	\lim_{m\to \infty}\limsup_{n\to \infty}\Big|\|u_n\|_{H_n}-\|u\|_{H}\Big|
	&\leq \lim_{m\to \infty}\limsup_{n\to \infty}\|\Phi_n\widetilde{u}_m-u_n\|_{H_n}+ \lim_{m\to \infty}\limsup_{n\to \infty}\Big|\|\Phi_n\widetilde{u}_m\|_{H_n}-\|u\|_{H}\Big|\\
	&=\lim_{m\to \infty}\Big|\|\widetilde{u}_m\|_{H}-\|u\|_{H}\Big|=0.
	\end{align*}
	
	\noindent $(iii)$ This is a straightforward consequence of the definition.
	
	\vspace{2mm}
	\noindent $(iv)$ From $(ii)$ and $(iii)$ the result is obtained as follows
	\begin{align*}
	\lim_{n\to \infty} (u_n,v_n)_{H_n} 
	&= \lim_{n\to \infty}\frac{1}{4} \Big( \|u_n+v_n\|^2_{H_n}-\|u_n-v_n\|^2_{H_n}\Big)\\
	&= \frac{1}{4} \Big( \|u+v\|^2_{H}-\|u-v\|^2_{H}\Big)= (u,v)_{H}.
	\end{align*}

	\noindent $(v)$ Since $\|u_n-v_n\|_{H_n}\xrightarrow{n\to \infty}0$ and $(u_n)_n$ converges to $u$ then $\|\widetilde{u}_m-u\|_{H}\xrightarrow{m\to \infty}0$ and
	\begin{align*}
	\lim_{m\to \infty}\limsup_{n\to \infty}\|\Phi_n\widetilde{u}_m-v_n\|_{H_n}\leq \lim_{m\to \infty}\limsup_{n\to \infty}\|\Phi_n\widetilde{u}_m-u_n\|_{H_n}+ \lim_{m\to \infty}\limsup_{n\to \infty}\|u_n-v_n\|_{H_n}= 0
	\end{align*}
	which shows that $(v_n)_n$ converges to $u$. 
	
		\medskip 
		
	\noindent $(vi)$ By density there exists $\widetilde{u}_m\in \mathscr{C}$ with $ \|\widetilde{u}_m-u\|_{H}\to 0$ it suffices to take 
	$u_n= \Phi_n\widetilde{u}_n\in H_n$ such that $\|\Phi_n\widetilde{u}_n-u_n\|_{H_n}=0$ for every $n\geq1$.
	
		\medskip
		
	\noindent $(vii)$ Note that 
	\begin{align*}
	\|u_n-v_n\|_{H_n}\leq \|\Phi_n\widetilde{u}_m-u_n\|_{H_n}+ \|\Phi_n\widetilde{v}_m-v_n\|_{H_n} +\|\Phi_n\widetilde{u}_m-\Phi_n\widetilde{v}_m\|_{H_n}.
	\end{align*}
	The claim readily follows since $\|\widetilde{u}_m-u\|_{H}\xrightarrow{m\to \infty}0$, $\|\widetilde{v}_m-u\|_{H}\xrightarrow{m\to \infty}0$ and we have 
	\begin{align*}
	\lim_{m\to \infty}\limsup_{n\to \infty}\|\Phi_n\widetilde{u}_m-\Phi_n\widetilde{v}_m\|_{H_n}=\lim_{m\to \infty}\limsup_{n\to \infty} \|\Phi_n(\widetilde{u}_m-\widetilde{v}_m)\|_{H_n} = \lim_{m\to \infty}\|\widetilde{u}_m-\widetilde{v}_m\|_{H}=0.
	\end{align*}
	
\end{proof}

\bigskip

\begin{lemma}
	Assume that  the sequence of Hilbert spaces $(H_n)_n$	converges to the Hilbert space H. A sequence $(u_n)_n$ with $u_n\in H_n$ converges to $u$ if and only if $(u_n,w_n)_{H_n}\xrightarrow{n\to \infty}(u,w)_{H}$ for every sequence 
	$(w_n)_n$, $w_n \in H_n$  converging weakly  to $w\in H$.
\end{lemma}
\begin{proof}
	The forward implication is patently a consequence of the definition of the weak convergence. Reciprocally assume $(u_n,w_n)_{H_n}\xrightarrow{n\to \infty}(u,w)_{H}$ for every sequence 
	$(w_n)_n$, $w_n \in H_n$  converging weakly  to $w\in H$. From Proposition \ref{prop:varying-hilbert-property} $(iv)$ we know that strong convergence implies the weak convergence, then choosing $(w_n)_n$ strongly converging to $w$ implies that $(u_n)$ is in particular  converging weakly  to $u$. Thus taking $w_n=u_n$ implies $\|u_n\|_{H_n} \xrightarrow{n\to \infty}\|u\|_{H}.$ Next by density we have $\|\widetilde{u}_m -u\|_{H}\xrightarrow{m\to \infty}0$ with $\widetilde{u}_m\in \mathscr{C}$. To conclude, it suffices to show that 
	\begin{align*}
	\lim_{m\to \infty}\limsup_{n\to \infty}\|\Phi_n\widetilde{u}_m-u_n\|_{H_n}=0.
	\end{align*}
	To do this, let us fix $m\geq 1$. We claim that $(\Phi_n\widetilde{u}_m)_n$ is  converging weakly  to 
$\widetilde{u}_m$. Indeed, for a sequence $(w_n)_n$, $w_n \in H_n$ strongly converging to $w\in H$ there exists $(\widetilde{w}_j)_j $ such that $\|\widetilde{w}_j-w\|_{H}\xrightarrow{j\to \infty}0$ and 
	\begin{align*}
	\lim_{j\to \infty}\limsup_{n\to \infty}\|\Phi_n\widetilde{w}_j-w_n\|_{H_n}= 0.
	\end{align*}
	We have
	\begin{align*}
	\big|(\Phi_n\widetilde{u}_m, w_n )_{H_n} - (\widetilde{u}_m, w )_{H} \big| 
	&\leq \big|(\Phi_n\widetilde{u}_m, \Phi_n\widetilde{w}_j)_{H_n} - (\widetilde{u}_m, w )_{H} \big| + \big|(\Phi_n\widetilde{u}_m, w_n-\Phi_n\widetilde{w}_j )_{H_n} \big| \\
	&\leq \big|(\Phi_n\widetilde{u}_m, \Phi_n\widetilde{w}_j)_{H_n} - (\widetilde{u}_m, w )_{H} \big| + \|\Phi_n\widetilde{u}_m\|_{H_n} \|w_n-\Phi_n\widetilde{w}_j \|_{H_n}.
	\end{align*}
	The claim follows once we show that $(\Phi_n\widetilde{u}_m, w_n )_{H_n} \xrightarrow{n\to \infty} (\widetilde{u}_m, w )_{H}$. Indeed, let us see that 
	\begin{align*}
	\lim_{j\to \infty}\limsup_{n\to \infty} (\Phi_n\widetilde{u}_m, \Phi_n\widetilde{w}_j)_{H_n} =\lim_{j\to \infty} (\widetilde{u}_m, \widetilde{w}_j)_{H} = (\widetilde{u}_m, w )_{H}. 
	\end{align*} 
	As $(\Phi_n\widetilde{u}_m)_n$ is  converging weakly  to $\widetilde{u}_m$ then in particular we have $(u_n,\Phi_n\widetilde{u}_m)_{H_n}\xrightarrow{n\to \infty}(u,\widetilde{u}_m)_{H}$.
	Finally, combining the foregoing steps yields
	\begin{align*}
	\lim_{m\to \infty}\limsup_{n\to \infty}\|\Phi_n\widetilde{u}_m-u_n\|^2_{H_n}&=
	\lim_{m\to \infty}\limsup_{n\to \infty} \Big(\|\Phi_n\widetilde{u}_m\|^2_{H_n}+\|u_n\|^2_{H_n} -2 (u_n,\Phi_n\widetilde{u}_m)_{H_n}\Big)\\
	&= \lim_{m\to \infty}\Big(\|\widetilde{u}_m\|^2_{H_n}+\|u\|^2_{H} -2 (u,\widetilde{u}_m)_{H}\Big)=0.
	\end{align*} 
\end{proof} 

\bigskip

\noindent Now we derive what can be viewed as a  weakly sequential compactness result in connection with the spaces $(H_n)_n$ and $H$. We believe that our proof here, although inspired by \cite{KS03}, is more comprehensible and uses rudimentary arguments.
\begin{theorem}\cite[Lemma 2.2 $\&$ 2.3]{KS03}\label{thm:varying-weak-convergence}
	Assume the sequence of separable Hilbert spaces $(H_n)_n$ converges to a separable Hilbert space H. Let $(u_n)_n$ be a sequence with $u_n\in H_n$. 
	\begin{enumerate}[$(i)$]
		\item If $(u_n)_n$ weakly converges to some $u\in H$ then the sequence $(\|u_n\|_{H_n})_n$ is bounded and we have $$\|u\|_H\leq \liminf\limits_{n\to\infty} \|u_n\|_{H_n}.$$ 
		\item Conversely, if the sequence $(\|u_n\|_{H_n})_n$ is bounded then there is a subsequence $(u_{n_k})_k$  converging weakly  to some $u\in H$. 
	\end{enumerate}
\end{theorem}

\medspace

\begin{proof} $(i)$ Assume $\sup\limits_{n\geq1} \|u_n\|_{H_n}=\infty$ then one is able to construct a subsequence $(n_k)_k$ such that $\|u_{n_k}\|_{H_{n_k}}\geq 2^k$ for all $k\geq1$. We have the right to set $v_k = \frac{u_{n_k}}{2^k \|u_{n_k}\|_{H_{n_k}}}$. It clearly appears that $\|v_{k}\|_{H_{n_k}} = 2^{-k}\xrightarrow{k\to \infty}0$. In virtue of Proposition \ref{prop:varying-hilbert-property}, $(v_k)_k$ converges to $0\in H$. On the other hand, it is also clear that the subsequence $(u_{n_k})_k$ weakly converges to $u$. By the definition of weak convergence it follows that $(u_{n_k},v_k)_{H_{n_k}}\xrightarrow{k\to \infty}0$ which contradictory to the fact that 
	\begin{align*}
	(u_{n_k},v_k)_{H_{n_k}}= \Big(u_{n_k},\frac{u_{n_k}}{2^k \|u_{n_k}\|_{H_{n_k}}}\Big)_{H_{n_k}} = \frac{\|u_{n_k}\|_{H_{n_k}}}{2^k }\geq 1.
	\end{align*}
	In conclusion, our initial assumption was wrong and hence $(\|u_n\|_{H_n})_n$ must be bounded. 
	
	\medspace
	
	\noindent Let us prove the liminf inequality. By invoking once more Proposition \ref{prop:varying-hilbert-property}, there exists $v_n\in H_n$ such that $(v_n)_n$ converging to $u$ and we also have $\|v_n\|_{H_n} \xrightarrow{n\to \infty}\|u\|_{H}$. By definition of weak convergence it follows that $(u_n,v_n)_{H_{n}}\xrightarrow{n\to \infty}(u,u)_{H}$. For $\varepsilon>0$ and certain $j_0\geq 1$ we have $\|v_n\|_{H_n}\leq \|u\|_{H}+\varepsilon$ for all $n\geq N$. Thus by Cauchy Schwarz inequality 
	\begin{align*}
	|(u_n,v_n)_{H_{n}}|\leq \|u_n\|_{H_n} \|v_n\|_{H_n} \leq \|u_n\|_{H_n} (\|u\|_{H} +\varepsilon)\quad\text{for all}\quad n\geq N.
	\end{align*} 
	Passing to the liminf as $n \to \infty $ and letting $\varepsilon\to0$ after, yields
	\begin{align*}
	\|u\|_{H} \leq \liminf_{n\to\infty}\|u_n\|_{H_n}.
	\end{align*} 
	
	\medskip
	
	\noindent $(ii)$ Assume  that $(\|u_n\|_{H_n})_n$ is bounded. Let $(\varphi_k)_k$ be an orthonormal basis of $H$. By density,
	for each $k\geq 1$ there exists a sequence $(\varphi_{k,m})_m\subset \mathscr{C}$ such that
	\begin{align}\label{eq:hpik-approx}
	\|\varphi_{k,m}-\varphi_k\|_{H}\leq 2^{-m}\quad\text{for every }\quad m\geq 1.
	\end{align}
	Observing that $ \lim\limits_{m\to \infty}\lim\limits_{n\to \infty}\|\Phi_n\varphi_{k,m}\|_{H_n}= \lim\limits_{m\to \infty}\|\varphi_{k,m}\|_{H} =1$ it is not difficult to obtain that for each $k\geq1$, the doubling sequence of real numbers $\big((u_n, \Phi_n\varphi_{k,m})_{H_n}\big)_{n,m}$ is bounded in $\R$. We may find a cluster point $a_k$, two suitable subsequences $(n_p)_p$ and $( m_j)_j$ and another sequence $(a_k(m_j))_j$ obeying the following rule: if we set $ a_k(n_p,m_j): =(u_{n_p} \Phi_{n_p}\varphi_{k,m_j})_{H_{n_p}},$ then
	\begin{align}
	a_k(m_j)= \lim_{p\to \infty} a_k(n_p, m_j) \quad\text{and }
	\quad \Big| a_k(m_j)-a_k\Big|\leq 2^{-j} . \label{eq:exists-ak}
	\end{align}
	In particular, we have 
	\begin{align*}
	a_k= \lim_{j\to \infty}\lim_{p\to \infty}a_k(n_p,m_j) =\lim_{j\to \infty}\lim_{p\to \infty}(u_{n_p}, \Phi_{n_p}\varphi_{k,m_j})_{H_{n_p}} \quad\text{for every }\quad k\geq 1.
	\end{align*}
	\noindent Let us  show the convergence of the series $\sum\limits_{k= 1}^\infty|a_k|^2$. Since $\|\varphi_{k,m}-\varphi_k\|_{H}\xrightarrow{m\to \infty}0$ it follows that
	\begin{align*}
	\lim_{m\to \infty}\lim_{n\to \infty} (\Phi_n\varphi_{k,m}, \Phi_n\varphi_{i,m})_{H_n} =\lim_{m\to \infty} (\varphi_{k,m}, \varphi_{i,m})_{H} =(\varphi_{k},\varphi_{i})_{H} = \delta_{k,i}. 
	\end{align*} 
	Further, recalling that $a_k(n_p,m_j): =(u_{n_p}, \Phi_{n_p}\varphi_{k,m_j})_{H_{n_p}}$ the following holds for every $\ell\in \mathbb{N}$.
	\begin{align*}
	0&\leq\lim_{j\to \infty}\liminf_{p\to \infty} \Big\| u_{n_p} - \sum_{k= 1}^\ell a_k(n_p,m_j)\Phi_{n_p}\varphi_{k,m_j} \Big\|^2_{H_{n_p}}\\
	& =\lim_{j\to \infty}\liminf_{p\to \infty} \left[\| u_{n_p} \|^2_{H_{n_p}} \hspace{-2ex} - 2\sum_{k= 1}^\ell| a_k(n_p,m_j)|^2  + \hspace{-0.6ex}
	\sum_{k= 1}^\ell\sum_{i= 1}^\ell a_k(n_p,m_j)a_i(n_p,m_j) (\Phi_{n_p}\varphi_{k,m_j}, \Phi_{n_p}\varphi_{i,m_j})_{H_{n_p}} \right]\\
	&= \liminf_{p\to \infty}\| u_{n_p} \|^2_{H_{n_p}} - 2\sum_{k= 1}^\ell| a_k|^2 
	+ \sum_{k= 1}^\ell\sum_{i= 1}^\ell a_k a_i \delta_{k,i}\\
	&=  \liminf_{p\to \infty}\| u_{n_p} \|^2_{H_{n_p}} - \sum_{k= 1}^\ell| a_k|^2 . 
	\end{align*}
	As a result, if we put $u:=\sum\limits_{k= 1}^\infty a_k\varphi_k$ then it is clear that $u\in H$ since letting $\ell\to \infty$, the above implies 
	\begin{align*}
	\|u\|^2_H=\sum_{k= 1}^\infty| a_k|^2\leq \liminf_{p\to \infty}\| u_{n_p} \|^2_{H_{n_p}}<\infty.
	\end{align*}
	
	\noindent Next, we show that $(u_{n_p})_p$ weakly converges to $u$. Let $(v_n)_n$, with $v_n\in H_n$ be a sequence strongly converging to some $v\in H$. That is there exists $(\widetilde{v}_m)_m\subset \mathscr{C}$ such that 
	\begin{align*}
	\lim_{m\to \infty}\|\widetilde{v}_m-v\|_{H}= 0\quad\text{and}\quad \lim_{m\to \infty}\limsup_{n\to \infty}\|\Phi_n\widetilde{v}_m-v_n\|_{H_n}= 0. 
	\end{align*}
	\noindent On the other hand, since $(\varphi_k)_k$ is an orthonormal basis of $H$ we know that 
	\begin{align*}
	v=\sum_{k= 1}^\infty b_k\varphi_k\quad\text{and}\quad \|v\|^2_H= \sum_{k= 1}^\infty |b_k|^2\quad\text{with}\quad b_k =(u, v)_{H}.
	\end{align*}
	To conclude, we show that $(u_{n_p},v_{n_p})_{H_{n_p}}\xrightarrow{p\to \infty} (u, v)_{H}.$ We put $ v^*_j= \sum\limits_{k= 1}^j b_k\varphi_{k,m_j}\in \mathscr{C}. $ From  \eqref{eq:wapprox} we obtain
	\begin{align*}
	\| v^*_j-v \|^2_H
	&= \Big\|\sum_{k= 1}^jb_k(\varphi_{k,m_j}-\varphi_{k})\Big\|^2_H+\sum_{k= j+1}^\infty |b_k|^2\\
	&\leq \Big(\sum_{k= 1}^j |b_k|^2\Big)\Big(\sum_{k= 1}^j\|\varphi_{k,m_j}-\varphi_{k}\|^2_H\Big)+ \sum_{k= j+1}^\infty |b_k|^2\\
	&\leq j4^{-m_j}\|v\|^2_H + \sum_{k= j+1}^\infty |b_k|^2\xrightarrow{j\to \infty}0. 
	\end{align*}
	It clearly follows that $\| v^*_j-\widetilde{v}_j \|^2_H\xrightarrow{j\to \infty}0 $ so that we have 
	\begin{align*}
	\lim_{j\to \infty}\lim_{p\to \infty}\Big| (u_{n_p}, \Phi_{n_p}v^*_{j}-\Phi_{n_p}\widetilde{v}_{j})_{H_{n_p}}\Big| \leq C\lim_{j\to \infty}\lim_{p\to \infty}\|\Phi_{n_p}(\widetilde{v}_{j}-v^*_{j})\|_{H_{n_p}}= C\lim_{j\to \infty} \|v^*_{j}-\widetilde{v}_{j}\|_{H} =0.
	\end{align*}
	Meanwhile recalling that $(u,v)_H = \sum\limits_{k= 1}^\infty a_kb_k$ using the above procedure we find that 
	\begin{align*}
	\lim_{j\to \infty}\lim_{p\to \infty} \Big|(u_{n_p},\Phi_{n_p} v^*_{j})_{H_{n_p}} -(u,v)_{H} \Big| 
	&\leq \lim_{j\to \infty} \sum_{k= 1}^j |b_k|\lim_{p\to \infty} |( u_{n_p},\Phi_{n_p} \varphi_{k,m_j})_{H_{n_p}} -a_k|+ \Big| \sum_{k= j+1}^\infty a_kb_k\Big|\\
	&= \lim_{j\to \infty} \sum_{k= 1}^j |b_k| |a_k(m_j) -a_k|+ \Big| \sum_{k= j+1}^\infty a_kb_k\Big|\\
	&\leq \lim_{j\to \infty} j2^{-j} \|v\|_H + \Big| \sum_{k= j+1}^\infty a_kb_k\Big|=0.
	\end{align*}

	\noindent Finally putting together the previous we obtain 
	\begin{align*}
	\limsup_{p\to \infty} \Big|(u_{n_p},v_{n_p})_{H_{n_p}} -(u,v)_{H} \Big| &\leq 
	\lim_{j\to \infty}\limsup_{p\to \infty} \Big| (u_{n_p},\Phi_{n_p} v^*_j)_{H_{n_p}} -(u,v)_{H} \Big| \\&+ \lim_{j\to \infty}\limsup_{p\to \infty} \Big|(u_{n_p},\Phi_{n_p}(v^*_{j}-\widetilde{v}_{j}))_{H_{n_p}} \Big|\\&
	+ \lim_{j\to \infty}\limsup_{p\to \infty}\Big|(u_{n_p},\Phi_{n_p} \widetilde{v}_{j} -v_{n_p})_{H_{n_p}} \Big|=0.
	\end{align*}
	%
\end{proof}

\section{Mosco convergence of nonlocal to local quadratic forms}
Here we establish the convergence in the sense of Mosco of a certain class of nonlocal quadratic
forms to some elliptic forms of gradient type. We essentially deal with Dirichlet forms. 
\begin{definition} Let $(X,\mathcal{B}, m)$ be a $\sigma$-finite measure space, where is $X$ is a locally compact separable metric space, $\mathcal{B}$ is the Borel set of $X$ and $m$ is a Radon measure on $X$ with full support. A Dirichlet form on $L^2(X,m)$ is a couple $(\mathcal{E} , \mathcal{D}(\mathcal{E}))$ satisfying $D_1-D_4$. Moreover $(\mathcal{E} , \mathcal{D}(\mathcal{E}))$  is called to be a  regular Dirichlet form on $L^2(X,m)$ if in addition $D_5$ holds true. 
	\begin{enumerate}[$D_1.$]
	\item  $\mathcal{D}(\mathcal{E})$ is dense in $L^2(X,m)$.
	\item $\mathcal{E}: \mathcal{D}(\mathcal{E})\times \mathcal{D}(\mathcal{E})\to \R$ is a positive, (i.e.  $\mathcal{E}(u,u)\geq0$ for all $u\in \mathcal{D}(\mathcal{E})$)  symmetric bilinear form. 
	
	\item $(\mathcal{E} , \mathcal{D}(\mathcal{E}))$ is closed if  the space $( \mathcal{D}(\mathcal{E}),
	 \|\cdot\|_{\mathcal{D}(\mathcal{E})}) $ is a Hilbert space. Here $\|\cdot\|_{\mathcal{D}(\mathcal{E})}$ represents the graph norm defined by $\|u\|^2_{\mathcal{D}(\mathcal{E})}:= \mathcal{E}_1(u,u)= \|u\|^2_{L^2(X,m) }+ \mathcal{E}(u,u).$

\item $(\mathcal{E} , \mathcal{D}(\mathcal{E}))$ is Markovian if  $u\in \mathcal{D}(\mathcal{E})$ and  $v= (u\lor 0)\land 1$ (called the normal contraction of $u$)  then $v\in  \mathcal{D}(\mathcal{E})$ and $\mathcal{E}(v,v)\leq \mathcal{E}(u,u)$. 
\item $(\mathcal{E} , \mathcal{D}(\mathcal{E}))$ is regular if it possesses a \textit{core} i.e. there is a subset $\mathcal{C}$ of $(\mathcal{D}(\mathcal{E})\cap C_c(X)$ such that $\mathcal{C}$ is dense in $( \mathcal{D}(\mathcal{E}), \|\cdot\|_{\mathcal{D}(\mathcal{E})})$ and $\mathcal{C}$ is dense in $( C_c(X), \|\cdot\|_{L^\infty(X)}).$
	\end{enumerate}
\end{definition}

\medskip

\noindent Referring to \cite{MR12, Fukushima-O-T}, any Dirichlet form is uniquely associated to a Generator $(A, \mathcal{D}(A))$, a strong  semigroup of contractions $(T_t)_{t\geq0}$ and a resolvent $(G_\lambda)_{\lambda>0}$ respectively characterized as follows. 

\noindent $\bullet$ Generator
\begin{enumerate}[$G_1.$]
	\item$ (A, \mathcal{D}(A))$ is a negative self-adjoint operator whose domain  $\mathcal{D}(A)$ is dense in $L^2(X,m)$, 
	\item $ \mathcal{D}(\mathcal{E})= \mathcal{D}(\sqrt{-A})$ and $ \mathcal{E}(u,v) = (\sqrt{-A}u, \sqrt{-A}v)_{L^2(X,m)}$ for all $u,v \in \mathcal{D}(\mathcal{E})$. 
	\item $A$ is closed, i.e. $( \mathcal{D}(A), \|\cdot\|_{\mathcal{D}(A)}) $ is a Hilbert space where $\|u\|^2_{\mathcal{D}(A)}= \|u\|^2_{L^2(X,m) }+\|Au\|^2_{L^2(X,m) }.$
	
	\item$\mathcal{D}(A)\subset  \mathcal{D}(\mathcal{E})$ and $ \mathcal{E}(u,v) = (-Au, v)_{L^2(X,m)}$ for all $u\in \mathcal{D}(A),\, v\in \mathcal{D}(\mathcal{E})$. 
\end{enumerate}
\noindent  $\bullet$ Semigroup
\begin{enumerate}[$S_1.$]
	\item$(T_t)_{t>0}$  is a semigroup on $L^2(X,m)$: $T_t$ is linear on $L^2(X, m)$ and  $T_sT_t=T_{t+s}$ for all $s,t>0$.
	\item The contraction property holds, i.e. $\|T_tu\|_{L^2(X,m)}\leq \|u\|_{L^2(X,m)}$, $t>0 $, $u\in L^2(X,m)$. 
	\item $(T_t)_{t>0}$ is strongly continuous, i.e.  $\|T_tu-u\|_{L^2(X,m)}\xrightarrow{t\to 0}0$. 
	\item  $T_tu  = e^{-At} u,$  $t>0$, $u\in L^2(X,m)$. $Au= \lim\limits_{t\to 0}\frac{T_t u-u}{t}$, $u\in \mathcal{D}(A)$.
\end{enumerate}
 \noindent  $\bullet$ Resolvent
\begin{enumerate}[$R_1.$]
	\item$(G_\lambda)_{\lambda>0}$  is a resolvent: $G_\lambda$ is linear on $L^2(X, m)$ and  $G_\lambda- G_\mu + (\lambda-\mu)G_\lambda G_\mu=0$, $\lambda, \mu>0$.
	\item The contraction property holds, i.e. $\|\lambda G_\lambda u\|_{L^2(X,m)}\leq \|u\|_{L^2(X,m)}$, $\lambda>0 $, $u\in L^2(X,m)$. 
	\item $(G_\lambda)_{\lambda>0}$ is strongly continuous, i.e.  $\|\lambda G_\lambda u-u\|_{L^2(X,m)}\xrightarrow{\lambda\to \infty}0$. 
	\item $G_\lambda u = \int_0^\infty e^{-\lambda t} T_tu$  $\lambda>0$, $u\in L^2(X,m)$. $Au=  \lambda u-G_\lambda^{-1} u$, $u\in \mathcal{D}(A)$.
\end{enumerate}

\noindent Let us recall the notion of Mosco convergence and  Gamma convergence  of quadratic forms on $L^2$- spaces.

%
\begin{definition}[Mosco convergence and Gamma convergence ]\label{def:mosco}
	Assume  $(\mathcal{E}^n, \mathcal{D}(\mathcal{E}^n))_{n\in \mathbb{N}}$ and $(\mathcal{E}, \mathcal{D}(\mathcal{E}))$  are  quadratic forms with dense domains in $L^2(X,m)$. One says that $(\mathcal{E}^n, \mathcal{D}(\mathcal{E}^n))_{n\in \mathbb{N}}$ converges in $L^2(X,m)$  in the Mosco sense (resp. in the Gamma sense)  to $(\mathcal{E}, \mathcal{D}(\mathcal{E}))$ if the following two conditions $M_1$ and $M_2$ (resp. $M_1$ and $M'_2$) are satisfied.
	
	\medskip

	\noindent \textbf{$M_1.$ (limsup):} For every $u\in L^2(X,m)$ there exists a sequence $(u_n)_n$ in $ L^2(X,m)$ such that   $u_n\in  \mathcal{D}(\mathcal{E}^n)$,  $u_n\to u$ (read $u_n$ strongly converges to $u$) in  $ L^2(X,m)$ and 
	\[\limsup_{n\to \infty} \mathcal{E}^n(u_n,u_n) \leq \mathcal{E}(u,u). \]
\noindent	\textbf{$M_2.$ (liminf):} For every sequence, $(u_n)_n$   with   $u_n\in  \mathcal{D}(\mathcal{E}^n)$ and every $u\in \mathcal{D}(\mathcal{E})$ such that   $u_n \rightharpoonup u$ (read $u_n$ weakly converges to $u$) in  $ L^2(X,m)$ we have, 
	\[\mathcal{E}(u,u)\leq  \liminf_{n\to \infty}  \mathcal{E}^n(u_n,u_n).\] 

\noindent \textbf{$M'_2.$ (liminf):} For every sequence, $(u_n)_n$   with   $u_n\in  \mathcal{D}(\mathcal{E}^n)$ and every $u\in \mathcal{D}(\mathcal{E})$ such that   $u_n \to u$ (read $u_n$ strongly converges to $u$) in  $ L^2(X,m)$ we have, 
\[\mathcal{E}(u,u)\leq  \liminf_{n\to \infty}  \mathcal{E}^n(u_n,u_n).\] 

\end{definition}

\begin{remark}
	$(i)$ It is worth emphasizing  that combining the $\limsup$ and $\liminf$ conditions,   the  $\limsup$ condition is equivalent to the existence of a sequence $(u_n)_n$ in $ L^2(X,m)$ such that   $u_n\in  \mathcal{D}(\mathcal{E}^n)$,  $u_n\to u$  in $ L^2(X,m)$ and 
	\[\lim_{n\to \infty} \mathcal{E}^n(u_n,u_n)=\mathcal{E}(u,u). \]
	$(ii)$ It is clear that  the convergence in the sense of Mosco implies Gamma convergence. The converse is true provided that the asymptomatic compactness holds, see Proposition  \ref{prop:Mosoc-vs-Gamma} below. 
	
\noindent $(iii)$ We adopt the convention  that for a given quadratic form $\big(\mathcal{E}, \mathcal{D}(\mathcal{E})\big)$, we have 
	$\mathcal{E}(u,u)= \infty$ whenever $u \not\in  \mathcal{D}(\mathcal{E})$.
	
	\noindent $(iv)$ Note that the Mosco convergence on Banach spaces in general is defined in \cite{Tol10}. 
\end{remark}

\begin{proposition}[Mosco vs Gamma ]\label{prop:Mosoc-vs-Gamma}
	Assume $(\mathcal{E}^n, \mathcal{D}(\mathcal{E}^n))_{n}$ converges in Gamma sense to  $(\mathcal{E}, \mathcal{D}(\mathcal{E}))$. Then $(\mathcal{E}^n, \mathcal{D}(\mathcal{E}^n))_{n}$ converges in Mosco sense to  $(\mathcal{E}, \mathcal{D}(\mathcal{E}))$ if $(\mathcal{E}^n, \mathcal{D}(\mathcal{E}^n))_{n\in \mathbb{N}}$ is asymptotically compact, i.e.  for any sequence $(u_n)_n$ such that $u_n\in \mathcal{D}(\mathcal{E}^n)$ and $\liminf\limits_{n\to \infty} \big(\mathcal{E}^n(u_n, u_n)+ \|u_n\|_{L^2(X,m)} \big)<\infty$ has  a strongly convergent subsequence. 
\end{proposition}
\begin{proof}
The proof is immediate. 
\end{proof}

\noindent To a closed form $\mathcal{E}$  corresponds a semigroup  $(T_t)_t$, Generator $(G_\lambda)_\lambda$  and  a stochastic process $(X_t)_t$.  According to \cite{Fukushima-O-T}, if $(X_t)_t$ is a Markov process (resp. a Hunt process) then $\mathcal{E}$ is a Dirichlet form (resp. a regular Dirichlet form).  The Mosco convergence relates the convergence of Markov processes and the convergence of their corresponding Dirichlet forms. 

\begin{theorem}[\cite{KS03,Mos94,Kol06}]\label{thm:mosco-characterisation}
Let	$(\mathcal{E}^n, D(\mathcal{E}^n))$ and $(\mathcal{E}, D(\mathcal{E}))$ be closed Dirichlet forms.  The following are equivalent:
	\begin{enumerate}[$(i)$]
		\item $(\mathcal{E}^n, D(\mathcal{E}^n))_n$ Mosco converges to $(\mathcal{E}, D(\mathcal{E})).$
		\item $(G^n_\lambda)_n$ strongly converges to $G_\lambda$ for every $\lambda$.
		\item $(T^n_t)_n$ strongly converges to $T_t$ for every $t$.
	\end{enumerate}
	
	\noindent Moreover,  if $(\mathcal{E}^n, D(\mathcal{E}^n))_n$ Mosco converges to $(\mathcal{E}, D(\mathcal{E}))$ then  $(X^n)_n$  converges to $X$ in finite dimensional distribution. 
\end{theorem}

\medskip

\noindent Next, we explain for which sequences of nonlocal quadratic forms we can prove convergence to a classical local gradient form.  Let us recall our standing set-up in this section. Let $(\nu_\alpha)_{\alpha\in(0,2)}$ be a family of L\' evy radial functions  approximating the Dirac measure at the origin, i.e. for every $\alpha, \delta > 0$
	\begin{align}\label{eq:assumption-nu-alpha}
	\begin{split}
	\nu_\alpha\geq 0,\,\,\text{ is radial}, \quad \int_{\mathbb{R}^d}	(1\land |h|^2)\nu_\alpha (h)\d h=1, \quad \lim_{\alpha\to 2}\int_{|h|>\delta}	\nu_\alpha(h)\d h=0\,.
	\end{split}
	\end{align}                      
	Moreover, we assume that $h \mapsto \nu_\alpha(h)$ is almost decreasing, i.e., for some $c \geq 1$ and all $x,y$ with $|x| \leq |y|$ we have $	\nu_\alpha(y) \leq c \, \nu_\alpha(x)$. Note that all examples of $(\nu_\alpha)_\alpha$  from Section \ref{sec:approx-dirac-mass} apply here as well. Moreover it is important to stress that  $(\nu_\alpha)_{0<	\alpha<2}$ defined as above generalizes the  set-up of  \cite{FGKV19}.

\medskip

\noindent  Given $\alpha \in (0,2)$,$J^\alpha: \mathbb{R}^d\times \mathbb{R}^d \setminus \operatorname{diag} \to [0, \infty]$ and sufficiently smooth functions $u,v: \mathbb{R}^d\to \mathbb{R}$, define 
\begin{align}
\mathcal{E}^{\alpha}_{\Omega}(u,v) &=  \iil_{\Omega \Omega} \big(u(y)-u(x)\big) \big(v(y)-v(x)\big) J^\alpha(x,y)\d x \, \d y\,, \label{eq:inner-form} \\
\mathcal{E}^{\alpha}(u,v) &= \iil_{(\Omega^c\times \Omega^c)^c} \big(u(y)-u(x)\big) \big(v(y)-v(x) J^\alpha(x,y) \d x \, \d y. \label{eq:ext-form}
\end{align}

\noindent  Here, for sequence  $(J^\alpha)_{0<\alpha<2}$ of  positive symmetric kernels $J^\alpha: \mathbb{R}^d\times \mathbb{R}^d \setminus \operatorname{diag} \to [0, \infty]$ we set-up  the following conditions:

\medskip

\begin{itemize}
	\item[(E)] There exists a constant $\Lambda\geq 1$ such that for every $\alpha\in  (0,2)$ and all $x,y \in \mathbb{R}^d$, with $0<|x-y|\leq 1$ 
	\begin{align}\label{eq:elliptic-condition}\tag{$E$}
	\Lambda^{-1} \nu_\alpha (x-y) \leq J^\alpha(x,y) &\leq  \Lambda \nu_\alpha (x-y) .
	\end{align}
	\item[(L)] For every $\delta >0$
	\begin{align} \label{eq:integrability-condition}\tag{$L$}
	\lim_{\alpha \to 2^-}\sup_{x\in \mathbb{R}^d} \int_{|h| > \delta} J^\alpha(x,x+h)dh=0.
	\end{align}
		\item[(I)] For each $\alpha \in (0,2)$ the kernel $J^\alpha$ is translation invariant, i.e., for every $h \in \mathbb{R}^d$
		\begin{align}\label{eq:translation-invariance}\tag{I}
		J^\alpha(x+h, y+h) = J^\alpha(x, y) .
		\end{align}
		
		\item[(G-E)] There exists a constant $\Lambda\geq 1$ such that for every $\alpha\in  (0,2)$ and all $x,y \in \mathbb{R}^d$, with $x\neq y$ 
		\begin{align}\label{eq:global-elliptic-condition}\tag{$G$-$E$}
		\Lambda^{-1} \nu_\alpha (x-y) \leq J^\alpha(x,y) &\leq  \Lambda \nu_\alpha (x-y) .
		\end{align}
\end{itemize}%
\noindent Clearly \eqref{eq:global-elliptic-condition} implies \eqref{eq:elliptic-condition} and \eqref{eq:integrability-condition}. Finally, let us define the limit object, which is a local quadratic form of gradient type. Given $x \in \R^d$ and $\delta > 0$, we define the symmetric matrix $A(x) = (a_{ij}(x))_{1\leq i,j\leq d}$ by
\begin{align}\label{eq:coef-matrix}
a_{ij}(x) = \lim_{\alpha\to 2^{-}} \int_{B_\delta}  h_ih_j J^\alpha(x,x+h)dh
\end{align}
and for $u,v \in H^{1}(\Omega)$ the corresponding bilinear form by 
\begin{align*}
\mathcal{E}^A(u,v):=  \il_{\Omega} \big( A(x)\nabla u(x)\cdot \nabla v(x) \big)  \d x \,. 
\end{align*}

\medskip

\begin{remark} Let us discuss the assumption on the family $(J^\alpha)_\alpha$. 
	$(i)$ Under conditions \eqref{eq:elliptic-condition} and \eqref{eq:integrability-condition} the expression $\int_{B_\delta}  h_ih_j J^{\alpha_n}(x,x+h)dx$ converges for a suitable subsequence of $(\alpha_n)$. The existence of the limit in \eqref{eq:coef-matrix} poses an implicit condition on the family $(J^\alpha)_\alpha$. 
$(ii)$ \eqref{eq:elliptic-condition} and \eqref{eq:integrability-condition} ensure that the quantity $a_{ij}(x)$  does not  depend on the choice of $\delta$ and is bounded as a function in $x$, i.e. for all $\delta, r>0$, 
\begin{align*}
a_{ij}(x) = \lim_{\alpha\to 2^{-}} \int_{B_\delta}  h_ih_j J^\alpha(x,x+h)dh=  \lim_{\alpha\to 2^{-}} \int_{B_r}  h_ih_j J^\alpha(x,x+h)dh.
\end{align*}
$(iii)$ Under condition \eqref{eq:translation-invariance} the functions $a_{ij}(x)$, $1\leq i,j\leq d,$ are constant in $x$. 
$(iv)$ Regarding Proposition \ref{prop:elliptic-matrix-chap1}, which asserts that 
	\begin{align*}
d^{-1}\Lambda^{-1}|\xi|^2\leq \langle A(x) \xi, \xi \rangle \leq d^{-1}\Lambda |\xi|^2, \quad\text{ for every } x, \xi \in \R^d \,,
\end{align*}
 condition \eqref{eq:elliptic-condition} is a sufficient condition for what can be seen as nonlocal version of the classical ellipticity condition for second order operators in divergence form. 
$(v)$ Condition \eqref{eq:integrability-condition} ensures that long-range interactions encoded by $J^\alpha(x,y)$ vanish as $\alpha\to 2^-$. As a result, for some $\alpha_0\in (0,2)$ we have  
\begin{align*} 
\kappa_0 = \sup_{\alpha \in(\alpha_0,2)}\sup_{x\in \mathbb{R}^d} \int_{|h|>1} J^\alpha(x,x+h)dh<\infty.
\end{align*}
Hence conditions \eqref{eq:elliptic-condition} and \eqref{eq:integrability-condition} imply the following uniform L\'evy integrability type property:
\begin{align}\label{eq:consequence-integrability}
\sup_{\alpha \in(\alpha_0,2)}\sup_{x\in \mathbb{R}^d} \int_{\mathbb{R}^d} (1\land |h|^2)J^\alpha(x,x+h) \d h<\infty \,.
\end{align}
\end{remark}

\vspace{2mm}

\noindent Let us begin with a simple but important observation.  
\begin{proposition} Assume $\Omega\subset \R^d$ is a sufficiently smooth open bounded set or $\Omega=\R^d$. Let $\nu$ be a unimodal L\'evy kernel. Let  $\nu_\alpha$, $J^\alpha$ and $A$ be as above. 
	\begin{enumerate}[$(i)$]
	\item The following are regular Dirichlet forms on $L^2(\Omega)$:  $(\mathcal{E}^A, H^1(\Omega))$, $(\mathcal{E}^A, H^1_0(\Omega))$, $(\mathcal{E}_\Omega, \HnuOm)$, $(\mathcal{E}_\Omega, H_{\nu,0}(\Omega))$, $(\mathcal{E}, \VnuOmO)$, $(\mathcal{E}^\alpha_\Omega, \HnuOma)$, $(\mathcal{E}^\alpha_\Omega, H_{\nu_\alpha,0}(\Omega))$ and  $(\mathcal{E}^\alpha, \VnuOmO)$. 
	\item $(\mathcal{E}, \VnuOm \cap L^2(\R^d))$ and  $(\mathcal{E}^\alpha, \VnuOma\cap L^2(\R^d))$ are regular Dirichlet forms on $L^2(\R^d)$. 
	 
	\item The form $(\mathcal{E}, \VnuOm )$ is a regular Dirichlet form on $L^2(\R^d, \widetilde{\nu})$. The form $(\mathcal{E}^\alpha, \VnuOma )$ is a regular Dirichlet form on $L^2(\R^d, \widetilde{\nu}_\alpha)$. In particular for $J^\alpha (x,y) = C_{d,\alpha}|x-y|^{-d-\alpha}$ then the form $(\mathcal{E}^\alpha, V^{\alpha/2}(\Omega|\R^d))$ is a regular Dirichlet form on $L^2(\R^d, (1+|x|)^{-d-\alpha})$.
	\end{enumerate}
\end{proposition}

\medskip 

\begin{proof}
For $(iii)$ the regularity of the forms is a consequence  of Theorem \ref{thm:density} and Lemma \ref{lem:natural-norm-on-V}. The remaining details are easy to establish using Theorem \ref{thm:density}, Theorem \ref{thm:density-bis} and Theorem \ref{thm:density-zero-outside}. 
\end{proof}

\medskip

\begin{proposition}
	Let $\alpha_0\in (0,2)$ be as in \eqref{eq:consequence-integrability}.  The quadratic forms $(\mathcal{E}^{\alpha}_{\Omega}(\cdot, \cdot), H_{\nu_\alpha} (\Omega)) $ and $\big(\mathcal{E}^{\alpha}(\cdot, \cdot), V_{\nu_\alpha}(\Omega|\mathbb{R}^d)\cap L^2(\mathbb{R}^d)\big)$  are well defined for every $\alpha\in (\alpha_0, 2)$. 
\end{proposition}

\medskip 
\begin{proof}
	Let $\alpha\in (\alpha_0, 2)$  and $u \in H_{\nu_\alpha}(\Omega)$.  By the assumption \eqref{eq:elliptic-condition} and relation \eqref{eq:consequence-integrability} we have 
	
	\begin{align*}
	\mathcal{E}^\alpha_{\Omega}(u,u) &= \hspace*{-2ex}\iil_{\Omega\Omega \cap \{|x-y|\leq 1\}} (u(x)-u(y))^2 J^\alpha(x,y)\d x\d y+ \hspace*{-2ex} \iil_{\Omega\Omega \cap \{|x-y|>1\}} (u(x)-u(y))^2 J^\alpha(x,y)\d x\d y\\
	&\leq   \Lambda \hspace*{-2ex}\iil_{\Omega\Omega \cap \{|x-y|\leq 1\}} (u(x)-u(y))^2 \nu_\alpha(x-y)\d x\d y+ 4 \il_{\Omega} u^2(x)\d x\il_{|x-y|>1} J^\alpha(x,y)\d y\\
	&\leq   \Lambda \iil_{\Omega\Omega } (u(x)-u(y))^2 \nu_\alpha(x-y)\d x\d y+ 4 \kappa_0 \il_{\Omega} u^2(x)d x \leq (\Lambda+4\kappa_0)\|u\|^2_{H_{\nu_\alpha}(\Omega)}<\infty\, .
	\end{align*}
	Now if $u \in V_{\nu_\alpha}(\Omega|\mathbb{R}^d) $ then, from the above  we  deduce  $\mathcal{E}^\alpha_{\Omega}(u,u)<\infty$. By the same argument we obtain 
	\begin{align*}
	&\iil_{\Omega\Omega^c } (u(x)-u(y))^2 J^\alpha(x,y)\d x\d y\\&\leq \Lambda \hspace*{-2ex}\iil_{\Omega\Omega^c \cap \{|x-y|\leq 1\}} (u(x)-u(y))^2 \nu_\alpha(x-y)\d x\d y+ 2\hspace*{-2ex} \iil_{\Omega\Omega^c \cap \{|x-y|>1\}} (u^2(x)+u^2(y))J^\alpha(x,y)\d x\d y\\
	&\leq   \Lambda \iil_{\Omega\Omega^c\cap \{|x-y|\leq 1\}} (u(x)-u(y))^2 \nu_\alpha(x-y)\d x\d y+ 2 \kappa_0 \il_{\Omega} u^2(x)d x+ 2 \kappa_0 \il_{\Omega^c} u^2(x)d x\\
	&\leq  \Lambda \iil_{\Omega\Omega^c} (u(x)-u(y))^2 \nu_\alpha(x-y)\d x\d y+ 2 \kappa_0 \il_{\mathbb{R}^d} u^2(x)dx<\infty\, .
	\end{align*}
	Finally, we obtain 
	\begin{align*}
	\mathcal{E}^\alpha(u,u) &= \mathcal{E}^\alpha_{\Omega}(u,u) + 2\iil_{\Omega\Omega^c } (u(x)-u(y))^2 J^\alpha(x,y)\d x\d y<\infty\, . 
	\end{align*}

\end{proof}

\medskip


\begin{lemma}\label{lem:cont-qua-form}
	Let  $ D \subset \mathbb{R}^d$ be an $H^1$-extension domain.  Assume $J^\alpha $ satisfies \eqref{eq:elliptic-condition} and \eqref{eq:integrability-condition} and let   $\alpha_0\in  (0,2)$be as in \eqref{eq:consequence-integrability}.  Then, there exists a constant $C:= C(D, \Lambda, d, \alpha_0)$ such that for every $u \in H^1(D)$ and every $\alpha\in  (\alpha_0,2)$ we have 
	\begin{align*}
	\mathcal{E}^{\alpha}_D(u,u) \leq C\| u\|^2_{H^1(D)}.
	\end{align*}
	
\end{lemma}

\medskip 

\begin{proof}
	From the symmetry of  $J^\alpha(x,y)$ and \eqref{eq:consequence-integrability} we have the following estimates
	\begin{align*}
		\iil_{D\times  D \cap \{|x-y|\geq 1\}} (u(x)-u(y))^2 J^\alpha(x,y)\d x \, \d y& \leq 2\il_{D}u^2(x) \d x \il_{ |x-y|\geq 1} J^\alpha(x,y) \, \d y
		\leq 2\kappa_0\|u\|^2_{L^2(D)}.
	\end{align*}
	By Lemma \ref{lem:boundedness-limsup} we get 
	\begin{align*}
	\iil_{D\times  D \cap \{|x-y|\leq 1\}} (u(x)-u(y))^2 \nu_{\alpha}(x-y)\d x \, \d y
	&\leq C\| u\|^2_{H^1(D)}.
	\end{align*}
	Combining the above estimates along with  the condition \eqref{eq:elliptic-condition} we get $ \mathcal{E}^{\alpha}_D(u,u) \leq C\| u\|^2_{H^1(D)}. $

\end{proof}

\medskip

\noindent The following result is reminiscent of Theorem \ref{thm:charact-w1p}. 
\begin{theorem}\label{thm:quadratic-convergence-BBM}
	Let $D\subset \mathbb{R}^d$  be a $H^1$-extension domain (or $D=\R^d$). Then, under assumptions \eqref{eq:elliptic-condition} and \eqref{eq:integrability-condition}, for all $u\in  H^{1}(D)$ we have
	\begin{align}\label{eq:quadratics-limit}
	\lim_{\alpha \to 2^-}\iil_{D D} (u(x)-u(y))^2 J^\alpha(x,y)\d x \, \d y = \int_{D}( A(x)\nabla u(x)\cdot\nabla u(x) )\d x.
	\end{align}
 In particular, if $C_{d,\alpha }$ is the normalization constant of the fractional Laplacian then we have  
	\begin{align}\label{eq:quadratics-limit-special}
	&\lim_{\alpha \to 2^-} \iil_{D D} (u(x)-u(y))^2 \nu_\alpha(x-y)\d x \, \d y =\frac{1}{d} \int_{D} |\nabla u(x) |^2 \d x\\
	&\lim_{\alpha \to 2^-} \frac{C_{d,\alpha }}{2} \iil_{D D} (u(x)-u(y))^2 |x-y|^{-d-\alpha}\d x \, \d y =\int_{D} |\nabla u(x) |^2 \d x.\notag
	\end{align}
\end{theorem}

\vspace{2mm}
\begin{proof}  Lemma \ref{lem:cont-qua-form} suggests that  it suffices to prove \eqref{eq:quadratics-limit}  for   $u$ in a  dense subset of $H^{1}(D)$. For instance, let us choose $u\in C_c^2(\overline{D}) $
	\begin{align} 
	\iil_{D\times  D \cap \{|x-y|\geq 1\}} (u(x)-u(y))^2 J^\alpha(x,y)\d x \, \d y 
	\leq  4\il_{D}u^2(x) \d x \il_{ |x-y|\geq 1} J^\alpha(x,y) \, \d y \xrightarrow{\alpha\to 2^-} 0. 
	\end{align}
	\noindent Now, we consider the mapping $F:D\times (0,2)\to \mathbb{R}$ with 
	$$F(x,\alpha):= \il_{|x-y|\leq 1} (u(x)-u(y))^2 J^\alpha(x,y) \d y.$$ 
	By Taylor expansion we have $u(y)-u(x) = \nabla u (x)\cdot(y-x)+ r_1(x,y)|x-y|^2$ and therefore,  we can write 
	$(u(y)-u(x))^2 = (\nabla u (x)\cdot(y-x))^2+ r(x,y)|x-y|^3$
with bounded remainders $r(x,y)$  and $r_1(x,y)$. Hence, $F(x,\alpha)$ becomes
	\begin{align*}
	F(x,\alpha)&= \il_{|x-y|\leq 1}  [\nabla u(x)\cdot (y-x)]^2 J^\alpha(x,y) \d y+ R(x,\alpha)\,
	\end{align*}
	with
	\begin{align*}
	|R(x,\alpha)| &:=\Big|\,  \int_{|x-y|\leq 1} r(x,y) |x-y|^3 J^\alpha(x,y) \d y\Big| 
	\leq   C  \int_{|h|\leq 1 }|h|^3\nu_{\alpha}(h)\d h\xrightarrow{\alpha \to 2^-}0.
	\end{align*}
	Here, we have applied \eqref{eq:elliptic-condition} and Remark \ref{rem:asymp-nu}. Finally, we obtain
	\begin{align*}
	\lim_{\alpha\to 2^-} F(x,\alpha) &= \lim_{\alpha\to 2^-} \il_{|x-y|\leq 1}  [\nabla u(x)\cdot (y-x)]^2  J^\alpha(x,y) \d y \\
	%
	%
	&= \sum_{0\leq i,j\leq d}   \partial_i u(x)  \partial_j u(x)  \lim_{\alpha\to 2^-} \il_{|x-y|\leq 1} ( y_i-x_i)(y_j-x_j) J^\alpha(x,y) \d y \\
	&= \sum_{0\leq i,j\leq d}  a_{ij} (x)\partial_i u(x)  \partial_j u(x) = \langle A (x)\nabla u(x), \nabla u(x)\rangle.
	\end{align*}
	Applying the dominated convergence theorem yields
	\begin{align*}
	\lim_{\alpha \to2^-} \hspace*{-3ex} \iil_{D\times D\cap \{|x-y|\leq1 \}} \hspace*{-2ex} (u(x)-u(y))^2J^\alpha(x,y)\d x \, \d y = \lim_{\alpha \to2^-} \il_{D} F(x,\alpha)\d x = \il_{D} \langle A(x)\nabla u(x), \nabla u(x) \rangle \ \d x.
	\end{align*}
	
\noindent The case $J^\alpha(x,y)= \nu_{\alpha}(x-y)$ follows from Theorem \ref{thm:BBM-result} since $K_{d,2}= \frac{1}{d}$.  In fact, the cases $J^\alpha(x,y)= \nu_{\alpha}(x-y)$ and $J^\alpha(x,y)=  \frac{C_{d,\alpha}}{2}|x-y|^{-d-\alpha}$ follow from Example \ref{Ex:J-guys-singular}.

%
%
\end{proof}

\medskip

\begin{corollary}\label{cor:alpha-convergence-Hilbert}
	Let $\Omega\subset \R^d$ be a $H^1$-extension domain. Then, under assumptions \eqref{eq:elliptic-condition} and \eqref{eq:integrability-condition}, for all $u\in  H^{1}(D)$ the family of Hilbert spaces $\big(\HnuOma, \|\cdot\|_\alpha\big)_\alpha$ converges to the Hilbert space  $\big(H^1(\Omega), \|\cdot\|_A\big)$  as $\alpha\to 2$ in the sense of Definition \ref{def:convergence-hilbert-space}. Here we consider the norms
	\begin{align*}
	\|u\|^2_\alpha&= \|u\|^2_{L^2(\Omega)} + \mathcal{E}_\Omega^\alpha(u,u)\quad\text{and}\quad
	\|u\|^2_A = \|u\|^2_{L^2(\Omega)} + \mathcal{E}^A(u,u).
	\end{align*}
\end{corollary}
\medskip

\begin{proof}
	Note that Lemma \ref{lem:radial-monotone} and  the regularity of $\Omega$ imply that $H^1(\Omega)$ is a dense subspace of $\HnuOma$. Thus result follows since by Theorem \ref{thm:quadratic-convergence-BBM} $\|u\|_\alpha\xrightarrow{\alpha\to 2}\|u\|_A$ for every $u\in H^1(\Omega)$. 
\end{proof}

\medskip

\begin{theorem} \label{thm:mosco-liminf-inequality} Let  $\Omega\subset \R^d$ be a $H^1$-extension domain. 
	Let $(\alpha_n)_n$ with $\alpha_n\in (0,2)$ be a sequence such that $\alpha_n \xrightarrow{n \to \infty}2$.  Assume  $(u_n)_n$ with $u_n\in H_{\nu_{\alpha_n}}(\Omega)$  is a sequence converging in  $L^2(\Omega)$ to $u\in H^1(\Omega)$. Under the assumptions \eqref{eq:elliptic-condition} and \eqref{eq:integrability-condition} we have 
\begin{align}\label{eq:liminf-ineq}
\il_{\Omega} (A(x)\nabla u(x)\cdot\nabla u(x)) \d x  \leq \liminf\limits_{n \to \infty} \iil_{\Omega\Omega} (u_{n}(x) -u_{n}(y))^2 J^{\alpha_n}(x,y)\d  x\d y.
\end{align}
\end{theorem}

\vspace{2mm}

\begin{proof}
	The claim clearly holds if $\liminf\limits_{n\to \infty}\mathcal{E}^{\alpha_n}_\Omega(u_n,u_n)= \infty$. Assume $\liminf\limits_{n\to \infty}\mathcal{E}^{\alpha_n}_\Omega(u_n,u_n)<\infty$. Then  the sequence of real numbers $(\|u_n\|_{\alpha_n})_n$ has a  bounded subsequence. By Corollary \ref{cor:alpha-convergence-Hilbert}  we know that the sequence of  Hilbert spaces $\big(H_{\nu_{\alpha_n}}, \|\cdot\|_{\alpha_n}\big)_n$ converges to the Hilbert space  $\big(H^1(\Omega), \|\cdot\|_A\big)$. Therefore, in view of Theorem \ref{thm:varying-weak-convergence}, there exists a subsequence of $(u_n)_n$ that we still denote  by $u_n$  converging weakly  (in the sense of Definition \ref{def:varying-convergence}) to some $u'\in H^1(\Omega)$ that is for every sequence $(v_n)_n$ with $v_n\in H_{\alpha_n}(\Omega)$ strongly converging to $v\in H^1(\Omega)$ we have  
	\begin{align*}
	\big(u_n, v_n\big)_{  H_{\alpha_n}(\Omega)}\xrightarrow{n\to \infty} \big(u', v\big)_{  H^1(\Omega)}. 
	\end{align*}
	Moreover, we have 
	\begin{align*}
	\|u'\|_A^2\leq \liminf_{n\to \infty}\|u_n\|^2_{  H_{\alpha_n}(\Omega)}. 
	\end{align*}
	In particular, for $v\in C_c^\infty(\Omega)$ the constant sequence $v_n=v$ strongly converges to $v$.  Thus since $(u_n)_n$ converges to u in $L^2(\Omega)$ we have, 
	\begin{align*}
	\big(u, v\big)_{L^2(\Omega)} +\lim_{n\to \infty}\mathcal{E}^{\alpha_n}_\Omega(u_n,v) = \lim_{n\to \infty} \big(u_n, v_n\big)_{  H_{\alpha_n}(\Omega)}= \big(u', v\big)_{ L^2(\Omega)}+ \mathcal{E}^A(u',v).   
	\end{align*}
	Wherefrom, we get $\big(u, v\big)_{L^2(\Omega)} = \big(u', v\big)_{L^2(\Omega)} $ for all $v\in C_c^\infty(\Omega)$. This implies that $u=u'$ a.e. on $\Omega$. Finally, since $\|u_n-u\|_{L^2(\Omega)}\xrightarrow{n \to \infty}0$, the above inequality  yields $\mathcal{E}^A(u,u)\leq \liminf\limits_{n\to \infty}\mathcal{E}^{\alpha_n}_\Omega(u_n,u_n).$ Indeed,
	\begin{align*}
	\|u\|^2_{L^2(\Omega)} + \mathcal{E}^A(u,u)= \|u\|^2_A&\leq \liminf_{n\to \infty}\|u_n\|^2_{  H_{\alpha_n}(\Omega)}\\
	&= \liminf_{n\to \infty}\|u_n\|^2_{ L^2(\Omega)}+ \liminf_{n\to \infty}\mathcal{E}^{\alpha_n}_\Omega
	(u_n,u_n) \\
	&= \|u\|^2_{L^2(\Omega)} + \liminf_{n\to \infty}\mathcal{E}^{\alpha_n}_\Omega(u_n,u_n). 
	\end{align*}
\end{proof}

\noindent As a consequence of Theorem \ref{thm:quadratic-convergence-BBM} and Theorem \ref{thm:mosco-liminf-inequality} the following is true. 

\begin{corollary}
 Let $\Omega\subset \mathbb{R}^d$ be a $H^1$-extension domain (or $\Omega=\R^d$).  Assume \eqref{eq:elliptic-condition} and \eqref{eq:integrability-condition}.
 Then the quadratic forms $(\mathcal{E}^\alpha_{\Omega}(\cdot, \cdot), H_{\nu_\alpha}( \Omega) )_{\alpha}$  and converges to  $( \mathcal{E}^A(\cdot,\cdot), H^{1}( \Omega) )$ in the Gamma sense in $L^2(\Omega)$ as $\alpha\to 2^-$.
\end{corollary}
\medskip

\noindent Finally, we now are in a  position to prove the  main result of this section. 
\begin{theorem}\label{thm:Mosco-convergence}
	Let $\Omega\subset \mathbb{R}^d$ be an open bounded set with a Lipschitz continuous boundary. Assume \eqref{eq:elliptic-condition} and \eqref{eq:integrability-condition}.
	Then the two families of  quadratic forms $(\mathcal{E}^\alpha_{\Omega}(\cdot, \cdot), H_{\nu_\alpha}( \Omega) )_{\alpha}$  and  
	$( \mathcal{E}^\alpha(\cdot, \cdot),V_{\nu_\alpha}( \Omega|\mathbb{R}^d) )_{\alpha}$ both converge to  $( \mathcal{E}^A(\cdot,\cdot), H^{1}( \Omega) )$ in the Mosco sense in $L^2(\Omega)$ as $\alpha\to 2^-$.
\end{theorem}

\medskip

\begin{proof}Note that by Theorem \ref{thm:density} $C_c^\infty(\mathbb{R}^d)$ is dense in $ \VnuOma $  and by Theorem \ref{thm:density-zero-outside} $C_c^\infty(\overline{\Omega})$ is dense in $\HnuOma $. Hence it follows that $\VnuOma$ and $\HnuOma$  are dense in  $L^2(\Omega)$.  We proof the $\limsup$ and the $\liminf$ conditions  separately. 
	
	\medskip	 
	
	\par \textbf{Limsup:}
	Let $u\in L^2(\Omega)$, if  $u \not\in H^1(\Omega)$ then the $\limsup$ statement holds true since $ \mathcal{E}^A(u,u)=\infty$. Now let  $u \in H^1(\Omega)$. By identifying  $u$ to one of its   extensions $\overline{u}\in  H^1(\mathbb{R}^d)$, for  the sake of simplicity we can always assume that $u \in H^1(\mathbb{R}^d) $. On the one hand, Theorem \ref{thm:quadratic-convergence-BBM} shows that $\lim\limits_{\alpha\to 2^{-}} \mathcal{E}^\alpha_{\Omega}(u,u)=\mathcal{E}^A(u,u)$. On the other hand,  we have 
	\begin{align*}
	\mathcal{E}^\alpha(u,u)=  \mathcal{E}^\alpha_{\Omega}(u,u) + 2\iil_{\Omega\Omega^c} (u(x)-u(y))^2J^\alpha(x,y)\d x\d y.
	\end{align*}
Since $\Omega$ is Lipschitz, adapting the proof of Theorem \ref{thm:collapsing-accross-boundary} by applying Theorem \ref{thm:quadratic-convergence-BBM} we find that
	\begin{align}\label{eq:extern-lim}
	\iil_{\Omega\Omega^c} (u(x)-u(y))^2J^\alpha(x,y)\d x\d y\to 0\quad \text{as $\alpha\to 2^-$.}
	\end{align}
Thus, we conclude that for $u\in H^1(\Omega)$,
	\[\limsup_{\alpha \to 2^{-}} \mathcal{E}_\Omega^\alpha(u,u)= \limsup_{\alpha \to 2^{-}} \mathcal{E}^\alpha(u,u)=  \mathcal{E}^A(u,u). \]
Thus, choosing the  constant sequence  $u_\alpha= u$ for all $ \alpha \in (0,2)$ we are provided with  the $\limsup$ condition for  both forms $( \mathcal{E}^\alpha_{\Omega}, \HnuOma )_{\alpha }$  and $( \mathcal{E}^\alpha, \VnuOma)_{\alpha}$.
	
		\medskip
	\par \textbf{Liminf}:  Let $u, u_n\in L^2(\Omega)$ be  such  that   $u_n  \rightharpoonup u$ in $L^2(\Omega)$. Necessarily, $(u_n)_{n}$ is bounded in $ L^2(\Omega)$. 
	Let $( \alpha_n)_n$ be a sequence in $(0,2)$ such that $ \alpha_n \to 2^-$ as $n\to \infty$.  If   $\liminf\limits_{n \to \infty} \mathcal{E}_\Omega^{\alpha_n}(u_n,u_n) =\infty$ then, 
	\[  \mathcal{E}^A(u,u)\leq \liminf_{n \to \infty} \mathcal{E}_\Omega^{ \alpha_n }(u_n,u_n) = \liminf_{n \to \infty} \mathcal{E}^{ \alpha_n }(u_n,u_n) =\infty.  \]
	Assume $\liminf\limits_{n \to \infty} \mathcal{E}_\Omega^{ \alpha_n }(u_n,u_n)<\infty$ then according to Theorem \ref{thm:asymp-compactness} the sequence $(u_n)_n$ has a subsequence  (which we again denote by $u_n$) converging  in $L^2(\Omega)$ to some $\widetilde{u}\in H^1(\Omega)$. Consequently, as  $u_n \rightharpoonup u$ it readily follows that $(u_n)_n$ converges strongly to $u$ in $L^2(\Omega)$. Therefore, taking into account that $u\in H^1(\Omega)$, the desired liminf inequality   is an immediate consequence of Theorem \ref{thm:mosco-liminf-inequality}. 
\end{proof}

\medskip
 \noindent The next result is a variant of  Theorem \ref{thm:Mosco-convergence} with $H^1(\Omega)$ replaced by $H^1_0(\Omega)$.

\begin{theorem}\label{thm:Mosco-convergence-bis}
	Let $\Omega\subset \mathbb{R}^d$ be an open bounded set with a continuous boundary. Assume \eqref{eq:elliptic-condition} and \eqref{eq:integrability-condition}.
	Then the two families of  quadratic forms $\big(\mathcal{E}^\alpha_{\Omega}(\cdot, \cdot),\overline{C_c^\infty(\Omega)}^{ H_{\nu_\alpha}( \Omega)} \big)_{\alpha}$  and  
	$\big( \mathcal{E}^\alpha(\cdot, \cdot),V^\Omega_{\nu_\alpha}( \Omega|\mathbb{R}^d) \big)_{\alpha}$ both converge to  $( \mathcal{E}^A(\cdot,\cdot), H_0^{1}( \Omega) )$ in the Mosco sense in $L^2(\Omega)$ as $\alpha\to 2^-$.
\end{theorem}

\noindent The result relies on the density of $\overline{C_c^\infty(\Omega)}^{ H_{\nu_\alpha}( \Omega)}$ resp. $V^\Omega_{\nu_\alpha}( \Omega|\mathbb{R}^d)$. The density of the first space is trivial. The density of the second space is formulated in Theorem \ref{thm:density-bis}. Apart from the density issue, the details of the proof are the same as in the proof of  Theorem \ref{thm:Mosco-convergence}.

\section{Convergence of Dirichlet and Neumann problems}

The main purpose of this section is to prove the convergence of weak solutions, eigenvalues and eigenfunctions  of the  nonlocal Dirichlet and Neumann problems the local ones. Here again we assume that  $(\nu_\alpha)_{\alpha\in(0,2)}$ satisfies the condition \eqref{eq:assumption-nu-alpha} and  the family$(J^\alpha)_\alpha$ of symmetric kernels either satisfies the conditions \eqref{eq:elliptic-condition} and \eqref{eq:integrability-condition} or the global elliptic condition \eqref{eq:global-elliptic-condition}. Define  the operators $\mathscr{L}_\alpha$ and $\mathscr{N}_\alpha$ by
\begin{align*}
\mathscr{L}_\alpha u(x)&= \pv  \int_{\R^d}(u(x)- u(y)) J^\alpha(x,y)\d y
~~\text{and}~~
\mathscr{N}_\alpha u(x)= \int_{\Omega}(u(x)- u(y))J^\alpha(x,y)\d y. 
\end{align*}
Recall  that we note
\begin{align*}
\mathcal{E}^\alpha(v,v) &=\frac{1}{2} \iint\limits_{(\Omega^c\times\Omega^c)^c} (v(x)-v(y))^2 J^\alpha(x, y)\d x\d y\,\quad \text{for all}~~v\in \VnuOma\\
\mathcal{E}^A(v,v)&=  \il_{\Omega} \big( A(x)\nabla v(x)\cdot\nabla v(x) \big)  \d x \,\quad \text{for all}~~v\in H^1(\Omega)\,.
\end{align*}
\noindent We remind that $A$ is the elliptic  symmetric matrix satisfying
$d^{-1}\Lambda^{-1}|\xi|^2\leq \langle A(x) \xi, \xi \rangle \leq d^{-1}\Lambda |\xi|^2, $ for every $x, \xi \in \R^d \,$ given  by (see \eqref{eq:coef-matrix}) 
\begin{align*}
a_{ij}(x) = \lim_{\alpha\to 2} \int_{B_\delta}  h_ih_j J^\alpha(x,x+h)dh. 
\end{align*}
%


\noindent   We also introduce  the outwards normal derivative of a function $v$ on $\partial \Omega$ with respect to the matrix $A$ (also called co-normal derivative of $v$)  which is defined for $x\in \partial \Omega$ by $$ \frac{\partial v}{\partial n_A}(x) = A(x) \nabla v(x) \cdot n(x).$$
\noindent We start with the following preparation result.

\medskip

\begin{lemma}\label{lem:colapsing-to-boundary}
Let $\Omega\subset \mathbb{R}^d$ be  Lipschitz, open and  bounded.  Let $\varphi \in C^2_b(\mathbb{R}^d)$ and $v\in \VnuOma$.  With the conditions \eqref{eq:assumption-nu-alpha},  \eqref{eq:elliptic-condition},  \eqref{eq:integrability-condition} and \eqref{eq:translation-invariance} in force, the following assertions hold. 
	\begin{enumerate}[$(i)$]
		\item There exist $\alpha_0\in (0,2)$ and a constant $C>0$ independent of $\alpha$ such that 
		\begin{align*}
		\sup_{\alpha\in (\alpha_0,2)}\Big|\int_{\Omega^c}\mathscr{N}_\alpha \varphi(y)v(y)\d y\Big|\leq C\|\varphi\|_{C^2_b(\mathbb{R}^d)} \|v\|_{\VnuOma}\,.
		\end{align*}
		\item Assume $v\in H^1(\mathbb{R}^d)$ then
		\begin{align*}
		\lim_{\alpha\to 2}\int_{\Omega^c} \mathscr{N}_\alpha \varphi(y)v(y)\d y =\int_{\partial\Omega}   \frac{\partial \varphi}{\partial n_A}(x) v(x)\d\sigma(x)\,.
		\end{align*} 
	\end{enumerate}
In particular if  $J^\alpha(x,y) = C_{d, \alpha} |x-y|^{-d-\alpha}$ then $\frac{\partial \varphi}{\partial n_A}(x)= \frac{\partial \varphi}{\partial n}(x)$. Here,  $C_{d, \alpha} $ is the normalization constant given by the formula \eqref{eq:explicit-consta-Cds}.
\end{lemma}

\medskip 

\begin{proof}
	The condition \eqref{eq:translation-invariance} gives $J^\alpha(x,x+h) = J^\alpha(x,x-h)$ which implies that
	\begin{align*}
	\mathscr{L}_\alpha u(x)&= \pv  \int_{\R^d}(u(x)- u(y)) J^\alpha(x,y)\d y= -\frac12\int_{\R^d}(u(x+h)+ u(x-h) -2u(x)) J^\alpha(x,x+h)\d h. 
	\end{align*}
Wherefrom, since $A$ is constant, it is easy to show that (see Proposition \ref{prop:elliptic-matrix-chap1})
\begin{align*}
\lim_{\alpha\to 2}\mathscr{L}_\alpha u(x) =- \operatorname{tr}(A(\cdot)D^2 u)(x) = -\operatorname{div}(A(\cdot)\nabla u)(x). 
\end{align*}
Combing the estimates \eqref{eq:consequence-integrability}, \eqref{eq:second-difference} and $|\varphi(x+h) -\varphi(x)|\leq 2\|\varphi\|_{C^1_b(\R^d)} (1\land |h|)$ we have 
	\begin{align*}
	|	\mathscr{L}_\alpha \varphi|\leq 4 \kappa\|\varphi\|_{C^2_b(\mathbb{R}^d)}\quad\text{and}\quad \mathcal{E}^\alpha(\varphi, \varphi)\leq 4 \kappa  |\Omega|\|\varphi\|_{C^2_b(\mathbb{R}^d)} \quad\text{for all }~~\alpha\in (\alpha_0,2)\,
	\end{align*}
where 
\begin{align*}
\kappa= \sup_{\alpha \in(\alpha_0,2)}\sup_{x\in \mathbb{R}^d} \int_{\mathbb{R}^d} (1\land |h|^2)J^\alpha(x,x+h) \d h<\infty \,.
\end{align*}
	\noindent Therefore, since the linear mapping $ v\mapsto \mathcal{E}^\alpha(\varphi, v) - \int_{\Omega}	\mathscr{L}_\alpha \varphi(x)v(x)\d x$ is continuous on $\VnuOma$, the Green-Gauss formula \eqref{eq:green-gauss-nonlocal} is applicable for $\varphi\in C_b^2(\mathbb{R}^d)$ and $v\in \VnuOma$ and hence  the above estimates imply $(i)$ as follows
	
	\begin{align*}
	\Big|\int_{\Omega^c} 	\mathscr{N}_\alpha \varphi(y)v(y)\d y\Big|&= \Big| \mathcal{E}^\alpha(\varphi, v) - \int_{\Omega}	\mathscr{L}_\alpha  \varphi(x)v(x)\d x\Big|\\
	&\leq \mathcal{E}^\alpha(\varphi, \varphi)^{1/2} \mathcal{E}^\alpha(v,v)^{1/2} +\|	\mathscr{L}_\alpha \varphi\|_{L^2(\Omega)} \|v\|_{L^2(\Omega)}\\
	&\leq C\|\varphi\|_{C^2_b(\mathbb{R}^d)}\|v\|_{\VnuOma}\,.
	\end{align*}
	\noindent We known that $	\mathscr{L}_\alpha \varphi (x)\xrightarrow{\alpha\to 2} -\operatorname{div}(A(\cdot)\nabla \varphi)(x)$ for all $x\in \mathbb{R}^d$ and since $|\mathscr{L}_\alpha\varphi|\leq 4 \kappa\|\varphi\|_{C^2_b(\mathbb{R}^d)}$, the Lebesgue dominated convergence Theorem yields
	\begin{align*}
	\int_{\Omega}\mathscr{L}_\alpha  \varphi(x)v(x)\d x \xrightarrow{\alpha\to 2} \-\int_{\Omega}\operatorname{div}(A(\cdot)\nabla \varphi) (x)v(x)\d x\,.
	\end{align*}
	On the other hand, Theorem \ref{thm:collapsing-accross-boundary} combined with Theorem \ref{thm:quadratic-convergence-BBM} implies that 
		\begin{align*}
	\mathcal{E}^\alpha(\varphi, v) \xrightarrow{\alpha\to 2} \int_{\Omega} (A(x)\nabla \varphi(x)\cdot\nabla v(x))\d x.
	\end{align*}
	
	\noindent Finally from the foregoing and the classical Green-Gauss formula we obtain $(ii)$ as follows
	\begin{align*}
	\lim_{\alpha\to 2}\int_{\Omega^c} \mathscr{N}_\alpha  \varphi(y)v(y)\d y &=
	\lim_{\alpha\to 2} \mathcal{E}^\alpha(\varphi, v)-\lim_{\alpha\to 2} \int_{\Omega} 
	\mathscr{L}_\alpha  \varphi(x)v(x)\d x\\
	&= \int_{\Omega}(A(x)\nabla \varphi(x)\cdot\nabla v(x))\d x - \int_{\Omega} \operatorname{div}(A(\cdot)\nabla \varphi)(x) v(x)\d x\\
	&=\int_{\partial\Omega} \frac{\partial \varphi}{\partial n_A}(x) v(x)\d\sigma(x)\,.
	\end{align*}
The case $J^\alpha(x,y) = C_{d, \alpha} |x-y|^{-d-\alpha}$ follows since $A= (\delta_{ij})_{1\leq i,j\leq d}$ see  from Example \ref{Ex:J-guys-singular}.

\end{proof}

\medskip

\noindent The following result is crucial for our purpose. 
\begin{theorem}\label{thm:asympt-conv-subsequence}
 Let $\Omega\subset \mathbb{R}^d$ be open, bounded and connected with Lipschitz boundary. Assume the conditions \eqref{eq:assumption-nu-alpha} and \eqref{eq:global-elliptic-condition} hold. Let $u_\alpha \in  \VnuOma$ such that 
 \begin{align*}
 \sup_{\alpha\in (\alpha_*,2)} \|u_\alpha\|_{L^2(\Omega)}+ \mathcal{E}^{\alpha}(u_\alpha, u_\alpha)<\infty
 \end{align*}
 \noindent Then there exist $u\in H^1(\Omega)$ and a subsequence $\alpha_n\xrightarrow{n\to \infty} 2$ such that $\|u_{\alpha_n} -u\|_{L^2(\Omega)} \xrightarrow{n\to \infty}0$ and 
\begin{align*}
\mathcal{E}^{\alpha_n}(u_{\alpha_n}, v)\xrightarrow{n\to \infty} \mathcal{E}^{A}(u, v) = \int_{\Omega} (A(x)\nabla u (x)\cdot\nabla v(x))\d x\,\quad \text{for all $v\in H^1(\R^d)$}. 
\end{align*}
The same is true for $\Omega=\R^d$ with the exception that the strong convergence on $L^2_{\operatorname{loc}}(\R^d)$, i.e. for every compact  $K\subset \R^d$ we have $\|u_{\alpha_n} -u\|_{L^2(K)} \xrightarrow{n\to \infty}0$. 
\end{theorem}
\medskip

\begin{proof}
	\noindent On $\HnuOma$ and $H^1(\Omega)$ we define the scalar products, 
\begin{align*}
\big(w ,\, v \big)_{H_{\nu_{\alpha}} (\Omega)}&= \int_{\Omega} w(x)v(x)\d x+ \iint\limits_{\Omega\Omega} (w(x)-w(y))(v(x)-w(y))J^\alpha(x,y)\d x\d y\\
\big(w, v \big)_{H^1 (\Omega)} &= \int_{\Omega} w(x)v(x)\d x+ \int_{\Omega} (A(x)\nabla w (x)\cdot\nabla v(x))\d x\,.
\end{align*}
\noindent In view of Theorem \ref{thm:quadratic-convergence-BBM},  for both cases $\Omega\neq \R^d$ and $\Omega= \R^d$, we have $ \|v\|_{\HnuOma}\xrightarrow{\alpha\to 2} \|v\|_{H^1(\Omega)}$ for all $v \in H^1(\Omega)$. In other words, 
the family of Hilbert spaces $(\HnuOma)_\alpha$ converges to the Hilbert space $H^1(\Omega)$ as $\alpha\to 2$ in the sense of Definition \ref{def:convergence-hilbert-space}.  Whence we deduce from Theorem \ref{thm:varying-weak-convergence}, there exists $u\in H^1(\Omega)$ and a subsequence $\alpha_n \xrightarrow{n\to \infty}2$ such that
\begin{align*}
\lim_{n\to \infty}\big(u_{\alpha_n}, v \big)_{H_{\nu_{\alpha_n}} (\Omega)} = \big(u, v \big)_{H^1 (\Omega)} \quad\text{for all $v\in H^1(\Omega)$}.
\end{align*}
\noindent On the other hand, for the case $\Omega\neq \R^d$ Theorem \ref{thm:asymp-compactness} asserts there exist a further subsequence that we still denote by $(\alpha_n)_n$ and a function $ u'\in H^1(\Omega)$ such that $\|u_{\alpha_n}-u'\|_{L^2(\Omega)} \xrightarrow[]{n\to \infty}0$. Instead, if $\Omega=\R^d$, then fact that $\|u_{\alpha_n} -u'\|_{L^2(K)}\xrightarrow{n\to \infty}0$ for every compact $K\subset \R^d$ is due to Theorem \ref{thm:asympto-local-compact}. In either case, it is not difficult by taking the test function $v\in C_c^\infty(\R^d)$ to show that $u'= u$ a.e in $\Omega$. In both cases, from this we get the weak convergence in $L^2(\Omega) $ that is, 
\begin{align*}
\lim_{n\to \infty}\big(u_{\alpha_n}, v \big)_{L^2(\Omega)} = \big(u, v \big)_{L^2 (\Omega)} \quad\text{for all $v\in H^1(\Omega)$}.
\end{align*}
Therefore, this and the above weak convergence imply  that
\begin{align}\label{eq:weak-con-semi-omega}
\iint\limits_{\Omega\Omega} (u_{\alpha_n}(x)-u_{\alpha_n}(y))(v(x)-v(y))J^{\alpha_n}(x,y)\d x\d y\xrightarrow{n\to \infty} \int_{\Omega} (A(x)\nabla u'(x)\cdot\nabla v(x))\d x\,.
\end{align}
\noindent The proof is thus complete for case $\Omega=\R^d$. Next, if $\Omega\neq \R^d$ then by the uniform boundedness in \eqref{eq:uniform-boundedness} together with \eqref{eq:global-elliptic-condition} we find that 
\begin{align*}
\iint\limits_{\Omega\Omega^c} \big|(u_{\alpha_n}(x)-u_{\alpha_n}(y))(v(x)-v(y))\big|&J^{\alpha_n}(x,y)\d x\d y\\& \leq C \iint\limits_{\Omega\Omega^c} (v(x)-v(y))^2 \nu_{\alpha_n}(x-y)\d x\d y\xrightarrow[]{n\to \infty}0\,
\end{align*}
where the convergence follows from Theorem \ref{thm:collapsing-accross-boundary}. This combined with \eqref{eq:weak-con-semi-omega} implies
\begin{align*}
\mathcal{E}^{\alpha_n}(u_{\alpha_n}, v) \xrightarrow{n\to \infty} \mathcal{E}^A(u, v)= \int_{\Omega} (A(x)\nabla u (x)\cdot\nabla v(x))\d x\,.
\end{align*}
\end{proof}

\noindent Here is our first convergence theorem concerning  weak solutions. 
\begin{theorem}[\textbf{Convergence of weak solution I}]\label{thm:convergence-solutioin-I} Let $\Omega\subset \mathbb{R}^d$ be open, bounded and connected with Lipschitz boundary.Let $(f_\alpha)_\alpha$ be a family  converging weakly  to some $f$ in $L^2(\Omega)$ as $\alpha\to2$. For $\varphi \in C^2_b(\mathbb{R}^d)$ define $g_\alpha=  \mathscr{N}_\alpha\varphi$ and $g=\frac{\partial \varphi}{\partial n_A}$.  Assume the conditions \eqref{eq:assumption-nu-alpha} and \eqref{eq:global-elliptic-condition} hold. Suppose $u_\alpha \in  \VnuOma^\perp$ is a weak solution of the nonlocal Neumann problem $\mathscr{L}_\alpha u= f_\alpha$ on $\Omega$ and $ \mathscr{N}_\alpha u= g_\alpha$  on $\Omega^c$ that is we have 
	\begin{align*}
	\mathcal{E}^\alpha(u_\alpha, v) = \int_{\Omega}f_\alpha(x) v(x)\d x + \int_{\Omega^c}g_\alpha(x) v(x)\d x\quad\text{for all }~~v\in \VnuOma^\perp\,.
	\end{align*}
	
	\noindent  Let $u\in H^1(\Omega)^\perp$ be the unique weak solution in $H^1(\Omega)^\perp$ of the  Neumann problem  $-\operatorname{div}(A(\cdot)\nabla u)=  f$ on $\Omega$ and $\dfrac{\partial u}{\partial n_A} =g$ on $\partial\Omega$ that is we have
	%
	\begin{align*}
	\mathcal{E}^A(u, v)= \int_{\Omega } f(x) v(x)\d x + \int_{\partial\Omega } g(x)v(x)\d\sigma(x)\quad\text{for all }~~v\in H^1(\Omega)^\perp\,. 
	\end{align*}
Under the condition \eqref{eq:translation-invariance} or else, if  $g_\alpha=g=0$, then $\|u_\alpha-u\|_{L^2(\Omega)}\xrightarrow[]{\alpha\to 2}0$, i.e. $(u_\alpha)_\alpha$ converges to $u$ in $L^2(\Omega)$.  Moreover, we have the weak convergence  $ \mathcal{E}^{\alpha}(u_{\alpha}, v) \xrightarrow{\alpha\to 2} \mathcal{E}^A(u, v)$ for all $v\in H^1(\R^d)$. 
\end{theorem}

\vspace{1mm}

\begin{proof}
	In virtue of Theorem \ref{thm:robust-poincare} combined with condition \eqref{eq:global-elliptic-condition} there exist $\alpha_0\in (0,2)$ and  a constant $B>0$ depending only on $\alpha_0, \Omega$ and $d$ such that for all $v \in L^2(\Omega)^\perp$ and all $\alpha\in (\alpha_0, 2)$
	\begin{align*}
	\|v\|^2_{L^2(\Omega)}\leq C\mathcal{E}^\alpha(v,v)\quad \text{with $C=\Lambda B$}. 
	\end{align*}
	Therefore, up to relabelling the constant $C$, for all $v \in L^2(\Omega)^\perp$ and all $\alpha\in (\alpha_0, 2)$ we have 
	\begin{align}\label{eq:uniform-coercivity}
	\|v\|^2_{\VnuOma}\leq C\mathcal{E}^\alpha(v,v). 
	\end{align}
	
	\noindent On one hand, by the weak convergence of $(f_\alpha)_\alpha$ we may assume $\sup\limits_{\alpha\in (0,2)}\|f_\alpha\|_{L^2(\Omega)}<\infty$. This together with the definition of $u_\alpha$ along with Lemma \ref{lem:colapsing-to-boundary} $(i)$ we get 
	\begin{align*}
	\mathcal{E}^\alpha(u_\alpha,u_\alpha) &= \int_{\Omega} f_\alpha(x) u_\alpha (x)\d x + \int_{\Omega^c} g_\alpha(y) u_\alpha (y)\d y\\
	&\leq \|u_\alpha\|_{\VnuOma}(\|f_\alpha\|_{L^2(\Omega)} + \|\varphi\|_{C_b^2(\mathbb{R}^d)})\\
	&\leq C \|u_\alpha\|_{\VnuOma},
	\end{align*}
	\noindent for some constant $C>0$ independent of $\alpha$ and $\alpha\in (\alpha_0, 2)$.  Combining this with \eqref{eq:uniform-coercivity} then for a generic constant $C>0$ independent of $\alpha$ we have the following uniform boundedness
	\begin{align}\label{eq:uniform-boundedness}
	\|u_\alpha\|_{\HnuOma}\leq \|u_\alpha\|_{\VnuOma}\leq C\qquad\text{for all } \alpha\in (\alpha_0, 2)\,.
	\end{align}
	\noindent In accordance to Theorem \ref{thm:asympt-conv-subsequence}, there exist a subsequence $(\alpha_n)_n$ such that $\alpha_n\to 2$ and a function $ u'\in H^1(\Omega)$ such that $\|u_{\alpha_n}-u'\|_{L^2(\Omega)} \xrightarrow[]{n\to \infty}0$ and 
	\begin{align*}
	\mathcal{E}^{\alpha_n}(u_{\alpha_n}, v) \xrightarrow{n\to \infty} \mathcal{E}^A(u', v)= \int_{\Omega} (A(x)\nabla u' (x)\cdot\nabla v(x))\d x\,\quad\text{for all $v\in H^1(\R^d)$}.
	\end{align*}
\noindent In addition we have $u'\in H^1(\Omega)^\perp$ since $u_{\alpha}\in \VnuOma^\perp$ for all $\alpha\in (0,2)$.  The proof will be completed if we show that $u=u'$. To this end, we fix $v\in H^1(\Omega)^\perp$, given that $\Omega$ has a Lipschitz boundary we let $\overline{v}\in H^1(\mathbb{R}^d)$ be an extension of $v$. We know that $\overline{v}\in\VnuOma^\perp$ for all $\alpha\in (0,2)$. Thus by definition of $u_{\alpha_n} $ it follows  that 
	\begin{align*}
	\mathcal{E}^{\alpha_n}(u_{\alpha_n}, \overline{v}) = \int_{\Omega}f_{\alpha_n}(x) v (x)\d x+ \int_{\Omega^c}g_{\alpha_n}(y) \overline{v}(y)\d y\,.
	\end{align*}
Now, given that $(f_{\alpha_n})_n$ weakly converges to $f$ in $L^2(\Omega)$, letting $n \to \infty$ in the above yields 
	\begin{align*}
	\mathcal{E}^A(u', v)= \int_{\Omega}f(x) v (x)\d x+ \int_{\partial\Omega}g(x)v(x)\d \sigma(x)\,
	\end{align*}
	if  $g_\alpha=g=0$ or  by applying Lemma \ref{lem:colapsing-to-boundary} $(ii)$, if condition \eqref{eq:translation-invariance} holds. This holds true for any $v\in H^1(\Omega)^\perp$. By uniqueness we have $u=u'$ on $\Omega$ since $u'\in H^1(\Omega)^\perp$ is the weak solution of the corresponding local Neumann problem. 
	
\noindent The same reasoning is true for any sequence $(\alpha_n)_n$ with $\alpha_n\to2$. Thereupon, the uniqueness of $u\in H^1(\Omega)^\perp$ implies that $\|u_\alpha-u\|_{L^2(\Omega)}\xrightarrow{\alpha\to 2}0$ and 	$\mathcal{E}^{\alpha}(u_{\alpha}, v) \xrightarrow{\alpha\to 2} \mathcal{E}^A(u, v)$  for all $v\in H^1(\R^d)$.
\end{proof}

\vspace{2mm}

 \noindent The general case can be captured as follows. 
\begin{theorem}[\textbf{Convergence of weak solution II}]\label{thm:convergence-solutioin-II} 
Under the assumptions of Theorem \ref{thm:convergence-solutioin-I} assume that 
\begin{align*}
\int_\Omega f_\alpha(x)\d x + \int_{\Omega^c}g_\alpha(y)\d y= \int_\Omega f(x)\d x + \int_{\partial\Omega}g(y)\d y=0\qquad\text{for all $\alpha\in(0,2)$}.  
\end{align*}
Assume $w_\alpha \in  \VnuOma$ is a weak solution of the nonlocal Neumann problem $\mathscr{L}_\alpha u= f_\alpha$ on $\Omega$ and $ \mathscr{N}_\alpha u= g_\alpha$  on $\Omega^c$, i.e.
\begin{align}\label{eq:variational-nonlocal-neumann-alpha}
\mathcal{E}^\alpha(w_\alpha, v) = \int_{\Omega}f_\alpha(x) v(x) + \int_{\Omega^c}g_\alpha(x) v(x)\quad\text{for all }~~v\in \VnuOma\,.
\end{align}

\noindent Assume that the family $ (c_\alpha)_\alpha$ with $c_\alpha = \hbox{$\fint_\Omega $} w_\alpha (x)\d x$ is bounded.  Then there exist $w\in H^1(\Omega)$ and a sequence $\alpha_n\to2$ such that $\|w_{\alpha_n} -w\|_{L^2(\Omega)}\xrightarrow{n\to \infty}0$, where $w\in H^1(\Omega)$
is a weak solution of the  Neumann problem  $-\operatorname{div}(A(\cdot)\nabla  u) = f$ on $\Omega$ and $\dfrac{\partial u}{\partial n_A} =g$ on $\partial\Omega$, i.e.
\begin{align}\label{eq:variational-limit-neumann}
\mathcal{E}^A(u,v)= \int_{\Omega } f(x) v(x)\d x + \int_{\partial\Omega } g(x)v(x)\d\sigma(x)\quad\text{for all }~~v\in H^1(\Omega)\,. 
\end{align}

\noindent Conversely, if $w\in H^1(\Omega)$ solves \eqref{eq:variational-limit-neumann} then there exists $w_\alpha\in \VnuOma$ solving \eqref{eq:variational-nonlocal-neumann-alpha} such that we have $\|w_{\alpha} -w\|_{L^2(\Omega)}\xrightarrow{\alpha\to2}0$ and $ \mathcal{E}^{\alpha}(w_{\alpha}, v) \xrightarrow{\alpha\to 2} \mathcal{E}^A(w, v)$  for all $H^1(\R^d)$. 
\end{theorem}

\begin{proof} 
	 It suffices to observe that  $w_\alpha = u_\alpha+ c_\alpha$ and $w=u+ c$ where $ c_\alpha = \mbox{$\fint_\Omega$} w_\alpha (x)\d x$ , $ c =\mbox{$\fint_\Omega$} w(x)\d x$ and the functions $u_\alpha\in \VnuOma^\perp$, $u\in H^1(\Omega)^\perp$ are uniquely determined as in Theorem \ref{thm:convergence-solutioin-I} such that $\|u_\alpha-u\|_{L^2(\Omega)}\xrightarrow{\alpha\to 2}0$. 
Thus, the forwards statement follows since by assumption $(c_\alpha)_{\alpha\in (\alpha_0, 2)}$ is a bounded family of real numbers and thus has a converging subsequence. For the converse it suffices to take $w_\alpha = u_\alpha+c$ with $c \in \R$. 
	 
\end{proof}

\noindent We will  need the following lemma to establish the convergence of Dirichlet problem. 

\begin{lemma}\label{lem:robust-poincare-J-alpha}  Let $\Omega\subset  \R^d$ be open and bounded with continuous boundary. Assume the conditions \eqref{eq:assumption-nu-alpha} and  \eqref{eq:elliptic-condition} hold. Then there exist $\alpha_*\in (0,2)$  and a constant $B>0$ independent of $\alpha$ such that for all $\alpha\in (\alpha_*, 2)$  we have 
	\begin{align*}
	\|u\|^2_{L^2(\Omega)}&\leq  2B \Lambda\, \mathcal{E}^\alpha(u,u) \quad \text{for all}~u\in V^\Omega_{\nu_\alpha}(\Omega|\R^d).
	\end{align*}
	
\end{lemma}

\medskip 

\begin{proof}
	Since  $\Omega$ has continuous boundary, by Theorem \ref{thm:density-zero-outside} it is sufficient to prove the inequality for $u\in C_c^\infty(\Omega)$.  In virtue of Theorem \ref{thm:robust-poincare-friedrichs} there exist $\alpha_0\in(0,2)$ and a constant $B>0$ independent of $\alpha$ such that for every $\alpha\in (\alpha_0, 2)$ we have 
	\begin{align*}
	\|u\|^2_{L^2(\Omega)}\leq B \iil_{ (\Omega^c\times \Omega^c)^c}(u(x) -u(y))^2\nu_\alpha(x-y)\d x\d y.
	\end{align*}
From the condition \eqref{eq:assumption-nu-alpha} there exists $\alpha_1\in (0,2)$ such that for all $\alpha\in (\alpha_1, 2)$ we have 
	\begin{align*}
	\int_{|h|\geq 1}\nu_\alpha(h)\d h \leq \frac{1}{16B}. 
	\end{align*}
	
	Let $\alpha_*=\max(\alpha_0, \alpha_1)$, thus for $\alpha\in (\alpha_*, 2)$, employing the above and the condition \eqref{eq:elliptic-condition} yields
	\begin{align*}
	\|u\|^2_{L^2(\Omega)}&\leq B \iil_{ (\Omega^c\times \Omega^c)^c}(u(x) -u(y))^2\nu_\alpha(x-y)\d x\d y \\
	&\leq B\Lambda  \iil_{ (\Omega^c\times \Omega^c)^c }(u(x) -u(y))^2 \mathbbm{1}_{B_1}(x-y) J^\alpha(x,y)\d x\d y+ 8B\|u\|^2_{L^2(\Omega)}\int_{|h|\geq 1}\nu_\alpha(h)\d h \\
	&\leq B\Lambda\, \mathcal{E}^\alpha(u,u) + \frac{1}{2}\|u\|^2_{L^2(\Omega)}.  
	\end{align*}
	\noindent Thereupon we get $\|u\|^2_{L^2(\Omega)} \leq 2B \Lambda\, \mathcal{E}^\alpha(u,u). $
	
\end{proof}
\vspace{1mm}

\begin{theorem}[\textbf{Convergence of weak solution III}]\label{thm:convergence-solutioin-III} 
	Let $\Omega\subset \mathbb{R}^d$ be open, bounded and connected with continuous boundary. Let  $(f_\alpha)_\alpha$ be a family  converging weakly  to some $f$ in $L^2(\Omega)$ as $\alpha\to2$. Let $g\in H^1(\R^d)$.
Under the conditions \eqref{eq:assumption-nu-alpha}, \eqref{eq:elliptic-condition} and \eqref{eq:integrability-condition}  assume $u_\alpha \in  \VnuOma$ is a weak solution of the nonlocal Dirichlet problem $\mathscr{L}_\alpha u= f_\alpha$ on $\Omega$ and $ u= g$  on $\Omega^c$, i.e. 
	$u_\alpha- g\in V^\Omega_{\nu_\alpha}(\Omega|\R^d)$ and
	\begin{align*}
	\mathcal{E}^\alpha(u_\alpha, v) = \int_{\Omega}f_\alpha(x) v(x)\quad\text{for all }~~v\in V^\Omega_{\nu_\alpha}(\Omega|\R^d) \,.
	\end{align*}
	
	\noindent  Then $\|u_{\alpha} -u\|_{L^2(\Omega)}\xrightarrow{\alpha\to2}0$ where $u\in H^1(\Omega)$ is the unique  weak solution of the Dirichlet problem  $\operatorname{div}(A(\cdot) \nabla) u= f$ on $\Omega$ and $u=g$ on $\partial\Omega$, i.e. we have  $u-g \in H^1_0(\Omega)$ and 
	\begin{align*}
	\mathcal{E}^A(u, v) = \int_{\Omega } f(x) v(x)\d x\quad\text{for all }~~v\in H^1_0(\Omega)\,. 
	\end{align*}
Moreover, we have the weak convergence  $ \mathcal{E}^{\alpha}(u_{\alpha}, v) \xrightarrow{\alpha\to 2} \mathcal{E}^A(u, v)$ for all $v\in H^1(\R^d)$. 
\end{theorem}

\begin{proof}
Since, $ u_\alpha-g\in V^\Omega_{\nu_\alpha}(\Omega|\R^d)$, in view of Lemma \ref{lem:robust-poincare-J-alpha} there exist $\alpha_*\in(0,2)$ and a constant $B>0$ independent of $\alpha$ such that for every $\alpha\in (\alpha_*, 2)$ we have 
\begin{align*}
\|u_\alpha-g\|^2_{L^2(\Omega)}\leq 2B\Lambda\,\mathcal{E}^\alpha(u_\alpha-g, u_\alpha-g)
\end{align*}
which implies  that 
\begin{align*}
\|u_\alpha-g\|^2_{\VnuOma}\leq (2B\Lambda+1)\mathcal{E}^\alpha(u_\alpha-g, u_\alpha-g).
\end{align*}
In accordance to the definition of $u_\alpha$,  taking $C= \sqrt{2B\Lambda+1}$, the following holds
\begin{align*}
	|\mathcal{E}^\alpha(u_\alpha-g, u_\alpha-g)|&\leq  |\mathcal{E}^\alpha(u_\alpha, u_\alpha-g)|+| \mathcal{E}^\alpha(g, u_\alpha-g) |\\
	&= \Big| \int_{\Omega}f_\alpha(x) (u_\alpha(x)-g(x))\d x\Big|+ |\mathcal{E}^\alpha(g,u_\alpha-g) |\\
	&\leq \|f_\alpha\|_{L^2(\Omega)}  \|u_\alpha-g\|_{L^2(\Omega)}+ \mathcal{E}^\alpha(g,g)^{1/2} \mathcal{E}^\alpha(u_\alpha- g,u_\alpha-g)^{1/2}\\
	&\leq C\mathcal{E}^\alpha(u_\alpha- g,u_\alpha-g)^{1/2}\big( \|f_\alpha\|_{L^2(\Omega)} + \mathcal{E}^\alpha(g,g)^{1/2}\big). 
	\end{align*}
That is $	\mathcal{E}^\alpha(u_\alpha- g,u_\alpha-g)^{1/2}\leq C \big( \|f_\alpha\|_{L^2(\Omega)} + \mathcal{E}^\alpha(g,g)^{1/2}\big).$ Further, given that $u_\alpha-g=0$ on $\Omega^c$,  we get $\|u_\alpha-g\|_{H_{\nu_\alpha}(\R^d)} = \|u_\alpha-g\|_{\VnuOma}.$ Hence we have 
\begin{align*}
\|u_\alpha-g\|_{H_{\nu_\alpha}(\R^d)} &= \Big( \|u_\alpha-g\|^2_{L^2(\Omega)}+ \mathcal{E}^\alpha(u_\alpha- g,u_\alpha-g)\Big)^{1/2}\\
&\leq \sqrt{2}C\mathcal{E}^\alpha(u_\alpha- g,u_\alpha-g)^{1/2}\\
&\leq 2C^2 \big( \|f_\alpha\|_{L^2(\Omega)} + \mathcal{E}^\alpha(g,g)^{1/2}\big). 
\end{align*}
Wherefrom, we find that 
\begin{align*}
\|u_\alpha\|_{H_{\nu_\alpha}(\R^d)} &\leq \|g\|_{H_{\nu_\alpha}(\R^d)}+ \|u_\alpha-g\|_{H_{\nu_\alpha}(\R^d)}\\
&\leq 2C^2 \big( \|f_\alpha\|_{L^2(\Omega)} +  \|g\|_{L^2(\R^d)}+ 2\mathcal{E}^\alpha(g,g)^{1/2}\big).
\end{align*}
On other the hand, under the conditions \eqref{eq:elliptic-condition} and \eqref{eq:integrability-condition} the estimate \eqref{eq:consequence-integrability} implies the following uniform L\'evy integrability type property:
\begin{align*}
\kappa_*= \sup_{\alpha \in(\alpha_*, 2)}\sup_{x\in \mathbb{R}^d} \int_{\mathbb{R}^d} (1\land |h|^2)J^\alpha(x,x+h) \d h<\infty \,\quad\text{with $\alpha_*>\alpha_0$.}
\end{align*}
\noindent Next, exploiting the estimate \eqref{eq:levy-p-estimate} we get 
\begin{align*}
\mathcal{E}^\alpha(g,g)\leq \iil_{ \R^d \R^d} (g(x) -g(y))^2 J^\alpha(x, x+h)\d h\d x\leq 4 \kappa_*\|g\|^2_{H^1(\R^d)} <\infty \,.
\end{align*}
\noindent By the weakly convergence of $(f_\alpha)_\alpha$ we may assume $\sup\limits_{\alpha\in (0,2)}\|f_\alpha\|_{L^2(\Omega)}<\infty$. Finally, we have shown that 
\begin{align}\label{eq:uniform-boundedness-D}
\|u_\alpha\|_{H_{\nu_\alpha}(\R^d)} \leq 2C^2 \big( \sup\limits_{\alpha\in (0,2)}\|f_\alpha\|_{L^2(\Omega)} +\|g\|_{L^2(\R^d)}+8 \kappa_*\|g\|^2_{H^1(\R^d)}\big):=M \quad\text{for all $\alpha\in (\alpha_*, 2)$}
\end{align}
where $M$ does not depend on $\alpha$. 
\noindent In view of Theorem \ref{thm:asympt-conv-subsequence}, there exist $u'\in H^1(\R^d)$ and a subsequence $\alpha_n \xrightarrow[]{n\to \infty}2$ such that for all compact $K\subset \R^d,$ $\|u_{\alpha_n}-u'\|_{L^2(K)} \xrightarrow{n\to \infty}0$ and for all $v\in H^1(\R^d)$
\begin{align}\label{eq:weak-con-semi-omega-D} 
\begin{split}
\iint\limits_{\R^d\R^d} (u_{\alpha_n}(x)-u_{\alpha_n}(y))(v(x)-v(y))J^{\alpha_n}(x,y)\d x\d y
\xrightarrow{n\to \infty} \int_{\R^d} (A(x)\nabla u'(x)\cdot\nabla v(x))\d x.
\end{split}
\end{align}
%
%
We claim that $ u'-g=0$ on $\Omega^c$. In fact for any compact $K\subset \Omega^c$ since  $u_{\alpha_n} -g =0$ on $\Omega^c$ we have  
$$\|u'-g\|_{L^2(K)}= \|u_{\alpha_n}-u'\|_{L^2(K)} \xrightarrow{n\to \infty}0.$$
 
 \noindent  Thus we have $u'-g \in H^1(\R^d)$ and $u'-g=0$ on $\Omega^c$ which means that $u'-g\in H^1_0(\Omega)$.  On the other hand,  we also have
$\|u_{\alpha_n}-u'\|_{L^2(\Omega)} \xrightarrow{n\to \infty}0$ since $\Omega$ is bounded.  
The proof will be completed if we  show that $u=u'$. To this end, we fix $v\in H^1_0(\Omega)$. We naturally assume that $v=0$ on $\Omega^c$ so that  $v\in H^1(\mathbb{R}^d)$.  The relation \eqref{eq:weak-con-semi-omega-D} implies
\begin{align*}
\mathcal{E}^{\alpha_n}(u_{\alpha_n}, v) \xrightarrow{n\to \infty} \mathcal{E}^A(u', v)= \int_{\Omega} (A(x)\nabla u' (x)\cdot\nabla v(x))\d x\,.
\end{align*}
\noindent We know that $v\in V^\Omega_{\nu_\alpha}(\Omega|\R^d)$ for all $\alpha\in (0,2)$. Thus by definition of $u_{\alpha_n} $ it follows  that 
\begin{align*}
\mathcal{E}^{\alpha_n}(u_{\alpha_n}, v) = \int_{\Omega}f_{\alpha_n}(x) v (x)\d x.
\end{align*}
Now, given that $(f_{\alpha_n})_n$ weakly converges to $f$ in $L^2(\Omega)$, letting $n \to \infty$ in the above gives 
\begin{align*}
\mathcal{E}^A(u', v)= \int_{\Omega}f(x) v (x)\d x\,.
\end{align*}
This holds true for any $v\in H^1_0(\Omega)$ and $u'-g\in H^1_0(\Omega)$. Thus, $u'\in H^1_0(\Omega)$ is the weak solution of the corresponding local Dirichlet problem. By uniqueness we have $u=u'$ on $\Omega$. 
The same reasoning is true for any sequence $(\alpha_n)_n$ with $\alpha_n\to2$. Thereupon, the uniqueness of $u\in H^1(\Omega)$ implies that $\|u_\alpha-u\|_{L^2(\Omega)}\xrightarrow{\alpha\to 2}0$ and $ \mathcal{E}^{\alpha}(u_{\alpha}, v) \xrightarrow{\alpha\to 2} \mathcal{E}^A(u, v)$ for all $v\in H^1(\R^d)$. 
  
\end{proof}

\medskip

\noindent Another result on the convergence of Dirichlet problems can be found in \cite{Voi17}. Next, we prove the convergence of nonlocal eigenvalues and eigenfunctions to the local ones.
\begin{theorem}[\textbf{Convergence of eigenvalues I}]\label{thm:convergence-eigenpair-I}
Let $\Omega\subset \mathbb{R}^d$ be open, bounded and connected with Lipschitz boundary. Assume that the conditions  \eqref{eq:global-elliptic-condition} and \eqref{eq:assumption-nu-alpha} are satisfied. 
For each $\alpha$, assume  that the operator $\mathscr{L}_\alpha$ has a family of normalized Neumann eigenpairs $(\mu_{\alpha,n}, \phi_{\alpha,n})_n$. That is for all $n\geq0$, $\phi_{\alpha,n} \in  \VnuOma$,  $\|\phi_{\alpha,n}\|_{L^2(\Omega)} =1$ the elements $(\phi_{\alpha,n})_n$ are mutually $L^2(\Omega)$-orthogonal, i.e. $(\phi_{\alpha,i}, \phi_{\alpha,k} )_{L^2(\Omega)}= 0$ if $i\neq k$, and 
	\begin{align*}
	\mathcal{E}^\alpha(\phi_{\alpha,n}, v) =\mu_{\alpha,n} \int_{\Omega}\phi_{\alpha,n} (x) v(x) \d x\quad\text{for all }~~v\in \VnuOma\,.
	\end{align*}
	\noindent Then for each $n\geq 0$, the eigenpair $(\mu_{\alpha,n}, \phi_{\alpha,n})_\alpha$ converges to $(\mu'_n, \phi'_n)$ in $\R\times L^2(\Omega)$ up to a subsequence as $\alpha \to 2$. To be more precise, $\mu_{\alpha_j, n}\xrightarrow{j\to \infty}\mu'_n$  in $\R$ and $\|\phi_{\alpha_j,n}-\phi'_n\|_{L^2(\Omega)}\xrightarrow{\alpha\to 2}0$  for some $\alpha_j\xrightarrow{j\to \infty}2$. Moreover the family $(\mu'_n, \phi'_n)_n$ is the sequence of the normalized Neumann eigenpairs of the $\operatorname{div}(A(\cdot)\nabla )$, i.e. $(\phi'_i,\phi'_k)_{L^2(\Omega)}=\delta_{i,k}$ and for each $n\geq 0$
	\begin{align*}
	\mathcal{E}^A(\phi'_n, v)= \mu'_n\int_{\Omega } \phi'_n(x) v(x)\d x \quad\text{for all }~~v\in H^1(\Omega)\,. 
	\end{align*}
\end{theorem}

\begin{proof}
First of all, we know  that $\phi_{\alpha,0}= \phi'_0= |\Omega|^{-1}$ and  $\mu_{\alpha,0}= \mu'_0=0$. By proceeding as in the proof of Theorem \ref{thm:convergence-sharp-const-poinc}, it can be shown that 
$\mu_{\alpha,1}\xrightarrow{\alpha\to 2}\mu'_1$. It can be shown that $\phi_{\alpha,1}$ converges to some $\phi'_1$ in $L^2(\Omega)^\perp$ up to a subsequence. More generally,  for fixed $n\geq 1$ assume there exists a sequence $\alpha_j:=\alpha_j(n)\xrightarrow{j\to \infty}2$ such  that for every $0\leq k\leq n-1$, $(\mu_{\alpha_j,k}, \phi_{\alpha_j(n),k})_n$ converges to $(\mu'_k, \phi'_k)$ in $\R\times L^2(\Omega)$ as $j\to \infty$.  Then by Rayleigh's quotient (see Chapter \ref{chap:IDEs}) we have 
\begin{align*}
\mu_{\alpha_j, n}= \mathcal{E}^{\alpha_j}(\phi_{\alpha_j, n}, \phi_{\alpha_j, n})= \min_{v\in V_{\alpha_j, n}}\Big\{\frac{\mathcal{E}^{\alpha_j}(v, v)}{\|v\|^2_{L^2(\Omega)}}\Big\}
\end{align*}

\noindent where $V_{\alpha_j, n}= \big\{ v\in V_{\alpha_j}(\Omega|\R^d): (\phi_{\alpha_j, k}, v )_{L^2(\Omega)}=0,~i=0,1, \cdots, n-1\big\}$. Observing that $V_{\alpha_j, n}\subset L^2(\Omega)^\perp$, it follows from the robust Poincar\'e inequality (see Corollary \ref{cor:robust-poincare}) that there exists a constant $C>0$  independent of $\alpha_j$ for which we have 
\begin{align*}
	\mu^{-1}_{\alpha_j,n}\leq C\quad\text{for all $j\geq 1$}. 
\end{align*}
\noindent That is, the sequence $(\mu^{-1}_{\alpha_jn})_j$ is bounded. Thus, without loss of generality we may assume that $\mu^{-1}_{\alpha_j,n}$ converges to some $\mu'^{-1}_{n}\geq 0$ as $j\to \infty$. On the other hand,  the sequence $(\mu_{\alpha_j,n})_j$ being bounded,  we have  
\begin{align*}
\|\phi_{\alpha_j,n}\|^2_{L^2(\Omega)} +  \mathcal{E}^{\alpha_j}(\phi_{\alpha_j, n},\phi_{\alpha_j, n})\leq  1+ \sup_{j\geq 1}\mu_{\alpha_j,n}:=M<\infty \quad\text{for all $j\geq 1$}. 
\end{align*}
Therefore, from Theorem \ref{thm:asympt-conv-subsequence} one is able to find a further subsequence  $(\alpha_j(n+1))_j$ which we denote $\alpha'_j =\alpha_j (n+1)$ of the sequence $(\alpha_j)_j$( recall $\alpha_j=\alpha_j(n)$) and $\phi'_n\in H^1(\Omega)$ such that
$\|\phi_{\alpha'_j,n}-\phi'_n\|_{L^2(\Omega)}\xrightarrow{j\to \infty}0$ and for all $v\in H^1(\R^d)$ we have 
\begin{align*}
\mathcal{E}^A(\phi'_n, v) 
&= \lim_{j\to \infty} \mathcal{E}^{\alpha'_j}(\phi_{\alpha'_j,n}, v)= \lim_{j\to \infty} \mu_{\alpha'_j,n} \int_{\Omega} \phi_{\alpha'_j,n} (x) v(x)\d x= \mu'_n \int_{\Omega}  \phi'_n (x) v(x)\d x.
\end{align*}
\noindent In particular, since $\partial \Omega$ is Lipschitz we have 
\begin{align*}
\mathcal{E}^A(\phi'_n, v) = \mu'_n \int_{\Omega}  \phi'_n (x) v(x)\d x\quad\text{for all $v\in H^1(\Omega)$}. 
\end{align*}
Moreover, for $0\leq i\leq k\leq n$ we have 
\begin{align*}
(\phi'_i,\phi'_k)_{L^2(\Omega)}= \lim_{j\to \infty}(\phi_{\alpha'_j, i}, \phi'_{\alpha'_j, k})_{L^2(\Omega)}=\delta_{i,k}. 
\end{align*}
By induction the result remains true for all $n\in \mathbb{N}_0$.
\end{proof}

\medskip

\noindent A similar reasoning leads to the convergence of eigenpairs associated with Dirichlet condition.
\begin{theorem}[\textbf{Convergence of eigenvalues II}]\label{thm:convergence-eigenpair-II}
	Let $\Omega\subset \mathbb{R}^d$ be open, bounded and connected.  Assume that the conditions  \eqref{eq:global-elliptic-condition} and \eqref{eq:assumption-nu-alpha} are satisfied. 
	For each fixed $\alpha$, assume  that the operator $\mathscr{L}_\alpha$ has a family of normalized Dirichlet eigenpairs $(\lambda_{\alpha,n}, \varphi_{\alpha,n})_n$. That is, for all $n\geq 1$, $\phi_{\alpha,n} \in   V^\Omega_{\nu_\alpha}(\Omega|\R^d)$,  $\|\varphi_{\alpha,n}\|_{L^2(\Omega)} =1$ the elements of $(\varphi_{\alpha,n})_n$ are mutually $L^2(\Omega)$-orthogonal, i.e. $(\varphi_{\alpha,i}, \varphi_{\alpha,k} )_{L^2(\Omega)}= 0$ if $i\neq k$ and 
	\begin{align*}
	\mathcal{E}^\alpha(\varphi_{\alpha,n}, v) =\lambda_{\alpha,n} \int_{\Omega}\varphi_{\alpha,n} (x) v(x) \d x\quad\text{for all }~~v\in V^\Omega_{\nu_\alpha}(\Omega|\R^d)\,.
	\end{align*}
	\noindent Then for each $n\geq 1$, the eigenpair $(\lambda_{\alpha,n}, \varphi_{\alpha,n})_\alpha$ converges to $(\lambda'_n, \varphi'_n)$ in $\R\times L^2(\Omega)$ up to a subsequence as $\alpha \to 2$. To be more precise, $\lambda_{\alpha_j, n}\xrightarrow{j\to \infty}\lambda'_n$  in $\R$ and $\|\varphi_{\alpha_j,n}-\phi'_n\|_{L^2(\Omega)}\xrightarrow{j\to \infty}0$ for some $\alpha_j\xrightarrow{j\to \infty}2$.  Moreover the family $(\lambda'_n, \varphi'_n)_n$ is the sequence of the normalized Dirichlet eigenpairs of the $\operatorname{div}(A(\cdot)\nabla )$, i.e. $(\varphi'_i,\varphi'_k)_{L^2(\Omega)}=\delta_{i,k}$ and for each $n\geq 1$
	\begin{align*}
	\mathcal{E}^A(\varphi'_n, v)= \lambda'_n\int_{\Omega } \varphi'_n(x) v(x)\d x \quad\text{for all }~~v\in H^1_0(\Omega)\,. 
	\end{align*}
\end{theorem}

\bigskip
\noindent Specializing the aforementioned results with $J^\alpha(x,y) = C_{d,\alpha}|x-y|^{-d-\alpha}$ that is we get $ \mathscr{L}_\alpha=(-\Delta)^{\alpha/2}$ and  $-\operatorname{div}(A(\cdot)\nabla) = -\Delta$ we get following corollaries. 

\begin{corollary}
	Let $\Omega\subset \mathbb{R}^d$ be  open bounded and connected with Lipschitz boundary. Let $(f_\alpha)_\alpha$ be a family  converging weakly  to some $f$ in $L^2(\Omega)$ as $\alpha\to2$. For $\varphi \in C^2_b(\mathbb{R}^d)$ define $g_\alpha=  \mathcal{N}_\alpha\varphi$ and $g=\frac{\partial \varphi}{\partial n}$.  Let $u_\alpha \in  V^{\alpha/2}(\Omega|\R^d)^\perp$ be the weak solution in $V^{\alpha/2}(\Omega|\R^d)^\perp$ of the nonlocal Neumann problem for the fractional Laplacian, i.e.  $u_\alpha$ is a weak solution of the following
	\begin{align*}
	\text{$(-\Delta)^{\alpha/2} u= f_\alpha$ on $~\Omega$ \quad and\quad $ \mathcal{N}_\alpha u= g_\alpha$  on $~\Omega^c$}. 
	\end{align*}
	\noindent  Then we have  $\|u_\alpha-u\|_{L^2(\Omega)}\xrightarrow{\alpha\to 2}0$ where $u\in H^1(\Omega)^\perp$ is the unique weak solution in $H^1(\Omega)^\perp$ of the  Neumann problem 
	\begin{align*}
	-\Delta u=f\, \, \text{ on ~~$\Omega$}\quad\text{ and}\quad  \frac{\partial u}{\partial n} =g \text{ on ~$\partial\Omega$}.
	\end{align*}	
	 Furthermore, let $(\mu_{\alpha,n}, \phi_{\alpha,n})_n$ be  the family of normalized Neumann eigenpairs of the fractional Laplacian $(-\Delta)^{\alpha/2}$. Let $(\mu'_n, \phi'_n)_n$ be the sequence of  normalized Neumann eigenpairs of the Laplacian $-\Delta$. Then for each $n\geq 0$, the eigenpair $(\mu_{\alpha,n}, \phi_{\alpha,n})_n$ converges to $(\mu'_n, \phi'_n)$ in $\R\times L^2(\Omega)$ up to a subsequence as $\alpha \to 2$.
\end{corollary}

\medskip 

\begin{corollary}
	Let $\Omega\subset \mathbb{R}^d$ be open, bounded and connected. Let $(f_\alpha)_\alpha$ be a family  converging weakly  to some $f$ in $L^2(\Omega)$ as $\alpha\to2$. Let $g\in H^1(\R^d)$.  Let $u_\alpha \in  V^{\alpha/2}(\Omega|\R^d)$ be the weak solution in $V^{\alpha/2}(\Omega|\R^d)$ of the nonlocal Dirichlet problem for the fractional Laplacian, i.e. $u_\alpha$ is a weak solution of the following
	\begin{align*}
	\text{$(-\Delta)^{\alpha/2} u= f_\alpha$ on $~\Omega$ \quad and\quad $ u= g$  on $~\Omega^c$}. 
	\end{align*}
	\noindent Then $\|u_\alpha-u\|_{L^2(\Omega)}\xrightarrow{\alpha\to 2}0$ where $u\in H^1(\Omega)$ be the unique weak solution in $H^1(\Omega)$ of the  Dirichlet problem 
	\begin{align*}
	-\Delta u=f\, \, \text{ on ~~$\Omega$}\quad\text{ and}\quad  u =g \text{ on ~$\partial\Omega$}.
	\end{align*}	
\noindent Furthermore, let $(\lambda_{\alpha,n}, \varphi_{\alpha,n})_n$ be  the family of normalized Dirichlet eigenpairs of the fractional Laplacian $(-\Delta)^{\alpha/2}$. Let $(\lambda'_n, \varphi'_n)_n$ be the sequence of  normalized Dirichlet eigenpairs of the Laplacian $-\Delta$. Then for each $n\geq 1$, the eigenpair $(\lambda_{\alpha,n}, \varphi_{\alpha,n})_n$ converges to $(\lambda'_n, \varphi'_n)$ in $\R\times L^2(\Omega)$ up to a subsequence as $\alpha \to 2$.
\end{corollary}

\medskip

\noindent Another special case is to consider $J^\alpha(x,y) = \frac1d\nu_\alpha (x,y)$ with $(\nu_\alpha)_\alpha$ satisfying the condition \eqref{eq:assumption-nu-alpha}. In this  case we also have $-\operatorname{div}(A(\cdot)\nabla) = -\Delta$.

\appendix 
\chapter{Lebesgue Spaces}\label{chap:complemnt-Lebesgue}
\section{Lebesgue spaces} \label{sec:lebesgue-space}

Let the triplet $(X, \mathcal{A}, \mu)$ be a measure space that is $X$ is a set, $\mathcal{A}$ is a $\sigma$-algebra on $X$ and $\mu: \mathcal{A}\to [0, \infty]$ is positive measure on $X $. Basic notions related to such a triplet are recorded in \cite{Alt16,bogachev2007volumeI,grafakos04} especially in \cite{EG15}. Our exposition here and in the next section abundantly relies upon these materials . 
\begin{definition}
	For $1\leq p\leq \infty$, the space $L^p(X)$ is the space of class of measurable functions $u:X\to \mathbb{R}$ such that $\|u\|_{L^p(X)}<\infty$ with 
	\begin{align*}
	\|u\|_{L^p(X)}&= \Big(\int_X|u(x)|^p\,\d \mu(x)\Big)^{1/p}\quad\textrm{for}\quad 1\leq
	p<\infty\\
	\|u\|_{L^\infty(X)}&= \inf\{c>0: |u|\leq c~~\mu-\textrm{a.e on $X$}\}\quad\textrm{for}\quad p=\infty.
	\end{align*}
\end{definition}

\noindent Let us quote without proofs some significant results on $L^p(X)$. We start with H\"older inequality which appears to be the most important inequality on $L^p$-spaces. 

\begin{theorem}[H\"older inequality]
	For $1\leq p,p'\leq \infty$, such that $ p+p'=pp'$ then for all measurable functions $u,v: X\to \mathbb{R}$ the following inequalities $\|uv\|_{L^1(X)}\leq \|u\|_{L^p(X)}\|v\|_{L^{p'}(X)}$ holds true that is 
	\begin{align*}
	\int_X|u(x)v(x)|\d \mu(x)\leq \Big(\int_X|u(x)|^p\,\d \mu(x)\Big)^{1/p} \Big(\int_X|u(x)|^{p'}\,\d \mu(x)\Big)^{1/p'}. 
	\end{align*}
	\noindent More generally assume $r,p_1\cdots, p_n\in [1, \infty]$ satisfy the relation 
	$$\frac{1}{p_1}+\cdots + \frac{1}{p_n}= \frac{1}{r}.$$
	Then for all measurable functions $u_1,\cdots,u_n: X\to \mathbb{R}$ we have 
	\begin{align*}
	\|u_1\cdots u_n\|_{L^r(X)}\leq \|u_1\|_{L^{p_1}(X)}\cdots \|u_n\|_{L^{p_n}(X)}.
	\end{align*}
	
	\noindent In particular, $u_1\cdots u_n\in L^r(X)$ once $u_i\in L^{p_i}(X), ~i=1,\cdots,n$. 
\end{theorem}

\medskip

\noindent One of the most influential inequalities on the  integration theory is the Jensen's inequality. 
\begin{theorem}[Jensen inequality]
	Let $u: X\to \mathbb{R}$ be measurable and let $\varphi: \mathbb{R}\to [0, \infty)$ be a convex function such that $\varphi\circ u\in L^1(X)$. Assume $\mu(X)<\infty$ then,
	\begin{align*}
	\varphi\left(\frac{1}{\mu(X)}\int_X u(x)\,\d \mu(x)\right) \leq \frac{1}{\mu(X)}\int_X \varphi\circ u(x)\,\d \mu(x).
	\end{align*}
\end{theorem}

\bigskip

\noindent Actually, it is possible to derive  H\"{o}lder inequality from Jensen's inequality.
Exploiting the H\"older inequality we are able to prove the Minkowski inequality. 

\begin{theorem}[Minkowski inequality]
	For $1\leq p\leq \infty$, then for all measurable functions $u,v: X\to \mathbb{R}$ one has the following triangle inequality
	\begin{align*}
	\|u+v\|_{L^p(X)}\leq \|u\|_{L^p(X)}+\|v\|_{L^{p}(X)}. 
	\end{align*}
\end{theorem}

\medskip

\noindent The Minkowski inequality shows that $L^p(X)$ is a normed space. We have the following.
\begin{theorem}[Riesz-Fischer]
	For $1\leq p\leq \infty$, the space $L^p(X)$ is a Banach space under the norm $u\mapsto\|u\|_{L^p(X)}$.
\end{theorem}

\bigskip

\noindent The next very important result is the monotone convergence theorem due
to Lebesgue and Beppo Levi.

\begin{theorem}[Monotone convergence or Beppo Levi]
	Assume $(u_n)_n$ is a sequence of nonnegative measurable functions on $X$ such that $0\leq u_n\leq u_{n+1}$ a.e for each $n\geq 1$ then
	\begin{align*}
	\int_X\lim_{n\to\infty} u_n(x)\,\d \mu(x)= \lim_{n\to\infty} \int_Xu_n(x)\,\d \mu(x).
	\end{align*}
\end{theorem}

\medskip

\noindent An alternative to the monotone convergence theorem when the monotonicity is violated is the so called Fatou's lemma. 
\begin{theorem}[Fatou's lemma]
	Assume $(u_n)_n$ is a sequence of nonnegative measurable functions on $X$ then
	\begin{align*}
	\int_X\liminf_{n\to\infty} u_n(x)\,\d \mu(x)\leq \liminf_{n\to\infty} \int_X u_n(x)\,\d \mu(x).
	\end{align*}
\end{theorem}

\medskip

\noindent Of course to complete the discussion on convergence one needs to recalls the Lebesgue Dominated convergence theorem but we defer this to the next section which will be derived as an immediate consequence of the Vitali convergence Theorem (cf. Theorem \ref{thm:lebesgue-dominated}).

\begin{theorem}[Fubini Theorem]
	Let $(X_1, \mathcal{A}_i, \mu_i)~i=1,2$ be two $\sigma$-finite measurable spaces and let $\big(X_1\times X_2, \mathcal{A}_1\otimes\mathcal{A}_2, \mu_1\otimes \mu_2\big)$ be the corresponding product space. Assume $v:X_1\times X_2\to \mathbb{R}$ is $\mathcal{A}_1\otimes\mathcal{A}_2$-measurable then\footnote{Fubini theorem for nonnegative functions is often known as the Fubini-Tonelli theorem} 
	\begin{relsize}{-.9}
	\begin{align*}
	\hspace{-1.5ex} \il_{X_1\times X_2} |v(x_1, x_2)|\,\d \mu_1\hspace{-0.5ex}\otimes\hspace{-0.5ex} \mu_2(x_1,x_2)\hspace{-0.2ex} = \il_{X_1} \Big(\hspace{-1ex} \il_{X_2}\hspace{-0.5ex} |v(x_1, x_2)|\,\d \mu_1(x_1)\Big) \d \mu_2(x_2)\hspace{-0.2ex} = \il_{X_2}\Big(\hspace{-1ex} \il_{X_1} |v(x_1, x_2)|\,\d \mu_1(x_2)\Big) \d \mu_2(x_1)
	\end{align*}
\end{relsize}
	additionally if $v\in L^1(X_1\times X_2)$ then, 
	\begin{relsize}{-.9}
	\begin{align*}
	 \il_{X_1\times X_2}v(x_1, x_2)\,\d \mu_1\otimes \mu_2(x_1,x_2)\hspace{-0.5ex} = \il_{X_1} \Big(\hspace{-0.5ex} \il_{X_2} v(x_1, x_2)\,\d \mu_1(x_1)\Big) \d \mu_2(x_2)\hspace{-0.5ex} = \il_{X_2} \Big(\hspace{-0.5ex} \il_{X_1} v(x_1, x_2)\,\d \mu_1(x_2)\Big) \d \mu_2(x_1).
	\end{align*}
	\end{relsize}
\end{theorem}

\section{The Vitali convergence theorem }\label{sec:vitali-theorem}

It is the aim of this section to establish the Vitali convergence theorem which provides necessary and sufficient conditions for convergence in $L^p$-spaces. Throughout, we assume that 
$(X, \mathcal{A}, \mu)$ is a measure space and by "measurable" or "integrable" we refer to this triplet. We start with some basics. 

\begin{definition}
	A sequence $(u_{n})_{n\in \mathbb{N}}$ of measurable functions is said to converge in measure to $u$ if for all $\varepsilon>0$
	\begin{align*}
	\lim_{n\to \infty}\mu\{x\in X:~~|u_n(x)-u(x)|\geq \varepsilon\}=0. 
	\end{align*}
	We said that $(u_{n})_{n\in \mathbb{N}}$ is Cauchy in measure if for all $\varepsilon>0$
	\begin{align*}
	\lim_{m,n\to \infty}\mu\{x\in X:~~|u_m(x)-u_n(x)|\geq \varepsilon\}=0.
	\end{align*}
\end{definition}

\medskip

\noindent Let us see the connections between the convergence in measure and other types of convergence. 
\begin{proposition}\label{prop:convergence-inmeasure}
	Let $u$ and $(u_{n})_{n\in \mathbb{N}}$ be measurable functions. The following statements hold true. 
	\begin{enumerate}[$(i)$]
		\item Suppose $\mu(X)<\infty$. If the sequence $(u_n)_n$ converges pointwise $\mu$-almost everywhere to $u$ then it converges in measure to $u$. 
		\item If the sequence $(u_n)_n$ converges in $L^p(X)~~(1\leq p<\infty)$ to $u$ then it converges in measure to $u$.
		\item If the sequence $(u_n)_n$ converges in measure to $u$ then there is a subsequence $(u_{n_j})_j$ of $(u_n)_n$ which converges pointwise $\mu$-almost everywhere to $u$.
	\end{enumerate}
	
\end{proposition}

\medskip

\begin{proof} Without any loss of generality we assume that $u=0$. \newline 
	$(i)$ Given that $u_n\to 0$ $\mu$-a.e., for any $\varepsilon>0$ the measurable set $E= \{x\in X:~\limsup\limits_{n\to \infty}|u_n(x)|>\varepsilon\}$ has measure zero and 
	contains the measurable sets $B_n= \bigcup\limits_{k\geq n} \{|u_k|\geq \varepsilon\}$ with $B_{n+1}\subset B_n$. Since $\mu(X)<\infty$, by monotonicity of the measure $\mu$ we get 
	\begin{align*}
	\lim\limits_{n\to \infty}\mu(\{x\in X:~|u_n(x)|>\varepsilon\} )\leq \lim\limits_{n\to \infty}\mu(\bigcup\limits_{k\geq n} \{|u_k|\geq \varepsilon\} ) = \mu(\bigcap\limits_{n\geq 1} B_n)\leq \mu(E)=0. 
	\end{align*}
	\noindent $(ii)$ Assume $|u_n\|_{L^p(X)}\xrightarrow{n\to \infty}0$ then the convergence in measure is a straightforward consequence of Chebyshev’s inequality, since for all $\varepsilon>0$ we have
	\begin{align*}
	\lim\limits_{n\to \infty}\mu(\{x\in X:~|u_n(x)|>\varepsilon\}) \leq \lim\limits_{n\to \infty} \frac{\|u_n\|_{L^p(X)}}{\varepsilon} = 0. 
	\end{align*}
	\noindent $(iii)$ Assume $u_n$ converges in measure to $0$ then for every $\varepsilon= \frac{1}{2^k}$ we capable to gradually construct $(n_k)_k$ with $n_k< n_{k+1}$ such that 
	\begin{align*}
	\mu(\{x\in X:~|u_n(x)|>\frac{1}{2^k}\}) \leq \frac{1}{2^k}\qquad\textrm{
		for all }\quad n\geq n_k.
	\intertext{Now let }
	E=\bigcap\limits_{n\geq 1}\bigcup\limits_{k\geq n} \{x\in X:~|u_{n_k}(x)|>\frac{1}{2^k}\} =\bigcap\limits_{n\geq 1} E_n.
	\end{align*}
	\noindent Observing that $E_{n+1}\subset E_n$ and $\mu(X)<\infty$ then the monotonicity of the measure $\mu$ yields
	\begin{align*}
	\mu(E) = \lim_{n\to \infty} \mu(E_n) \leq \lim_{n\to \infty}\sum_{k=n}^{\infty} \mu(\{|u_{n_k}|>\frac{1}{2^k}\}) \leq \lim_{n\to \infty}\sum_{k=n}^{\infty}\frac{1}{2^k} = \lim_{n\to \infty}\frac{1}{2^{n-1}}=0.
	\end{align*}
	\noindent If $x\in X\setminus E$ then there exists $N_x\geq 1$ such that $x\not\in \{x\in X:~|u_{n_k}(x)|>\frac{1}{2^k}\}$ for all $k\geq N_x$ i.e $|u_{n_k}(x)|\leq \frac{1}{2^k}$ for all $k\geq N_x$ which means that $u_{n_k}(x)\xrightarrow{k\to \infty}0$. Hereby providing the subsequence sought for.
\end{proof}

\medskip

\noindent Setwise the convergence in measure on a finite measure space may be defined by means of a topology. 
\begin{proposition}
	Assume the measure $\mu$ is finite, i.e. $\mu(X)<\infty$. Let $(u_{n})_{n\in \mathbb{N}}$ be a sequence of measurable functions. A necessary and sufficient condition for $(u_n)_n$ to converge in measure is that $\rho(u_n)\xrightarrow{n\to\infty}0$ where
	\begin{align*}
	\rho(u)= \int_X\frac{|u(x)|}{1+|u(x)|}\,\d \mu(x). 
	\end{align*}
\end{proposition}

\medskip

\begin{proof}
	First and foremost, it clearly appears that for each $\varepsilon>0$ and every $n\geq 1$ we have 
	\begin{align*}
	A_n= \{x\in X:~|u_n(x)|\geq \varepsilon\}= \Big\{x\in X:~\frac{|u_n(x)|}{1+|u_n(x)|}\geq\frac{\varepsilon}{1+\varepsilon}\Big\}.
	\end{align*}
	Utilizing this remark, entails on the  one hand that $\mu(A_n)\leq \frac{1+\varepsilon}{\varepsilon}\rho(u_n)$ and on the other hand that 
	\begin{align*}
	\rho(u_n)&= \int_{A_n} \frac{|u_n(x)|}{1+|u_n(x)|}\,\d \mu(x)+ \int_{X\setminus A_n}\frac{|u_n(x)|}{1+|u_n(x)|}\,\d \mu(x)\\
	& \leq \int_{A_n}1\,\d \mu(x)+ \int_{X\setminus A_n}\frac{\varepsilon}{1+\varepsilon }\,\d \mu(x)= \mu_(A_n)+ \varepsilon\mu(X). 
	\end{align*}
	The claim readily follows by letting $n\to \infty$ and after $\varepsilon\to 0$. 
	
\end{proof}

\medskip

\begin{remark}
	It clearly appears that the convergence in measure induces a topology on the set of class of measurable
	functions on $X$ which is metrisable relatively to the distance function 
	\begin{align*}
	\rho(u,v) := \rho(u-v) = \int_X\frac{|u(x)-v(x)|}{1+|u(x)-v(x)|}\,\d \mu(x).
	\end{align*}
\end{remark}

\noindent The convergence in measure induces a topology on the space of class of measurable functions. The following theorem states its completeness under this topology. For a proof we refer to \cite{bogachev2007volumeI, grafakos04}.
\begin{theorem}
	A sequence $(u_n)_n$ of measurable functions on $X$ that is Cauchy in measure converges in measure up to a subsequence. 
\end{theorem}

\bigskip

\begin{theorem}[Ergorov's Theorem]\label{thm:Ergorov}
	Assume that $E\subset X$ is a measurable subset with finite measure, i.e. $\mu(E)<\infty$. If $(u_n)_n$ is a sequence of measurable functions converging almost everywhere on $A$ to a function $u$ then for every $\varepsilon>0$ there is $A\subset E$ such that $\mu(A)<\varepsilon$ and the convergence of $(u_n)_n$ to $u$ holds uniformly on $E\setminus A$. 
\end{theorem}

\medskip

\begin{proof}
	Let us fix $k\geq 1$ and consider the measurable set
	\begin{align*}
	A_k= \bigcap_{n\geq 1}A_{n,k} \qquad\text{with}\quad A_{n,k}= \bigcup_{m\geq n}\{x\in E:~ |u_m(x)-u(x)|\geq \frac{1}{k}\}. 
	\end{align*}
	\noindent Let $x\in A_k$, it turns out that for all $n\geq 1$ there is $m\geq n$ such that $|u_m(x)-u(x)|\geq \frac{1}{k}$. 
	In other words, $x$ belongs to a set of measure zero given that $(u_n)_n$ converges to $u$ almost everywhere on $E$. Hence it is easy to recognise that $\mu(A_k)=0$. 
	Moreover $\mu(E)<\infty$ and $A_{n+1,k}\subset A_{n,k}$ so that by the monotonicity of the measure $\mu$, $\lim\limits_{n \to \infty}\mu(A_{n,k})= \mu(A_k) =0$. Now for fixed $\varepsilon>0$ there is $n_k$ such that $\mu(A_{n_k,k})<\frac{\varepsilon}{2^k}.$ 
	Therefore setting 
	$$A= \bigcup\limits_{k\geq 1}A_{n_k,k} = \bigcup\limits_{k\geq 1}\bigcup\limits_{m\geq n_k}\{x\in E:~ |u_m(x)-u(x)|\geq \frac{1}{k}\}$$ 
	one obtains $\mu(A)<\varepsilon$ and the uniform convergence $(u_n)_n$ on $E\setminus A$ follows since for each $k\geq 1$
	\begin{align*}
	\sup_{x\in E\setminus A} |u_m(x)-u(x)| < \frac{1}{k}\qquad\textrm{for all}\quad m\geq n_k.
	\end{align*}
\end{proof}

\medskip

\noindent It is significant to mention that the conclusion of Ergorov's theorem may fail if the finiteness of $\mu(E)$ is violated. Indeed on $\mathbb{R}$, the sequence $u_n= \mathds{1}_{[n,n+1]}$ converges pointwise everywhere to $u=0$ on $\mathbb{R}$ and does not uniformly since for each $n\geq 1$ $\sup\limits_{x\in \mathbb{R}}|u_n(x)| =1$. 

\bigskip

\noindent Let us now introduce some important notions such as the equiintegrability and the tightness. 
\begin{definition}
	A set $\mathcal{F}$ of measurable functions on $X$ is called uniformly integrable (or equiintegrable) if: for every $\varepsilon>0$ there exists $\delta>0$ such that for a measurable set $E\in \mathcal{A}$, if $\mu(E)<\delta$, then
	\begin{align*}
	\int_E |u(x)|\,\d \mu(x)<\varepsilon\qquad\textrm{for all}\quad u\in \mathcal{F}.
	\end{align*}
	This can be fairly written in short form as 
	\begin{align*}
	\lim_{\mu(E)\to 0} \sup_{u\in \mathcal{F}} \int_E |u(x)|\,\d \mu(x)= 0.
	\end{align*}
	More generally, $\mathcal{F}$ is called $p$-uniformly integrable (or $p$-equiintegrable) for some $1\leq p<\infty$ if
	\begin{align*}
	\lim_{\mu(E)\to 0} \sup_{u\in \mathcal{F}} \int_E |u(x)|^p\,\d \mu(x)= 0.
	\end{align*}
\end{definition}

\bigskip

\noindent The next result gives an equivalent definition of the equiintegrability. 
\begin{proposition}
	Let $\mathcal{F}$ be a set of measurable functions. The following assertions are true. 
	
	\begin{enumerate}[$(i)$]
		\item A necessary and sufficient condition for $\mathcal{F}$ to be uniformly integrable is that 
		\begin{align*}
		\lim_{\mu(E)\to 0} \sup_{u\in \mathcal{F}} \Big|\int_E u(x)\,\d \mu(x)\Big|= 0.
		\end{align*}
		\item If in addition, $\mathcal{F}$ is a bounded subset of $L^p(X)$ then $\mathcal{F}$ is uniformly integrable if and only if 
		\begin{align*}
		\lim_{R\to \infty} \sup_{u\in \mathcal{F}} \il_{\{|u| \geq R\}}| u(x)|^p\,\d \mu(x)= \inf_{R\geq 0} \sup_{u\in \mathcal{F}} \il_{\{|u| \geq R\}}| u(x)|^p\,\d \mu(x)= 0.
		\end{align*}
		
	\end{enumerate}
\end{proposition}

\medskip

\begin{proof} $(i)$ It is clear by definition that the uniform integrability implies
	\begin{align*}
	\lim_{\mu(E)\to 0} \sup_{u\in \mathcal{F}} \Big|\int_E u(x)\,\d \mu(x)\Big|= 0.
	\end{align*}
	For the converse, let $\varepsilon$ and $\delta$ be such that if $ \mu(E)<\delta$, 
	\begin{align*}
	\Big|\int_E u(x)\,\d \mu(x)\Big|<\varepsilon/2\qquad\textrm{for all}\quad u \in \mathcal{F}.
	\end{align*}
	Observe that for every $u\in \mathcal{F}$ and fix a measurable set $E$ with $\mu(E)<\delta $, the sets $E_1=E\cap \{u\geq 0\}$ and $E_2=E\cap \{u< 0\}$ are measurable too. Hence, the desired claim follows from the above since $u $ is arbitrarily chosen and $\mu(E_1),\mu(E_2)<\delta$ and
	\begin{align*}
	\int_E|u(x)| \,\d \mu(x)&= \int_{E\cap \{u\geq 0\}} \hspace*{-5ex}u(x) \,\d \mu(x) + \int_{E\cap \{u< 0\}} \hspace*{-5ex}-u(x) \,\d \mu(x)\\
	&=\Big|\int_{E_1}u(x) \,\d \mu(x)\Big| + \Big|\int_{E_2} u(x) \,\d \mu(x)\Big|<\varepsilon.
	\end{align*}
	$(ii)$ Assume $\mathcal{F}$ is bounded in $L^p(X)$ and let $M= \sup\limits_{u\in \mathcal{F}}\|u\|_{L^p(X)}$. For $u\in \mathcal{F}$ and every $R>0$ we have
	\begin{align*}
	\mu\{|u|>R\} = \frac{1}{R^p}\il_{\{|u|>R\}}|u(x)|\,\d \mu(x)\leq \frac{1}{R^p}\int_X|u(x)|\,\d \mu(x) \leq \frac{M}{R^p}. 
	\end{align*}
	This forces $\mu\{|u|>R\} \to 0$ as $R\to \infty$ and hence if $\mathcal{F}$ is uniformly integrable then
	\begin{align*}
	\lim_{R\to \infty} \sup_{u\in \mathcal{F}} \il_{\{|u| \geq R\}}| u(x)|\,\d \mu(x)= \inf_{R\geq 0} \sup_{u\in \mathcal{F}} \il_{\{|u| \geq R\}}| u(x)|\,\d \mu(x)= 0.
	\end{align*}
	The converse holds straightforwardly since for  each measurable set $E$ with a small measure, the following is true. 
	\begin{align*}
	\sup_{u\in \mathcal{F}} \int_E |u(x)|\,\d \mu(x) 
	&\leq \sup_{u\in \mathcal{F}} \il_{E\cap \{|u|< R\}} |u(x)|\,\d \mu(x)+ \sup_{u\in \mathcal{F}} \il_{E\cap \{|u|\geq R\}} |u(x)|\,\d \mu(x)\\
	& \leq M\mu(E)+ \sup_{u\in \mathcal{F}} \il_{\{|u|\geq R\}} |u(x)|\,\d \mu(x).
	\end{align*}
	
\end{proof}

\bigskip

\noindent We now introduce the tightness which an is accurate expedient concept preventing a family of measurable function from escaping to infinity.
\begin{definition}
	Let $1\leq p<\infty$. A set $\mathcal{F}$ of measurable functions on $X$ is called to be $p$-tight (or simply tight for $p=1$) if for every $\varepsilon>0$ there exists measurable set $E$ with a finite measure i.e $ \mu(E)<\infty$ such that 
	\begin{align*}
	\int_{X\setminus E} |u(x)|^p\,\d \mu(x)<\varepsilon\qquad\textrm{for all}\quad u \in \mathcal{F}.
	\end{align*}
	This can be written in compact form as 
	\begin{align*}
	\inf_{\stackrel{\mu(E) <\infty}{E~\textrm{measurable}} } \sup_{u\in \mathcal{F}} \int_{X\setminus E} |u(x)|^p\,\d \mu(x)= 0.
	\end{align*}
	
\end{definition}

\noindent When $X=\mathbb{R}^d$ is equipped with the Lebesgue measure the $p$-tightness $\mathcal{F}$ is equivalent to 
\begin{align*}
\lim_{R\to \infty}\sup_{u\in \mathcal{F}}\il_{\{|x|\geq R\}}|u(x)|^p\,\d x =0. 
\end{align*}

\medskip

\noindent Warning: Some authors define equiintegrability as uniform integrability plus 
tightness. Of course if $\mu(X)<\infty$ then the tight is gratuitous.

\bigskip

\begin{remark} Let $u\in L^p(X)$ then the set $\mathcal{F}=\{u\}$ is $p$-tight and $p$-uniformly integrable. Indeed, since $u$ is integrable, it emanates form the monotone convergence theorem that for $\varepsilon>0$ there is $n\geq 1$ large enough such that 
	\begin{align*}
	\int_X |u(x)|^p\mathds{1}_{\{|u|\leq \frac{1}{n}\}}(x)\,\d \mu(x)<\varepsilon
	\intertext{and as well }
	\int_X |u(x)|^p\mathds{1}_{\{|u|>n\}}(x)\,\d \mu(x)<\varepsilon/2. 
	\end{align*}
	On one hand, the tightness follows from the first estimate since letting $A_n= \{x\in X:~|u(x)|>\frac{1}{n}\}$ the Chebyshev’s inequality yields
	\begin{align*}
	\mu(A_n)= \mu(\{x\in X:~|u(x)|>\frac{1}{n}\}) \leq n\|u\|_{L^p(X)}<\infty
	\end{align*}
	and we have 
	\begin{align*}
	\il_{X\setminus A_n} |u(x)|^p \,\d \mu(x) =\int_X |u(x)|^p\mathds{1}_{\{|u|\leq \frac{1}{n}\}}(x)\,\d \mu(x)<\varepsilon.
	\end{align*}
	On the other hand, for any measurable set $E$ with $\mu(E)<\delta= \varepsilon/2n^p$ the second inequality leads to the uniform integrability through the following. 
	\begin{align*}
	\int_E|u(x)|^p\,\d \mu(x)&= \int_E |u(x)|^p\mathds{1}_{\{|u|\leq n\}}(x)\,\d \mu(x)+ \int_E|u(x)|^p\mathds{1}_{\{|u|>n\}}(x)\,\d \mu(x)\\
	&\leq n^p\mu(E)+ \int_X|u(x)|^p\mathds{1}_{\{|u|>n\}}(x)\,\d \mu(x)<\varepsilon. 
	\end{align*}
	Correspondingly, upon the above observation any finite subset $\{u_1,\cdots,u_N\}$ of $L^p(X)$ is also $p$-tight and $p$-uniformly integrable. 
	
\end{remark}

\bigskip

\noindent We now are in a position to state the Vitali convergence theorem. 
\begin{theorem}[Vitali convergence theorem]\label{thm:vitali}
	Let $u$ and $(u_n)_n$ be a sequence of measurable functions. Then $(u_n)_n$ converges to $u$ in $L^p(X)$ with $1\leq p<\infty$ if and only if $(u_n)_{n}$ is $p$-uniformly integrable, $p$-tight and converges in measure to $u$. 
	
\end{theorem}

\medskip

\begin{proof}
	Assume $(u_n)_n$ converges to $u$ in $L^p(X)$ then the convergence in measure follows from to Proposition \ref{prop:convergence-inmeasure}. Note that for every measurable set $A\subset X$ we have 
	\begin{align*}
	\int_A |u_n(x)|^p\,\d \mu(x)\leq 2^{p-1}\int_A |u(x)|^p\,\d \mu(x)+2^{p-1}\int_X |u_n(x)-u(x)|^p\,\d \mu(x). 
	\end{align*}
	\noindent Given that for each $N\geq 1$ the set $\{u,u_1, \cdots, u_N\}$ is $p-$uniformly integrable and $p$-tight as a finite subset of $L^p(X)$ so is the sequence $(u_n)_n$ since $\|u_n-u\|_{L^p(x)}\xrightarrow{n\to \infty}0$. 
	
	\noindent Let us now show that the converse holds true. Again from Proposition \ref{prop:convergence-inmeasure} we can assume that $(u_n)_n$ converges almost everywhere to $u$ on $X$. Now for fixed $\varepsilon>0$ the $p$-uniform integrability and the $p$-tightness imply that there exist  $\delta>0$ and a measurable set $E\subset X$ with $\mu(E)<\infty$ both depending solely on $\varepsilon$ such that 
	\begin{align*}
	&\sup_{n\geq 1}\int_{X\setminus E}|u_n(x)|^p\,\d \mu(x)<\varepsilon \quad\textrm{and}\quad \sup_{n\geq 1}\int_{A}|u_n(x)|^p\,\d \mu(x)<\varepsilon \quad\textrm{for all}\quad \mu(A)<\delta.
	\end{align*}
	\noindent Given the pointwise convergence, with the aid of the Fatou's lemma the above estimates respectively imply 
	\begin{align*}
	&\int_{X\setminus E}|u(x)|^p\,\d \mu(x)\leq \liminf_{n\to\infty}\int_{X\setminus E}|u_n(x)|^p\,\d \mu(x)<\varepsilon 
	\intertext{and}
	&\int_{A}|u(x)|^p\,\d \mu(x)\leq \liminf_{n\geq 1}\int_{A}|u_n(x)|^p\,\d \mu(x)<\varepsilon \quad\textrm{for all}\quad \mu(A)<\delta.
	\end{align*}
	Combining these estimates with the previous ones gives 
	\begin{align}
	\hspace{-1ex}\sup_{n\geq 1}\int_{X\setminus E}\hspace{-3ex}|u_n(x)-u(x)|^p\,\d \mu(x)< 2^p\varepsilon ~~\textrm{and}~~ \sup_{n\geq 1}\int_{A}|u_n(x)-u(x)|^p\,\d \mu(x)<2^p\varepsilon, ~\forall~ \mu(A)<\delta\label{eq:uniform-integrability-tight-eps}.
	\end{align}
	Since $\mu(E)<\infty$, for the choice of $\delta>0$ as above the Ergorov's Theorem \ref{thm:Ergorov} reveals that there is a measurable set $A\subset E$ with $\mu(A)<\delta$ such that $(u_n)_n$ converges uniformly to $u$ on $E\setminus A$, i.e. $$\sup\limits_{x\in E\setminus A}|u_n(x)-u(x)|^p\xrightarrow{n\to \infty}0.$$ For the particular choices of $E$ and $A$ with $\mu(A)<\delta$ in combination with the estimates in \eqref{eq:uniform-integrability-tight-eps} we get
	\begin{align*}
	\int_{X}|u_n(x)-u(x)|^p\,\d \mu(x)&\leq \int_{X\setminus E}|u_n(x)-u(x)|^p\,\d \mu(x)+ \int_{A}|u_n(x)-u(x)|^p\,\d \mu(x)\\
	&+ \int_{E\setminus A}|u_n(x)-u(x)|^p\,\d \mu(x)\\
	&\leq 2^p\varepsilon +2^p\varepsilon + \mu(E\setminus A)\sup_{x\in E\setminus A}|u_n(x)-u(x)|^p.
	\end{align*}
	\noindent The expected convergence in $L^p(X)$ thus follows by letting $n\to \infty$ and $\varepsilon\to 0$ successively.
\end{proof}

\bigskip

\begin{remark}
	Note that the tightness cannot be completely dropped. Indeed on $\mathbb{R}$, if we consider $u_n= \mathds{1}_{[n,n+1]}$ then $(u_n)_n$ converges pointwise everywhere to $u=0$ and is uniformly integrable since for each $n\geq 1$ and every measurable set $E$ with $|E|<\varepsilon$, 
	$$\int_{E}u_n dx = |E\cap [n,n+1]|\leq |E|<\varepsilon.$$
	\noindent Nevertheless, 
	\begin{align*}
	1= \lim_{n \to \infty}\int_{X}u_n(x)\,\d x \neq\int_{X} \lim_{n \to \infty}u_n(x) \,\d x=0. 
	\end{align*}
\end{remark}

\bigskip

\noindent We now visit some consequences of  Vitali's convergence theorem. 
\begin{corollary}\label{cor:consequence-vitali1}
	Let $(u_n)_n$ be a sequence of $L^p(X)$ with $(1\leq p<\infty)$ converging almost everywhere to $u\in L^p(X)$. Then $(u_n)_n$ converges to $u$ in $L^p(X)$ if and only if $(|u_n|^p)_n$ converges to $|u|^p$ in $L^1(X)$.
\end{corollary}
\begin{proof}
	The result plainly springs from Theorem \ref{thm:vitali} based on the observation that $(u_n)_n$ is $p$-uniformly integrable and $p$-tight if and only if $(|u_n|^p)_n$ is uniformly integrable and tight. 
\end{proof}

\medskip

\noindent A deeper version of the above result is the following. 
\begin{corollary}
	Let $(u_n)_n$ be a sequence of $L^p(X)$ with $(1\leq p<\infty)$ converging almost everywhere to $u\in L^p(X)$. Then $(u_n)_n$ converges to $u$ in $L^p(X)$ if and only if $\|u_n\|_{L^p(X)} \xrightarrow{n \to \infty}\|u\|_{L^p(X)}$. 
\end{corollary}

\medskip

\begin{proof} The forward implication is obvious hence let us prove the converse statement.
	Assume $\|u_n\|_{L^p(X)} \to\|u\|_{L^p(X)}$ as $ {n \to \infty}$. Put $g_n =\max(0, |u|^p-|u_n|^p)$. We have $g_n\xrightarrow{n\to \infty}0$ pointwise almost everywhere, $(g_n)_n$ uniformly integrable and tight since $0\leq g_n\leq |u|^p$ for each $n\geq 1$ and $|u|^p\in L^1(X)$. Wherefore, Theorem \ref{thm:vitali} implies 
	\begin{align*}
	\int_X \max(0, |u|^p-|u_n|^p)\d \mu(x) =\int_Xg_n(x)\d \mu(x)\xrightarrow{n \to \infty}0. 
	\end{align*}{}
	\noindent Taking into account the identity $ |u_n|^p-|u|^p= h_n-g_n$ with $h_n=\max(0, |u_n|^p-|u|^p)$ and $ g_n= \max(0, |u|^p-|u_n|^p)$, the assumption also entails 
	\begin{align*}
	\lim_{n \to \infty} \int_X h_n(x)\d \mu(x)= \lim_{n \to \infty} \big(\|u_n\|^p_{L^p(X)} -\|u\|^p_{L^p(X)}\big) - \lim_{n \to \infty} \int_Xg_n(x)\d \mu(x) =0. 
	\end{align*}
	Finally, since $| |u_n|^p-|u|^p|=h_n+g_n= \max(0, |u_n|^p-|u|^p)+ \max(0, |u|^p-|u_n|^p)$ we get
	\begin{align*}
	\lim_{n \to \infty} \int_X ||u_n(x)|^p-|u|^p(x)|\,\d \mu(x)=0. 
	\end{align*}
	That is $(|u_n|^p)_n$ converges to $|u|^p$ in $L^1(X)$ thus by Corollary \ref{cor:consequence-vitali1} $(u_n)_n$ converges to $u$ in $L^p(X)$. 
	
	\medskip
	
	\noindent Alternatively one may consider  the sequence of positive functions $(2^{p-1}(|u_n|^p+|u|^p) -|u_n-u|^p)_n$ which converges almost everywhere to $2^p|u|$. Therefore, in view of Fatou's lemma and the assumption we get
	\begin{align*}
	2^{p} \|u\|^p_{L^p(X)} &\leq \liminf_{n\to\infty} \int_X\big[2^{p-1}(|u_n|^p+|u|^p) -|u_n-u|^p\big]\,\d \mu(x)\\
	&= 2^{p} \|u\|^p_{L^p(X)} -\limsup_{n\to \infty}\int_X|u_n-u|^p\big)\,\d \mu(x). 
	\end{align*}
	The result follows since we have shown that $\limsup\limits_{n\to \infty}\|u_n-u\|_{L^p(X)} =0$.
\end{proof}

\medskip

\noindent The special case $p=1$ provides the following well-known and established Scheff\'e lemma. 
\begin{corollary}[Scheff\'e Lemma, \cite{Wil91}, p.55]\label{cor:Scheffe-lemma}
	Let $(u_n)_n$ be a sequence of $L^1(X)$ converging almost everywhere to $u\in L^1(X)$. Then $(u_n)_n$ converges to $u$ in $L^1(X)$ if and only if $$\int_X|u_n(x)|\,\d \mu(x) \xrightarrow{n \to \infty}\int_X|u(x)|\,\d \mu(x).$$ 
\end{corollary}

\noindent A typical application of Vitali’s theorem is provided by the next result. 
\begin{corollary}
	Assume $\mu(X)<\infty$ and $1<p<\infty$. Let $(u_n)_n$ be a bounded sequence of $L^p(X)$ 
	converging almost everywhere to $u$. Then $(u_n)_n$ converges to $u$ in $L^r(X)$ for all $1\leq r<p$. 
\end{corollary}

\begin{proof}
	The $r$-tightness $(u_n)_n$ obviously holds true and the $r$-uniform integrability follows from H\"{o}lder inequality 
	\begin{align*}
	\sup_{n\geq 1} \int_E |u_n(x)|^r\,\d \mu(x)\leq \mu(E)^{1-\frac{r}{p}} \sup_{n\geq 1} \Big(\int_X |u_n(x)|^p\,\d \mu(x)\Big)^{r/p}\leq C \mu(E)^{1-\frac{r}{p}}\xrightarrow{\mu(E)\to 0}0.
	\end{align*}
	
\end{proof}

\medskip

\noindent The celebrated Lebesgue dominated convergence theorem appears as an immediate consequence of Vitali's convergence theorem. 
\begin{theorem}[Lebesgue dominated convergence theorem]\label{thm:lebesgue-dominated}
	Let $(u_n)_n$ be a bounded sequence of $L^p(X)$ with $(1\leq p<\infty)$ converging almost everywhere to $u$. Assume that  there exists $g\in L^p(X)$ such that 
	$$\sup_{n\geq 1}|u_n(x)|\leq g(x) \qquad\textrm{for almost every }\quad x\in X $$ then $u\in L^p(X)$ and 
	\begin{align*}
	\lim_{n\to \infty}\int_X|u_n(x)-u(x)|^p\,\d \mu(x)=0.
	\end{align*}
\end{theorem}

\medskip

\begin{proof}
	The function $g\in L^p(X)$ is $p$-uniformly integrable and $p$-tight so is $(u_n)_n$ because $|u_n|\leq g$ for all $n\geq 1$ and thus the result follows from Theorem \ref{thm:vitali}. 
\end{proof}

\noindent Another consequence of the Vitali theorem is Brezis-Lieb lemma \cite[Section 4.5]{Bre10} is the following. 
\begin{theorem}[Brezis-Lieb lemma] 
Let $(u_n)_n$ be a bounded sequence of $L^p(X)$ with $(1<p<\infty)$ converging almost everywhere to $u$.
Then 
\begin{align*}
\|u\|^p_{L^p(X)} = \lim_{n\to \infty}\big\{ \|u_n\|^p_{L^p(X)} -\|u_n-u\|^p_{L^p(X)}\big\}.
\end{align*}
In particular, if $\|u_n\|_{L^p(X)}\xrightarrow{n\to\infty} \|u\|_{L^p(X)}$ then $\|u_n-u\|_{L^p(X)}\xrightarrow{n\to\infty}0.$
\end{theorem}

\begin{proof}
  Note that boundedness plus almost everywhere convergence imply the weak convergence. Then $|u_n-u|\rightharpoonup0$ (weakly) in $L^p(X)$ and  $|u_n-u|^{p-1} \rightharpoonup0$ (weakly) in $L^{p'}(X)$. One concludes by noticing that There exists a constant $C>0$ such that for all $a,b\in \R$ (taking $a= u_n-u$, $b=u$),
  \begin{align*}
  	||a+b|-|a|-|b||\leq C(|b|^{p-1}a+ |a|^{p-1}b).
  \end{align*}
 
\end{proof}
%
%
%
%
%
%
%
%
%

\addcontentsline{toc}{chapter}{Bibliography}

\newcommand{\etalchar}[1]{$^{#1}$}

\end{document}